\documentclass[11pt]{amsart}
\usepackage[a4paper,margin=25mm]{geometry}
\usepackage{amsmath,amssymb,amsthm,mathtools,mathrsfs}
\usepackage{booktabs,longtable,array,enumitem}
\usepackage{needspace,placeins}
\usepackage[T1]{fontenc}
\usepackage{lmodern,microtype}
\usepackage{xcolor}
\usepackage{tikz}
\usetikzlibrary{arrows.meta,positioning,calc}
\usepackage[colorlinks=true,linkcolor=blue!45!black,citecolor=blue!45!black,urlcolor=blue!45!black]{hyperref}
\usepackage{orcidlink}
\hypersetup{pdftitle={$G$-crossed braided tensor categories of twisted modules},pdfauthor={Shun Xu}}
\AtBeginDocument{\setlength{\abovedisplayskip}{8pt plus 2pt minus 3pt}\setlength{\belowdisplayskip}{8pt plus 2pt minus 3pt}\setlength{\abovedisplayshortskip}{4pt plus 2pt minus 2pt}\setlength{\belowdisplayshortskip}{6pt plus 2pt minus 2pt}}
\setlist[enumerate]{label=\textup{(\roman*)},leftmargin=*,itemsep=3pt}
\setlist[itemize]{leftmargin=*,itemsep=3pt}
\newcolumntype{P}[1]{>{\raggedright\arraybackslash}p{#1}}
\numberwithin{equation}{section}
\newtheorem{theorem}{Theorem}[section]
\newtheorem{proposition}[theorem]{Proposition}
\newtheorem{lemma}[theorem]{Lemma}
\newtheorem{corollary}[theorem]{Corollary}
\theoremstyle{definition}
\newtheorem{definition}[theorem]{Definition}
\theoremstyle{remark}

\newcommand{\C}{\mathbb C}
\newcommand{\R}{\mathbb R}
\newcommand{\Z}{\mathbb Z}
\newcommand{\N}{\mathbb N}
\newcommand{\cO}{\mathcal O}
\newcommand{\cC}{\mathcal C}
\newcommand{\cT}{\mathcal T}
\newcommand{\cJ}{\mathcal J}
\newcommand{\cQ}{\mathcal Q}
\newcommand{\cF}{\mathcal F}
\newcommand{\cI}{\mathcal I}

\newcommand{\bo}{\boxtimes}
\newcommand{\wt}{\operatorname{wt}}
\newcommand{\Res}{\operatorname{Res}}
\newcommand{\Hom}{\operatorname{Hom}}
\newcommand{\End}{\operatorname{End}}
\newcommand{\Irr}{\operatorname{Irr}}
\newcommand{\Span}{\operatorname{span}}

\newcommand{\Aut}{\operatorname{Aut}}
\newcommand{\Fitt}{\operatorname{Fitt}}
\newcommand{\id}{\operatorname{id}}
\newcommand{\ee}{\varepsilon}
\newcommand{\vval}{\operatorname{v}}
\newcommand{\scrD}{\mathscr D}

\newcommand{\diff}{\mathrm{diff}}
\newcommand{\can}{\mathrm{can}}
\newcommand{\PBW}{\mathrm{PBW}}

\title[$G$-crossed tensor categories]{$G$-crossed braided tensor categories of twisted modules}
\author[Shun Xu]{Shun Xu~\orcidlink{0009-0006-8080-8107}}
\address{School of Mathematical Sciences, Anhui University, Hefei 230601, Anhui, China}
\email{shunxu@ahu.edu.cn}
\date{}

\begin{document}
\begin{abstract}
Let $V$ be a strongly rational vertex operator algebra and $G\leq\Aut(V)$ a finite group. We prove Huang's conjecture that the categories of $g$-twisted $V$-modules, for $g\in G$, form a $G$-crossed braided fusion category. This category is faithfully $G$-graded, with tensor products representing the complete twisted intertwining operators of Du--Huang, the canonical pullback $G$-action, and crossed braiding. In particular, the associativity, coherence and rigidity required by the conjecture all hold.
\end{abstract}
\maketitle

\section{Introduction}
\subsection{Twisted modules and Huang's conjecture}
The representation theory of a vertex operator algebra governs the chiral sectors of a conformal field theory. For a finite automorphism group, orbifold theory relates three kinds of objects: modules for the original vertex operator algebra, modules for the fixed-point algebra, and twisted modules attached to group elements. The orbifold computations of Dijkgraaf, Vafa, Verlinde and Verlinde already show that the twisted sectors are not organized by an ordinary braided tensor category: the sector label changes when insertion points are exchanged \cite{DVVV}. The expected structure is therefore a $G$-graded tensor category equipped with a $G$-action and a crossed braiding.

The algebraic foundations of twisted modules were established by Dong, Li and Mason. They developed the finite-order twisted-module theory, introduced twisted Zhu algebras, and proved the correspondence between simple modules for those algebras and simple admissible twisted modules \cite{DLM}. Their modular invariance theorem for twisted trace functions provided further finiteness results under rationality and cofiniteness hypotheses \cite{DLMModular}. These works supply the basic representation-theoretic framework for orbifold sectors.

The tensor theory for untwisted modules was developed by Huang and Lepowsky \cite{HL}, and later extended by Huang, Lepowsky and Zhang to suitable logarithmic settings \cite{HLZ}. In this theory tensor products represent intertwining operators at prescribed insertion points, while analytic continuation produces the associativity and braiding isomorphisms. Huang's differential-equation and regular-singularity results provide the convergence and extension properties needed for the construction \cite{HuangDifferential}, and his rigidity theorem yields a modular tensor category under strong rationality hypotheses \cite{HuangRigidity}.

On the orbifold side, Kirillov described $G$-equivariant modular categories and pointed to twisted modules as a fundamental example \cite[Example~5.5]{KirillovG}. He also asked for a direct construction of the monoidal structure in terms of twisted modules themselves \cite[Remark~5.6]{KirillovG}. Huang formulated the corresponding finite-group conjecture \cite[Conjecture~6.2]{HuangProblems}: if $V$ is strongly rational and $G\leq\Aut(V)$ is finite, then the twisted sectors should form a $G$-crossed tensor category.

The fixed-point approach starts from a suitable category of $V^G$-modules. Carnahan and Miyamoto proved regularity for fixed points of a finite-order automorphism, and consequently for finite solvable groups \cite{CarnahanMiyamoto}. Given an ambient vertex tensor category of $V^G$-modules containing $V$, McRae constructed the crossed category and proved that its equivariantization recovers that ambient tensor category \cite{McRae}.

Tan established differential equations and convergence for particular twisted-intertwiner configurations with an untwisted first input and a fixed twisted second input \cite{TanDifferential}.

The present paper follows a direct twisted-module route. Du and Huang introduced twisted intertwining operators, proved their skew-symmetry and contragredient isomorphisms, and constructed universal binary tensor products together with crossed commutativity maps under suitable hypotheses \cite{DH}. They restated the finite-group conjecture in this analytic language \cite[Conjecture~1.2]{DH}. Passing from their binary constructions to a $G$-crossed braided fusion category requires four further steps: finite-dimensional control of the relevant Ward quotients, regular singularity near collisions, reconstruction of the intermediate twisted modules from collision data, and coherence of the resulting associativity and crossed braiding isomorphisms.

\subsection{The main result}
Recall that $V$ is \emph{strongly rational} if it is simple, of CFT type, self-contragredient, rational and $C_2$-cofinite. Here CFT type means $V=\bigoplus_{n\ge 0}V_n$ with $V_0=\C\mathbf 1$, rationality means complete reducibility of admissible untwisted modules, and $C_2$-cofiniteness means
\[
 \dim_\C V/\Span\{u_{-2}v:u,v\in V\}<\infty.
\]
For $g\in G$, let $\cC_g$ be the category of grading-restricted generalized $g$-twisted $V$-modules, and set $\cC_G:=\bigoplus_{g\in G}\cC_g$. The grading restriction means that the generalized $L(0)$-eigenspaces are finite-dimensional and their real parts are bounded below. We write $\bo$ for the tensor product to be constructed. The categorical conventions for $G$-crossed braided categories are recalled in Section~\ref{sec:categoricalconventions}, and the Du--Huang complete twisted intertwining operators are given in Definition~\ref{def:completeintertwiner}.

\begin{theorem}\label{thm:main}
Let $V$ be strongly rational and let $G\leq\Aut(V)$ be finite. The category $\cC_G$ has a natural structure of a faithfully $G$-graded, $G$-crossed braided fusion category. Each sector $\cC_g$ is nonzero. Its tensor product sends $\cC_g\times\cC_h$ to $\cC_{gh}$ and represents the complete twisted intertwining operators of Definition~\ref{def:completeintertwiner}. The unit is $V$, the action is the strict pullback $Y_{T_kX}(v,x)=Y_X(k^{-1}v,x),\ T_kT_l=T_{kl},$ and for $X\in\cC_g$ the crossed braiding is $c_{X,Y}:X\bo Y\longrightarrow T_gY\bo X.$ The associativity isomorphisms have the defining product--iterate identity on their original overlap, and the braiding has the Du--Huang upper-half-plane normalization. The pentagon, triangle, both crossed hexagons and action equivariance hold. The category is rigid; a dual of a $g$-twisted object lies in the $g^{-1}$-sector.
\end{theorem}

The theorem establishes Huang's finite-group twisted-module tensor-category conjecture \cite[Conjecture~6.2]{HuangProblems} directly in the category of twisted modules, with tensor products representing the complete class of Du--Huang intertwining operators.

\subsection{Outline of the proof}
The construction begins with the finite-order representation theory of $V$. Strong rationality implies $g$-regularity for each $g\in G$, so the categories $\cC_g$ are finite semisimple and their simple objects are ordinary twisted modules (Proposition~\ref{prop:regularmodules}). This permits the use of restricted duals and genuine homogeneous intermediate-weight sums. It also reduces the construction of tensor products and the verification of identities to finitely many simple types, without presupposing an associativity theorem for products of twisted intertwining operators.

The first analytic input is a finite differential Ward quotient. We realize twisted mode actions as residues of equivariant currents on a based, framed Hurwitz cover, retaining the local roots and the action on the output restricted dual. Finite-state current generation, cyclic spanning, and control of principal parts make the seed Ward quotient finite over the configuration space (Proposition~\ref{prop:seedfinite}). Saturating its relations under the placement derivatives gives a finite-rank vector bundle with an integrable connection (Theorem~\ref{thm:finiteconn}). Original three-point intertwining-operator functionals satisfy these Ward identities. Their factorization through the finite quotient bounds fusion rules (Proposition~\ref{prop:fusion}), and semisimplicity then gives the universal $P(z)$-tensor products (Theorem~\ref{thm:binary}). The units and the crossed commutors are fixed by the module vertex operators and the Du--Huang skew transformations, with their stated logarithm branches.

The next step concerns analytic continuation of the original, rather than merely formal, multi-input channels. We extend the prescribed framed covers across collision divisors by balanced admissible covers and compare the relative-weight integral current lattices. This establishes regular singularity of the canonical Ward connection (Theorem~\ref{thm:canonicalregular}); its logarithmic extension controls the local exponents and gives convergent expansions (Theorem~\ref{thm:formalconv}). For actual products and iterates we separately identify the four original three-point sources, their inherited roots, the output restricted dual, and the common comparison paths (Lemma~\ref{lem:E0actualsourcecomparison}). The corresponding full-state Ward identities and their differential saturation (Theorem~\ref{thm:fullsmallward}), together with the covariant continuation of the original E0 functionals (Theorem~\ref{thm:E0covariantcontinuation}), connect the finite connection to normally convergent original weight sums. Thus a formal horizontal solution is not substituted for the originally specified product or iterate.

Collision expansions yield coefficient kernels with a common discrete lower support (Proposition~\ref{prop:collisionuniformsupport}). We propagate $V$-insertions on the normalized components and use the complete residue identities, including the output contragredient action, to construct the intermediate weak twisted modules (Theorem~\ref{thm:partialnodeweakmodule}). Finite-dimensionality and twisted regularity identify their ordinary summands and recover complete three-point intertwining operators at the node (Proposition~\ref{prop:partialnodevertices}). The remaining comparison is coefficientwise: balanced global sections realize prescribed node jets (Lemma~\ref{lem:matchingnodejetlift}), and the resulting mode identities identify the two node actions on each collision coefficient (Theorem~\ref{thm:matchedcoefficientmodes}). These identifications give factorization by actual iterates with the original physical-coordinate normalizations (Theorem~\ref{thm:partialnodeexactspan}), rather than only a formal decomposition of the collision kernel.

For three inputs, exact factorization and the binary universal property give a natural epimorphism between the two bracketings (Lemma~\ref{lem:productepi}). The finite-orbit argument, using the crossed commutors and the finite grading, upgrades it to an isomorphism (Lemma~\ref{lem:finiterotation} and Theorem~\ref{thm:smallassoc}). The product--iterate equality on its positive overlap determines the associator, including its parallel-transport normalization. The unit triangle and the strong monoidal $G$-action are checked against the original universal operators (Proposition~\ref{prop:unittriangle} and Proposition~\ref{prop:action}). For the pentagon, the five four-input channels are compared on one based cover and one common analytic domain; normal convergence permits the original two-grade sums to be compared before their equality is transferred to module morphisms (Theorem~\ref{thm:criticalpentagon}). For the two crossed hexagons, the actual exchange paths, full translation covariance, root phases and prescribed skew normalization identify the respective routes (Theorem~\ref{thm:criticalhexagons}).

Finally, restricted contragredients supply nonzero pairings compatible with the $G$-action (Lemma~\ref{lem:contragredientpairing}). The equivariantization $\cC_G^G$ is finite semisimple and braided, with a strong monoidal forgetful functor; the pairings give the required non-negligibility of its simple objects (Lemma~\ref{lem:equivariantization}). Applying the Etingof--Penneys rigidity criterion there and transferring duals through the forgetful functor proves rigidity of $\cC_G$ (Proposition~\ref{prop:rigidity}). Together with the preceding coherence results and nonvanishing of each sector, this proves Theorem~\ref{thm:main}.

\subsection{Proof dependencies}
The first diagram displays the analytic results leading to exact iterate factorization, and the second records their use in the construction of the crossed braided structure and in the rigidity argument. In particular, regular singularity, the full-state Ward identities for the original E0 channels, and coefficientwise node matching enter separately into the factorization theorem.

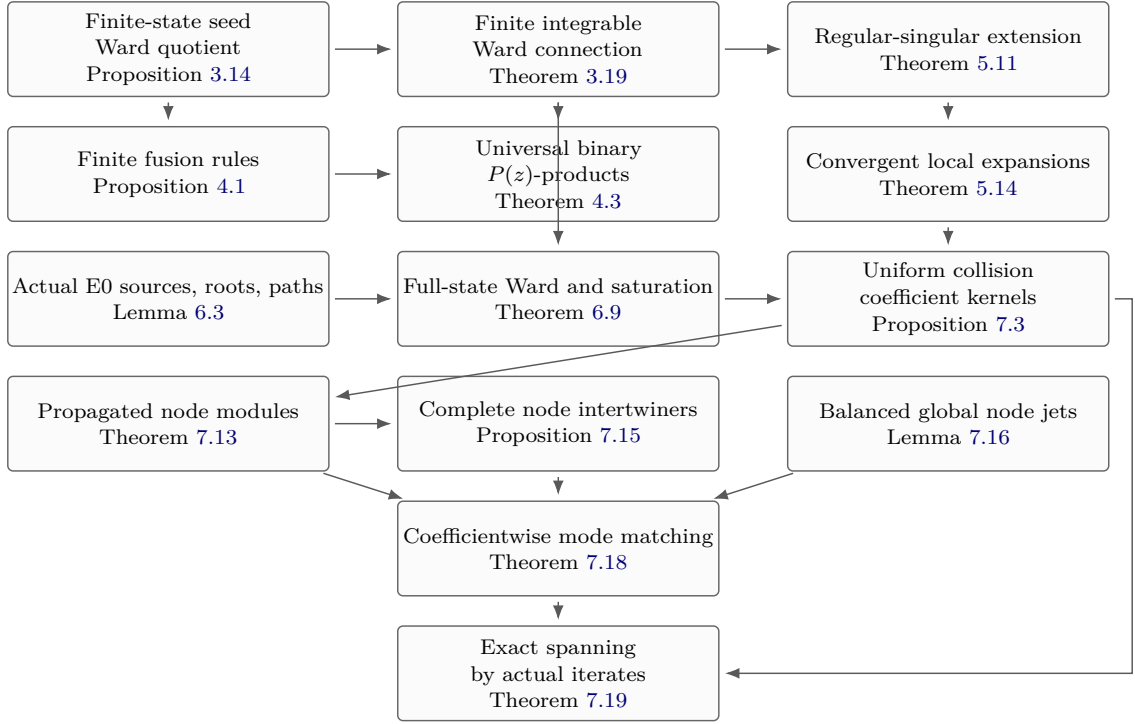
\begin{figure}[!htbp]
\centering
\begin{tikzpicture}[
  >=Latex, line width=.55pt,
  box/.style={draw=black!60, rounded corners=2pt, fill=gray!3,
    align=center, text width=41mm, minimum height=12.5mm, inner sep=2pt,
    font=\scriptsize},
  arr/.style={->,draw=black!65,shorten <=2pt,shorten >=2pt}]
\node[box] (seed) at (0,0) {Finite-state seed Ward quotient\\Proposition~\ref{prop:seedfinite}};
\node[box] (conn) at (5.15,0) {Finite integrable Ward connection\\Theorem~\ref{thm:finiteconn}};
\node[box] (reg) at (10.30,0) {Regular-singular extension\\Theorem~\ref{thm:canonicalregular}};
\node[box] (fusion) at (0,-1.65) {Finite fusion rules\\Proposition~\ref{prop:fusion}};
\node[box] (binary) at (5.15,-1.65) {Universal binary $P(z)$-products\\Theorem~\ref{thm:binary}};
\node[box] (conv) at (10.30,-1.65) {Convergent local expansions\\Theorem~\ref{thm:formalconv}};
\node[box] (e0) at (0,-3.30) {Actual E0 sources, roots, paths\\Lemma~\ref{lem:E0actualsourcecomparison}};
\node[box] (ward) at (5.15,-3.30) {Full-state Ward and saturation\\Theorem~\ref{thm:fullsmallward}};
\node[box] (kern) at (10.30,-3.30) {Uniform collision coefficient kernels\\Proposition~\ref{prop:collisionuniformsupport}};
\node[box] (node) at (0,-4.95) {Propagated node modules\\Theorem~\ref{thm:partialnodeweakmodule}};
\node[box] (vertex) at (5.15,-4.95) {Complete node intertwiners\\Proposition~\ref{prop:partialnodevertices}};
\node[box] (jet) at (10.30,-4.95) {Balanced global node jets\\Lemma~\ref{lem:matchingnodejetlift}};
\node[box] (match) at (5.15,-6.60) {Coefficientwise mode matching\\Theorem~\ref{thm:matchedcoefficientmodes}};
\node[box] (exact) at (5.15,-8.25) {Exact spanning by actual iterates\\Theorem~\ref{thm:partialnodeexactspan}};
\draw[arr] (seed)--(conn);
\draw[arr] (conn)--(reg);
\draw[arr] (seed)--(fusion);
\draw[arr] (fusion)--(binary);
\draw[arr] (conn)--(binary);
\draw[arr] (reg)--(conv);
\draw[arr] (e0)--(ward);
\draw[arr] (conn)--(ward);
\draw[arr] (conv)--(kern);
\draw[arr] (ward)--(kern);
\draw[arr] (kern)--(node);
\draw[arr] (node)--(vertex);
\draw[arr] (node)--(match);
\draw[arr] (vertex)--(match);
\draw[arr] (jet)--(match);
\draw[arr] (match)--(exact);
\draw[arr] (kern.east) -- ++(3mm,0) |- (exact.east);
\end{tikzpicture}
\caption{Analytic control and node reconstruction. The original channel Ward identities, the regular-singular expansion, and the balanced node-mode comparison jointly support exact iterate factorization.}
\label{fig:analyticdependencies}
\end{figure}

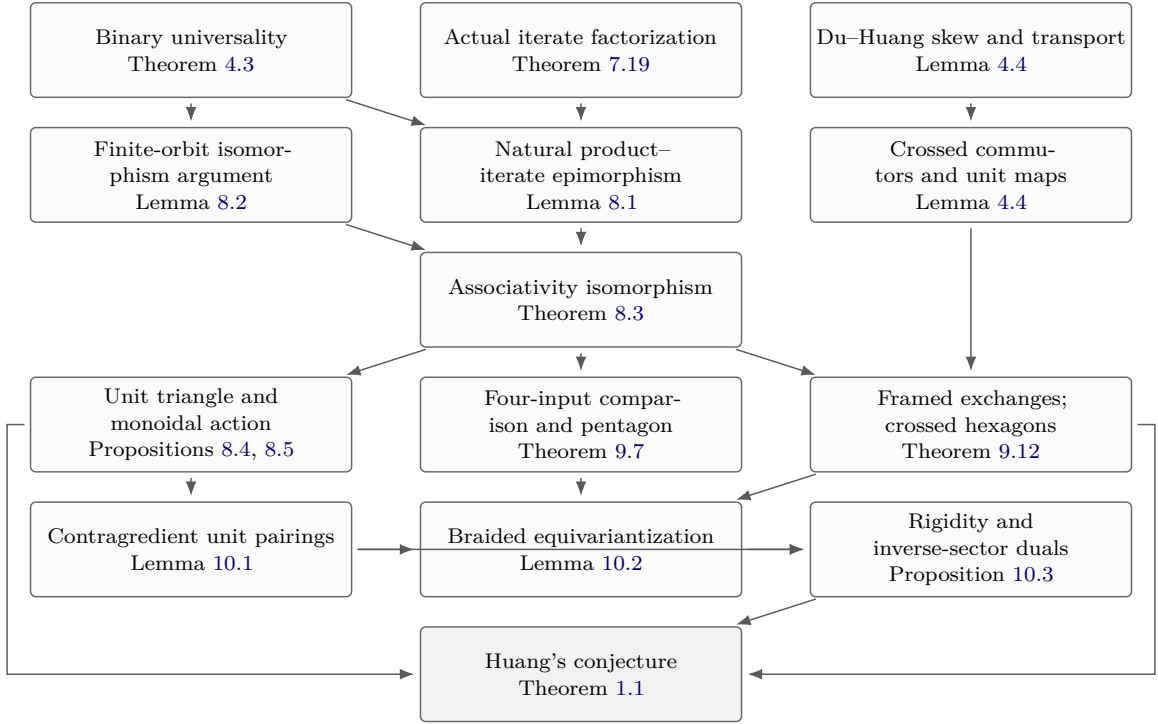
\begin{figure}[!htbp]
\centering
\begin{tikzpicture}[
  >=Latex,line width=.55pt,
  box/.style={draw=black!60,rounded corners=2pt,fill=gray!3,
    align=center,text width=41mm,minimum height=12.5mm,inner sep=2pt,font=\scriptsize},
  conclusion/.style={box,fill=gray!10},
  arr/.style={->,draw=black!65,shorten <=2pt,shorten >=2pt}]
\node[box] (binary) at (0,0) {Binary universality\\Theorem~\ref{thm:binary}};
\node[box] (exact) at (5.15,0) {Actual iterate factorization\\Theorem~\ref{thm:partialnodeexactspan}};
\node[box] (skew) at (10.30,0) {Du--Huang skew and transport\\Lemma~\ref{lem:units}};
\node[box] (rotation) at (0,-1.65) {Finite-orbit isomorphism argument\\Lemma~\ref{lem:finiterotation}};
\node[box] (epi) at (5.15,-1.65) {Natural product--iterate epimorphism\\Lemma~\ref{lem:productepi}};
\node[box] (comm) at (10.30,-1.65) {Crossed commutors and unit maps\\Lemma~\ref{lem:units}};
\node[box] (assoc) at (5.15,-3.30) {Associativity isomorphism\\Theorem~\ref{thm:smallassoc}};
\node[box] (unit) at (0,-4.95) {Unit triangle and monoidal action\\Propositions~\ref{prop:unittriangle}, \ref{prop:action}};
\node[box] (penta) at (5.15,-4.95) {Four-input comparison and pentagon\\Theorem~\ref{thm:criticalpentagon}};
\node[box] (hexa) at (10.30,-4.95) {Framed exchanges; crossed hexagons\\Theorem~\ref{thm:criticalhexagons}};
\node[box] (pair) at (0,-6.60) {Contragredient unit pairings\\Lemma~\ref{lem:contragredientpairing}};
\node[box] (eqvt) at (5.15,-6.60) {Braided equivariantization\\Lemma~\ref{lem:equivariantization}};
\node[box] (rigid) at (10.30,-6.60) {Rigidity and inverse-sector duals\\Proposition~\ref{prop:rigidity}};
\node[conclusion] (main) at (5.15,-8.25) {Huang's conjecture\\Theorem~\ref{thm:main}};
\draw[arr] (binary)--(rotation);
\draw[arr] (exact)--(epi);
\draw[arr] (binary)--(epi);
\draw[arr] (skew)--(comm);
\draw[arr] (rotation)--(assoc);
\draw[arr] (epi)--(assoc);
\draw[arr] (assoc)--(unit);
\draw[arr] (assoc)--(penta);
\draw[arr] (assoc)--(hexa);
\draw[arr] (comm)--(hexa);
\draw[arr] (unit)--(pair);
\draw[arr] (pair)--(eqvt);
\draw[arr] (penta)--(eqvt);
\draw[arr] (hexa)--(eqvt);
\draw[arr] (eqvt)--(rigid);
\draw[arr] (pair)--(rigid);
\draw[arr] (rigid)--(main);
\draw[arr] (unit.west) -- ++(-3mm,0) |- (main.west);
\draw[arr] (hexa.east) -- ++(3mm,0) |- (main.east);
\end{tikzpicture}
\caption{From exact factorization to the $G$-crossed braided fusion structure. Coherence is proved using original channel identities before the rigidity argument via equivariantization.}
\label{fig:categoricaldependencies}
\end{figure}

\FloatBarrier
\section{Preliminaries}\label{sec:setup}
\subsection{Vertex operator algebras and twisted modules}\label{sec:moduleprelim}
Put $\N=\Z_{\geq0}$ and $\C^*=\C\setminus\{0\}$. For a formal Laurent series, $\Res_x f(x)=[x^{-1}]f(x)$; set $\delta(x)=\sum_{j\in\Z}x^j$ and $\delta_{i,j}=1$ if $i=j$, $0$ otherwise. We use $\binom{\eta}{0}=1$ and $\binom{\eta}{j}=\frac{1}{j!}\prod_{k=0}^{j-1}(\eta-k)$ for $j\geq1$. The binomial expansion of $(x-y)^\eta$ is in nonnegative powers of $y$.

A vertex operator algebra $(V,Y,\mathbf1,\omega)$ has $V=\bigoplus_{d\in\Z}V_d$, $\dim V_d<\infty$, $V_d=0$ for $d\ll0$, $\mathbf1\in V_0$, $\omega\in V_2$, and
$Y(u,x)=\sum_{n\in\Z}u_nx^{-n-1}\in\End(V)[[x,x^{-1}]]$ with $u_nv=0$ for $n\gg0$. The vacuum, creation and Jacobi identities are $Y(\mathbf1,x)=\id$, $Y(u,x)\mathbf1\in u+xV[[x]]$, and
\begin{align*}
&x_0^{-1}\delta\!\left(\frac{x_1-x_2}{x_0}\right)Y(u,x_1)Y(v,x_2)
-x_0^{-1}\delta\!\left(\frac{x_2-x_1}{-x_0}\right)Y(v,x_2)Y(u,x_1)\\
&\hspace{10mm}=x_2^{-1}\delta\!\left(\frac{x_1-x_0}{x_2}\right)Y(Y(u,x_0)v,x_2).
\end{align*}
For $Y(\omega,x)=\sum_{n\in\Z}L(n)x^{-n-2}$,
\[
[L(m),L(n)]=(m-n)L(m+n)+\frac{c_V}{12}(m^3-m)\delta_{m+n,0}\id,
\qquad L(0)|_{V_d}=d\id,
\]
and $Y(L(-1)u,x)=\partial_xY(u,x)$. Here $c_V\in\C$ and $\wt(u)=d$ for $u\in V_d$. A $V$-submodule is a subspace stable under all $u_n$; $V$ is simple if its adjoint module is simple. The group $\Aut(V)$ preserves $\mathbf1$, $\omega$ and $Y$. The CFT-type condition is $V_0=\C\mathbf1$, $V_d=0$ for $d<0$. Define $C_2(V)=\Span\{u_{-2}v:u,v\in V\}$; $C_2$-cofiniteness is $\dim_\C(V/C_2(V))<\infty$.

Fix $g\in\Aut(V)$ with $N=\operatorname{ord}(g)<\infty$ and put
\[
V^{(r)}=\{u\in V:gu=e^{2\pi ir/N}u\},\qquad
V=\bigoplus_{r=0}^{N-1}V^{(r)}.
\]
A \emph{weak $g$-twisted module} $(M,Y_M^g)$ is a vector space with $Y_M^g(\mathbf1,x)=\id_M$ and, for $u\in V^{(r)}$,
$Y_M^g(u,x)=\sum_{n\in r/N+\Z}u_nx^{-n-1}$, where $u_nw=0$ for $n\gg0$ and each $w\in M$. The twisted Jacobi identity may equivalently be written as the following mode identity \cite[Section~3]{DLM}:
for $u\in V^{(r)}$, $v\in V^{(s)}$, $m\in r/N+\Z$, $n\in s/N+\Z$, $p\in\Z$, and $w\in M$,
\begin{align}\label{eq:twistedBorcherds}
\sum_{j\geq0}\binom mj(u_{p+j}v)_{m+n-j}w
&=\sum_{j\geq0}(-1)^j\binom pj\bigl(u_{m+p-j}v_{n+j}w\notag\\
&\hspace{34mm}-(-1)^pv_{n+p-j}u_{m+j}w\bigr).
\end{align}
On the left, $u_{p+j}v\in V^{(r+s\bmod N)}$ is computed in $V$; every remaining mode acts on $M$. Truncation makes each sum finite on $w$. Write $Y_M^g(\omega,x)=\sum_{j\in\Z}L_M(j)x^{-j-2}$. Then the Virasoro and translation relations hold on $M$, and for $u\in V_d$,
$[L_M(0),u_n]=(d-n-1)u_n$. A submodule of $M$ is stable under all twisted modes; $M$ is irreducible if it has no proper nonzero submodule, and cyclic if $M=\Span\{u^{(1)}_{n_1}\cdots u^{(k)}_{n_k}w:k\geq0\}$ for some $w\in M$. Morphisms commute with all modes.

For $\lambda\in\C$ put $M_{[\lambda]}=\bigcup_{r\geq1}\ker(L_M(0)-\lambda)^r$. A \emph{generalized $g$-twisted module} is a weak $g$-twisted module with $M=\bigoplus_{\lambda\in\C}M_{[\lambda]}$. It is \emph{lower bounded} if $M_{[\lambda]}=0$ for $\Re\lambda<B$ for some $B\in\R$, and \emph{grading restricted} if moreover $\dim M_{[\lambda]}<\infty$ for every $\lambda$. It is \emph{ordinary} if $L_M(0)$ is diagonalizable; then $M_\lambda=M_{[\lambda]}$. In particular, $u_nM_{[\lambda]}\subseteq M_{[\lambda+d-n-1]},\; u\in V_d,\; \wt(M)=\{\lambda:M_{[\lambda]}\ne0\}$.
For a generalized homogeneous vector $0\ne w\in M_{[\lambda]}$ write $\wt(w)=\lambda$. The canonical projections are $\pi^M_\lambda\colon M\to M_{[\lambda]}$. For grading-restricted $M$ set
$M'=\bigoplus_\lambda M_{[\lambda]}^*$ and $\overline M=\prod_\lambda M_{[\lambda]}$, with pairing $\langle w',w\rangle=w'(w)$, extended componentwise to $M'\times\overline M$. A weight-preserving $f\colon M\to N$ induces $\overline f\colon\overline M\to\overline N$. The restricted dual separates $\overline M$.

The contragredient $M'$ has the $g^{-1}$-twisted vertex operator \cite[Section~3]{DLM}
\begin{equation}\label{eq:contragredientprelim}
\bigl\langle Y_{M'}^{g^{-1}}(u,x)w',w\bigr\rangle
=\bigl\langle w',Y_M^g(e^{xL(1)}(-x^{-2})^{L(0)}u,x^{-1})w\bigr\rangle.
\end{equation}
For $u\in V_d$, $(-x^{-2})^{L(0)}u=(-1)^dx^{-2d}u$ and $L(1)^ru=0$ for $r\gg0$. Morphisms are transposed contravariantly, $f'\colon N'\to M'$, and the evaluation map yields $M\simeq M''$. For $V$, self-contragredience is equivalent to a nondegenerate invariant bilinear form $(\,,\,)$ satisfying $(Y(u,x)v,w)=(v,Y(e^{xL(1)}(-x^{-2})^{L(0)}u,x^{-1})w)$.

An \emph{admissible $g$-twisted module} has an additional grading $M=\bigoplus_{j\in N^{-1}\N}M(j)$, $u_nM(j)\subseteq M(j+d-n-1)$ for $u\in V_d$, with $M(t)=0$ for $t<0$. \emph{$g$-rationality} means semisimplicity of this admissible category; \emph{$g$-regularity} means that every weak $g$-twisted module is an algebraic direct sum of ordinary irreducible $g$-twisted modules. Rationality and regularity mean the respective assertions for $g=1$. \emph{Strong rationality} means that $V$ is simple, of CFT type, rational, $C_2$-cofinite and self-contragredient. Henceforth $V$ is strongly rational and $G\leq\Aut(V)$ is finite. Its invariant form is normalized by $(\mathbf1,\mathbf1)=1$.

Let $\cC_g$ be the category of grading-restricted generalized $g$-twisted modules, with
\[
\Hom_V(X,Y)=\{f\in\Hom_\C(X,Y):fY_X^g(u,x)=Y_Y^g(u,x)f\}.
\]
Let $\cC_G=\bigoplus_{g\in G}\cC_g$ denote the external direct sum, with $\Hom_{\cC_G}((g,X),(h,Y))=0$ for $g\ne h$. An object of $\cC_G$ is a finite direct sum of labeled objects $(g,X)$. For $k\in G$ define $T_k(g,X)=(kgk^{-1},\phi_kX),\; Y_{\phi_kX}(u,x)=Y_X^g(k^{-1}u,x)$.
On underlying vector spaces and morphisms $T_k$ is the identity; as functors $T_kT_l=T_{kl}$ and $T_1=\id$ \cite[Section~2]{DH}.

To compare this formal convention with the analytic convention of \cite[Section~2]{DH}, for $z\in\C^*$ and $p\in\Z$ set $\ell_p(z)=\log|z|+i\arg z+2\pi ip$, $0\leq\arg z<2\pi$. For a vector space $U$, write $U\{x\}[\log x]$ for formal series with complex exponents and polynomial logarithmic coefficients, with finite logarithmic degree for every fixed input. A series $F(x)=\sum_{n\in\C}\sum_{j=0}^{K_n}a_{n,j}x^{-n-1}(\log x)^j$ has the $p$-th branch $F^p(z)$ obtained from $x^{-n-1}=e^{(-n-1)\ell_p(z)}$, $\log x=\ell_p(z)$ whenever the branch is defined. For a finite-order $g$, the analytic twisted-module equivariance axiom is $(Y_M^g)^{p+1}(gu,z)=(Y_M^g)^p(u,z),\; p\in\Z,\ z\in\C^*$.
The analytic definition also requires its duality, $L(0)$-grading and $L(-1)$-derivative axioms; equivariance alone does not imply the twisted Jacobi identity.

\begin{lemma}\label{lem:fieldlogs}
If $g^N=1$, every finite-degree logarithmic coefficient in a $g$-twisted vertex operator satisfying the preceding equivariance condition is independent of $\log x$. For $u\in V^{(r)}$, its nonzero modes have indices in $r/N+\Z$.
\end{lemma}
\begin{proof}
Let $gu=e^{2\pi ir/N}u$, and fix $w\in M$ and $n\in\C$. Write the coefficient of $x^{-n-1}$ in $Y_M^g(u,x)w$ as $P_n(\log x)\in M[\log x]$. Equivariance gives $e^{2\pi i(r/N-n-1)}P_n(t+2\pi i)=P_n(t)$.
Since $P_n(t+2\pi i)$ and $P_n(t)$ have the same leading coefficient, $P_n\ne0$ implies $e^{2\pi i(r/N-n-1)}=1$, hence $n\in r/N+\Z$. Now $P_n(t+2\pi i)=P_n(t)$; comparison of the next-to-leading coefficients forces $\deg P_n=0$.
\end{proof}

For grading-restricted analytic generalized twisted modules, the duality axiom of \cite[Section~2]{DH}, together with the finite-order reduction above, gives the formal twisted Jacobi identity \eqref{eq:twistedBorcherds} by the product/iterate coefficient comparison of \cite[Theorem~2.4 and its proof]{HuangTw}. Conversely, the classical twisted Jacobi identity gives the three convergent expansions of the same finite-branched rational matrix coefficients, hence the analytic duality axiom. Indeed, after multiplication by suitable powers of $z_1$, $z_2$ and $z_1-z_2$, the two-point matrix coefficients against $M'$ are finite Laurent sums in $z_1^{1/N},z_2^{1/N}$; expansion in the three chambers gives the required identities. The argument uses finite-support restricted-dual vectors and the generalized weight decomposition, but not semisimplicity of $L_M(0)$. Thus the two notions coincide in the grading-restricted finite-order setting.

Write $\Irr(\cC)$ for the set of isomorphism classes of simple objects of a category $\cC$, and $\operatorname{Ob}(\cC)$ for its objects.

\begin{proposition}\label{prop:regularmodules}
For every $g\in G$, $\cC_g\ne0$, and every weak $g$-twisted module is a direct sum of ordinary irreducible $g$-twisted modules. Moreover, $|\Irr(\cC_g)|<\infty$, and each grading-restricted generalized or cyclic weak $g$-twisted module is a finite direct sum of ordinary irreducibles.
\end{proposition}
\begin{proof}
By \cite[Theorem~4.5]{ABD}, $V$ is regular. The hypotheses of \cite[Corollary~5.26]{CarnahanMiyamoto} hold: $V$ is simple, regular, nonnegatively graded and has a nondegenerate invariant form, while $\operatorname{ord}(g)<\infty$. That corollary gives $g$-regularity, including complete reducibility of \emph{unrestricted weak} $g$-twisted modules. Finiteness of irreducible types follows from $g$-rationality and \cite[Theorem~8.1]{DLM} (cf.\ \cite[Theorem~2.6(2)]{DR}).

Suppose $W=\bigoplus_{i\in I}S_i$ is grading restricted, with every $S_i$ ordinary irreducible. For each isomorphism class $[S]$, choose $\lambda_S\in\wt(S)$ and let $m_S=|\{i\in I:S_i\simeq S\}|$. Then $m_S\dim S_{\lambda_S}\leq\dim W_{[\lambda_S]}<\infty$.
There are finitely many $[S]$, hence $|I|=\sum_{[S]}m_S<\infty$. If $W=\langle w\rangle$ is cyclic, then $w\in\bigoplus_{i\in I_0}S_i$ for a finite subset $I_0\subset I$. Since the latter is a submodule, $W=\bigoplus_{i\in I_0}S_i$. Finally, \cite[Theorem~9.1]{DLMModular} applies because $V$ is simple and $C_2$-cofinite, and supplies a nonzero irreducible ordinary $g$-twisted module; thus $\cC_g\ne0$.
\end{proof}

Choose a set $\mathcal S_g\subset\operatorname{Ob}(\cC_g)$ containing one representative of each class in $\Irr(\cC_g)$ and put $\mathcal S=\bigsqcup_{g\in G}\mathcal S_g$. Let $M$ be any simple ordinary $g$-twisted module and $N_M=\operatorname{ord}(g)$. If $0\ne w\in M_\lambda$, then the submodule it generates has weights in $\lambda+N_M^{-1}\Z$; irreducibility and lower boundedness imply the existence of $h_M\in\lambda+N_M^{-1}\Z$ with
\[
M_{h_M}\ne0,\qquad M=\bigoplus_{a\in N_M^{-1}\N}M_{h_M+a},\qquad h_{M'}=h_M.
\]
For $T=M_0'\otimes M_1\otimes\cdots\otimes M_m$ with $M_i\in\mathcal S$, take $D\in\N_{>0}$ divisible by every $N_{M_i}$. On a homogeneous pure tensor $w_0'\otimes w_1\otimes\cdots\otimes w_m$, with $w_0'\in(M_0')_{\lambda_0}$ and $w_i\in(M_i)_{\lambda_i}$, define the nonnegative integer
\[
\ee(w_0'\otimes w_1\otimes\cdots\otimes w_m)
=D\sum_{i=0}^m(\lambda_i-h_{M_i})\in\N.
\]
Here $\lambda_0$ is the weight of $w_0'$, and $M_0'$ has the same weight spectrum as $M_0$. The quantity $\ee$ is an energy grading of homogeneous tensors, not a linear functional on $T$.

For $M\in\mathcal S_g$, $q\in\C^*$ and a fixed logarithm $\ell(q)$, let $q^{L_M(0)-h_M}=e^{\ell(q)(L_M(0)-h_M)}$. For $b\in V_d\cap V^{(r)}$ and $n\in r/N_M+\Z$, $q^{L_M(0)-h_M}b_nq^{-(L_M(0)-h_M)}=q^{d-n-1}b_n,\; q^a=e^{a\ell(q)}$.
This follows from $b_nM_\lambda\subseteq M_{\lambda+d-n-1}$; for finite direct sums the formula is applied on each simple summand with its own lowest weight.

\subsection{Analytic branches and twisted intertwining operators}
For $r\in\Z_{\geq1}$ put $\mathcal M_r=F_r(\C^*)=\{z=(z_i)_{i=1}^r\in(\C^*)^r:z_i\ne z_j\ (i\ne j)\}$.
A multivalued analytic function on $\mathcal M_r$ is an analytic germ continued along paths in $\mathcal M_r$; a maximal continuation has a germ at each endpoint on every admissible path. For a base point $z^0$, $\pi_1(\mathcal M_r,z^0)$ is the group of endpoint-fixed homotopy classes of loops; continuation induces monodromy on germs. A braid is a homotopy class of collision-free point motions; it is pure if each labeled endpoint equals its starting point. Below $\sigma_i$ is the clockwise half-exchange of the $i$-th and $(i+1)$-st labeled points (a lower half-turn of their difference); products of braid symbols are read in chronological order. See \cite[Section~2]{DH}.

For two variables, let $\mathcal M_2^0$ denote the cut domain obtained from $\mathcal M_2$ by cuts along $z_1\in\R_{>0}$, $z_2\in\R_{>0}$ and $z_1-z_2\in\R_{>0}$, with boundary values approached from the corresponding upper half-planes \cite[Section~2]{DH}. Its preferred logarithms are $\ell_0(z_1)$, $\ell_0(z_2)$ and $\ell_0(z_1-z_2)$ as defined in Section~\ref{sec:moduleprelim}. Write $f^e$ for the preferred branch on $\mathcal M_2^0$. With a retained logarithm $L$, $z^\lambda:=e^{\lambda L},\; \bigl(e^{\pm\pi i}z,\log(e^{\pm\pi i}z)\bigr):=(e^{L\pm\pi i},L\pm\pi i)$.
All branch substitutions below refer to the indicated logarithms, not solely to the complex-number arguments.

We use three absolute convergence chambers $\mathcal A_1=\{|z_1|>|z_2|>0\}$, $\mathcal A_2=\{|z_2|>|z_1|>0\}$ and $\mathcal A_3=\{|z_2|>|z_1-z_2|>0\}$. Their preferred comparison regions are
\begin{align*}
 \mathcal D_1&=\mathcal A_1\cap\mathcal M_2^0
   \cap\{|\arg(z_1-z_2)-\arg z_1|<\pi/2\},\\
 \mathcal D_2&=\mathcal A_2\cap\mathcal M_2^0
   \cap\{-3\pi/2<\arg(z_1-z_2)-\arg z_2<-\pi/2\},\\
 \mathcal D_3&=\mathcal A_3\cap\mathcal M_2^0
   \cap\{|\arg z_1-\arg z_2|<\pi/2\}.
\end{align*}
Thus branch comparisons specify both inequalities of absolute values and actual logarithms.

Let $r\geq1$, $z=(z_1,\ldots,z_r)$ be a row vector, $A\in\operatorname{GL}_r(\C)$, $\beta\in\C^r$, and
\[
\eta=zA-\beta,\qquad z=(\eta+\beta)A^{-1},\qquad
\delta\in\{0,\infty\}^r,\qquad
\xi_i=\begin{cases}\eta_i,&\delta_i=0,\\ \eta_i^{-1},&\delta_i=\infty.\end{cases}
\]
For $\epsilon=(\epsilon_i)\in\R_{>0}^r$ put $\Delta^\times_\epsilon=\{\xi:0<|\xi_i|<\epsilon_i\ (1\leq i\leq r)\}$. Define $\xi^\lambda=\prod_i\xi_i^{\lambda_i}$ for $\lambda\in\C^r$ and $(\log\xi)^m=\prod_i(\log\xi_i)^{m_i}$ for $m\in\Z_{\geq0}^r$, using a fixed branch of each logarithm. The corner $\eta=\delta$ is \emph{component-isolated} if $\Delta^\times_\epsilon$ lies in the analytic domain of $f((\eta+\beta)A^{-1})$ for some $\epsilon$. It is \emph{regular} if there are finite $\Lambda\subset\C^r$, $B\in\Z_{\geq0}$ and coefficient functions $F^b_{\lambda,m}$ holomorphic at $\xi=0$ such that each branch $b$ near that corner has the absolutely convergent expansion
\begin{equation}\label{eq:analyticregularcorner}
 f^b((\eta+\beta)A^{-1})
 =\sum_{\lambda\in\Lambda}\sum_{\substack{m\in\Z_{\geq0}^r\\0\leq m_i\leq B}}
 \xi^\lambda(\log\xi)^m F^b_{\lambda,m}(\xi),
\end{equation}
where the exponent set $\Lambda$ and bound $B$ are common to all branches $b$; the holomorphic functions and convergence radii may depend on $b$. This is the local power-series reformulation of the exponent sets $D_i=\bigcup_j(\alpha_{ij}+\Z_{\geq0})$ (at $0$) or $D_i=\bigcup_j(\alpha_{ij}-\Z_{\geq0})$ (at $\infty$) in \cite[Definition~2.6]{DH}. For a single-valued meromorphic $f$, a pole of order $a_i$ along $\xi_i=0$ means $\xi_i^{a_i}f$ is holomorphic and $a_i\in\Z_{\geq0}$ is minimal.

Let $D$ be a branch domain, and write $K\Subset D$ for a compact subset. For scalar weight contractions indexed by a set $\mathfrak W$ of internal weight tuples, \emph{absolute convergence} is $\sum_{\boldsymbol\mu\in\mathfrak W}|a_{\boldsymbol\mu}(z)|<\infty$. \emph{Normal convergence} means uniform convergence of this series on every compact $K$ in the branch domain, equivalently,
\[
\forall K\Subset D\;\forall\epsilon>0\;\exists F\Subset\mathfrak W:
\quad\sup_{z\in K}\sum_{\boldsymbol\mu\in\mathfrak W\setminus F}
|a_{\boldsymbol\mu}(z)|<\epsilon,
\]
where $F\Subset\mathfrak W$ means that $F$ is finite. Each $a_{\boldsymbol\mu}$ is the full finite-dimensional contraction at a fixed weight tuple. A domain $U\subset\C^r$ is \emph{complete Reinhardt} if $z\in U$ and $|w_i|\leq|z_i|$ for all $i$ imply $w\in U$. On such a domain, holomorphic Taylor expansions at $0$ converge absolutely and normally on compact subsets; see Lemma~\ref{lem:Reinhardt}.

We give the operator definition used throughout the paper. Fix $X\in\cC_g$, $Y\in\cC_h$ and $Z\in\cC_{gh}$. A typed operator or channel has its input and output modules, their ordered sector labels, and its bracketing specified; at each vertex the output label is the ordered product of the input labels.

\begin{definition}\label{def:completeintertwiner}
A twisted intertwining operator of type $\binom Z{X\ Y}$, in the sense of \cite[Definition~2.8]{DH} for generalized twisted modules without chosen group-element actions, is a linear map
\[
 \mathcal Y:X\otimes Y\longrightarrow Z\{x\}[\log x],
 \quad
 \mathcal Y(a,x)b=\sum_{n\in\C}\sum_{j=0}^K
                 \mathcal Y_{n,j}(a)b\,x^{-n-1}(\log x)^j
\]
with the following four properties (the coefficient truncation and branch conditions are part of the definition).
\begin{enumerate}
\item For $a\in X$, $b\in Y$, fixed $n\in\C$ and $j$,
$\mathcal Y_{n+l,j}(a)b=0$ for all sufficiently large integers $l\geq0$.

\item For every $u\in V$, $a\in X$, $b\in Y$ and $z'\in Z'$
there is a maximally continued multivalued function $f$ on $\mathcal M_2$, with preferred branch $f^e$, such that
 \begin{align*}
S_1&=\sum_\lambda
  \langle z',Y_Z(u,z_1)\pi^Z_\lambda\mathcal Y(a,z_2)b\rangle,\\
S_2&=\sum_\lambda
  \langle z',\mathcal Y(a,z_2)\pi^Y_\lambda Y_Y(u,z_1)b\rangle,\\
S_3&=\sum_\lambda
  \langle z',\mathcal Y(\pi^X_\lambda
Y_X(u,z_1-z_2)a,z_2)b\rangle
 \end{align*}
converge absolutely on $\mathcal A_1,\mathcal A_2,\mathcal A_3$, respectively. Their sums equal $f^e$ on $\mathcal D_1,\mathcal D_2,\mathcal D_3$, respectively.

\item For each $k\in\Z_{\geq3}$, states
$u_1,\ldots,u_{k-1}\in V$, and $a,b,z'$ as above, the series
 \begin{align*}
&\sum_{\lambda_1,\ldots,\lambda_{k-1}}
 \big\langle z',
Y_Z(u_1,z_1)\pi^Z_{\lambda_1}\cdots\\[-2pt] &\hspace{20mm}\cdots Y_Z(u_{k-1},z_{k-1})\pi^Z_{\lambda_{k-1}}
                         \mathcal Y(a,z_k)b\big\rangle
 \end{align*}
converges absolutely for $|z_1|>\cdots>|z_k|>0$ and has maximal multivalued analytic continuation $f_k$ on $\mathcal M_k$. For each $z_k\in\C^*$ and each choice
\[
(M_{ij})_{1\leq i<j\leq k-1}\in\Z_{>0}^{\binom{k-1}{2}},\qquad
x^{M_{ij}}Y(u_i,x)u_j\in V[[x]]\quad(1\leq i<j\leq k-1),
\]
every component-isolated corner $\eta=(z_1,\ldots,z_{k-1})-\beta=\delta$, $\beta\in\C^{k-1}$, $\delta\in\{0,\infty\}^{k-1}$, of
\[
\prod_{1\leq i<j\leq k-1}(z_i-z_j)^{M_{ij}}
 f_k(z_1,\ldots,z_{k-1},z_k)
\]
is regular in the sense of \eqref{eq:analyticregularcorner}. The quantifier agrees with \cite[Definition~2.8(3)]{DH}: the latter also specifies $M_{ij}$ for $i\geq j$, which do not occur in the product; for every prescribed set $(M_{ij})_{i<j}$, those unused entries can be chosen sufficiently large by the lower truncation of $Y(u_i,x)u_j$.

\item The formal derivative satisfies $\frac{d}{dx}\mathcal Y(a,x)b =\mathcal Y(L(-1)a,x)b.$ \end{enumerate}
\end{definition}

The four conditions are those of \cite[Definition~2.8]{DH}; the restriction $Z\in\cC_{gh}$ is our ordered-sector convention. We call (ii) the one-$V$ duality axiom and (iii) the multiple-$V$ analytic-extension axiom. Put $\cI\binom Z{X\ Y}=\cI^Z_{X,Y}$.
The fusion-rule multiplicity is $\dim_\C\cI^Z_{X,Y}$, whose finiteness is proved later. In each contraction, $\pi^M_\lambda$ is the projector from Section~\ref{sec:moduleprelim}; the summand at a fixed weight lies in a finite-dimensional weight space. Only finite sums of homogeneous components are needed for each algebraic input.

Since $g\omega=h\omega=\omega$, one-$V$ duality with $u=\omega$, followed by the usual residue calculation for $z_1^{m+1}$, gives, for $m\in\Z_{\geq-1}$,
\begin{equation}\label{eq:IOvircommutator}
[L_Z(m),\mathcal Y(a,t)]_{Y,Z}
:=L_Z(m)\mathcal Y(a,t)-\mathcal Y(a,t)L_Y(m)
=\sum_{j=0}^{m+1}\binom{m+1}{j}t^{m+1-j}\mathcal Y(L_X(j-1)a,t).
\end{equation}
Indeed, setting $\xi=z_1-t$, one has
\[
\operatorname{Res}_{\xi}\,(t+\xi)^{m+1}Y_X(\omega,\xi)a
=\sum_{j=0}^{m+1}\binom{m+1}{j}t^{m+1-j}L_X(j-1)a.
\]
The commutator is defined here between a map $Y\to Z$ and the respective Virasoro operators. In particular,
\begin{align*}
L_Z(-1)\mathcal Y(a,t)-\mathcal Y(a,t)L_Y(-1)&=\mathcal Y(L_X(-1)a,t)=\partial_t\mathcal Y(a,t),\\
L_Z(0)\mathcal Y(a,t)-\mathcal Y(a,t)L_Y(0)&=\mathcal Y(L_X(0)a,t)+t\partial_t\mathcal Y(a,t).
\end{align*}
This derivation uses only the one-$V$ duality and the conformal-vector modes, independently of product--iterate associativity for general twisted intertwining operators (cf.\ \cite[Definition~2.8 and the remarks following it]{DH}).

\begin{lemma}\label{lem:IOlogs}
Let $X,Y,Z$ be ordinary, $a\in X_\alpha$, $b\in Y_\beta$, $z'\in(Z')_\gamma$, and $\mathcal Y\in\cI\binom Z{X\ Y}$. Then $\langle z',\mathcal Y(a,t)b\rangle =C_{z',a,b}\,t^{\gamma-\alpha-\beta},\; C_{z',a,b}\in\C$,
with no nonconstant powers of $\log t$. Consequently, $\pi^Z_\gamma\mathcal Y_{n,j}(X_\alpha)Y_\beta=0\; \text{unless}\; j=0,\; -n-1=\gamma-\alpha-\beta$.
\end{lemma}
\begin{proof}
Pairing the $m=0$ instance of \eqref{eq:IOvircommutator} with $z'\in(Z')_\gamma$ and applying it to $a\in X_\alpha$, $b\in Y_\beta$ yields $\left(t\partial_t-(\gamma-\alpha-\beta)\right) F(t)=0,\; F(t):=\langle z',\mathcal Y(a,t)b\rangle$.
Write $F(t)=\sum_{n\in\C}t^{-n-1}P_n(\log t)$, $P_n\in\C[s]$, with finite logarithmic degree. Comparing the formally independent powers $t^{-n-1}$ gives $P_n'(s)+\bigl(-n-1-\gamma+\alpha+\beta\bigr)P_n(s)=0$.
If $-n-1\ne\gamma-\alpha-\beta$, the polynomial differential operator $\partial_s+(-n-1-\gamma+\alpha+\beta)$ has trivial kernel; otherwise $P_n'=0$. Thus only the exponent $\gamma-\alpha-\beta$ remains, with logarithmic degree $0$. Since every $Z_\gamma$ is finite-dimensional and $Z'=\bigoplus_\gamma Z_\gamma^*$ separates weight components, the coefficient assertion follows.
\end{proof}

A formal map is \emph{weight compatible} if $\mathcal Y_{n,j}(X_\alpha)Y_\beta\subseteq Z_{\alpha+\beta-n-1}\; (n\in\C,\ j\geq0)$.
For ordinary $X,Y,Z$, Lemma~\ref{lem:IOlogs} proves this directly (with $j>0$ identically zero). All modules of $\cC_G$ are ordinary by Proposition~\ref{prop:regularmodules}.

For $\mathcal Y\in\cI\binom Z{X\ Y}$, $X\in\cC_g$, $Y\in\cC_h$, define the skew transforms on the underlying vector spaces by
\begin{equation}\label{eq:skewprelim}
(\Omega_\pm\mathcal Y)(b,t)a =e^{tL_Z(-1)}\mathcal Y(a,e^{\pm\pi i}t)b .
\end{equation}
The logarithm in the right side is $\log t\pm\pi i$. The transformations have types $\Omega_+\mathcal Y\in\cI\binom Z{Y\ T_{h^{-1}}X},\; \Omega_-\mathcal Y\in\cI\binom Z{T_gY\ X}$.
These are isomorphisms of operator spaces by \cite[Theorem~3.1 and the immediately following corollary]{DH}, for generalized twisted modules \emph{without} actions of the twisting automorphisms. Here the pullback $T_k$ is the functor $\phi_k$ of \cite[Section~2]{DH}, not the distinct $\varphi_k$ construction requiring an action of $k$ on the module. For a labeled object $W\in\cC_q$ put $\operatorname{tw}(W)=q$ (the label is retained in the external sum $\cC_G$). Then
\[
\operatorname{tw}(Y)\operatorname{tw}(T_{h^{-1}}X)
=h(h^{-1}gh)=gh,\qquad
\operatorname{tw}(T_gY)\operatorname{tw}(X)
=(ghg^{-1})g=gh.
\]
In particular $\Omega_-\Omega_+\mathcal Y=\mathcal Y=\Omega_+\Omega_-\mathcal Y$ with the corresponding source/target identifications.

The contragredient transformations are defined by
\begin{align*}
 \langle(A_\pm\mathcal Y)(a,t)z',b\rangle
&=\big\langle z',
    \mathcal Y(e^{tL_X(1)}e^{\pm\pi iL_X(0)}
t^{-2L_X(0)}a,t^{-1})b\big\rangle .
\end{align*}
Here $\log(t^{-1})=-\log t$. Their types are $A_+\mathcal Y\in\cI\binom{T_gY'}{X\ Z'},\; A_-\mathcal Y\in\cI\binom{Y'}{X\ T_{g^{-1}}Z'}$.
Theorem~3.3 of \cite{DH} proves that these operators satisfy the four analytic axioms of Definition~\ref{def:completeintertwiner} and have the displayed types. Their sector labels satisfy
\[
\operatorname{tw}(T_gY')=gh^{-1}g^{-1},\qquad
\operatorname{tw}(T_{g^{-1}}Z')=g^{-1}(gh)^{-1}g=g^{-1}h^{-1},
\]
and hence $g(gh)^{-1}=gh^{-1}g^{-1}$ and $g(g^{-1}h^{-1})=h^{-1}$. No action of $G$ on an individual twisted module is required for these transformations; only the pullbacks of Section~\ref{sec:moduleprelim} occur. For $a\in X_\lambda$, $e^{\pm\pi iL_X(0)}t^{-2L_X(0)}a =e^{\pm\pi i\lambda}t^{-2\lambda}a$.
Since $L_X(1)$ lowers weight by $1$ and $X$ is lower bounded, $e^{tL_X(1)}a$ is a finite sum. The $e^{tL_Z(-1)}$ in \eqref{eq:skewprelim} acts in $\overline Z$: for fixed output weight $\mu$, its projection involves finitely many inputs from the weights $\mu-j$, $j\in\N$.

For $X\in\cC_g$, the module field $Y_X^g$ gives an intertwining operator in $\cI\binom X{V\ X}$, with the multiple-$V$ property supplied by the twisted-module correlation-function rationality (cf.\ \cite[Sections~2--3]{DH}). Since $T_1=\id$, its $\Omega_-$-transform is the creation operator
\[
 \Gamma_X(w,t)v=(\Omega_-Y_X)(w,t)v
=e^{tL_X(-1)}Y_X(v,e^{-\pi i}t)w
           \in\cI\binom X{X\ V}.
\]
In particular $\Gamma_X(w,t)\mathbf1=e^{tL_X(-1)}w$.

\subsection{\texorpdfstring{$P(z)$}{P(z)} products, transport and channels}
Fix $g,h\in G$, $X\in\cC_g$, $Y\in\cC_h$, $Z\in\cC_{gh}$, $z\in\C^*$, and $L\in\C$ with $e^L=z$. For $\mathcal Y\in\cI\binom Z{X\ Y}$, its evaluation at $L$ is the linear map
\[
\operatorname{ev}_L(\mathcal Y)=I_{\mathcal Y,L}:X\otimes Y\longrightarrow\overline Z,
\qquad
I_{\mathcal Y,L}(a\otimes b)=\mathcal Y(a,x)b
\big|_{x^\rho=e^{\rho L},\,\log x=L}.
\]
The evaluation is defined componentwise in $\overline Z=\prod_{\delta}Z_\delta$. Indeed, for $a\in X_\alpha$, $b\in Y_\beta$, and $\delta\in\wt(Z)$, Lemma~\ref{lem:IOlogs} gives
\begin{equation}\label{eq:Pzinverse}
\pi^Z_\delta\mathcal Y(a,t)b
=e^{(\delta-\alpha-\beta)(\log t-L)}
\pi^Z_\delta I_{\mathcal Y,L}(a\otimes b).
\end{equation}
Thus $\operatorname{ev}_L$ is injective. Conversely, \eqref{eq:Pzinverse}, extended linearly over the finite homogeneous decompositions of $a,b$, recovers the unique operator whose evaluation is a given map in
$\operatorname{im}(\operatorname{ev}_L)$. Following
\cite[Section~4, definition of twisted $P(z)$-intertwining maps]{DH}, put
\[
\cI_{P(z),L}\binom Z{X\ Y}
:=\operatorname{im}\!\left(
\operatorname{ev}_L:\cI\binom Z{X\ Y}
\longrightarrow\Hom_\C(X\otimes Y,\overline Z)\right).
\]
A $P(z)$-intertwining map means an element of this image; no arbitrary linear map into $\overline Z$ is thereby declared to be an intertwining map. The notation $P(z)$ corresponds to the ordered marked points $(\infty,z,0)$.

\begin{definition}\label{def:Pzproduct}
A $P(z)$-product of $(X,Y)$ is $(U,I)$ with $U\in\cC_{gh}$ and $I\in\cI_{P(z),L}\binom U{X\ Y}$. A morphism $(U,I)\to(W,J)$ is $f\in\Hom_V(U,W)$ such that $\overline f\circ I=J$.
A $P(z)$-tensor product is an initial object in this category of products; when it exists, denote it by $(X\bo_{P(z)}Y,I^z_{X,Y})$.
\end{definition}

Equivalently, for every $W\in\cC_{gh}$, the map
\begin{equation}\label{eq:Pzuniversal}
\Hom_V(X\bo_{P(z)}Y,W)
\xrightarrow{\ \sim\ }
\cI_{P(z),L}\binom W{X\ Y},\qquad
f\longmapsto\overline f\circ I^z_{X,Y}
\end{equation}
is a bijection, natural in $W$ \cite[Section~4]{DH}. If $f:X\to\widetilde X$ and $q:Y\to\widetilde Y$ are module maps and the two $P(z)$-tensor products exist, their tensor-product map is
uniquely specified by $\overline{f\bo_{P(z)}q}\circ I^z_{X,Y} =I^z_{\widetilde X,\widetilde Y}\circ(f\otimes q)$. Hence $\id_X\bo_{P(z)}\id_Y=\id_{X\bo_{P(z)}Y}$ and
$(f_2f_1)\bo_{P(z)}(q_2q_1)
=(f_2\bo_{P(z)}q_2)(f_1\bo_{P(z)}q_1)$, by
\eqref{eq:Pzuniversal}. Existence in $\cC_{gh}$ is proved in
Theorem~\ref{thm:binary}; put $X\bo Y=X\bo_{P(1)}Y$, with $\log1=0$. For each chosen pair $(z,L)$, define the associated universal formal operator by $\mathcal Y^z_{X,Y}:=\operatorname{ev}_L^{-1}(I^z_{X,Y}) \in\cI\binom{X\bo_{P(z)}Y}{X\ Y}$.

Let $\gamma:[0,1]\to\C^*$ satisfy $\gamma(0)=z_1$, $\gamma(1)=z_2$. Fix the logarithm $L_2$ used for $I_{X,Y}^{z_2}$ and let $L_1^\gamma$ be its continuation from $z_2$ to $z_1$ along the reverse of $\gamma$; thus $e^{L_1^\gamma}=z_1$.
For the universal formal operator $\mathcal Y_{X,Y}^{z_2}\in\cI\binom{X\bo_{P(z_2)}Y}{X\ Y}$, define the parallel transport $\mathsf T_\gamma$ by
\begin{equation}\label{eq:Pztransport}
\begin{gathered}
\mathsf T_\gamma:X\bo_{P(z_1)}Y\xrightarrow{\ \sim\ }X\bo_{P(z_2)}Y,\\
\overline{\mathsf T_\gamma}\circ I_{X,Y}^{z_1}
=\operatorname{ev}_{L_1^\gamma}(\mathcal Y_{X,Y}^{z_2}).
\end{gathered}
\end{equation}
The right side is a $P(z_1)$-intertwining map, so
\eqref{eq:Pzuniversal} determines $\mathsf T_\gamma$.
The inverse-path construction gives $\mathsf T_{\gamma^{-1}}=\mathsf T_\gamma^{-1}$; for composable paths $\gamma_1,\gamma_2$,
$\mathsf T_{\gamma_2*\gamma_1}
=\mathsf T_{\gamma_2}\circ\mathsf T_{\gamma_1}$.
These statements use analytic continuation of the logarithm and the universal property, as in \cite[Section~4, parallel transport]{DH}. The resulting map depends only on the endpoint-fixed homotopy class.
In the common formal representative of Theorem~\ref{thm:binary}, $\mathsf T_\gamma=\id$ when $z_1,z_2>0$ and $\gamma$ is the positive-real interval, with its real logarithms.

For $X\in\cC_g$, the operators $Y_X^g\in\cI\binom X{V\ X}$ and $\Gamma_X=\Omega_-Y_X^g\in\cI\binom X{X\ V}$ from Section~\ref{sec:moduleprelim} and \eqref{eq:skewprelim} induce, by \eqref{eq:Pzuniversal}, the module maps
\begin{equation}\label{eq:Pzunits}
\begin{aligned}
l_X &: V\bo X\longrightarrow X,
&\overline{l_X}\circ I^1_{V,X}(v\otimes w)&=Y_X^g(v,1)w,\\
r_X &: X\bo V\longrightarrow X,
&\overline{r_X}\circ I^1_{X,V}(w\otimes v)&=\Gamma_X(w,1)v.
\end{aligned}
\end{equation}
Here $\Gamma_X(w,1)\mathbf1=e^{L_X(-1)}w$; Lemma~\ref{lem:units} proves $l_X$ and $r_X$ are isomorphisms.

For $X\in\cC_g$, $Y\in\cC_h$, and $z=e^L$, set $L_-=L+\pi i$, so $e^{L_-}=-z$. The target $T_gY\in\cC_{ghg^{-1}}$ gives $(ghg^{-1})g=gh$; moreover, $T_{g^{-1}}(T_gY)=Y$.
Define $\mathcal Z:=\operatorname{ev}_{L_-}^{-1}(I^{-z}_{T_gY,X})\in
\cI\binom{T_gY\bo_{P(-z)}X}{T_gY\ X}$.
By \cite[Section~3, skew symmetry, and Section~4, crossed commutativity]{DH},
$\Omega_+\mathcal Z\in
\cI\binom{T_gY\bo_{P(-z)}X}{X\ Y}$.
Thus \eqref{eq:Pzuniversal} gives an isomorphism
\begin{equation}\label{eq:Pzrawcrossed}
\begin{gathered}
\mathcal R^z_{X,Y}:X\bo_{P(z)}Y
\xrightarrow{\ \sim\ }T_gY\bo_{P(-z)}X,\\
\overline{\mathcal R^z_{X,Y}}\circ I^z_{X,Y}(a\otimes b)
=\operatorname{ev}_{L}(\Omega_+\mathcal Z)(a\otimes b)
=e^{zL_{T_gY\bo_{P(-z)}X}(-1)}
I_{T_gY,X}^{-z}(b\otimes a).
\end{gathered}
\end{equation}
The inverse is induced by $\Omega_-$ together with $\Omega_-\Omega_+=\id$; neither transformation requires an extra $G$-action on the underlying twisted modules. For $\gamma$ from $-1$ to $1$ through the upper half-plane,
$\log(-1)=\pi i$ is the backward continuation of $\log1=0$, and $c_{X,Y} :=\mathsf T_\gamma^{T_gY,X}\circ\mathcal R^1_{X,Y}: X\bo Y\xrightarrow{\ \sim\ }T_gY\bo X$. This is the branch normalization of \cite[Section~4]{DH}.

For the channel notation used below, let $\tau$ be a finite rooted binary tree; each internal vertex $v$ has ordered children $v_L,v_R$. Define its anchor recursively by $p(v)=p(v_R)$ and let
$r_v=p(v_L)-p(v_R)$; the anchor of each leaf is its marked position. Translate all positions so that $p(\operatorname{root}(\tau))=0$. On the configuration domain where every $r_v\ne0$, put
\[
R=r_{\operatorname{root}(\tau)},\qquad
u_v=\frac{r_v}{r_{\operatorname{parent}(v)}}
\quad(v\ne\operatorname{root}(\tau)),\qquad
r_v=R\!\prod_{w\in\operatorname{path}(\operatorname{root}(\tau),v)
\setminus\{\operatorname{root}(\tau)\}}\!u_w,
\]
where the product runs over the nonroot vertices on the unique root--$v$ path. The same chart is restricted further to pairwise distinct leaf positions; only its nonzero, collision-free locus is used.
Each internal edge carries its specified output/input module $M$, and contractions use $\pi^M_\mu:M\to M_\mu$, $\mu\in\wt(M)$. For a basis $\{e_{\mu,j}\}_j$ of $M_\mu$ and its dual $\{e_{\mu,j}^*\}_j$, the finite-rank projector is
\[
\pi^M_\mu=\sum_j e_{\mu,j}\otimes e_{\mu,j}^*
\quad\text{in }\End_\C(M),
\]
where $e_{\mu,j}^*$ is extended by zero off $M_\mu$. The full sum $\sum_{\mu\in\wt(M)}\pi^M_\mu$ is used only coefficientwise on $M$; for completed channel contractions its absolute and normal convergence is established in
Sections~\ref{sec:regularity} and \ref{sec:threeinputward}.

For $A\in\cC_g$, $B\in\cC_h$, $C\in\cC_k$, $M\in\cC_{hk}$, $N\in\cC_{gh}$, $D\in\cC_{ghk}$, let
\[
\mathcal Y_1\in\cI\binom D{A\ M},\quad
\mathcal Y_2\in\cI\binom M{B\ C},\quad
\mathcal Y_3\in\cI\binom D{N\ C},\quad
\mathcal Y_4\in\cI\binom N{A\ B}.
\]
For $a\in A$, $b\in B$, $c\in C$, and $d'\in D'$, the two formal contractions are
\begin{align*}
\Phi_P(z_1,z_2)&=
\sum_{\mu\in\wt(M)}
\langle d',\mathcal Y_1(a,z_1)\pi^M_\mu
\mathcal Y_2(b,z_2)c\rangle,\\
\Phi_I(z_1,z_2)&=
\sum_{\nu\in\wt(N)}
\langle d',\mathcal Y_3(
\pi^N_\nu\mathcal Y_4(a,z_1-z_2)b,z_2)c\rangle.
\end{align*}
Their positive-real comparison locus is $z_1>z_2>z_1-z_2>0$, with $\log z_1,\log z_2,\log(z_1-z_2)\in\R$. No product--iterate equality or convergence is assumed in these definitions.

\subsection{Monoidal and crossed braided categories}\label{sec:categoricalconventions}
A $\C$-linear category $\mathcal A$ has $\C$-vector spaces $\Hom_{\mathcal A}(X,Y)$ and $\C$-bilinear composition. We work with additive categories (finite biproducts). For functors $F,H:\mathcal A\to\mathcal B$, a natural transformation $\eta:F\Rightarrow H$ is a family $\eta_X:FX\to HX$ satisfying $H(f)\eta_X=\eta_YF(f)$ for all $f:X\to Y$; it is invertible if and only if every $\eta_X$ is invertible. Write $\operatorname{Nat}(F,H)$ for the set of natural transformations $F\Rightarrow H$. A bifunctor is biadditive if it is additive in each variable. An idempotent $p\in\End(P)$, $p^2=p$, splits if there are $i:X\to P$ and $r:P\to X$ with $ri=\id_X$ and $ir=p$. We call $X$ a retract (direct summand) of $P$. Cauchy completeness means finite biproducts and splitting of all idempotents.

For the $\C$-linear abelian categories used below, \emph{finite semisimple} means that $\Irr(\mathcal A)=\{[S_1],\ldots,[S_m]\}$ is finite, $\dim_{\C}\Hom_{\mathcal A}(X,Y)<\infty$ for all $X,Y$, and every object is a finite direct sum of simple objects. Schur's lemma gives $\Hom(S_i,S_j)=0$ for $i\ne j$ and $\End(S_i)\cong\C$. Every short exact sequence then splits, as does every idempotent. An object $X$ is indecomposable if $X\ne0$ and $X\cong A\oplus B$ implies $A=0$ or $B=0$. A morphism $f:X\to Y$ is epic if $\forall Z$ and $p,q:Y\to Z$, $pf=qf\Rightarrow p=q$; for the module categories here, this is equivalent to surjectivity.

A monoidal category $(\mathcal A,\bo,\mathbf u,a,l,r)$ consists of a tensor bifunctor $\bo$, a unit $\mathbf u$, and natural isomorphisms $a_{X,Y,Z}:(X\bo Y)\bo Z\xrightarrow{\sim}X\bo(Y\bo Z)$, $l_X:\mathbf u\bo X\xrightarrow{\sim}X$ and $r_X:X\bo\mathbf u\xrightarrow{\sim}X$. Composition is read from right to left; for $X,Y,Z,W\in\mathcal A$, the pentagon and triangle are
\begin{align*}
a_{X,Y,Z\bo W}a_{X\bo Y,Z,W} &=(\id_X\bo a_{Y,Z,W})a_{X,Y\bo Z,W} (a_{X,Y,Z}\bo\id_W),\\ (\id_X\bo l_Y)a_{X,\mathbf u,Y} &=r_X\bo\id_Y.
\end{align*}
For the module category the unit will be $V$. The symbol $\mathbf1\in V$ continues to denote its vacuum vector.

A left dual of $X$ is an object $X^*$ with $\operatorname{ev}_X:X^*\bo X\to\mathbf u$ and $\operatorname{coev}_X:\mathbf u\to X\bo X^*$ satisfying
\begin{align*}
r_X(\id_X\bo\operatorname{ev}_X)a_{X,X^*,X} (\operatorname{coev}_X\bo\id_X)l_X^{-1}&=\id_X,\\ l_{X^*}(\operatorname{ev}_X\bo\id_{X^*})a^{-1}_{X^*,X,X^*} (\id_{X^*}\bo\operatorname{coev}_X)r_{X^*}^{-1} &=\id_{X^*}.
\end{align*}
A right dual of $X$ is an object ${}^*X$ with $\operatorname{ev}^R_X:X\bo{}^*X\to\mathbf u$ and $\operatorname{coev}^R_X:\mathbf u\to{}^*X\bo X$, with
\begin{align*}
l_X(\operatorname{ev}^R_X\bo\id_X)a^{-1}_{X,{}^*X,X} (\id_X\bo\operatorname{coev}^R_X)r_X^{-1}&=\id_X,\\ r_{{}^*X}(\id_{{}^*X}\bo\operatorname{ev}^R_X)a_{{}^*X,X,{}^*X} (\operatorname{coev}^R_X\bo\id_{{}^*X})l_{{}^*X}^{-1} &=\id_{{}^*X}.
\end{align*}
Rigidity means that every object has a left and a right dual. A \emph{fusion category} over $\C$ is a rigid finite semisimple $\C$-linear monoidal category with $\End(\mathbf u)\cong\C$ and a $\C$-bilinear tensor bifunctor \cite[Definition~4.1.1]{EGNO}.

A strong monoidal functor $F:(\mathcal A,\bo,\mathbf u_{\mathcal A})\to(\mathcal B,\bo,\mathbf u_{\mathcal B})$ has natural isomorphisms $J_{X,Y}:FX\bo FY\xrightarrow{\sim}F(X\bo Y)$ and $j:\mathbf u_{\mathcal B}\xrightarrow{\sim}F\mathbf u_{\mathcal A}$, satisfying
\begin{align*}
F(a_{X,Y,Z})J_{X\bo Y,Z}(J_{X,Y}\bo\id) &=J_{X,Y\bo Z}(\id\bo J_{Y,Z})a_{FX,FY,FZ},\\ F(l_X)J_{\mathbf u,X}(j\bo\id)&=l_{FX},\\ F(r_X)J_{X,\mathbf u}(\id\bo j)&=r_{FX}.
\end{align*}
The maps $J$ are called tensorators. A functor is faithful if it is injective on each morphism space, and exact if it preserves short exact sequences. The opposite category reverses the arrows; a contravariant functor is an ordinary functor from that opposite category.

We use the braided $G$-crossed convention of \cite[Definition~8.24.1, (8.104)--(8.107)]{EGNO}, equivalently \cite[Definition~4.41]{DGNO}. We impose $T_kT_l=T_{kl}$ as \emph{underlying functors}, so their monoidal action compositors $\gamma_{k,l}:T_kT_l\Rightarrow T_{kl}$ are identities; in general the tensorators $J^k$ are \emph{not} identities.

A $G$-grading of a $\C$-linear monoidal category is $\cC_G=\bigoplus_{g\in G}\cC_g$, with $\Hom(X_g,Y_h)=0$ for $X_g\in\cC_g$, $Y_h\in\cC_h$, $g\ne h$, $\cC_g\bo\cC_h\subseteq\cC_{gh}$ and $\mathbf u\in\cC_1$. It is faithful if $\cC_g\ne0$ for all $g\in G$. A strict underlying action is a family of autoequivalences $T_k:\cC_G\to\cC_G$ with $T_kT_l=T_{kl}$ and $T_1=\id$. Each $T_k$ is strong monoidal, with $J^k_{X,Y}:T_kX\bo T_kY\xrightarrow{\sim}T_k(X\bo Y)$ and $j_k:\mathbf u\xrightarrow{\sim}T_k\mathbf u$ satisfying the equations above and
\begin{equation}\label{eq:Jcomp}
J^{kl}_{X,Y}=T_k(J^l_{X,Y})J^k_{T_lX,T_lY},\qquad j_{kl}=T_k(j_l)j_k,\quad J^1=\id,\quad j_1=\id .
\end{equation}
The grading condition is $T_k(\cC_g)=\cC_{kgk^{-1}}$ ($g,k\in G$).

A crossed braiding is a natural family $c_{X,Y}:X\bo Y\xrightarrow{\sim}T_gY\bo X$ ($X\in\cC_g$, $Y\in\cC_G$), extended additively in $X$. For all $k\in G$, equivariance is
\begin{equation}\label{eq:equiv}
T_k(c_{X,Y})J^k_{X,Y} =J^k_{T_gY,X}c_{T_kX,T_kY}.
\end{equation}
The common target is $T_k(T_gY\bo X)$, using $T_{kgk^{-1}}T_k=T_kT_g$. For $X\in\cC_g$, $Y\in\cC_h$, $Z\in\cC_G$, the two crossed hexagons are
\begin{align}
&(\id_{T_gY}\bo c_{X,Z})a_{T_gY,X,Z} (c_{X,Y}\bo\id_Z)
       \notag\\
&\quad=a_{T_gY,T_gZ,X} ((J^g_{Y,Z})^{-1}\bo\id_X)c_{X,Y\bo Z}a_{X,Y,Z},
                                        \label{eq:hex1}\\
&(c_{X,T_hZ}\bo\id_Y)a^{-1}_{X,T_hZ,Y} (\id_X\bo c_{Y,Z})
       \notag\\
&\quad=a^{-1}_{T_{gh}Z,X,Y}c_{X\bo Y,Z}a^{-1}_{X,Y,Z}.
                                        \label{eq:hex2}
\end{align}
In \eqref{eq:hex2}, $T_gT_h=T_{gh}$ identifies both codomains with $(T_{gh}Z\bo X)\bo Y$. Thus \eqref{eq:Jcomp}, \eqref{eq:equiv}, \eqref{eq:hex1} and \eqref{eq:hex2}, together with the monoidal functor identities, are precisely the specialization of \cite[Definition~8.24.1]{EGNO} to $\gamma_{k,l}=\id$. If $G=\{1\}$, $T_1=\id$, $J^1=j_1=\id$, then $c_{X,Y}:X\bo Y\to Y\bo X$ is the usual braiding.

For the strict underlying action, a $G$-equivariant object is $(X,(u_k)_{k\in G})$ with $u_k:T_kX\xrightarrow{\sim}X$, $u_1=\id_X$ and $u_{kl}=u_kT_k(u_l)$. A morphism $f:(X,u^X)\to(Y,u^Y)$ satisfies $fu_k^X=u_k^YT_k(f)$ for all $k\in G$. The equivariantization $\cC_G^G$ has these objects and morphisms, with forgetful functor $U(X,u)=X$ \cite[Section~2.7]{EGNO}.

In a Cauchy complete $\C$-linear \emph{braided} monoidal category $\mathcal B$ with $\dim\Hom(A,B)<\infty$ and $\End(\mathbf u)=\C$, an indecomposable $X$ is \emph{non-negligible} if $\exists Y\in\mathcal B$ such that $\mathbf u$ is a retract of $Y\bo X$; equivalently, $\exists\,i:\mathbf u\to Y\bo X,\ r:Y\bo X\to\mathbf u$ with $ri=\id_{\mathbf u}$. It has \emph{moderate growth} if $\exists n\in\Z_{>0}$ with $\dim_{\C}\End(X^{\bo n})<n!$. Under precisely these categorical hypotheses, \cite[Theorem~1.1]{EtingofPenneys} asserts
\[
  X\text{ indecomposable, non-negligible, of moderate growth}
  \quad\Longrightarrow\quad X\text{ rigid}.
\]
We apply this theorem only to $\mathcal B=\cC_G^G$, after the braided monoidal structure has been constructed.

For a finite semisimple monoidal category $\mathcal D$ with simples $S_1,\ldots,S_m$, the Grothendieck ring is $K_0(\mathcal D)=\bigoplus_{i=1}^m\Z[S_i]$, with $[X\oplus Y]=[X]+[Y]$ and $[X][Y]=[X\bo Y]$. Write $[X\bo S_j]=\sum_{i=1}^m (N_X)_{ij}[S_i]$; then $N_X\in\operatorname{Mat}_{m\times m}(\Z_{\ge0})$. These matrices are used in Section~\ref{sec:rigidity}.

For $A,B$ in a locally small category, Yoneda gives $\operatorname{Nat}(\Hom(B,-),\Hom(A,-))\xrightarrow{\sim}\Hom(A,B),\; \eta\longmapsto b:=\eta_B(\id_B)$.
Indeed, naturality yields $\eta_Z(f)=f b$ for every $f:B\to Z$, and conversely $b:A\to B$ gives $\eta_Z(f)=fb$. This fixes the variance and arrow direction used below.

\subsection{Covers, currents and principal parts}\label{sec:geometricconventions}
Put $\cO=\C[[x]]$ and $\Aut\cO=\{f(x)=\sum_{j\ge1}a_jx^j\in x\C[[x]]:a_1\in\C^*\}$.
The action on $\cO$ is by substitution; composition of the induced point maps has the opposite order. For $r\ge1$, the $r$-jet group is $\Aut\cO/\{f:f(x)\equiv x\pmod{x^{r+1}}\}$. An operator is nilpotent if some positive power is zero; it is unipotent if its difference from the identity is nilpotent. The subgroup $\Aut_+\cO=\{f\in\Aut\cO:f'(0)=1\}$ is prounipotent,
with Lie algebra $x^2\C[[x]]\partial_x$. On the graded state space, fix the coordinate-change convention $\mathrm dR_V(x^{j+1}\partial_x)=-L(j)$ ($j\ge0$) and
$R_V(ax)=a^{-L(0)}$ ($a\in\C^*$), as in \cite[Sections~3--4]{FS}. Write $f(x)=a_1u(x)$ with $u(x)=x+O(x^2)$ and $u(x)=\exp(\sum_{j\geq1}c_jx^{j+1}\partial_x)x$. The $c_j$ are determined successively by the coefficients of $u$. In the ring-composition convention, $R_V(f)=\exp\!\left(-\sum_{j\geq1}c_jL(j)\right)a_1^{-L(0)}.$ For the CFT-type VOA under consideration, put
$V_{\le d}=\bigoplus_{0\le r\le d}V_r$ ($d\in\N$). Then $L(j)V_r\subseteq V_{r-j}$ and $L(j)|_{V_{\le d}}=0$ for $j>d$;
thus $R_V$ restricts to $V_{\le d}$ and factors through a finite jet group. At an observation point $x$, a coordinate change $f$ uses the centered series $f_x(v)=f(x+v)-f(x)$. A homogeneous state is primary if $L(j)u=0$ for every $j>0$. Its coordinate action is a weight character. For other states the same representation includes the positive-Virasoro terms \cite{FS}.

A flat module is one whose tensor product preserves exact sequences; a flat family has local rings flat over the base. For the smooth relative curves considered here, the geometric fibers are smooth
one-dimensional curves and relative formal coordinates exist locally. An
\'{e}tale morphism is smooth of relative dimension zero; over $\C$ it is
locally a biholomorphism after analytification. A proper morphism is separated, of finite type and universally closed. A projective morphism admits a closed embedding into a projective space over its base. For a smooth relative curve $C/S$, let $\omega_{C/S}=\Omega^1_{C/S}$ be its canonical line bundle, abbreviated to $\omega$ when the base is fixed. The torsor $\operatorname{Coord}(C/S)$ has, over each point, the set of formal fiber coordinates; the coordinate group acts freely and transitively on that set. If $A=\bigoplus_{r} (A\cap V_r)\subset V$ is a finite-dimensional
$R_V(\Aut\cO)$-stable subspace, put $\mathscr V_A=\operatorname{Coord}(C/S)\times^{\Aut\cO}A$, where the associated-bundle convention uses the above $R_V$. The current bundle is $\mathscr V_A\otimes\omega$; since $R_V(ax)|_{V_d}=a^{-d}$,
its weight-$d$ associated graded factor has transition character
$a^{1-d}$, i.e. that of $\omega^{\otimes(1-d)}$ \cite{FS}. A gauge is a change of local frame, expressed by an invertible matrix; a Laurent gauge has Laurent-series entries with a finite pole bound.

The coarse marked curve is $\mathbb P^1$, with moving points $(p_1,\ldots,p_m)\in F_m(\C^*)$, a fixed input at $0$, and output at $\infty$. Fix a subgroup $H\leq G$. A principal $H$-cover of its punctured part is an étale covering with a left $H$-action that is free and transitive on every fiber. Its branched extension is an $H$-cover of the marked curve. A positive meridian is a counterclockwise loop about a marking, conjugated into
the chosen basepoint; at $\infty$ positivity is measured in $z^{-1}$. We fix its based path as follows. At a real reference configuration $p_1>\cdots>p_m>0$, choose a basepoint to the right of all inputs and a geometric system of based meridians in the \emph{input order}
\begin{equation}\label{eq:orderedmeridians}
 (\gamma_{p_1},\ldots,\gamma_{p_m},\gamma_0,\gamma_\infty),
 \qquad
 \gamma_{p_1}\cdots\gamma_{p_m}\gamma_0\gamma_\infty=1.
\end{equation}
Loop products are traversed from left to right. For a chosen unramified based lift $\widetilde b$, define $\rho(\gamma)\in H$ by requiring that the lift of $\gamma$ starting at $\widetilde b$ end at $\rho(\gamma)\widetilde b$; with the left deck action this makes $\rho$ a homomorphism. The reference paths are continued with the configuration. The subscript $0$ always denotes the geometric zero marking, not the first place in the ordered list. Thus a binary channel at $(z,0,\infty)$ uses $(\gamma_z,\gamma_0,\gamma_\infty)$ and labels $(g,h,(gh)^{-1})$.

For a disk containing consecutive inputs, its based positive boundary is the product of their meridians in this order. We retain the common boundary stem when the disk becomes a bubble. For the lower half-turn in which the right input passes below the left one, comparison of the transported stems with the terminal ordered system gives
\begin{equation}\label{eq:meridianhurwitz}
 (\gamma_i,\gamma_{i+1})\longmapsto
 (\gamma_i\gamma_{i+1}\gamma_i^{-1},\gamma_i),
 \qquad
 (a,b)\longmapsto(aba^{-1},a).
\end{equation}
Indeed, the stem of the input passing below is retained, while the other stem acquires its positive meridian and inverse on its two sides; their ordered product remains $\gamma_i\gamma_{i+1}$. The inverse half-turn gives $(a,b)\mapsto(b,b^{-1}ab)$. These conventions also specify the based paths on normalized components and along the exchange paths in Section~\ref{sec:coherence}. A change of a distinguished lift by $h\in H$ conjugates its label by $h$ and replaces its module by $T_hM$; a root is always continued from its initial choice.

The image $g\in H$ of a marked meridian generates the cyclic inertia group $I_p=\operatorname{Stab}_H(p)=\langle g\rangle$ at its chosen lift. Put $N=|I_p|$, $\zeta_N=e^{2\pi i/N}$ and
$\mu_N=\{\zeta\in\C^*:\zeta^N=1\}$. Fix a coarse parameter $x$ and a special root $u$ with $x=u^N$. We act on functions by $(hf)(P)=f(h^{-1}P)$, with the convention $gu=\zeta_N^{-1}u$.

The Riemann existence correspondence identifies finite topological covering spaces of a smooth complex algebraic variety with finite \'{e}tale algebraic covers \cite[Expos\'e XII, Theorem 5.1]{SGA1}. For a finite cover of a smooth proper curve punctured at finitely many points,
normalization extends the cover uniquely as a finite normal cover. In
characteristic zero its completed local map is $\C[[x]]\hookrightarrow
\C[[u]]$, $x\mapsto u^N$, with $N$ equal to the inertia order. An isomorphism fixing a based lift fixes the connected principal cover. A nilpotent thickening $X\subset X'$ is a closed inclusion defined by an ideal sheaf whose sufficiently high power is zero; it preserves the underlying topological space. For such a thickening, restriction is an equivalence between finite \'{e}tale covers of $X'$ and of $X$ \cite[Tag 0BQB]{Stacks}. Consequently, for an \'{e}tale $H$-cover of a constant punctured curve,
its pullback to any nilpotent thickening is the unique extension
\emph{up to the unique $H$-equivariant isomorphism restricting to the fixed
special-fiber identification}; a chosen based lift removes the remaining cover automorphisms.

A based lift, distinguished marking lifts and chosen special roots form a framing. A framed Hurwitz chart parametrizes these cover data over a configuration chart; the finite étale charts used here are constructed in Lemma~\ref{lem:framedcharts}. At other lifts the labels are conjugated and the assigned modules are the strict pullbacks $T_hM$. For a subgroup $K\leq H$ and a component $C_0$ with $K$-action, the induced cover is $H\times^K C_0=(H\times C_0)/\!\sim$, where $(h,c)\sim(hk^{-1},kc)$. Its states are induced by the same pullbacks.

A node has two branches, locally given by $xy=0$. An admissible cover of a family of nodal marked curves is étale away from markings and nodes and, after finite base extension, has at a node over a base parameter $t$ the equations $UV=t^a,\ x=U^N,\ y=V^N,\ xy=t^{aN}.$ The node inertia is balanced: $(U,V)\mapsto(\zeta U,\zeta^{-1}V)$ on points, $\zeta^N=1$. For a genus-zero coarse curve, stability means that every component has at least three special points, counting both markings and nodes. The balanced admissible-cover compactification of \cite[Sections~2--4]{ACV}
provides such local models for covers of stable pointed curves, after finite base extension; stability of the target and the balancing assumption are part of this assertion. A root base change of a parameter is $q=t^e$,
with positive integer $e$ and retained logarithms $\log q=e\log t$. A thickening to order $J$ is base change to $\C[[t]]/(t^{J+1})$.

A compactification of a variety $X$ is a proper variety $\overline X$ with an open immersion $X\hookrightarrow\overline X$. For an integral normal base $Y$ and a finite extension $K'/K(Y)$, normalization means the
integral closure of $\cO_Y$ in $K'$ (on each affine open). For a smooth variety and an ideal sheaf $I$, its blow-up is $\operatorname{Proj}\bigoplus_{n\geq0}I^n$. An arrangement of smooth subvarieties has clean intersections if its intersections are smooth and their tangent spaces are the intersections of the corresponding tangent spaces. A building set is a collection such that the minimal members containing each intersection meet transversely and have that intersection as their common intersection. The wonderful compactification is the iterated blow-up of these centers in an inclusion-compatible order. The specific cluster arrangement and its centers are given in Definition~\ref{def:clustercompactification}.

For a marked point with inertia $I_p=\langle g\rangle$, choose $x=u^N$ and a homogeneous eigenbasis $b\in V^{(r)}$ of the finite state frame. Express the pullback of an $I_p$-invariant meromorphic state-current section in this frame as
$\sigma=\sum_{r,b}f_{b,p}(u)b\,dx$ with $ f_{b,p}(u)\in u^r\C((u^N))$ for $gb=\zeta_N^rb$; thus the coefficients are Laurent series in $u$, \emph{not} necessarily
functions of $u^N$. The congruence follows from $g(f_{b,p}(u)b)=f_{b,p}(\zeta_N^{-1}u)\zeta_N^rb$. An admissible current has finite poles confined to the marked orbit divisors. The descended local action on a weak $g$-twisted module $M_i$ is
\begin{equation}\label{eq:localcurrent25}
R_i(\sigma)w=\frac1N\Res_u
\left(\sum_{r,b}f_{b,p}(u)Y_{M_i}^g(b,u^N)w\right)
Nu^{N-1}\,du.
\end{equation}
Indeed, if $f_{b,p}(u)=\sum_{k\equiv r\,(\mathrm{mod}\,N)}f_{b,k}u^k$, then \eqref{eq:localcurrent25} equals $\sum_{r,b}\sum_{n\in r/N+\Z}f_{b,Nn}b_nw$; for each $w$ the sum is finite by finite pole order and lower truncation.
The factor $N^{-1}$ converts the residue of a pulled-back differential to its residue in the coarse coordinate $x$. At the output the contragredient convention is \eqref{eq:contragredientprelim}, with coarse coordinate $x=1/z$.
For an $H$-invariant current, the source residue over an orbit $H/I_p$ has $|H|/N$ equal contributions on the cover, each $N$ times
the descended residue; hence its average $|H|^{-1}\sum_{q\in H/I_p}
\Res_q$ equals the action \eqref{eq:localcurrent25}.
The finite frame, finite pole order and twisted lower truncation ensure that $R_i(\sigma)w$ is a finite sum for every algebraic vector $w$.

For the tensor product $T=\bigotimes_iM_i$ of the marked modules over a base $S$, put $\cT=\cO_S\otimes_{\C}T$. For a global admissible current
$\sigma$, let $R(\sigma)=\sum_iR_i(\sigma)$ on $T$ and set $\mathscr W=\sum_{\sigma}\operatorname{Im}R(\sigma)\subseteq\cT$. The coinvariants and Ward functionals are, respectively,
$\mathscr Q=\cT/\mathscr W$ and $\operatorname{Hom}_{\cO_S}(\mathscr Q,\cO_S)$; on a fiber they satisfy $\Phi(R(\sigma)w)=0$ for all $\sigma,w$.
These are the scalar conformal blocks in the present conventions. The canonical finite current bundle and allowed pole divisors are specified in Section~\ref{sec:currents}.

A ring is Noetherian if every ideal is finitely generated; equivalently, its ideals satisfy the ascending-chain condition. For a ring $R$, finite type means finite generation as an $R$-module. A coherent sheaf is locally the sheaf associated to a finitely presented module: such a module has a presentation $R^p\to R^q\to M\to0$ with finite $p,q$. A vector bundle is a locally free sheaf of finite rank. A finite projective module is a direct summand of a finite free module; over an affine variety this is the module of sections of a vector bundle. Put $Q^\vee=\Hom_R(Q,R)$. For such $Q$ and any module $B$, choose $q_i\in Q$ and $\phi_i\in Q^\vee$ with $w=\sum_i\phi_i(w)q_i$. Then $Q^\vee\otimes B\simeq\Hom_R(Q,B),\ \phi\otimes b\mapsto(w\mapsto\phi(w)b),$ with inverse $f\mapsto\sum_i\phi_i\otimes f(q_i)$.

An effective Cartier divisor $D$ is locally the zero locus of one non-zero-divisor. Write $F(D)=F\otimes\cO(D)$. For a locally free sheaf $F$ on a smooth proper curve $C$ and effective
Cartier divisor $D$, the sheaf $F(D)/F$ is supported on $D$ and records principal parts. If $Q$ is a marked section with ideal $I_Q$, the sheaf $F/I_Q^{J+1}F$ records its $J$-jets. The exact sequence $0\to F\to F(D)\to F(D)/F\to0$ induces
$H^0(C,F(D))\twoheadrightarrow H^0(C,F(D)/F)$ whenever $H^1(C,F)=0$. Moreover, $H^1(C,F)^*\simeq H^0(C,F^\vee\otimes\omega)$, and for a principal part $\eta$ and $s\in H^0(C,F^\vee\otimes\omega)$ the Serre pairing
with its boundary is $\sum_{p\in\operatorname{supp}D}
\Res_p\langle\eta_p,s_p\rangle$ \cite[Chapter~III, Section~7]{Hartshorne}. For a projective family $p:C\to S$ and relatively ample divisor $A$, relative Serre vanishing gives $R^1p_*F(mA)=0$ for sufficiently large $m$. A line bundle is relatively ample if, locally on the base, some positive power defines a closed embedding into a projective space over that base. Its corresponding divisor is called relatively ample. These are the cohomology and lifting facts used below
\cite{Hartshorne}.

Two base-change statements will be used. For a projective curve family over a Noetherian ring $R$, a coherent sheaf $F$ and a flat $R$-algebra $B$, flat base change gives $H^0(C,F)\otimes_RB\simeq H^0(C_B,F_B)$ \cite[Tag 02KH]{Stacks}. Inverting a DVR uniformizer is this statement with $B$ its fraction field. For a projective flat curve family $p:C\to S$ over a Noetherian
base and a coherent $S$-flat sheaf $F$, the fiberwise condition $H^1(C_s,F_s)=0$ on a neighborhood of $s_0$ yields $R^1p_*F=0$ there and $p_*F$ locally free, with $(p_*F)\otimes_{\cO_S}\kappa(s)\simeq H^0(C_s,F_s)$;
formation of $p_*F$ then commutes with base change on that neighborhood
\cite[Chapter~III, Section~12]{Hartshorne}. We apply this after fixing all current pole orders and twisting by a sufficiently large relatively ample incoming divisor. The relevant sheaves on the smooth component families are vector bundles, hence flat. If $R$ is a DVR with uniformizer $t$, an $R$-module $M$ is flat if and
only if $t:M\to M$ is injective. This applies to a reflexive sheaf on a normal nodal total space only after checking that the total space is an integral $R$-scheme and $t$ is a non-zero-divisor on the sheaf.

For a complex algebraic variety $X$, its analytification $X^{\rm an}$ has its complex analytic topology and structure sheaf; an algebraic coherent sheaf $F$ induces $F^{\rm an}$. On a projective $X$, analytification is an equivalence on coherent sheaves and preserves their cohomology and morphisms. This is the GAGA correspondence \cite{SerreGAGA}. In particular, a meromorphic finite-bundle section with poles bounded by a fixed divisor $D$ is an analytic section of $F(D)$, and is algebraic by this correspondence. This is the algebraization used for propagated sections on products of proper component covers.

We use the following \emph{finite-rank} strong residue theorem
\cite[Theorem~6.1]{FS}. Let $C$ be a smooth connected projective
curve, $P\subset C$ a nonempty finite set, and $F$ a finite-rank algebraic vector bundle. For $\eta_p\in(F^\vee\otimes\omega)\otimes_{\cO_{C,p}}\C((x_p))$ with a finite negative part in each chosen local parameter $x_p$,
\begin{equation}\label{eq:strongresidue25}
\begin{aligned}
&\exists\eta\in H^0(C\setminus P,F^\vee\otimes\omega),\quad
\widehat\eta_p=\eta_p\quad(p\in P)\\
&\hspace{9mm}\Longleftrightarrow\quad
\sum_{p\in P}\Res_p\langle\eta_p,s_p\rangle=0\quad
\bigl(\forall s\in H^0(C\setminus P,F)\bigr).
\end{aligned}
\end{equation}
The pairing uses the local expansions of $s$. A relative analytic version is given in \cite[Theorem~3.41]{GuiZhangPartial}. Equation~\eqref{eq:strongresidue25} is applied here only to fixed
finite-rank state-current bundles on smooth proper components; it does not by itself assert an infinite-rank or nodal extension. The formal-Laurent formulation in \eqref{eq:strongresidue25}
requires no a priori convergence. Fix the effective pole divisor $D$ from the negative parts; put $E_D=H^0(C,F^\vee\otimes\omega(D))$, so $\dim_\C E_D<\infty$. For $J\ge0$, the locus $A_J\subseteq E_D$ of sections agreeing with
$\eta_p$ modulo $x_p^{J+1}$ for every $p$ is affine. Serre duality identifies the obstruction to $A_J\ne\varnothing$ with residue pairings against meromorphic sections of $F$ with poles supported on
$P$; the right side of \eqref{eq:strongresidue25} annihilates them. Thus $A_J\ne\varnothing$, $A_{J+1}\subseteq A_J$ and, because $\dim E_D<\infty$, the descending affine chain stabilizes.
Any $\eta\in\bigcap_{J\ge0}A_J$ has all prescribed Laurent coefficients, proving algebraicity and convergence of these germs. For a finite group $H$ in characteristic zero, the invariant functor
is exact by the projector $|H|^{-1}\sum_{h\in H}h$.

A module over a domain is torsion free if $rw=0$ with $r\ne0$ implies $w=0$. For $F^\vee=\mathcal Hom(F,\cO)$, reflexivity means that $F\to F^{\vee\vee}$ is an isomorphism. On a normal variety a reflexive sheaf is torsion free and is determined by its restriction off a codimension-two subset. A Noetherian domain is normal if it is integrally closed in its fraction field. Depth is the maximal length of a sequence in the maximal ideal that successively acts by non-zero-divisors. A local ring is regular if the minimal number of generators of its maximal ideal equals its dimension. A Noetherian ring satisfies $R_1$ when its localizations at primes of height at most one are regular, and $S_2$ when every localization has depth at least $\min(2,\dim)$. Serre's criterion says that a reduced Noetherian ring is normal precisely when it satisfies $R_1$ and $S_2$. A local ring is Cohen--Macaulay when its depth equals its dimension; the corresponding condition on all localizations implies $S_2$. A hypersurface local ring is a quotient of a regular local ring by one non-zero-divisor. Such rings are Cohen--Macaulay. For $R=L[[t]]$ ($L$ a characteristic-zero field) and $a\ge1$,
$A=R[[u,v]]/(uv-t^a)$ is a two-dimensional complete local domain. As a hypersurface it is $S_2$; its singular locus is contained in the closed point $(t,u,v)$ (empty if $a=1$), so it is $R_1$.
Consequently $A$ is normal by Serre's criterion. In particular $t$
is a non-zero-divisor on $A$, a fact needed for the flatness argument. A regular parameter system is a minimal set of generators of that maximal ideal. For a nonzero element $f$, its order is the largest $n$ with $f$ in the $n$th power of the maximal ideal; in a regular local ring its initial term is a homogeneous polynomial in the parameter classes.

A discrete valuation ring (DVR) is a one-dimensional local principal ideal domain. If $t$ generates its maximal ideal, every nonzero fraction has the form $t^nu$ with $u$ a unit; its valuation is $\operatorname v_t(t^nu)=n$. The completed characteristic-zero rings used below are $O=L[[t]]$, with fraction field $K=L((t))$ and residue field $L$. A prime divisor is an irreducible codimension-one subvariety; on a normal variety its local ring is a DVR. The $t$-adic completion is the inverse limit over $O/(t^n)$. For $L=\C((s))$, the same notation describes the additional formal Laurent parameter. We use Hensel's simple-root lifting: a root of a polynomial modulo $t$ whose derivative is nonzero lifts uniquely to a root in the complete DVR with that residue. An excellent ring is a Noetherian ring with geometrically regular formal fibers, universally catenary dimension theory, and an open regular locus in every finite-type algebra. Universally catenary means that saturated prime chains with fixed endpoints have the same length in every finite-type algebra. Geometric regularity is regularity after field extension. Rings essentially of finite type over $\C$ and their complete local rings have these properties.

We use two local analytic extension facts. In the ring $\C\{z_1,\ldots,z_n\}$ of convergent power-series germs, if the zero set of an ideal $I$ is contained in the zero set of $Q$, the analytic Nullstellensatz gives $Q\in\sqrt I$, hence $Q^b\in I$ for some $b$. If $A(z):\cO^p\to\cO^q$ has rank $q$ on
$\{Q\ne0\}$, the ideal $I_q(A)$ of maximal minors has zero locus contained in $\{Q=0\}$; thus $Q^b\in I_q(A)$ for some $b\ge0$. Writing $Q^b$ as a combination of these minors and using the
adjugate identities yields $B(z)$ over $\cO[Q^{-1}]$ such that
$A(z)B(z)=I_q$ and $Q^bB(z)$ has holomorphic entries. On a normal complex space, holomorphic functions and sections of a locally free sheaf extend uniquely across analytic subsets of codimension at least two. Thus a section meromorphic off such a subset extends meromorphically if its poles there can first be cleared by one fixed local divisor. We use this Hartogs form only after establishing the finite divisor orders. In particular, on a root polydisk with poles confined to its axes, multiplying by $\prod_i s_i^{b_i}$ for sufficiently large $b_i$ gives a holomorphic section on the whole polydisk.

For a finite presentation $R^p\xrightarrow{A}R^q\to M\to0$, $\Fitt_j(M)=I_{q-j}(A)$, with $I_k(A)=R$ for $k\le0$ and $I_k(A)=0$ for $k>\min(p,q)$. The ideals are presentation-independent
and satisfy $\Fitt_j(M\otimes_RB)=\Fitt_j(M)B$ for any ring map $R\to B$. On an open subset $U\subseteq\operatorname{Spec} R$, $M|_U$ is locally free of rank $r$ precisely when
$\Fitt_r(M)|_U=\cO_U$ and $\Fitt_{r-1}(M)|_U=0$. These facts enter the connection argument in Lemma~\ref{lem:connfree}.

\subsection{Connections and formal expansions}
Let $U\subset S$ be a smooth complex algebraic coordinate chart, $R=\cO_S(U)$, and $\operatorname{Der}_{\C}(R)$ the $R$-module of $\C$-derivations of $R$. A connection on an $R$-module $M$ is a $\C$-linear assignment $\delta\mapsto\nabla_\delta\in\End_{\C}(M)$ with
$\nabla_\delta(fw)=\delta(f)w+f\nabla_\delta w$ and $\nabla_{f\delta}=f\nabla_\delta$. It is integrable if $[\nabla_\delta,\nabla_\gamma]=\nabla_{[\delta,\gamma]}$.
For smooth $S$, these data are equivalent to a left $\scrD_S$-module structure, where $\scrD_S$ is generated locally by $R$ and $\operatorname{Der}_{\C}(R)$ with $[\delta,f]=\delta(f)$. A section $w$ is horizontal if $\nabla_\delta w=0$ for every $\delta$. On $M^\vee=\Hom_R(M,R)$, the dual connection is
$\bigl(\nabla^\vee_\delta\Phi\bigr)(w)=\delta(\Phi(w))-\Phi(\nabla_\delta w)$; hence $\nabla^\vee\Phi=0$ is equivalent to $\delta(\Phi(w))=\Phi(\nabla_\delta w)$.
If $M$ is free of rank $q$ on $U$, write $\nabla=d+A$ in a frame $e$, with $A\in\operatorname{Mat}_{q}(\Omega^1_U)$. For column coordinates $v$, row coordinates $\Phi$, and a new frame $e'=eB$, $B\in\operatorname{GL}_q(R)$, $dv+Av=0,\; d\Phi-\Phi A=0,\; A'=B^{-1}AB+B^{-1}dB$.
An analytic integrable connection on a finite-rank holomorphic bundle has parallel transport along paths, invariant under endpoint-fixed homotopy \cite[Chapter~I, Section~2]{Deligne}.

For the finite marking index set $I$, let $T=\bigotimes_{i\in I}M_i$ and $\cT=R\otimes_{\C}T$ as in Section~\ref{sec:geometricconventions}; write $z_i\in R$ for the moving coarse insertion coordinates and set $D_i=L(-1)^{(i)}$ on $T$. In fixed physical-coordinate frames define
\begin{equation}\label{eq:Tconn}
 \nabla_\delta(fw)=\delta(f)w+
 f\sum_{i\in I_{\mathrm{mov}}}\delta(z_i)D_iw.
\end{equation}
Here $I_{\mathrm{mov}}\subseteq I$ consists of the moving inputs; fixed coordinates satisfy $\delta(z_i)=0$. Since $[D_i,D_j]=0$ and $\delta\gamma(z_i)-\gamma\delta(z_i)=[\delta,\gamma](z_i)$, we have $[\nabla_\delta,\nabla_\gamma]=\nabla_{[\delta,\gamma]}$. The use of $D_i$ agrees with the $L(-1)$-derivative property of twisted intertwining operators \cite[Definition~2.8]{DH}. For any $R$-submodule $J\subseteq\cT$, put
\[
 J^\nabla=
 \sum_{m\geq0}\ \sum_{\delta_1,\ldots,\delta_m\in\operatorname{Der}_{\C}(R)}
 R\,\nabla_{\delta_1}\cdots\nabla_{\delta_m}J,
 \qquad \nabla( J^\nabla)\subseteq J^\nabla.
\]
It is the least connection-stable $R$-submodule containing $J$, and $\nabla$ induces a connection on $\cT/J^\nabla$. The canonical differential Ward quotient is specified in Definition~\ref{def:diff}; no finite-rank assumption is imposed on $\cT$.

Let $L\supseteq\C$ be a field, $O=L[[t]]$, $K=L((t))$, and let $Q_K$ be a finite-dimensional $K$-vector space with a connection for $t\partial_t$ ($\partial_t|_L=0$). A full lattice is an $O$-submodule $\Lambda_t\subset Q_K$ satisfying
$\Lambda_t\simeq O^{\dim_K Q_K}$ and $K\Lambda_t=Q_K$. Regular singularity at $t=0$ is equivalent to
\[
 \exists\Lambda_t:\quad
 \nabla_{t\partial_t}\Lambda_t\subseteq\Lambda_t,
 \qquad \operatorname{Res}_{\Lambda_t}(\nabla)
 =\nabla_{t\partial_t}\big|_{\Lambda_t/t\Lambda_t}
 \in\End_L(\Lambda_t/t\Lambda_t).
\]
The residue depends on the lattice; its eigenvalues modulo $\Z$ do not \cite[Chapter~II, Sections~1 and~3]{Deligne}. Let $X$ be smooth over $\C$ and $D=\bigcup_{i=1}^rD_i$ a simple normal-crossing divisor, locally $D=\{\xi_1\cdots\xi_r=0\}$ with coordinates $(\xi_1,\ldots,\xi_r,y_1,\ldots,y_s)$. Then
$\Omega^1_X(\log D)=\sum_i\cO_X\,d\xi_i/\xi_i+\sum_j\cO_X\,dy_j$ locally. A logarithmic extension of a finite-rank integrable connection $(Q,\nabla)$ on $X\setminus D$ is a locally free $\overline Q$ on $X$ satisfying
$\overline Q|_{X\setminus D}=Q$ and $\nabla\overline Q\subseteq\overline Q\otimes\Omega^1_X(\log D)$. For a regular integrable connection, Deligne's canonical logarithmic extension is characterized by
\[
 \operatorname{Spec}\!\left(\operatorname{Res}_{D_i}\nabla\right)
 \subseteq\{\lambda\in\C:0\leq\operatorname{Re}\lambda<1\},
 \qquad 1\leq i\leq r,
\]
where the strip selects one representative of each exponent class in $\C/\Z$ \cite[Chapter~II, Sections~5.3--5.6]{Deligne}. Any two logarithmic extensions are identified after allowing finite-order meromorphic poles along $D$.

For the global regularity criterion, let $U$ be a smooth complex algebraic variety, $(Q,\nabla)$ a \emph{finite-rank algebraic integrable} connection on $U$, and $\overline U$ a normal proper compactification. Then regularity along each codimension-one irreducible component of $\overline U\setminus U$ implies regularity after pullback to every smooth curve mapping to $U$, and hence at all divisorial valuations; in particular a logarithmic extension exists after choosing a smooth normal-crossing compactification \cite[Theorem~1.1.1, (ii)$\Rightarrow$(iii)]{Andre}\cite[Chapter~II, Sections~4--5]{Deligne}. This criterion is used only after the finite quotient has been proved locally free with integrable connection (Lemma~\ref{lem:connfree} and Theorem~\ref{thm:finiteconn}).

For $\widehat A=\C[[\xi_1,\ldots,\xi_r,y_1,\ldots,y_s]]$, a finite set $\Lambda\subset\C^r$, and $b\in\Z_{\geq0}$, let $\xi^\lambda=\prod_i\xi_i^{\lambda_i}$ and $(\log\xi)^\beta=\prod_i(\log\xi_i)^{\beta_i}$ be formal symbols. Set
\begin{equation}\label{eq:Nilsson}
 \mathscr N_{\Lambda,b}
 =\sum_{\lambda\in\Lambda}\ \sum_{0\leq\beta_i\leq b}
 \xi^\lambda(\log\xi)^\beta
 \widehat A[(\xi_1\cdots\xi_r)^{-1}].
\end{equation}
Here $\widehat A[(\xi_1\cdots\xi_r)^{-1}]=\bigcup_{n\geq0}(\xi_1\cdots\xi_r)^{-n}\widehat A$; thus every element has finitely many exponent classes modulo $\Z^r$, bounded logarithmic degree, and a finite pole bound in every $\xi_i$. Formally,
$\partial_{\xi_i}\xi^\lambda=\lambda_i\xi^{\lambda-e_i}$ and $\partial_{\xi_i}\log\xi_j=\delta_{ij}\xi_i^{-1}$, where $e_i\in\Z^r$ is the $i$th standard basis vector; $\partial_{y_j}$ acts coefficientwise on $\widehat A$.
For $O=L[[t]]$ and $\lambda\in\C$, the formal monomial line is $Oe_\lambda$ with
$\nabla_{t\partial_t}(f e_\lambda)=(t\partial_tf+\lambda f)e_\lambda$; $e_\lambda$ is denoted $t^\lambda$ in collision expansions. Finite direct sums retain the corresponding complex exponents.

\section{Canonical currents and the differential Ward quotient}\label{sec:currents}\label{sec:control}
\subsection{Finite state spaces, spanning words and physical currents}
For each $g\in G$, fix a set $B_g\subset V_{>0}=\bigoplus_{d>0}V_d$ of homogeneous $g$-eigenvectors whose images form a basis of $(V/C_2(V))_{>0}$. This is possible because $C_2(V)$ is graded and $g$-stable; $\dim(V/C_2(V))<\infty$ and $V_0=\C\mathbf1$. Put
\begin{equation}\label{eq:E}
 B=\{\mathbf1,\omega\}\cup\bigcup_{g\in G}\bigcup_{k\in G}kB_g,
 \qquad E=\Span\{L(j_1)\cdots L(j_s)b:b\in B,\ j_i\in\Z_{>0},\ s\in\N\}.
\end{equation}
The sets $B_g$ are used in the $g$-twisted spanning theorem below; the larger space $E$ is used for coordinate-stable currents. In particular, $kB_g$ is not asserted to consist of $g$-eigenvectors.

\begin{lemma}\label{lem:Efinite}
$E\subset V$ is finite dimensional, graded, $G$-stable, and $L(j)E\subseteq E$ for $j>0$. The $\Aut\cO$-action on $E$ factors through a finite jet quotient.
\end{lemma}
\begin{proof}
Take $d=\max\{\wt(b):b\in B\}$. Since $[L(0),L(j)]=-jL(j)$, 
\[
 E\subseteq\bigoplus_{q=0}^dV_q,
 \qquad L(j)|_E=0\quad(j>d).
\]
The finite-dimensionality follows from $\dim V_q<\infty$; the $L(j)$-stability follows directly from the definition. For $k\in G$, $kL(j)=L(j)k$ and $kB=B$, hence $kE=E$. On $E$, $L(0)$ is semisimple with integral eigenvalues, and the Lie algebra generated by $L(j)$, $j>0$, acts nilpotently and factors through $L(1),\ldots,L(d)$. Thus the origin-fixing coordinate action factors through the quotient of $\Aut\cO$ by the subgroup $\{\rho:\rho(t)\equiv t\pmod{t^{d+2}}\}$.
\end{proof}

In particular, $L(2)\omega=(c_V/2)\mathbf1\in E$. The associated stress-tensor term is evaluated in Section~\ref{app:stress}; $L(-1)E\subseteq E$ is neither asserted nor needed.

We apply \cite[Proposition~2.12]{Tan} to the positive-weight representative set $B_g$, not to the $G$-stable union $B$. Choose $M_\PBW\in\Z_{>0}$ such that
\[
 \bigoplus_{d\ge M_\PBW}V_d\subseteq C_2(V).
\]
For an ordinary $g$-twisted module $M$ generated by finitely many nonzero homogeneous vectors $u\in U_M$, choose $L_u\in\Z$ such that $b_{\alpha_b-n}u=0\; (b\in B_g,\ n\le L_u),\; gb=e^{2\pi i\alpha_b}b,\; 0\le\alpha_b<1$.
For $N\in\Z$ with $N\ge\max_{u\in U_M}L_u$, let $\mathcal W_N(M)$ consist of the vectors
\begin{equation}\label{eq:PBW}
 b^1_{\alpha_1-n_1}\cdots b^s_{\alpha_s-n_s}u,
 \qquad u\in U_M,\quad b^i\in B_g,\quad
 n_1\ge\cdots\ge n_s>L_u,\quad s\in\N,
\end{equation}
with the additional condition $\#\{i:n_i=n\}<M_\PBW\; \text{for each }n\le N$.
Here $\alpha_i=\alpha_{b^i}$. Tan's proposition asserts $M=\Span\mathcal W_N(M)$ for every such $N$. The repetition bound $M_\PBW-1$ is independent of $N$ and $u$. Homogeneous $u$ and $b^i$ make every word homogeneous, since for $b\in V_d$ and $v\in M_\lambda$,
$b_{\alpha_b-n}v\in M_{\lambda+d-\alpha_b+n-1}$. Every ordinary simple twisted module is generated by any nonzero homogeneous vector; finite direct sums therefore admit finite sets $U_M$.

\begin{lemma}\label{lem:lowdepth}
Let $R,N\in\Z$ satisfy $N\ge\max(R-1,\max_uL_u)$. Put $K_{R,N}(M)=\Span\{w\in\mathcal W_N(M):s=0\ \text{or}\ n_1<R\}$.
Then $\dim K_{R,N}(M)<\infty$. For finitely many ordinary twisted modules $M_i$, one can choose a common $N$ and a finite set of homogeneous tensor generators of $\bigotimes_iK_{R,N}(M_i)$.
\end{lemma}
\begin{proof}
For a nonempty word in $K_{R,N}(M)$ with terminal generator $u$, each $n_i$ belongs to the finite set $\{L_u+1,\ldots,R-1\}\subseteq(-\infty,N]$. Thus $s\le(M_\PBW-1)\max\{0,R-L_u-1\}$.
The sets $B_g,U_M$ are finite. Hence there are finitely many such words, and the assertion about a finite tensor product follows.
\end{proof}

\begin{lemma}\label{lem:descent}
Let $T=\bigotimes_{i=0}^mM_i$, with the $M_i$ ordinary simple twisted modules (including a restricted contragredient at $i=0$), $A$ a commutative $\C$-algebra, and $J\subseteq A\otimes_\C T$ an $A$-submodule. Fix $R\in\Z$ and $N\ge\max(R-1,\max_{i,u}L_u)$. Assume that every homogeneous tensor $w$ of words in the fixed families $\mathcal W_N(M_i)$ containing a mode of depth $\ge R$ satisfies
\begin{equation}\label{eq:descent-condition}
 1\otimes w\equiv\sum_{j=1}^{q}a_j\otimes v_j\pmod J,
 \qquad a_j\in A,\quad v_j\in T\ \text{homogeneous},\quad
 \ee(v_j)<\ee(w).
\end{equation}
Then $ (A\otimes_\C T)/J$ is generated as an $A$-module by the images of $\bigotimes_{i=0}^mK_{R,N}(M_i)$.
\end{lemma}
\begin{proof}
For $e\in\N$, induct on $e$ for all homogeneous tensors with $\ee=e$. Every homogeneous vector of $T$ is a finite linear combination of homogeneous tensors of words in the fixed $\mathcal W_N(M_i)$, since each spanning identity respects the conformal-weight grading. A tensor with no depth $\ge R$ belongs to $\bigotimes_iK_{R,N}(M_i)$. Otherwise apply \eqref{eq:descent-condition}. Each $v_j$ has energy $<e$, so the induction hypothesis applies to its class modulo $J$; after rewriting $v_j$ in the same families, every nonzero homogeneous term still has energy $\ee(v_j)$. This proves the claim by well-founded induction on $\ee\in\N$. The same argument applies to each of the finitely many simple external blocks of Proposition~\ref{prop:regularmodules}.
\end{proof}

Let $\mathfrak B=F_m(\C^*)=\{(p_1,\ldots,p_m)\in(\C^*)^m:p_i\ne p_j\ (i\ne j)\}$ and $\mathfrak U=F_{m+1}(\C^*)$, where the final coordinate $z$ varies in $\mathbb P^1\setminus\{0,p_1,\ldots,p_m,\infty\}$. Prescribe the local labels $g_i=\rho(\gamma_{p_i})$, $g_0=\rho(\gamma_0)$ and $g_\infty=\rho(\gamma_\infty)$, with
\[
 g_1\cdots g_mg_0g_\infty=1,
 \qquad \langle g_1,\ldots,g_m,g_0\rangle=H\le G,
\]
in the ordered meridian convention \eqref{eq:orderedmeridians}. We work over a smooth affine open of $\mathfrak B$ or its finite \'etale Hurwitz cover, retaining the physical coordinates $z,z-p_i,s=z^{-1}$.

\begin{lemma}\label{lem:framedcharts}
For these data there exists a smooth finite \'etale chart $S\to\mathfrak B$ carrying a smooth proper connected family of $H$-Galois covers of $\mathbb P^1$, together with a chosen lift of each marking and a formal parameter $u_i$ at that lift satisfying $x_i=u_i^{N_i}$, where $N_i=\operatorname{ord}(g_i)$ and $x_i\in\{z,z-p_j,z^{-1}\}$. An unramified based lift is retained as a separate part of the framing. Homotopic base paths relative to endpoints, with the same initial full framing, have identical final framed points of $S$.
\end{lemma}
\begin{proof}
Disjoint-disk isotopies of the marked points give local trivializations of the forgetful map $\mathfrak U\to\mathfrak B$, whose fiber is $\C\setminus\{0,p_1,\ldots,p_m\}$. The fiber and, inductively, $\mathfrak B$ are aspherical. Thus the homotopy sequence gives
\[
 1\longrightarrow\Pi\longrightarrow\pi_1(\mathfrak U)
 \longrightarrow\pi_1(\mathfrak B)\longrightarrow1,
 \qquad\Pi=\pi_1(\C\setminus\{0,p_1,\ldots,p_m\}).
\]
The continuous section $z(p)=2+\sum_i|p_i|$ satisfies $|z(p)|>|p_i|$ and splits this sequence. Fix the resulting semidirect product $\pi_1(\mathfrak U)=\Pi\rtimes\pi_1(\mathfrak B)$. The meridians define a surjection $\rho:\Pi\twoheadrightarrow H$ with the stated labels. Put $\Gamma=\{\gamma\in\pi_1(\mathfrak B):\rho\circ\operatorname{Ad}_{\gamma}^{-1}=\rho\}$.
Since $H$ is finite and $\Pi$ is finitely generated, $[\pi_1(\mathfrak B):\Gamma]<\infty$. The connected topological covering $\mathfrak B_\rho\to\mathfrak B$ corresponding to $\Gamma$ is algebraic finite \'etale by \cite[Expos\'e XII, Th\'eor\`eme~5.1]{SGA1}. Let $\mathfrak U_\rho=\mathfrak U\times_{\mathfrak B}\mathfrak B_\rho$. The homomorphism
\[
 \widehat\rho:\pi_1(\mathfrak U_\rho)=\Pi\rtimes\Gamma\longrightarrow H,
 \qquad(\eta,\gamma)\longmapsto\rho(\eta)
\]
is well defined because $\Gamma$ fixes $\rho$. Its kernel defines a connected principal $H$-cover of $\mathfrak U_\rho$, algebraic finite \'etale by the same theorem; the $H$-action algebraizes by full faithfulness.

Normalize $\mathbb P^1\times\mathfrak B_\rho$ in this finite cover. The normalization is finite (the base is excellent) and hence proper and projective over $\mathfrak B_\rho$. Locally analytically near each marked divisor $x_i=0$, the punctured covering is a disjoint union of models $x_i=u_i^{N_i},\; (u_i,t)\in\mathbb D\times T$,
with $T$ a simply connected analytic chart of $\mathfrak B_\rho$. Each normalization is $\mathbb D\times T$, so the normalized family is analytically smooth over $T$. The algebraic family is smooth by \cite[Expos\'e XII, Proposition~3.1(iv)]{SGA1}, and its reduced marking preimages are finite \'etale over $\mathfrak B_\rho$. Uniqueness of normal finite extension follows from \cite[Expos\'e XII, Proposition~5.3 and Th\'eor\`eme~5.4]{SGA1}.

Pass to the finite \'etale product of the reduced marking preimages to choose distinguished lifts $Q_i$. Set
\[
 \mathcal L_i=I_{Q_i}/I_{Q_i}^2,\qquad
 \delta_i=(\pi^*x_i)\bmod I_{Q_i}^{N_i+1}
 \in\Gamma(S,\mathcal L_i^{\otimes N_i}).
\]
Because the ramification index is $N_i$, $\delta_i$ is nowhere zero. The functor of $N_i$th roots $\ell_i^{\otimes N_i}=\delta_i$ is a finite \'etale $\mu_{N_i}$-torsor: in a local trivialization it has equation $a^{N_i}=c$ with $c\in\cO_S^\times$. After choosing $\ell_i$, take a local parameter $v_i$ with leading conormal class $\ell_i$. There is a unique formal unit $a_i(v_i)$ with $a_i(0)=1$ such that $\pi^*x_i=v_i^{N_i}a_i(v_i)$, and $u_i=v_i a_i(v_i)^{1/N_i},\; u_i\equiv v_i\pmod{v_i^2},\; u_i^{N_i}=\pi^*x_i$.
Uniqueness makes the formal roots compatible on overlaps. These finitely many torsors give the required finite \'etale chart $S\to\mathfrak B$.

The splitting section $z(p)$ retains an unramified based lift: the subgroup $\Gamma$ maps trivially under $\widehat\rho$, so it admits a continuous lift over $\mathfrak B_\rho$. This based lift is additional to the distinguished ramified lifts and roots; their mutual identification requires the initially chosen reference paths and is not inferred from a local root alone. On a connected component of $S$, the ordered inertia generators are locally constant as elements of $H$ (rather than only up to conjugacy), once the based paths are transported. Homotopy lifting for the finite covering $S\to\mathfrak B$ gives the final assertion.
\end{proof}

\begin{definition}\label{def:commonoutsetframing}
For each ordered triple of external sectors, fix the following \emph{single smooth outset framing} before choosing any pair of intertwining operators or any homogeneous coefficient.  In the positive real chamber
$\mathfrak P^+=\{(a,b)\in\mathbb R^2:a>b>0\}$ take $0<\rho_0<1/16$ and
\[
 p_*=(a_*,b_*)=\bigl(1+\rho_0/2,(1+\rho_0)/2\bigr),
 \qquad b(p)=2+a+b.
\]
On the connected framed Hurwitz component determined by the ordered monodromy, choose once a lift of the unramified point $b(p_*)$, the full right-direct/lower-comb homotopy classes of reference paths to $a_*,b_*,0,\infty$, and the corresponding distinguished marking lifts and roots of the centered coarse parameters.  Continue all these data together by path lifting; in particular the small loops have ordered words
$\gamma_A,\gamma_B,\gamma_C$ and the cluster stems have words $\gamma_B\gamma_C$ (product) and $\gamma_A\gamma_B$ (iterate).  On the positive rays use real logarithms, while the binary Du--Huang single-operator germs always retain their original preferred continuations and $\log1=0$.  These choices depend on the sector labels and the framed component, not on an intertwining operator or an internal grade.

Specifically, continue the same initial full framing along
\[
 (a(s),b(s))=\bigl(1+s\rho_0,(1-s)\rho_0+s\bigr),\quad 0\le s\le1.
\]
Both coordinates are positive and $a(s)-b(s)>0$; hence the corresponding lower-comb paths are isotopic in the ordered configuration space.  Their lifts define the product framing at $(1,\rho_0)$ with basepoint $3+\rho_0$, and the iterate framing at $(1+\rho_0,1)$ with basepoint $4+\rho_0$.  At either boundary take the actual admissible degeneration of this transported generic cover and inherit its external roots, based stems and logarithms; its two node roots must be balanced and are not separately reselected.  A finite ramified pullback used to express the node is always taken with the already chosen positive logarithm.  Other connected framed components admit the same type of independent outset choice.  For each binary degeneration, the auxiliary rooted three-point
reference covers used in Lemma~\ref{lem:singleinsertion} are obtained
jointly by restriction of this one outset cover, rather than by
two additional independent choices of a node based lift.  More
precisely, retain one lifted point of the punctured balanced collar;
its outside and bubble restrictions determine the two pointed component
covers and their common node stem.  The two reference node roots are
the restrictions of this one balanced rooted collar (expressed in the
standard three-point coordinates), and the reference node module on
the bubble is the restricted contragredient of the outside module,
with the strict deck pullback understood on both sides.  This is an
auxiliary \emph{geometric} choice made before selecting the pair of
intertwining operators: it does not select, renormalize or alter any
intertwining operator or its prescribed $z=1$, $\log1=0$ branch.
The existence and the comparison with the prescribed preferred
single-operator germs are checked in
Lemma~\ref{lem:E0actualsourcecomparison}; in particular the above
prescription is not an extra compatibility hypothesis on those germs.
A root/deck change or a braid later on is a \emph{transport of these
data}, not a declaration that the original coefficient formula
remains numerically unchanged.
\end{definition}

\begin{definition}\label{def:fiveoutsetframing}
For a fixed ordered four-input problem (the output is at infinity), use the
same smooth five-point Hurwitz cover before choosing a binary tree,
intertwining operators, or homogeneous grades.  On the positive chamber
\(\mathfrak P_5^+=\{z_A>z_B>z_C>z_E=0\}\), choose the point
\(p_5=(7,6,4,0)\), the unramified point
\(b_5(p)=2+z_A+z_B+z_C\), one lift of it, and the complete
right-direct/lower-comb paths to the four finite marks and infinity.  Lift
these paths from that one point, and fix the resulting marked lifts and a
special root of each physical centered local parameter (and of \(Z^{-1}\)
at infinity).  All these choices, and the positive real logarithms of the
configuration differences, are continued together along paths in
\(\mathfrak P_5^+\).  No new external root or based lift is selected at a
cherry.  The two node roots and their paired node module frames are
inherited from the \emph{one actual} balanced degeneration along its retained
lower-comb cluster stem.  At a boundary they are compared, rather than
silently equated, with the auxiliary three-input representatives of
Definition~\ref{def:commonoutsetframing}.

Three-input outset framings are reference representatives, not a requirement
that the same triple representative be literally the restriction of
every four-input extension.  A reference change acts on the modules, typed
twists, coefficients and Ward currents simultaneously by
Lemma~\ref{lem:relativecoordinates} and Proposition~\ref{prop:canonicalrealization}.
In particular a change of reference does not replace the original
intertwining operator \(\mathcal Y(\cdot,1)\cdot\), its \(\log1=0\) branch,
or an external algebraic vector by a newly defined operator or vector.
\end{definition}

\begin{lemma}\label{lem:jets}
Let $p:C\to S$ be a smooth projective family of curves with an action of a finite group $H$ over $S$, and let $\cF$ be an $H$-equivariant vector bundle. Let $A,Q$ be disjoint, $H$-invariant relative effective Cartier divisors, with $\cO_C(A)$ relatively ample. For every fixed $K\in\N$, locally on $S$ there is $m\in\N$, independent of $n\in\N$, for which
\begin{equation}\label{eq:jetmap}
 \bigl(p_*\cF(mA+nQ)\bigr)^H\longrightarrow
 \bigl(p_*(\cF(mA+nQ)/\cF(mA-KQ))\bigr)^H
\end{equation}
is surjective for all $n\in\N$.
\end{lemma}
\begin{proof}
Relative Serre vanishing for the projective morphism $p$ and the $p$-ample line bundle $\cO_C(A)$ gives, after replacing $S$ by a neighborhood if needed, an integer $m$ such that $R^1p_*\cF(mA-KQ)=0$.
For each $n\in\N$, the sequence $0\longrightarrow\cF(mA-KQ)\longrightarrow\cF(mA+nQ) \longrightarrow\cF(mA+nQ)/\cF(mA-KQ)\longrightarrow0$
is exact because $n+K\ge0$. Pushing forward and using the same vanishing gives surjectivity before taking $H$-invariants, simultaneously for every $n$. Since $|H|\in\C^\times$, the projector $|H|^{-1}\sum_{h\in H}h$ makes invariants exact, proving \eqref{eq:jetmap}. The quotient is supported on the relative divisor $(n+K)Q$; since $Q$ is relative Cartier and $\cF$ locally free, its direct image is locally free over $S$ (after shrinking if necessary). Hence \eqref{eq:jetmap} lifts arbitrary \emph{$H$-invariant} sections of the indicated quotient locally on $S$.
\end{proof}

At a sink marking of ramification order $e$, the reduced $H$-orbit divisor $A$ satisfies $eA=\pi^*(\{\text{sink}\}\times S)$, hence is relatively ample. For several target orbits one replaces $KQ$ by the fixed sum of nontarget jet divisors; the cohomology kernel remains independent of target pole depth.

\subsection{Ramified coordinate covariance and residue normalization}
Put $E_{\le d}=E\cap\bigoplus_{0\le j\le d}V_j$, $E_{<d}=E_{\le d-1}$, and $\operatorname{gr}_dE=E_{\le d}/E_{<d}$. Since $L(j)E_{\le d}\subseteq E_{\le d-j}$ for $j>0$, the coordinate action on $\operatorname{gr}_dE$ is $a_1^{-d}\id$ for $f(x)=a_1x+O(x^2)$. Hence the weight-$d$ associated graded of the current bundle is
\[
 \operatorname{gr}_d(\mathscr V_E\otimes\omega_{C/S})
 \simeq(\operatorname{gr}_dE)\otimes_{\C}\omega_{C/S}^{\otimes(1-d)}.
\]
For $z=t^N$, the induced tensor-density substitution on the punctured disk is $z^r(dz)^{1-d}=N^{1-d}t^{Nr+(N-1)(1-d)}(dt)^{1-d}$.
This is the associated graded formula, not the transformation law for a general $b\in E_{\le d}$: the latter also contains $L(j)b$ ($j>0$), as in the coordinate-change representation of \cite[Sections~3--4]{FS}.

\begin{lemma}\label{lem:gauge}
Let $N\ge1$. Over $\C[t,t^{-1}]$, the centered change from the cover coordinate $v$ to the coarse coordinate $w$, $w=f_t(v):=(t+v)^N-t^N$,
induces on $E$ a matrix $A_N(t)\in\operatorname{GL}(E\otimes_\C\C[t,t^{-1}])$. Both $A_N(t)$ and $A_N(t)^{-1}$ have finite Laurent entries and pole bounds independent of a current's mode index.
\end{lemma}
\begin{proof}
Write $f_t(v)=a_1(t)v+\sum_{k=2}^Na_k(t)v^k$, where
\[
 a_1(t)=Nt^{N-1},\qquad
 \frac{a_k(t)}{a_1(t)}=\frac1N\binom Nk t^{1-k}\quad(2\le k\le N).
\]
Thus the coefficients of every finite jet of $f_t$ and the inverse of its linear coefficient belong to $\C[t,t^{-1}]$. By Lemma~\ref{lem:Efinite}, $R_V|_E$ is an algebraic representation of a finite jet group; its scaling weights are integers and its positive-mode operators are nilpotent. Therefore $A_N(t)$ has entries in $\C[t,t^{-1}]$. The unipotent factor has determinant $1$ and the scaling factor has determinant a signed integral power of $a_1(t)$; thus $\det A_N(t)\in\C^\times t^\Z$. The adjugate formula gives $A_N(t)^{-1}\in\operatorname{Mat}(\dim E,\C[t,t^{-1}])$. Since $E$ is fixed, the finitely many Laurent valuations of these matrices have a common lower bound independent of the mode index.
\end{proof}

The matrix of Lemma~\ref{lem:gauge} compares frames only where $t\ne0$: the ramified map $z=t^N$ is not an invertible coordinate change at $t=0$ when $N>1$. If $K\ge0$ bounds the negative Laurent powers of the inverse state-current gauge, then, after the prescribed cover principal part is fixed, the regular cover coefficients through degree $K$ determine the negative Laurent part in the coarse frame. The required additional regular jets may be imposed using Lemma~\ref{lem:jets}. The bound $K$ is independent of the prescribed mode depth.

The state coordinate action uses the integral weights of $V$. For an ordinary module, Lemma~\ref{lem:relativecoordinates} gives an origin-fixing action using relative weights $L(0)-h_M$. The placement connection is defined in fixed physical module frames by \eqref{eq:Tconn}.

\begin{lemma}\label{lem:canonicalfinite}
Let $p:C\to S$ be a smooth family of curves with an $H$-action, and let $E$ be the full state space of \eqref{eq:E}. Then
\begin{equation}\label{eq:canonicalfinite}
F_E=\bigl(\operatorname{Coord}(C/S)
       \times^{\Aut\cO}E\bigr)\otimes\omega_{C/S}
\end{equation}
is an algebraic $H$-equivariant vector bundle of rank $\dim_\C E$. At a marking with $x=u^N$, the two frames on $u\ne0$ are related by a matrix in $\operatorname{GL}(E\otimes\C[u,u^{-1}])$ after incorporating the differential Jacobian; its Laurent pole bounds are independent of mode depth. These statements also hold on the relative smooth locus of a nodal curve over a $\C$-algebra DVR.
\end{lemma}
\begin{proof}
By Lemma~\ref{lem:Efinite}, $E$ is a finite-dimensional algebraic representation of a finite jet quotient of $\Aut\cO$. The associated bundle $\mathscr V_E$ was defined in Section~\ref{sec:geometricconventions}; it is locally free of rank $\dim E$. Its tensor product with $\omega_{C/S}$ has the same rank. Since $hE=E$ and $hL(j)=L(j)h$ for all $h\in H$, the deck action commutes with the coordinate action and defines an $H$-linearization.

At $u_0\ne0$, set $w=(u_0+v)^N-u_0^N$. Lemma~\ref{lem:gauge} supplies the state matrix and its Laurent inverse. On $u\ne0$ the differential transition is multiplication by $dx/du=Nu^{N-1}\in\C[u,u^{-1}]^\times$. Hence the full current matrix and its inverse are Laurent. The construction is polynomial in finite jets and in $N^{\pm1}$; it therefore commutes with base change over any $\C$-algebra, including a $\C$-algebra DVR, provided the relative curve is smooth at the observation point.
\end{proof}

\begin{lemma}\label{lem:relativecoordinates}
Let $M$ be an ordinary simple $g$-twisted module, $N=\operatorname{ord}(g)$, and $M=\bigoplus_{a\in N^{-1}\N}M_{h+a}$, where $h=h_M$ is as defined in Section~\ref{sec:moduleprelim}. Put $x=u^N$ and
\[
 \Aut_N\cO=\{\rho\in\Aut\C[[u]]:\rho(\C[[u^N]])=\C[[u^N]]\}
 =\{\rho:\rho(u)=u\,a(u^N),\ a(t)\in\C[[t]]^\times\}.
\]
It is the group of ring automorphisms (opposite to composition of the displayed point maps), with an exact sequence $1\to\mu_N\to\Aut_N\cO\xrightarrow{\rho\mapsto\rho(u)^N}\Aut\C[[x]]\to1$; the target uses $x=u^N$. There is an action $R_{M,h}$ of $\Aut_N\cO$ on $M$, algebraic on each finite-dimensional relative-weight cutoff, whose infinitesimal action is
\begin{equation}\label{eq:relativeautlie}
u^{Nk+1}\partial_u\longmapsto
 \begin{cases}
-N(L_M(0)-h\id),&k=0,\\ -NL_M(k),&k>0.
 \end{cases}
\end{equation}
In particular $R_{M,h}(a u)=a^{-N(L_M(0)-h\id)}$ for $a\in\C^\times$. If $f(x)=\rho(u)^N$ and $f_x(v)=f(x+v)-f(x)$, then
\begin{equation}\label{eq:relativefieldcov}
R_{M,h}(\rho) Y_M\bigl(R_V(f_x)^{-1}b,f(x)\bigr) R_{M,h}(\rho)^{-1}=Y_M(b,x).
\end{equation}
All fractional powers in the left-hand field are expanded using the chosen root $\rho(u)$ of $f(x)$. Here $R_V$ is the state action of Section~\ref{sec:geometricconventions}, with no lowest-weight shift. The identity holds coefficientwise in $\operatorname{End}(M)((u))$ on each fixed input, with the usual lower-truncation interpretation. The same construction applies blockwise to finite direct sums of ordinary simples. This is the lowest-weight-shifted version of the coordinate covariance in \cite[Section~4]{FS}.
\end{lemma}
\begin{proof}
The first displayed description of $\Aut_N\cO$ agrees with the special-coordinate group of \cite[Section~3]{FS} and follows from $\rho(u)^N\in\C[[u^N]]$: the unique $N$th root of $\rho(u)^N/u^N$ with prescribed nonzero constant term belongs to $\C[[u^N]]$. Its Lie algebra is $\prod_{k\ge0}\C v_k$, $v_k=u^{Nk+1}\partial_u$, with $[v_k,v_l]=N(l-k)v_{k+l}$. For $k,l\ge0$ the Virasoro central term vanishes. Replacing $L_M(0)$ by $L_M(0)-h\id$ preserves all brackets; hence \eqref{eq:relativeautlie} is a Lie algebra representation.

For $q\in N^{-1}\N$, put $M^{(\le q)}=\bigoplus_{0\le a\le q}M_{h+a}$. Then $\dim M^{(\le q)}<\infty$ and $L_M(k)M^{(\le q)}\subseteq M^{(\le q)}$ for $k\ge0$. Moreover $L_M(k)|_{M^{(\le q)}}=0$ for $k>q$, and every product of more than $\lfloor q\rfloor+1$ strictly positive $L_M(k)$ vanishes on $M^{(\le q)}$. Thus the positive Lie algebra integrates through a finite unipotent jet group. Since $N(L_M(0)-h\id)$ has eigenvalues in $\N$, the scaling operators $a^{-N(L_M(0)-h\id)}$ give an algebraic $\C^\times$-action. The brackets $[L_M(0),L_M(k)]=-kL_M(k)$ verify the semidirect-product relation. The compatible actions on all $M^{(\le q)}$ define $R_{M,h}$; no condition $Nh\in\Z$ is required. The shift is essential when $h$ is not in $N^{-1}\Z$.

We verify the field identity by the module Virasoro commutator. The twisted module Virasoro commutator, for $k\ge0$, is
\begin{align*}
[L(k),Y_M(b,x)] &=x^{k+1}\partial_xY_M(b,x)\\ &\quad+\sum_{j=0}^k\binom{k+1}{j+1} x^{k-j}Y_M(L(j)b,x).
\end{align*}
It follows from the twisted module Jacobi identity for the conformal vector, and the module translation identity. For $\eta(x)=\sum_{k\ge0}c_kx^{k+1}$ (interpreted coefficientwise in $\C[[x]]$), the commutator becomes
\begin{equation}\label{eq:relativeinfinitesimal}
 \left[\sum_{k\ge0}c_kL(k),Y_M(b,x)\right]
=\eta(x)\partial_xY_M(b,x) +\sum_{j\ge0}\frac{\eta^{(j+1)}(x)}{(j+1)!} Y_M(L(j)b,x).
\end{equation}
Subtracting $c_0h\id$ has no effect on this equation.

For $f(x)=x+\epsilon\eta(x)$, choose $\rho(u)=u+\epsilon\eta(u^N)/(Nu^{N-1})$, so that $\rho(u)^N=f(x)$ modulo $\epsilon^2$. The infinitesimal state action at the observation point is
\[
 R_V(f_x)=1-\epsilon\sum_{j\ge0}
 \frac{\eta^{(j+1)}(x)}{(j+1)!}L(j),
\]
where the sum is finite on each state. The module action is $1-\epsilon\bigl(\sum_{k\ge0}c_kL_M(k)-c_0h\id\bigr)$. Substitution in
\eqref{eq:relativefieldcov} shows that its
$\epsilon$-coefficient is zero by
\eqref{eq:relativeinfinitesimal}. Exponentiating proves
the formula for the prounipotent subgroup, coefficientwise in the formal variables. For $\rho(u)=au$ and $f(x)=a^Nx$, take homogeneous $b\in V_d$. The mode commutator implies $a^{-N(L_M(0)-h\id)}b_na^{N(L_M(0)-h\id)} =a^{-N(d-n-1)}b_n$.
The exponent $N(d-n-1)$ is integral, since $n\in N^{-1}\Z$. Together with $R_V(f_x)^{-1}b=a^{Nd}b$ and the branch substitution $x^{-n-1}\mapsto a^{-N(n+1)}x^{-n-1}$, this gives \eqref{eq:relativefieldcov} on every mode.
Scalings and the prounipotent subgroup generate the full origin-fixing group, proving the claimed formula.

For a root rotation $\zeta=e^{2\pi i/N}$, put $S=R_{M,h}(\zeta u)=e^{-2\pi i(L_M(0)-h\id)}$. Then $S^N=1$ and, for $gb=e^{2\pi ir/N}b$,
\[
 S\,Y_M^g(b,x)\,S^{-1}=Y_M^g(gb,x),\qquad
 Y_M^g(b,x)\big|_{u\mapsto\zeta u}=Y_M^g(g^{-1}b,x).
\]
Their product yields \eqref{eq:relativefieldcov} for $f(x)=x$. This also fixes the direction of the monodromy and agrees with the chosen root branch convention.

Each simple summand of a finite direct sum has its own lowest weight; apply the same construction on each summand. Module modes preserve these summands, so the field identity remains valid on their direct sum.
\end{proof}

The two coordinate actions are distinct: $R_V(ax)=a^{-L(0)}$ on states, whereas $R_{M,h}(au)=a^{-N(L_M(0)-h\id)}$ on a twisted module. Equation~\eqref{eq:relativefieldcov} compares them, but does not identify either action with the placement connection \eqref{eq:Tconn}, which is taken in fixed physical frames.

\begin{definition}\label{def:canonicalmarked}
Let $\pi:C\to\mathbb P^1\times S$ be the smooth framed $H$-cover of Lemma~\ref{lem:framedcharts}, $p:C\to S$ its structure map, and $\cF=F_E$ the $H$-equivariant finite-rank current bundle of \eqref{eq:canonicalfinite}. Its left linearization is $h\cdot(p_0,(b,v))=(hp_0,(hb,v\circ h^{-1})),\; h\in H$,
with the induced action on relative differentials, as in \cite[Section~5, equation~(5.2)]{FS}. For each marked orbit, let $\widetilde Q_i$ be its distinguished lift, $I_i=\langle g_i\rangle\le H$ its inertia group, $N_i=|I_i|$, and $u_i$ the chosen root coordinate satisfying $u_i(g_i p_0)=\zeta_i u_i(p_0),\; \zeta_i=e^{2\pi i/N_i},\; x_i=u_i^{N_i}$.
Here $x_i$ is the fixed physical local coordinate ($z-p_i$, $z$, or $z^{-1}$, as appropriate). The external factor at $\widetilde Q_i$ is a $g_i$-twisted module $M_i$, with $h_i:=h_{M_i}$; at the output orbit it is the prescribed contragredient module. Let $\mathcal B_i$ be a homogeneous $g_i$-eigenbasis of $E$, and put $g_i b=e^{2\pi i r_b/N_i}b$, $0\le r_b<N_i$.

On $u_i\ne0$, denote by $\iota_{x_i}(b)$ the state frame determined by $x_i$ at the variable observation point. The transition to an ordinary cover frame is the finite Laurent gauge of Lemma~\ref{lem:canonicalfinite}. Write the germ of an $H$-invariant meromorphic current in this frame as
\[
 s=\sum_{b\in\mathcal B_i}f_b(u_i)\,\iota_{x_i}(b)\,dx_i,
 \qquad f_b(u_i)\in\cO_S((u_i)),\quad
 f_b(\zeta_i u_i)=e^{2\pi i r_b/N_i}f_b(u_i).
\]
Define the normalized local action on $M_i$ by
\begin{equation}\label{eq:canonicalmarkedres}
 R_i(s)w=\frac1{N_i}\Res_{u_i=0}\left(
  \sum_{b\in\mathcal B_i}f_b(u_i)
  Y_{M_i}^{g_i}(b,u_i^{N_i})w\right)N_i u_i^{N_i-1}\,du_i.
\end{equation}
The coefficientwise finite interpretation is established in Proposition~\ref{prop:canonicalrealization}. For a global $H$-invariant current set $R(s)=\sum_iR_i(s)$ on the external tensor, with $R_i(s)$ acting on the $i$th factor. Coordinate changes use the representation $R_{M_i,h_i}$ of Lemma~\ref{lem:relativecoordinates}; the placement connection \eqref{eq:Tconn} remains in fixed physical module frames.
\end{definition}

\begin{proposition}\label{prop:canonicalrealization}
For each $i$ and $g_i$-eigenvector $b\in E$ with eigenvalue $e^{2\pi i r_b/N_i}$, the invariant punctured-disk current coefficient $u_i^{N_i n}$, $n\in r_b/N_i+\Z$, acts by $b_n$. More generally, the local action \eqref{eq:canonicalmarkedres} is well defined on every algebraic external tensor, is compatible with special-coordinate changes and deck transport, and has the exact mode expansion
\begin{equation}\label{eq:canonicalmoderes}
 R_i(s)w=\sum_{b\in\mathcal B_i}\sum_{n\in r_b/N_i+\Z}
 \bigl[u_i^{N_i n}\bigr]f_b(u_i)\,b_nw.
\end{equation}
The state-current gauge and its inverse have Laurent pole bounds independent of $n$. Invariant finite principal-part prescriptions modulo $\cF(D)$ in this local frame admit global lifts after the auxiliary pole allowance of Lemma~\ref{lem:smoothsymparts}.
\end{proposition}
\begin{proof}
Lemma~\ref{lem:canonicalfinite} constructs $\cF$, its $H$-linearization and the inverse finite Laurent current gauge. In the chosen $x_i$-frame, inertia leaves $x_i,dx_i$ fixed, so
\[
 s(g_ip_0)=g_i s(p_0)
 \ \Longleftrightarrow\ 
 f_b(\zeta_i u_i)=e^{2\pi i r_b/N_i}f_b(u_i)
 \ \Longleftrightarrow\ 
 f_b(u_i)\in u_i^{r_b}\cO_S((u_i^{N_i})).
\]
Consequently $[u_i^k]f_b=0$ unless $k\equiv r_b\pmod{N_i}$. For $Y_{M_i}^{g_i}(b,x)=\sum_{n\in r_b/N_i+\Z}b_nx^{-n-1}$, we have
\begin{equation}\label{eq:canonicalpuremode}
 \frac1{N_i}\Res_{u_i=0}
 u_i^{N_i n}Y_{M_i}^{g_i}(b,u_i^{N_i})N_i u_i^{N_i-1}\,du_i=b_n,
 \qquad n\in r_b/N_i+\Z.
\end{equation}
Indeed, the summand with index $m$ has exponent $N_i(n-m)-1$ in $u_i$. This proves \eqref{eq:canonicalmoderes} by linearity and also fixes the $1/N_i$ normalization.

For a fixed meromorphic $s$, write $f_b(u_i)=\sum_{k\ge-P_i}c_{b,k}u_i^k$ for some $P_i\in\N$, uniformly over the finite basis $\mathcal B_i$. Then $[u_i^{N_i n}]f_b=0$ if $N_i n<-P_i$. For each fixed $w\in M_i$, lower truncation gives $b_nw=0$ for $n\gg0$. Thus each inner sum in \eqref{eq:canonicalmoderes} is finite on $w$; the external tensor has finitely many factors. The coefficients vary algebraically over $S$, giving an $\cO_S$-linear action on local relative sections.

For a prescribed pure mode, form $u_i^{N_i n}\iota_{x_i}(b)\,dx_i$ on the punctured root disk. The inverse gauge of Lemma~\ref{lem:gauge}, including the Jacobian and every lower-weight term, has Laurent entries; hence this germ has a finite pole order in the ordinary cover frame. Its inertia condition is the congruence already verified. Conversely, a stabilizer-invariant finite jet at $\widetilde Q_i$ extends by deck transport to an $H$-invariant jet on the whole orbit. The assertion that this jet is the germ of a \emph{global} meromorphic current requires the cohomological lifting of Lemma~\ref{lem:smoothsymparts}, not merely the local gauge. The gauge's fixed Laurent bound also controls the number of regular cover jets needed to detect the base-frame negative part, independently of $n$.

Under a special-coordinate change $\rho$, the state-current gauge, its differential Jacobian and the module action $R_{M_i,h_i}(\rho)$ satisfy \eqref{eq:relativefieldcov}. Taking residues gives the conjugate action in the changed module frame; the scalar lowest-weight shift cancels in this conjugation. For $h\in H$, the pullback module $T_hM_i$ has $Y_{T_hM_i}(hb,x_i)=Y_{M_i}(b,x_i)$, while the distinguished root and state are transported by $h$. Thus the local normalized residues are preserved by deck transport, including the root rotations. The coordinate covariance agrees with \cite[Sections~4--6]{FS}.
\end{proof}

If a nodal family is identified on its generic fiber with the chosen framed cover, the finite coordinate bundle and the local actions agree under this identification, by base change of the associated bundle and the prescribed root and physical-coordinate identifications. This statement is confined to the \emph{generic smooth fiber}; it makes no assertion about an extension of the residue action across a node.

\begin{lemma}\label{lem:smoothsymparts}
Let $\widetilde Q_i$ be the distinguished lifts of Definition~\ref{def:canonicalmarked}. Henceforth, when discussing \emph{$H$-invariant divisors}, write
\[
 Q_i:=\sum_{h\in H/I_i}h\widetilde Q_i,
 \qquad Q:=\sum_iQ_i,\qquad \cF:=F_E.
\]
Thus $Q_i$ now denotes the \emph{reduced full orbit divisor}, not its distinguished point $\widetilde Q_i$. There exists $m\in\N$, locally on $S$ and independent of $i$ and $n\in\N$, such that, with $D=mQ$,
\begin{equation}\label{eq:smoothsymparts}
 \bigl(p_*\cF(D+nQ_i)\bigr)^H\twoheadrightarrow
 \bigl(p_*(\cF(D+nQ_i)/\cF(D))\bigr)^H
 \qquad(i,\ n\ge0).
\end{equation}
Moreover, the modes arising from local expansions of sections of $\cF(D)$ at the marked orbits have a common upper bound on conformal degree.
\end{lemma}
\begin{proof}
For every marking $i$, ramification of order $N_i=|I_i|$ gives the Cartier-divisor identity $N_iQ_i=\pi^*\bigl(\{x_i=0\}\times S\bigr)$.
Since $\mathcal O_{\mathbb P^1\times S}(\{x_i=0\}\times S)$ is relatively ample and $\pi$ is finite, $\mathcal O_C(N_iQ_i)$ is $p$-ample. Hence $\mathcal O_C(Q_i)$, and therefore $\mathcal O_C(Q)$, is $p$-ample. As $p$ is proper, $S$ is Noetherian and $\cF$ is coherent, relative Serre vanishing \cite[Tag~02O1]{Stacks} gives, after shrinking $S$ if necessary, $R^1p_*\cF(mQ)=0\; (m\ge m_0)$.
Fix such an $m$. For every $i$ and $n\in\N$, the exact sequence $0\longrightarrow\cF(D)\longrightarrow\cF(D+nQ_i) \longrightarrow\cF(D+nQ_i)/\cF(D)\longrightarrow0$
has the same left-hand term $\cF(D)$. Its pushforward is therefore surjective onto the principal-part quotient for all $i,n$. Each divisor is $H$-invariant; since $|H|\in\C^\times$, averaging $|H|^{-1}\sum_{h\in H}h$ is an exact projector. This proves \eqref{eq:smoothsymparts}. The quotient is supported on $nQ_i$, which is finite over $S$, and is locally free as an $\cO_S$-module after pushforward because $Q_i$ is a disjoint union of sections and $\cF$ is locally free there.

For the uniform degree bound, fix the finite Laurent state-current gauges from Lemma~\ref{lem:gauge}. The divisor $D$ gives, at every $i$, a common $P_i\in\N$ such that every base-frame coefficient of a section of $\cF(D)$ lies in $u_i^{-P_i}\cO_S[[u_i]]$. Formula \eqref{eq:canonicalmoderes} then permits only $n\ge-P_i/N_i$; for $b\in\mathcal B_i$, the corresponding mode has degree $\wt(b)-n-1$. Hence $\deg(b_n)\le\wt(b)+P_i/N_i-1 \le\max_{c\in\mathcal B_i}\wt(c)+P_i/N_i-1$.
The maximum over the finitely many marked orbits yields the required common bound, including every nonprimary lower-weight gauge term.
\end{proof}

For an invariant meromorphic one-form $\eta$ on the smooth cover, put
\[
 \Res_i^{\mathrm{base}}\eta:=\frac1{N_i}\Res_{\widetilde Q_i}\eta.
\]
Its residues are constant on the $H$-orbit of $\widetilde Q_i$, and $|H\widetilde Q_i|=|H|/N_i$. Consequently
\[
 \sum_{q\in H\widetilde Q_i}\Res_q\eta
 =\frac{|H|}{N_i}\Res_{\widetilde Q_i}\eta
 =|H|\Res_i^{\mathrm{base}}\eta.
\]
If $\eta$ has no poles outside these finitely many orbits, the residue theorem on the smooth proper cover yields $\sum_i\Res_i^{\mathrm{base}}\eta=0$.

\subsection{Finite generation and the ramified stress tensor}\label{app:stress}
Let $T=\bigotimes_{i=0}^mM_i$ be the external tensor of Section~\ref{sec:moduleprelim}, $\cT=\cO_S\otimes_\C T$, and $R(s)=\sum_iR_i(s)$ the normalized action of Definition~\ref{def:canonicalmarked}. Define the \emph{seed Ward submodule}
\[
 \cJ_E:=\sum_{s\in(p_*\cF(*Q))^H}\cO_S\,R(s)\cT\subseteq\cT,
 \qquad \cF(*Q):=\varinjlim_{r\ge0}\cF(rQ).
\]
Here the sum and the notation for relative sections are understood locally on the smooth framed chart $S$; $Q=\sum_iQ_i$ is the reduced full-orbit divisor of Lemma~\ref{lem:smoothsymparts}. Every section occurring in the sum has a finite pole bound, and $R(s)w$ is a finite sum for each $w\in T$ by Proposition~\ref{prop:canonicalrealization}.

\begin{proposition}\label{prop:seedfinite}
The $\cO_S$-module $\cT/\cJ_E$ is of finite type. More precisely, locally on $S$ it is generated by the images of $\bigotimes_{i=0}^m K_{R,N}(M_i)$ for a common choice of $R,N\in\Z$ as in Lemma~\ref{lem:lowdepth}.
\end{proposition}
\begin{proof}
Fix $A=mQ$ as in Lemma~\ref{lem:smoothsymparts}; in particular $R^1p_*\cF(A)=0$. By its degree estimate, there is $C_0\in\Z$ such that, at every marked orbit, any nonzero mode $c_\nu$ occurring in the local residue of a section of $\cF(A)$ satisfies
\begin{equation}\label{eq:boundedcorrection}
 \deg(c_\nu):=\wt(c)-\nu-1\le C_0,
 \qquad c\in\mathcal B_i,\quad \nu\in r_c/N_i+\Z.
\end{equation}
The bound is uniform in the section, the marked orbit and the external tensor on which its residue is evaluated. It incorporates all terms of the finite state-coordinate gauge, including $L(j)c$ and $\mathbf1$.

For each $i$, let $B_{g_i}$ and $\alpha_b$ be those of \eqref{eq:PBW}. Since $\bigcup_iB_{g_i}$ is finite, choose $R\in\Z$ so that $\Delta_i(b,n):=\wt(b)-\alpha_b+n-1>C_0,\; i=0,\ldots,m,\; b\in B_{g_i},\; n\ge R$.
Fix $i,b,n$ as above and write $\mu=\alpha_b-n$, $\Delta_*=\wt(b)-\mu-1$. At $\widetilde Q_i$ prescribe the invariant Laurent germ
\[
 s_{i,b,\mu}^{\mathrm{loc}}=u_i^{N_i\mu}\,\iota_{x_i}(b)\,dx_i,
 \qquad N_i\mu\in\Z,\quad N_i\mu\equiv r_b\pmod{N_i}.
\]
By Lemma~\ref{lem:gauge} and Proposition~\ref{prop:canonicalrealization}, this germ has finite pole order in the ordinary cover frame. For some $r\ge0$ its class modulo $\cF(A)$ thus determines an inertia-invariant principal part in $\cF(A+rQ_i)/\cF(A)$. Extend this class to the whole $H$-orbit by deck transport and prescribe zero at the other marked orbits. Work on one fixed affine open of $S$ on which the common divisor $A$ and its cohomology vanishing have been chosen. The pushforwards in \eqref{eq:smoothsymparts} are quasi-coherent, so taking sections over this affine open preserves that surjection for every $r$. Thus no sequence of mode-dependent shrinkings is needed. It gives a global invariant section $s_{i,b,\mu}\in\Gamma(C/S,\cF(A+rQ_i))^H$
with this principal part. In the distinguished root frames,
\begin{equation}\label{eq:warddeep}
 R_i(s_{i,b,\mu})=b_\mu+U_i,
 \qquad R_j(s_{i,b,\mu})=U_j\quad(j\ne i),
\end{equation}
where each $U_j$ is a locally finite sum of modes $c_\nu$ satisfying $\deg(c_\nu)\le C_0$. Indeed, subtraction of $s_{i,b,\mu}^{\mathrm{loc}}$ at $i$, and of zero at $j\ne i$, has germ in $\cF(A)$; the uniform bound \eqref{eq:boundedcorrection} applies. Formula \eqref{eq:canonicalpuremode} fixes the coefficient of $b_\mu$ to be $1$. These identities are identities of operators on individual module vectors, not assertions of finite sums in the unrestricted mode algebra.

Choose $N\ge\max(R-1,\max_{i,u}L_u)$ and a homogeneous PBW tensor $w=\bigotimes_jw_j$ in the fixed spanning families \eqref{eq:PBW} with some depth $n_1\ge R$. Because the depths are nonincreasing, we may choose such a factor $i$ with $w_i=b_{\alpha_b-n_1}v_i,\; v_i=b^2_{\alpha_2-n_2}\cdots b^s_{\alpha_s-n_s}u$,
where $s\ge1$ (and the empty tail is $u$). Put $v=\bigotimes_{j\ne i}w_j\otimes v_i$ with its factors restored to their original order. Then $w=b_\mu^{(i)}v$ and $\ee(w)=\ee(v)+D\Delta_*$, where $D\in\Z_{>0}$ is the common weight denominator in the definition of $\ee$ in Section~\ref{sec:moduleprelim}. Applying the seed relation $R(s_{i,b,\mu})v\equiv0\pmod{\cJ_E}$ and \eqref{eq:warddeep} gives
\begin{equation}\label{eq:actualdescent33}
 1\otimes w\equiv-\sum_{j=0}^{m}U_j^{(j)}v
 \pmod{\cJ_E},\qquad
 \ee(c_\nu^{(j)}v)-\ee(w)
 =D\bigl(\deg(c_\nu)-\Delta_*\bigr)<0
\end{equation}
for each nonzero homogeneous mode summand $c_\nu$ of $U_j$. For a fixed $v$, the sum in \eqref{eq:actualdescent33} is finite by lower truncation, and its coefficients lie in $\cO_S$ locally. Since $D\Delta_*,D\deg(c_\nu)\in\Z$, the strict inequality applies to the integer-valued energy $\ee\ge0$. Thus \eqref{eq:actualdescent33} verifies precisely \eqref{eq:descent-condition}, and Lemmas~\ref{lem:descent} and \ref{lem:lowdepth} give finite local generators. The finitely many simple external blocks are treated with a common maximum of the bounds and thresholds.
\end{proof}

\begin{lemma}\label{lem:singleinsertion}
Let $W_1\in\cC_g$, $W_2\in\cC_h$, $W_3\in\cC_{gh}$ be ordinary, let $\mathcal Y$ have type $\binom{W_3}{W_1\ W_2}$ and satisfy Definition~\ref{def:completeintertwiner}, and fix $z\in\C^*$ and the logarithmic branch used by $\mathcal Y$. For homogeneous $w_1\in W_1$, $w_2\in W_2$ and $w_3'\in W_3'$, the single-$V$ insertion matrix coefficients determine, on the framed connected $H$-cover of $(\mathbb P^1;z,0,\infty)$, $H=\langle g,h\rangle$, a meromorphic $H$-invariant dual section of $\mathscr V_E$, where $F_E=\mathscr V_E\otimes\omega_{C/S}$. Its only possible poles lie above $z$, $0$ and $\infty$; its three local germs are given by the $W_1$-, $W_2$- and contragredient $W_3'$-twisted vertex operators, respectively.
\end{lemma}
\begin{proof}
Set $g_z=g$, $g_0=h$ and $g_\infty=(gh)^{-1}$, so that $g_zg_0g_\infty=1$ in \eqref{eq:orderedmeridians}. Use the based paths obtained from the preferred three germs in \cite[Definition~2.8]{DH}, with this input order. By their two-point duality and analytic continuation, for $u\in V$ the three coefficient expansions
\begin{align*}
 F_\infty(x;u)&=\langle w_3',Y_{W_3}^{gh}(u,x)\mathcal Y(w_1,z)w_2\rangle,\\
 F_0(x;u)&=\langle w_3',\mathcal Y(w_1,z)Y_{W_2}^{h}(u,x)w_2\rangle,\\
 F_z(x;u)&=\langle w_3',\mathcal Y(Y_{W_1}^{g}(u,x-z)w_1,z)w_2\rangle
\end{align*}
represent branches of the same multivalued holomorphic function on $\mathbb P^1\setminus\{0,z,\infty\}$. Its continuation around a positive local meridian at $a\in\{z,0,\infty\}$ acts on the inserted state by $g_a^{-1}$, with the right-continuation and root conventions of \cite[Section~2]{DH}; at $\infty$ this meridian is clockwise in the affine coordinate $x$. Hence the state monodromy factors through $H$. The framed $H$-cover constructed in Lemma~\ref{lem:framedcharts} trivializes this monodromy (with the matching base paths and chosen lifts), and the continuation relations give its deck-equivariant transition rule.

For a homogeneous basis $u\in E$, the resulting functions are holomorphic on the punctured cover. The finite $\Aut\cO$-representation on $E$ (Lemma~\ref{lem:Efinite}) and the state-field coordinate covariance of \cite[Sections~3--5]{FS} identify their transition matrices with those of the dual state bundle. Thus they give a holomorphic $H$-invariant section of $\mathscr V_E^\vee$ away from the marked orbits; no assertion is made here for the full infinite-dimensional state bundle.

At $0$ and $z$, $F_0$ and $F_z$ have finite negative power bounds in the corresponding root coordinates by lower truncation of the two input twisted vertex operators, together with the convergence in \cite[Definition~2.8]{DH}. At $\infty$, use $s=x^{-1}$ and the opposite-vertex formula \eqref{eq:contragredientprelim}:
\[
 \langle Y_{W_3'}^{(gh)^{-1}}(u,s)w_3',w_3\rangle
 =\langle w_3',Y_{W_3}^{gh}\bigl(e^{sL(1)}(-s^{-2})^{L(0)}u,s^{-1}\bigr)w_3\rangle.
\]
For $u\in V_d$, the inserted state on the right equals $(-1)^d s^{-2d}\sum_{j=0}^{d}\frac{s^j}{j!}L(1)^j u$,
so the sum is finite; here $L(1)^{d+1}u=0$. Lower truncation on $W_3'$ bounds the pole at $s=0$. Taking a maximum over the finite homogeneous basis of $E$ gives a common pole bound. The Laurent gauges of Lemma~\ref{lem:gauge} and their inverses have fixed finite pole orders, so the same conclusion holds in the ordinary cover frames. The section extends meromorphically across all marked orbits. Since $\mathscr V_E$ is algebraic of finite rank on the smooth projective cover, GAGA \cite{SerreGAGA} identifies this meromorphic section with an algebraic meromorphic section of $\mathscr V_E^\vee$.
\end{proof}

\begin{proposition}\label{prop:binaryward}
For $W_1,W_2,W_3,\mathcal Y$ as in Lemma~\ref{lem:singleinsertion}, the three-point functional $\Phi_{\mathcal Y,z}(w_3'\otimes w_1\otimes w_2) :=\langle w_3',\mathcal Y(w_1,z)w_2\rangle$
using the fixed branch of $\log z$, annihilates the fiber $\cJ_E\otimes_{\cO_S}\kappa(z)$ on the standard framed cover at $z$. Here $\kappa(z)=\C$.
\end{proposition}
\begin{proof}
Fix $s\in\Gamma(C,\cF(*Q))^H$ and homogeneous external vectors. Let $\lambda\in\Gamma(C,\mathscr V_E^\vee(*Q))^H$ be the dual state section from Lemma~\ref{lem:singleinsertion}. Since $s\in\Gamma(C,\mathscr V_E\otimes\omega_{C/S}(*Q))^H$, the evaluation $\eta:=\langle\lambda,s\rangle$ is an $H$-invariant meromorphic one-form on the smooth proper cover $C$. Its poles lie only above $z,0,\infty$, and its local Laurent residues, using \eqref{eq:canonicalmarkedres} and the three branch germs in that lemma, satisfy $\Res_{\widetilde Q_i}\eta=N_i\,\Phi_{\mathcal Y,z}(R_i(s)w),\; i\in\{z,0,\infty\}$.
Thus the classical residue theorem, $|H\widetilde Q_i|=|H|/N_i$ and the exact normalization of \eqref{eq:canonicalpuremode} yield
\[
 0=\frac1{|H|}\sum_{q\in C}\Res_q\eta
 =\sum_{i\in\{z,0,\infty\}}\frac1{N_i}\Res_{\widetilde Q_i}\eta
 =\Phi_{\mathcal Y,z}(R(s)w).
\]
For every fixed $z$, linearity gives $\Phi_{\mathcal Y,z}(\cJ_E\otimes_{\cO_S}\kappa(z))=0$; the identity varies holomorphically on each chosen branch domain. The use of the residue theorem is restricted to this smooth proper finite cover and finite-rank $F_E$, as in \cite[Sections~5--6]{FS}.
\end{proof}

For $z=f(t)=t^N$ on the punctured root disk $t\ne0$, define $T_z(z)=Y(\omega,z)$ and let $T_t(t)$ denote its coordinate transform, with
\[
 T_t(t)=(f'(t))^2T_z(f(t))+\frac{c_V}{12}\{f,t\},
 \qquad \{f,t\}=\frac{f'''(t)}{f'(t)}-\frac32\left(\frac{f''(t)}{f'(t)}\right)^2
 =\frac{1-N^2}{2t^2}.
\]
The corresponding unnormalized pullback of $z^rT_z(z)\,dz$ is
\begin{align}\label{eq:stress}
 f^rT_z(f)\,df
 &=\frac{f^r}{f'}\left(T_t-\frac{c_V}{12}\{f,t\}\right)dt\notag\\
 &=\frac1N t^{Nr-N+1}T_t(t)\,dt
 +\frac{c_V(N^2-1)}{24N}t^{Nr-N-1}\,dt.
\end{align}
The second term is scalar and corresponds to the vacuum component $\mathbf1\in E$: indeed, $L(2)\omega=\frac{c_V}{2}\mathbf1$. It has a residue only when $r=1$, in which case the \emph{unnormalized} Schwarzian residue is $c_V(N^2-1)/(24N)$. When the normalized local action \eqref{eq:canonicalmarkedres} is used, its additional factor $1/N$ must also be applied. Equation \eqref{eq:stress} is a coordinate identity on $t\ne0$; it is not an invertible change of formal coordinate at $t=0$ for $N>1$.

\subsection{The differential quotient and horizontal channel functionals}\label{sec:diff}

Let $U\subset S$ be a smooth framed configuration chart, $A=\cO_S(U)$, and $\cT(U)=A\otimes_\C T$. Write $D_i=L_{M_i}(-1)^{(i)}$, with $D_i=0$ for fixed output coordinates, and use the placement connection
\eqref{eq:Tconn}. For $\delta,\gamma\in\operatorname{Der}_\C(A)$,
\[
 [D_i,D_j]=0,\qquad
 [\nabla_\delta,\nabla_\gamma]
 =\sum_{i\in I_{\mathrm{mov}}}
 \bigl(\delta\gamma(z_i)-\gamma\delta(z_i)\bigr)D_i
 +[\delta,\gamma]
 =\nabla_{[\delta,\gamma]}.
\]
The last equality uses the action of $[\delta,\gamma]$ on the scalar coefficients. It remains valid after pullback to the finite \emph{\'etale} framed Hurwitz chart. Under a change of physical module frame $e'=eB$,
$\nabla'=B^{-1}\nabla B$; the residue covariance of Proposition~\ref{prop:canonicalrealization} identifies the transformed seed Ward relations. Thus the constructions below are compatible on overlaps where the stated frame changes are defined.

\begin{definition}\label{def:diff}
For the seed Ward module $\cJ_E$ of Section~\ref{app:stress}, define
\begin{equation}\label{eq:diffward34}
 \cJ_{\diff}:=\sum_{r\ge0}\ \sum_{\delta_1,\ldots,\delta_r
 \in\operatorname{Der}_\C(A)}
 A\,\nabla_{\delta_1}\cdots\nabla_{\delta_r}\cJ_E,
 \qquad \cQ:=\cT/\cJ_{\diff},\qquad
 \cJ_E^{\diff}:=\cJ_{\diff}.
\end{equation}
Here $r=0$ means $\cJ_E$. The formula is interpreted sheafwise on the framed atlas.
\end{definition}

For $j\in\cJ_E$, $f\in A$, and $\delta,\delta_1,\ldots,\delta_r
\in\operatorname{Der}_\C(A)$,
\[
 \nabla_\delta\bigl(f\nabla_{\delta_1}\cdots
 \nabla_{\delta_r}j\bigr)
 =\delta(f)\nabla_{\delta_1}\cdots\nabla_{\delta_r}j
 +f\nabla_\delta\nabla_{\delta_1}\cdots\nabla_{\delta_r}j
 \in\cJ_{\diff}.
\]
Hence $\nabla\cJ_{\diff}\subseteq\cJ_{\diff}$, and $\cJ_{\diff}$ is the smallest $\nabla$-stable $A$-submodule of $\cT$ containing $\cJ_E$. In particular, $\nabla$ induces an integrable connection $\overline\nabla$ on $\cQ$.

\begin{lemma}\label{lem:connfree}
Let $X$ be a smooth variety over $\C$ and let $M$ be a coherent $\cO_X$-module with an integrable connection. Then $M$ is locally free of finite rank. If $N\subset M$ is connection-stable, the induced connection on $M/N$ is integrable.
\end{lemma}
\begin{proof}
The first assertion is \cite[Proposition~8.8]{Katz}. Fix $x\in X$ and put $A=\cO_{X,x}$, $I_j=\Fitt_j(M_x)$. Since $A$ is regular Noetherian and $M_x$ is finite over $A$,
$M_x$ is finitely presented. For $\delta\in\operatorname{Der}_\C(A)$, put $A_\epsilon=A[\epsilon]/(\epsilon^2)$ and $\sigma_0(a)=a$, $\sigma_\delta(a)=a+\epsilon\delta(a)$. The Leibniz rule gives an $A_\epsilon$-linear isomorphism
\[
 A_\epsilon\otimes_{A,\sigma_\delta}M_x
 \overset{\sim}{\longrightarrow}
 A_\epsilon\otimes_{A,\sigma_0}M_x,
 \qquad 1\otimes m\longmapsto
 1\otimes m+\epsilon\otimes\nabla_\delta m,
\]
with inverse obtained by replacing $\epsilon$ by $-\epsilon$. By the base-change formula for Fitting ideals, $\sigma_\delta(I_j)A_\epsilon=\sigma_0(I_j)A_\epsilon$; comparison modulo $\epsilon^2$ gives
$\delta(I_j)\subseteq I_j$. At a regular local $\C$-algebra,
\[
 \bigl(\forall\delta\in\operatorname{Der}_\C(A),\
 \delta(I)\subseteq I\bigr)
 \quad\Longrightarrow\quad I\in\{0,A\}.
\]
Indeed, if $0\ne I\subsetneq A$, choose $f\in I$ of minimal positive $\mathfrak m_x$-adic order. In regular local parameters, a coordinate derivation lowers the order of the nonzero initial homogeneous form of
$f$, contradicting $\delta(f)\in I$; characteristic zero is used here. Consequently each $I_j$ is $0$ or $A$. If $r=\dim_{A/\mathfrak m_x}(M_x/\mathfrak m_x M_x)$, then $\Fitt_r(M_x)=A$ and $\Fitt_{r-1}(M_x)=0$.
A minimal presentation $A^p\to A^r\to M_x\to0$ has all entries of the relation matrix in $\Fitt_{r-1}(M_x)=0$, so $M_x\simeq A^r$. For $r=0$, Nakayama's lemma gives $M_x=0$. Finally, if $N$ is stable then
\[
 \bigl([\overline\nabla_\delta,\overline\nabla_\gamma]
 -\overline\nabla_{[\delta,\gamma]}\bigr)[m]
 =\bigl[([\nabla_\delta,\nabla_\gamma]
 -\nabla_{[\delta,\gamma]})m\bigr]=0.
\]
\end{proof}

\begin{theorem}\label{thm:finiteconn}
The canonical $\cQ$ is a finite-rank locally free $\cO_S$-module with integrable connection. Locally, it is generated by the images of the finite low-depth tensors of Proposition~\ref{prop:seedfinite}.
\end{theorem}
\begin{proof}
On $U$, Proposition~\ref{prop:seedfinite} gives a finite set $\mathcal K\subset T$ for which
\[
 \sum_{w\in\mathcal K}A[w]=\cT/\cJ_E,
 \qquad \cT/\cJ_E\twoheadrightarrow\cQ.
\]
Thus $\cQ$ is coherent, since $A$ is Noetherian. Equation
\eqref{eq:diffward34} and the calculation preceding Lemma~\ref{lem:connfree}
give an integrable quotient connection. Apply that lemma (equivalently,
\cite[Proposition~8.8]{Katz}) to obtain local freeness.
\end{proof}

For a rooted finite tree $\tau$ with typed intertwining operators $\mathcal Y_v$, let $r_v$ be its vertex displacements and set $u_v=r_v/r_{\operatorname{parent}(v)}$ for every nonroot vertex.
After splitting any finite direct sum of intermediate modules, fix one ordinary simple module on each internal edge $e$, with $W_e=\bigoplus_{n\ge0}(W_e)_{h_e+n/N_e}$ and projections $\pi_{e,n}:W_e\to(W_e)_{h_e+n/N_e}$.
A \emph{coefficientwise finite formal channel} is the contraction of these operators over the $\pi_{e,n}$, with homogeneous external vectors and a homogeneous restricted-dual output, such that each
coefficient in the independent variables $u_v$ is a finite sum. For the ordinary modules used here this finiteness follows directly: by the $L(0)$-weight identity of Lemma~\ref{lem:IOlogs}, the exponent
of $u_e$ fixes $h_e+n/N_e$. Thus a fixed exponent tuple fixes all internal weights; each weight space is finite dimensional, and the finite tree has only finitely many intermediate simple summands.
These are formal coefficient identities; no convergence of the associated complex-variable series is assumed.

On a chosen algebraic channel chart, let $\mathscr R$ be a differential $\cO_S$-module in which these coefficientwise expressions are defined. One may use the corresponding branch of an iterated Laurent-series
field in the independent root ratios, adjoining the finitely many initial generalized-power monomials. Rational coordinate functions act by their expansions, and derivations act coefficientwise by the
formal chain rule. Write $\partial_\delta$ for this action of $\delta\in\operatorname{Der}_\C(A)$; hence $\partial_\delta(fF)=\delta(f)F+f\partial_\delta F$.

\begin{lemma}\label{lem:channelderivative}
Let $\Phi_\tau:\cT\to\mathscr R$ be the $\cO_S$-linear extension of such a channel, with fixed output puncture. Then
\begin{equation}\label{eq:channeldiff34}
 \partial_\delta\Phi_\tau(w)
 =\Phi_\tau(\nabla_\delta w),
 \qquad \delta\in\operatorname{Der}_\C(A),\quad w\in\cT.
\end{equation}
\end{lemma}
\begin{proof}
For a vertex $\mathcal Y$ of type $\binom{W_3}{W_1\ W_2}$, the $L(-1)$-derivative axiom and the Virasoro-field duality of \cite[Definition~2.8]{DH} give
\[
 \partial_x\mathcal Y(a,x)
 =\mathcal Y(L_{W_1}(-1)a,x)
 =L_{W_3}(-1)\mathcal Y(a,x)
 -\mathcal Y(a,x)L_{W_2}(-1).
\]
All $L(-1)$ operators in this formula are those of the indicated modules. Apply the formal chain rule to each vertex displacement. At an internal edge, the two appearances of the corresponding
$L(-1)$ operator cancel after the homogeneous contraction. Induction on the number of vertices gives, for a constant external state $w\in T$,
\[
 \partial_\delta\Phi_\tau(w)
 =\sum_{i\in I_{\mathrm{mov}}}
 \delta(z_i)\Phi_\tau(D_iw).
\]
Every coefficient is a finite contraction by the preceding weight argument, so the differentiation and cancellation are coefficientwise legitimate. For $f\in A$, the Leibniz rule gives
$\partial_\delta\Phi_\tau(fw)
 =\Phi_\tau(\delta(f)w+f\sum_i\delta(z_i)D_iw)$,
which is \eqref{eq:channeldiff34}.
\end{proof}

\begin{proposition}\label{prop:factor}
Let $\mathscr R$ be as above, and let $\Phi:\cT\to\mathscr R$ be a $\cO_S$-linear formal channel satisfying \eqref{eq:channeldiff34}. Assume the \emph{full seed Ward identity in $\mathscr R$}:
\begin{equation}\label{eq:seedformal34}
 \Phi\bigl(R(s)w\bigr)=0\quad\text{in }\mathscr R
 \qquad\bigl(s\in(p_*\cF(*Q))^H,\ w\in\cT\bigr).
\end{equation}
Then $\Phi(\cJ_{\diff})=0$, and it determines a horizontal section $\overline\Phi\in\Gamma(\cQ^\vee\otimes_{\cO_S}\mathscr R)$: $\overline\Phi\circ(\cT\twoheadrightarrow\cQ)=\Phi,\; \nabla^\vee_\delta\overline\Phi=0$.
\end{proposition}
\begin{proof}
By \eqref{eq:seedformal34}, $\Phi(\cJ_E)=0$. For $j\in\cJ_E$ and $r\ge0$, iterating
\eqref{eq:channeldiff34} gives
\[
 \Phi(\nabla_{\delta_1}\cdots\nabla_{\delta_r}j)
 =\partial_{\delta_1}\cdots\partial_{\delta_r}\Phi(j)=0.
\]
Hence \eqref{eq:diffward34} implies $\Phi(\cJ_{\diff})=0$. By Theorem~\ref{thm:finiteconn}, $\cQ$ is finite projective on $U$; therefore
\[
 \cQ^\vee\otimes_A\mathscr R
 \overset{\sim}{\longrightarrow}\Hom_A(\cQ,\mathscr R),
 \qquad \ell\otimes F\longmapsto
 \bigl([w]\mapsto\ell([w])F\bigr).
\]
The induced map corresponds to $\overline\Phi$, and
\[
 (\nabla^\vee_\delta\overline\Phi)([w])
 =\partial_\delta\overline\Phi([w])
 -\overline\Phi(\overline\nabla_\delta[w])=0.
\]
\end{proof}

Identity \eqref{eq:seedformal34} is an identity of \emph{formal series}, not merely an equality at one specialized value of the placement variables. If it has first been proved only after
localizing $\mathscr R$, its descent requires injectivity of the localization map; otherwise the conclusion is valid only in the localized module. The actual multi-input Ward identities used to
verify \eqref{eq:seedformal34} are established separately in Section~\ref{sec:threeinputward}; they do not follow from Definition~\ref{def:diff} or Theorem~\ref{thm:finiteconn}.

\section{Fusion rules and binary tensor products}\label{sec:binary}
\subsection{Finite fusion and the universal binary tensor product}
Let $g,h\in G$, $X\in\cC_g$, $Y\in\cC_h$, and $\mathcal S_{gh}$ be the simple representatives fixed in Section~\ref{sec:moduleprelim}. For $Z\in\cC_{gh}$ put
\[
 \cI^Z_{X,Y}:=\cI\binom Z{X\ Y},\qquad
 \cI^Z_{P(z),L;X,Y}:=\cI_{P(z),L}\binom Z{X\ Y},
 \quad z\in\C^*,\quad e^L=z.
\]
By \eqref{eq:Pzinverse}, evaluation
$\operatorname{ev}_L:\cI^Z_{X,Y}\xrightarrow{\sim}
\cI^Z_{P(z),L;X,Y}$ is a linear bijection.

\begin{proposition}\label{prop:fusion}
For $X\in\cC_g$, $Y\in\cC_h$, and $Z\in\cC_{gh}$, $\dim_{\C}\cI^Z_{X,Y}<\infty$. More precisely, if $X,Y,Z$ are simple and $\tilde z$ is a framed point over $(\infty,z,0)$, then, for $T=Z'\otimes X\otimes Y$,
\begin{equation}\label{eq:fusionbound41}
 \cI^Z_{X,Y}\hookrightarrow
 \Hom_{\C}\bigl((\cT/\cJ_E)\otimes_{\cO_S}\kappa(\tilde z),\C\bigr),
 \qquad \kappa(\tilde z)=\C.
\end{equation}
\end{proposition}
\begin{proof}
First suppose $X,Y,Z$ are simple. For $\mathcal Y\in\cI^Z_{X,Y}$, put
\[
 \Phi_{\mathcal Y,\tilde z}(z'\otimes x\otimes y)
 :=\big\langle z',\operatorname{ev}_L(\mathcal Y)(x\otimes y)\big\rangle,
 \qquad z'\in Z',\ x\in X,\ y\in Y.
\]
Propositions~\ref{prop:canonicalrealization} and
\ref{prop:binaryward} give
$\Phi_{\mathcal Y,\tilde z}(\cJ_E\otimes\kappa(\tilde z))=0$; hence \eqref{eq:fusionbound41} is well defined and linear. If $\Phi_{\mathcal Y,\tilde z}=0$, then for $x\in X_\alpha$, $y\in Y_\beta$, and
$\gamma\in\wt(Z)$, Lemma~\ref{lem:IOlogs} yields
\[
 \pi^Z_\gamma\mathcal Y(x,t)y
 =e^{(\gamma-\alpha-\beta)(\log t-L)}
  \pi^Z_\gamma\operatorname{ev}_L(\mathcal Y)(x\otimes y)=0.
\]
The restricted dual $Z'=\bigoplus_\gamma Z_\gamma^*$ separates $\overline Z=\prod_\gamma Z_\gamma$, so $\mathcal Y=0$. Proposition~\ref{prop:seedfinite} makes the fiber in
\eqref{eq:fusionbound41} finite dimensional, proving the assertion
for simple modules.

For arbitrary $X,Y,Z$, Proposition~\ref{prop:regularmodules} gives finite decompositions into ordinary simples. Projection onto, and inclusion of, the simple summands give a finite direct-sum
isomorphism of the corresponding intertwining-operator spaces; all axioms of Definition~\ref{def:completeintertwiner} are preserved under these compositions. The simple case therefore implies the assertion.
\end{proof}

\begin{lemma}\label{lem:Schur}
For $S,T\in\mathcal S_g$,
\[
 \Hom_V(S,T)\cong
 \begin{cases}\C,&S=T,\\0,&S\ne T.\end{cases}
\]
\end{lemma}
\begin{proof}
For $f\in\End_V(S)$, choose $\lambda\in\wt(S)$ with $S_\lambda\ne0$. Since $[f,L_S(0)]=0$ and $\dim S_\lambda<\infty$, there exist $a\in\C$ and $0\ne v\in\ker(f-a\id)\cap S_\lambda$. Simplicity gives $\ker(f-a\id)=S$, hence $f=a\id$.
For $S\ne T$, a nonzero map $S\to T$ would be an isomorphism, a contradiction.
\end{proof}

\begin{theorem}\label{thm:binary}
For $X\in\cC_g$, $Y\in\cC_h$, define
\begin{equation}\label{eq:Pz}
 U_{X,Y}:=\bigoplus_{S\in\mathcal S_{gh}}
 (\cI^S_{X,Y})^*\otimes_{\C}S.
\end{equation}
Then, for every $(z,L)$ with $e^L=z\ne0$, the pair $(U_{X,Y},I^z_{X,Y})$ defined below is a $P(z)$-tensor product in $\cC_{gh}$, with a natural bijection
\begin{equation}\label{eq:binaryrepresent41}
 \Hom_V(U_{X,Y},W)\xrightarrow{\sim}
 \cI^W_{P(z),L;X,Y},\qquad
 f\longmapsto\overline f\circ I^z_{X,Y}
 \quad(W\in\cC_{gh}).
\end{equation}
Moreover, $(X,Y)\mapsto U_{X,Y}$ extends to a bifunctor $\cC_g\times\cC_h\to\cC_{gh}$.
\end{theorem}
\begin{proof}
Write $H_S:=\cI^S_{X,Y}$. By Proposition~\ref{prop:fusion} and finiteness of $\mathcal S_{gh}$, all summands of
\eqref{eq:Pz} have finite multiplicity, so $U_{X,Y}\in\cC_{gh}$.
For each $S$, choose a basis $(\mathcal Y_{S,a})_a$ of $H_S$ and its dual basis $(e_{S,a})_a$ of $H_S^*$. Set
\begin{equation}\label{eq:binarycanonical41}
 \begin{aligned}
 \mathcal Y^{\can}_{X,Y}(x,t)y
   &:=\sum_{S\in\mathcal S_{gh}}\sum_a
       e_{S,a}\otimes\mathcal Y_{S,a}(x,t)y,\\
 I^z_{X,Y}&:=\operatorname{ev}_L(\mathcal Y^{\can}_{X,Y})
     :X\otimes Y\longrightarrow\overline{U_{X,Y}}.
 \end{aligned}
\end{equation}
The tensor $\sum_a e_{S,a}\otimes\mathcal Y_{S,a}$ is the canonical element of $H_S^*\otimes H_S$; thus \eqref{eq:binarycanonical41} is basis independent. Its sum is finite and each summand satisfies
Definition~\ref{def:completeintertwiner}; hence $\mathcal Y^{\can}_{X,Y}\in\cI^{U_{X,Y}}_{X,Y}$ and $I^z_{X,Y}\in\cI^{U_{X,Y}}_{P(z),L;X,Y}$.

Let $W\in\cC_{gh}$ and put $M_S(W):=\Hom_V(S,W)$. Proposition~\ref{prop:regularmodules} and Lemma~\ref{lem:Schur} give the natural isomorphism
\[
 \bigoplus_{S\in\mathcal S_{gh}}M_S(W)\otimes S
 \xrightarrow{\sim}W,\qquad m\otimes s\longmapsto m(s).
\]
By additivity of the typed intertwining-operator axioms,
\begin{equation}\label{eq:binaryHom41}
 \begin{aligned}
 \Hom_V(U_{X,Y},W)
 &\cong\bigoplus_S\Hom_{\C}(H_S^*,M_S(W))\\
 &\cong\bigoplus_SM_S(W)\otimes H_S
 \cong\cI^W_{X,Y}
 \xrightarrow[\operatorname{ev}_L]{\sim}
 \cI^W_{P(z),L;X,Y}.
 \end{aligned}
\end{equation}
The middle isomorphism sends $m\otimes\mathcal Y$ to $m\circ\mathcal Y$. On $e_{S,a}\otimes s$, the first isomorphism is induced by the $V$-map $f$; substitution in \eqref{eq:binarycanonical41} shows that
\eqref{eq:binaryHom41} sends $f$ precisely to
$\overline f\circ I^z_{X,Y}$. This proves
\eqref{eq:binaryrepresent41}, including uniqueness and
naturality in $W$, i.e. Definition~\ref{def:Pzproduct}.

For $f:X\to X'$ in $\cC_g$ and $q:Y\to Y'$ in $\cC_h$, the operator $\mathcal Y^{\can}_{X',Y'}\circ(f\otimes q)$ is of type $\binom{U_{X',Y'}}{X\ Y}$. The universal property gives a unique $V$-map
$ f\bo_{P(z)}q:U_{X,Y}\to U_{X',Y'}$ with $\overline{f\bo_{P(z)}q}\circ I^z_{X,Y} =I^z_{X',Y'}\circ(f\otimes q)$. Identity and composition laws follow from uniqueness. The construction uses semisimplicity and
Proposition~\ref{prop:fusion}; it does not invoke the additional existence assumptions in \cite[Section~4]{DH}.
\end{proof}

For different $(z,L)$ the representing object \eqref{eq:Pz} is unchanged, whereas its evaluation map is
\eqref{eq:binarycanonical41} with the indicated logarithm.
Write $X\bo Y:=U_{X,Y}=X\bo_{P(1)}Y$ ($\log1=0$), and $\bo_t:=\bo_{P(t)}$. For $z_1,z_2\in\R_{>0}$ and positive-real $\gamma:z_1\to z_2$, the logarithms continue within $\R$, and \eqref{eq:Pztransport} gives $\mathsf T_\gamma=\id_{U_{X,Y}}$.
This is the positive-real normalization of
\cite[Section~4, parallel transport]{DH}; for general paths,
$\mathsf T_\gamma$ retains its winding dependence.

\subsection{Operator-pair spaces and elementary structure maps}
Fix $g,h,k\in G$, $A\in\cC_g$, $B\in\cC_h$, $C\in\cC_k$ and $D\in\cC_{ghk}$. Define
\begin{align*}
\mathsf P_{ABC}^{D}
&:=\bigoplus_{M\in\mathcal S_{hk}}
 \cI\binom D{A\ M}\otimes_{\C}\cI\binom M{B\ C},
 \\
\mathsf I_{ABC}^{D}
&:=\bigoplus_{N\in\mathcal S_{gh}}
 \cI\binom D{N\ C}\otimes_{\C}\cI\binom N{A\ B}.
\end{align*}
Proposition~\ref{prop:fusion} and $|\mathcal S_{hk}|,|\mathcal S_{gh}|<\infty$ imply $\dim_\C\mathsf P_{ABC}^{D},\dim_\C\mathsf I_{ABC}^{D}<\infty$. Theorem~\ref{thm:binary}, applied twice, and the natural
simple-summand decomposition of its proof give
\begin{equation}\label{eq:pairHom42}
\mathsf P_{ABC}^{D}\simeq
\Hom_V(A\bo(B\bo C),D),\qquad
\mathsf I_{ABC}^{D}\simeq
\Hom_V((A\bo B)\bo C,D).
\end{equation}
More explicitly, the first identification uses $B\bo C=\bigoplus_M(\cI^M_{B,C})^*\otimes M$ and
\[
 \cI^D_{A,B\bo C}
 \simeq\bigoplus_M\cI^M_{B,C}\otimes_\C\cI^D_{A,M};
\]
the second is obtained similarly. These are algebraic identifications of operator spaces; no convergence of an infinite product--iterate contraction is asserted here.

For $q\in G$, $X\in\cC_g$ and $Y\in\cC_h$, regard the same formal series $\mathcal Y^{\can}_{X,Y}(x,t)y$ as an operator of type
\[
 \binom{T_q(X\bo Y)}{T_qX\ T_qY}.
\]
Indeed, $Y_{T_qW}(v,t)=Y_W(q^{-1}v,t)$, $q\omega=\omega$, and each defining analytic identity for the pulled-back operator is the corresponding identity for the original operator with $v_i$ replaced by $q^{-1}v_i$; see \cite[Sections~2--3]{DH}.
Theorem~\ref{thm:binary} therefore defines the \emph{action tensorator}
\begin{equation}\label{eq:J}
 J^q_{X,Y}:T_qX\bo T_qY\xrightarrow{\ \sim\ }T_q(X\bo Y),
 \qquad
 \overline{J^q_{X,Y}}\,I^1_{T_qX,T_qY}=I^1_{X,Y},
\end{equation}
where the right side takes values in $\overline{T_q(X\bo Y)}$ on the unchanged underlying vector space. Applying the same construction with $q^{-1}$ gives $(J^q_{X,Y})^{-1} =T_q\bigl(J^{q^{-1}}_{T_qX,T_qY}\bigr)$.
Uniqueness in \eqref{eq:Pzuniversal} yields
\[
 J^1_{X,Y}=\id_{X\bo Y},\qquad
 J^{pq}_{X,Y}=T_p(J^q_{X,Y})\circ J^p_{T_qX,T_qY}
 \quad(p,q\in G).
\]
The unit is $V\in\cC_1$; put $ j_q:V\to T_qV$, $j_q(v)=q^{-1}v$. Then
\[
 Y_{T_qV}(u,t)j_q(v)
 =q^{-1}Y_V(u,t)v=j_q(Y_V(u,t)v),
 \quad j_1=\id_V,\quad j_{pq}=T_p(j_q)\circ j_p.
\]
In particular $T_q$ is a strict action on the underlying categories, not a strict monoidal action with $J^q=\id$.

For $X\in\cC_g$, $Y\in\cC_h$, the raw crossed commutativity map is \eqref{eq:Pzrawcrossed}; fix the upper-half-plane path $\gamma:-1\to1$ specified there and put
\begin{equation}\label{eq:c}
 c_{X,Y}:=\mathsf T^{T_gY,X}_\gamma\circ\mathcal R^1_{X,Y}
 :X\bo Y\longrightarrow T_gY\bo X.
\end{equation}
Thus the positive half-exchange sends $(g,h)$ to $(ghg^{-1},g)$. The isomorphism assertion is proved below.

\begin{lemma}\label{lem:units}
For $X\in\cC_g$ the morphisms $l_X,r_X$ of
\eqref{eq:Pzunits} are isomorphisms. For $X\in\cC_g$,
$Y\in\cC_h$, both $\mathcal R^1_{X,Y}$ and $c_{X,Y}$ are isomorphisms.
\end{lemma}
\begin{proof}
Let $Z\in\cC_g$ and $\mathcal Y\in\cI\binom Z{V\ X}$. By the $L(-1)$-derivative property,
\[
 \partial_t\mathcal Y(\mathbf1,t)=
 \mathcal Y(L(-1)\mathbf1,t)=0,
 \qquad f_{\mathcal Y}(x):=\mathcal Y(\mathbf1,t)x.
\]
The one-$V$ duality of Definition~\ref{def:completeintertwiner}, first with $\mathbf1$ and then with $Y_V(v,s)\mathbf1=e^{sL(-1)}v$, gives
\[
f_{\mathcal Y}Y_X^g(v,t)=Y_Z^g(v,t)f_{\mathcal Y},
 \qquad
 \mathcal Y(v,t)x=Y_Z^g(v,t)f_{\mathcal Y}(x).
\]
The second identity also uses $\mathcal Y(e^{sL(-1)}v,t)=e^{s\partial_t}\mathcal Y(v,t)$; both identities are obtained in the common analytic germs and then coefficientwise. Conversely, for $f\in\Hom_V(X,Z)$, $Y_Z^g(v,t)f(x)$ is an intertwining
operator of type $\binom Z{V\ X}$, since $f$ is a module morphism and the module vertex operator satisfies the analytic twisted-module identities. Consequently,
\[
 \cI\binom Z{V\ X}\xrightarrow{\ \sim\ }\Hom_V(X,Z),
 \quad\mathcal Y\longmapsto f_{\mathcal Y},
 \quad f\longmapsto Y_Z^g(\cdot,t)f(\cdot).
\]
Together with Theorem~\ref{thm:binary}, this is a natural bijection $\Hom_V(V\bo X,Z)\simeq\Hom_V(X,Z)$ induced by $l_X$; its inverse is induced by $Y_X^g$, hence $l_X$ is an isomorphism.

For the right unit, the skew isomorphism
\cite[Theorem~3.1]{DH} has the special case
\[
 \Omega_+: \cI\binom Z{X\ V}
 \xrightarrow{\ \sim\ }\cI\binom Z{V\ X},
 \qquad \Omega_+\Gamma_X=Y_X^g,
\]
since the second input has twist $1$ and $\Gamma_X=\Omega_-Y_X^g$. Thus $\cI\binom Z{X\ V}\simeq\Hom_V(X,Z)$, with the canonical map corresponding to $\Gamma_X$. The same universal argument shows that $r_X$ is an
isomorphism. No identification of $T_gV$ with $V$ is required for this step.

Finally, put $W=T_gY\bo_{P(-1)}X$. By \cite[Theorem~3.1]{DH}, the mutually inverse transformations $\Omega_+: \cI\binom W{T_gY\ X} \xrightarrow{\ \sim\ }\cI\binom W{X\ Y},\; \Omega_-=(\Omega_+)^{-1}$ are well typed because $T_{g^{-1}}T_gY=Y$.
Applied to the universal $P(-1)$-operator and its inverse transform, they give inverse representing morphisms by
\eqref{eq:Pzuniversal}. Hence $\mathcal R^1_{X,Y}$ is
invertible. So is $c_{X,Y}$, since $\mathsf T_\gamma^{-1}=\mathsf T_{\gamma^{-1}}$ by \eqref{eq:Pztransport}. The logarithms for both half-turns are those of \eqref{eq:skewprelim} and
\eqref{eq:Pzrawcrossed}.
\end{proof}

\begin{lemma}\label{lem:generation}
For $X\in\cC_g$, $Y\in\cC_h$, the homogeneous coefficients of $\mathcal Y^{\can}_{X,Y}$ span $X\bo Y$. For every finite rooted binary tree, the iterated algebraic canonical coefficients span its representing module.
In particular, a module morphism out of either module is zero if it vanishes on all the indicated coefficients.
\end{lemma}
\begin{proof}
Fix $(z,L)$ and put $U=X\bo Y$. For a $P(z)$-product $(W,I)$ define
\[
 \lambda_{I,w'}(x\otimes y):=
 \langle w',I(x\otimes y)\rangle,
 \quad w'\in W',\qquad
 \Lambda_{X,Y}^z:=
 \Span_{\C}\{\lambda_{I,w'}\}\subseteq(X\otimes Y)^*,
\]
where $(W,I)$ ranges over the $P(z)$-products in $\cC_{gh}$. In \cite[Section~4, Proposition~4.3, Theorems~4.6 and~4.8]{DH}, $\Lambda_{X,Y}^z$ is equipped with its $(gh)^{-1}$-twisted contragredient action, and its
restricted dual represents the $P(z)$-product functor, provided that the stated closure and double-dual hypotheses hold. Here they hold: Proposition~\ref{prop:regularmodules} gives finite semisimplicity and finite simple spectra,
Proposition~\ref{prop:fusion} gives finite fusion spaces, and $W''\simeq W$ for grading-restricted $W$ by Section~\ref{sec:moduleprelim}. Thus uniqueness of representing objects, Theorem~\ref{thm:binary}, and the defining pairing
give an isomorphism
\begin{equation}\label{eq:dualspan42}
 U'\xrightarrow{\ \sim\ }\Lambda_{X,Y}^z,
 \qquad
 w'\longmapsto\bigl(x\otimes y\mapsto
 \langle w',I^z_{X,Y}(x\otimes y)\rangle\bigr).
\end{equation}
In particular the displayed map is injective.

For $\mu\in\wt(U)$, let $E_\mu$ be the linear span of $\pi^U_\mu I^z_{X,Y}(x\otimes y)$ for homogeneous $x\in X$, $y\in Y$. If $E_\mu\ne U_\mu$, finite dimensionality gives $0\ne w'\in U_\mu^*$ with
$w'(E_\mu)=0$. Extending it by zero on the other weight spaces yields $w'\in U'$ whose image in
\eqref{eq:dualspan42} is zero, a contradiction. Hence
$E_\mu=U_\mu$ for every $\mu$. By Lemma~\ref{lem:IOlogs},
\[
 \pi_\mu^U\mathcal Y^{\can}_{X,Y}(x,t)y
 =e^{(\mu-\wt(x)-\wt(y))(\log t-L)}
 \pi_\mu^U I^z_{X,Y}(x\otimes y).
\]
Therefore the formal homogeneous coefficients span $U$.

For a finite rooted binary tree, interpret an \emph{iterated algebraic coefficient} as follows: choose a homogeneous output coefficient at each internal vertex and apply the next binary operator to the resulting homogeneous vectors.
Every selected coefficient is a finite algebraic vector; no infinite internal-weight sum is taken. The claim follows by induction on the number of internal vertices. If the two subtree coefficient sets span $A$ and $B$, bilinearity
expresses each homogeneous coefficient of the canonical operator $A\otimes B\to(A\bo B)\{t\}$ as a finite linear combination of full-tree coefficients. The binary case shows that these span $A\bo B$. Finally, a linear map
vanishing on a spanning set is zero.
\end{proof}

\section{Regularity of the canonical differential quotient}
\label{sec:regularity}\label{sec:canonicalboundary}
\subsection{Actual covers at the boundary}
Let $U=F_m(\C^*)$ and let $S\to U$ be a connected component of the finite framed Hurwitz cover of Definition~\ref{def:canonicalmarked}. Write $K(S)$ for its function field. All boundary valuations below
are divisorial valuations in characteristic zero.

\begin{definition}\label{def:clustercompactification}
For $A\subseteq\{1,\ldots,m\}$, put $\Delta_A=\{p_i=p_j:i,j\in A\}$ for $|A|\ge2$, and $\Delta_{a,A}=\{p_i=a:i\in A\}$ for $a\in\{0,\infty\}$ and $A\ne\varnothing$. Let
\[
 \mathcal G=\{\Delta_A:|A|\ge2\}\cup
 \{\Delta_{0,A},\Delta_{\infty,A}:A\ne\varnothing\}
 \subseteq (\mathbb P^1)^m.
\]
The equations are interpreted in the appropriate affine coordinate near $0$ or $\infty$. Nonempty intersections are indexed by disjoint maximal equality blocks, each with at most one anchor in
$\{0,\infty\}$. The minimal members of $\mathcal G$ containing such an intersection impose independent linear equations in these local coordinates; hence $\mathcal G$ is a building set in the sense of
\cite[Definition~2.1 and Theorems~1.2--1.3]{Li}.
Let $\overline U_{\rm cl}$ be the iterated blow-up of members of $\mathcal G$ of codimension at least two in an inclusion-compatible order. The divisorial members of $\mathcal G$ are retained as boundary
divisors. The cited wonderful-compactification theorem gives
\[
 \overline U_{\rm cl}\text{ smooth and proper},\qquad
 \overline U_{\rm cl}\setminus U\text{ a simple normal crossing divisor}.
\]
Set $\overline S_{\rm cl}=\operatorname{Norm}_{K(S)}(\overline U_{\rm cl})$, the normalization in $K(S)$.
\end{definition}

Write $\overline U=\overline U_{\rm cl}$ and $\overline S=\overline S_{\rm cl}$ throughout this section.

\begin{lemma}\label{lem:clustercharts}
Let $D\subseteq\overline U\setminus U$ be a boundary prime divisor and $L=\C(D)$. On a dense open of $D$, a transversal has completed base $L[[q]]$ and one of the following forms, after extending $L$ if necessary:
\begin{enumerate}
\item $p_i=z_*+q\xi_i$ for $i\in A$, where $A$ is a free finite
cluster ($|A|\ge2$) or a cluster anchored at $0$ ($z_*=0$). The $\xi_i$ are pairwise distinct; in the anchored case the fixed marking has $\xi_0=0$ and $\xi_i\ne0$ for moving markings.
\item $p_i=q^{-1}\xi_i$ for $i\in A$, where $A$ is anchored at
$\infty$, with distinct $\xi_i\in L^\times$. Put $s=z^{-1}$, $\xi=qz$ and $\eta=\xi^{-1}$; then $s=q\eta$.
\end{enumerate}
All other markings have distinct, nonspecial limits. Blowing up the coarse total curve at the collision produces a two-component nodal family with local equation $xy=q$. Each component meets an original
marked section (among $0,\infty,p_1,\ldots,p_m$). Under $q=t^e$ the node has equation $xy=t^e$.
\end{lemma}
\begin{proof}
At the generic point of $D$ exactly one cluster divisor occurs. For a finite cluster take the nonzero collision scale $q$ and $\xi=(z-z_*)/q$; near the node use $x=z-z_*$ and $y=\xi^{-1}$, so $xy=q$.
For an infinity cluster perform the construction in $s_i=p_i^{-1}$; near its node use $x=s$ and $y=\xi= qz$, again giving $xy=q$. At the generic point the nonzero cluster ratios and their pairwise
differences are units. A free cluster contains at least two moving markings; an anchored cluster contains its fixed marking and at least one moving marking. Its complement contains the other fixed marking.
Thus both components meet the original marked divisor. Base change $q=t^e$ gives the final equation.
\end{proof}

\begin{lemma}\label{lem:primes}
For $\nu:\overline S\to\overline U$, every codimension-one point $\eta\in\overline S\setminus S$ maps to a codimension-one point of $\overline U\setminus U$. Moreover, $\mathcal O_{\overline S,\eta}$ is a DVR and $\nu^{-1}(U)\simeq S$.
\end{lemma}
\begin{proof}
The normalization is finite because $\overline U$ is a variety over $\C$ \cite[Tag~0BXR]{Stacks}. Both varieties have dimension $m$; for the finite dominant morphism $\nu$, the image of an irreducible
closed set has the same dimension. Thus a prime divisor maps to a prime divisor, and normality makes its height-one local ring a DVR. The finite normal cover $S\to U$ is the normalization of $U$ in
$K(S)$; uniqueness of integral closure gives $\nu^{-1}(U)\cong S$.
\end{proof}

\begin{lemma}\label{lem:pureramification}
Let $\operatorname{char}L=0$ and let $K'/L((q))$ be a finite extension with its $q$-adic valuation. After a further finite extension, there
exist a finite extension $L''/L$, a uniformizer $t$, and $e\ge1$ such that $\mathcal O_{K''}=L''[[t]],\; q=t^e,\; (t\partial_t)|_{L''}=0$. If $L=\C((s))$, after a further finite coefficient-field extension
there is a uniformizer $z$ with $s=z^f$, $f\ge1$.
\end{lemma}
\begin{proof}
The complete discretely valued field $K'$ has a finite separable residue-field extension $L'/L$. Lift a primitive element of $L'/L$ by Hensel's lemma applied to its separable minimal polynomial over
$L\subset L[[q]]$. Thus $L'\subset\mathcal O_{K'}$ and $\mathcal O_{K'}\simeq L'[[\pi]]$ for a uniformizer $\pi$. Write $q=\pi^e\epsilon(\pi)$ with $\epsilon(\pi)\in L'[[\pi]]^\times$. After adjoining an $e$th root
of $\epsilon(0)$, the binomial series supplies $\epsilon(\pi)^{1/e}\in L''[[\pi]]^\times$. Set $t=\pi\epsilon(\pi)^{1/e}$; then $q=t^e$, and the inverse formal coordinate change gives $\mathcal O_{K''}=L''[[t]]$.
Since $L''$ is the chosen coefficient field, $t\partial_t$ vanishes on $L''$. If $L=\C((s))$, the same argument for the finite extension $L''/\C((s))$ gives $s=z^f$ (the residue field is $\C$).
Further integral powers of $t,z$ may clear a finite collection of
\emph{rational} relative-weight denominators; arbitrary complex
lowest weights remain formal monomial factors, not integral powers.
\end{proof}

\begin{lemma}\label{lem:actualcoverextension}
Let $\eta_S$ be a codimension-one boundary point of $\overline S$ lying above a boundary point $\eta_U$ of $\overline U$, and put $R=\widehat{\mathcal O}_{\overline S,\eta_S}$, $K=\operatorname{Frac}(R)$. Let $C_\eta\to\mathbb P^1_{K}$
be the pullback of the actual connected framed $H$-cover of Definition~\ref{def:canonicalmarked} along $\operatorname{Spec}K\to S$. The coarse family is the base change of the family in
Lemma~\ref{lem:clustercharts} along the finite extension of completed boundary DVRs. After a finite DVR extension $R\subset R'$, there is an admissible $H$-cover $C\to X$ with $X/R'$ the prescribed two-component coarse
family, $C/R'$ projective and nodal, and $C\times_{R'}\operatorname{Spec}K'\cong C_\eta\otimes_K K'$, where $K'=\operatorname{Frac}(R')$. Every irreducible component of $C_0$ meets the inverse image of an original marked section. The generic
$H$-action, inertia and framing are preserved under the chosen generic-fiber identification.
\end{lemma}
\begin{proof}
Write $\Sigma$ for the divisor of original markings. If a component of the coarse special fiber has fewer than three special points, add disjoint auxiliary sections in the relative smooth locus,
generically avoiding the branch locus of $C_\eta$; their inertia is trivial. After this stabilization $(X,\Sigma_{\rm st})/R$ is a stable pointed family, and the generic framed cover, branched only
at the original markings, defines a representable twisted $H$-cover of its smooth generic fiber.

The balanced twisted-cover stack is proper over the stack of stable pointed curves \cite[Corollary~3.0.5]{ACV}. Its valuative criterion, after a finite DVR extension, yields a balanced twisted $H$-cover
whose coarse pointed curve is the specified stabilized family. The equivalence with admissible $H$-covers is
\cite[Theorem~4.3.2]{ACV}; the source is projective and nodal by
\cite[Lemma~2.2.1]{ACV}. Forgetting the auxiliary markings does not
alter $X$ or the cover, and the original marking divisor $\Sigma$ remains distinguished from the auxiliary divisor.

The admissible map $C_0\to X_0$ is finite, and each irreducible component of $C_0$ maps onto an irreducible component of $X_0$. Every coarse component meets $\Sigma_0$ by Lemma~\ref{lem:clustercharts}, hence so does every source component.
At each marked orbit the admissible local form, in suitable completed root coordinates, is $x=u^{N_i}$ with $N_i=|I_i|$; its source point is smooth over $R'$. The marking lifts of the generic framed object
extend, after a further finite DVR extension, through their finite proper preimage schemes. The isomorphism on the generic fiber is part of the extension data; it does not assert a canonical trivialization of the special fiber.
\end{proof}

\begin{lemma}\label{lem:nodalnormal}
Let $R$ be an excellent characteristic-zero DVR and $p:C\to\operatorname{Spec}R$ a flat nodal-curve family of finite type with smooth generic fiber. Then the total space $C$ is normal
(componentwise if disconnected). The assertion remains true after any finite extension of DVRs.
\end{lemma}
\begin{proof}
At a point of the relative smooth locus the total space is regular. At a closed source node, after an étale extension splitting the node,
its complete local ring has the form $\widehat{\mathcal O}_{C,x}\simeq \widehat{R'}[[u,v]]/(uv-c),\; 0\ne c\in\mathfrak m_{R'}$, by \cite[Tag~0H8M]{Stacks}. The generic fiber is smooth, so $c\ne0$;
writing $c=\epsilon t^a$ and rescaling one variable reduces to $uv=t^a$ with $a\ge1$. This hypersurface is Cohen--Macaulay, hence satisfies $(S_2)$. Its only possible nonregular point is the
closed node $u=v=t=0$ (of codimension two); consequently it satisfies $(R_1)$. By Serre's criterion it is normal. Excellence and faithful flatness of completion allow descent of normality to the local ring.
The same argument applies after finite DVR base change. In particular, $uv=t^a$ for $a>1$ is generally singular but normal; cf.
\cite[Tag~0CDC]{Stacks}.
\end{proof}

\begin{lemma}\label{lem:directextension}
In the setting of Lemma~\ref{lem:nodalnormal}, put $U=C^{\rm sm}_{/R}$ and $j:U\hookrightarrow C$. For every finite-rank vector bundle $F^\circ$ on $U$, $F=j_*F^\circ$ is the unique coherent reflexive extension of $F^\circ$ to $C$;
it is $R$-flat and restricts to $F^\circ$ on $U$. An $H$-linearization of $F^\circ$ extends uniquely to $F$. No local freeness of $F$ at a source node is asserted. After any finite DVR base change the same
construction applies on the new normal total space to the pulled-back bundle on its relative smooth locus.
\end{lemma}
\begin{proof}
The closed complement $C\setminus U$ consists of finitely many special-fiber nodes, each of codimension two in the normal total surface. On each integral component, the depth at these points is
at least two. Thus the extension theorem for reflexive sheaves
\cite[Tag~0EBJ]{Stacks} gives
\[
 F=j_*F^\circ,\qquad F\simeq F^{\vee\vee},\qquad
 j^*F\simeq F^\circ.
\]
Normality implies $F$ is torsion-free. Since all components of $C$ dominate $\operatorname{Spec}R$ and $C/R$ is flat, the DVR uniformizer is a non-zero-divisor on $F$. Thus $F$ is torsion-free
as an $R$-module and hence $R$-flat. The same uniqueness extends $H$-linearization maps and their cocycle identities. For $R\to R'$, apply Lemma~\ref{lem:nodalnormal} to $C_{R'}$ and
repeat the reflexive-extension construction there. This construction is \emph{not} asserted to commute with arbitrary base change across nodes; for the canonical coordinate bundle its restriction to the
relative smooth locus is the ordinary base change of the associated bundle by Lemma~\ref{lem:canonicalfinite}.
\end{proof}

\subsection{Integral currents and a weighted logarithmic lattice}
\begin{lemma}\label{lem:symparts}
Let $O$ be a DVR with $|H|\in O^\times$, and let $p:C\to\operatorname{Spec}O$ be a projective, flat family of nodal curves equipped with an $O$-linear action of the finite group $H$.
Let $F$ be an $H$-equivariant coherent $\cO_C$-module. Suppose that $Q_1,\ldots,Q_s$ are pairwise disjoint, $H$-invariant, effective Cartier divisors, finite and flat over $O$, contained in the relative
smooth locus, and that $F$ is locally free near $Q=\sum_iQ_i$. For each geometric irreducible component $C_\alpha$ of the closed fiber, assume $\deg\bigl(\cO_C(Q)|_{C_\alpha}\bigr)>0$.
Then there is an $m\ge0$, independent of $i$ and $n$, such that, with $D=mQ$, for all $i$ and $n\ge0$, $H^0\bigl(C,F(D+nQ_i)\bigr)^H \twoheadrightarrow H^0\bigl(C,F(D+nQ_i)/F(D)\bigr)^H$.
\end{lemma}
\begin{proof}
By the componentwise degree criterion
\cite[Tag~0AYQ, Lemma~33.44.15]{Stacks},
$\cO_C(Q)$ is ample on the geometric closed fiber and hence $p$-ample over the DVR. Relative Serre vanishing
\cite[Tag~02O1]{Stacks} gives $m_0\ge0$ such that
\[
 R^1p_*F(mQ)=0\qquad(m\ge m_0).
\]
Since $\operatorname{Spec}O$ is affine, $H^1(C,F(mQ))=0$. For $D=mQ$, local freeness of $F$ along $Q_i$ gives, for every $n\ge0$, the exact sequence
\[
 0\longrightarrow F(D)\longrightarrow F(D+nQ_i)
 \longrightarrow F(D+nQ_i)/F(D)\longrightarrow0.
\]
Indeed the inclusion is an isomorphism off $Q_i$ and is multiplication by a non-zero-divisor in a local trivialization near $Q_i$. The cohomology sequence therefore yields the asserted surjection
before invariants. As $|H|\in O^\times$, the Reynolds projector
\[
 \mathsf{Av}_H=|H|^{-1}\sum_{h\in H}h
\]
shows that $(\cdot)^H$ is exact on $O[H]$-modules. Thus the surjection remains surjective after taking invariants. For an $H$-orbit of markings, a stabilizer-invariant jet at one
lift specifies an $H$-invariant jet on the complete orbit, to which the same surjection applies.
\end{proof}

For $b\in V_d$ homogeneous and $z-z_*=q\xi$, the scaling part of the state-field covariance gives
\begin{equation}\label{eq:linearcurrent}
 \xi^nY_\xi(b,\xi)\,d\xi
 =q^{d-n-1}(z-z_*)^nY_z(b,z)\,dz,
 \qquad \Delta=d-n-1.
\end{equation}
Here $n$ is an allowed mode index (possibly fractional in a chosen twisted root coordinate); the powers are interpreted on the chosen finite ramified base. The residue mode has conformal-weight shift
$\Delta$, and its physical coefficient is $q^\Delta$. The $L(0)$-homogeneity used here concerns the state $b$; relative module weights enter only through the corresponding mode indices. Apply the identity to each homogeneous summand of a finite
Laurent coordinate gauge, without identifying the full gauge with a scalar dilation.

\begin{lemma}\label{lem:weightednormal}\label{prop:filteredderivative}
\label{lem:signedweights}
Let $O=L[[t]]$, $K=L((t))$, and let $(\cQ_K,\nabla)$ be a finite-dimensional differential module over $K$ arising from a finite tensor product of ordinary twisted-module blocks as in
Lemma~\ref{lem:descent}. Fix $R,N$ and the homogeneous PBW spanning families $\mathcal W_N(M_i)$ of that lemma. Choose $u_1,\ldots,u_\ell$ as a basis of $\bigotimes_iK_{R,N}(M_i)$ homogeneous in the individual external
weights. Let $\nu$ be an integer-valued function of this external weight tuple. Thus $\nu$ is unchanged by every weight-homogeneous PBW rewrite. Suppose that every PBW tensor $w$ containing a leftmost
mode of depth $\ge R$ admits a finite relation in $\cQ_K$,
\begin{equation}\label{eq:B}
 [w]=\sum_b c_b[w_b],\qquad
 \ee(w_b)<\ee(w),\qquad
 \vval_t(c_b)\ge\nu(w_b)-\nu(w),
\end{equation}
with $c_b\in K$. Suppose moreover that for every such homogeneous PBW tensor, including the low-depth tensors,
\begin{equation}\label{eq:filteredderivative}
 \nabla_{t\partial_t}[w]=\sum_a d_a[w_a],\qquad
 \vval_t(d_a)\ge\nu(w_a)-\nu(w),
\end{equation}
where each sum is finite and $d_a\in K$. Then every homogeneous PBW tensor $w$ has a finite expansion
\begin{equation}\label{eq:normal}
 [w]=\sum_{j=1}^\ell a_j[u_j],\qquad
 \vval_t(a_j)\ge\nu(u_j)-\nu(w),
\end{equation}
and
\begin{equation}\label{eq:weightedmember}
 L_\nu=\sum_{j=1}^\ell O\,t^{\nu(u_j)}[u_j],\qquad
 t^{\nu(w)}[w]\in L_\nu,
 \qquad \nabla_{t\partial_t}L_\nu\subseteq L_\nu.
\end{equation}
Thus $L_\nu$ is a finite free full $O$-lattice in $\cQ_K$; no lower bound on $\nu$ is required.
\end{lemma}
\begin{proof}
Apply strong induction on $\ee(w)\in\N$. For a low-depth word, expand it in the $u_j$ of the same external weight tuple; the coefficients lie in $\C$ and have valuation $0$,
while all terms have the same $\nu$. Otherwise apply \eqref{eq:B} and expand each $w_b$ in the fixed PBW families. Each term has the same energy and
$\nu$ as $w_b$; therefore the induction hypothesis applies and $\vval_t(c_ba_{bj}) \ge\nu(w_b)-\nu(w)+\nu(u_j)-\nu(w_b) =\nu(u_j)-\nu(w)$. The non-Archimedean valuation inequality preserves this bound after
summing over $b$, proving \eqref{eq:normal} and the membership in
\eqref{eq:weightedmember}.
The classes $[u_j]$ span $\cQ_K$ by Lemma~\ref{lem:descent}; thus $K L_\nu=\cQ_K$. As a finitely generated $O$-submodule of a $K$-vector space, $L_\nu$ is torsion-free, hence free over the DVR $O$. For any homogeneous $w$,
\[
 \begin{aligned}
 \nabla_{t\partial_t}\bigl(t^{\nu(w)}[w]\bigr)
 &=\nu(w)t^{\nu(w)}[w]
 +\sum_a t^{\nu(w)-\nu(w_a)}d_a\,
             t^{\nu(w_a)}[w_a]
 \in L_\nu,
 \end{aligned}
\]
by \eqref{eq:filteredderivative} and the membership already proved. Taking $w=u_j$ gives stability of the finite generating set and hence of $L_\nu$.
\end{proof}

\begin{lemma}\label{lem:DVRdescent}
Let $K_0=L_0((q))$ with $\operatorname{char}L_0=0$, and let $K/K_0$ be a finite field extension. If a finite-dimensional differential module $(M,\nabla)$ over $K_0$ becomes regular singular over $K$, then it is regular singular over $K_0$.
\end{lemma}
\begin{proof}
Put $O_0=L_0[[q]]$, let $O_K$ be its integral closure in $K$, and set $\delta=q\partial_q$. The extension is finite separable, $O_K$ is a finite DVR over $O_0$, and $\delta$ extends uniquely to $K$. In equicharacteristic zero this
extension is logarithmic: $\delta(O_K)\subseteq O_K$. Indeed, after a finite extension of the residue-field constants, $q=u(t)t^e$ with $u(t)\in L[[t]]^\times$, and $\delta(t)/t =\bigl(e+t\,u'(t)/u(t)\bigr)^{-1}\in O_K^\times$.
Regular singularity over $K$ therefore supplies a full $O_K$-lattice $\Lambda\subset M_K=M\otimes_{K_0}K$ with $\nabla_\delta\Lambda\subseteq\Lambda$. Define the $K_0$-linear map
\[
 \operatorname{Tr}_M
 =\id_M\otimes\operatorname{Tr}_{K/K_0}:M_K\to M,
 \qquad \Lambda_0=\operatorname{Tr}_M(\Lambda).
\]
Since $\Lambda$ is finite over $O_K$ and $O_K$ is finite over $O_0$, $\Lambda_0$ is a finite $O_0$-module. It is torsion-free and spans $M$ over $K_0$, because $\operatorname{Tr}_{K/K_0}\ne0$ for a finite
separable extension. Moreover the trace commutes with the unique extension of $\delta$:
\[
 \nabla_\delta\operatorname{Tr}_M(x)
 =\operatorname{Tr}_M(\nabla_\delta x),\qquad
 \nabla_\delta\Lambda_0\subseteq\Lambda_0.
\]
Hence $\Lambda_0$ is a full logarithmically stable $O_0$-lattice, which is the formal regular-singular lattice criterion
\cite[Chapter~II]{Deligne}.
\end{proof}

\begin{theorem}\label{thm:canonicalregular}
Let $S$ be the framed Hurwitz chart of Lemma~\ref{lem:framedcharts}, $\cQ$ the finite locally free differential Ward quotient of Theorem~\ref{thm:finiteconn}, and $\nabla$ its algebraic integrable
connection. Then $(\cQ,\nabla)$ is regular singular. For a smooth proper compactification $\widetilde S$ with normal-crossing boundary, it admits an algebraic logarithmic extension on $\widetilde S$.
\end{theorem}
\begin{proof}
Write $U=F_m(\C^*)$, let $\overline U$ be the wonderful model of Definition~\ref{def:clustercompactification}, and set $\overline S=\operatorname{Norm}_{K(S)}(\overline U)$. By
\cite[Theorems~1.2--1.3]{Li} and Lemmas~\ref{lem:clustercharts}
and \ref{lem:primes}, $\overline S$ is normal and proper, and the codimension-one points of $\overline S\setminus S$ lie over precisely the free finite, $0$-anchored and $\infty$-anchored cluster divisors.
The possible intersections of cluster divisors have codimension $\ge2$ on $\overline U$ and do not constitute additional boundary primes. We prove regularity at every boundary prime of $\overline S$.

Fix such a prime $\eta$, put
\[
 O_0=\widehat{\mathcal O}_{\overline S,\eta},\qquad
 K_0=\operatorname{Frac}(O_0),\qquad
 (M_0,\nabla)=\cQ_{\eta_S}\otimes_{K(S)}K_0.
\]
Here $\eta_S$ is the generic point of $S$; the tensor product is taken using $K(S)\hookrightarrow K_0$. Lemmas~\ref{lem:actualcoverextension} and \ref{lem:pureramification} give, after a finite extension of
complete DVRs, $O=L[[t]]\supset O_0$, $q=t^e$, $e\ge1$, with $\delta=t\partial_t$ and $\delta|_L=0$, together with an extension of the actual generic framed $H$-cover to a projective nodal
admissible $H$-cover $p:C\to\operatorname{Spec}O$. Here $L$ includes the tangential generic parameters, which are $\delta$-constant. Every irreducible component of $C_0$ meets an
original marked orbit. In particular, if $P_i$ denotes the reduced complete $H$-orbit of the $i$th extended marking and $P=\sum_iP_i$, then
\[
 \deg(\cO_C(P)|_{C_\alpha})>0
 \qquad(C_\alpha\subset C_0\text{ geometrically irreducible}).
\]
The $P_i$ are disjoint invariant Cartier divisors in the relative smooth locus. Since $\operatorname{char}L=0$, $|H|\in O^\times$.

Let $F^\circ$ be the finite state-current bundle of Lemma~\ref{lem:canonicalfinite} on $C^{\mathrm{sm}}_{/O}$ and put $F=j_*F^\circ$, with $j:C^{\mathrm{sm}}_{/O}\hookrightarrow C$. By Lemmas~\ref{lem:nodalnormal} and \ref{lem:directextension},
$F$ is a coherent $H$-equivariant reflexive sheaf, $O$-flat and locally free near $P$. Apply Lemma~\ref{lem:symparts} to $Q_i=P_i$, $Q=P$, and $F$. There is an $m\ge0$ for which
\begin{equation}\label{eq:511parts}
 \Gamma(C,F(mP+nP_i))^H\twoheadrightarrow
 \Gamma\bigl(C,F(mP+nP_i)/F(mP)\bigr)^H
 \qquad(i,\ n\ge0).
\end{equation}
This is a statement about the \emph{extended integral current sheaf}; no regular-singular property of $M_0$ has been used.

Choose component coordinates at the marked sections. The physical coarse coordinates $x_{i,\mathrm{phys}}$ and their component coordinates $x_{i,\mathrm{comp}}$ satisfy
\begin{equation}\label{eq:511scales}
 x_{i,\mathrm{phys}}=q^{a_i}x_{i,\mathrm{comp}},\qquad
 \begin{array}{c|c|c}
 \text{boundary}&a_i\ (i\text{ an input})&a_0\ (\text{output at }\infty)\\ \hline
 \text{free finite cluster }A&1\ (i\in A),\ 0\ (i\notin A)&0\\
 \text{$0$-anchored cluster }A&1\ (i\in A),\ 0\ (i\notin A)&0\\
 \text{$\infty$-anchored cluster }A&-1\ (i\in A),\ 0\ (i\notin A)&1
 \end{array}
\end{equation}
with $a_i\in\{-1,0,1\}$. In the second row the fixed input at $0$ also has $a_i=1$; in the third row the component output coordinate is $s/q$, where $s=z^{-1}$. These are centered coordinates:
for an escaping input $x_{i,\mathrm{comp}}=q(z-p_i)$, while for a finite bubble input $x_{i,\mathrm{comp}}=(z-p_i)/q$. After further finite ramification, $ea_i/N_i\in\Z$ for each inertia
order $N_i$. Choose the component roots with the phases inherited from the fixed generic framing. The identities
\[
 x_{i,\mathrm{comp}}=u_{i,\mathrm{comp}}^{N_i},\qquad
 u_{i,\mathrm{phys}}=t^{ea_i/N_i}u_{i,\mathrm{comp}}
\]
then hold on the generic marked formal disks. In particular, they have no extra $t$-dependent unit. When $a_i<0$, only the component root is asserted to be integral over $O$.

We verify both hypotheses of Lemma~\ref{lem:signedweights}. Use the fixed finite $L(j)$-stable state space $E$ and the fixed homogeneous representatives $B_{g_i}$ of \eqref{eq:PBW}.
The full state-current coordinate matrices and their inverses at $P_i$ are integral in component root coordinates by Lemma~\ref{lem:canonicalfinite}; their Laurent root orders are bounded by constants depending on $E$ and $N_i$, not on the
prescribed mode depth. Consequently there is a common $C_0\in\Z$ such that for any $s\in\Gamma(C,F(mP))$ and every nonzero homogeneous mode $b_\mu$ occurring in its component-coordinate residue, $\Delta(b_\mu):=\wt b-\mu-1\le C_0$.
Indeed, near $P_i$ the fixed pole bound of $F(mP)$ and the finite Laurent bounds of its two frame changes bound the negative powers of $u_{i,\mathrm{comp}}$; the corresponding mode degree tends to
$-\infty$ as the root exponent increases. All coefficients of these modes belong to $O$. The bound is independent of $s$, $i$ and the external tensor on which its residue is evaluated.

Choose $R$ so that $\Delta_*:=\wt b-\alpha_b+n-1>\max\{0,C_0\}\; (b\in B_{g_i},\ n\ge R,\ \text{all }i)$. Fix a high-depth mode $b_{\alpha_b-n}$ at orbit $P_i$. Its pure component-coordinate current gives an invariant principal
part modulo $F(mP)$ on $P_i$, with zero principal parts at every other $P_j$; the inertia eigenvalue condition of Definition~\ref{def:canonicalmarked} is satisfied. For a sufficiently
large $k$ (depending on $n$), \eqref{eq:511parts} lifts this principal part to $s_{i,b,n}\in\Gamma(C,F(mP+kP_i))^H$. The normalized residue formula of Proposition~\ref{prop:canonicalrealization} gives, in \emph{component}
coordinates, $R_i(s_{i,b,n})=b_{\alpha_b-n}+U_i,\; R_j(s_{i,b,n})=U_j\; (j\ne i)$, where each homogeneous mode term $b'_{\mu'}$ of each $U_j$ obeys $\Delta'=\wt b'-\mu'-1\le C_0$ and has a coefficient
$\beta_{j,b',\mu'}(t)\in O$. The coefficient of the target mode is $1$ by the $1/N_i$ normalized orbit residue. The series $U_j$ act by finite sums on each fixed tensor, by lower truncation; no
uniformly finite sum of operators on the whole module is claimed. The generic fiber of $s_{i,b,n}$ is an invariant meromorphic current on the original cover; therefore its Ward relation vanishes in
$M_0\otimes_{K_0}K$, by the definition of the differential quotient.

Fix a homogeneous PBW tensor $w=b_{\alpha_b-n}^{(i)}v$ with the leftmost high-depth factor removed, as in the proof of Proposition~\ref{prop:seedfinite}. The current relation in the
differential quotient, first in component frames and then in the physical frames \eqref{eq:511scales}, gives
\begin{equation}\label{eq:511descent}
 [w]=\sum_{j,b',\mu'}c_{j,b',\mu'}(t)
                  [b_{\mu'}'^{(j)}v],\qquad
 c_{j,b',\mu'}(t)
 =-\beta_{j,b',\mu'}(t)\,
 t^{e(a_j\Delta'-a_i\Delta_*)}.
\end{equation}
Only finitely many terms act nontrivially on $v$. The scaling exponent is precisely \eqref{eq:linearcurrent}, applied to each homogeneous mode; it does not omit the positive-Virasoro
terms, which have already been included in the bounded integral component gauges. Let $D>0$ clear all relative grading denominators, so that the energy $\ee$ of Section~\ref{sec:moduleprelim}
is in $\N$. For $w_{j,b',\mu'}=b_{\mu'}'^{(j)}v$, $\ee(w_{j,b',\mu'})-\ee(w) =D(\Delta'-\Delta_*)<0$. For each fixed simple external block, put
\begin{equation}\label{eq:511nu}
 \nu(w)=e\sum_{k=0}^m a_k(\wt w_k-h_k),
\end{equation}
where $h_k$ is the lowest $L(0)$-weight on its $k$th factor. Our finite ramification makes $\nu(w)\in\Z$, since $\wt w_k-h_k\in N_k^{-1}\Z$. Because the tensor $v$ is the same in the target and in every correction term,
\begin{equation}\label{eq:511valuation}
 \begin{aligned}
 \nu(w_{j,b',\mu'})-\nu(w)
     &=e(a_j\Delta'-a_i\Delta_*),\\
 \vval_t(c_{j,b',\mu'})
     &=\vval_t(\beta_{j,b',\mu'})
       +e(a_j\Delta'-a_i\Delta_*)
       \ge\nu(w_{j,b',\mu'})-\nu(w).
 \end{aligned}
\end{equation}
Homogeneous PBW rewrites preserve the individual external weights, and hence $\nu$. Equations \eqref{eq:511descent}--\eqref{eq:511valuation} are exactly hypothesis \eqref{eq:B}, with a common threshold $R$.

It remains to verify \eqref{eq:filteredderivative} for the
original placement connection, not for a coordinate-gauge
replacement. In fixed physical module frames, the chosen transverse derivation $\delta=t\partial_t$ has
\begin{equation}\label{eq:511derivative}
 \nabla_\delta[w]
 =\sum_{i=1}^m(\delta p_i)[L(-1)^{(i)}w],
 \qquad \delta|_L=0.
\end{equation}
Here all moving inputs are included, whereas the fixed input at $0$ and the output at $\infty$ contribute no placement derivative. For $p_i=z_*+t^e\xi_i$ (finite bubble), $\vval_t(\delta p_i)\ge e$ if $a_i=1$;
for $p_i=t^{-e}\xi_i$ (escaping input), $\vval_t(\delta p_i)=-e$ and $a_i=-1$; for an unscaled leg, either $\delta p_i=0$ or $\vval_t(\delta p_i)\ge0=e a_i$. Since $[L(0),L(-1)]=L(-1)$, every nonzero term $w_i'=L(-1)^{(i)}w$ satisfies
\[
 \nu(w_i')-\nu(w)=ea_i,
 \qquad
 \vval_t(\delta p_i)\ge\nu(w_i')-\nu(w).
\]
Rewrite each $w_i'$ into the fixed PBW families; this is a finite weight-homogeneous rewrite and does not change $\nu(w_i')$. Thus \eqref{eq:511derivative} verifies the coefficientwise hypothesis \eqref{eq:filteredderivative}.

By Lemma~\ref{lem:signedweights}, the extended differential module $M_0\otimes_{K_0}K$ has a full $O$-lattice $\Lambda$ with $\nabla_\delta\Lambda\subseteq\Lambda$. The lowest weights $h_i$
contribute only scalar monomial factors; the integral lattice uses the relative weights in \eqref{eq:511nu}. Since $q=t^e$ and $\delta|_L=0$, one has $q\partial_q=e^{-1}\delta$ on $K$ (with the chosen coefficient
field); thus $\nabla_{q\partial_q}\Lambda\subseteq\Lambda$ and regular singularity follows over $K$. Lemma~\ref{lem:DVRdescent} then gives regular singularity at $\eta$ over $K_0$.

Finally, $\cQ$ is a finite-rank algebraic vector bundle with integrable connection on the smooth variety $S$ by Theorem~\ref{thm:finiteconn}. We have checked every boundary prime of the \emph{fixed normal proper} $\overline S$. Therefore
\cite[Theorem~1.1.1, (ii)$\Rightarrow$(i)]{Andre}
applies to $S$ and $(\cQ,\nabla)$; it also covers valuations exceptional over $\overline S$. The existence of a logarithmic extension on a smooth proper normal-crossing compactification follows from Deligne's regular-singular extension theorem
\cite[Chapter~II, Sections~4--5]{Deligne}.
\end{proof}

\subsection{Formal expansions and convergence of the original weight sums}

Let $\tau$ be a finite rooted binary tree with root vertex $v_0$, internal-vertex set $\mathcal V(\tau)$ and leaf set $I$. For $v\in\mathcal V(\tau)$ write $I_v\subseteq I$ for its
leaf descendants. Attach a typed twisted intertwining operator to each internal vertex and fix homogeneous $w_i\in(W_i)_{\alpha_i}$ $(i\in I)$ and $w_0'\in(W_0')_\gamma$. Every internal edge carries an ordinary grading-restricted module;
if it is not simple, decompose it into its finitely many simple summands. At an internal vertex $v$ let $r_v$ denote its nonzero formal displacement. Put
\[
 R=r_{v_0},\qquad
 u_v=\frac{r_v}{r_{\operatorname{par}(v)}}
 \quad(v\in\mathcal V(\tau)\setminus\{v_0\}),
\]
and regard $R$ and the $u_v$ as independent formal variables. The internal sums are taken over homogeneous projections of the intermediate modules, in the indicated tree order.

\begin{lemma}\label{lem:treesupport}
For fixed simple summands on the internal edges, let $\mu_v=h_v+n_v/N_v$, $n_v\in\Z_{\ge0}$, be the weight on the edge from $v\ne v_0$ to $\operatorname{par}(v)$; an absent homogeneous space contributes zero. There are coefficients
$C_{\boldsymbol n}\in\C$ such that the original projected matrix coefficient has the coefficientwise identity
\begin{equation}\label{eq:treeweightseries}
 \Phi_\tau
 =R^{\gamma-\sum_{i\in I}\alpha_i}
  \sum_{\boldsymbol n\in\Z_{\ge0}^{\mathcal V(\tau)\setminus\{v_0\}}}
  C_{\boldsymbol n}
  \prod_{v\ne v_0}u_v^{\,h_v+n_v/N_v-\sum_{i\in I_v}\alpha_i}.
\end{equation}
Each $C_{\boldsymbol n}$ is a finite contraction of finite-dimensional weight spaces. Summing over the finitely many simple choices yields a finite sum of series of the form
\[
 R^{\gamma-\sum_i\alpha_i}
 u^{\lambda}F(u_1^{1/N},\ldots,u_d^{1/N}),
 \qquad F\in\C[[u_1^{1/N},\ldots,u_d^{1/N}]],
\]
where $d=|\mathcal V(\tau)|-1$ and $N$ is a common multiple of the finitely many $N_v$. The expression has finite exponent classes, no logarithms and a well-defined coefficient for every monomial. No analytic convergence is asserted.
\end{lemma}
\begin{proof}
By Lemma~\ref{lem:IOlogs}, an intertwining-operator coefficient with homogeneous input weights $\beta_1,\beta_2$ and output weight $\beta_0$ is proportional to $r_v^{\beta_0-\beta_1-\beta_2}$, without logarithms.
Taking the product over $v\in\mathcal V(\tau)$ and cancelling every intermediate weight once with sign $+$ and once with sign $-$ gives the exact identity
\[
 \prod_{v\in\mathcal V(\tau)}
 r_v^{\mu_v-\mu_{v,1}-\mu_{v,2}}
 =R^{\gamma-\sum_{i\in I}\alpha_i}
 \prod_{v\ne v_0}
 u_v^{\mu_v-\sum_{i\in I_v}\alpha_i},
\]
where $\mu_{v_0}=\gamma$, and $\mu_{v,j}$ denotes the weight of the $j$th incoming edge (or the corresponding external leaf). For each tuple $\boldsymbol n$, all intermediate homogeneous
spaces have finite dimension. Their canonical projections and contractions consequently define $C_{\boldsymbol n}$ by a finite sum, even though the unprojected intermediate spaces are infinite dimensional. The map
$\boldsymbol n\mapsto(n_v/N_v)_{v\ne v_0}$ is injective. There are finitely many simple choices, each with a fixed lower weight $h_v$ and a discrete nonnegative $N_v^{-1}$-grading.
These observations prove \eqref{eq:treeweightseries} and its stated support properties.
\end{proof}

For the next assertion, put $A=\C\{q_1,\ldots,q_r,y_1,\ldots,y_s\}$, $\widehat A=\C[[q_1,\ldots,q_r,y_1,\ldots,y_s]]$ and $Q_D=q_1\cdots q_r$. An \emph{adapted ratio corner} is a chart with boundary $Q_D=0$ in which the variables occurring in
\eqref{eq:treeweightseries} satisfy, after a finite root base change,
\begin{equation}\label{eq:adaptedratio}
 u_v=q^{m_v}\epsilon_v(q,y),\qquad
 m_v\in\Z_{\ge0}^r\setminus\{0\},\quad
 \epsilon_v\in A^\times,\qquad
 q^{m_v}=\prod_{j=1}^r q_j^{m_{vj}}.
\end{equation}
In particular, $u_v=q_v$ is adapted. The nonzero-row condition is essential for coefficientwise substitution of an arbitrary formal power series in the $u_v$: substitution at $u_v(0,y)\ne0$ is not
available before convergence has been proved. Choose the branches of $\epsilon_v^\lambda$ by their nonzero constant terms, formally near the corner. After extracting the finitely many base monomials and making a
common root substitution, \eqref{eq:adaptedratio} maps every series of Lemma~\ref{lem:treesupport} into some $\mathscr N_{\Lambda,b}$ of \eqref{eq:Nilsson}, with $b=0$. Indeed, for fixed $\boldsymbol k\in\Z_{\ge0}^r$, a contribution from
$\boldsymbol n\in\Z_{\ge0}^d$ requires
\[
 \sum_v n_v m_{vj}\le N k_j\quad(1\le j\le r),\qquad
 \sum_v n_v\sum_{j=1}^r m_{vj}\le N\sum_{j=1}^rk_j.
\]
Since $\sum_jm_{vj}\ge1$ for all $v$, only finitely many $\boldsymbol n$ satisfy these inequalities. The units $\epsilon_v^{h_v+\cdots}$ have well-defined formal Taylor expansions after choosing a value of their constant terms.
The finite power series in the transverse variables are interpreted in the $(q,y)$-adic completion.

\begin{lemma}\label{lem:jointsupport}
Let $U$ be an adapted ratio corner, let $(\cQ,\nabla)$ be the finite locally free algebraic differential Ward quotient on $U^*=U\setminus\{Q_D=0\}$, and suppose it has a logarithmic extension $\overline{\cQ}$ near the corner.
Let $\mathbb K$ be a common formal coefficient extension admitting the injections of $\widehat A[Q_D^{-1}]$ and of the finitely many generalized-power monomials under consideration. Suppose that
\[
 \Phi\in\operatorname{Hom}_{\mathcal O(U^*)}
 (\cQ|_{U^*},\mathbb K),\qquad \nabla^\vee\Phi=0,
\]
and that every evaluation $\Phi([w])$, for a constant external state $w$ and its associated section $[w]$ of the Ward quotient, has the coefficientwise joint support of Lemma~\ref{lem:treesupport} after \eqref{eq:adaptedratio}.
Then the coordinate row of $\Phi$ in a local logarithmic frame of $\overline{\cQ}$ belongs to $\mathscr N_{\Lambda,b}^{\oplus\operatorname{rk}\cQ}$ for some finite $\Lambda\subset\C^r$, $b\ge0$.
\end{lemma}
\begin{proof}
Work on a sufficiently small algebraic affine neighbourhood of the corner and write $s=\operatorname{rk}\cQ$. Since $\cQ$ is a coherent quotient of the free module on external tensor states,
Noetherianity yields finitely many constant states $w_1,\ldots,w_p$ whose images $\sigma_j=[w_j]$ generate on a possibly smaller $U^*$. Choose an algebraic logarithmic extension by Theorem~\ref{thm:canonicalregular} and
\cite[Chapter~II, Proposition~5.4]{Deligne}.
Each $\sigma_j$ is a rational section of this extension, regular off the coordinate boundary. Thus some common $a\ge0$ makes $Q_D^a\sigma_j$ holomorphic at the corner. In a local logarithmic frame form the $s\times p$ holomorphic
matrix $B$ with columns $Q_D^a\sigma_j$. Then $\operatorname{rank}B(q,y)=s\; (Q_D(q)\ne0)$. The maximal-minor ideal $I_s(B)\subset A$ therefore satisfies $V(I_s(B))\subseteq V(Q_D)$. By the analytic Nullstellensatz,
$Q_D^c\in I_s(B)$ for some $c\ge0$. In particular, for suitable $c_J\in A$ and $s\times s$ minors $B_J$,
\[
 Q_D^c=\sum_Jc_J\det B_J,\qquad
 C=Q_D^{-c}\sum_Jc_J\iota_J\operatorname{adj}(B_J)
 \in\operatorname{Mat}_{p\times s}(A[Q_D^{-1}]),\qquad
 BC=I_s.
\]
Here $\iota_J:\C^s\hookrightarrow\C^p$ inserts the indexed coordinates. If $h$ is the row of $\Phi$, then $h=(hB)C,\; (hB)_j=Q_D^a\Phi([w_j])$. Each entry of $hB$ lies in a finite Nilsson module by the
hypothesis and \eqref{eq:adaptedratio}. Multiplication by $C\in A[Q_D^{-1}]$ preserves finite exponent classes, finite coordinate pole orders and the uniform logarithmic-degree bound. Thus $h\in\mathscr N_{\Lambda,b}^{\oplus s}$, for suitable
$\Lambda,b$. The equation $\nabla^\vee\Phi=0$ is an equality of formal coefficients in the chosen injection into $\mathbb K$; no analytic convergence has been used.
\end{proof}

\begin{theorem}\label{thm:formalconv}
Let $\Delta^{r+s}$ be a sufficiently small complex polydisk with coordinates $q=(q_i)_{i=1}^r$ and $y=(y_j)_{j=1}^s$. Put $D=\{q_1\cdots q_r=0\}$. Let $(Q,\nabla)$ be a finite-rank integrable holomorphic connection on $\Delta^{r+s}\setminus D$
admitting a locally free logarithmic extension over $\Delta^{r+s}$. Let $f$ be a formal horizontal section of $Q$ or $Q^\vee$, expressed in a meromorphic analytic frame whose transition to a
logarithmic frame belongs to $\operatorname{GL}_{\operatorname{rk}Q}(A[Q_D^{-1}])$. If its coordinates lie in a finite module $\mathscr N_{\Lambda,b}$ of \eqref{eq:Nilsson}, then $f$ is the
formal expansion of an analytic horizontal section on any sufficiently small simply connected branch sector of $\Delta^{r+s}\setminus D$. Its generalized-power coefficient series converge absolutely and locally normally after a suitable
finite collection of exponent classes is factored out.
\end{theorem}
\begin{proof}
First choose Deligne's canonical logarithmic extension, with the residue strip fixed in Section~\ref{sec:setup}. Any logarithmic extension in the hypothesis and this canonical one have the same meromorphic localization along $D$; after shrinking, the identity on the punctured polydisk is represented by a matrix and its inverse in $A[Q_D^{-1}]$. Thus the change from the given meromorphic frame to a canonical-extension frame preserves the finite Nilsson condition.
By Deligne's local normal-crossing classification
\cite[Chapter~II, Proposition~5.4 and Remark~5.5(ii)]{Deligne},
in a suitable holomorphic frame of \emph{this canonical extension} the connection has the form
\begin{equation}\label{eq:delignelocalnormal}
 d-\sum_{i=1}^rR_i\frac{dq_i}{q_i},\qquad
 R_i\in\operatorname{Mat}_n(\C),\quad
 [R_i,R_j]=0,\quad n=\operatorname{rk}Q.
\end{equation}
This is a local statement near the normal-crossing point, with contractible transverse parameters $y$; it uses the chosen canonical extension, normal crossings and integrability. The comparison with the original logarithmic lattice is only required to be meromorphic with finite boundary pole orders, not holomorphic in a resonant lattice. Hence the transformed formal horizontal section
$\widetilde f$ still has finitely many exponent classes, bounded logarithmic degrees and finite coordinate pole orders. Choose branches of $\log q_i$. Since $R_i$ commute, the matrix
\[
 P(q)=\prod_{i=1}^rq_i^{R_i}
 =\prod_i\exp(R_i\log q_i)
\]
is a finite sum of terms $q^\lambda(\log q)^\beta$ with constant matrix coefficients, where the logarithmic degrees are bounded by the nilpotence indices of the Jordan blocks of the $R_i$. For $c=P(q)^{-1}\widetilde f$, horizontality and
\eqref{eq:delignelocalnormal} give
\[
 \partial_{q_i}c=0\quad(1\le i\le r),\qquad
 \partial_{y_j}c=0\quad(1\le j\le s)
\]
as formal identities in a finite Nilsson module. To identify the kernel of these derivations, group monomials by
classes in $\C^r/\Z^r$. For one variable, $\partial_q\bigl(q^{\lambda+m}P_0(\log q)\bigr) =q^{\lambda+m-1} \bigl((\lambda+m)P_0+P_0'\bigr)(\log q)$. For a polynomial $P_0$, the identity $(\lambda+m)P_0+P_0'=0$ forces $P_0=0$ unless
$\lambda+m=0$, in which case $P_0$ is constant. The Laurent lower bound permits coefficient comparison in each exponent class. Repeating for $q_1,\ldots,q_r$ and comparing the ordinary coefficients of the $y_j$ gives $c\in\C^n$. Thus
\[
 \widetilde f=P(q)c,\qquad c\in\C^n.
\]
This is an analytic horizontal section on each fixed branch sector, and its formal expression is the prescribed one. The inverse logarithmic-frame gauge is holomorphic after multiplication by a sufficiently high power of $Q_D$.
Its convergent Laurent expansion, multiplied by the finite power--logarithm expression $P(q)c$, gives absolutely and locally normally convergent coefficient series on smaller branch sectors.
The dual connection has the same properties, with residue matrices $-R_i^{\mathsf T}$.
\end{proof}

\begin{lemma}\label{lem:Reinhardt}
Let $\Omega\subset\C^r$ be a complete Reinhardt domain with $0\in\Omega$, and let $f\in\mathcal O(\Omega)$. Put
\[
 a_{\boldsymbol n}
 =\frac{1}{\boldsymbol n!}
  (\partial_u^{\boldsymbol n}f)(0),\qquad
 \boldsymbol n\in\Z_{\ge0}^r.
\]
Here $\boldsymbol n!=\prod_i n_i!$ and $\partial_u^{\boldsymbol n}=\prod_i\partial_{u_i}^{n_i}$. Then
\begin{equation}\label{eq:Reinhardtnormal}
 f(u)=\sum_{\boldsymbol n\in\Z_{\ge0}^r}
 a_{\boldsymbol n}u^{\boldsymbol n},\qquad
 \sum_{\boldsymbol n\in\Z_{\ge0}^r}
 \sup_{u\in K}|a_{\boldsymbol n}u^{\boldsymbol n}|<\infty
 \quad(K\Subset\Omega).
\end{equation}
\end{lemma}
\begin{proof}
For $p\in\Omega$, choose $\rho_i>|p_i|$, $\rho_i>0$, so that $\overline P_\rho=\{u:|u_i|\le\rho_i\}\Subset\Omega$; this is possible by openness and completeness. Let $M_\rho=\max_{\overline P_\rho}|f|$. The polydisk Cauchy formula yields
\[
 a_{\boldsymbol n}
 =\frac{1}{(2\pi i)^r}
 \int_{|\zeta_1|=\rho_1}\!\cdots\!
 \int_{|\zeta_r|=\rho_r}
 \frac{f(\zeta)}{\prod_i\zeta_i^{n_i+1}}
 \,d\zeta_r\cdots d\zeta_1,\qquad
 |a_{\boldsymbol n}|\le M_\rho\rho^{-\boldsymbol n}.
\]
For $0<\sigma_i<\rho_i$,
\[
 \sum_{\boldsymbol n\ge0}
 \sup_{u\in\overline P_\sigma}
 |a_{\boldsymbol n}u^{\boldsymbol n}|
 \le M_\rho\prod_{i=1}^r(1-\sigma_i/\rho_i)^{-1}<\infty.
\]
The series equals $f$ near $0$, hence on $P_\rho$ by the identity theorem. A compact $K\Subset\Omega$ is contained in a finite union of such smaller polydisks $P_{\sigma}$, giving
\eqref{eq:Reinhardtnormal}.
\end{proof}

\begin{theorem}\label{thm:globalratio}
Let $\Omega\subset\C^r$ be a connected complete Reinhardt domain, $0\in\Omega$, and set
\[
 D_i=\Omega\cap\{u_i=0\},\qquad
 D=\bigcup_{i=1}^rD_i,\qquad
 \Omega^*=\Omega\setminus D.
\]
Let $(Q,\nabla)$ be a finite-rank integrable holomorphic connection on $\Omega^*$, regular singular along $D$, with Deligne canonical logarithmic extension $(\overline Q,\overline\nabla)$ on $\Omega$. Write $\overline Q(*D)$ for the sheaf of
meromorphic sections of $\overline Q$ with poles supported on $D$. Let $\widehat\Phi$ be a formal horizontal section of $Q^\vee$ at $0$ whose coordinates in a local logarithmic frame belong to $\mathscr N_{\Lambda,b}$ of \eqref{eq:Nilsson}.
For $1\le j\le m$, suppose $q_j\in\Gamma(\Omega^*,Q),\; q_j\in\Gamma(\Omega,\overline Q(*D))$. Then $\widehat\Phi$ determines a multivalued analytic horizontal section $\Phi^{\mathrm{an}}$ of $Q^\vee$ on $\Omega^*$.
For any simply connected branch domain $B\subset\Omega^*$ with fixed branches of $\log u_i$, the scalar $f_j=\Phi^{\mathrm{an}}(q_j)|_B$ admits an expansion
\begin{equation}\label{eq:globalratioexpansion}
 f_j(u)=\sum_{\lambda\in\Lambda_j}
 \sum_{0\le\beta_i\le b_j}
 u^\lambda(\log u)^\beta
 \sum_{\boldsymbol n\in\Z_{\ge0}^r}
 c_{j,\lambda,\beta,\boldsymbol n}u^{\boldsymbol n-\boldsymbol a},
\end{equation}
where $\Lambda_j\subset\C^r$ is finite, $b_j\in\Z_{\ge0}$, $\boldsymbol a\in\Z_{\ge0}^r$, and each inner power series converges absolutely and normally on compact subsets of $\Omega$. In particular,
\eqref{eq:globalratioexpansion} converges absolutely and normally
on compact subsets of $B$ and has the prescribed formal germ at $0$.
\end{theorem}
\begin{proof}
Choose $\boldsymbol b\in\R_{>0}^r$ with $\overline P_{\boldsymbol b}\Subset\Omega$. For $u\in\Omega^*$, $r_i=|u_i|$, put
\[
 c_i(r_i)=\frac{b_ir_i}{b_i+r_i},\qquad
 0<c_i(r_i)\le r_i,\quad c_i(r_i)<b_i.
\]
Coordinatewise radial interpolation from $r_i$ to $c_i(r_i)$ stays in $\Omega^*$ by completeness; interpolation from $c_i(r_i)$ to $b_i/2$ stays in $P_{\boldsymbol b}\cap(\C^*)^r$. Both homotopies preserve the arguments, and their final map is
$u\mapsto(b_i u_i/(2|u_i|))_i$, which restricts to the identity on the torus $T=\{|u_i|=b_i/2\}$. Thus
\begin{equation}\label{eq:omegapi1}
 \Omega^*\simeq T\cong(S^1)^r,\qquad
 \pi_1(\Omega^*)\cong\Z^r.
\end{equation}
Let $M_i\in\operatorname{GL}_n(\C)$, $n=\operatorname{rk}Q$, be the monodromies corresponding to the standard generators, so $[M_i,M_k]=0$. On their simultaneous generalized
eigenspaces, choose logarithms and use the finite nilpotent series $\log(I+N)=\sum_{k\ge1}(-1)^{k+1}N^k/k$.
This gives commuting $R_i\in\operatorname{Mat}_n(\C)$ with $e^{2\pi iR_i}=M_i,\; [R_i,R_k]=0,\; \operatorname{Spec}(-R_i)\subset \{z\in\C:0\le\Re z<1\}$. The connection
\[
 (Q_0,\nabla_0)=
 \left(\mathcal O_{\Omega^*}^{\oplus n},
 d-\sum_{i=1}^rR_i\frac{du_i}{u_i}\right)
\]
has the monodromies $M_i$, and its standard extension $(\mathcal O_\Omega^{\oplus n},\overline\nabla_0)$ has residues $-R_i$. By the local-system/flat-connection correspondence and \eqref{eq:omegapi1}, there is an isomorphism
of flat holomorphic bundles
\[
 (Q,\nabla)\xrightarrow{\sim}(Q_0,\nabla_0)
 \quad\text{on }\Omega^*.
\]
The functorial uniqueness of Deligne's canonical extension
\cite[Chapter~II, Proposition~5.4]{Deligne} extends this
isomorphism to the logarithmic extensions on $\Omega$. In particular, the corresponding local and global changes of logarithmic frame are holomorphic and invertible across $D$.

In the standard frame, the dual horizontal equation is
\[
 d\phi+\phi\sum_{i=1}^rR_i\frac{du_i}{u_i}=0.
\]
By Theorem~\ref{thm:formalconv}, the prescribed formal germ has a unique analytic realization near $0$. Equivalently, coefficientwise comparison in the finite Nilsson module gives
\begin{equation}\label{eq:globalphimodel}
 \widehat\Phi=c\prod_{i=1}^ru_i^{-R_i},\qquad
 \Phi^{\mathrm{an}}=c\prod_{i=1}^ru_i^{-R_i},
 \quad c\in(\C^n)^*,
\end{equation}
where the second expression is taken with the chosen logarithms. It analytically continues along every path in $\Omega^*$.

Under the global logarithmic-frame isomorphism, each $q_j$ becomes a vector $v_j\in\mathcal O(\Omega^*)^{\oplus n}$ meromorphic along $D$. Every $D_i$ is a nonempty, connected, smooth irreducible hypersurface: it is invariant under radial
contraction to $0$. Since all poles are supported on $D_1\cup\cdots\cup D_r$, the polar divisor of the finitely many components of the $v_j$ has a finite bound $\sum_i a_iD_i$ for some $a_i\in\Z_{\ge0}$, independent of $j$.
Normality of $\Omega$ (or Hartogs extension across the codimension-two intersections) gives
\[
 h_j(u)=u^{\boldsymbol a}v_j(u)
 \in\mathcal O(\Omega)^{\oplus n},\qquad
 h_j(u)=\sum_{\boldsymbol n\ge0}h_{j,\boldsymbol n}u^{\boldsymbol n}.
\]
By Lemma~\ref{lem:Reinhardt},
\begin{equation}\label{eq:globalcauchy}
 \sum_{\boldsymbol n\ge0}
 \sup_{u\in K}\|h_{j,\boldsymbol n}u^{\boldsymbol n}\|<\infty
 \qquad(K\Subset\Omega).
\end{equation}

The commuting matrices $R_i$ have a simultaneous generalized eigenspace decomposition. Hence, for a finite $\Lambda'\subset\C^r$ and some $b\ge0$,
\[
 \prod_{i=1}^ru_i^{-R_i}
 =\sum_{\lambda\in\Lambda'}
 u^\lambda\!
 \sum_{0\le\beta_i\le b}
 C_{\lambda,\beta}(\log u)^\beta,
 \qquad C_{\lambda,\beta}\in\operatorname{Mat}_n(\C).
\]
Combining this with \eqref{eq:globalphimodel} and $v_j=u^{-\boldsymbol a}h_j$ yields
\eqref{eq:globalratioexpansion}, with
$\Lambda_j\subseteq\Lambda'$ and coefficients obtained from $cC_{\lambda,\beta}h_{j,\boldsymbol n}$. The estimate \eqref{eq:globalcauchy} proves absolute normal convergence on every compact subset of a fixed branch $B$:
$u^\lambda$, $u^{-\boldsymbol a}$ and $(\log u)^\beta$ are bounded there. Near $0$, uniqueness of the formal horizontal solution identifies this expansion with the original formal evaluation $\widehat\Phi(q_j)$.
\end{proof}

\begin{proposition}\label{prop:algebraictests}
Let $Z$ be a smooth complex algebraic variety, $D=\bigcup_{i=1}^rD_i\subset Z$ a simple normal-crossing divisor, and $U=Z\setminus D$. Let $(F,\nabla)$ be a finite-rank algebraic integrable connection on $U$, regular
singular along $D$. Suppose that
\[
\begin{aligned}
 &\Omega\subset Z^{\mathrm{an}},\quad
 \Omega\subset\C^r,\quad 0\in\Omega,
 \quad\Omega\text{ complete Reinhardt},\\
 &D\cap\Omega=\{u_1\cdots u_r=0\}\cap\Omega,
\end{aligned}
\]
where the second inclusion denotes the specified analytic coordinate identification and $U\cap\Omega=\Omega^*$. Let $\overline F$ be Deligne's canonical logarithmic extension on $Z$ \cite[Chapter~II, Proposition~5.4]{Deligne}.
For $q_1,\ldots,q_m\in\Gamma(U,F)$ one has
\begin{equation}\label{eq:algebraictestmeromorphic}
 q_j|_{\Omega^*}\in
 \Gamma\bigl(\Omega,\overline F^{\mathrm{an}}(*D)\bigr),
 \qquad 1\le j\le m.
\end{equation}
Consequently, if $\widehat\Phi$ is a formal horizontal section of $F^\vee$ at $0$ satisfying the finite-Nilsson hypothesis of Theorem~\ref{thm:globalratio}, then
\begin{equation}\label{eq:algebraictestnormal}
 \sum_{\lambda,\beta}\sum_{\boldsymbol n\ge0}
 \sup_{u\in K}\bigl|
 u^\lambda(\log u)^\beta
 c_{j,\lambda,\beta,\boldsymbol n}
 u^{\boldsymbol n-\boldsymbol a}\bigr|<\infty,
 \qquad K\Subset B\subset\Omega^*,
\end{equation}
for the finite exponent--logarithm expansion of $\Phi^{\mathrm{an}}(q_j)$ on each simply connected branch $B$ (the fixed logarithms being used throughout).

In particular, suppose $U$ lies in a smooth algebraic Ward chart, $F$ is the algebraic pullback of $\cQ$ along $U$, and $\cO_U\otimes T\twoheadrightarrow F$ is the induced algebraic quotient map on all of $U$. Then each fixed
$w\in T$ defines an algebraic test section $q_w=[w]$.
\end{proposition}
\begin{proof}
By regular singularity along $D$, the algebraic connection extends to a locally free logarithmic connection $\overline F$ with the canonical residue interval
\cite[Chapter~II, Proposition~5.4]{Deligne}.
Fix $p\in D$ and an algebraic neighbourhood $W\subset Z$ on which $\overline F$ is free. The two bundles agree over $W\cap U$; hence, in this frame, each $q_j$ is a vector of
regular functions on $W\cap U$. Since $\cO_Z(W\cap U)\subset\C(Z),\; q_j\in \C(Z)^{\oplus\operatorname{rk}F}$, these are rational functions on $W$. Their denominators have no irreducible factors outside $D\cap W$, because $q_j$ is
regular on $W\cap U$. After shrinking $W$ at $p$, the normal-crossing parameters $u_i$ give, for some $\boldsymbol a\in\Z_{\ge0}^r$,
\[
 \Bigl(\prod_{i=1}^ru_i^{a_i}\Bigr)q_j
 \in\Gamma(W,\overline F),\qquad 1\le j\le m.
\]
Equivalently $q_j\in\overline F(*D)$; analytification and restriction to $\Omega$ yield
\eqref{eq:algebraictestmeromorphic}. The canonical
extension used here is the same canonical extension as in Theorem~\ref{thm:globalratio}; no arbitrary analytic gauge with uncontrolled growth is inserted.

The hypotheses of Theorem~\ref{thm:globalratio} now hold for $(F^{\mathrm{an}}|_{\Omega^*},\nabla)$ and the finite collection $\{q_j\}$. That theorem and Lemma~\ref{lem:Reinhardt} give
\eqref{eq:algebraictestnormal}. For the final assertion,
$1\otimes w$ is an algebraic section of $\cO_U\otimes T$, so its image under the stipulated algebraic quotient map, and then its algebraic pullback, belong to $\Gamma(U,F)$. No regularity conclusion for a nonalgebraic test section is used.
\end{proof}

For the following proposition, a rooted typed tree is written in the anchored convention of Lemma~\ref{lem:treesupport}: $r_v$ is the displacement of the operator at the internal vertex $v$, the output anchor is the anchor of its second
input, and $R=r_\rho,\; u_v=r_v/r_{\operatorname{par}(v)}\; (v\ne\rho)$. The ratios $u_v$ are independent; no additional, configuration-dependent change of these variables is implicit.

\begin{proposition}\label{prop:gradesnormal}
Fix homogeneous external vectors of weights $a_i$, a homogeneous terminal restricted-dual vector of weight $a_{\mathrm{out}}$, and an ordinary finite typed tree $\tau$. For each fixed choice of simple internal modules,
assume that its original grade-projected channel is a formal horizontal germ in the anchored ratios $u_v$ and that the corresponding scalar evaluation satisfies Theorem~\ref{thm:globalratio} (in particular its
meromorphic-test-section hypothesis), or Proposition~\ref{prop:algebraictests}, on the specified complete Reinhardt domain $\Omega$. Then, for every branch domain $B\subset\Omega^*$ and $K\Subset B$, the original intermediate-weight sum
satisfies
\begin{equation}\label{eq:originalgradesnormal}
 \sum_{\boldsymbol\mu}\sup_{u\in K}
 \left|C_{\boldsymbol\mu}
 R^{a_{\mathrm{out}}-\sum_i a_i}
 \prod_{v\ne\rho}
 u_v^{\mu_v-\sum_{i\in A_v}a_i}\right|<\infty.
\end{equation}
Here each $C_{\boldsymbol\mu}$ is the contraction of the finite-dimensional intermediate weight spaces before its absolute value is taken, and $R\ne0$ is fixed (uniformity also holds for $R$ in a compact set on a fixed branch).
For finite direct sums of ordinary internal modules,
\eqref{eq:originalgradesnormal} holds after the finite sum
over simple choices.
\end{proposition}
\begin{proof}
First fix simple internal modules. For every internal vertex $v$ let $v_{\mathrm L}$ and $v_{\mathrm R}$ be its two ordered children and $A_v$ its set of descendant leaves. Write $w\preceq v$ if $w$ lies on the root-to-$v$ path, and put
$\mu_\rho=a_{\mathrm{out}}$ and $\mu_i=a_i$ on the leaf edges. The homogeneity of ordinary intertwining operators (Lemma~\ref{lem:IOlogs}) gives the vertex exponent $e_v=\mu_v-\mu_{v_{\mathrm L}}-\mu_{v_{\mathrm R}}$. Writing
$r_v=R\prod_{\rho\ne w\preceq v}u_w$ and telescoping the exponents on each subtree gives
\[
 \sum_{v\in\mathcal V(\tau)}e_v
 =a_{\mathrm{out}}-\sum_i a_i,\qquad
 \sum_{w\succeq v}e_w
 =\mu_v-\sum_{i\in A_v}a_i\quad(v\ne\rho).
\]
Thus the contribution of the fixed internal weights is exactly
\[
C_{\boldsymbol\mu}
 R^{a_{\mathrm{out}}-\sum_i a_i}
 \prod_{v\ne\rho}
 u_v^{\mu_v-\sum_{i\in A_v}a_i}.
\]
All internal homogeneous spaces are finite-dimensional; $C_{\boldsymbol\mu}$ is therefore a finite contraction. This identity is coefficientwise and does not presuppose convergence (cf. Lemma~\ref{lem:treesupport}).

Write
\[
 \mu_v=h_v+n_v/N_v,\quad n_v\in\Z_{\ge0},\qquad
 \kappa_v=h_v-\sum_{i\in A_v}a_i,
\]
assigning coefficient zero when the indicated grade is empty. Choose $N\in\Z_{>0}$ divisible by every $N_v$, set $k_v=N/N_v$ and $u_v=\xi_v^N$. Set
\[
 \widetilde\Omega=\{\xi\in\C^{|\mathcal V(\tau)|-1}:
 (\xi_v^N)_v\in\Omega\},\qquad
 \widetilde\Omega^*=\widetilde\Omega\cap(\C^*)^{|\mathcal V(\tau)|-1}.
\]
Then
\[
 \widetilde\Omega\text{ complete Reinhardt},\quad
 0\in\widetilde\Omega,\qquad
 R^{-a_{\mathrm{out}}+\sum_i a_i}
 \xi^{-N\boldsymbol\kappa}\widehat\Phi_\tau
 =\sum_{\boldsymbol n\ge0}C_{\boldsymbol n}
 \prod_{v\ne\rho}\xi_v^{k_vn_v}.
\]
The exponents $\boldsymbol n\mapsto(k_vn_v)_{v\ne\rho}$ are injective. In particular the series on the right is an ordinary formal Taylor series without negative powers or logarithms, and each coefficient is exactly one original
weight-tuple contraction.

Let $f$ be the analytic scalar evaluation supplied by the hypothesis, and initially on a branch put
\[
 H(\xi)=R^{-a_{\mathrm{out}}+\sum_i a_i}
 \xi^{-N\boldsymbol\kappa}
 f((\xi_v^N)_{v\ne\rho}).
\]
Near $0$ the formal germ of $H$ is the displayed Taylor series. By Theorem~\ref{thm:globalratio}, $f$ is an analytic continuation of this germ, with finitely many exponent classes, bounded logarithmic degrees and finite coordinate
pole orders. The local Taylor germ has trivial coordinate monodromy. Uniqueness of analytic continuation makes the monodromy of $H$ trivial around the generators of $\pi_1(\widetilde\Omega^*)\cong\Z^{|\mathcal V(\tau)|-1}$
(cf. \eqref{eq:omegapi1}). Hence $H$ is single-valued on $\widetilde\Omega^*$.

The finite-power logarithmic expansion of Theorem~\ref{thm:globalratio}, together with this trivial monodromy, shows that $H$ has only integral Laurent powers and no logarithmic terms along the coordinate divisors:
indeed, the generalized-power--logarithm coefficients with nontrivial coordinate monodromy must vanish. Thus $H$ extends meromorphically to $\widetilde\Omega$ with poles supported on the coordinate divisors. Its principal parts
vanish in a neighbourhood of $0$, by the Taylor germ. Since every $\widetilde\Omega\cap\{\xi_v=0\}$ is connected, analytic continuation along that divisor makes all these principal parts vanish globally. Hartogs extension across
their intersections then gives
\[
 H\in\mathcal O(\widetilde\Omega),\qquad
 H(\xi)=\sum_{\boldsymbol n\ge0}
 C_{\boldsymbol n}\prod_{v\ne\rho}\xi_v^{k_vn_v}.
\]
Lemma~\ref{lem:Reinhardt} yields
\begin{equation}\label{eq:gradescauchy}
 \sum_{\boldsymbol n\ge0}\sup_{\xi\in\widetilde K}
 \left|C_{\boldsymbol n}
 \prod_{v\ne\rho}\xi_v^{k_vn_v}\right|<\infty
 \quad(\widetilde K\Subset\widetilde\Omega).
\end{equation}
On a compact $K\Subset B$, choose the coordinatewise $N$th roots consistently on $B$. Their values lie in a compact $\widetilde K\Subset\widetilde\Omega^*$; the fixed factors $R^{a_{\mathrm{out}}-\sum_i a_i}$ and
$\xi^{N\boldsymbol\kappa}$ are bounded there. Combining \eqref{eq:gradescauchy} with the injectivity of $\boldsymbol n\mapsto(k_vn_v)_v$ gives
\eqref{eq:originalgradesnormal} without interchanging any
infinite sums before absolute convergence is proved. Finally, finitely many simple internal choices contribute a finite sum of bounds of the same form.
\end{proof}

\section{Actual product and iterate Ward identities}\label{sec:threeinputward}
\subsection{One-edge contraction on the actual cover}

Retain the objects $A\in\cC_g$, $B\in\cC_h$, $C\in\cC_k$ and $D\in\cC_{ghk}$ from Section~\ref{sec:setup}, together with simple ordinary $M\in\cC_{hk}$ and $N\in\cC_{gh}$. Put $T=D'\otimes A\otimes B\otimes C$. For homogeneous
$d'\in D'$, $a\in A$, $b\in B$ and $c\in C$, write
\[
 (\lambda_D,\lambda_A,\lambda_B,\lambda_C)
 =(\wt(d'),\wt(a),\wt(b),\wt(c)),\qquad
 \Delta=\lambda_D-\lambda_A-\lambda_B-\lambda_C.
\]
For $U\in\{M,N\}$, with lowest weight $h_U$, define $\Lambda_U=\{\gamma\in N_U^{-1}\N:U_{h_U+\gamma}\ne0\},\; N_U=\operatorname{ord}(\operatorname{tw}(U))$. For $\gamma\in\Lambda_U$, choose a basis $\{u_{\gamma,j}\}_{j=1}^{r_{U,\gamma}}$ of
$U_{h_U+\gamma}$ and its dual basis $\{u_{\gamma,j}^{*}\}_{j=1}^{r_{U,\gamma}}\subset U'$, and put
\begin{equation}\label{eq:gradecoevaluation61}
 \mathscr E_{\gamma}^{U}
 =\sum_{j=1}^{r_{U,\gamma}}
 u_{\gamma,j}\otimes u_{\gamma,j}^{*}
 \in U_{h_U+\gamma}\otimes U_{h_U+\gamma}^{*},\qquad
 \Gamma_U(\rho)=\sum_{\gamma\in\Lambda_U}
 \rho^{\gamma}\mathscr E_{\gamma}^{U}.
\end{equation}
The tensor $\mathscr E_{\gamma}^{U}$ corresponds to $\id_{U_{h_U+\gamma}}$ under
$U_{h_U+\gamma}\otimes U_{h_U+\gamma}^{*}
\simeq\End_{\C}(U_{h_U+\gamma})$; hence it is
basis independent. The sum $\Gamma_U$ is \emph{formal}; it is not asserted to lie in $U\otimes U'$. If $e\in\N_{>0}$ is divisible by $N_M$ and $N_N$, then $e\Lambda_U\subset\N$, so $\Gamma_U(t^e)$ is a
coefficientwise defined power series in $t$. This substitution does not require $e h_U\in\Z$.

For $\mathcal Y_1\in\cI\binom D{A\ M}$ and $\mathcal Y_2\in\cI\binom M{B\ C}$, put
\begin{align*}
 \ell_P(d',a,m)&=\left.\langle d',\mathcal Y_1(a,x)m\rangle\right|_{x=1},
 &r_P(m',b,c)&=\left.\langle m',\mathcal Y_2(b,x)c\rangle\right|_{x=1}.
\end{align*}
Here evaluation means $x^\lambda\mapsto1$ with $\log1=0$ \emph{after pairing with the restricted dual}; for homogeneous inputs each paired series has one power $x^\lambda$ and no logarithms by Lemma~\ref{lem:IOlogs}.
For algebraic inputs only finitely many homogeneous terms occur. Both evaluations are therefore finite sums after pairing with the restricted dual. Define the formal scalar contraction
\begin{equation}\label{eq:productcoeval61}
 \Psi_P(\rho)(d'\otimes a\otimes b\otimes c)
 =\sum_{\gamma\in\Lambda_M}\rho^\gamma
   \sum_{j=1}^{r_{M,\gamma}}
     \ell_P(d',a,m_{\gamma,j})
     r_P(m_{\gamma,j}^*,b,c).
\end{equation}
For $\mathcal Y_3\in\cI\binom D{N\ C}$ and $\mathcal Y_4\in\cI\binom N{A\ B}$, define
\[
 \ell_I(d',n,c)=\left.\langle d',\mathcal Y_3(n,x)c\rangle\right|_{x=1},
 \qquad
 r_I(n',a,b)=\left.\langle n',\mathcal Y_4(a,x)b\rangle\right|_{x=1},
\]
and
\begin{equation}\label{eq:iteratecoeval61}
 \Psi_I(\rho)(d'\otimes a\otimes b\otimes c)
 =\sum_{\gamma\in\Lambda_N}\rho^\gamma
   \sum_{j=1}^{r_{N,\gamma}}
     \ell_I(d',n_{\gamma,j},c)
     r_I(n_{\gamma,j}^*,a,b).
\end{equation}
Every inner sum in \eqref{eq:productcoeval61} and
\eqref{eq:iteratecoeval61} is finite and basis independent.
These definitions use only the single-operator properties of
\cite[Definition~2.8]{DH}; they do not invoke convergence of
products of two twisted intertwining operators.

\begin{lemma}
\label{lem:actualonecontraction}
With $R>0$ and $0<\rho<1$, $z_1=R$, $z_2=R\rho$ for the product, and $z_2=R$, $z_1=R(1+\rho)$ for the iterate, the formal internal-grade contractions satisfy
\begin{align*}
 \Phi_P(R,R\rho)(d'\otimes a\otimes b\otimes c)
 &=R^\Delta\rho^{h_M-\lambda_B-\lambda_C}
   \Psi_P(\rho)(d'\otimes a\otimes b\otimes c),
 \\
 \Phi_I(R(1+\rho),R)(d'\otimes a\otimes b\otimes c)
 &=R^\Delta\rho^{h_N-\lambda_A-\lambda_B}
   \Psi_I(\rho)(d'\otimes a\otimes b\otimes c).
\end{align*}
The equalities hold coefficientwise in the relative grades; they extend to branches with $\ell(R\rho)=\ell(R)+\ell(\rho)$ and, respectively, $\ell(z_1-z_2)=\ell(R)+\ell(\rho)$. No analytic convergence of either left side is assumed.
\end{lemma}
\begin{proof}
For $\gamma\in\Lambda_M$ and $m\in M_{h_M+\gamma}$, Lemma~\ref{lem:IOlogs} gives
\begin{align*}
 \langle d',\mathcal Y_1(a,R)m\rangle
 &=R^{\lambda_D-\lambda_A-h_M-\gamma}\ell_P(d',a,m),\\
 \langle m',\mathcal Y_2(b,R\rho)c\rangle
 &=(R\rho)^{h_M+\gamma-\lambda_B-\lambda_C}
    r_P(m',b,c)
 \quad(m'\in M_{h_M+\gamma}^{*}).
\end{align*}
Consequently,
\[
 R^{\lambda_D-\lambda_A-h_M-\gamma}
 (R\rho)^{h_M+\gamma-\lambda_B-\lambda_C}
 =R^\Delta\rho^{h_M-\lambda_B-\lambda_C}\rho^\gamma.
\]
Inserting the finite coevaluation $\mathscr E_{\gamma}^{M}$ yields precisely the $\gamma$-th summand of \eqref{eq:productcoeval61}. The same calculation with $N$, $\mathcal Y_3$ and $\mathcal Y_4$ yields
\[
 (R\rho)^{h_N+\gamma-\lambda_A-\lambda_B}
 R^{\lambda_D-h_N-\gamma-\lambda_C}
 =R^\Delta\rho^{h_N-\lambda_A-\lambda_B}\rho^\gamma.
\]
The prescribed logarithms justify the displayed power identities; finite-dimensional duality proves basis independence. Both equalities are coefficientwise, before summing over $\gamma$.
\end{proof}

For the algebraic tensor $T$, define the weight-diagonal operators
\[
 \begin{aligned}
 D_P(\rho)&=\id_{D'}\otimes\id_A\otimes
 \rho^{-L_B(0)}\otimes\rho^{-L_C(0)},\\
 D_I(\rho)&=\id_{D'}\otimes\rho^{-L_A(0)}\otimes
 \rho^{-L_B(0)}\otimes\id_C,
 \end{aligned}
 \qquad
 \rho^{\alpha}=e^{\alpha\ell(\rho)}.
\]
Write $\Phi_P(\rho):=\Phi_P(1,\rho)$ and $\Phi_I(\rho):=\Phi_I(1+\rho,1)$. Then the preceding lemma becomes, as an equality of formal scalar functionals on $T$,
\begin{equation}\label{eq:physicalscaling61}
 \Phi_P(\rho)=\rho^{h_M}\Psi_P(\rho)\circ D_P(\rho),
 \qquad
 \Phi_I(\rho)=\rho^{h_N}\Psi_I(\rho)\circ D_I(\rho).
\end{equation}
All operators in $D_P,D_I$ act on finite sums of weight vectors. In particular, if a current changes an external weight, its action must be composed with the full operators $D_P,D_I$ in
\eqref{eq:physicalscaling61}, not with their values on a fixed
homogeneous tensor.

\begin{lemma}
\label{lem:canonicalnodecurrent}
Let $\star\in\{P,I\}$ and let $s$ be an $H$-invariant canonical seed current on the actual generic framed cover of the corresponding three-input chart. After a finite extension of the collision DVR, the following hold.
\begin{enumerate}
\item The restriction of $s$ to the normalization of each component
of the special fiber has a coefficientwise expansion in meromorphic three-point currents on its actual based vertex cover. There is $m=m(s)\geq0$ such that $t^m s$ extends to the coherent current
bundle with its prescribed external poles.
\item At a source node of inertia $N$, with a compatible choice of
root coordinates, one has
\begin{equation}\label{eq:actualbalancednode62}
  uv=t^a,\qquad x=u^N,\qquad y=v^N,\qquad
  xy=\rho=t^{aN},\qquad a\in\N_{>0},
\end{equation}
and the two node-current frames differ by the constant operator $\Sigma=e^{L(1)}(-1)^{L(0)+1}$, including all the $L(1)$-terms.
\item For an invariant paired one-form $\eta$, the orbit-normalized
node action is
\begin{equation}\label{eq:normalizednode62}
  \frac{1}{|H_v|}\sum_{p\in H_v\widetilde Q}\Res_p\eta
  =\frac1N\Res_{\widetilde Q}\eta,
\end{equation}
where $H_v$ is the deck group of the distinguished connected vertex cover, $\widetilde Q$ is its selected node lift and $N=|\operatorname{Stab}_{H_v}(\widetilde Q)|$.
\end{enumerate}
The statements concern each fixed generic current $s$; no pole bound independent of $s$ is asserted.
\end{lemma}
\begin{proof}
For the product put $Z=z/R$ and $(z_1,z_2)=(R,R\rho),\; x=Z,\; \xi=Z/\rho,\; y=\xi^{-1}$. For the iterate put $(z_1,z_2)=(R(1+\rho),R)$ and $x=Z-1,\; \xi=(Z-1)/\rho,\; y=\xi^{-1}$.
In either case $xy=\rho$. Each normalized coarse component has three special points, with coordinates $0,1,\infty$. Work over the completed DVR of a boundary point of the actual Hurwitz cover, rather than
of the unframed configuration base. By
\cite[Corollary~3.0.5 and Theorem~4.3.2]{ACV}, after finite DVR
extension its generic balanced twisted $H$-cover extends across the stable marked boundary. The equivalence with admissible $H$-covers
\cite[Theorem~4.3.2]{ACV} gives the local node charts; the total
source is projective and nodal by \cite[Lemma~2.2.1]{ACV}. The chosen generic lift over the auxiliary unramified basepoint extends, after a further finite extension, by the valuative criterion
for a finite cover; compatible inertia-root phases are fixed after adjoining the finitely many required roots. Thus the extended cover is compared with the original \emph{framed} generic cover, rather
than with a newly selected special-fiber cover.

Let $O=L[[t]]$, $K=L((t))$ after this extension. The completed local rings at a coarse node and a source node initially have forms
\[
 O[[X,Y]]/(XY-\epsilon t^{aN}),\qquad
 O[[U,V]]/(UV-\epsilon' t^a),\qquad
 X=U^N,\quad Y=V^N,
\]
with $\epsilon,\epsilon'\in O^\times$ satisfying $\epsilon=(\epsilon')^N$ after a choice of compatible local charts. The coarse rational functions $x,y$ above define the same two
oriented branches and the same smoothing divisor as $X,Y$. Comparing their divisors in the normal completed coarse local ring gives $x/X,y/Y\in O[[X,Y]]/(XY-\epsilon t^{aN})^\times$.
After finite residue-field extension these units have unique $N$th roots with specified residues (Hensel's lemma, $N\in O^\times$). Absorb them into $U,V$, and absorb the remaining smoothing unit into
the opposite root and, if needed, the base uniformizer. This gives
\eqref{eq:actualbalancednode62}. It also fixes the residual
$\mu_N$-ambiguity by the inherited marked roots. The inertia action is $u\mapsto\zeta u$, $v\mapsto\zeta^{-1}v$; in particular neither branch is independently re-framed. The ramification exponent is fixed more precisely as follows.  Write $n=\operatorname{ord}(h)$ for the inertia order of the cluster-boundary label $h=g_Bg_C$ (product) or $h=g_Ag_B$ (iterate), and $N_X=\operatorname{ord}(g_X)$ for an external marked leg $X$.  Let $m_{\star}$ be the order of the finite permutation monodromy of the selected full framed Hurwitz sheet around the given boundary.  Enlarge the already allowed DVR base extension, if necessary, to $\rho=t^{e_{\star}}$ with
\begin{equation}\label{eq:E0sufficientbasechange}
 n\mid e_{\star},\qquad m_{\star}\mid e_{\star},\qquad
 \begin{cases}N_B,N_C\mid e_P,&\star=P,\\
 N_A,N_B\mid e_I,&\star=I.
 \end{cases}
\end{equation}
Thus $a=e_{\star}/n$ is a positive integer, and the scaled external root parameters $\rho^{1/N_X}=t^{e_{\star}/N_X}$ are single-valued on the same root disk.  For example one may take the least common multiple of the displayed integers, or a positive multiple needed for other finite torsors.  Enlarging $e_{\star}$ only pulls back the given rooted family: on $t>0$, $\log\rho=e_{\star}\log t$, and no marked lift, root phase or original intertwining-operator coefficient is replaced.  Lowest weights remain in their specified generalized powers and need not become integral after this base change.

On $D(u)\cap D(v)$ one has $y=\rho/x$. Its centered Taylor coordinate change is
\[
 \tau\longmapsto
 \frac{\rho}{x+\tau}-\frac{\rho}{x}
 =\frac{-\rho x^{-2}\tau}{1+\tau/x}.
\]
For $b\in V_d$, the inverse state transformation is $(-\rho x^{-2})^{L(0)}e^{-x^{-1}L(1)}b$. Since $L(1)^{d+1}b=0$, it is a finite sum. Its $L(1)^jb/j!$-term, including the differential, satisfies
\begin{equation}\label{eq:balancednodeframe}
 x^{d-1}(-\rho x^{-2})^{d-j}(-x^{-1})^j\,dx
 =(-1)^{d+1}y^{d-j-1}\,dy.
\end{equation}
Indeed, $dx=-\rho y^{-2}dy$, $x=\rho/y$ and all powers of $\rho$ cancel. Consequently the balanced state-current transition is
\[
 \Sigma=e^{L(1)}(-1)^{L(0)+1},\qquad
 (-1)^{L(0)}e^{L(1)}=e^{-L(1)}(-1)^{L(0)},\qquad
 \Sigma^2=\id_E.
\]
Here $L(0)$ has integer eigenvalues on $V$, $E$ is stable under $L(1)$, and the state automorphisms commute with $L(0),L(1)$. The coordinate change is fractional-linear, so its Schwarzian
vanishes. At ramified external markings the observation gauges, including vacuum contributions, are composed with the full coordinate-change cocycle; the calculation above is for the underlying state-current frame.

Put $B=O[[u,v]]/(uv-t^a)$. It is a normal two-dimensional Cohen--Macaulay hypersurface (Serre's $R_1+S_2$), and $D(u)\cup D(v)$ is the complement of its closed node. The two balanced trivializations, with transition $\Sigma$, yield a free
$B$-module $B\otimes_L E$ on this punctured spectrum (change the second frame by $\Sigma$). The coherent current bundle is reflexive there; hence the depth-two extension equivalence
\cite[Tag~0EBJ]{Stacks} identifies its reflexive extension across
the codimension-two node with $B\otimes_L E$. This conclusion uses reflexivity and the actual frame cocycle, not just the normality of $B$.

The normalized special-fiber components carry the \emph{induced} based covers. Write $h_{AB}=g_Ag_B$ and $h_{BC}=g_Bg_C$. In the input order of \eqref{eq:orderedmeridians}, their complete three-point tuples are
\[
 \begin{array}{c|cc}
 &\text{outside component}&\text{bubble component}\\ \hline
 P&(g_A,h_{BC},(g_Ah_{BC})^{-1})
   &(g_B,g_C,h_{BC}^{-1})\\
 I&(h_{AB},g_C,(h_{AB}g_C)^{-1})
   &(g_A,g_B,h_{AB}^{-1})
 \end{array}
\]
To obtain these tuples, transport the reference stems from the positive real collision chart and retain the boundary stem of the consecutive cluster. Its meridian is $\gamma_B\gamma_C$ for $P$ and $\gamma_A\gamma_B$ for $I$. On the bubble the outgoing meridian is the inverse boundary, whereas on the outside component it is the incoming boundary. Thus the two node labels are exactly $h_{BC},h_{BC}^{-1}$ or $h_{AB},h_{AB}^{-1}$ in these based frames, and every displayed triple has product one. No exchange of stems is made during either collision.
The based finite-\'etale classification identifies these with the binary Ward covers having the same monodromy and selected lifts: equality of their based monodromy homomorphisms makes the comparison independent of the path used to lift a point. The selected component basepoints are obtained from the inherited boundary stems. If another distinguished component is used, its labels and states are simultaneously changed by $(g,M)\mapsto(kgk^{-1},T_kM)$; the same deck change is applied on both sides of a node. The comparison extends to the marked Kummer disks using the specified roots, and to the proper curves by finite-cover
algebraicity. At finite external markings the centered physical coordinates are respectively $R$ and $R\rho$ times the component coordinates; at infinity the factor is $R^{-1}$. Their root powers
are fixed by the retained logarithms and framing. A change of the internal root by a deck rotation acts on $M$ by $S=e^{-2\pi i(L_M(0)-h_M)}$ and on $M'$ by $(S^{-1})^t$; it fixes $\mathscr E^M_\gamma$ in \eqref{eq:gradecoevaluation61}.
No complex ground-weight factor is discarded; its separate scaling is given by Lemma~\ref{lem:actualonecontraction}.

Fix $s$ as in the statement. Let $D$ be an $H$-invariant horizontal pole divisor, disjoint from the source nodes, bounding all generic poles of $s$, and let $F(D)$ be the corresponding coherent extension
of the finite current bundle on the proper nodal source $p:C_O\to
\operatorname{Spec}O$. Properness and flat base change give
\[
 (p_*F(D))\otimes_OK\simeq H^0(C_K,F(D)_K),
 \qquad p_*F(D)\ \text{finite over }O.
\]
Thus $t^m s$ extends for some $m=m(s)\geq0$. Since $|H|\in O^\times$, averaging preserves $H$-invariance. At each source node, the balanced coefficients of $s$ belong to $t^{-m}B$.
This bound is uniform over the finitely many nodes \emph{for this} $s$.

To compare the coefficients with the fixed three-point component currents, set $O_J=O/(t^{J+1})$, let $U_v=\mathbb P^1\setminus\{0,1,\infty\}$, and denote the selected connected component of the special vertex cover above
$U_v$ by $C_v^\circ$. Over $U_v\times\operatorname{Spec}O_J$ the extra moving marked section is absent. For example, the product outside equation and its bubble analogue have factors $Z-\rho=Z(1-\rho/Z),\; 1-\rho\xi$,
which are units because $\rho$ is nilpotent in $O_J$; the iterate factors are $(Z-1)(1-\rho/(Z-1))$ and $1+\rho\xi$. The actual cover is therefore finite \'etale on this open. Restriction
along the nilpotent thickening is an equivalence on finite-\'etale covers \cite[Tag~0BQB]{Stacks}. The selected special-fiber comparison consequently has a unique compatible lift
\[
 \iota_{v,J}:C_v^\circ\times\operatorname{Spec}O_J
 \xrightarrow{\sim}(C_{O_J}|_{U_v})^{\mathrm{selected}},
\]
for every $J\geq0$, where the full $H$-cover is obtained by induction from the connected stabilizer $H_v$. Uniqueness yields compatibility in $J$ and deck equivariance. At the deleted markings
and on the two punctured node branches use the tame Kummer presentations. With $N_i\in O^\times$, an $N_i$th root of a unit with prescribed reduction lifts uniquely, and an $N_i$th root of
unity reducing to $1$ is $1$. These facts extend the given comparison over the ramified disks, with their chosen root phases. On $D(u)$ and $D(v)$ the relations $v=t^a/u$ and $u=t^a/v$
identify the comparisons with the two restrictions of the node frame. Thus all the comparisons agree on overlaps, including the state-current cocycle and the full module observation gauges.

It follows that at each external marked orbit, modulo $t^{J+1}$, $R_{i,\mathrm{rel}}(s) =R_{i,v}(\iota_{v,J}^{*}s)$, where $R_{i,v}$ is precisely the orbit-normalized selected-lift
residue of \eqref{eq:canonicalmarkedres}. The physical-to-component weight dilations are applied subsequently, as in Lemma~\ref{lem:physicalcontraction}; they are not absorbed into this identity.

Finally, an element of $t^{-m}B$ has an expansion $\sum_{r,s,\ell\geq0}c_{rs\ell}t^{\ell-m}u^rv^s$, subject to $uv=t^a$. On the $u$ branch its $(r,s,\ell)$-term equals $t^{as+\ell-m}u^{r-s}$;
for a fixed coefficient $t^j$, negative $u$-powers require $s\leq\lfloor(j+m)/a\rfloor$. Interchanging $r,s$ gives the same result on the $v$ branch. The finite coordinate gauges add a
fixed pole shift, so every coefficient has finite pole order at the node. \emph{Global meromorphicity} of this coefficient on the fixed proper vertex cover follows in addition from restricting the
\emph{global} section $t^ms\in H^0(C_O,F(D))$ and applying the
compatible identifications $\iota_{v,J}$; finite Laurent pole order alone would not imply this extension. Its other poles are confined to the external marked orbits. Consequently every coefficient is a
bona fide three-point Ward current on the corresponding actual vertex cover.

For the orbit normalization let $\widetilde Q$ be a node lift. An $H_v$-invariant paired one-form has equal residues at the $|H_v|/N$ points of $H_v\widetilde Q$. Therefore
\eqref{eq:normalizednode62} holds, and for an admissible mode
$n\in\frac rN+\Z$ the selected-root normalization reads
\[
 \frac1N\Res_{u=0}
 u^{Nn}Y_M^g(b,u^N)\,Nu^{N-1}du=b_n,
 \qquad b\in V^{(r)}.
\]
Induction from $H_v$ to $H$ shows that the global $H$-average of the full node orbit equals this distinguished-component average; source nodes over the same coarse node are not counted again.
The node module carries only the indicated cyclic inertia twist; no additional $H$-action on it is assumed.
\end{proof}

\begin{lemma}\label{lem:E0actualsourcecomparison}
With the outset data of Definition~\ref{def:commonoutsetframing}, the component-based comparisons in Lemma~\ref{lem:canonicalnodecurrent} represent the original Du--Huang three-point source germs used in \eqref{eq:productcoeval61} and \eqref{eq:iteratecoeval61}.  For $R=1$ the physical and component coordinates are
\[
\begin{array}{c|c|c|c}
 &\text{outside }Z&\text{node }x&\text{bubble }\xi\\\hline
 P&Z&Z&Z/\rho\\
 I&Z&Z-1&(Z-1)/\rho
\end{array},\qquad y=\xi^{-1},\quad xy=\rho.
\]
In particular, for the iterate the \emph{outside} original operator $\mathcal Y_3(n,1)c$ is inserted at $Z=1$; the centered $x=Z-1$ describes only its node neighborhood, while $\mathcal Y_4(a,1)b$ is evaluated at $\xi=1$.  The product has original operators $\mathcal Y_1(a,1)m$ and $\mathcal Y_2(b,1)c$ at the corresponding outside and bubble insertion points.
\end{lemma}
\begin{proof}
We construct the vertex functionals \emph{from the four given
intertwining operators}; in particular we do not identify arbitrary
solutions of a three-point Ward system.  Apply
Lemma~\ref{lem:singleinsertion} to each $\mathcal Y_i$ with its
unchanged $z=1$, $\log 1=0$ branch, obtaining the four dual current
sections $\lambda_i^{\mathrm{DH}}$ on their specified rooted
three-point covers.  The purpose of the following comparison is to
pull these very sections, and not newly chosen three-point blocks, to
the normalized components of the actual degeneration.

Fix $0<\rho<1/16$.  The two smooth physical configurations and
based points are
\[
 (1,\rho,0),\quad b_P=3+\rho;
 \qquad (1+\rho,1,0),\quad b_I=4+\rho.
\]
Continue the single outset lift, lower-comb stems, marked roots
and logarithms of Definition~\ref{def:commonoutsetframing} to these
points before taking either collision.  If $h_{BC}=g_Bg_C$ and
$h_{AB}=g_Ag_B$, the retained cluster stems have ordered words
$\gamma_B\gamma_C$ and $\gamma_A\gamma_B$.  The resulting vertex
data are as follows; the pairs in the fourth column give the
ordered finite input markings in the normalized coordinate.
\[
\begin{array}{c|c|c|c|c}
 &\text{original IO}&\text{component coordinate}
 &\text{normalized inputs}&\text{ordered inertia triple}\\\hline
 P_o&\mathcal Y_1&Z&(A:1,M:0)
 &(g_A,h_{BC},(g_Ah_{BC})^{-1})\\
 P_b&\mathcal Y_2&\xi=Z/\rho&(B:1,C:0)
 &(g_B,g_C,h_{BC}^{-1})\\
 I_o&\mathcal Y_3&Z&(N:1,C:0)
 &(h_{AB},g_C,(h_{AB}g_C)^{-1})\\
 I_b&\mathcal Y_4&\xi=(Z-1)/\rho&(A:1,B:0)
 &(g_A,g_B,h_{AB}^{-1})
\end{array}
\]
In particular, the outside iterate insertion is at $Z=1$, not at
$x=Z-1=0$; $x$ is only its centered node parameter.  These data
follow by restricting the given based Hurwitz cover along the
common boundary stem, with no permutation of meridians.

The following path comparison includes the one-$V$ observation
variable.  On a standard three-point curve with marks $1,0,\infty$
choose a right basepoint $B_0>1$ and $0<\epsilon<1/4$.  The path to
$1$ is the real interval from $B_0$ to $1$; the lower path to $0$
consists of the real interval to $1+\epsilon$, the arc
\[
 w(\theta)=1+\epsilon e^{-i\theta},\qquad0\le\theta\le\pi,
 \tag{*}
\]
and the real interval from $1-\epsilon$ to $0$.  The output path
runs to $+\infty$ on the right.  Moving $B_0$ in $(1,\infty)$ is a
homotopy relative to the marked endpoints, since all the resulting
initial intervals avoid $0,1$.

We verify that the lower arc $(*)$ gives the fixed preferred
continuation in Definition~\ref{def:completeintertwiner}, rather than
merely an admissible geometric stem.  Denote the additional
$V$-insertion position by $w$ and the intertwining-operator position
by $z$.  Put $\alpha=3\pi/2$ and
\[
 (w_\theta,z_\theta)=
 \bigl(e^{i\alpha}(1+\epsilon e^{-i\theta}),
       e^{i\alpha}\bigr),\qquad 0\le\theta\le\pi.
\]
Throughout this path, $|w-z|=\epsilon<|z|$ and
$|\arg w-\arg z|<\pi/2$, so it remains in the original preferred
$\mathcal D_3$ of Definition~\ref{def:completeintertwiner}.  Its
initial point belongs also to $\mathcal D_1$, and its final point to
$\mathcal D_2$: at the latter point the \emph{continuous} difference
of preferred arguments is $-\pi$, not $+\pi$.  Consequently the three
single-insertion expansions $F_\infty,F_z,F_0$ are continuations of
one and the same original $\lambda_i^{\mathrm{DH}}$ along
this arc, by the defining one-$V$ duality of $\mathcal Y_i$.
To move this comparison to $z=1$ use the explicit homotopy
\[
 H(t,\theta)=\bigl(e^{i(1-t)\alpha}
 (1+\epsilon e^{-i\theta}),\ e^{i(1-t)\alpha}\bigr),
 \quad 0\le t\le1,
\]
which stays in $wz(w-z)\ne0$.  Its uninterrupted logarithms are
\[
 L_z=i(1-t)\alpha,\quad
 L_w=\operatorname{Log}_{\mathrm{right}}(1+\epsilon e^{-i\theta})
      +i(1-t)\alpha,\quad
 L_{w-z}=\log\epsilon+i(1-t)\alpha-i\theta.
\]
Thus at $t=1$ it is precisely the lower arc $(*)$, with
$L_z=0$ and $L_{w-z}=\log\epsilon-i\theta$.  The rotation is used
to compare the \emph{observation-current} germs; it neither
redefines $\mathcal Y_i$ nor changes its final $z=1$, $\log1=0$
coefficient.  In particular the sign of the bypass is fixed by the
Du--Huang preferred regions, not chosen to force a Ward identity.

To compare these standard paths with the inherited ones, take the
common boundary point at $Z=\sqrt\rho$ in the product and
$Z=1+\sqrt\rho$ in the iterate.  It lies strictly between the two
scales in either case.  For the product, the lower-comb stem from
$b_P$ passes below $A=1$ and approaches $Z=\sqrt\rho$; in the
outside normalization this becomes the lower path from the right to
the node $0$.  On the bubble its image is
$\xi=\rho^{-1/2}>1$, so sliding this positive right basepoint in
$(1,\infty)$ gives the standard paths to $B=1$ and $C=0$.
For the iterate the stem from $b_I$ reaches $Z=1+\sqrt\rho$
from the right of the whole $(A,B)$ cluster: it becomes the direct
right-hand path to the outside node at $Z=1$, while its bubble image
$\xi=\rho^{-1/2}>1$ gives the direct path to $A=1$ and the lower
path to $B=0$.  In all four cases the homotopies are compositions of
these positive real basepoint slides and the lower arc $(*)$;
there is no added loop about $0$ or $1$.  Because the \emph{full}
cluster-boundary stem, rather than just its meridian word, is retained,
the comparisons preserve its inherited unramified based lift and its
noncommutative deck label.  This proves the four required
endpoint-fixed \emph{observation-path} comparisons.  The operator
insertion itself stays at normalized $z=1$ throughout, with
$\log z=0$.

The special roots and node comparison give the same
\emph{logged} transport.  The actual marking roots on either bubble
are obtained from the fixed physical roots by
$Z=\rho\xi$ or $Z-1=\rho\xi$; the centered roots on its two
external legs therefore acquire precisely the factors
$\exp(\ell(\rho)/N_X)$.  These are continuous from the original
positive logarithm and are not reset.  The two outside coordinates
are the global $Z$, as in the table.  At a node put
$x=Z$ in the product and $x=Z-1$ in the iterate, and let
$y=\xi^{-1}$.  On the same retained annulus point
$x=y=\sqrt\rho$ we have
\[
 L_x=L_y=\tfrac12\ell(\rho),\qquad
 L_x+L_y=\ell(\rho),\qquad xy=\rho.
\]
Write $\zeta=\exp(2\pi i/n)$.  For the jointly induced three-point reference covers of
Definition~\ref{def:commonoutsetframing}, take the two node roots
$u_{\mathrm{DH}}=\exp(L_x/n)$ and
$v_{\mathrm{DH}}=\exp(L_y/n)$ \emph{together} on this same
annulus: starting with the actual root pair, apply the unique
balanced $\mu_n$ rotation making the outside root positive; the
bubble root then has the prescribed positive value because its product
with the outside root is $\exp(\ell(\rho)/n)$.
This is the simultaneous auxiliary realization of the already fixed
DH preferred continuation, not a new analytic continuation or a
separate choice of the two roots.  The DH intertwining-operator
definition fixes the preferred analytic branch but does not itself
specify an independent root of the \emph{internal node} on each of
two distinct three-point covers.  In the actual balanced chart
$x=u^n$, $y=v^n$, $uv=t^a$ and $\rho=t^{an}$ after the prescribed
positive DVR base change.  Hence there are unique residue classes
$r,s\pmod n$ such that
\[
 u=\zeta^r u_{\mathrm{DH}},\quad
 v=\zeta^s v_{\mathrm{DH}},\qquad
 \zeta^{r+s}=\frac{uv}{u_{\mathrm{DH}}v_{\mathrm{DH}}}=1.
\]
Thus $s=-r\pmod n$ is a \emph{consequence} of the original
balanced relation, not a choice of a new node root.  The associated
relative-weight frame operators are
$S^r=\exp(-2\pi i r(L_M(0)-h_M))$ and $(S^{-r})^t$
(on the iterate node replace $M$ by $N$).  The same retained boundary
stem supplies the noncommutative deck comparison on both legs, with
the appropriate strict module pullbacks $T_d$.  It cannot generate
an additional independent node rotation.

More explicitly, let $J_{C,i}$ be the actual rooted/deck/coordinate
comparison from each component leg to its original reference
leg. A deck change $d$ first changes its sector to the strict pullback
$T_dW_i$; the common based comparison then identifies the three
corresponding typed legs with those of the original DH triple. We
write $J_{C,i}$ for this complete, typed frame map, rather than
silently identifying $T_dW_i$ with $W_i$ as twisted modules.
Transporting $\lambda_i^{\mathrm{DH}}$ along that explicit
comparison gives, on actual leg vectors $v_1,v_2$ and output
restricted dual $v_3'$, exactly
\[
 \lambda_C^{\mathrm{act}}(v_3',v_1,v_2)
 =\left\langle (J_{C,3}^{-1})^t v_3',\,
     \mathcal Y_i(J_{C,1}v_1,1)J_{C,2}v_2\right\rangle.
\]
The intertwining operator on the right is literally the fixed
$\mathcal Y_i$ from \eqref{eq:productcoeval61} or
\eqref{eq:iteratecoeval61}, with its $\log1=0$ branch; neither an
operator nor its matrix coefficient has been redefined. The typed
frame maps also carry the current states according to
Lemma~\ref{lem:relativecoordinates}.
The operators $J_{C,i}$ are the actual maps supplied by the based
cover comparison, the specified Kummer root maps, and
Lemma~\ref{lem:relativecoordinates}; the output always uses the
inverse transpose. At an output marking the local coordinate is
$s=Z^{-1}$ on an outside component, or $s=\xi^{-1}$ on a bubble;
the contragredient observation-current state change is, for
$b\in V_d$,
\[
 e^{sL(1)}(-s^{-2})^{L(0)}b
 =(-1)^d s^{-2d}\sum_{j=0}^{d}\frac{s^j}{j!}L(1)^jb.
\]
It is this full finite polynomial, together with the inverse
transpose on the output, that is transported; replacing it by the
leading $(-1)^d s^{-2d}$ term would not give the same dual section.  Under the initial reference identifications of
Definition~\ref{def:commonoutsetframing}, this formula is exactly the
original corresponding $\ell_P,r_P,\ell_I$ or $r_I$.  In other
framed representatives it is the displayed \emph{covariant} formula,
not an assertion that the untransported scalar expression is invariant.

We verify the node pairing with the \emph{jointly induced}
three-point reference data; no independent DH node-root reference is
hidden in Definition~\ref{def:completeintertwiner}.  That definition
fixes the operators and their preferred analytic continuations, not
an additional choice of two unramified Hurwitz lifts or node roots.
Lemma~\ref{lem:singleinsertion} realizes the resulting fixed
three-point dual sections on a cover equipped with auxiliary based
paths and marked roots.  These are precisely the auxiliary choices
which Definition~\ref{def:commonoutsetframing} makes jointly at the
outset.

For $U=M$ (product) or $U=N$ (iterate), let $U_a$ be the typed
outside node leg of the actual normalized cover and $U_a'$ its paired
bubble output leg.  They have inertia $h=\operatorname{tw}(U)$ and
$h^{-1}$, respectively; the pairing $U_a'\times U_a\to\C$ is the
restricted contragredient pairing.  On the physical punctured collar
choose \emph{one} lifted point $\widetilde q$ above
$x=y=\sqrt\rho>0$ (after a harmless finite positive base extension
if needed).  Restrict its based stems to the outside and bubble and
then normalize the two components.  As these are restrictions of one
lift, any change of this lift by a deck element $d\in H$ acts on both
legs by the same strict pullback, with types
$T_dU_a\in\cC_{dhd^{-1}}$ and
$(T_dU_a)'=T_d(U_a')\in\cC_{dh^{-1}d^{-1}}$; there is no
identification $T_dU\simeq U$ as $V$-modules when their twist types
differ.  The reference component covers are the pointed three-point
covers obtained from these restrictions in the coordinates $Z$,
$\xi=Z/\rho$ (product), or $Z$, $\xi=(Z-1)/\rho$
(iterate).  Their ordered meridians and preferred observation paths
are exactly those checked above.  Based-cover uniqueness identifies
the corresponding \emph{pointed covers}; it is not used to identify
conformal blocks.  Instead we pull back the four dual sections
\emph{already defined} by the fixed $\mathcal Y_i$ through
Lemma~\ref{lem:singleinsertion}.  Since the insertion in each
normalized three-point coordinate remains $1$ with $\log1=0$, this
operation changes no coefficient $\ell_P,r_P,\ell_I,r_I$.

The root and frame comparison is as follows.  Use the common
balanced equations
\[
 x=u^n,\qquad y=v^n,\qquad uv=\tau,
 \qquad\tau^n=\rho,
 \qquad n=\operatorname{ord}(h),
\]
and the positive logarithms retained from $\widetilde q$.  The
outside and bubble reference roots $u_r,v_r$ are \emph{the images of
the same} rooted collar in these two exact component coordinates;
thus $u_r^n=x$, $v_r^n=y$ and $u_rv_r=\tau$.  This last equality is
a relation in the jointly transported cover, not an equality
postulated for independent root torsors.  Any comparison of the
actual roots with this joint reference preserves the same coarse
$x,y$ and the same lifted plumbing parameter, so it is necessarily
\[
 (u,v)=(\zeta^r u_r,\zeta^{-r}v_r),
 \qquad\zeta=e^{2\pi i/n},\qquad r\in\Z/n\Z.
\]
There are no positive-coordinate jets in this comparison: $x=Z$
(or $Z-1$), $y=\xi^{-1}$ are precisely the reference node
coordinates.  By Lemma~\ref{lem:relativecoordinates}, on
$U_{h_U+\gamma}$ the outside root gauge is
$S_r=e^{-2\pi i r(L_U(0)-h_U)}$ and acts by
$e^{-2\pi i r\gamma}$.  The bubble output is the restricted
contragredient of this same node leg, so the actual opposite
root and the same typed deck restriction act by $(S_r^{-1})^t$.
In the common reference types the two comparison maps therefore have
the legitimate directions
\[
 A=J_{o,U}:U_a\longrightarrow U,
 \qquad B=J_{b,U'}:U_a'\longrightarrow U',
 \qquad B=(A^{-1})^t.
 \tag{\ensuremath{\mathrm{E0\text{-}pair}}}
\]
Here the equality is taken after the common strict $T_d$
pullback on both legs; prior to that operation $A$ and $B$ have
different twisting-sector targets and cannot be composed as
untyped module automorphisms.  Equivalently, the square expressing
$\langle B\varphi,A m\rangle_{U'=U^*}
 =\langle\varphi,m\rangle_{U_a'}$ commutes for every $m,\varphi$.
The equality follows from the \emph{explicit} annulus automorphism
$(u,v)\mapsto(\zeta^ru,\zeta^{-r}v)$ and the common deck restriction,
not from an assertion that any two independently chosen reference
frames are automatically paired.

Both the strict pullback and $S_r$ preserve $L(0)-h_U$; thus $A$
commutes with the relative-weight projectors, and the same is true of
$B$.  Choose a basis $\{e_{\gamma,j}\}$ of
$(U_a)_{h_U+\gamma}$ and its restricted-dual basis.  Then
$\{Ae_{\gamma,j}\}$ and
$\{(A^{-1})^t e_{\gamma,j}^*\}$ are precisely dual bases
in the original abstract module $U$ and $U'$.  Therefore,
coefficient by coefficient,
\[
 (A\otimes B)\mathscr E_\gamma^{U,a}
 =\sum_j Ae_{\gamma,j}\otimes(A^{-1})^t e_{\gamma,j}^*
 =\mathscr E_\gamma^{U,\mathrm{ref}},
 \qquad
 B_\gamma A_\gamma^t=\id_{U_{h_U+\gamma}^*}.
\]
The displayed identities apply to both $U=M$ and $U=N$; no
$L(j>0)$ operator occurs in the \emph{node module} comparison since
the two exact node coarse coordinates have not changed.  The
$L(1)$ terms in the output opposite-current formula above act on
the observed $V$-state and remain fully present there.

To check the fixed coefficients rather than only a Ward solution,
substitute these two dual bases directly into the grade-$\gamma$
inner sums of \eqref{eq:productcoeval61} and
\eqref{eq:iteratecoeval61}.  In the jointly induced normalized
external frames they become, respectively,
\[
 \sum_j\ell_P(d',a,Ae_{\gamma,j})
            r_P((A^{-1})^t e_{\gamma,j}^*,b,c),
 \qquad
 \sum_j\ell_I(d',Ae_{\gamma,j},c)
            r_I((A^{-1})^t e_{\gamma,j}^*,a,b),
\]
which are exactly the original $\gamma$-th finite sums, by
basis independence.  If one uses a different external framing, the
four external arguments (including the output restricted dual) must
be transported by the corresponding $J_{C,i}$; the naked coefficient
is not asserted to be invariant.  The physical-to-normalized input
scalings are already the specified $D_P(\rho),D_I(\rho)$, and the
nonintegral lowest-weight factor is the unchanged
$\rho^{h_U}=\exp(h_U\ell(\rho))$ of
\eqref{eq:physicalscaling61}.  Thus neither an IO nor its chosen
$\log1=0$ branch has been renormalized.

If independently chosen three-point node
references $A_0:U_a\to U_r$ and $B_0:U_a'\to U_r'$ are used instead,
the corresponding finite sum contains the genuine typed compensator
$K_\gamma=B_{0,\gamma}A_{0,\gamma}^{t}$ (with the external gauges).
The present argument does not claim $K_\gamma=1$ for arbitrary such
independent references.  Those auxiliary choices are not imposed by
Definition~\ref{def:completeintertwiner} or
Lemma~\ref{lem:singleinsertion}; the common outset prescription
specifies the paired representatives before either original IO is
applied.  No E0/E1 Ward, associativity, convergence of two-IO
products, or coherence is used.
 \end{proof}

\begin{lemma}\label{lem:nodewindow}
Fix a current $\sigma$ and $m=m(\sigma)$ as in Lemma~\ref{lem:canonicalnodecurrent}, and let $uv=t^a$, $x=u^N$, $y=v^N$, $e=aN$ at a source node. For a homogeneous $b\in E\cap V_d$ and a term $t^{-m}t^\ell u^rv^s b$ of its balanced frame, set
\[
 r,s,\ell\in\N,\qquad k=r-s,\qquad
 n=d-1+\frac{k}{N},\qquad \delta=e(d-n-1)=-ak.
\]
The contribution of this term and a grade-$p$ summand of $\Gamma_M(t^e)$ to the coefficient of $t^L$ is zero unless
\begin{equation}\label{eq:nodewindow}
 \begin{gathered}
 p=e\gamma\in\N,\qquad
 a\max(0,k)\le p\le L+m-a\max(0,-k),\\
 |k|\le\left\lfloor\frac{L+m}{a}\right\rfloor.
 \end{gathered}
\end{equation}
In particular, for every $L\in\Z$ the sum over modes, relative grades, and monomials is finite. Its two incident-node actions, including every $L(1)$ correction, cancel coefficientwise in $M\otimes M'$.
\end{lemma}
\begin{proof}
Write $M=\bigoplus_{\gamma\in\Lambda_M}M_{h_M+\gamma}$, $\Lambda_M\subseteq N^{-1}\N$, and
\[
 P_p=\mathscr E^M_{p/e},\qquad
 \Gamma_M(t^e)=\sum_{p\in e\Lambda_M}t^pP_p,
 \qquad P_p=0\quad(p\notin e\Lambda_M).
\]
Thus $e\Lambda_M\subseteq\N$. The $u$-branch expression is $t^{as+\ell-m}u^kb\,x^{d-1}dx$, whose selected-root residue extracts $b_n$ by \eqref{eq:canonicalpuremode}. Since $b_n:M_{h_M+p/e}\longrightarrow M_{h_M+(p+\delta)/e},\; \delta=-ak$,
its action vanishes for $p+\delta<0$. Its $t$-degree in the weighted contraction equals $p+as+\ell-m$. Consequently $p\ge a\max(0,k),\; s\ge\max(0,-k),\; p+as+\ell=L+m$. These relations imply \eqref{eq:nodewindow}; if $L+m<0$, no
term exists. For fixed $L$, the allowed integers $k,p$ are finite, and $as+\ell=L+m-p$ leaves only finitely many pairs $(s,\ell)\in\N^2$ and $r=s+k\in\N$. Since $E$ is finite
dimensional, the same bound applies to a complete element of $t^{-m}\bigl(O[[u,v]]/(uv-t^a)\bigr)\otimes E$.

Write $b_n^t$ for the algebraic transpose of $b_n$. Finite homogeneous duality gives, in $M\otimes M'$,
\begin{equation}\label{eq:node63transpose}
 (b_n\otimes\id)P_p
 =(\id\otimes b_n^t)P_{p+\delta},\qquad
 (b_n\otimes\id)\Gamma_M(t^e)
 =t^{-\delta}(\id\otimes b_n^t)\Gamma_M(t^e).
\end{equation}
Both identities are understood coefficientwise with the zero-grade convention. By the contragredient formula
\cite[Section~2]{DH}\cite[Section~3]{DLM}, for $b\in V_d$,
\begin{equation}\label{eq:node63opposite}
 b_n^t=(-1)^d\sum_{j=0}^{d}\frac1{j!}
 \bigl(L(1)^jb\bigr)^{M'}_{2d-j-n-2}.
\end{equation}
Here $d\in\N$, $L(1)^{d+1}b=0$, and $M'$ is the $g^{-1}$-twisted restricted contragredient. By
\eqref{eq:balancednodeframe} the paired $v$-branch is
\[
 t^{ar+\ell-m}v^{-k}
 \sum_{j=0}^{d}\frac{(-1)^{d+1}}{j!}
 (L(1)^jb)y^{d-j-1}dy.
\]
Its $j$-th selected-root residue extracts $(L(1)^jb)^{M'}_{d-j-1-k/N}$, and
\[
 d-j-1-\frac{k}{N}=2d-j-n-2,
 \qquad ar=as+ak=as-\delta.
\]
Hence the $v$-branch action is exactly $-t^{-\delta}(\id\otimes b_n^t)$ times the $u$-branch scalar coefficient. Combining this with \eqref{eq:node63transpose} gives $\bigl[(b_n\otimes\id) -t^{-\delta}(\id\otimes b_n^t)\bigr] \Gamma_M(t^e)=0$.
The minus-side grade is $p+\delta$; its total $t$-degree is $(p+\delta)+ar+\ell-m=p+as+\ell-m$, so its coefficient is in the same finite window. Finally, $n\in\alpha+\Z$ for $gb=e^{2\pi i\alpha}b$ precisely when
$k/N\in\alpha+\Z$. As $gL(1)=L(1)g$, all states $L(1)^jb$ have the same inertia character, while the opposite modes act in the inverse-twisted module. Summing over the finitely many node orbits and a finite homogeneous basis of $E$ proves the
assertion.
\end{proof}

\begin{lemma}\label{lem:physicalcontraction}
For the product and iterate of \eqref{eq:productcoeval61} and
\eqref{eq:iteratecoeval61}, the node contractions of
Lemma~\ref{lem:nodewindow} and the three-point Ward identities recover the original four-point seed Ward identities in the physical coordinates. The equalities hold coefficientwise on the
outset rooted collision branches of Definition~\ref{def:commonoutsetframing}, with their actual component comparison of Lemma~\ref{lem:E0actualsourcecomparison}. In addition, both formal contractions satisfy the independent $L(-1)$-derivative identities in their two moving insertion variables.
\end{lemma}
\begin{proof}
Let $a,b,c,d'$ have weights $\lambda_A,\lambda_B,\lambda_C,\lambda_D$, respectively, and put $\Delta=\lambda_D-\lambda_A-\lambda_B-\lambda_C$. For a product summand of intermediate weight $h_M+\gamma$,
\begin{align*}
 &R^{\lambda_D-\lambda_A-h_M-\gamma}
 (R\rho)^{h_M+\gamma-\lambda_B-\lambda_C}\\
 &\hspace{2em}=
 R^\Delta\rho^{h_M-\lambda_B-\lambda_C}\rho^\gamma.
\end{align*}
For an iterate summand of weight $h_N+\gamma$ the corresponding identity is
\[
 (R\rho)^{h_N+\gamma-\lambda_A-\lambda_B}
 R^{\lambda_D-h_N-\gamma-\lambda_C}
 =R^\Delta\rho^{h_N-\lambda_A-\lambda_B}\rho^\gamma.
\]
These are identities of individually finite homogeneous contractions, assuming the compatible logarithms in Lemma~\ref{lem:actualonecontraction}. Put $(M_P,M_I)=(M,N)$. In particular, at $R=1$,
\begin{equation}\label{eq:node64physical}
 \Phi_\star(\rho)=\rho^{h_{M_\star}}
 \Psi_\star(\rho)\circ D_\star(\rho),
\end{equation}
where $D_P,D_I$ are the weight-diagonal operators defined in
\eqref{eq:physicalscaling61} and its preceding display.

We compare the \emph{geometrically defined} external
residues.  Fix the given generic $H$-invariant seed current $s$;
its restrictions to the actual normalized components are those
of Lemma~\ref{lem:canonicalnodecurrent}, not newly chosen Ward
currents.  In the product use the coarse outside coordinate $Z$,
and on the bubble put $x_B=Z-\rho$, $x_C=Z$ and
$\xi_B=\xi-1$, $\xi_C=\xi$, where $Z=\rho\xi$.
In the iterate, the outside coordinate is still $Z$, with its
original IO inserted at $Z=1$; on the bubble put
$x_A=Z-(1+\rho)$, $x_B=Z-1$ and
$\xi_A=\xi-1$, $\xi_B=\xi$, where
$Z-1=\rho\xi$.  Thus for each bubble external marking
$x_i=\rho\xi_i$ exactly.  On the other external marking and
at infinity the physical and outside coordinates coincide for
$R=1$.

Let $b\in E\cap V_d$ be homogeneous, and let
$n\in r_b/N_i+\Z$ be an allowed twisted mode at such a
bubble marking.  In its centered physical state-current frame
consider a coefficient $c_{i,b,n}(t)x_i^n b\,(dx_i)^{1-d}$.
This density notation records the exact scaling of the state-current
frame under a \emph{linear} coordinate change: since
$x_i=\rho\xi_i$, the induced component coefficient is
$\rho^{n+1-d}c_{i,b,n}(t)$.  No positive Virasoro mode occurs
in this linear transition, even when $L(1)b\ne0$; the complete
positive-mode transition at the \emph{node} is instead the
balanced transition \eqref{eq:balancednodeframe} already used
in Lemma~\ref{lem:nodewindow}.  To check the orbit
normalization explicitly, for the selected ramified root
$u_i^{N_i}=x_i$ one has
\[
 \frac1{N_i}\Res_{u_i=0}
 u_i^{N_i n}Y_{M_i}^{g_i}(b,u_i^{N_i})
 N_i u_i^{N_i-1}du_i=b_n.
\]
Under the \emph{inherited} root change
$u_i=\rho^{1/N_i}\widetilde u_i$ the same formula in the
component root has coefficient
$\rho^{n+1-d}c_{i,b,n}(t)$; the root powers use the fixed
$\ell(\rho)$ and hence produce no unrecorded root phase.
This calculation is compatible, modulo every $t^{J+1}$, with the
actual marked-lift comparisons $\iota_{v,J}$ in the proof of
Lemma~\ref{lem:canonicalnodecurrent}.  In particular it
compares the same current $s$, including all Laurent
coefficients of its extension, and does not infer the component
action merely from abstract covariance.

Put $\Delta_{b,n}=d-n-1$.  The original twisted-module
commutator gives, on a vector in any homogeneous relative-weight
space,
\begin{equation}\label{eq:node64conjugation}
 [L(0),b_n]=\Delta_{b,n}b_n,\qquad
 \rho^{-L(0)}b_n\rho^{L(0)}
 =\rho^{-\Delta_{b,n}}b_n
 =\rho^{n+1-d}b_n.
\end{equation}
Indeed, if $w_i$ has weight $\lambda$ then both sides of
$\rho^{-L(0)}b_nw_i
=\rho^{n+1-d}b_n\rho^{-L(0)}w_i$
are $\rho^{-(\lambda+\Delta_{b,n})}b_nw_i$.
Consequently, writing $R_{i,\mathrm{phys}}(s)$ for the
orbit-normalized residue of \eqref{eq:canonicalmarkedres} on
the actual external marking, and $R_{i,\mathrm{vert}}(s)$
for the restriction of that same current to its normalized
component, one obtains
\[
 R_{i,\mathrm{vert}}(s)
   =\rho^{-L_i(0)}R_{i,\mathrm{phys}}(s)
    \rho^{L_i(0)}
 \quad(i=B,C\text{ for }P;\ i=A,B\text{ for }I)
\]
coefficientwise on algebraic inputs.  For the unscaled outside
input and the output $D'$ the comparison is the identity.
The latter statement also uses the prescribed \emph{restricted
contragredient} frame: an input frame map $G$ induces
$(G^{-1})^t$ on its dual; for $R=1$ the output map itself is
$G=\id$, and the internal output-dual map is already paired
with the node input by Lemma~\ref{lem:E0actualsourcecomparison}.
These observations prove the independent geometric--algebraic
comparison
\begin{equation}\label{eq:node64externalcomparison}
 D_\star(\rho)R_{\mathrm{phys}}(s)w
 =R_{\mathrm{vert}}(s)D_\star(\rho)w,
 \qquad \star\in\{P,I\},
\end{equation}
for each finite coefficient of the rooted formal expansion.
The equality is not asserted as an unqualified continuous-operator
identity on an infinite product completion.  The scalar powers
are kept on the chosen branch, including nonintegral external
lowest weights.

The Ward cancellation also follows directly, without invoking
\eqref{eq:node64externalcomparison}.  After multiplying $s$ by
the $t^{m(s)}$ supplied by Lemma~\ref{lem:canonicalnodecurrent},
restrict its extension to the two actual normalized components.
Apply Proposition~\ref{prop:binaryward} separately to the
outside and bubble original three-point functionals,
identified by Lemma~\ref{lem:E0actualsourcecomparison}.
At a fixed coefficient $t^L$, Lemma~\ref{lem:nodewindow}
restricts the intermediate grades and all node Laurent indices
to its finite window.  Pair the two three-point Ward equalities
with the \emph{unchanged} degree-$\gamma$ coevaluation, weighted
by $t^{e\gamma}$.  The incident node terms are respectively
$(b_n\otimes\id)\Gamma_{M_\star}(t^e)$ and
$-t^{-\delta}(\id\otimes b_n^t)\Gamma_{M_\star}(t^e)$,
including the $L(1)^j$-terms; they cancel by
\eqref{eq:node63transpose}--\eqref{eq:node63opposite}.
The surviving external terms are exactly the restrictions of $s$
just computed, hence give
\begin{equation}\label{eq:node64ward}
 \Psi_\star(\rho)
  \bigl(R_{\mathrm{vert}}(s)D_\star(\rho)w\bigr)=0,
 \qquad
 D_\star(\rho)R_{\mathrm{phys}}(s)w
 =R_{\mathrm{vert}}(s)D_\star(\rho)w
\end{equation}
in each finite homogeneous coefficient.  As every $t^L$ involves
only finitely many current monomials, node modes and weights, no
analytic convergence or exchange of infinite sums is used;
localizing the extension by $t^{-m(s)}$ recovers the original
current.  Finally \eqref{eq:node64physical} gives, for the
original physical four-point Ward action,
\[
 \Phi_\star(\rho)(R_{\mathrm{phys}}(s)w)
 =\rho^{h_{M_\star}}
 \Psi_\star(\rho)
   (R_{\mathrm{vert}}(s)D_\star(\rho)w)=0.
\]
The internal ground-weight scalar $\rho^{h_{M_\star}}$
is retained.  For a further common scale $R$, the output
restricted dual transforms by the inverse transpose of the
input scaling, namely $R^{L_{D'}(0)}$ versus $R^{-L(0)}$
on each input.  All the external centered-coordinate changes
are affine dilations; their Schwarzian is zero.  This does not
discard the ramified-root Schwarzian already included in the
full state-current gauge of \eqref{eq:canonicalmarkedres}.

For the derivative statement, use the $L(-1)$-derivative property and translation commutator of a single twisted intertwining operator \cite[Definition~2.8]{DH}:
\[
\partial_x\mathcal Y(a,x)=\mathcal Y(L(-1)a,x),
 \quad
 L(-1)\mathcal Y(a,x)-\mathcal Y(a,x)L(-1)
 =\mathcal Y(L(-1)a,x).
\]
For a product, differentiating each of the two insertion coordinates gives respectively the insertion of $L(-1)a$ and $L(-1)b$, with the intermediate projections held fixed.
For an iterate let $q=z_1-z_2$, and write $F=\langle d',\mathcal Y_3(\mathcal Y_4(a,q)b,z_2)c\rangle$. Then, coefficientwise,
\begin{align*}
 \partial_{z_1}F
 &=\langle d',\mathcal Y_3(
    \mathcal Y_4(L(-1)a,q)b,z_2)c\rangle,\\
 \partial_{z_2}F
 &=\langle d',\mathcal Y_3(
    L(-1)\mathcal Y_4(a,q)b
    -\mathcal Y_4(L(-1)a,q)b,z_2)c\rangle\\
 &=\langle d',\mathcal Y_3(
    \mathcal Y_4(a,q)L(-1)b,z_2)c\rangle.
\end{align*}
For each prescribed relative power, the $L(0)$-weight equation fixes the intermediate grade; insertion of $L(-1)$ shifts that grade by one. Thus each displayed coefficient is a finite
homogeneous contraction, and the differentiation is performed coefficientwise in the formal series. The chain rule and the Leibniz rule then give \eqref{eq:channeldiff34} on the physical placement
chart. The substitution $\rho=t^e$ is injective on the relative grade support. A root rotation acts as $S=e^{-2\pi i(L_M(0)-h_M)}$ on $M$ and $(S^{-1})^t$ on $M'$, fixing each $\mathscr E^M_\gamma$; the iterate uses $N$ in place
of $M$. Thus the identities respect the retained framing.
\end{proof}

\begin{theorem}\label{thm:smallwardautomatic}
Fix the single outset product--iterate framing of Definition~\ref{def:commonoutsetframing}, and let $\star\in\{P,I\}$.  In its preferred collision-root formal neighborhood, the original contraction $\widehat\Phi_\star$ of \eqref{eq:productcoeval61} or \eqref{eq:iteratecoeval61}, normalized by \eqref{eq:physicalscaling61}, is a coefficientwise well-defined formal scalar functional.  For every canonical seed current $s$ pulled back to this formal neighborhood and every external algebraic tensor $w$,
\begin{equation}\label{eq:actualfourward65}
 \widehat\Phi_\star\bigl(R(s)w\bigr)=0
 \quad\text{coefficientwise in the chosen root parameter}.
\end{equation}
For each placement derivation $\delta$ the independent formal derivative identity and differential saturation give
\begin{equation}\label{eq:actualfourhorizontal65}
 \partial_\delta\widehat\Phi_\star(w)
 =\widehat\Phi_\star(\nabla_\delta w),\qquad
 \widehat\Phi_\star(\cJ_{\diff})=0.
\end{equation}
Consequently this is a horizontal \emph{formal} functional on the completed finite differential quotient $\widehat{\cQ}$.  The same initial rooted choice is used for every typed intertwining operator and every homogeneous coefficient.  No coefficientwise statement about the \emph{untransported} original series on arbitrary other framed charts is claimed here; its covariant analytic continuation is established in Theorem~\ref{thm:E0covariantcontinuation} below.
\end{theorem}
\begin{proof}
The geometric differential quotient and its connection were constructed without a four-point Ward hypothesis in Theorem~\ref{thm:finiteconn}; their regular singularity was established in Theorem~\ref{thm:canonicalregular}.  Fix $\star$, a current $s$, and homogeneous $w$.  The source of the two original three-point functionals is identified on the \emph{same actual based cover} by Lemma~\ref{lem:E0actualsourcecomparison}.  The selected nodal extension of Lemma~\ref{lem:canonicalnodecurrent} uses the base change \eqref{eq:E0sufficientbasechange}; it does not change the external rooted/logged sources.  After multiplying by $t^{m(s)}$, restrict the given current to the two normalized components.  For every coefficient of $t^L$, Lemma~\ref{lem:nodewindow} gives the finite window
\begin{equation}\label{eq:actualwindow65}
 \begin{gathered}
 p=e\gamma\in\N,\qquad
 a\max(0,k)\le p\le L+m-a\max(0,-k),\\
 |k|\le\lfloor(L+m)/a\rfloor,\qquad
 p+as+\ell=L+m,
 \end{gathered}
\end{equation}
where $e=an$ with $n$ the node inertia order.  Thus applying the two three-point Ward identities of Proposition~\ref{prop:binaryward} involves only finitely many states, intermediate grades and Laurent monomials.  All node terms, including the $L(1)^j$ corrections and their inverse-transpose actions, cancel by \eqref{eq:node63transpose} and \eqref{eq:node63opposite}.  As the component germs are the original $\ell_\star,r_\star$, the remaining external action satisfies
\begin{equation}\label{eq:relativeward65}
 \Psi_\star(\rho)\bigl(R_{\mathrm{vert}}(s)D_\star(\rho)w\bigr)=0.
\end{equation}
Here $R_{\mathrm{vert}}$ is obtained from the \emph{actual restrictions} of $s$, not defined by conjugating its physical action.  The full coordinate law of Lemma~\ref{lem:relativecoordinates} and the mode identity
$\rho^{-L(0)}v_n\rho^{L(0)}=\rho^{-(\wt(v)-n-1)}v_n$
apply separately to each $L(1)^jv$ state and to the output restricted contragredient.  Therefore, in these finite contractions, Lemma~\ref{lem:physicalcontraction} gives
\[
 D_\star(\rho)R_{\rm phys}(s)
 =R_{\rm vert}(s)D_\star(\rho),\qquad
 \widehat\Phi_\star(R_{\rm phys}(s)w)
 =\rho^{h_{M_\star}}\Psi_\star(R_{\rm vert}(s)D_\star w)=0.
\]
The factors $\rho^{h_M},\rho^{h_N}$ remain on the chosen logarithmic branch, whether or not the lowest weights are integral.  Injectivity of the finite root substitution on the original lower-bounded relative-grade support descends this coefficientwise equality from the $t$-expansion.  Lemma~\ref{lem:physicalcontraction} separately proves the formal $L(-1)$ derivative identities; Lemma~\ref{lem:channelderivative} and Proposition~\ref{prop:factor} then yield \eqref{eq:actualfourhorizontal65}.  Neither convergence nor associativity nor an all-chart Ward assertion has entered the argument.
\end{proof}

\begin{corollary}\label{cor:smallconvergenceautomatic}
For any two actual twisted intertwining operators with ordinary intermediate module, the original intermediate-grade contractions converge absolutely and normally on compact subsets of the respective physical chambers
\[
 \mathcal A_P=\{(z_1,z_2):|z_1|>|z_2|>0\},\qquad
 \mathcal A_I=\{(z_1,z_2):|z_2|>|z_1-z_2|>0\},
\]
on every logarithmic branch obtained by continuing the outset physical/rooted data within these chambers. Their sums are the evaluations of analytic horizontal sections of the canonical differential Ward quotient on the corresponding branches; no assertion is made that an untransported original grade sum converges outside its physical chamber.
\end{corollary}
\begin{proof}
Fix a simple intermediate module and homogeneous external and restricted-dual vectors. Put
\begin{equation}\label{eq:ratiochoice66}
 (R,q)=
 \begin{cases}
 (z_1,z_2/z_1),&\star=P,\\
 (z_2,(z_1-z_2)/z_2),&\star=I.
 \end{cases}
\end{equation}
Then $R\ne0$, $0<|q|<1$ and the third moving-puncture difference is $R(1-q)$ or $R(1+q)$, respectively. Thus the configuration is nonsingular throughout the punctured unit disk.
With $R=1$ let $S_\star\to
\mathbb A^1\setminus\{0,\epsilon_\star\}$ be the finite
\'{e}tale framed Hurwitz cover of Lemma~\ref{lem:framedcharts},
where $(\epsilon_P,\epsilon_I)=(1,-1)$. Its analytic pullback to $0<|q|<1$ has finite permutation monodromy. Choose $d\in\N_{>0}$ divisible both by the order of this monodromy and by the denominators of the relative grades of the fixed simple
intermediate module; set $q=\xi^d$. On the punctured $\xi$-disk the finite cover splits into single-valued analytic sheets. This operation removes the geometric permutation monodromy; it
does \emph{not} remove the arbitrary complex lowest-weight power $q^{h_{M_\star}}$ or its specified logarithm.

Let $Z$ be the normalization of the pulled-back algebraic Hurwitz curve, restricted to the irreducible component containing the selected analytic sheet. Its morphism to the $\xi$-line is finite. Over $0<|\xi|<1$ the chosen sheet is an
analytic section. It extends across $\xi=0$: on an affine neighbourhood of $0$, integral generators of the finite algebra satisfy monic equations with holomorphic coefficients and hence
are locally bounded along the section; the removable-singularity theorem extends each generator. These extensions retain the monic relations and yield a section of $Z^{\mathrm{an}}$ above $\xi=0$.
At its point $p$, the normal curve $Z$ is smooth. If $\xi=\varepsilon(s)s^e$ in its local parameter $s$, with $\varepsilon(0)\ne0$, the existence of a section implies $e=1$. Consequently the selected analytic sheet, with its point
$p$, is a disk with coordinate $\xi$. Since the original framed cover is \'{e}tale away from $q=0,\epsilon_\star$, it has no other branching on $|\xi|<1$.

The pulled-back finite-rank Ward connection is algebraic and integrable on this smooth algebraic curve away from its boundary; Theorem~\ref{thm:canonicalregular} and stability of regular
singularity under finite pullback give its regular singularity at $p$ (cf.\ \cite[Chapter~II, Proposition~5.4]{Deligne}). Each fixed external tensor yields an algebraic test section on the smooth
framed chart. It extends meromorphically in the canonical logarithmic extension at $p$ by Proposition~\ref{prop:algebraictests}. Theorem~\ref{thm:smallwardautomatic}, now applied \emph{only to this preferred collision-root neighborhood}, and Proposition~\ref{prop:factor} give the formal
horizontal identity independently of any all-chart Ward conclusion. By \eqref{eq:ratiochoice66} and Lemma~\ref{lem:treesupport}, its normalized grade series has a fixed fractional-power shift and only nonnegative integral powers of $\xi$ after the chosen root substitution. Hence the adapted
ratio condition of \eqref{eq:adaptedratio} and the finite-Nilsson condition of Lemma~\ref{lem:jointsupport} hold at $p$. Theorem~\ref{thm:formalconv} identifies this formal series with a
convergent local horizontal germ. Theorem~\ref{thm:globalratio} and Proposition~\ref{prop:algebraictests} continue that germ throughout the rooted disk and give absolute normal convergence of its finite exponent-class Taylor series there.

More precisely, the homogeneity calculation of Lemma~\ref{lem:actualonecontraction} gives
\[
\Phi_\star(R,q)=R^{\Delta}
 q^{\kappa_\star}\sum_{\gamma\in\Lambda_{M_\star}}
 C_{\star,\gamma}\,q^\gamma,\qquad
 \kappa_P=h_M-\lambda_B-\lambda_C,\quad
 \kappa_I=h_N-\lambda_A-\lambda_B.
\]
Each $C_{\star,\gamma}$ is the original finite-dimensional contraction of the grade-$\gamma$ space. Under $q=\xi^d$, the map $\gamma\mapsto d\gamma\in\N$ is injective. Thus distinct relative grades are distinct Taylor exponents after the fixed
factor $q^{\kappa_\star}$ has been extracted. Uniqueness of formal coefficients identifies the coefficients of the analytic expansion with $C_{\star,\gamma}$, not with regrouped sums.
Proposition~\ref{prop:gradesnormal}, or equivalently the Cauchy estimates of Lemma~\ref{lem:Reinhardt} on the $\xi$-disk, now gives
\[
\sum_{\gamma\in\Lambda_{M_\star}}
 \sup_{\xi\in K}|C_{\star,\gamma}\,\xi^{d\gamma}|<\infty
 \qquad(K\Subset\{0<|\xi|<1\}).
\]
The fixed factors $R^\Delta q^{\kappa_\star}$ are bounded on compact subsets of a chosen logarithmic branch of the physical chamber. Reinstating them, and summing the finitely many simple
intermediate choices, proves the asserted normal convergence on compact subsets of $\mathcal A_\star$. The resulting analytic functions have exactly the original coefficientwise Taylor germs of Theorem~\ref{thm:smallwardautomatic}; their extension to other framed charts is the separate covariant assertion of Theorem~\ref{thm:E0covariantcontinuation}.
\end{proof}

\subsection{The full-state differential quotient}

Fix an affine smooth framed chart $S=\operatorname{Spec}A$ and the actual based $H$-cover $p:C\to S$ of Section~\ref{sec:geometricconventions}. Let $\cT=A\otimes_{\C}T$ with the physical connection $\nabla$ of
\eqref{eq:Tconn}. Since $V$ is of CFT type, put
\[
 V_{\le D}=\bigoplus_{n=0}^{D}V_n\quad(D\in\N),\qquad
 D_0\in\N,\quad E\subseteq V_{\le D_0}.
\]
For $D\ge D_0$, the finite-dimensional $G$-module $V_{\le D}$ is stable under $L(0)$ and $L(j)$ for $j>0$, since $L(j)V_n\subseteq V_{n-j}$. Thus the coordinate-change action of
\cite[Sections~3--5]{FS} restricts to $V_{\le D}$, including the
finite operator $e^{L(1)}(-1)^{L(0)+1}$. Write $\cF_D$ for the corresponding finite-rank current bundle; the inclusions $V_{\le D}\hookrightarrow V_{\le D+1}$ induce $\cF_D\hookrightarrow\cF_{D+1}$. Let $P$ be the reduced sum of
all marked $H$-orbits, and define
\[
 \begin{aligned}
 \mathscr S_D&=\Gamma(C,\cF_D(*P))^H,\\
 \cJ_D&=\sum_{s\in\mathscr S_D}\sum_{w\in\cT}
                  A\,R_D(s)w\subseteq\cT,\\
 \cJ_D^{\diff}
  &=\sum_{r\ge0}\sum_{\delta_1,\ldots,\delta_r\in
                       \operatorname{Der}_{\C}(A)}
     A\,\nabla_{\delta_1}\cdots\nabla_{\delta_r}\cJ_D,
 \end{aligned}
\]
where $R_D$ is the orbit-normalized residue action of Definition~\ref{def:canonicalmarked} (equivalently
\eqref{eq:canonicalmarkedres}), and $r=0$ denotes $\cJ_D$.
Here $\Gamma(C,\cF_D(*P))=\bigcup_{m\ge0}\Gamma(C,\cF_D(mP))$; no current involves infinitely many state weights. Set
\begin{equation}\label{eq:fullstatequotient62}
 \begin{gathered}
 \cJ_{\infty}^{\diff}
   =\bigcup_{D\ge D_0}\cJ_D^{\diff},\qquad
 \cQ_{\infty}=\cT/\cJ_{\infty}^{\diff},\\
 \cJ_E^{\diff}\subseteq\cJ_D^{\diff}
 \subseteq\cJ_{D+1}^{\diff},\qquad
 \cQ=\cT/\cJ_E^{\diff}\twoheadrightarrow\cQ_\infty.
 \end{gathered}
\end{equation}
Every $\cJ_D^{\diff}$ and their union are $\nabla$-stable by the Leibniz rule, so $\nabla$ induces an integrable connection on their quotients whenever these quotients are formed.

\begin{lemma}\label{lem:fullbinarywardfixedD}
Fix $D\ge D_0$ and an original ordinary twisted intertwining operator
$\mathcal Y$ of type $\binom{W_3}{W_1\ W_2}$, with its already fixed
$z=1$, $\log1=0$ germ (or its specified logarithmic continuation).
On its actual smooth rooted three-point $H$-cover, the original
three-point matrix-coefficient functional annihilates the normalized
residue action of every $H$-invariant current in
$H^0(C_3,\cF_D(*P_3))$.  The assertion includes all finite
$L(j>0)$ coordinate corrections and the restricted contragredient
at infinity; no convergence of a four-point contraction is involved.
\end{lemma}
\begin{proof}
The coordinate-change representation on the state space
$V_{\le D}$ is finite: $L(j)V_r\subseteq V_{r-j}$ for $j>0$,
and $L(j)$ vanishes on $V_r$ when $j>r$.  Thus for an origin-fixing
change of formal coordinate the $\Aut\mathcal O$ action on
$V_{\le D}$ is an algebraic finite-jet matrix. It commutes with
$H\le\Aut(V)$ because the latter preserves the conformal vector.
Its dual and the Jacobian define the finite-rank bundle
$\mathscr V_D^\vee$ and the current bundle
$\cF_D=\mathscr V_D\otimes\omega$ on the smooth cover.

For every $b\in V_{\le D}$ use the \emph{same fixed}
$\mathcal Y$ to define the three single-insertion coefficients
$F_z(x;b)$, $F_0(x;b)$, $F_\infty(x;b)$ of
Lemma~\ref{lem:singleinsertion}.  The Du--Huang single-insertion
analytic duality identifies their continuation along the fixed
rooted paths; this step does not reconstruct an intertwining
operator from Ward solutions.  Their $H$-monodromy and the complete
finite-jet coordinate matrices give an $H$-invariant holomorphic
section $\lambda_D$ of $\mathscr V_D^\vee$ away from the three
marked orbits.  Since $V_{\le D}$ has a finite homogeneous basis,
the lower truncation of the input twisted vertex operators gives
one finite lower Laurent bound at the two finite markings.  At
infinity, for $b\in V_d$ the complete opposite insertion uses
\[
 e^{xL(1)}(-x^{-2})^{L(0)}b
 =(-1)^d x^{-2d}\sum_{j=0}^{d}
   \frac{x^j}{j!}L(1)^j b.
\]
The sum terminates, and lower truncation on the restricted
contragredient $W_3'$ provides a finite bound there as well.
The finite state-current coordinate gauges have finite Laurent
orders and therefore preserve meromorphicity at the markings.
Consequently $\lambda_D$ extends as an algebraic meromorphic
section on the entire proper three-point cover (by GAGA).

For any invariant $s\in H^0(C_3,\cF_D(*P_3))$, the pairing
$\eta=\langle\lambda_D,s\rangle$ is a global invariant
meromorphic one-form, with no poles away from the marked orbits.
In a distinguished $N_i$th-root coordinate $x_i=u_i^{N_i}$,
the coefficient of each permitted twisted mode $b_n$ is exactly
\[
 \frac1{N_i}\operatorname{Res}_{u_i=0}
   u_i^{N_i n}Y^{{g_i}}_{W_i}(b,u_i^{N_i})
                     N_i u_i^{N_i-1}\,du_i=b_n.
\]
The complete $L(1)$ terms and the output inverse transpose are
already included in the current/frame and opposite-field formulas.
Each marked $H$-orbit has $|H|/N_i$ points, hence the usual residue
theorem gives
\[
 0=\frac1{|H|}\sum_{x\in C_3}\operatorname{Res}_x\eta
   =\sum_i\Phi_{\mathcal Y,1}(R_{D,i}(s)w).
\]
Linearity treats nonhomogeneous $s$.  This proves the stated Ward
identity for each fixed $D$; neither Proposition~\ref{prop:binaryward}
on $E$ alone nor a Ward identity for an infinite state bundle is
being used as a substitute.
\end{proof}

\begin{theorem}\label{thm:fullsmallward}
Let $\star\in\{P,I\}$, and let $\widehat\Phi_\star$ be the
\emph{preferred-outset formal} original product or iterate of
Theorem~\ref{thm:smallwardautomatic}, with its fixed rooted/logarithmic
collision branch. For every fixed $D\ge D_0$,
\begin{equation}\label{eq:fullseedward62}
 \Phi_\star(\cJ_D)=0,\qquad
 \Phi_\star(\cJ_D^{\diff})=0,\qquad
 \Phi_\star(\cJ_\infty^{\diff})=0.
\end{equation}
Consequently $\widehat\Phi_\star$ induces a horizontal formal
functional on the completion of $\cQ_\infty$ at this collision root,
using the finite regular-singular lattice of the prior finite quotient.
The assertion is coefficientwise in the original physical
normalization \emph{on the preferred outset branch}; arbitrary framed
charts are treated by Theorem~\ref{thm:E0covariantcontinuation}.
No convergence theorem for a $V_{\le D}$-dependent quotient is
assumed.
\end{theorem}
\begin{proof}
Fix $D$, an invariant $s\in\mathscr S_D$, and an algebraic
homogeneous external tensor $w$. The current space is genuinely
finite dimensional in the state direction: the conformal grading is
finite dimensional at each weight, and every $L(j>0)$ preserves
$V_{\le D}$.  The origin-fixing coordinate matrices (including
$L(1)^j$, the differential Jacobian and the ramified root charts)
are consequently finite matrices.  Therefore $\cF_D$ is a
finite-rank algebraic current bundle.  The $H$ action, the
marking-root constructions and the inclusions
$\cF_D\hookrightarrow\cF_{D+1}$ respect those matrices.

Use the actual admissible degeneration supplied by
Lemma~\ref{lem:canonicalnodecurrent}, not a separately chosen
three-point cover. After its allowed finite DVR and Kummer extension,
put $O=L[[t]]$, $K=L((t))$.  At each source node the complete local
ring is
\[
 B=O[[u,v]]/(uv-t^a),\qquad x=u^N,\quad y=v^N.
\]
The underlying state-current transition computed in
\eqref{eq:balancednodeframe} is, on $V_{\le D}$,
\[
 \Sigma_D=e^{L(1)}(-1)^{L(0)+1}\big|_{V_{\le D}},\qquad
 \Sigma_D b=(-1)^{d+1}\sum_{j=0}^d
                 \frac{L(1)^jb}{j!}\quad(b\in V_d).
\]
In particular no term leaves $V_{\le D}$.  The transition is
constant in the balanced state-current frames, including its
Jacobian; it commutes with the inertia action.  Changing the second
trivialization by $\Sigma_D$ identifies the bundle on the punctured
node with the restriction of $B\otimes V_{\le D}$.  The hypersurface
$B$ is normal and Cohen--Macaulay of dimension two.  Thus the
reflexive-extension equivalence across its closed codimension-two
point (\cite[Tag~0EBJ]{Stacks}) extends this bundle there as a
coherent, locally free node stalk.  Away from nodes the prescribed
finite Kummer marking charts and the original smooth current bundle
supply a coherent sheaf on the entire proper total cover; gluing
with these node extensions gives $\cF_{D,O}$.
This construction concerns the sheaf on the actual total cover,
not the restriction of a reflexive sheaf under a possibly nonflat
special-fiber base change.

Choose an invariant effective external marking divisor $P$ bounding
all generic poles of $s$.  For a suitable finite $m$, $s$ lies in
$H^0(C_K,\cF_{D,O}(mP)_K)$.  Properness and coherence make
$p_*\cF_{D,O}(mP)$ a finite $O$-module.  Since base change here is localization to the generic fiber, flatness of localization gives
\[
 (p_*\cF_{D,O}(mP))\otimes_OK
   \simeq H^0(C_K,\cF_{D,O}(mP)_K).
\]
It follows that $t^{m_s}s$ is a global coherent section for some
$m_s=m_s(s,D)\ge0$.  Averaging by $|H|^{-1}\sum_{h\in H}h$
retains invariance.  No uniform exponent in $s$ or $D$ is required.

Take the coefficient of $t^L$ after the extension.  Near a node a
term in the integral finite-state current has the form
$t^{-m_s+\ell}u^rv^s b$ with $r,s,\ell\ge0$, $b\in V_d$, $d\le D$.
For this paragraph $s$ in the exponent denotes an integer, not the
current itself. Set
$k=r-s$, $n=d-1+k/N$, $e=aN$, and $p=e\gamma$ for
the intermediate relative grade $\gamma\ge0$.  The mode degree is
$d-n-1=-k/N$, so the source and target grade restrictions give
$p\ge a\max(0,k)$.  Since $p+as+\ell=L+m_s$ and
$s\ge\max(0,-k)$, we obtain
\begin{equation}\label{eq:fullnodewindow62}
 \begin{gathered}
 a\max(0,k)\le p\le L+m_s-a\max(0,-k),\\
 |k|\le\left\lfloor\frac{L+m_s}{a}\right\rfloor,\qquad
 p+as+\ell=L+m_s.
 \end{gathered}
\end{equation}
Thus at fixed $(D,L)$ there are only finitely many integers
$k,p,r,s,\ell$, and only finitely many homogeneous state basis
vectors.  All $L(1)^jb$ are in $V_{\le D}$, and the internal
restricted-grade spaces are finite dimensional.  Formula
\eqref{eq:node63opposite} gives all of their transposed node terms,
including the change of index and every $L(1)^j$ correction.
Lemma~\ref{lem:nodewindow} therefore cancels the two node residues
in this same $t^L$ coefficient; this uses neither an infinite
sum interchange nor a bound common to all $D$.

The thickening comparison in
Lemma~\ref{lem:canonicalnodecurrent} transports the coefficientwise
restrictions of $t^{m_s}s$ to meromorphic currents on the actual
normalized smooth three-point component covers.  Its proof is
statewise and its finite-jet/frame constructions above preserve
$V_{\le D}$. Lemma~\ref{lem:E0actualsourcecomparison} identifies
the \emph{fixed original} three-point IO germs in those frames.
Apply Lemma~\ref{lem:fullbinarywardfixedD} on each component,
coefficient by coefficient, to the original outside/bubble
functionals.  The sum of their external and node residues is zero;
after node cancellation only the external residues remain.
The already established \emph{geometric} physical-coordinate
comparison of Lemma~\ref{lem:physicalcontraction} acts on each
$b\in V_{\le D}$ by the identical finite matrices and on output
restricted dual by inverse transpose.  It therefore sends precisely
the remaining component action to the original physical
$R_D(s)$ action, without introducing another formal functional.
The factors $\rho^{h_M}$ or $\rho^{h_N}$ and the logarithmic branch
are retained. Hence
$\Phi_\star(R_D(s)w)=0$ in $\mathscr R_\star$ after the faithful
finite root substitution.  $A$-linearity proves
$\Phi_\star(\cJ_D)=0$.

The independent formal derivative identity of
Lemma~\ref{lem:channelderivative}, also calculated directly from the original
IO $L(-1)$-axiom in the proof of
Theorem~\ref{thm:smallwardautomatic} rather than from a full-state Ward assumption,
reads for every placement derivation $\delta$
\[
 \partial_\delta\Phi_\star(fw)
  =\delta(f)\Phi_\star(w)
     +f\Phi_\star\!\left(\sum_i\delta(z_i)L(-1)^{(i)}w\right)
  =\Phi_\star(\nabla_\delta(fw)),\qquad f\in A.
\]
The coefficient module $\mathscr R_\star$ is the fixed
lower-bounded iterated Laurent/root series module: differentiating
an individual exponent only shifts a finite grade, and every fixed
coefficient of the original contraction is finite dimensional.
Multiplication by expanded coefficients of $A$ and formal
Leibniz differentiation preserve this lower-bounded series module;
no summation of infinitely many contributions to one coefficient is
needed.  For any $j\in\cJ_D$ and any finite list of derivations,
iteration gives
\[
 \Phi_\star(\nabla_{\delta_1}\cdots\nabla_{\delta_r}j)
  =\partial_{\delta_1}\cdots\partial_{\delta_r}
     \Phi_\star(j)=0.
\]
The definition of $\cJ_D^{\diff}$ includes their $A$-multiples,
so $\Phi_\star(\cJ_D^{\diff})=0$.  The embeddings of the finite
state bundles preserve currents and their residue actions, hence
$\cJ_D^{\diff}\subseteq\cJ_{D+1}^{\diff}$, and each ideal is
$\nabla$-stable by the Leibniz rule.  Every element of
$\cJ_\infty^{\diff}=\bigcup_D\cJ_D^{\diff}$ belongs to \emph{one}
such finite level; therefore its image is zero.  Since the finite
root pullback of the selected collision coefficient field is
injective, no vanishing established only in a nonfaithful quotient
of coefficient rings has been used.

The preceding equations first give the unambiguous algebraic
factorization $\cQ_\infty\to\mathscr R_\star$ and its formal
horizontality.  To make the completion assertion precise without
using the later Theorem~\ref{thm:fullstabilization}, work at the
chosen collision-boundary generic point after the same finite
ramification.  By the prior Theorems~\ref{thm:finiteconn} and
\ref{thm:canonicalregular}, the finite quotient $\cQ$ has a
finite-rank regular-singular $O$-lattice $\Lambda$ in its generic
fiber.  Take its image $\Lambda_\infty$ in the generic quotient
$\cQ_\infty\otimes_AK$.  It is a finite, torsion-free
$O$-module (hence free), and is stable under the induced
$t\partial_t$ connection; this construction uses only the
algebraic quotient, not any stabilization theorem.
Choose finite $O$-generators of $\Lambda_\infty$ and lifts in
$\Lambda$. Their original formal coefficient functions, evaluated
by $\Phi_\star$, each have finite lower $t$-valuation after a
fixed initial generalized-power factor is extracted.  A single
lower bound works for these finitely many lifts.  Thus the induced
$O$-linear functional is continuous from the $t$-adic lattice
topology to the lower-bounded complete formal coefficient module.
It extends uniquely to $\widehat{\Lambda_\infty}$, or equivalently
to the collision-root completion defined by this lattice.
Its derivative equation proves horizontality there.  This is not
an assertion of finite local freeness of $\cQ_\infty$ on all of $S$;
that assertion belongs to Theorem~\ref{thm:fullstabilization}.
\end{proof}

\begin{lemma}\label{lem:E0formalframedtransport}
Let \(c\) be a specified rooted/logged comparison from an actual smooth
three-input component to a preferred E0 representative, \emph{within the same
positive lower/direct based-path homotopy class}.  It consists of the
relative-endpoint homotopy of the based paths, all marked-root and
output-root comparisons, and the resulting typed deck pullbacks and
special-coordinate changes.  In particular it has no nontrivial pure-braid
analytic continuation; such continuation is governed by the flat
connection, not by a finite state-frame gauge.  After a finite Kummer coefficient-field
extension, let \(J_c\) act on every input by its ordered deck and special
coordinate action, and on the restricted output dual by the inverse
transpose.  Let \(j_c\) be the induced comparison of invariant finite-state
meromorphic currents.  For each \(D\ge D_0\), each finite-pole current
\(\sigma\) and algebraic external tensor \(w\),
\begin{equation}\label{eq:E0formaltransport}
 J_cR_D^{\rm act}(\sigma)w
   =R_D^{\rm ref}(j_c\sigma)J_cw,
 \qquad
 (\Phi_{\rm ref}\circ J_c)(R_D^{\rm act}(\sigma)w)=0.
\end{equation}
The second identity holds for the original product or iterate E0
formal functional, on each finite relative-weight coefficient, including
all formal derivative Ward relations.  The comparison can change a module
\(M\) into \(T_dM\): the arrows in the first formula are typed, and no
isomorphism \(M\cong T_dM\) is presumed.  On successive comparisons
\(J_{c_2c_1}=J_{c_2}J_{c_1}\), with the intermediate module types retained.
\end{lemma}
\begin{proof}
In a centered root parameter the current action is the normalized residue
\eqref{eq:canonicalmarkedres}.  For every local state \(b\),
\eqref{eq:relativefieldcov} intertwines its complete field, including the
positive Virasoro modes \(L(j)\), and the current differential is changed by
its actual Jacobian.  Thus the two normalized residues are conjugate,
coefficientwise, proving the first identity at one marking.  For a deck
comparison \(d\), put
\(Y_{T_dM}(b,z)=Y_M(d^{-1}b,z)\); its canonical underlying-vector-space
arrow intertwines \(Y_M(b,z)\) with \(Y_{T_dM}(db,z)\).
Summing over the marked orbits gives the first identity on the tensor
product.  At infinity the dual arrow is \((J_D^{-1})^*\), so the same
calculation includes the full opposite and \(L(1)\) terms, without
assuming the output vector fixed.  Coordinate changes preserving a marked
point act on states through \(L(j)\), \(j\ge0\), and deck pullback preserves
conformal weight.  Hence \(V_{\le D}\) is preserved, as are finite pole
orders and the finite truncations in the residue; each equality lies in a
finite Kummer extension of the original Laurent coefficient module.

Theorem~\ref{thm:fullsmallward} says
\(\Phi_{\rm ref}R_D^{\rm ref}=0\) for every \(D\), without convergence.
Compose this equality with the first one.  For derivatives use the
coefficient-field gauge relation
\(\nabla^{\rm act}=J_c^{-1}\nabla^{\rm ref}J_c\), obtained by differentiating
the same coordinate/current identities and the complete frame \(J_c\);
the derivative of \(J_c\) must not be discarded.  Thus the differential
Ward ideal is carried into the differential Ward ideal.  Finally path
lifting, root-coordinate substitution and the strict definition of deck
pullback respect ordered composition.  Consequently the induced arrows
obey the stated composition law.  Neither analytic continuation nor a
five-input Ward identity is used.
\end{proof}

\begin{theorem}\label{thm:fullstabilization}
There exists $D_*\ge D_0$ on the fixed affine chart such that
\begin{equation}\label{eq:fullstabilization62}
 \cJ_\infty^{\diff}=\cJ_{D_*}^{\diff},\qquad
 \cQ_\infty\cong\cQ/\mathcal K,
 \quad\mathcal K:=\cJ_{D_*}^{\diff}/\cJ_E^{\diff}.
\end{equation}
The module $\cQ_\infty$ is finite locally free on $S$, and its induced integrable connection is regular singular in the sense of Theorem~\ref{thm:canonicalregular}. No uniform cutoff $D_*$ is asserted across different affine charts.
\end{theorem}
\begin{proof}
Let $\pi:\cT\twoheadrightarrow\cQ$ be the quotient map and put
\[
 K_D=\pi(\cJ_D^{\diff})
     =\cJ_D^{\diff}/\cJ_E^{\diff}\subseteq\cQ,
 \qquad K=\bigcup_{D\ge D_0}K_D.
\]
The inclusions $E\subseteq V_{\le D}$ and the compatible current-bundle
inclusions imply $\cJ_E^{\diff}\subseteq\cJ_D^{\diff}
\subseteq\cJ_{D+1}^{\diff}$, so all displayed quotients and inverse
images are well defined.  The smooth affine coordinate ring $A$ is
Noetherian and $\cQ$ is finite over $A$ by
Theorem~\ref{thm:finiteconn}.  Thus its submodule $K$ is finitely
generated. Choose $K=\sum_{i=1}^{r}Ak_i$ with $k_i\in K_{D_i}$, and
take $D_*\ge\max(D_0,D_1,\ldots,D_r)$ (take $D_*=D_0$ if $K=0$).
Since the $K_D$ are increasing, $K=K_{D_*}$, and for every $D$,
\[
 \pi^{-1}(K_D)=\cJ_D^{\diff}:
\]
indeed, if $\pi(x)=\pi(j)$ with $j\in\cJ_D^{\diff}$, then
$x-j\in\ker\pi=\cJ_E^{\diff}\subseteq\cJ_D^{\diff}$.
Taking inverse images gives \eqref{eq:fullstabilization62}.
This $D_*$ is chosen only on the fixed affine chart.

Each $\cJ_D^{\diff}$ is stable under all placement derivations by
its definition and the Leibniz rule. Hence $K$ is stable under the
induced connection on $\cQ$, and the formulas
\[
 \nabla^\infty_\delta[q]=[\overline\nabla_\delta q],\qquad
 \nabla^\infty_\delta(f[q])
   =\delta(f)[q]+f\nabla^\infty_\delta[q],\qquad
 [\nabla^\infty_\delta,\nabla^\infty_\gamma]
   =\nabla^\infty_{[\delta,\gamma]}
\]
are well defined.  Thus $\cQ_\infty=\cQ/K$ is a finite coherent
$A$-module with integrable algebraic connection.  Since $S$ is smooth
over $\C$, \cite[Proposition~8.8]{Katz}, or
Lemma~\ref{lem:connfree}, makes it locally free, including at points
where a priori it might have torsion.

\Needspace{6\baselineskip}
For regularity, work on one irreducible component of $S$ and write
$F=\C(S)$.  Let $\overline S$ be the \emph{fixed normal proper}
compactification of Theorem~\ref{thm:canonicalregular}.  For any
codimension-one point $\eta$ of $\overline S\setminus S$ put
$O_\eta=\widehat{\mathcal O}_{\overline S,\eta}$ and
$K_\eta=\operatorname{Frac}(O_\eta)$; this is a complete DVR because
$\overline S$ is normal.  The generic quotient
\[
 0\longrightarrow K\otimes_A F
 \longrightarrow\cQ\otimes_A F
 \overset{\pi_F}{\longrightarrow}\cQ_\infty\otimes_A F
 \longrightarrow0
\]
is exact by flat localization.  It is an exact sequence of differential
modules: stability under $\operatorname{Der}_{\C}(A)$ extends by
localization to rational derivations of $F$, and the quotient map is
horizontal.  Extending scalars along $F\hookrightarrow K_\eta$ gives
the horizontal surjection
\[
 \pi_\eta:Q_\eta:=\cQ\otimes_AK_\eta
 \twoheadrightarrow
 Q_{\infty,\eta}:=\cQ_\infty\otimes_AK_\eta.
\]
In particular, its kernel is stable under the extended logarithmic
derivation $\delta_\eta=t\partial_t$, with the chosen coefficient
field held constant. This is a statement about the generic
differential quotient; it does not assert an exact base change for a
reflexive integral extension of a current bundle. The algebraic
equality $\cJ_\infty^{\diff}=\cJ_{D_*}^{\diff}$ also remains true
under this flat field extension.

By Theorem~\ref{thm:canonicalregular}, the original differential
module $Q_\eta$ has a full $O_\eta$-lattice $\Lambda_\eta$ with
$\nabla_{\delta_\eta}\Lambda_\eta\subseteq\Lambda_\eta$.
Define its image \emph{inside the quotient vector space}
\[
 \Lambda_{\infty,\eta}:=\pi_\eta(\Lambda_\eta)
 \subseteq Q_{\infty,\eta}.
\]
It is generated over $O_\eta$ by the images of finitely many
generators of $\Lambda_\eta$.  As a submodule of the
$K_\eta$-vector space $Q_{\infty,\eta}$ it is torsion-free; this
does not require that $\pi_\eta$ extend as a surjection between
previously chosen integral models.  Since
$K_\eta\Lambda_\eta=Q_\eta$ and $\pi_\eta$ is surjective,
$K_\eta\Lambda_{\infty,\eta}=Q_{\infty,\eta}$.
A finite torsion-free module over the DVR $O_\eta$ is free, so
\[
 \operatorname{rank}_{O_\eta}\Lambda_{\infty,\eta}
 =\dim_{K_\eta}Q_{\infty,\eta}.
\]
The connection on the quotient is well defined for $\delta_\eta$
because its kernel is stable, and horizontality gives, for every
$\lambda\in\Lambda_\eta$,
\[
 \nabla^\infty_{\delta_\eta}\bigl(\pi_\eta(\lambda)\bigr)
 =\pi_\eta\bigl(\nabla_{\delta_\eta}\lambda\bigr)
 \in\Lambda_{\infty,\eta}.
\]
The Leibniz rule extends this inclusion to the full image lattice;
thus $\Lambda_{\infty,\eta}$ is a full logarithmically stable DVR
lattice, also when the quotient has rank zero.  If a lattice for
$Q_\eta$ was first obtained after a finite Kummer extension, the
DVR-descent step of Lemma~\ref{lem:DVRdescent}, already used in
Theorem~\ref{thm:canonicalregular}, supplies the lattice over this
$O_\eta$ before taking its image.  No descent of an arbitrarily
chosen integral quotient is assumed.

We have established regularity at each boundary prime of the
same fixed normal proper compactification.  The hypotheses of
\cite[Theorem~1.1.1]{Andre} now hold: $S$ is smooth in characteristic
zero and $(\cQ_\infty,\nabla^\infty)$ is a finite-rank algebraic
integrable connection. Its codimension-one criterion therefore gives
regular singularity globally, including along exceptional divisorial
valuations on other models; a separate uniform cutoff $D_*$ at such
valuations is unnecessary.

Finally, for the finite projective module $\cQ_\infty$ and every
$A$-module $\mathscr R_\star$, evaluation gives
\[
 \cQ_\infty^\vee\otimes_A\mathscr R_\star
 \overset{\sim}{\longrightarrow}
 \operatorname{Hom}_A(\cQ_\infty,\mathscr R_\star),
 \qquad f\otimes r\longmapsto(q\mapsto f(q)r).
\]
Theorem~\ref{thm:fullsmallward} already constructs the original
formal functionals on the algebraic quotient (and its specified
collision-root completion) and proves their derivative identity;
this identification merely places those same functionals in the
dual connection.  It does not appeal to the later analytic
continuation or convergence theorems.
\end{proof}

\begin{theorem}\label{thm:E0covariantcontinuation}
For either original product or iterate from Definition~\ref{def:commonoutsetframing}, the normally convergent function of Corollary~\ref{cor:smallconvergenceautomatic} has a unique continuation as a horizontal analytic functional along paths in the framed universal cover of its collision-free component.  On an arbitrary resulting framed chart its \emph{covariant local expression}, with transported marked roots, based lifts, deck-module frames, logarithms and the inverse-transpose output restricted dual, annihilates the local seed Ward ideal and the full-state ideal $\cJ_\infty^{\diff}$.  The original untransported grade series is asserted to converge only in its own product or iterate chamber.  A closed pure braid is allowed to act by nontrivial monodromy, possibly mixing channel solutions.
\end{theorem}
\begin{proof}
We distinguish the preferred series, its analytic germ, and the analytic
continuation of that germ.  Write
$\cT_\alpha=\bigotimes_i W_{i,\alpha}$, with the output factor equal
to the \emph{restricted} contragredient $D_\alpha'$, and let
$\mathscr S_{D,\alpha}$ denote the actual global, finite-pole,
$H$-invariant $V_{\le D}$-currents on the corresponding framed cover.
All assertions about an infinite external tensor below are first
assertions on its algebraic direct sum, on which a specified vector
has only finitely many homogeneous components.

Choose a simply connected relatively compact subchamber $U$ in the
original product or iterate convergence domain, with the outset
based roots and logarithms.  For each algebraic external tensor $w$,
Corollary~\ref{cor:smallconvergenceautomatic} gives normal convergence
of the \emph{fixed original} IO contraction $\Phi_\star(w)$ on $U$.
For a fixed $D$, a finite-pole current $s$, and such $w$,
Proposition~\ref{prop:canonicalrealization} makes $R_D(s)w$ a
finite sum of algebraic external tensors.  Therefore
$\Phi_\star(R_D(s)w)$ is an analytic function on $U$; its formal
coefficient series is zero by Theorem~\ref{thm:fullsmallward},
including the common root substitution used there.  Analytic
uniqueness of the convergent expansion gives
$\Phi_\star(R_D(s)w)=0$ on $U$.  The same argument applies to
every generator of $\cJ_D^{\diff}$: the finite iterated
$L(-1)$-derivative identities of Lemma~\ref{lem:channelderivative}
are identities of normally convergent analytic functions, and
termwise derivatives converge on smaller compact subsets by Cauchy's
formula.  Multiplication by a regular chart coefficient and
Leibniz differentiation introduce no infinite sum at any fixed
algebraic input.  Hence this particular analytic germ kills
$\cJ_\infty^{\diff}=\bigcup_D\cJ_D^{\diff}$ and satisfies
$\delta\Phi_\star(w)=\Phi_\star(\nabla_\delta w)$ for every
placement derivation.  It thus descends to a horizontal analytic
functional on the already finite locally free $\cQ_\infty|_U$;
locally a finite set of representatives of a frame of $\cQ_\infty$
suffices to verify holomorphy.  In particular, it is the original
IO germ, not an arbitrarily selected solution of the Ward system.
No analytic statement about the original grade series outside $U$
is used.

On an overlap of two genuinely identified framed Hurwitz charts
$\alpha,\beta$, denote by
$j_{\alpha\beta,D}:\mathscr S_{D,\beta}\xrightarrow{\sim}
\mathscr S_{D,\alpha}$ the expression of the \emph{same global
current} under the cover, root and special-coordinate comparison.
The allowed pole orbits are transported, not changed; finite state
jets preserve $V_{\le D}$ since $L(j)V_{\le D}\subseteq V_{\le D}$
for $j\ge0$.  At an input marking the typed arrow
$C_{i,\alpha\beta}:W_{i,\beta}\to W_{i,\alpha}$ is obtained by
first taking the strict deck pullback $W_{i,\beta}\to
T_{d_i}W_{i,\beta}$ and then applying the relative-weight
special-coordinate operator $R_{T_{d_i}W_{i,\beta},h_i}(\rho_i)$;
its target is the module named by the actual $\alpha$ framing.
The order is part of this formula: it does not identify distinct
twisting sectors or commute two deck elements.

At the output the local physical parameter is $s=z^{-1}$, and the
assigned module is $D'$, with the inverse inertia.  Define
\emph{on the restricted dual itself} the typed comparison
$C_{\infty,\alpha\beta}:D_\beta'\to D_\alpha'$ by the same
strict pullback followed by
$R_{D_\alpha',h_{D'}}(\rho_\infty)$, using its own output
special-coordinate parameter.  Every positive Virasoro operator
$L_{D'}(j)$, $j>0$, lowers relative grade, so this and its inverse
act by finite sums on each algebraic vector of $D'$.  This is
essential: the unrestricted transpose of an origin-fixing
positive-mode exponential on $D$ generally takes $D'$ into its
\emph{completion}, not into $D'$ itself.  To express the correct
output arrow by inverse transpose, let $\widehat D$ be the full
graded completion paired with $D'$ and put
$G_D:=(C_{\infty,\alpha\beta}^{-1})^t$ on that completed
pairing space.  Then, as a well-typed map on restricted duals,
\[
 C_{\infty,\alpha\beta}=(G_D^{-1})^t,
 \qquad
 \langle C_{\infty,\alpha\beta}d',G_Dv\rangle_\alpha
    =\langle d',v\rangle_\beta.
\]
The complete opposite-field formula
$Y_{D'}(b,s)=Y_D(e^{sL(1)}(-s^{-2})^{L(0)}b,s^{-1})^t$
identifies this output action with the original contragredient
residue convention, including all $L(1)^jb$ and the derivative
Jacobian.  It does \emph{not} identify
$G_D$ with $R_{D,h_D}(\rho_\infty)$ on the uncompleted primal.
Put $B_{\alpha\beta}=\bigotimes_i C_{i,\alpha\beta}$, including
this output arrow.  Both $B_{\alpha\beta}$ and its inverse, and
their derivatives in chart parameters, are well defined on
$\cT_\beta$ coefficientwise in every algebraic vector.

For each $s_\beta\in\mathscr S_{D,\beta}$ the centered
identity \eqref{eq:relativefieldcov} transports the full state
field and its actual Jacobian.  At the output apply that identity
\emph{to $D'$}, and the opposite-field formula just displayed.
For a root rotation $u\mapsto\zeta u$ the input factor is
$\zeta^{-N(L(0)-h)}$ and the fractional mode eigencharacter is
its reciprocal in the normalized root residue; at the output use
the inverse-inertia action on $D'$.  Deck changes satisfy
$Y_{T_dW}(db,x)=Y_W(b,x)$, and their chronological products are
typed.  Thus Proposition~\ref{prop:canonicalrealization} gives
for the \emph{same actual global current}
\begin{equation}\label{eq:E0actualgauge}
 B_{\alpha\beta}R_\beta(s_\beta)
   =R_\alpha(j_{\alpha\beta,D}s_\beta)B_{\alpha\beta},
 \qquad
 B_{\alpha\beta}\cJ_{D,\beta}=\cJ_{D,\alpha}.
\end{equation}
The reverse inclusion uses the inverse current and frame arrows;
in particular a local residue calculation alone would not imply
this equality without the global-current comparison.

In a local column trivialization write $\nabla_\alpha=d+A_\alpha$.
A section with $\alpha$-coordinates $B_{\alpha\beta}w_\beta$
has $\beta$-connection
\begin{equation}\label{eq:E0differentialgauge}
 \nabla_\beta=B_{\alpha\beta}^{-1}\nabla_\alpha B_{\alpha\beta},
 \qquad
 A_\beta=B_{\alpha\beta}^{-1}A_\alpha B_{\alpha\beta}
   +B_{\alpha\beta}^{-1}dB_{\alpha\beta}.
\end{equation}
Here the placement term $\sum_i dz_i L_i(-1)$ is the connection
in the original fixed physical frames; a variable special-coordinate
frame contributes the displayed derivative of $B$ and is not
silently given the same fixed-frame matrix.  Equation
\eqref{eq:E0differentialgauge} implies for every derivation and
$j\in\cJ_{D,\beta}$ that
$B\nabla_{\beta,\delta}j=\nabla_{\alpha,\delta}(Bj)$.
Using \eqref{eq:E0actualgauge}, the Leibniz rule and induction on
the number of derivatives, we therefore obtain
\begin{equation}\label{eq:E0idealframe}
 B_{\alpha\beta}\cJ_{D,\beta}^{\diff}
  =\cJ_{D,\alpha}^{\diff},\qquad
 B_{\alpha\beta}\cJ_{\infty,\beta}^{\diff}
  =\cJ_{\infty,\alpha}^{\diff}.
\end{equation}
No common integer $D_*$ over different charts is required.
These are equalities after restriction to the overlap and, if the
coordinate gauges are only analytic, after analytification.
Consequently $B$ descends to an isomorphism of the full-state
quotients with their connections.

For three successive actual comparisons, the ring-map convention
for special coordinates, their finite-jet actions and strict
$T_{d_1}T_{d_2}=T_{d_1d_2}$ give
$C_{i,\alpha\gamma}=C_{i,\alpha\beta}C_{i,\beta\gamma}$
with the intermediate twisted sector retained.  The same follows
at infinity from the \emph{restricted-dual} operators and the
contragredient pairing, hence
$B_{\alpha\gamma}=B_{\alpha\beta}B_{\beta\gamma}$.
This is a \emph{coordinate-gauge cocycle}; it is not analytic
parallel transport along a path.  The framed atlas therefore
glues the finite locally free integrable bundle
$(\cQ_\infty,\nabla^\infty)$ and its analytic dual.  By local
existence and uniqueness for holomorphic integrable connections,
the initial \emph{specific} analytic horizontal germ extends along
every collision-free framed path.  It depends only on the
endpoint-fixed homotopy class and is single valued on the framed
universal cover.  In an endpoint chart its coordinates are
$F_\beta(w)=F_\alpha(B_{\alpha\beta}w)$ for a pure overlap
change, with the horizontally continued $F_\alpha$ used when a
path is involved.  Equations \eqref{eq:E0idealframe} and
\eqref{eq:E0differentialgauge} give
\[
 F_\beta(\cJ_{\infty,\beta}^{\diff})=0,
 \qquad
 \delta F_\beta(w)=F_\beta(\nabla_{\beta,\delta}w).
\]
For a closed pure braid, the endpoint coordinate gauge may be
trivial while the flat parallel transport is not; its monodromy
acts on the horizontal solution space and preserves the full-state
Ward quotient.  In particular there is no claim that the
untransported original grade sum converges outside its initial
chamber or that a pure braid has trivial monodromy.
\end{proof}

\section{Node reconstruction and actual iterate factorization}\label{sec:reconstruction}
\subsection{The three-point Ward converse}
Fix simple ordinary modules $X\in\cC_g$, $Y\in\cC_h$ and $Z\in\cC_{gh}$, and the based three-point $H$-cover of $(\mathbb P^1;\infty,1,0)$, where $H=\langle g,h\rangle$. Put $T_3=Z'\otimes X\otimes Y$. For $D\geq0$ let
$\cJ_{D,3}\subseteq T_3$ be generated by the normalized residues of $H$-invariant meromorphic currents with states in $V_{\leq D}$ and poles confined to the three marked orbits. Define
\[
 \cJ_{\infty,3}=\bigcup_{D\geq0}\cJ_{D,3},\qquad
 \mathfrak W_{XY}^{Z}=(T_3/\cJ_{\infty,3})^*.
\]
Choose $D_0$ with $E\subseteq V_{\leq D_0}$. The three-point seed relations of Proposition~\ref{prop:seedfinite} are then contained in $\cJ_{D_0,3}$, so that proposition implies $\dim_\C\mathfrak W_{XY}^{Z}<\infty$.
Here the residue normalization and the based-root identification are those of Proposition~\ref{prop:canonicalrealization}; differential saturation is not imposed.

\begin{definition}\label{def:singleinsertionoperators}
For ordinary $X,Y,Z$ define $\widetilde{\cI}\binom Z{X\ Y}$ to be the space of maps $\mathcal Y:X\otimes Y\to Z\{z\}[\log z]$ satisfying items (i), (ii) and (iv) of Definition~\ref{def:completeintertwiner}, with the same preferred
branches $\mathcal A_i,\mathcal D_i$, and $\pi^Z_\gamma\mathcal Y_{n,j}(X_\alpha)Y_\beta=0\; \text{if }\gamma\ne\alpha+\beta-n-1$. No condition (iii) of that definition is imposed. The $L(0)$-commutator calculation in Lemma~\ref{lem:IOlogs}
uses only the one-$V$ duality axiom (with $u=\omega$); hence $\mathcal Y_{n,j}=0$ for $j>0$ in this ordinary setting. Consequently
\[
 \cI\binom Z{X\ Y}\subseteq
 \widetilde{\cI}\binom Z{X\ Y}.
\]
\end{definition}

\begin{proposition}\label{prop:singlewardconverse}
Evaluation on the specified branch of $z=1$ induces an isomorphism $\operatorname{ev}_1:\widetilde{\cI}\binom Z{X\ Y} \xrightarrow{\ \sim\ }\mathfrak W_{XY}^{Z}$. For $\ell\in\mathfrak W_{XY}^{Z}$, $x\in X_\alpha$,
$y\in Y_\beta$ and $z'\in (Z')_\gamma$, the inverse is
\begin{equation}\label{eq:singlewardinverse}
 \big\langle z',\mathcal Y_\ell(x,z)y\big\rangle
 =e^{(\gamma-\alpha-\beta)\log z}\ell(z',x,y),
\end{equation}
where $\log1=0$ and the logarithm is continued along the prescribed based paths.
\end{proposition}
\begin{proof}
For homogeneous $x\in X_\alpha$, $y\in Y_\beta$ and each
$\gamma\in\wt(Z)$, finite dimensionality of $Z_\gamma$ gives a unique
$v_\gamma(x,y)\in Z_\gamma$ such that
$\langle z',v_\gamma(x,y)\rangle=\ell(z',x,y)$ for every
$z'\in Z_\gamma^*$.  Set
\[
 \mathcal Y_\ell(x,z)y
   =\sum_{\gamma\in\wt(Z)}v_\gamma(x,y)
                   z^{\gamma-\alpha-\beta}.
\]
Only the coefficient in $Z_\gamma$ is used in each subsequent
pairing: evaluation at $z=1$ does not mean summing the completed
output in $Z$.  Since $Z$ is ordinary and grading-restricted, its
weights are bounded below and lie in finitely many cosets of
$N_Z^{-1}\mathbb Z$; hence the coefficient with
$n=\alpha+\beta-\gamma-1$ vanishes for sufficiently large $n$ in
each allowed coset.  This is the required lower truncation, and the
same coefficient formula is linear on algebraic $X\otimes Y$.  No
logarithmic coefficient has been introduced.  We establish the one-$V$ analytic axiom directly.

Fix homogeneous $z'\in(Z')_\gamma$, $x\in X_\alpha$,
$y\in Y_\beta$, a cutoff $D$, and $u\in V_{\le D}$.
On the specified based $H$-cover of $(\mathbb P^1;\infty,1,0)$
consider the following \emph{formal, coefficientwise} families:
\[
\begin{aligned}
 F_\infty(w;u)&=\sum_{\delta\in\wt(Z)}
  \langle z',Y_Z(u,w)v_\delta(x,y)\rangle,\\
 F_0(w;u)&=\ell(z',x,Y_Y(u,w)y),\\
 F_1(w;u)&=\ell(z',Y_X(u,w-1)x,y).
\end{aligned}
\]
The notation for the second and third families means application of
$\ell$ to each coefficient, not evaluation of a convergent series.
For $u\in V_d$ the coefficient of $v_\delta$ in the first family
has the single possible mode $u_n$ with
$n=\delta+d-\gamma-1$, and therefore its $w$-power is
$w^{\gamma-\delta-d}$.  Thus the apparently infinite completed
output is never acted upon as a whole: every coefficient in the
local parameter $s=w^{-1}$ is a finite-dimensional matrix
coefficient, and the lower bound on $\delta$ gives a finite negative
part in $s$.  More invariantly, the output local expression is
written with $Y_{Z'}$ and the opposite-vertex coordinate
transformation
\[
 e^{sL(1)}(-s^{-2})^{L(0)}u
  =(-1)^d\sum_{j=0}^{d}\frac{s^{j-2d}}{j!}L(1)^ju.
\]
This is a finite sum, and lower truncation on the algebraic vector
$z'$ gives the same finite pole bound.  At $0$ and $1$, lower
truncation of $Y_Y(u,w)y$ and $Y_X(u,w-1)x$ gives finite negative
parts in their respective root parameters.  A homogeneous basis of
the finite-dimensional $V_{\le D}$ has a common bound.  All
coordinate/root gauges and their inverses act on this truncation
through finite matrices of Laurent series with finite pole order;
therefore they preserve the finite-negative-part property.

These assignments, for every $u\in V_{\le D}$, are precisely a
collection $\tau_{D,p}$ of \emph{complete} formal local sections of
$\mathscr V_D^\vee$ at all marked lifts $p\in Q$.  At the output,
the restricted contragredient and its full $L(1)$ correction define
the infinity component; the inertia congruences of the three
original twisted fields and the chosen roots define all other
lifts by $\tau_{D,hp}=h\tau_{D,p}$.  This uses neither a Ward
converse nor convergence of the reconstructed intertwining operator.

The theorem in \cite[Theorem~6.1]{FS} must here be applied to the
finite-rank bundle
$\mathscr E_D=\mathscr V_D\otimes\omega_C$:
its local dual one-form
$\mathscr E_D^\vee\otimes\omega_C$ is exactly
$\mathscr V_D^\vee$, while the global test sections are
$H^0(C\setminus Q,\mathscr E_D)$.  In particular the theorem
compares all Laurent coefficients, not merely principal
parts.  For an invariant global test current $a$, the actual
orbit-normalized residue formula of
Proposition~\ref{prop:canonicalrealization} yields
\[
 \sum_{p\in Q}\operatorname{Res}_p
       \langle\tau_{D,p},a_p\rangle
   =|H|\,\ell\bigl(R_D(a)(z'\otimes x\otimes y)\bigr)=0.
\]
Indeed, if $N_i$ is the inertia order over marking $i$, there are
$|H|/N_i$ lifts, each with residue $N_i$ times its normalized
marked action; neither fractional modes nor the infinity leg
produce an additional factor.  If $a$ is not invariant, the
$H$-equivariance of the entire local collection implies that
the total residue pairing with $a$ equals that with
$|H|^{-1}\sum_{h\in H}ha$, so it is again zero.  This verifies
the strong residue criterion for every required global test
section, not only for invariant ones.

The finite-rank Strong Residue Theorem now supplies a section
\[
 \lambda_D\in H^0(C\setminus Q,\mathscr V_D^\vee)
\]
whose complete germs are exactly $\tau_{D,p}$.  More explicitly,
fix a common pole bound $r$ for these formal germs.  Requiring a
section of the finite-dimensional space
$H^0(C,\mathscr V_D^\vee(rQ))$ to match their coefficients
through any prescribed finite jet is an affine linear problem.
The Serre-dual obstruction is a residue pairing against a global
section of $\mathscr V_D\otimes\omega_C(*Q)$, and is zero by
the criterion above.  Hence every finite-jet problem is solvable;
the resulting nested nonempty affine subspaces of this one finite-
dimensional space stabilize and have a common point.  This proves
that all Taylor coefficients, not merely principal parts, agree.
Its poles at $Q$ are finite by the already established negative-power
bounds; it has no poles away from $Q$.  This use of the theorem does not infer
an analytic germ directly from an arbitrary formal power series:
existence of the algebraic global section forces the specified
formal series to be its convergent local expansions.  Uniqueness
follows since $C$ is connected and a meromorphic section with zero
complete germ at one marked lift vanishes identically.  Thus the
$H$-average of $\lambda_D$ equals $\lambda_D$ itself.  For
$D\le D'$, the finite-jet state-bundle inclusions commute with
the prescribed field germs; restriction of $\lambda_{D'}$ has the
same complete germs as $\lambda_D$, so uniqueness proves
$\lambda_{D'}|_{\mathscr V_D}=\lambda_D$.  No infinite-rank
strong residue theorem is being used.

The chosen finite $H$-cover has no branch points away from
$0,1,\infty$.  Each $\lambda_D$ is meromorphic on that cover with
poles confined to the lifts of these three points.  Consequently
its three marked Laurent expansions converge absolutely and
locally normally in the root-coordinate domains lying over
$|w|>1$, $|w|<1$, and $|w-1|<1$, respectively, with the analytic
branches selected by the original based paths.  They are the
three formal series above, coefficient for coefficient.  Thus the
sums over projected intermediate weights appearing in
Definition~\ref{def:completeintertwiner}(ii) converge absolutely
in its $\mathcal A_1,\mathcal A_2,\mathcal A_3$ and represent the
analytic continuation of the same meromorphic section on
its selected rooted cover.  For arbitrary $z_2\ne0$ substitute
$w=z_1/z_2$; the $L(0)$-homogeneity of each of the three module
fields gives the same factor
$\exp((\gamma-\alpha-\beta-d)\log z_2)$.
Its logarithm is continued from $\log1=0$ along the fixed based
paths, and the root/deck comparisons are exactly those used to
construct $\tau_D$.  The resulting multivalued function on
$\mathcal M_2$ therefore has the preferred restrictions on the
original $\mathcal D_i$.  This proves the complete one-$V$ axiom,
including the required analytic convergence, without using the
multiple-$V$ axiom (iii).

Use the invariant global Virasoro current corresponding to the
Möbius vector field $w\partial_w$.  In local physical coordinates
$\xi=w-1$ at $1$, $w$ at $0$, and $s=w^{-1}$ at infinity, its
vector fields are
$(1+\xi)\partial_\xi$, $w\partial_w$, and
$-s\partial_s$.  Their canonical marked residues, including the
restricted-dual opposite action at infinity, give the single
Ward identity
\[
 0=\ell(z',L_X(-1)x,y)
   +(\alpha+\beta-\gamma)\ell(z',x,y).
\]
The changes of physical coordinates here are Möbius, so their
Schwarzian derivatives vanish; the root-coordinate factors are
already included in the normalized residue action and do not
supply a further scalar.  Hence
\begin{equation}\label{eq:singlewardLminusone}
 \ell(z',L_X(-1)x,y)
   =(\gamma-\alpha-\beta)\ell(z',x,y).
\end{equation}
Since $L_X(-1)x$ has weight $\alpha+1$, differentiating the
reconstructed monomial gives
\[
 \partial_z\langle z',\mathcal Y_\ell(x,z)y\rangle
  =(\gamma-\alpha-\beta)z^{\gamma-\alpha-\beta-1}
       \ell(z',x,y)
  =\langle z',\mathcal Y_\ell(L_X(-1)x,z)y\rangle.
\]
The restricted dual separates the output weight spaces, so the
formal derivative axiom holds.  All defining requirements of
$\widetilde{\cI}\binom Z{X\ Y}$ have now been checked.

Conversely, let $\mathcal Y$ satisfy only items (i), (ii), (iv)
and weight compatibility.  The proof of
Lemma~\ref{lem:singleinsertion} for a single inserted state
uses exactly the analytic convergence and continuation in item
(ii), the finite $L(1)$ opposite-field formula at infinity, and
finite state-coordinate covariance; it never inserts a second
independent $V$-state.  Repeating that proof for every fixed
$V_{\le D}$ supplies a meromorphic $H$-equivariant
$\mathscr V_D^\vee$ section with its three prescribed complete
germs.  Pairing it with any invariant global current and applying
the ordinary residue theorem, with the same $|H|/N_i$ orbit
multiplicities, proves
$\operatorname{ev}_1(\mathcal Y)(\cJ_{D,3})=0$
for every $D$.  Thus evaluation lies in $\mathfrak W_{XY}^Z$.
The $L(0)$-commutator calculation of
Lemma~\ref{lem:IOlogs} uses only item (ii) and ordinary grading,
so each paired coefficient is exactly a constant times
$z^{\gamma-\alpha-\beta}$ without nonconstant logarithms.
Consequently reconstruction after evaluation reproduces every
weight component of $\mathcal Y$; evaluation after reconstruction
is $\ell$ by construction.  The two maps are inverse.
\end{proof}

For the collision construction choose simple ordinary modules $A\in\cC_{g_A}$, $B\in\cC_{g_B}$, $C\in\cC_{g_C}$ and $D\in\cC_{g_Ag_Bg_C}$, and an actual product channel $F$ whose intermediate object is a finite direct sum.
A finite linear combination of actual product channels is handled linearly. For arbitrary external objects use the finite semisimple decompositions of Proposition~\ref{prop:regularmodules}; the finite
union of their exponent sets and the largest logarithm bound give common estimates. The converse just proved applies on each simple summand; the extra-field axiom is not inferred here.

Fix the positive branch and the collision chart $z_2=1,\; z_1=1+q,\; q=t^e,\; e\in\Z_{>0},\; t>0$. Set $T_{\rm in}=A\otimes B$ and $T_{\rm out}=D'\otimes C$, using the vector-space permutation only to regroup the external
factors; the operator order and the based paths are unchanged. For homogeneous $a\in A$, $b\in B$ and $o\in T_{\rm out}$ define
\begin{equation}\label{eq:collisionnormalizedkernel}
 \widehat F(t;o;a,b)=
 e^{e(\wt a+\wt b)\log t}\,
 F(1+t^e,1;o;a,b).
\end{equation}
The prefactor includes the lowest weights, which need not be rational. An iterate through an ordinary intermediate module $N$ with weights $h_N+\gamma$ ($\gamma\geq0$) has individual normalized
powers $t^{e(h_N+\gamma)}$. This last statement is a homogeneous coefficient identity as in Lemma~\ref{lem:actualonecontraction}, not an assertion of convergence for a previously unconstructed iterate.

\begin{proposition}\label{prop:collisionuniformsupport}
Let $F$ be a fixed actual product channel with the external modules and collision normalization of
\eqref{eq:collisionnormalizedkernel}. There exist a finite subset
$\Sigma\subset\C$ and $L\in\N$, depending on $F$ and the chosen collision branch but not on $o,a,b$, such that
\begin{equation}\label{eq:collisioncoefficientkernels}
 \widehat F(t;o;a,b)
 =\sum_{\sigma\in\Sigma}\sum_{n\ge0}\sum_{l=0}^L
 t^{\sigma+n}(\log t)^l
 \ell_{\sigma+n,l}(o\otimes a\otimes b),
 \qquad
 \ell_{\sigma+n,l}\in (T_{\rm out}\otimes T_{\rm in})^*.
\end{equation}
The sets $\sigma+\N$, $\sigma\in\Sigma$, are pairwise disjoint. In particular, $\inf\{\Re\beta:\ell_{\beta,l}\ne0 \text{ for some }0\le l\le L\}>-\infty$ whenever $\widehat F\ne0$. The equality in
\eqref{eq:collisioncoefficientkernels} is an equality of
branchwise Nilsson expansions; the coefficient functionals are ordinary algebraic linear functionals, with no continuity assumption on the infinite-dimensional external tensor product.
\end{proposition}
\begin{proof}
Put $O=\C[[t]]$, $K=\C((t))$ and $\delta=t\partial_t$. On homogeneous $w=o\otimes a\otimes b$, $a\in A_\alpha$, $b\in B_\beta$, set $\nu(w)=e\bigl((\alpha-h_A)+(\beta-h_B)\bigr)\in\N$.
The finite ramification $q=t^e$ is chosen so that the relative weight differences are integral. By Lemma~\ref{lem:weightednormal} and the boundary estimates of Theorem~\ref{thm:canonicalregular},
the finite-rank quotient $Q_K$ has a full $O$-lattice $\Lambda\subset Q_K$ with $\nabla_\delta\Lambda\subseteq\Lambda,\; t^{\nu(w)}[w]\in\Lambda$. Here $[w]$ denotes the image of the external algebraic tensor.
The horizontal full-state quotient of Theorem~\ref{thm:fullstabilization} (if denoted $Q_{\infty,K}$ in its construction) is a differential quotient of $Q_K$. Let $\pi:Q_K\twoheadrightarrow Q_{\infty,K}$ and put
\[
 \Lambda_\infty=\pi(\Lambda),\qquad
 \nabla^\infty_\delta\Lambda_\infty\subseteq\Lambda_\infty,
 \qquad K\Lambda_\infty=Q_{\infty,K}.
\]
The image is a finite torsion-free $O$-module, hence a full free lattice over the DVR $O$. In particular, the passage to full-state Ward relations does not require a separate uniform cutoff for the external vectors.

It is not necessary, and has not been proved, that $\Lambda_\infty$ is generated by weighted algebraic tensors. Instead put
\[
 M_\nu:=\sum_{\substack{w\in T\\w\text{ homogeneous}}}
 O\,t^{\nu(w)}\pi([w])\ \subseteq\Lambda_\infty.
\]
Since $O$ is Noetherian and $\Lambda_\infty$ is a finite free $O$-module, $M_\nu$ is a finite $O$-module. By the definition of this sum, finitely many of its generators can be chosen in the \emph{actual specified form}
$g_j=t^{\nu(w_j)}\pi([w_j])$ with homogeneous algebraic $w_j$; choose $1\le j\le s$ so that $M_\nu=\sum_j Og_j$. Every $t^{\nu(w)}\pi([w])$ belongs to this $M_\nu$, so
\begin{equation}\label{eq:collisionlatticereduction73}
 t^{\nu(w)}\pi([w])=\sum_{j=1}^s a_j(t;w)g_j,
 \qquad a_j(t;w)\in O.
\end{equation}
This is finite for each $w$, and its fixed generators are independent of $w$. No equality $M_\nu=\Lambda_\infty$ is asserted or needed.

We also check the \emph{analytic} meaning of these coefficients. On the chosen collision disk let $R=\C\{t\}$ and $\widehat R=O$. The algebraic quotient and its specified algebraic tensor classes are meromorphic in an algebraic local frame near the puncture. Thus each of the finite set $g_j$, and each fixed $t^{\nu(w)}\pi([w])$, lies in $t^{-b}R^r$ after choosing a suitable common $b$ for these finitely many vectors. The $R$-module $M_R=\sum_jRg_j\subseteq t^{-b}R^r$ has completion $M_R\otimes_R\widehat R=\sum_jOg_j$, as a submodule, by flatness of $R\to\widehat R$. Therefore the membership in \eqref{eq:collisionlatticereduction73}, checked over $O$, implies membership already in $M_R$: apply faithful flatness to the image of $t^{\nu(w)}\pi([w])$ in $t^{-b}R^r/M_R$. In particular the $a_j(t;w)$ may be chosen in $\C\{t\}$, not merely in $\C[[t]]$. The number $b$ and these coefficient functions may depend on $w$; their \emph{pole-free} character and the finite set of $g_j$ do not. This distinction will give genuinely convergent, rather than only formal, Nilsson expansions.

By Theorem~\ref{thm:fullsmallward} on the preferred collision completion, together with Theorem~\ref{thm:E0covariantcontinuation} for the actual analytic germ at the present framed placement, $F$ annihilates the full-state Ward relations and is horizontal. Consequently it is a horizontal functional on $Q_{\infty,K}$. The regular-singular lattice $\Lambda_\infty$ and the
regular-singular local theory
\cite[Chapter~II, Proposition~5.4]{Deligne} give, for the
finitely many $g_j$, expansions of the form
\[
 F(g_j)=\sum_{\lambda\in\Xi_j}\sum_{l=0}^{L_j}
 t^\lambda(\log t)^l f_{j,\lambda,l}(t),\qquad
 f_{j,\lambda,l}(t)\in\C[[t]],
 \quad |\Xi_j|<\infty.
\]
These are analytic regular-singular (Frobenius--Nilsson) expansions on the fixed logarithmic branch, with the power series $f_{j,\lambda,l}$ convergent near $t=0$; convergence follows from the analytic regular-singular connection, not from a formal $O$-lattice alone. If a generator vanishes under $F$, its expansion is zero. Set $\Xi=\bigcup_j\Xi_j$ and $L=\max_jL_j$ (take $L=0$ when all vanish). Equation~\eqref{eq:collisionlatticereduction73} with the \emph{convergent} coefficient functions just proved produces only nonnegative integer exponent shifts, never increases $L$, and gives locally normally convergent expansions for each fixed algebraic $w$. Finally the normalization in
\eqref{eq:collisionnormalizedkernel} is exactly
\begin{equation}\label{eq:collisionnormalfactor73}
 \widehat F(t;w)
 =t^{e(h_A+h_B)}F\bigl(t^{\nu(w)}\pi([w])\bigr),
\end{equation}
where both powers use the prescribed branch. Hence the possible exponents belong to the finite union $e(h_A+h_B)+\Xi+\N$.

Partition this finite union by $\C/\Z$. Within each congruence class the finite set of initial exponents is totally ordered by integer difference; choose its least element $\sigma$. Initial exponents in other classes have disjoint progressions, and the union of all progressions in the given class is contained in $\sigma+\N$. Taking one $\sigma$ per class produces a finite set $\Sigma$ whose progressions $\sigma+\N$ are pairwise disjoint. In particular, although the coefficients depend on $w$, the exponents and logarithmic degree are bounded using a fixed finite set of generators of $M_\nu$. Coefficient extraction in
\eqref{eq:collisionnormalfactor73} defines
$\ell_{\sigma+n,l}(w)$ linearly for homogeneous $w$, and by finite linear extension for every algebraic tensor. The finite uniform sets $\Sigma$ and $\{0,\ldots,L\}$ give the asserted common lower bound and complete the proof.
\end{proof}

\subsection{Component currents and finite-rank collision coefficients}
Keep the actual admissible cover constructed in Lemma~\ref{lem:actualcoverextension}. At each source node write $uv=t^{a_j},\ a_j>0,$ with the fixed compatible coarse and root coordinates. For a finite state cutoff $V_{\leq d}$ use the balanced current frame of \eqref{eq:balancednodeframe}, including its entire $L(1)$ correction and the involution $e^{L(1)}(-1)^{L(0)+1}$. On the normalization we measure a current's order at a node by the orders of its coefficients in this frame. The component order is measured after the opposite-branch state transform.

For $k\geq1$, let $\mathscr G_{{\rm in},k}$ be the union over finite state cutoffs of invariant component currents on the inner normalized cover, with arbitrary finite poles at its actual external marked orbits and with balanced order at least $k$ at every node lift. Define $\mathscr G_{{\rm out},k}$ similarly. Each component current includes its $H$-induced translates. Its residue action is normalized by the same orbit and inertia factors as the original seed.

\begin{lemma}\label{lem:truncatedcomponentlift}
Fix $J\in\N$ and $k\ge1$ such that $a_jk\ge J+1\; \text{at every source node }j$. For $\sigma\in\mathscr G_{{\rm in},k}$ with a finite state cutoff $V_{\le d}$, there are an integer $m=m(\sigma,J,d)\ge0$ and an $H$-invariant current
\[
 s\in H^0\bigl(C_O,\mathscr F_{d,O}(mP)\bigr)^H,
 \qquad O=\C[[t]],
\]
whose restriction to the generic family is an actual global current. Here $P$ is the full external marked-orbit divisor and $\mathscr F_{d,O}$ is the coherent, $O$-flat finite-state current
sheaf on the chosen admissible family, in the balanced node frames. With the component-coordinate residue action (before undoing the external relative-weight dilations),
\begin{equation}\label{eq:truncatedcurrentaction}
 R_{\rm rel}(s)\equiv R_{\rm in}(\sigma)
 \pmod{t^{J+1}\operatorname{Hom}_{\C}(T,T)[[t]]},
 \qquad
 R_{\rm rel}(s)=R_{\rm in}(\sigma)
 +\sum_{r\ge J+1}t^rR_r.
\end{equation}
The equality is coefficientwise on $T=T_{\rm out}\otimes
T_{\rm in}$; every residue coefficient on a fixed algebraic tensor
is a finite sum. Exchanging the components gives the analogous outer-current statement, with zero inner component jet.
\end{lemma}
\begin{proof}
Fix $d,J,k$ and put $O_J=O/(t^{J+1})$.  We check the
actual comparison at a node; a comparison of finite-etale
covers on a punctured three-point component alone does not supply
a filtration statement at its missing node.

For either collision, the coarse coordinates are the actual rational
coordinates $x=Z$, $y=\xi^{-1}$ (product), or $x=Z-1$,
$y=\xi^{-1}$ (iterate), with $xy=\rho$.  At a selected source node,
after the already permitted base change and unit-root choices in
Lemma~\ref{lem:canonicalnodecurrent}, these are \emph{equalities}
\[
 xy=t^{aN},\qquad x=u^N,\qquad y=v^N,\qquad uv=t^a,
\]
where $N$ is the node inertia order.  In particular, the $y$ used
on the actual punctured node is the same $y$ as on the fixed
normalized three-point component; it is not changed by a new
coordinate choice.

Write $v_0$ for the \emph{selected} normalized-component root,
$v_0^N=y$.  By the explicit product/iterate rational charts, the
coarse open $U_v=\mathbb P^1\setminus\{0,1,\infty\}$ away from the
node and external marking disks has constant coordinate over $O_J$:
additional moving-marking equations differ from the removed factors
by units, since $\rho$ is nilpotent.  Hence the actual cover and the
constant normalized cover are finite etale over
$U_v\times\operatorname{Spec}O_J$.  Tag~0BQB of \cite{Stacks}
lifts their fixed, based special-fiber isomorphism uniquely to the
$H_v$-equivariant comparison $\iota_{v,J}$; the comparisons commute
with reduction in $J$ and with induction to the full $H$-orbit.
The root check is independent of this etale statement.  Restrict
$\iota_{v,J}$ to the punctured $y$-disk.  On the selected inertia
sheet the two actual finite-etale rings are
\[
 A_0=O_J((y))[v_0]/(v_0^N-y),\qquad
 A_1=O_J((y))[v]/(v^N-y).
\]
The ring on the right is the indicated summand of the actual node
collar $B_J[v^{-1}]$, $B_J=B/t^{J+1}B$.
Because the comparison is over the same coarse coordinate
$y$, its pullback takes $v_0$ to $\eta v$ with $\eta^N=1$.
The inherited selected root and based lift give $\eta\equiv1\pmod t$.
Since $N$ is invertible, $(\eta-1)(1+\cdots+\eta^{N-1})=0$ implies
$\eta=1$.  Thus
\begin{equation}\label{eq:truncatedrootcomparison}
 \iota_{v,J}^{\sharp}(y)=y,\qquad
 \iota_{v,J}^{\sharp}(v_0)=v
 \quad\text{on the selected punctured node sheet}.
\end{equation}
The same calculation with $u_0^N=x$ holds on the other branch;
for other sheet choices the prescribed constant root phase and the
strict deck pullback are applied \emph{together}, with no
identification of differently twisted modules.  This establishes
precisely the extension to the punctured Kummer collars that is not
provided by Tag~0BQB on its own.

The finite-state current frame is associated to the actual coarse
coordinate and its selected root.  In the frame in which the node
order was defined immediately before this lemma, the transition
between the two oriented branches is the entire constant
operator $\Sigma_d=e^{L(1)}(-1)^{L(0)+1}|_{V_{\le d}}$,
by \eqref{eq:balancednodeframe}; its $L(1)$ expansion is finite,
$\Sigma_d^2=1$, and it commutes with the inertia action.
Equation~\eqref{eq:truncatedrootcomparison}, together with the
complete coordinate representation, therefore identifies the
normalized and thickened balanced coefficients without a negative
Laurent root gauge.  In particular, if the chosen component current
has at the node the full expansion
$\sigma=\sum_{r\ge k}c_rv_0^r$ in the balanced frame, its actual
transport on the punctured $v$-branch is exactly
\begin{equation}\label{eq:truncatedcollartransport}
 \widetilde\sigma_J=\sum_{r\ge k}c_rv^r
 \in v^kO_J[[v]]\otimes V_{\le d}.
\end{equation}
The symbols $c_r$ here include the already prescribed
opposite-branch $\Sigma_d$ transform; it is not applied a second
time.  The external marked-root comparisons use the same fixed
coarse component coordinates and their tame Kummer roots, so
$\widetilde\sigma_J$ has precisely the transported external germs
of $\sigma$, with no additional poles away from the original
marked orbits.

This is a substantive restriction on allowed gauges: the punctured
unit $1+t/v$ is not a comparison of roots both satisfying $v^N=y$,
even though it reduces to $1$ when $t=0$.  Nor may a matrix in
$\operatorname{GL}_d(O_J((v)))$ be inserted independently of the
current-coordinate cocycle.  A regular \emph{node}-frame correction
$G\in\operatorname{GL}_d(B_J)$, if used instead, is also harmless:
$v^kG$ is zero on $D(u)$ when $ak\ge J+1$, even though writing
$u=t^a/v$ may lower the displayed exponent of $v$ while adding
compensating powers of $t$.  Thus one needs the actual rooted
comparison (or this weaker integral weighted property), not a
blanket assertion about arbitrary punctured gauges.

Set the outer section to zero.  At every node the element
\[
 f(v)=\sum_{r\ge k}c_rv^r\in B_J\otimes V_{\le d}
\]
is a \emph{regular} completed-node section whose restrictions
are the prescribed two sections: it equals
\eqref{eq:truncatedcollartransport} on $D(v)$, while on $D(u)$
\[
 f(t^a/u)=\sum_{r\ge k}c_rt^{ar}u^{-r}=0
 \quad\text{in }B_J[u^{-1}]\otimes V_{\le d},
\]
since $ar\ge ak\ge J+1$.  We have constructed a section of the
whole nonreduced node ring, not inferred gluing merely from
agreement on $D(u)\cap D(v)$ (which can be empty over $O_J$).
The construction respects $\Sigma_d$, opposite root rotation and
$H$-transport, and is applied simultaneously on all node orbits.

Let $A$ be the Noetherian local ring of $C_{O_J}$ at one node and
$M$ its finite current module.  Choose an algebraic local element $f$ whose image in the completed
node has \emph{linear term} $u+v$ modulo the square of the maximal
ideal. Such an element exists by the isomorphism of cotangent spaces
under completion; equality to the entire formal series $u+v$ is
not asserted. Its special-fiber reduction vanishes only at the node
on a sufficiently small local neighborhood, so $D(f)$ is the
punctured local curve; $f$ avoids the minimal branch primes and is
a non-zero-divisor.
The already defined punctured component sections give an element
$q\in M_f$.  The completed gluing above says that its image in
$\widehat M_f$ lies in $\widehat M$.  Completion
$A\to\widehat A$ is faithfully flat by
\cite[Tag~00MC]{Stacks}; because $f$ is $M$-regular, the natural
injection $M\hookrightarrow M_f$ gives
\[
 (M_f/M)\otimes_A\widehat A
 \simeq\widehat M_f/\widehat M.
\]
The class of $q$ becomes zero on the right, hence was zero in
$M_f/M$.  Thus $q\in M$.  This proves genuine local algebraic
extension on the nonreduced $O_J$-curve.  Since there are finitely
many nodes and the extensions agree on their punctured overlaps,
they glue with the component sections to a global $H$-invariant
section
\begin{equation}\label{eq:truncatedsectionmod73}
 \bar s_J\in H^0(C_{O_J},
 \mathscr F_{d,O_J}(m_0P_{O_J}))^H
\end{equation}
for some finite external pole bound $m_0$.

By Lemmas~\ref{lem:clustercharts} and
\ref{lem:actualcoverextension}, every irreducible special-fiber
source component maps finitely onto a coarse component meeting an
original external marking.  Hence the actual Cartier divisor $P$
has positive degree on each geometric irreducible special-fiber
component; $\mathcal O(P)$ is relatively ample for this proper
relative nodal curve.  The finite-state sheaf is coherent and
$O$-flat by its actual smooth and balanced free frames, so relative
Serre vanishing \cite[Tag~02O1]{Stacks} gives, after increasing
$m\ge m_0$, $R^1p_*\mathscr F_{d,O}(mP)=0$.
Moreover multiplication by $t^{J+1}$ is injective, and therefore
\[
0\to\mathscr F_{d,O}(mP)\xrightarrow{t^{J+1}}
\mathscr F_{d,O}(mP)\to
\mathscr F_{d,O_J}(mP_{O_J})\to0
\]
is exact.  Since $\operatorname{Spec}O$ is affine, its long
cohomology sequence gives surjectivity of $H^0$ reduction.  Lift
$\bar s_J$ and average by $|H|^{-1}\sum_{h\in H}h$ (defined in
characteristic zero) to obtain the required invariant global
$s$ with $s\bmod t^{J+1}=\bar s_J$.  No bound uniform in
$\sigma,J,d$ is used.

At each actual external marked lift, compare the complete retained
root coordinate and finite-state coordinate frame with the fixed
component chart.  The comparison in the first paragraph preserves
all Laurent coefficients modulo $t^{J+1}$, not just a finite
principal part, and on the other component the transported section
is zero.  Thus, before reversing any external relative-weight
dilation, the marked normalized residue definition
\eqref{eq:canonicalmarkedres}, including $1/N_i$ and the entire
finite $L(1)$ opposite-field transform on an output restricted
dual, gives for every external leg $i$
\[
 R_{i,\mathrm{rel}}(s)\equiv
 \begin{cases}R_{i,\mathrm{in}}(\sigma),&i\text{ inner},\\
 0,&i\text{ outer},\end{cases}
 \pmod{t^{J+1}}.
\]
The deck map is the retained strict pullback in its actual twisting
sector; this comparison does not identify $T_hW$ with $W$ by fiat.
A fixed state cutoff, bounded marked poles and twisted-field lower
truncation show that, on any fixed algebraic tensor, only finitely
many modes contribute to a given $t$-coefficient.  The reduction
therefore takes place in the coefficientwise operator module of
\eqref{eq:truncatedcurrentaction} and does not require completion of $\operatorname{End}(T)$.  Summing the finite orbit-normalized marked
residues proves that equation.  Interchanging the two normalized
components, and using $u$ in place of $v$ in the same local argument,
proves the outer-current assertion.
\end{proof}

The finite-order lifting follows the truncated sewing-current method of \cite[Section~2]{GuiZhangSF}; unlike the untwisted setting there, the balanced node frame and tame-root covariance are verified explicitly above.

On either normalized component $\widetilde C_\bullet$, let $P_\bullet$ be the reduced divisor of its external marked orbits and $Q_\bullet$ the reduced divisor of all its node lifts. Both are
$H$-stable and $P_\bullet\cap Q_\bullet=\varnothing$. Let $F_{E,\bullet}$ be the finite-rank canonical current bundle restricted to that component, with its retained root coordinates. For $k\ge1$ define the coherent subsheaf
\[
 F_{E,\bullet,k}:=F_{E,\bullet}(-kQ_\bullet)
 =\ker\!\left(F_{E,\bullet}\longrightarrow
   \bigoplus_{p\in Q_\bullet}
   (F_{E,\bullet})_p/\mathfrak m_p^k(F_{E,\bullet})_p\right).
\]
Here vanishing is measured in the balanced node frames, so the notation $(-kQ_\bullet)$ is understood in those frames. Since $\widetilde C_\bullet$ is smooth, this sheaf is locally free;
its restriction to a neighborhood of $P_\bullet$ is $F_{E,\bullet}$. For $\bullet\in\{{\rm in},{\rm out}\}$, write $\Gamma_{\bullet,k}:=\mathscr G_{\bullet,k}T_\bullet,\; B_{\bullet,k}:=T_\bullet/\Gamma_{\bullet,k}$.
These are $\C$-vector-space quotients; the currents and their actions are taken on the actual normalized component with the inherited monodromy and based roots.

\begin{lemma}\label{lem:partialcomponentfinite}
For every $k\in\Z_{\ge1}$, $\dim_\C B_{{\rm in},k}<\infty,\; \dim_\C B_{{\rm out},k}<\infty$.
\end{lemma}
\begin{proof}
Fix a component and $k$. It suffices to use the invariant currents valued in $E\subset V$: their action space is contained
in $\Gamma_{\bullet,k}$. Let $p:\widetilde C_\bullet\to
\operatorname{Spec}\C$ denote the proper smooth component cover.
Every connected component of $\widetilde C_\bullet$ meets $P_\bullet$ (Lemma~\ref{lem:clustercharts}); hence $\mathcal O_{\widetilde C_\bullet}(P_\bullet)$ is ample. Since $F_{E,\bullet,k}$ is coherent, Serre vanishing gives an integer
$m=m(E,k)\ge0$ such that
\begin{equation}\label{eq:componentserre75}
 H^1\!\left(\widetilde C_\bullet,
   F_{E,\bullet,k}(mP_\bullet)\right)=0.
\end{equation}
The $H$-action preserves all the divisors and the balanced frames; $|H|^{-1}\sum_{h\in H}h$ is defined over $\C$.

Fix an external marked orbit $H\widetilde p_i$, a homogeneous $b\in B_{g_i}\subset E$ and an admissible mode $\mu=\alpha_b-n$, with $n$ sufficiently large. Denote by $D_i=\sum_{p\in H\widetilde p_i}p$ its orbit divisor.
For $r\ge0$ there is an exact sequence
\[
 0\longrightarrow F_{E,\bullet,k}(mP_\bullet)
 \longrightarrow F_{E,\bullet,k}(mP_\bullet+rD_i)
 \longrightarrow\mathscr P_{i,r}\longrightarrow0,
\]
where $\mathscr P_{i,r}$ is the finite-length principal-part sheaf supported on $D_i$. By \eqref{eq:componentserre75},
\begin{equation}\label{eq:componentjet75}
 H^0\!\left(\widetilde C_\bullet,
   F_{E,\bullet,k}(mP_\bullet+rD_i)\right)^H
 \twoheadrightarrow H^0(\mathscr P_{i,r})^H .
\end{equation}
The invariant prescribed principal part is the pure mode germ $u_i^{N_i\mu}\iota_{x_i}(b)\,dx_i$, induced over $H\widetilde p_i$, and is zero at the other external orbits. As $F_{E,\bullet,k}=F_{E,\bullet}$ near $P_\bullet$ and the
root-coordinate gauge is a finite Laurent matrix (Lemma~\ref{lem:gauge}), this germ defines an element of $H^0(\mathscr P_{i,r})^H$ for some finite $r$. Choose an invariant lift $s_{i,b,\mu}$ by
\eqref{eq:componentjet75}; it belongs to
$\mathscr G_{\bullet,k}$ by construction.

The correction bound must be computed after subtraction of the \emph{entire prescribed principal part}. All remaining germs belong to the fixed space $F_{E,\bullet,k}(mP_\bullet)$, independently of $n$ and of the auxiliary pole order $r$. The finite inverse root-coordinate gauge has a fixed lower Laurent exponent $-P_i$ at each marking, and $E$ has a finite homogeneous basis. Thus any correction mode $c_\nu$ has $N_i\nu\ge -P_i$ and
\[
 \deg(c_\nu)=\operatorname{wt}(c)-\nu-1
 \le \max_{c\in\mathcal B_i}\operatorname{wt}(c)+P_i/N_i-1.
\]
Taking the maximum over the finitely many markings gives $C_k$. At the output marking the current acts on the \emph{actual restricted contragredient} $D'$ using its $g^{-1}$-twisted modes. The opposite-field $L(1)^j$ terms are first included in the finite gauge, and each such state has its true weight $\operatorname{wt}(c)-j$. Its modes obey the same conformal-weight commutator on $D'$, so the same finite \emph{upper} degree bound applies, without identifying output and input mode labels or reversing the energy inequality.

Its external residue operators have the form
\begin{equation}\label{eq:componentward75}
 R_i(s_{i,b,\mu})=b_\mu+U_i,
 \qquad R_j(s_{i,b,\mu})=U_j\quad(j\ne i),
\end{equation}
where the germs of the correction terms belong to $F_{E,\bullet,k}(mP_\bullet)$ near $P_\bullet$. In particular there is an integer $C_k$, depending on $m,E,k$
and the finitely many marked-root gauges but not on $n$, such that $\deg(c_\nu):=\operatorname{wt}(c)-\nu-1\le C_k$ for every homogeneous mode summand $c_\nu$ of any $U_j$. The same bound applies to the full finite $L(j)$ corrections
in the root-coordinate gauges. This is the argument of
\eqref{eq:boundedcorrection}, now applied to
$F_{E,\bullet,k}(mP_\bullet)$.

Choose $R_k\in\Z$ satisfying $\Delta(b,n):=\operatorname{wt}(b)-\alpha_b+n-1>C_k,\; b\in B_{g_i},\; n\ge R_k$ for every external module on the component. Fix $N_{\rm PBW}\ge\max(R_k-1,\max_uL_u)$ and use the
$g_i$-eigenvector spanning families of
\cite[Proposition~2.12]{Tan}, as fixed in
\eqref{eq:PBW}. If a tensor word contains a leftmost
$b_{\alpha_b-n}^{(i)}$ with $n\ge R_k$, write $w=b_\mu^{(i)}v$. By applying the Ward relation for $s_{i,b,\mu}$ to $v$, \eqref{eq:componentward75} gives
\[
 [w]= -\sum_j[U_j^{(j)}v]\quad\text{in }B_{\bullet,k},
 \qquad
 \ee(c_\nu^{(j)}v)-\ee(w)
 =D\bigl(\deg(c_\nu)-\Delta(b,n)\bigr)<0
\]
for every nonzero homogeneous correction term; $D$ is a common positive integer clearing the relative-weight denominators. The sum on each fixed $v$ is finite by lower truncation.
Lemma~\ref{lem:descent}, applied over $A=\C$, reduces the quotient to the image of the finite-dimensional space $\bigotimes_i K_{R_k,N_{\rm PBW}}(M_i)$ by Lemma~\ref{lem:lowdepth}. Thus $\dim B_{\bullet,k}<\infty$.
The same argument applies to the other normalized component.
\end{proof}

For $\beta\in\Sigma+\N$ define
\[
\begin{aligned}
 J(\beta)&=\max\bigl(\{0\}\cup
  \{n\in\N:\ \beta=\sigma+n,\ \sigma\in\Sigma\}\bigr),\\
 k(\beta)&=\max\left\{1,\ \max_j
  \left\lceil\frac{J(\beta)+1}{a_j}\right\rceil\right\}.
\end{aligned}
\]
The maximum is finite, and $a_jk(\beta)\ge J(\beta)+1$ for every source node $j$. Put $\ell_{\beta,l}=0$ if $\beta\notin\Sigma+\N$ or $l\notin\{0,\ldots,L\}$.

\begin{theorem}\label{thm:finitecollisionkernels}
For every $\beta\in\Sigma+\N$ and $0\le l\le L$,
\[
 \ell_{\beta,l}\bigl(\Gamma_{{\rm out},k(\beta)}
             \otimes T_{\rm in}
       +T_{\rm out}\otimes\Gamma_{{\rm in},k(\beta)}\bigr)=0,
\]
so that
\begin{equation}\label{eq:collisionfactorrank76}
 \ell_{\beta,l}\in
 \bigl(B_{{\rm out},k(\beta)}\otimes
       B_{{\rm in},k(\beta)}\bigr)^*
 \cong B_{{\rm out},k(\beta)}^*\otimes
       B_{{\rm in},k(\beta)}^*.
\end{equation}
In particular, $\ell_{\beta,l}$ has finite tensor rank as a bilinear form on $T_{\rm out}\times T_{\rm in}$.
\end{theorem}
\begin{proof}
Put $J=J(\beta)$ and $k=k(\beta)$. On a homogeneous $w=o\otimes a\otimes b\in T_{\rm out}\otimes T_{\rm in}$ set
\[
 \mathsf D_t(w)=t^{e((\operatorname{wt}a-h_A)
                         +(\operatorname{wt}b-h_B))}w,
 \qquad c_t=t^{e(h_A+h_B)}.
\]
The exponent defining $\mathsf D_t$ is in $\N$, and the collision normalization \eqref{eq:collisionnormalizedkernel} is
\begin{equation}\label{eq:collisiongauge76}
 \widehat F(t;w)=c_t F(t;\mathsf D_t w).
\end{equation}
For the physical Ward action $R_{\rm phys}(s)$ we \emph{calculate} its relative-frame action using the geometric component roots. If $b\in V_d$ and $n$ is an allowed, possibly fractional, twisted mode, put $\Delta_{b,n}=d-n-1$. On $A$ and $B$ (and only on these two legs) the physical-to-component linear scale is $x=t^e\xi$. The complete current density transformation, including the Jacobian and the fixed $1/N_i$ root-orbit normalization, carries the physical $b_n$ term to $t^{-e\Delta_{b,n}}b_n$. Directly, for a homogeneous input $v$ of relative weight $\gamma$,
\[
 (\mathsf D_t^{-1}b_n\mathsf D_t)v
 =t^{e\gamma-e(\gamma+\Delta_{b,n})}b_nv
 =t^{-e\Delta_{b,n}}b_nv.
\]
The shift $\Delta_{b,n}$ is in $e^{-1}\Z$ for the selected common ramification; hence every displayed exponent is an integer. For a nonprimary state the finite coordinate transformation replaces $b$ by the finitely many $L(1)^jb$; each term is treated with \emph{its own} weight $d-j$ and allowed mode index. The exact root comparisons, continuous lowest-weight branches and typed noncommutative deck pullbacks are those of the actual component-coordinate action in Lemma~\ref{lem:truncatedcomponentlift}; the relative lowest weights cancel inside the conjugation. On the unscaled $D'$ and $C$ legs the dilation is the identity. The $D'$ action is the restricted-contragredient opposite field, including the finite $L(1)$ correction, not a transpose of an arbitrary positive-mode operator on a completion.

Consequently, after this residue calculation, one obtains the correctly typed equality
$R_{\rm rel}(s)=\mathsf D_t^{-1}R_{\rm phys}(s)\mathsf D_t$
in the common fraction field, and for the currents lifted in Lemma~\ref{lem:truncatedcomponentlift} it has the integral coefficientwise expansion \eqref{eq:truncatedcurrentaction}. This is not a definition of an unrelated component action and is not an unconditional assertion of integrality for arbitrary meromorphic currents. Theorem~\ref{thm:fullsmallward} and
\eqref{eq:collisiongauge76} therefore imply
\begin{equation}\label{eq:collisionrelativeward76}
 \widehat F(t;R_{\rm rel}(s)w)
 =c_tF(t;R_{\rm phys}(s)\mathsf D_t w)=0
\end{equation}
for each actual invariant current $s$ and algebraic $w$. Here a formal operator series is evaluated coefficientwise; $\mathsf D_t$ is diagonal on homogeneous vectors and need not be replaced by a scalar when a current changes their weights.

Choose $\sigma\in\mathscr G_{{\rm in},k}$. By Lemma~\ref{lem:truncatedcomponentlift}, there exists an actual global invariant $s$ with
\[
R_{\rm rel}(s)=R_{\rm in}(\sigma)
  +\sum_{j\ge J+1}t^jS_j,
 \qquad S_j\in\operatorname{Hom}_{\C}(T,T),
\]
where for any fixed algebraic input $w$, each $S_jw$ is a finite algebraic sum and the congruence is understood in the relative frames. Apply \eqref{eq:collisionrelativeward76} to $w$:
\begin{equation}\label{eq:collisionannihilate76}
 \widehat F(t;R_{\rm in}(\sigma)w)
   =-\sum_{j\ge J+1}t^j\widehat F(t;S_jw).
\end{equation}
By Proposition~\ref{prop:collisionuniformsupport}, every $\widehat F(t;S_jw)$ has exponents in the same $\Sigma+\N$ and logarithmic degree at most $L$, independently of $j,w$. If $t^j\widehat F(t;S_jw)$ had a term
$t^\beta(\log t)^l$, then $\beta-j=\sigma_0+n,\; \sigma_0\in\Sigma,\ n\in\N\; \Longrightarrow\; j\le\beta-\sigma_0\le J$, contradicting $j\ge J+1$. This calculation does not evaluate an infinite sum of endomorphisms on the algebraic tensor product. Fix $w$ and first form each of the finite algebraic tensors $S_jw$. For each exponent class and exponent $\beta$, the condition $\beta-j\in\Sigma+\N$ implies
\[
 0\le j\le \max_{\sigma\in\Sigma}(\Re\beta-\Re\sigma),
\]
so only finitely many $j$ can affect that coefficient. Each such $S_jw$ is a finite sum of homogeneous algebraic tensors and therefore has the same $\Sigma,L$ bounds. The logarithmic factors change only their logarithmic degree, not their exponent; no cancellation across different $\C/\Z$ classes is used. Thus the sum in \eqref{eq:collisionannihilate76} is a legitimate coefficientwise sum. As every possible contribution to $t^\beta(\log t)^l$ would require $j\le J$, none of the terms with $j\ge J+1$ contributes. Therefore $\ell_{\beta,l}(R_{\rm in}(\sigma)w)=0$. The outer-current lifting gives $\ell_{\beta,l}(R_{\rm out}(\sigma')w)=0$ for every $\sigma'\in\mathscr G_{{\rm out},k}$.

Over $\C$ the tensor product is exact, so
\[
 \frac{T_{\rm out}\otimes T_{\rm in}}
 {\Gamma_{{\rm out},k}\otimes T_{\rm in}
      +T_{\rm out}\otimes\Gamma_{{\rm in},k}}
 \cong B_{{\rm out},k}\otimes B_{{\rm in},k}.
\]
Both factors are finite dimensional by Lemma~\ref{lem:partialcomponentfinite}, and taking duals yields
\eqref{eq:collisionfactorrank76}.
\end{proof}

For $o\in T_{\rm out}$ and $\beta\in\Sigma+\N$ put
\[
 \lambda_{\beta,l,o}(a,b)=\ell_{\beta,l}(o;a,b),\qquad
 K_{\beta,l}=\operatorname{Span}_{\C}
       \{\lambda_{\beta,l,o}:o\in T_{\rm out}\},\qquad
 K_\beta=\sum_{l=0}^L K_{\beta,l},
\]
with $\lambda_{\beta,L+1,o}=0$. By Theorem~\ref{thm:finitecollisionkernels}, $K_{\beta,l}\subseteq B_{{\rm in},k(\beta)}^*,\; \dim_\C K_\beta<\infty$.

\begin{lemma}\label{lem:collisionenergy}
Define $\mathcal E\in\operatorname{End}_\C(T_{\rm in}^*)$ by
\[
 (\mathcal E\lambda)(a,b)
 =\lambda(L_A(-1)a,b)
  +\lambda(L_A(0)a,b)+\lambda(a,L_B(0)b).
\]
Then, for $\beta\in\Sigma+\N$, $0\le l\le L$ and $o\in T_{\rm out}$,
\begin{equation}\label{eq:collisionenergychain}
 (\mathcal E-\beta/e)\lambda_{\beta,l,o}
    =\frac{l+1}{e}\lambda_{\beta,l+1,o}.
\end{equation}
Consequently
\[
 \mathcal E K_\beta\subseteq K_\beta,
 \qquad
 (\mathcal E-\beta/e)^{L+1}K_\beta=0.
\]
\end{lemma}
\begin{proof}
Take $a\in A_\alpha$, $b\in B_{\beta'}$ and $o\in T_{\rm out}$, and fix the branch $\log q=e\log t$ for $q=t^e$, $t>0$. Set, initially for $q>0$,
\[
 \widehat F_q(q;o;a,b)
   :=q^{\alpha+\beta'}F(1+q,1;o;a,b),
 \qquad
 \widehat F(t;o;a,b)=\widehat F_q(t^e;o;a,b).
\]
The $L(-1)$-derivative property of the actual product intertwining operators \cite[Definition~2.8]{DH}, with $z_2=1$ fixed, gives $\partial_q F(1+q,1;o;a,b) =F(1+q,1;o;L_A(-1)a,b)$. Since $L_A(-1)a\in A_{\alpha+1}$ and
$L_A(0)a=\alpha a$, $L_B(0)b=\beta'b$,
\begin{align*}
 q\partial_q\widehat F_q(q;o;a,b)
 &=q^{\alpha+\beta'+1}
       F(1+q,1;o;L_A(-1)a,b)
     +(\alpha+\beta')\widehat F_q(q;o;a,b)\\
 &=\widehat F_q(q;o;L_A(-1)a,b)
   +\widehat F_q(q;o;L_A(0)a,b)
   +\widehat F_q(q;o;a,L_B(0)b).
\end{align*}
Thus $t\partial_t\widehat F(t;o;a,b)$ equals $e$ times the sum of the three normalized evaluations on the last line. By Proposition~\ref{prop:collisionuniformsupport}, both sides have a common locally normally convergent Nilsson expansion
with exponent set $\Sigma+\N$ and logarithmic degree $\le L$. Since $t\partial_t\bigl(t^{\beta}(\log t)^l\bigr) =\beta t^{\beta}(\log t)^l +l\,t^{\beta}(\log t)^{l-1}$,
coefficient uniqueness in the chosen branch gives $\beta\lambda_{\beta,l,o} +(l+1)\lambda_{\beta,l+1,o} =e\,\mathcal E\lambda_{\beta,l,o}$. This is \eqref{eq:collisionenergychain}. It holds for arbitrary
$a,b$ by their finite homogeneous decompositions. Iterating
\eqref{eq:collisionenergychain} proves the two final assertions.
\end{proof}

Let $\pi_\bullet:\widetilde C_\bullet\to C_\bullet$ be either inherited normalized $H$-cover, with the retained based roots and physical-coordinate identifications. Denote by $P_\bullet$ the
reduced divisor of its external marked lifts and by $Q_\bullet$ the reduced divisor of its node lifts; they are disjoint and $H$-stable. Every connected component of $\widetilde C_\bullet$ meets $P_\bullet$ by Lemma~\ref{lem:clustercharts}.
Let $T_\bullet$ be the algebraic tensor product of the ordinary incoming modules on that component. For $d\in\N$ let $F_{d,\bullet}$ denote its finite-rank current bundle with states in $V_{\le d}$ and the balanced node frames of
Lemma~\ref{lem:canonicalnodecurrent}. For $k\ge1$ set
\[
 \mathscr G_{\bullet,d,k}
 :=H^0\!\left(\widetilde C_\bullet,
       F_{d,\bullet}(*P_\bullet)(-kQ_\bullet)\right)^H,
 \qquad
 \mathscr G_{\bullet,k}:=\bigcup_{d\ge0}\mathscr G_{\bullet,d,k},\qquad
 \mathscr G_{\bullet,k}T_\bullet
 :=\sum_{s\in\mathscr G_{\bullet,k}}R_\bullet(s)T_\bullet.
\]
Here $(-kQ_\bullet)$ denotes vanishing of each \emph{balanced root-frame coefficient} to order $\ge k$; $R_\bullet$ is the orbit-normalized current action, and every summand acts through a finite state cutoff. Define
\begin{equation}\label{eq:partialnodedualspace}
 B_k:=B_{\bullet,k}=T_\bullet/\mathscr G_{\bullet,k}T_\bullet,
 \qquad
 \mathscr D_k:=B_k^*\subseteq T_\bullet^*,
 \qquad
 \mathscr D:=\bigcup_{k\ge1}\mathscr D_k.
\end{equation}
For a fixed component $\bullet$, the inclusions $\mathscr D_k\subseteq\mathscr D_{k+1}$ follow from
$\mathscr G_{\bullet,d,k+1}\subseteq
\mathscr G_{\bullet,d,k}$. This is the same full-state quotient $B_{\bullet,k}$ used in Theorem~\ref{thm:finitecollisionkernels}; $B_k$ abbreviates it when the component is fixed. Lemma~\ref{lem:partialcomponentfinite} already proves $\dim_\C B_k<\infty$: its proof uses only $E$-valued currents, whose relations are contained in the full-state relation space. No additional quotient is being formed here.
Write $\mathscr G_k=\mathscr G_{\bullet,k}$ when the component is fixed. Write $\mathscr D_{\rm in}$ and $\mathscr D_{\rm out}$ for the spaces obtained from the two components. At the inner outgoing node
the chosen based path has label $h^{-1}$, and at the outer node it has label $h$, where $h=g_Ag_B$; these are the inverse inertia labels fixed by the balanced covering construction in Lemma~\ref{lem:canonicalnodecurrent}.

Let $N=\operatorname{ord}(h)$ and choose a node root $z=s^N$. For $f(z)\in\C((z))$, the change-of-variable formula is
\begin{equation}\label{eq:normalizednoderesidue}
 \operatorname{Res}_{z=0}f(z)\,dz
 =\frac1N\operatorname{Res}_{s=0}
       f(s^N)Ns^{N-1}\,ds.
\end{equation}
For an $H$-invariant paired one-form $\eta$ and a node lift $\widetilde Q$ with stabilizer of order $N$,
\[
 \frac1{|H|}\sum_{p\in H\widetilde Q}
            \operatorname{Res}_{p}\eta
 =\frac1N\operatorname{Res}_{\widetilde Q}\eta.
\]
The scalar paired integrand is invariant, so
\eqref{eq:normalizednoderesidue} applies after descent even when
individual state fields have fractional modes.

For homogeneous $u\in V_d$ in the $r$-th eigenspace of the node inertia, the physical balanced frame is $e_u=\iota_z(u)z^{d-1}dz$, with the full opposite-branch $L(1)$-transformation. An invariant root coefficient has the form
\[
 f_u(s)\in s^r\C((s^N)),\qquad
 s^j e_u\text{ is allowed only if }j\equiv r\pmod N.
\]
Its monomial $s^j e_u=z^{d-1+j/N}\iota_z(u)\,dz$ extracts the mode $u_{d-1+j/N}$ by the normalized residue formula; in particular $\deg u_{d-1+j/N}=-j/N,\; j\ge k\ \Longrightarrow\ \deg u_{d-1+j/N}\le-k/N$.
The subsequence $k=N(a+1)$, $a\in\N$, is cofinal. For such $k$ the allowed node coefficients in the $r$-eigenspace are precisely $s^{N(a+1)+r}\C[[s^N]]\,e_u =z^{d+a}s^r\C[[z]]\,\iota_z(u)\,dz$. This identifies the root-order filtration with the
weight-dependent outgoing-point filtration of
\cite[Section~3.1, definition of the outgoing filtration]{GuiZhangPartial}
(in the untwisted case $r=0$); the inertial congruence and full balanced frame are retained here.

\subsection{Free-point propagation and the intermediate twisted modules}
\begin{lemma}\label{lem:partialequivariantresidue}
Fix a normalized component $\widetilde C_\bullet$, $k\ge1$, $\lambda\in\mathscr D_k$, $w\in T_\bullet$ and $d\ge0$. Put
\[
 \mathscr F_{d,k}:=F_{d,\bullet}(-kQ_\bullet),\qquad
 \mathscr F_{d,k}^{\dagger}
 :=\mathcal Hom_{\cO_{\widetilde C_\bullet}}
       (\mathscr F_{d,k},\omega_{\widetilde C_\bullet})
 \simeq F_{d,\bullet}^{\vee}\otimes
       \omega_{\widetilde C_\bullet}(kQ_\bullet).
\]
The calibrated incoming fields $Y_i^{g_i}(u,z_i)$, $u\in V_{\le d}$, determine $H$-compatible formal germs $\tau_{d,p}$ of $\mathscr F_{d,k}^{\dagger}$ at every $p\in P_\bullet$. They are the Laurent expansions of a unique
\emph{$H$-invariant} section
\begin{equation}\label{eq:partialpropagatedsection78}
 \tau_d\in H^0\!\left(\widetilde C_\bullet\setminus P_\bullet,
             \mathscr F_{d,k}^{\dagger}\right)^H.
\end{equation}
This section is meromorphic at $P_\bullet$ and, in a local $F_{d,\bullet}^{\vee}\otimes\omega$ frame, has poles of order at most $k$ along $Q_\bullet$. For $d\le d'$ the restriction of $\tau_{d'}$ to $V_{\le d}$ is $\tau_d$.
Consequently evaluation at $p\notin P_\bullet\cup Q_\bullet$ provides a well-defined linear functional $P\lambda(u,p;w)$, $u\in V$, with the inherited finite-order monodromy. The construction does not assert parameter
holomorphy for varying $p$ without the additional finite-jet argument below.
\end{lemma}
\begin{proof}
All bundles in the statement are on the \emph{smooth, proper} normalized component cover. The divisor $Q_\bullet$ is reduced, $H$-stable and disjoint from $P_\bullet$; hence $\mathscr F_{d,k}$ and $\mathscr F_{d,k}^{\dagger}$ are
finite-rank vector bundles. In particular the node condition is encoded in the test-current bundle itself:
\begin{equation}\label{eq:partialnodepole78}
 H^0\bigl(\widetilde C_\bullet,
        \mathscr F_{d,k}(*P_\bullet)\bigr)^H
   =\mathscr G_{\bullet,d,k},\qquad
 \mathscr F_{d,k}^{\dagger}
   =F_{d,\bullet}^{\vee}\otimes\omega(kQ_\bullet).
\end{equation}
Here $\omega=\omega_{\widetilde C_\bullet}$. Thus regularity of a section of $\mathscr F_{d,k}^{\dagger}$ at the node means
\emph{at most} a $k$-th order pole in the unfiltered dual
frame; it does not mean regularity in that latter frame.

Let $p\in P_\bullet$ lie over the incoming marking $i$. For $u\in V_{\le d}$, $v_i$ the $i$-th tensor factor of $w$, and $z_i=s_i^{N_i}$ the chosen physical/root coordinates, use the original $g_i$-twisted field
$Y_i^{g_i}(u,z_i)v_i$ and the calibrated finite Laurent coordinate gauge of Lemma~\ref{lem:gauge} to define the local formal dual germ $\tau_{d,p}$. More precisely, its pairing with a local current is the corresponding local
field-residue pairing with $\lambda$; the physical residue is converted to a root residue by
\eqref{eq:normalizednoderesidue}. For each $u$ and $w$, lower
truncation gives $Y_i^{g_i}(u,z_i)v_i\in M_i((z_i^{1/N_i}))$, with a finite negative part. As $V_{\le d}$ is finite dimensional and the calibrated gauge has only finitely many Laurent entries, $\tau_{d,p}$ has a finite pole order in a
finite-rank root trivialization. The $g_i$-eigenvalue congruences and the actual based-root identifications supply $\tau_{d,hp}=h\tau_{d,p}$ for $h\in H$. No monodromy-free interpretation of an individual fractional field is made.

For any $s\in\mathscr G_{\bullet,d,k}$, the definitions of the orbit-normalized actions give the \emph{unnormalized} residue identity
\begin{equation}\label{eq:partialresidue78}
 \sum_{p\in P_\bullet}
     \Res_p\langle\tau_{d,p},s\rangle
       =|H|\,\lambda(R_\bullet(s)w)=0.
\end{equation}
All local residues on the left are finite algebraic sums: $s$ has finite-order poles, $d<\infty$, and lower truncation applies to the fixed vector $w$. The factor $|H|$ follows from summing the full $H$-orbit, including its inertia
stabilizers, rather than from treating one root lift as an unramified point.

Let now $s\in H^0(\widetilde C_\bullet,
\mathscr F_{d,k}(*P_\bullet))$ be arbitrary and put
$\operatorname{Av}_H(s)=|H|^{-1}\sum_{h\in H}hs$. Since the collection $\{\tau_{d,p}\}_{p\in P_\bullet}$ is $H$-compatible, its total residue pairing with $s$ equals that with $\operatorname{Av}_H(s)$. Thus
\eqref{eq:partialresidue78} holds for all test currents.
Apply the finite-rank strong residue theorem
\cite[Theorem~6.1]{FS} (equivalently its bundle formulation in
\cite[Section~3, strong residue theorem]{GuiZhangPartial})
to $\mathscr F_{d,k}$ and to these incoming formal dual germs. The theorem yields a section of $\mathscr F_{d,k}^{\dagger}$ on $\widetilde C_\bullet\setminus P_\bullet$ with precisely the
prescribed Laurent expansions. The finite negative parts show that it is meromorphic at the incoming marks. Averaging the result over $H$ leaves those germs unchanged and gives
\eqref{eq:partialpropagatedsection78}. This is an application
of the \emph{finite-rank} residue theorem, not of an infinite-rank or untwisted propagation statement.

Every connected component of $\widetilde C_\bullet$ meets $P_\bullet$ by Lemma~\ref{lem:clustercharts}. A meromorphic dual section with zero complete Laurent germ at an incoming point therefore vanishes identically on that connected
component. The maps $F_{d,\bullet}\hookrightarrow F_{d',\bullet}$ respect the balanced filtration and the calibrated incoming fields. Restriction of $\tau_{d'}$ to $F_{d,\bullet}$ has the same
incoming germs as $\tau_d$, so uniqueness proves their compatibility. Their evaluations define $P\lambda$ on the algebraic direct union $V=\bigcup_dV_{\le d}$, without any infinite sum over $d$.

The following parameter criterion will be used for the finite-rank construction. Let $B$ be a smooth parameter open, and suppose a finite-rank bundle $\mathcal A$ on $\widetilde C_\bullet\times B$ already includes a
\emph{fixed relative effective polar divisor}, finite flat
over $B$; the incoming divisor $P_\bullet\times B$ is fixed and meets every connected component. Suppose the prescribed incoming Laurent germs have a common finite negative bound and depend holomorphically, coefficient by
coefficient, on $b\in B$. Assume that for every $b$ the pointwise residue criterion supplies a section of $\mathcal A_b(*P_\bullet)$ with those germs.

After shrinking $B$ about a chosen point, choose $m\gg0$ so that relative Serre vanishing and cohomology and base change yield
\[
 R^1\operatorname{pr}_{B*}\mathcal A(mP_\bullet)=0,
 \qquad
 \mathcal H:=\operatorname{pr}_{B*}\mathcal A(mP_\bullet)
 \text{ is locally free},
\]
with fiber $H^0(\widetilde C_\bullet,\mathcal A_b(mP_\bullet))$. The pole bound may be enlarged once to contain the pointwise sections locally on $B$; equivalently one assumes a uniform
such $m$ before applying this criterion. Since each connected component meets $P_\bullet$, the complete incoming jets separate $\mathcal H_b$; finite dimensionality gives an $r\ge1$ such that the length-$r$ jet map is injective at the
chosen base point. By semicontinuity, after further shrinking,
\[
 \operatorname{jet}_r:\mathcal H\longrightarrow
 \operatorname{pr}_{B*}\!
 \left(\mathcal A(mP_\bullet)|_{r(P_\bullet\times B)}\right)
\]
is an injective bundle map with locally free cokernel. The holomorphic $r$-jet vector of the prescribed germs lies fiberwise in its image; hence it has a unique holomorphic preimage in $\mathcal H$. By pointwise uniqueness and
existence, that preimage has all the prescribed germs. Such local families glue by uniqueness. In particular, pointwise strong-residue existence alone is not being used to claim parameter holomorphy without a fixed pole bound.
\end{proof}

\begin{lemma}\label{lem:partialpropagationward}
Let $k\ge1$, $\lambda\in\mathscr D_k$, and let $p$ be a free $H$-orbit disjoint from $P_\bullet\cup Q_\bullet$. The functional $P\lambda(\,\cdot\,,p;\,\cdot\,)$ of Lemma~\ref{lem:partialequivariantresidue} satisfies all Ward
relations of the pointed component with an additional \emph{untwisted} $V$-module at $Hp$, without changing the outgoing level $k$. The construction can be iterated. For two free orbits and a common
physical coordinate $\xi=p_1-p_2$ in a simply connected disk disjoint from the other marked orbits,
\begin{equation}\label{eq:partialpropagationope}
 P^2\lambda(u,p_1;v,p_2;w)
   =\sum_{m\in\mathbb Z}\xi^{-m-1}
       P\lambda(u_mv,p_2;w)
   =P\lambda(Y(u,\xi)v,p_2;w).
\end{equation}
The Laurent series in \eqref{eq:partialpropagationope} converges normally on compact subsets of a sufficiently small punctured disk; its entire regular part ($m<0$) is included. The opposite ordering
represents the same meromorphic functional away from all orbital diagonals, with the corresponding $H$-transported states and frames.
\end{lemma}
\begin{proof}
We distinguish the finite-rank Ward calculation, the family construction, and the two-free-point comparison. The untwisted propagation and double-propagation statements of
\cite[Section~3, ``Propagation of partial conformal blocks'' and
``Double propagation'']{GuiZhangPartial} provide a comparison of the method, but are not applied as results for twisted modules here.

\noindent(i)\ 
Fix $d,k,w,v$ and a free point $p_0$. On an $H$-free parameter neighborhood $B$ of $p_0$ write
$\Gamma=\sum_{h\in H}\Gamma_h\subset
 \widetilde C_\bullet\times B$, where
$\Gamma_h=\{(x,p):x=hp\}$. The graphs are pairwise disjoint on $B$, and $\Gamma$ is a relative effective Cartier divisor. If
$\sigma\in H^0(\widetilde C_\bullet,
 F_{d,\bullet}(*P_\bullet+KHp_0)(-kQ_\bullet))^H$,
choose $m\gg0$ so that
\[
 R^1\operatorname{pr}_{B*}\mathcal A_{d,k,K,m}=0,
 \quad
 \mathcal A_{d,k,K,m}
 :=\operatorname{pr}_1^*F_{d,\bullet}(-kQ_\bullet)
       \otimes\mathcal O(mP_\bullet+K\Gamma).
\]
Here $mP_\bullet$ is relatively ample on every connected component, and $\mathcal A_{d,k,K,m}$ is a vector bundle, flat over $B$. After shrinking $B$, cohomology and base change makes
$\operatorname{pr}_{B*}\mathcal A_{d,k,K,m}$ locally free and its fiber at $p_0$ contains $\sigma$. Extend $\sigma$ in a local frame, and average over $H$; this yields a family $\sigma(p)$ of actual
currents of \emph{the same node level $k$}. For a fixed finite set of test currents, one $d,K,m$ suffices. Near a collision with an old incoming point, each algebraic family is meromorphic after the
finite root base change. Multiplication by some parameter power makes its coefficients integral. No holomorphic dependence of $P\lambda$ is assumed in this construction.

\noindent(ii)\ 
Work first with a common unramified local coordinate $z$, with $p$ represented by $q$. For $u\in V$ and the test coefficient $(z-q)^{-1}$, the two expansions are
\begin{align*}
 \iota_{z,q}(z-q)^{-1}
 &=\sum_{j\ge0}q^jz^{-j-1},&
 \iota_{q,z}(z-q)^{-1}
 &=-\sum_{j\ge0}q^{-j-1}z^j.
\end{align*}
Consequently $R_>(q)=\sum_{j\ge0}q^ju_{-j-1}$ and $R_<(q)=-\sum_{j\ge0}q^{-j-1}u_j$ act on the original incoming module. The actual module Jacobi identity,
\eqref{eq:twistedBorcherds} with $g=1$ and $p=-1$, gives
\[
R_>(q)Y_i(v,q)w-Y_i(v,q)R_<(q)w
       =Y_i(u_{-1}v,q)w.
\]
This equality is coefficientwise: in the first term the lower bound of $Y_i(v,q)w$ bounds $j$ at each fixed power; in the second $u_jw=0$ for $j\gg0$.

At a marking of inertia order $N_i$, write $g_iu=e^{2\pi i\alpha}u$, $0\le\alpha<1$, and fix its root/physical branches. For $K\ge1$ set $b_{K,j}=\binom{K+j-1}{j}$. The equivariant coefficient $z^\alpha(z-q)^{-K}$ has expansions
\begin{align*}
 \iota_{z,q}\bigl(z^\alpha(z-q)^{-K}\bigr)
  &=\sum_{j\ge0}b_{K,j}q^jz^{\alpha-K-j},\\
 \iota_{q,z}\bigl(z^\alpha(z-q)^{-K}\bigr)
  &=(-1)^K\sum_{j\ge0}b_{K,j}q^{-K-j}z^{\alpha+j}.
\end{align*}
Thus $R_>(q)=\sum_{j\ge0}b_{K,j}q^ju_{\alpha-K-j}$, $R_<(q)=(-1)^K\sum_{j\ge0}b_{K,j}q^{-K-j}u_{\alpha+j}$. The coefficient identity \eqref{eq:twistedBorcherds}, now with $m=\alpha$, $p=-K$, yields
\begin{align*}
 &R_>(q)Y_i^{g_i}(v,q)w
       -Y_i^{g_i}(v,q)R_<(q)w\\
 &\qquad=\sum_{\ell\ge0}\binom{\alpha}{\ell}
       q^{\alpha-\ell}Y_i^{g_i}(u_{\ell-K}v,q)w.
\end{align*}
Indeed $(-1)^j\binom{-K}{j}=b_{K,j}$ and $(-1)^{-K}=(-1)^K$. The sum on the right is finite by $u_{\ell-K}v=0$ for $\ell\gg0$. In the first term, its fixed $q$-coefficient uses only finitely many $j$ by lower
truncation of $Y_i^{g_i}(v,q)w$; in the second term, $u_{\alpha+j}w=0$ for $j\gg0$. Substitution $q=\tau^{N_i}$ makes each asserted coefficient a finite sum in the root Laurent series. Integral Laurent factors and the
finite positive-Virasoro corrections in the calibrated current gauge give the same assertion by linearity and finite shifts.

Apply these identities to each finite state-frame coefficient of the meromorphic current family from (i). The coefficients in its incoming and moving-point Laurent expansions have common
lower bounds after a single parameter-power multiplication. For a fixed parameter coefficient, only finitely many current modes and input vectors occur. Since the original functional annihilates $\mathscr G_{\bullet,k}T_\bullet$ \emph{on every
algebraic input tensor}, the transferred expression vanishes. In the physical/root normalization of
\eqref{eq:normalizednoderesidue}, the $N_i$ free lifts
approaching one inertia orbit cancel the factor $1/N_i$. We obtain the identity of complete root Laurent germs
\begin{equation}\label{eq:propagatedwardidentity}
 (P\lambda)(R_p(\sigma)v;w)
       +\sum_i(P\lambda)(v;R_i(\sigma)w)=0.
\end{equation}
Both sides are meromorphic in $p$ on the indicated finite cover. Their common zero germ at an old incoming orbit implies the identity everywhere on that connected component. Thus $P\lambda(\cdot,p;\cdot)$ annihilates every additional-point
Ward relation of node level $k$; it is again an element of the corresponding partial Ward dual. The argument uses only the existing twisted \emph{incoming} modules, not a conjectural node module.

\noindent(iii)\ 
For a fixed list of algebraic state vectors and a finite state cutoff, lower truncation gives finite pole bounds at each old marking and at every added free orbit; the associated divisors
are finite sums of the fixed markings and the graphs $\Gamma_h$. On an $H$-free parameter open they are relative Cartier divisors. Apply the finite-jet criterion at the end of Lemma~\ref{lem:partialequivariantresidue} to this specified
finite-rank bundle, after enlarging the incoming divisor by $mP_\bullet$. The $H^1$-vanishing in (i), flatness and base change apply. The prescribed jets are obtained from finitely many actual incoming module coefficients and are holomorphic in
the parameters. Injectivity of the finite jet map and pointwise existence therefore give a unique holomorphic family. The same argument works after adjoining further copies of the genuine vacuum module $V$ at free orbits, without changing $k$.

\noindent(iv)\ 
Fix $p_2$ and consider the pointed component with the additional
actual incoming $V$-module at $p_2$. Applying the now
established one-point propagation to $u$ at $p_1=p_2+\xi$ prescribes, \emph{by the definition of its incoming germ},
\[
 Y_V(u,\xi)v=\sum_{m\in\mathbb Z}
           (u_mv)\xi^{-m-1}.
\]
For each $m\in\mathbb Z$, its propagated scalar coefficient is therefore exactly $P\lambda(u_mv,p_2;w)$. This statement includes $m<0$; in particular, it does \emph{not} follow only
from the singular terms $m\ge0$. The finite-rank residue theorem of Lemma~\ref{lem:partialequivariantresidue}, applied after the additional insertion, identifies this complete formal germ with
the Laurent expansion of an actual meromorphic dual section. Since $u_mv=0$ for $m\gg0$, its principal part is finite, and its regular Taylor part converges normally on every smaller disk not meeting another marked orbit. This proves
\eqref{eq:partialpropagationope}, with all coefficients and the
asserted convergence, and does not presuppose that OPE as a prior analytic identity. Other lifts are obtained by $H$-transport.

We prove the uniqueness used to interchange the insertions. At an added unramified free point $p$, prescribe the local current principal part $u\,d\xi/\xi$ and outgoing
vanishing $-kQ_\bullet$. By the finite-rank principal-parts sequence and $H^1(F_{d,\bullet}(mP_\bullet-kQ_\bullet))=0$ for $m\gg0$, this part lifts to an $H$-equivariant global current, with other poles only at the old incoming points.
On the vacuum, its residue is $u_{-1}\mathbf1=u$, whereas the regular local terms act by $b_j\mathbf1=0$ for $j\ge0$. Hence a partial Ward functional with a free $V$-insertion is determined by its restriction to
$\mathbf1$ at that insertion. Restricting a partial Ward functional to $\mathbf1$ deletes the free insertion: currents regular at that point act trivially on the vacuum, and the old Ward relations remain. Two iterated propagations in opposite
orders are partial Ward functionals of level $k$, and their restriction to both vacua is $\lambda$. Apply the just-proved uniqueness successively at the two free orbits; they coincide.
Joint holomorphy from (iii) identifies the two germs with the same meromorphic continuation on the collision-free locus.
\end{proof}

For a pair of states define their singular OPE order and the finite orbit bound by
\[
\begin{split}
p(u,v)&=\max\bigl(\{0\}\cup\{j+1:j\geq0,\ u_jv\neq0\}\bigr),\\
 \overline M_H(u,v)&=\max\bigl(\{1\}\cup
                 \{p(\gamma u,\delta v):\gamma,\delta\in H\}\bigr).
\end{split}
\]
The maximum is finite by state truncation and finiteness of $H$. We use $\overline M_H$ to bound the transported diagonal poles; the multipliers in the operator condition remain the original $M_{ij}$.

\begin{lemma}\label{lem:actualcomponentcorner}
Let $\pi_i:C_i\to\mathbb P^1$ be finite morphisms with $C_i$ smooth projective, put $C=\prod_{i=1}^nC_i$, and write $\C(C)$ for its rational-function field. Fix $f\in\C(C)$ and $P\in\C[z_1,\ldots,z_n]$, and set
$F(z)=P(z)f$ on the branches of $\pi=\prod_i\pi_i$. For $A\in\operatorname{GL}_n(\C)$, $\beta\in\C^n$, put
\[
 \eta=zA-\beta,\qquad z=(\eta+\beta)A^{-1},\qquad
 \xi_i=\begin{cases}\eta_i,&\delta_i=0,\\
                    \eta_i^{-1},&\delta_i=\infty,
       \end{cases}\quad\delta\in\{0,\infty\}^n.
\]
Write $\C\{\xi\}=\C\{\xi_1,\ldots,\xi_n\}$ for the ring of convergent power-series germs at $\xi=0$. If $\eta=\delta$ is component-isolated for the \emph{full set of branches} of $F$, then there exist
$N_i\in\Z_{>0}$, $b_i\in\Z_{\geq0}$, and, for each branch $a$, $H_{a,\mathbf r}\in\C\{\xi_1,\ldots,\xi_n\}$, $0\leq r_i<N_i$, such that
\begin{equation}\label{eq:actualcornerexpansion}
 F^a\bigl((\eta+\beta)A^{-1}\bigr)
 =\sum_{\substack{0\leq r_i<N_i\\1\leq i\leq n}}
   \prod_{i=1}^n\xi_i^{(r_i-b_i)/N_i}
                  H_{a,\mathbf r}(\xi).
\end{equation}
The series in \eqref{eq:actualcornerexpansion} converge absolutely on sufficiently small cut punctured polydisks. Thus $\eta=\delta$ is regular in the sense of
\cite[Definition~2.6]{DH}, with logarithm degree $0$ and common
finite exponent sets for all branches. The same conclusion holds for any finite family and, branchwise, for $c_aF^a$ with $c_a\in\C$. The constants need not be compatible with analytic continuation, so this last assertion
concerns the local expansions only. The common exponent sets may depend on the chosen family and on $P$.
\end{lemma}
\begin{proof}
Put $\Delta^\times=\{\xi:0<|\xi_i|<\epsilon_i\}$, shrinking $\epsilon$ so that every branch is holomorphic on $\Delta^\times$. By GAGA \cite{SerreGAGA}, $f$ is algebraic; therefore $F$ satisfies a monic equation over $\C(\xi)$,
\begin{equation}\label{eq:cornerpolynomial}
 F^d+\sum_{j=0}^{d-1}a_j(\xi)F^j=0,
 \qquad a_j\in\C(\xi_1,\ldots,\xi_n).
\end{equation}
Here $\eta_i=\xi_i$ or $\xi_i^{-1}$, so the affine-linear substitution does not change algebraicity. The coefficients $a_j$ are elementary symmetric functions of all conjugate branches of $F$, counted with multiplicity. Since these
branches are holomorphic on $\Delta^\times$, the $a_j$ have no poles there. After shrinking $\Delta$, the polar divisors of the rational functions $a_j$ are consequently supported on
$\bigcup_i\{\xi_i=0\}$. Hence, for some $\mathbf m\in\Z_{\geq0}^n$,
\[
 \xi^{\mathbf m}a_j\in\C\{\xi\}
 \quad(0\leq j<d),\qquad
 \xi^{\mathbf m}=\prod_i\xi_i^{m_i}.
\]
The analytic continuation of the finitely many branches gives finite monodromy on $\Delta^\times$, whose fundamental group is $\Z^n$. Choose $N_i\geq1$ divisible by the order of the monodromy of its $i$th standard generator. The substitutions
\[
 \xi_i=s_i^{N_i}\qquad(1\leq i\leq n)
\]
then make every branch single-valued and holomorphic on the punctured $s$-polydisk. Equation~\eqref{eq:cornerpolynomial}, the preceding coefficient bounds, and the elementary bound for the roots of a monic polynomial give a common
$\mathbf b\in\Z_{\geq0}^n$ such that
\[
 \sup_{0<|s_i|<\epsilon_i^{1/N_i}}
 \left|s^{\mathbf b}F^a(s_1^{N_1},\ldots,s_n^{N_n})\right|
 <\infty\qquad\text{for all branches }a.
\]
(The choice $b_i=N_im_i$ is sufficient after enlarging the $m_i$ if necessary.) The bounded functions extend across the coordinate axes by the removable-singularity theorem, applied successively, and are holomorphic on the whole $s$-polydisk.
Write their convergent Taylor series as
\[
 s^{\mathbf b}F^a(s^{\mathbf N})
   =\sum_{\mathbf q\in\Z_{\geq0}^n}c_{a,\mathbf q}s^{\mathbf q}.
\]
Group the exponents uniquely as $q_i=N_i\ell_i+r_i$, $0\leq r_i<N_i$. Then
\[
 H_{a,\mathbf r}(\xi)
  =\sum_{\boldsymbol\ell\in\Z_{\geq0}^n}
       c_{a,\mathbf N\boldsymbol\ell+\mathbf r}
       \xi^{\boldsymbol\ell}\in\C\{\xi\},
\]
which proves \eqref{eq:actualcornerexpansion} and its absolute convergence. At $\delta_i=\infty$, $\xi_i=\eta_i^{-1}$ turns $\xi_i^{(r_i-b_i)/N_i}\C\{\xi_i\}$ into a series with exponents in $(b_i-r_i)/N_i-\Z_{\geq0}$, as required by
\cite[Definition~2.6]{DH}. Thus one may take the common
exponent set
\[
 \Lambda=\prod_{i=1}^n
  \left\{\frac{r-b_i}{N_i}:0\leq r<N_i\right\},
 \qquad L=0.
\]
There are finitely many branches and finitely many functions in a specified family. Taking common multiples of the $N_i$ and increasing the $b_i$ gives common exponent sets; multiplication
by $c_a\in\C$ changes only the coefficients. No uniform bound over unbounded families of polynomials $P$ is asserted.
\end{proof}

\begin{proposition}
\label{prop:partialallfields}
Iterating the preceding construction gives $P^n\lambda$ for every finite number of free $V$ inputs. For fixed states and incoming vectors these are meromorphic finite-state dual sections on the product of the proper component covers, with poles only at incoming orbits, outgoing orbits and orbital diagonals. The poles along a $V$--$V$ diagonal have integer order. The original local ordered expansions converge absolutely on their prescribed radial chambers.

For a full three-point Ward functional of actual ordinary modules, the homogeneous inverse
\eqref{eq:singlewardinverse} therefore satisfies the complete
additional-field requirement for complete twisted intertwining operators \cite[Definition~2.8(3)]{DH}.

\end{proposition}
\begin{proof}
First iterate the \emph{pointwise} Ward construction of Lemma~\ref{lem:partialpropagationward}. It uses only actual incoming modules, including each added copy of $V$, and never a not-yet-constructed node module. Its vacuum-uniqueness argument makes these pointwise propagations independent of their order. For a full three-point functional take $Q$ empty and include the actual contragredient module at infinity among the incoming factors.

We prove parameter holomorphy by induction on the number $n$ of free fields. Fix $\lambda\in\mathscr D_k$, all original input vectors, the states at the first $n-1$ free points, and a finite cutoff $V_{\leq d}$ for the last state. Let $B$ be an algebraic parameter open where those first $n-1$ orbits are disjoint from one another and from $P,Q$. In the last variable the pole bounds are already fixed, independently of $b\in B$. At an original incoming orbit this is lower truncation on the fixed incoming vector, uniformly on a finite basis of $V_{\leq d}$ and its finite coordinate corrections. At a previously added free orbit the incoming module is $V$ with its fixed vector, so the same assertion follows from state truncation. More explicitly, if the last test state has weight at most $d$ and an earlier vacuum input has homogeneous weight $d_j$, then its $m$-th mode on that fixed input vanishes for $m\geq d+d_j$ by the nonnegative grading of $V$. Thus an integer $M_j=\max(1,d+d_j)$ bounds the corresponding OPE pole, uniformly as the earlier free positions vary on the collision-free base. At an original incoming orbit one takes the maximum of the finitely many lower-truncation bounds on the fixed incoming vector over a homogeneous basis of $V_{\leq d}$ and over the finitely many root/deck and positive-mode coordinate corrections. These operations lower the state weight and preserve finite dimensionality, so the resulting integer $c_i$ is independent of the parameter values. At $Q$ it is precisely the dual bound of the fixed balanced level $k$. Include all $H$-translates. Thus there are integers $c_i,M_j\geq0$, depending on these fixed data but not on the free positions, such that the last-variable dual section belongs fiberwise to
\[
\mathcal A_{d,k}=\operatorname{pr}_1^*\bigl(
 F_{d,\bullet}^{\vee}\otimes\omega_{\widetilde C_\bullet}(kQ_\bullet)\bigr)
 \otimes\cO\left(\sum_i c_i(P_i\times B)
                  +\sum_{j<n}\sum_{\gamma\in H}M_j\Gamma_{j,\gamma}\right).
\]
Here $F_{d,\bullet}$ is the finite-rank current bundle of Lemma~\ref{lem:partialequivariantresidue}, and $F_{d,\bullet}^{\vee}\otimes\omega(kQ_\bullet)$ is precisely its node-filtered dual bundle in~\eqref{eq:partialnodepole78}; the
node bound is therefore included once, not added as an unspecified extra pole. The graphs are $x=\gamma p_j$. After restricting to the $H$-free collision-free parameter open $B$,
they are disjoint relative Cartier divisors and avoid $P_\bullet\cup Q_\bullet$. Consequently $\mathcal A_{d,k}$ is a finite-rank vector bundle, flat over $B$, and the incoming divisor $P_\bullet\times B$ is
relatively ample on every connected component.

The preceding bounds show, before choosing any auxiliary $m$, that the \emph{pointwise} section belongs to $H^0(\widetilde C_\bullet,(\mathcal A_{d,k})_b)$ for every $b\in B$: the original incoming poles are bounded by the $c_i$, the moving poles by the $M_j$, and the node order by $k$. Hence its membership in $H^0(\widetilde C_\bullet,(\mathcal A_{d,k}(mP_\bullet))_b)$ for every $m\geq0$ follows from the established pointwise bounds. Apply the finite-jet criterion at the end of Lemma~\ref{lem:partialequivariantresidue} to $\mathcal A_{d,k}$: increase the incoming twist by $mP_\bullet$ so that the required relative $R^1\operatorname{pr}_{B*}$ vanishes,
and choose a finite incoming jet order $r$ for which the resulting jet map is fiberwise injective. After shrinking $B$ if necessary, cohomology and base change make its source locally free, and its
image is a holomorphic subbundle. Each specified incoming jet is a previous $(n-1)$-point propagation evaluated on one of finitely many actual mode vectors: the state cutoff, the fixed pole orders
and $r$ make this list finite. These jets are holomorphic by induction. Pointwise existence from Lemma~\ref{lem:partialequivariantresidue} places the jet family in the image subbundle. Inverting the jet map onto its image gives a
holomorphic section whose complete incoming germs coincide, by pointwise uniqueness, with the pointwise propagation. This proves joint holomorphy away from the stated divisors. It uses only pole
bounds in the \emph{last} free variable, established at this induction stage; the polar bounds in the earlier variables are obtained below by permuting them into the last position.

Now use the already established pointwise permutation symmetry to put any chosen variable last. Its incoming and outgoing axis orders have the same position-independent bounds. At a pair diagonal the full Laurent OPE of Lemma~\ref{lem:partialpropagationward} gives the integer bound $\overline M_H(u_i,u_j)$ in the uncorrected physical frames. For the finite positive-Virasoro frame corrections, one must not assert that this same minimum OPE order remains sufficient. If $u_i\in V_{d_i}$ and $u_j\in V_{d_j}$ are homogeneous, all states appearing in their finite coordinate/deck corrections have weights at most $d_i$ and $d_j$, respectively. Since $V_{<0}=0$, their $m$-th products vanish for $m\geq d_i+d_j$. Thus the explicit auxiliary bound
\[
 M_{ij}^{\mathrm{aux}}:=\max\bigl(1,\overline M_H(u_i,u_j),d_i+d_j\bigr)
\]
controls every orbit-pair pole in these frames; it is fixed independently of the free positions. It is an auxiliary convergence bound and does \emph{not} replace any exponent $M_{ij}$ prescribed in Definition~\ref{def:completeintertwiner}(iii). The negative OPE terms are finite and their coefficients are lower-point propagations, already holomorphic by induction. There are only finitely many axis and orbital-pair divisors. Choose their fixed maximum orders for the data. Multiplication by these divisor equations removes every pole at its generic point and extends the section away from intersections. The remaining locus has codimension at least two in the smooth product of curves; Hartogs extension in local vector-bundle frames removes it. Consequently $P^n\lambda$ is meromorphic on the whole proper product with one fixed finite polar divisor. GAGA then makes its finite-state scalar evaluations algebraic \cite{SerreGAGA}. All constants here may depend on $n,k$, the finite state cutoffs and the fixed algebraic vectors. They do not depend on the free positions; no bound uniform over all vectors or all $n$ has been used.

We check convergence of the original ordered sums. On a marked disk containing no other marked orbit, use a common root coordinate for its cyclic inertia. After clearing the axis orders and pair diagonal factors, the scalar is holomorphic on the product of these root disks. Its multiple Taylor series converges absolutely by Cauchy's estimate. For a normalized full three-point block at $\infty,1,0$, the output disk has physical coordinate $x=z^{-1}$ and radius $|x|<1$. A common root substitution $x_i=s_i^N$ trivializes all its finite covering branches. Clear the diagonal factors using the finite auxiliary orbit bounds $M_{ij}^{\mathrm{aux}}$. The resulting function $H(s)$ has only coordinate poles, so $s^bH(s)$ is holomorphic on $|s_i|<1$. For smaller polyradii $r_i<\rho_i<1$, put $M=\sup_{|s_i|=\rho_i}|s^bH(s)|$. Its Taylor coefficients satisfy
\[
 |A_{\boldsymbol n}|\leq M\prod_i\rho_i^{-n_i},
 \qquad
 \sum_{\boldsymbol n\geq0}|A_{\boldsymbol n}|\prod_i r_i^{n_i}
 \leq M\prod_i(1-r_i/\rho_i)^{-1}.
\]
The finitely many negative powers and root exponent classes give a convergent multiple Laurent series before restoring the pair factors.

On the chamber $|z_1|>\cdots>|z_r|>1$ one has $|x_1|<\cdots<|x_r|<1$. For $i<j$, restoring a cleared pole of integer order $M_{ij}^*=M_{ij}^{\mathrm{aux}}$ uses
\[
 (z_i-z_j)^{-M_{ij}^*}
 =x_i^{M_{ij}^*}\sum_{a\geq0}
   \binom{M_{ij}^*+a-1}{a}(x_i/x_j)^a.
\]
Each geometric series converges absolutely there. On a compact subchamber, cover by finitely many small product neighborhoods on which an upper bound for $|x_i|$ is smaller than a lower bound for $|x_j|$ whenever $i<j$. The preceding Cauchy estimate and the geometric bounds then give absolute normal convergence of their finite product with the Laurent series. For a fixed resulting exponent, the coefficient convolution is finite, not merely conditionally summable. To see this already for two variables, suppress a fixed Kummer residue class and write $H=\sum_{m_1\geq-b_1,\,m_2\geq-b_2}h_{m_1,m_2}x_1^{m_1}x_2^{m_2}$. Restoring one factor of order $M$ makes the coefficient of $x_1^{e_1}x_2^{e_2}$ equal to
\[
 \sum_{0\leq a\leq e_1+b_1-M}
   \binom{M+a-1}{a}h_{e_1-M-a,\,e_2+a},
\]
where an empty range contributes zero. On the Kummer cover each fixed residue class gives the identical bound after dividing the corresponding integral exponents by the root order. With more variables, the smallest-index variable bounds all geometric indices leaving it; inductively each subsequent variable bounds its outgoing indices after those entering it have been fixed. Thus each target multidegree has finitely many contributors. This is an ordered finite coefficient calculation. The incoming and exchanged chambers use the same calculation in their centered coordinates.

Finally compare these coefficients with the prescribed weight sums. For homogeneous inserted states $u_1,\ldots,u_r$ and final dual weight $\lambda_Z$, let $n_i$ be the mode index of the $i$th field in the outside-to-inside order. Its exponent in $z_i$ is $-n_i-1$, and its degree is $\wt u_i-n_i-1$. The weight immediately after the first $j$ outside fields have been removed is $\mu_j=\lambda_Z-\sum_{i=1}^j(\wt u_i-n_i-1)$, for $1\leq j\leq r$. This triangular formula determines all projected intermediate weights from the mode indices, and its differences recover every $n_i$. Hence one Laurent coefficient is precisely the finite-dimensional contraction at one original weight tuple. The prescribed germs, the complete OPE and the finite coefficient convolution identify it with that contraction before convergence is used. Absolute normal convergence of this Laurent expansion is therefore absolute normal convergence of the original sums over $\mu_1,\ldots,\mu_r$. Finite simple-block decompositions add finitely many such estimates.

For homogeneous $u_i\in V_{d_i}$, $a\in X_\alpha$, $b\in Y_\beta$, and $z'\in (Z')_\gamma$, the original $L(0)$-commutator formula gives, on compatible logarithm branches,
\[
 f_k(z_1,\ldots,z_k)
 =z_k^{\gamma-\alpha-\beta-\sum_{i<k}d_i}
   f_k(z_1/z_k,\ldots,z_{k-1}/z_k,1).
\]
The power of $z_k$ is retained as a fixed complex power; it is not claimed to be algebraic in $z_k$. For each fixed $z_k\in\C^*$, the dependence on $z_1,\ldots,z_{k-1}$ is an
algebraic multivalued scalar on the finite product cover, with the specified branches. The preceding estimate yields absolute normal convergence of the original sums in $|z_1|>\cdots>|z_k|>0$.

Now fix \emph{any} $(M_{ij})_{1\le i<j<k}$ allowed by
\cite[Definition~2.8(3)]{DH}, so that
$z^{M_{ij}}Y_V(u_i,z)u_j\in V[[z]]$. The product $P_M(z)=\prod_{1\le i<j<k}(z_i-z_j)^{M_{ij}}$ is a polynomial in the free positions. The scalar $P_M(z)f_k(z_1,\ldots,z_{k-1},z_k)$, for fixed $z_k$,
remains algebraic on that finite cover, irrespective of whether the chosen exponents equal the auxiliary pole bounds $M_{ij}^*=M_{ij}^{\mathrm{aux}}$ used in the convergence proof.
Lemma~\ref{lem:actualcomponentcorner} therefore applies to every component-isolated corner of the \emph{complete collection of its branches}. No replacement of the prescribed $M_{ij}$ by larger integers is made in this regularity argument.

The finite set of deck branches and the finitely many transported state vectors allow one to take common multiples of the ramification indices and common pole bounds in Lemma~\ref{lem:actualcomponentcorner}. Thus the exponent sets and
logarithm degree $0$ are common to all these branches for the
fixed inputs and multipliers. Changing the branch of the
fixed power of $z_k$ multiplies its value by a scalar and does not alter those exponents in $z_1,\ldots,z_{k-1}$. Finite linearity extends the assertion to arbitrary algebraic vectors. This proves
the original radial absolute convergence and the full quantifier over permitted $M_{ij}$ and component-isolated corners in
\cite[Definition~2.8(3)]{DH}; the derivative axiom was proved
independently in Proposition~\ref{prop:singlewardconverse}.
\end{proof}

\begin{corollary}
\label{cor:completewardconverse}
Under the standing hypotheses on $V$ and $G$, the three-point Ward space and the intertwining-operator spaces satisfy $\mathfrak W_{XY}^{Z} \cong\widetilde{\cI}\binom Z{X\ Y} =\cI\binom Z{X\ Y}.$ \end{corollary}
\begin{proof}
Proposition~\ref{prop:singlewardconverse} proves
$\operatorname{ev}_1:\widetilde{\cI}\binom Z{X\ Y}
\xrightarrow{\sim}\mathfrak W_{XY}^Z$ and the independent
$L(-1)$-derivative identity. By Proposition~\ref{prop:partialallfields}, its inverse satisfies item~(3) of Definition~\ref{def:completeintertwiner}; it already
satisfies items~(1), (2), (4) by definition. Hence $\widetilde{\cI}\binom Z{X\ Y}\subseteq\cI\binom Z{X\ Y}$. The reverse inclusion is the evaluation-to-Ward map established in Proposition~\ref{prop:binaryward} and
Proposition~\ref{prop:singlewardconverse}, proving the claimed isomorphism.
\end{proof}

\subsection{Residue identities and the full node action}
Fix a node label $g$ of order $N$, a based root $z=s^N$, and $u\in V_d$ with $gu=e^{2\pi i r/N}u$ ($0\le r<N$). Put
$\mathcal P_{\lambda,u}(z;w)
=P\lambda(\iota_z(u),z;w)$. For
$m\in r/N+\Z$, define
\begin{equation}\label{eq:partialnodemode}
 \begin{aligned}
 (u_m\lambda)(w)
 &:=\Res_z z^m\mathcal P_{\lambda,u}(z;w)\,dz\\
 &:=\frac1N\Res_s
       s^{Nm}\mathcal P_{\lambda,u}(s^N;w)
       Ns^{N-1}\,ds,
 \qquad \lambda\in\mathscr D,
 \quad w\in T_\bullet .
 \end{aligned}
\end{equation}
The second line specifies the residue of a fractional-power series; its normalization is \eqref{eq:normalizednoderesidue}. In particular, if
$\mathcal P_{\lambda,u}(z;w)
 =\sum_{j\in r/N+\Z}c_jz^{-j-1}$,
then $(u_m\lambda)(w)=c_m$. The monodromy convention
$\mathcal P_{\lambda,u}(e^{2\pi i}z;w)
 =\mathcal P_{\lambda,g^{-1}u}(z;w)$
forces these indices. Extend \eqref{eq:partialnodemode} linearly in $u$ and write $Y_{\mathscr D}(u,z)=\sum_m u_mz^{-m-1}$; the action on $\mathscr D$ and lower truncation are established in Theorem~\ref{thm:partialnodeweakmodule}.

\begin{theorem}\label{thm:partialnodeweakmodule}
The maps \eqref{eq:partialnodemode} are independent of the finite-state cutoff and compatible with the specified root-coordinate transport. For homogeneous $u\in V_d$ and every $k\ge1$,
\begin{equation}\label{eq:nodelevelbounds713}
 \begin{gathered}
 u_m\mathscr D_k=0
 \quad\left(m\in\frac rN+\Z,\quad
             m\ge d-1+\frac{k}{N}\right),\\
 u_m\mathscr D_k\subseteq
 \mathscr D_{k+\max\{0,N(d-m-1)\}}
 \quad\left(m\in\frac rN+\Z\right).
 \end{gathered}
\end{equation}
Consequently $Y_{\mathscr D}$ makes $\mathscr D$ a weak $g$-twisted $V$-module. More precisely, for eigenvectors $u,v$, $m\in r_u/N+\Z$, $n\in r_v/N+\Z$ and $p\in\Z$,
\begin{equation}\label{eq:nodeborcherds713}
 \begin{aligned}
 &\sum_{j\ge0}\binom mj
       (u_{p+j}v)_{m+n-j}\\
 &\quad=\sum_{j\ge0}(-1)^j\binom pj
       \left(u_{m+p-j}v_{n+j}
           -(-1)^pv_{n+p-j}u_{m+j}\right)
 \end{aligned}
\end{equation}
(as an identity on $\mathscr D$); the sums are finite on every fixed vector. Moreover $Y_{\mathscr D}(\mathbf1,z)=\id_{\mathscr D}$ and
$Y_{\mathscr D}(L(-1)u,z)
 =\partial_zY_{\mathscr D}(u,z)$.
In particular, for $k=N(a+1)$, $a\in\N$, $u_m\mathscr D_k=0$ whenever $m\ge d+a$.
\end{theorem}
\begin{proof}
All assertions are verified first on a fixed $\mathscr D_k$. For $\lambda\in\mathscr D_k$ the $k$-vanishing condition on the balanced current frame $e_u=\iota_z(u)z^{d-1}\,dz$, together with Lemma~\ref{lem:partialequivariantresidue}, gives the
root Laurent bound
\begin{equation}\label{eq:nodepolebound713}
 \mathcal P_{\lambda,u}(s^N;w)
 \in s^{1-k-Nd}\C[[s]],
\end{equation}
with only the character $-r$ terms retained. Thus the lowest possible root exponent is $1-k-Nd$. The residue in \eqref{eq:partialnodemode} selects $s^{-N(m+1)}$; for $m\ge d-1+k/N$, $-N(m+1)\le-k-Nd<1-k-Nd$.
This proves the first assertion of \eqref{eq:nodelevelbounds713} without using any node-module identity. It also gives the lower truncation of $Y_{\mathscr D}(u,z)\lambda$. Uniqueness of finite-rank propagation in
Lemma~\ref{lem:partialequivariantresidue}, the compatibility of its restrictions $V_{\le D}\subset V_{\le D'}$, and the based coordinate covariance imply the stated choice independence.

For stability, let $v\in V_e$ be a $g$-eigenvector, $n\in r_v/N+\Z$, and put
\[
 \delta=e-n-1\in\tfrac1N\Z,\qquad
 k'=k+\max\{0,N\delta\}\in\N_{>0}.
\]
Take $s\in\mathscr G_{k'}$. A homogeneous node term of $s$ is $u z^q\,dz$ with $u\in V_d$ and $q\ge d-1+k'/N$. When its free-point residue is moved across the insertion $v$, the full local OPE \eqref{eq:partialpropagationope} gives
\[
\Res_{\zeta=z}\zeta^q
        Y(u,\zeta-z)v\,d\zeta
   =\sum_{j\ge0}\binom qj z^{q-j}u_jv.
\]
The sum is finite since $u_jv=0$ for $j\gg0$. For each nonzero summand,
\[
 \begin{split}
 &(q+n-j)-\bigl(\wt(u_jv)-1\bigr)\\
 &\quad=q+n-d-e+2
 \ge \frac{k'}{N}-\delta\\
 &\quad=\frac{k}{N}+\max\{0,\delta\}-\delta
 \ge\frac{k}{N},
 \qquad \wt(u_jv)=d+e-j-1.
 \end{split}
\]
The first bound of \eqref{eq:nodelevelbounds713}, already proved for $\lambda\in\mathscr D_k$, therefore annihilates its node residue; the mode-index congruence follows from $g(u_jv)=e^{2\pi i(r_u+r_v)/N}u_jv$.
Incoming residues vanish by the complete one-free-point Ward identity \eqref{eq:propagatedwardidentity}. For an arbitrary $s$ one works in its finite state cutoff: only finitely many node terms contribute to the fixed residue,
by \eqref{eq:nodepolebound713}, and the incoming fields are lower truncated on the fixed algebraic tensor. Hence
\[
 (v_n\lambda)(R_\bullet(s)w)=0
 \quad(s\in\mathscr G_{k'},\ w\in T_\bullet),
 \qquad v_n\lambda\in\mathscr D_{k'}.
\]
This proves the second assertion of
\eqref{eq:nodelevelbounds713}, and therefore defines
$Y_{\mathscr D}$ on the algebraic union $\mathscr D$.

We prove the two iterated mode identities directly from the original
propagation, before using the twisted Jacobi identity. For
$\lambda\in\mathscr D_k$ set
\[
 \mathcal P_{u,v}(\zeta,z;w)
 :=P^2\lambda(\iota_\zeta(u),\zeta;
                  \iota_z(v),z;w).
\]
Lemma~\ref{lem:partialpropagationward} gives this particular double
propagation with both added modules equal to the genuine untwisted
vacuum module $V$. Proposition~\ref{prop:partialallfields} supplies its
joint meromorphicity on the actual rooted configuration cover, not an
arbitrarily chosen solution of its Ward equations. In the physical
coordinate $\xi=\zeta-z$ its \emph{complete} diagonal Laurent germ is
\begin{equation}\label{eq:nodephysicalope713}
 \mathcal P_{u,v}(z+\xi,z;w)
 =\sum_{j\in\Z}\xi^{-j-1}
          P\lambda(\iota_z(u_jv),z;w),
 \qquad 0<|\xi|\ll|z|.
\end{equation}
This is the complete incoming germ at the second free $V$-insertion
in~\eqref{eq:partialpropagationope}; in particular $j<0$ is not
discarded. All state-coordinate corrections are transported through
the actual based physical frames before this formula is read. One must
not identify two arbitrary root trivializations when taking its
coefficients.

We verify the first nested residue. Fix an allowed $n$ and define on
the $|\zeta|>|z|$ chamber the coefficientwise single-variable object
\[
 Q_n(\zeta;w):=\Res_z z^n
             \mathcal P_{u,v}(\zeta,z;w)\,dz,
\]
where every fractional residue means exactly the normalized rooted
residue of~\eqref{eq:normalizednoderesidue}. This is well-defined:
for a fixed finite state cutoff, the double propagation is meromorphic
with finite polar orders at both axes and at their orbital diagonals.
On the chosen $|z|<|\zeta|$ annulus, the expansion of each diagonal
pole is a power series in $z/\zeta$ with a finite lower $z$-exponent.
For a specified $z$-exponent, only finitely many terms from each
such expansion and the lower-bounded $z$-Laurent series can contribute.
Consequently $Q_n$ is a meromorphic rooted section in $\zeta$ with
poles only at the original incoming orbits and at the node; it has
some finite node pole order $K(n,u,v,w,d)$. We do not claim that
this order is the original $k$.

Let $p_i$ be any old incoming marking. On disjoint rooted disks about
$p_i$ and the node, expand the prescribed complete incoming germ of
$\mathcal P_{u,v}$ against every local state $b\in V_{\le d}$.
Normal convergence on smaller disjoint annuli (or the corresponding
finite coefficient calculation just described) gives, for every
incoming mode coefficient, the identity
\begin{align*}
 &[\text{incoming mode}]\,Q_n(\zeta;w)\\
 &\quad=\Res_z z^n
 [\text{incoming mode}]\,\mathcal P_{u,v}(\zeta,z;w)\,dz\\
 &\quad=[\text{incoming mode}]\,
          P(v_n\lambda)(\iota_\zeta(u),\zeta;w).
\end{align*}
The last equality is the \emph{definition}
$ (v_n\lambda)(w_i')=\Res_z z^nP\lambda(v,z;w_i')\,dz$
for each algebraic vector $w_i'$ obtained by applying the actual old
incoming mode to $w$. It uses no node-module Jacobi identity.
The stability assertion already proved above gives
$v_n\lambda\in\mathscr D_{k'}$ with
$k'=k+\max\{0,N(\wt(v)-n-1)\}$. Choose a finite $K$ at least this $k'$
and the actual pole order of $Q_n$ at the node (and increase the state
cutoff if necessary). Both $Q_n$ and the independently constructed
$P(v_n\lambda)(u,\zeta;w)$ then belong to the same finite-rank
dual bundle $F_{d,\bullet}^{\vee}\otimes\omega(KQ)$ and have
\emph{identical complete incoming germs at every old marking}.
Uniqueness in Lemma~\ref{lem:partialequivariantresidue}, component by
component, gives equality of these meromorphic rooted sections.
Taking the outer $\zeta^m$ coefficient proves the first identity below.
The same argument with the two variables interchanged, first
extracting $\zeta^m$ in $|\zeta|<|z|$, uses
$u_m\lambda\in\mathscr D_{k+\max\{0,N(\wt(u)-m-1)\}}$
and proves the second identity. Thus
\begin{equation}\label{eq:nestednodemodes713}
 \begin{aligned}
 u_mv_n\lambda
  &=\Res_{|\zeta|>|z|}\Res_z
       \zeta^m z^n\mathcal P_{u,v}(\zeta,z;\cdot)
       \,dz\,d\zeta,\\
 v_nu_m\lambda
  &=\Res_{|z|>|\zeta|}\Res_\zeta
       \zeta^m z^n\mathcal P_{u,v}(\zeta,z;\cdot)
       \,d\zeta\,dz .
 \end{aligned}
\end{equation}
Here and throughout, the two residues are successive
\emph{normalized rooted} residues, not residues of an undefined
infinite-dimensional tensor product. All comparisons are made after
fixing a test vector, a finite state cutoff, and the two specified
mode indices. This also proves the compatibility of nested extraction
with the selected root paths and strictly typed deck pulls.

We compute the diagonal normalization rather than suppressing the
root branches. Put $\epsilon=e^{2\pi i/N}$,
$\zeta=s_1^N$, $z=s_2^N$, and write
$gu=\epsilon^{r_u}u$, $gv=\epsilon^{r_v}v$.
The root rotations in the two free variables multiply the propagated
scalar by $\epsilon^{-r_u}$ and $\epsilon^{-r_v}$, respectively.
The corresponding factors from $s_1^{Nm}$ and $s_2^{Nn}$ are
$\epsilon^{Nm}=\epsilon^{r_u}$ and
$\epsilon^{Nn}=\epsilon^{r_v}$; the Jacobians
$N s_i^{N-1}ds_i$ are invariant under a full root rotation.
Thus the \emph{paired} two-form
$\zeta^mz^n\mathcal P_{u,v}(\zeta,z;w)\,d\zeta\,dz$
is single valued around either physical coordinate axis, on the
specified analytic branch. Its physical diagonal $\zeta=z$
has exactly the root preimages
$s_1=\epsilon^b s_2$ ($0\le b<N$). Transport the $b$-th
local root germ to $b=0$ using the actual cyclic inertia action;
the state factor is $\epsilon^{-br_u}$ while the mode multiplier
is $\epsilon^{bNm}$, whose product is
\[
 \epsilon^{b(Nm-r_u)}=1.
\]
For example the pullback of $d\zeta/(\zeta-z)$ has residue $1$ at
each of the $N$ simple root diagonals:
\[
 \Res_{s_1=\epsilon^b s_2}
 \frac{N s_1^{N-1}\,ds_1}{s_1^N-s_2^N}=1.
\]
After the transported eigencharacters cancel, the $N$ local root
residues are equal. The one physical $\zeta$-residue is the sum of
these root residues divided by $N$; therefore its diagonal
coefficient is $N/N=1$. The second normalized residue in the $z$
variable similarly equals its physical residue. No additional $N$
or $1/N$ remains. The cyclic inertia calculation is local to each
node; arbitrary global deck arrows are composed in their actual
chronological order, with source and target transported sectors,
not commuted or silently identified.

Now fix $p\in\Z$ and pair the two-meromorphic propagation with
$(\zeta-z)^p\zeta^mz^n\,d\zeta\,dz$. On the respective radial
annuli the \emph{different} binomial expansions are
\[
 \begin{aligned}
 (\zeta-z)^p
 &=\sum_{j\ge0}(-1)^j\binom pj\zeta^{p-j}z^j,
 &&|z|<|\zeta|,\\
 (\zeta-z)^p
 &=(-1)^p\sum_{j\ge0}(-1)^j\binom pj
                       \zeta^jz^{p-j},
 &&|\zeta|<|z|.
 \end{aligned}
\]
Using the just-proved nested identities, their iterated residues are,
respectively,
\[
 \sum_{j\ge0}(-1)^j\binom pj
       u_{m+p-j}v_{n+j}\lambda,
 \qquad
 (-1)^p\sum_{j\ge0}(-1)^j\binom pj
       v_{n+p-j}u_{m+j}\lambda.
\]
The difference of these two radial residues equals the sum of the
physical diagonal residues by the residue theorem, with the
normalization computed above. Put $\xi=\zeta-z$ and apply the
\emph{complete} local OPE~\eqref{eq:nodephysicalope713}. The formal
expansion
$\zeta^m=(z+\xi)^m=\sum_{j\ge0}\binom mj z^{m-j}\xi^j$
gives exactly
\[
 \Res_z z^n\Res_{\xi=0}
   (z+\xi)^m\xi^p
   \sum_{q\in\Z}\xi^{-q-1}P\lambda(u_qv,z;\cdot)
    \,d\xi\,dz
 =\sum_{j\ge0}\binom mj(u_{p+j}v)_{m+n-j}\lambda.
\]
The equality is again coefficientwise. For fixed $p$ its diagonal
sum stops when $u_{p+j}v=0$; each radial sum stops because the
\emph{rightmost} mode $u_{m+j}$ or $v_{n+j}$ kills $\lambda$ for
$j\gg0$, by the already proved first bound in
\eqref{eq:nodelevelbounds713}. Finite cutoff and lower truncation
make all remaining scalar coefficients finite. Thus contour
comparison, equivalently the equality of three complete Laurent-germ
residues, proves~\eqref{eq:nodeborcherds713} for every integer $p$,
including $p<0$, and for every $w\in T_\bullet$.
Since the evaluations on all $w$ separate $T_\bullet^*$, this is
an operator identity on the actual algebraic union $\mathscr D$.
It has precisely the fractional-mode and sign convention of
\eqref{eq:twistedBorcherds}; no twisted module has been assumed.

For the vacuum, propagation with a free $\mathbf1$-input gives $P\lambda(\mathbf1,z;w)=\lambda(w)$, and hence $\mathbf1_m=\delta_{m,-1}\id_{\mathscr D}$. Taking $v=\mathbf1$, $p=-2$, $n=-1$ in
\eqref{eq:nodeborcherds713}, and using
$u_{-2}\mathbf1=L(-1)u$, $u_{-1}\mathbf1=u$ and $u_j\mathbf1=0$ for $j\ge0$, gives $(L(-1)u)_{m-1}+m u_{m-2}=u_{m-2},\; (L(-1)u)_q=-q\,u_{q-1}$. Therefore $Y_{\mathscr D}(L(-1)u,z)=\partial_zY_{\mathscr D}(u,z)$.
For the conformal state put $L_{\mathscr D}(q)=\omega_{q+1}$.
In the $p=0$ specialization of~\eqref{eq:nodeborcherds713}
use $\omega_0\omega=L(-1)\omega$,
$\omega_1\omega=2\omega$, $\omega_2\omega=0$ and
$\omega_3\omega=(c_V/2)\mathbf1$. With field mode indices
$m+1,n+1$ this gives
\[
\begin{aligned}
 [L_{\mathscr D}(m),L_{\mathscr D}(n)]
 &=-(m+n+2)\omega_{m+n+1}
   +2(m+1)\omega_{m+n+1}\\
 &\quad+\frac{c_V}{2}\binom{m+1}{3}
                    \mathbf1_{m+n-1}\\
 &=(m-n)L_{\mathscr D}(m+n)
       +\frac{c_V}{12}(m^3-m)\delta_{m+n,0}\id_{\mathscr D}.
\end{aligned}
\]
Together with~\eqref{eq:nodelevelbounds713}, linearity, the eigenmode
congruence and the vacuum identity, these are the weak $g$-twisted-module axioms. No grading restriction or
semisimplicity of $\mathscr D$ has been used.
\end{proof}

\begin{proposition}\label{prop:partialnodefiniteordinary}
Put $g_{\rm in}=h^{-1}$ and $g_{\rm out}=h$, where $h=g_Ag_B$.  For $\bullet\in\{{\rm in},{\rm out}\}$ there are multiplicities $m_{\bullet,S}\in\N$ such that
\begin{equation}\label{eq:partialnodefinite714}
 \mathscr D_\bullet\cong
 \bigoplus_{S\in\mathcal S_{g_\bullet}}
 S^{\oplus m_{\bullet,S}},\qquad
 \sum_{S\in\mathcal S_{g_\bullet}}
 m_{\bullet,S}\dim_\C S_{h_S}\leq\dim_\C B_{\bullet,1}<\infty.
\end{equation}
On $\mathscr D_{\rm in}$, the conformal mode constructed in Theorem~\ref{thm:partialnodeweakmodule} satisfies
\begin{equation}\label{eq:partialnodeenergy}
 (L_{\mathscr D_{\rm in}}(0)\lambda)(a,b)
 =\lambda\bigl((L_A(-1)+L_A(0))a,b\bigr)
  +\lambda\bigl(a,L_B(0)b\bigr)
 =(\mathcal E\lambda)(a,b).
\end{equation}
For the collision coefficients of Proposition~\ref{prop:collisionuniformsupport},
\begin{equation}\label{eq:partialnodelogfree714}
 \ell_{\beta,l}=0\qquad
 (\beta\in\Sigma+\N,\;1\leq l\leq L).
\end{equation}
\end{proposition}
\begin{proof}
Fix $\bullet\in\{{\rm in},{\rm out}\}$ and write $g=g_\bullet$ and
$N=\operatorname{ord}(g)$.  By Theorem~\ref{thm:partialnodeweakmodule}
and Proposition~\ref{prop:regularmodules} there is an algebraic internal
direct sum $\mathscr D_\bullet=\bigoplus_{i\in I_\bullet}M_i$ of ordinary
irreducible $g$-twisted modules; the index set is not assumed finite.
Let $\Omega_i=(M_i)_{h_i}\ne0$ be the lowest-weight space.  For
$u\in V_d\cap V_g^{(r)}$ and $j\ge1$ with $j\equiv r\pmod N$, the
actual balanced root monomial $s^je_u$, with
$e_u=\iota_z(u)z^{d-1}dz$ and $z=s^N$, has normalized node residue
$u_{d-1+j/N}$, not $Nu_{d-1+j/N}$; see
\eqref{eq:normalizednoderesidue} and \eqref{eq:partialnodemode}.
The node weak-module commutator gives
\[
 [L_{M_i}(0),u_{d-1+j/N}]=-\frac jN u_{d-1+j/N}.
\]
Because $M_i$ is ordinary and has no weight below $h_i$,
\begin{equation}\label{eq:partialground714}
 u_{d-1+j/N}\Omega_i=0 \qquad (j\ge1,\ j\equiv r\pmod N).
\end{equation}
This applies to the entire finite-jet
balanced frame: its $L(a>0)$ corrections are already expressed in
that frame, preserve finite-state truncations, and each resulting
monomial has its own allowed inertia index.  No identification of
different deck-twisted sectors is made.

We relate these node modes to the \emph{full} global Ward
annihilator defining $B_{\bullet,1}$.  Take $\lambda\in\Omega_i\subset
\mathscr D_\bullet$; choose $K$ with $\lambda\in\mathscr D_K$.
For an arbitrary $q\in\mathscr G_{\bullet,1}$, work in the finite
state cutoff containing $q$ and pair the meromorphic current $q$
with the actual propagated dual section of
Lemma~\ref{lem:partialequivariantresidue}.  The ordinary global residue
theorem, including the $H$-orbit normalization, reads
\[
 \lambda(R_\bullet(q)w)
   =-\sum_{\text{node orbits}}\sum_{u,j\ge1}
          c_{u,j}(q)\,(u_{d_u-1+j/N}\lambda)(w),
 \qquad w\in T_\bullet.
\]
Here the local coefficients $c_{u,j}(q)$ are those of $q$ in the
actual balanced root frames, with the node-orbit sum understood
with its stabilizer factors.  The sum is finite for fixed $w$:
$q$ has a finite state cutoff and the propagated dual section has
finite node pole order; hence only finitely many positive root
powers can contribute a residue.  If the cover has several node
orbits the same calculation is made on each, with its actual node
sector.  Every term on the right vanishes by the lowest-weight
calculation, so $\lambda(\mathscr G_{\bullet,1}T_\bullet)=0$.
Thus
\[
 \bigoplus_{i\in I_\bullet}\Omega_i\ \subseteq\
 \mathscr D_{\bullet,1}=B_{\bullet,1}^*.
\]
The sum remains direct in the algebraic dual $T_\bullet^*$.
Lemma~\ref{lem:partialcomponentfinite} gives
\[
 \sum_{i\in I_\bullet}\dim\Omega_i
 \le\dim B_{\bullet,1}<\infty.
\]
Since every $\Omega_i\ne0$, the index set is finite.  Replacing each
$M_i$ by its representative in $\mathcal S_g$ proves
\eqref{eq:partialnodefinite714}, and $L_{\mathscr D_\bullet}(0)$
is diagonalizable on the entire finite direct sum.

For the energy identity, use the actual normalized inner component
with coarse physical coordinate $x$ and the marked locations
$1,0,\infty$.  The global vector field $\vartheta=x\partial_x$ has
precisely the following local expansions:
\[
 \begin{array}{c|c|c}
 \text{marked point}&\text{physical coordinate}&\vartheta\\ \hline
 1&\xi=x-1&(1+\xi)\partial_\xi\\
 0&x&x\partial_x\\
 \infty&z=x^{-1}&-z\partial_z.
 \end{array}
\]
For $T_y(y)=Y(\omega,y)=\sum_{n\in\mathbb Z}L(n)y^{-n-2}$,
residues of $T_y(y)f(y)dy$ extract $L(-1)$ for $f=1$ and
$L(0)$ for $f=y$.  Hence the physical local residues at these
three points act respectively as
\begin{equation}\label{eq:partialstressres714}
 L_A(-1)+L_A(0),\qquad L_B(0),\qquad
 -L_{\mathscr D_{\rm in}}(0).
\end{equation}
At the last point the sign comes from $-z\partial_z$; its mode
is precisely the $\omega_1$ mode of
\eqref{eq:partialnodemode}, acting on the constructed
$h^{-1}$-twisted node module.  The original opposite-field
identification at $x=\infty$ is taken before this physical residue;
it is not the transpose of an untransported input mode.  The
physical changes $x\mapsto x-1$ and $x\mapsto x^{-1}$ are M\"obius,
so introduce no Schwarzian term.

The ramified root changes are not M\"obius.
At a marking of inertia order $N_i$, set $y=s^{N_i}$.
Equation~\eqref{eq:stress} gives the exact identity
\[
 T_s(s)=N_i^2s^{2N_i-2}T_y(s^{N_i})
     -\frac{c_V(N_i^2-1)}{24}s^{-2}.
\]
For the Euler residue, after the original orbit factor $1/N_i$,
\begin{align*}
 \frac1{N_i}\Res_s yT_y(y)\,dy
 &=\frac1{N_i^2}\Res_s sT_s(s)\,ds
       +\frac{c_V(N_i^2-1)}{24N_i^2}\\
 &=\Res_y yT_y(y)\,dy.
\end{align*}
In the second line we substituted the displayed transformation:
\[
 \Res_s sT_s(s)\,ds=N_i^2\Res_y yT_y(y)\,dy
 -c_V(N_i^2-1)/24.
\]
(The residues on the right are in their
respective normalized physical-mode convention.)  Thus the vacuum
Schwarzian term is \emph{retained and canceled} when the physical
mode is recovered; there is no additional lowest-weight or
central-charge shift.  The translation part $f(y)=1$ has no
Schwarzian residue, since the scalar term in
\eqref{eq:stress} is proportional to $s^{-N_i-1}ds$.
The same computation is performed at the output node, with its
actual inverse inertia and restricted-contragredient field.

To use the global Ward relation without assuming the desired energy
identity, fix $\lambda\in\mathscr D_{{\rm in},K}$ and pair its
finite-state propagated meromorphic dual section with the global
conformal current determined by $\vartheta$.  Unlike a
$(-KQ)$ test current, this particular current need not vanish at
the node: the meromorphic paired one-form is nevertheless defined,
and the global residue theorem includes its node residue.  That
residue is the \emph{definition} of the node conformal mode via
\eqref{eq:partialnodemode}; the two old-mark residues use the
original twisted-module Virasoro fields.  With the actual based
root frames and the normalizations just computed, their sum is
zero.  Therefore for every algebraic $a,b$,
\[
 (L_{\mathscr D_{\rm in}}(0)\lambda)(a,b)
 =\lambda((L_A(-1)+L_A(0))a,b)
    +\lambda(a,L_B(0)b)
 =(\mathcal E\lambda)(a,b).
\]
Finite homogeneous decomposition proves the formula for arbitrary
$a,b$, hence establishes both \eqref{eq:partialnodeenergy} and the
preservation of $\mathscr D_{\rm in}$ by $\mathcal E$.

For every $\beta,l,o$, Theorem~\ref{thm:finitecollisionkernels}
puts $\lambda_{\beta,l,o}$ in $\mathscr D_{\rm in}$.
Lemma~\ref{lem:collisionenergy} and the energy identity give, with
$H=L_{\mathscr D_{\rm in}}(0)$ and $c=\beta/e$,
\begin{equation}\label{eq:partiallogchain714}
 (H-c)\lambda_{\beta,l,o}
   =\frac{l+1}{e}\lambda_{\beta,l+1,o},
 \qquad 0\le l\le L,\qquad\lambda_{\beta,L+1,o}=0.
\end{equation}
Since $H$ is semisimple, its generalized $c$-eigenspace equals
its $c$-eigenspace (and both are zero if $c$ is not an eigenvalue).
The chain itself gives
$(H-c)^{L+1}\lambda_{\beta,0,o}=0$, so
$(H-c)\lambda_{\beta,0,o}=0$ and the \emph{$l=0$ equation} yields
$\lambda_{\beta,1,o}=0$.  Inductively, the equation at $l=j-1$
then gives $\lambda_{\beta,j,o}=0$ for every $j\ge1$.
This argument also covers repeated simple summands having the same
weight; no distinctness of the eigenvalues is needed.
Finally $o,a,b$ range over all algebraic tensors, so their
matrix coefficients separate $\ell_{\beta,l}$, establishing
\eqref{eq:partialnodelogfree714} for every $\beta$ and every $l>0$.
\end{proof}

\begin{proposition}\label{prop:partialnodevertices}
Set $h=g_Ag_B$, $U=\mathscr D_{\rm in}'$, and fix the based logarithm $\log1=0$. There are natural isomorphisms
\begin{equation}\label{eq:partialvertexrepresent715}
 U\cong A\bo_{P(1)}B,\qquad
 \Hom_V(S,\mathscr D_{\rm out})
 \xrightarrow{\ \sim\ }\cI\binom D{S\ C}
 \quad(S\in\mathcal S_h).
\end{equation}
They are induced by the evaluation maps at $1$; their inverse operators satisfy all four conditions of Definition~\ref{def:completeintertwiner}, including the original additional-field requirement of \cite[Definition~2.8]{DH}.
\end{proposition}
\begin{proof}
By Proposition~\ref{prop:partialnodefiniteordinary}, $\mathscr D_{\rm in}$ and $\mathscr D_{\rm out}$ are finite ordinary twisted modules. In particular
\[
 \mathscr D_{\rm in}=\bigoplus_{\mu} (\mathscr D_{\rm in})_\mu,
 \qquad U=\bigoplus_\mu(\mathscr D_{\rm in})_\mu^*,
 \qquad \overline U=\prod_\mu U_\mu,
\]
with finite-dimensional weight spaces. The algebraic tensor product $T_{\rm in}=A\otimes B$ and the evaluation pairing $\mathscr D_{\rm in}\subseteq T_{\rm in}^*$ are those of
\eqref{eq:partialnodedualspace}.
For $a\in A$, $b\in B$, define $j(a,b)\in\overline U$ by
\begin{equation}\label{eq:nodej715}
 \langle\lambda,j(a,b)\rangle=\lambda(a\otimes b)
 \quad(\lambda\in\mathscr D_{\rm in}).
\end{equation}
Its component in $U_\mu$ is the restriction of evaluation to $(\mathscr D_{\rm in})_\mu$; no finite weight support for $j(a,b)$ is asserted. For homogeneous $a\in A_\alpha$, $b\in B_\beta$ put
\begin{equation}\label{eq:nodeY715}
 \mathcal Y^{\rm in}(a,z)b
  =\sum_{\mu\in\wt(U)}
     z^{\mu-\alpha-\beta}\pi^U_\mu j(a,b),
 \qquad \log1=0.
\end{equation}
Since $U$ is lower bounded and has only finitely many relative-weight cosets, this is lower truncated. For every $\lambda\in(\mathscr D_{\rm in})_\mu$, its scalar coefficient is
\begin{equation}\label{eq:nodeYscalar715}
 \langle\lambda,\mathcal Y^{\rm in}(a,z)b\rangle
    =z^{\mu-\alpha-\beta}\lambda(a\otimes b).
\end{equation}

The following finite-direct-sum compatibility will be used below.  The
three-point converse in Proposition~\ref{prop:singlewardconverse} is
stated for simple ordinary modules.  If the output is a finite ordinary
direct sum $Z=\bigoplus_{i=1}^q Z_i$, then $Z'=\bigoplus_i Z_i'$ and
the global three-point current action preserves each output summand.
Consequently, with the original based cover and all its roots fixed,
\[
 \mathfrak W_{AB}^{Z}=\bigoplus_{i=1}^q\mathfrak W_{AB}^{Z_i},
 \qquad
 \cI\binom Z{A\ B}=\bigoplus_{i=1}^q\cI\binom{Z_i}{A\ B}.
\]
The first equality is obtained by applying every Ward generator to
each output projection; the second follows by projecting the four
Du--Huang axioms to the finitely many summands, and conversely taking
finite sums of their normally convergent matrix coefficients on
compact subchambers.  The finitely many prescribed corner exponent
sets and logarithm bounds may be united.  The same argument, separately
in any one of the three factors, applies to finite direct sums of
inputs.  Hence the two earlier reconstruction results apply to $U$
and to any finite ordinary test module $Z$; no extension to unrestricted
infinite direct sums is used.

We verify the analytic conditions before invoking any universal property. If $s$ is a global invariant current on the based three-point component cover, with state support in $V_{\le d}$
and a finite pole divisor at the marked orbits, its node Laurent expansion has finitely many contributing modes on each algebraic input. The node field of Theorem~\ref{thm:partialnodeweakmodule}
is defined by the same orbit-normalized residues as the external current actions of Definition~\ref{def:canonicalmarked}. For fixed $w=a\otimes b$ and a finite state cutoff $d$, let
$\tau_{\lambda,w,d}$ be the finite-rank propagated dual section from Lemma~\ref{lem:partialequivariantresidue}. At the node, its Laurent coefficients are the residues defining the node modes in Theorem~\ref{thm:partialnodeweakmodule}.  More explicitly, write a fixed
finite-state current in the balanced root frame at $Q$ as
$s_Q=\sum_{u,j}c_{u,j}s^j e_u$ with $e_u=\iota_z(u)z^{d_u-1}dz$
and $z=s^N$.  The coefficient of its node residue against the
germ $\tau_{\lambda,w,d}$ is
$\sum_{u,j}c_{u,j}(u_{d_u-1+j/N}\lambda)(w)$,
with $j$ restricted by the actual inertia character and with the
factor $1/N$ of~\eqref{eq:normalizednoderesidue}.  Only finitely
many pairs $(u,j)$ contribute: $u$ has finite state support,
$s_Q$ has a finite Laurent principal part at $Q$, and $\tau$ has a finite pole bound
there.  At a translated lift use its transported $H$-state and
root, rather than identifying two different twisted sectors.
Thus the orbit-normalized local pairing satisfies
\[
 (R_Q(s)\lambda)(w)
 =\frac1{|H|}\sum_{p\in HQ}
       \Res_p\langle\tau_{\lambda,w,d},s\rangle.
\]
At the two incoming marked orbits the same pairing gives $\lambda(R_P(s)w)$, with the contragredient opposite-field coordinate transform at $\infty$. The sum of all residues of this meromorphic one-form on the smooth projective cover is
zero. Hence
\begin{equation}\label{eq:nodeevaluationward}
 (R_Q(s)\lambda)(a\otimes b)
 +\lambda\bigl(R_P(s)(a\otimes b)\bigr)=0,
 \qquad \lambda\in\mathscr D_{\rm in}.
\end{equation}
Here $R_Q(s)$ is the actual node-mode action, and $R_P(s)$ is the sum of the incoming actions. The calculation takes place in a single finite state cutoff; increasing that cutoff leaves the Laurent coefficients unchanged by the compatibility
in Lemma~\ref{lem:partialequivariantresidue}. No current with an infinite state support or unrestricted infinite node mode sum is used. The functional
\[
 E:\mathscr D_{\rm in}\otimes A\otimes B\longrightarrow\C,
 \qquad E(\lambda\otimes a\otimes b)=\lambda(a\otimes b),
\]
annihilates all three-point global Ward relations by
\eqref{eq:nodeevaluationward}; each algebraic $\lambda$ is a
finite sum of homogeneous components. Thus $E$ belongs to the full three-point Ward space $\mathfrak W_{AB}^{U}$ of Section~\ref{sec:reconstruction}. Proposition~\ref{prop:singlewardconverse}
therefore proves all three single-$V$ branches and their absolute original projected sums for \eqref{eq:nodeY715}. Proposition~\ref{prop:partialallfields}, applied to the same full three-point Ward functional, gives the original finite
strings of additional $V$ fields, the radial weight sums and the component-isolated regularity for every allowed collection $(M_{ij})$; the based branches are unchanged. These are items
(ii) and (iii) of Definition~\ref{def:completeintertwiner}. For $\lambda\in(\mathscr D_{\rm in})_\mu$, the stress-current identity \eqref{eq:partialnodeenergy} yields $\lambda(L_A(-1)a\otimes b) =(\mu-\alpha-\beta)\lambda(a\otimes b)$.
By \eqref{eq:nodeYscalar715}, this is precisely
$\partial_z\mathcal Y^{\rm in}(a,z)b
=\mathcal Y^{\rm in}(L_A(-1)a,z)b$, so item (iv) holds.
Thus
\[
\mathcal Y^{\rm in}\in\cI\binom U{A\ B},
 \qquad \operatorname{ev}_0(\mathcal Y^{\rm in})=j,
\]
where $\operatorname{ev}_0$ denotes evaluation at $z=1$ with $\log1=0$ as in \eqref{eq:Pzinverse}.

For an arbitrary $Z\in\cC_h$ define
\[
\Theta_Z:\Hom_V(U,Z)\longrightarrow\cI\binom Z{A\ B},
 \qquad f\longmapsto\overline f\circ\mathcal Y^{\rm in}.
\]
This is well-defined and natural in $Z$. We construct its inverse without assuming \eqref{eq:partialvertexrepresent715}. For $\mathcal Y\in\cI\binom Z{A\ B}$ set
\begin{equation}\label{eq:nodeiota715}
 (\iota_{\mathcal Y}z')(a\otimes b)
 :=\langle z',\mathcal Y(a,1)b\rangle,
 \qquad z'\in Z'.
\end{equation}
Fix a finite sum of homogeneous vectors $z'\in Z'$ and let $R(z')$ be the greatest relative weight, measured above the lowest weights of the finitely many simple summands containing its components. Choose $k\ge1$ with $k/N>R(z')$, where
$N=\operatorname{ord}(h)$. Every balanced node coefficient of a current in $\mathscr G_{{\rm in},k}$ has allowed root order $j\ge k$ and hence acts on $Z'$ through a mode of degree $-j/N<-R(z')$. It annihilates $z'$. The complete three-point
Ward identity of $\mathcal Y$ (Proposition~\ref{prop:binaryward}, with the contragredient node label) consequently gives
\[
 \iota_{\mathcal Y}(z')\bigl(\mathscr G_{{\rm in},k}
                   (A\otimes B)\bigr)=0,
 \qquad \iota_{\mathcal Y}(z')\in\mathscr D_{{\rm in},k}.
\]
The level $k$ is allowed to depend on $z'$, as required by the algebraic union $\mathscr D_{\rm in}=\bigcup_k\mathscr D_{{\rm in},k}$.
The actual Virasoro $L(0)$ relation also checks the grading of this
map.  If $z'\in(Z')_\gamma$, $a\in A_\alpha$ and
$b\in B_\beta$, the $L(0)$ commutator and the original
$L(-1)$ derivative identity for $\mathcal Y$ at $z=1$ imply
\[
 \langle z',\mathcal Y(L_A(-1)a,1)b\rangle
 =(\gamma-\alpha-\beta)
 \langle z',\mathcal Y(a,1)b\rangle.
\]
Together with~\eqref{eq:partialnodeenergy} this gives
$L_{\mathscr D_{\rm in}}(0)\iota_{\mathcal Y}(z')
=\gamma\iota_{\mathcal Y}(z')$.  Hence $\iota_{\mathcal Y}$
is homogeneous-weight preserving, so its restricted transpose below
is defined on the algebraic restricted dual; no full algebraic dual
or completed-dual transpose is substituted.

For the promised mode comparison, fix $u\in V_d$ and $w=a\otimes b$.
In the \emph{actual output-node root coordinate} $z=s^N$ let
\[
 G_{u,z',w}(z)
 :=\big\langle Y_{Z'}^{h^{-1}}(u,z)z',
              \mathcal Y(a,1)b\big\rangle,
\]
understood coefficientwise, with the original contragredient
opposite-field and full finite $L(1)$ gauge at the output.  This
is the node expansion of the original one-$V$ three-point duality
function.  Its complete germs at both old incoming marks are exactly
those prescribed for the one-free-point propagation of
$\iota_{\mathcal Y}(z')$: for every external mode coefficient,
they agree by the product and iterate expansions of that same
one-$V$ duality function, with unchanged based roots and deck stems.
For example, if $v_q^A a$ is the coefficient of the incoming
field $Y_A(v,\xi)a$ at $1$, the compared coefficient is exactly
\[
 \big\langle Y_{Z'}^{h^{-1}}(u,z)z',
       \mathcal Y(v_q^Aa,1)b\big\rangle;
\]
the other old input gives
$\langle Y_{Z'}(u,z)z',\mathcal Y(a,1)v_q^Bb\rangle$.
The original one-$V$ duality applies separately to every algebraic
$v_q^Aa$ and $v_q^Bb$; their full coefficient strings, including
regular coefficients, are the prescribed incoming germs.  This
comparison uses no prospective node-mode intertwining identity.
To justify the asserted global meromorphic comparison rather than
postulate it, use the reverse direction of
Proposition~\ref{prop:singlewardconverse}: the original complete
$\mathcal Y$ determines its actual full three-point Ward
functional.  Its finite-state strong-residue reconstruction is the
meromorphic dual section whose output Laurent germ is $G$, with the
same original based paths.  The propagated side is meromorphic by
Lemma~\ref{lem:partialequivariantresidue}.  Both therefore are
finite-pole rooted dual sections; a common finite cutoff and a
common finite node level $K$ contain them.
Complete-germ uniqueness in
Lemma~\ref{lem:partialequivariantresidue} therefore identifies
their whole Laurent series, not only their singular parts:
\[
 \mathcal P_{\iota_{\mathcal Y}z',u}(z;a,b)
   =G_{u,z',a\otimes b}(z).
\]
Taking the \emph{normalized} coefficient of $z^{-m-1}$ yields,
for every allowed $m\in r/N+\mathbb Z$,
\[
 \bigl(u_m\iota_{\mathcal Y}z'\bigr)(a,b)
 =\langle u_m^{Z'}z',\mathcal Y(a,1)b\rangle
 =\bigl(\iota_{\mathcal Y}(u_m^{Z'}z')\bigr)(a,b).
\]
On every fixed vector these coefficients and their root-orbit sums
are finite; the $h^{-1}$ node and $Z'$ have the same twisted label.
This proves the mode intertwining before transposition, for
all algebraic $a,b$, without relying on a later matching lemma.
Thus $\iota_{\mathcal Y}:Z'\to\mathscr D_{\rm in}$ is a $V$-homomorphism. Restricted contragredience of the finite ordinary objects supplies a $V$-map
\[
 f_{\mathcal Y}=\iota_{\mathcal Y}':U=\mathscr D_{\rm in}'
            \longrightarrow (Z')'\cong Z.
\]
Its evaluation satisfies, for homogeneous $a,b,z'$,
\[
 \langle z',\overline{f_{\mathcal Y}}j(a,b)\rangle
 =\langle\iota_{\mathcal Y}(z'),j(a,b)\rangle
 =\langle z',\mathcal Y(a,1)b\rangle.
\]
The weight-compatibility formula of Lemma~\ref{lem:IOlogs} then identifies the entire operators: $\Theta_Z(f_{\mathcal Y})=\mathcal Y$. Conversely,
\eqref{eq:nodej715} gives
$\iota_{\Theta_Z(f)}=f':Z'\to\mathscr D_{\rm in}$, so $f_{\Theta_Z(f)}=f$. Hence $\Theta_Z$ is a natural bijection. Evaluation at $1$ identifies it with $f\mapsto\overline f\circ j$, the universal map of
Definition~\ref{def:Pzproduct}. The uniqueness of the $P(1)$-tensor product in Theorem~\ref{thm:binary} gives the first isomorphism in \eqref{eq:partialvertexrepresent715}, with its canonical evaluation map.

For the outer component take $S\in\mathcal S_h$ and $f\in\Hom_V(S,\mathscr D_{\rm out})$. Evaluation gives $E_f(d'\otimes s\otimes c)=f(s)(d'\otimes c),\; d'\in D',\ s\in S,\ c\in C$. The full outer-node Ward identity, the equivariance of $f$,
and the finite homogeneous decomposition of $s$ make $E_f$ a full three-point Ward functional on $D'\otimes S\otimes C$. Proposition~\ref{prop:singlewardconverse} reconstructs a
single-insertion operator from $E_f$; Proposition~\ref{prop:partialallfields} verifies all its additional-$V$ conditions, for the prescribed multipliers, branches and original projected sums. On homogeneous
$d'\in(D')_\gamma$, $s\in S_\alpha$, $c\in C_\beta$ its matrix coefficient is explicitly
\begin{equation}\label{eq:outerY715}
 \langle d',\mathcal Y_f(s,z)c\rangle
    =z^{\gamma-\alpha-\beta}f(s)(d'\otimes c),
\end{equation}
which also verifies truncation. The outer stress-current Ward relation, with marked points $\infty,1,0$, gives $f(L_S(-1)s)(d'\otimes c) =(\gamma-\alpha-\beta)f(s)(d'\otimes c)$, and \eqref{eq:outerY715} proves the derivative condition.
Consequently $\mathcal Y_f\in\cI\binom D{S\ C}$.

Conversely, for $\mathcal Y\in\cI\binom D{S\ C}$ define
\begin{equation}\label{eq:outeriota715}
 f_{\mathcal Y}(s)(d'\otimes c)
  :=\langle d',\mathcal Y(s,1)c\rangle.
\end{equation}
For every algebraic $s\in S$, its finite relative-weight support supplies $k$ such that all balanced node coefficients of order $\ge k$ annihilate $s$. The full three-point Ward identity makes $f_{\mathcal Y}(s)\in\mathscr D_{{\rm out},k}$.
Fix $u\in V_d$, $s_0\in S$, and algebraic $d',c$.
The original complete one-$V$ duality for $\mathcal Y$ provides,
in the \emph{physical coordinate of the outer incoming node at $1$},
the whole local germ
\[
  G^{\rm out}_{u,s_0,d',c}(z)
  =\langle d',\mathcal Y(Y_S^{h}(u,z)s_0,1)c\rangle,
\]
including its regular Taylor part.  Each coefficient is a finite
matrix coefficient on the algebraic input $s_0$.  The independently
propagated germ
$\mathcal P_{f_{\mathcal Y}(s_0),u}(z;d'\otimes c)$
has the same complete germs at the old output and $C$ marking:
these follow from the other two original one-$V$ regions of
the same $\mathcal Y$.  For instance, at the old $C$-mark
and for an incoming coefficient $v_q^C c$ the compared scalar is
$\langle d',\mathcal Y(Y_S(u,z)s_0,1)v_q^Cc\rangle$;
at the old output it is obtained by the full opposite-field
operator acting on $d'$.  All these are equal as \emph{complete}
formal incoming coefficients by the fixed original operator's
one-$V$ duality, with the finite extra-field expansion supplied
when necessary by Proposition~\ref{prop:partialallfields}.
The reverse Ward direction of Proposition~\ref{prop:singlewardconverse}
identifies the matrix coefficients of this original $\mathcal Y$
with an actual meromorphic three-point dual section; the propagation
on the other side is meromorphic by
Lemma~\ref{lem:partialequivariantresidue}.  Both objects have finite node pole
order after choosing a common cutoff and level $K$, since
$f_{\mathcal Y}(s_0)\in\mathscr D_{{\rm out},k}$ and the original
twisted field $Y_S(u,z)s_0$ is lower truncated.
By the complete-germ uniqueness in
Lemma~\ref{lem:partialequivariantresidue} they agree as rooted
meromorphic sections.  Extracting each normalized node coefficient
therefore gives the actual typed identity
\[
 (u_m^{\mathscr D_{\rm out}}f_{\mathcal Y}(s_0))(d',c)
 = f_{\mathcal Y}(u_m^S s_0)(d',c).
\]
The required root phases and $1/N$ normalization occur on both
sides through the same based node coordinate; when a deck
arrow changes the sector it is applied to both sides in the fixed
chronological order.  Hence $f_{\mathcal Y}$ is a genuine
$V$-homomorphism, rather than merely a linear map into $T_{\rm out}^*$. The evaluation identities
\eqref{eq:outerY715}--\eqref{eq:outeriota715} and
Lemma~\ref{lem:IOlogs} show that $\mathcal Y\mapsto f_{\mathcal Y}$ and $f\mapsto\mathcal Y_f$ are mutually inverse and natural in $S$. This proves the second identification in
\eqref{eq:partialvertexrepresent715}.  If $a:S_1\to S$ is a
$V$-homomorphism in the $h$-sector, then
$E_{f\circ a}(d',s,c)=E_f(d',a(s),c)$, so the outer
bijection is natural contravariantly in $S$; the inner bijection is
natural covariantly in $Z$ by composition of weight-preserving
maps.  The two directions recover the \emph{complete} formal
intertwining operators because their fixed homogeneous matrix
coefficients have the same value at $1$ and the forced weight power
$z^{\gamma-\alpha-\beta}$, not merely identical full-Ward
annihilators.
\end{proof}

\subsection{Balanced node jets and the coefficient kernel}
Fix the actual admissible source family $p:C_O\to\operatorname{Spec}O$, $O=\C[[t]]$, of Lemma~\ref{lem:canonicalnodecurrent}.  At a chosen node lift use
\begin{equation}\label{eq:nodejetlocalring}
 B=O[[s_+,s_-]]/(s_+s_--t^a),\qquad
 x=s_+^N,\quad y=s_-^N,\quad q=t^e,\quad e=aN.
\end{equation}
Here $a,N\in\Z_{>0}$ and $N$ is the inertia order.  In the balanced current frame of \eqref{eq:balancednodeframe} the completed
finite-state sheaf is $\widehat{\mathscr F}_{d,O}=B\otimes_\C
V_{\le d}$. We identify this common free module using the plus-branch balanced frame. Thus a coefficient state $b$ in the common frame is represented by $b$ on the plus branch and by $\Sigma_db$ on the minus branch, where $\Sigma_d=e^{L(1)}(-1)^{L(0)+1}|_{V_{\le d}}$. This convention applies for both signs of $k$. Fix $b\in V_d$ homogeneous for the inertia action, and
$k\in\Z$ such that the coefficient below is inertia-invariant. Put
\begin{equation}\label{eq:balancedprescribedjet}
 \sigma_k(b)=
 \begin{cases}
    s_+^k b,&k\ge0,\\
    s_-^{-k}b,&k<0.
 \end{cases}
\end{equation}
The state-current transition to the opposite branch is the full $L(1)$-corrected balanced transition, with the opposite inertia character.  In particular, on the punctured overlap,
\begin{equation}\label{eq:balancedjetbranch}
 s_+^k=t^{ak}s_-^{-k}\quad(k\ge0),\qquad
 s_-^{-k}=t^{-ak}s_+^k\quad(k<0).
\end{equation}
The monomials in \eqref{eq:balancedprescribedjet} belong to the integral node module; neither identity asserts the integrality of the displayed negative power on its own.

\begin{lemma}\label{lem:matchingnodejetlift}
Let $J\in\Z_{\ge0}$ and choose $K\in\Z_{>0}$ satisfying $K>|k|,\; aK\ge J+1$. There exist $m\ge0$ and an $H$-invariant section $s\in H^0(C_O,\mathscr F_{d,O}(mP))^H$ whose jet at the chosen node, in the balanced frame, equals
$\sigma_k(b)$ modulo
\begin{equation}\label{eq:matchingjetideal}
 I_{J,K}=(t^{J+1},s_+^K,s_-^K)\subset B.
\end{equation}
At the other node lifts the prescribed jets are its $H$-translates. Here $P$ is the actual full external marked-orbit divisor and $\mathscr F_{d,O}$ is the $O$-flat coherent finite-state current sheaf of Lemma~\ref{lem:truncatedcomponentlift}.
In the \emph{balanced relative component frames}, the node germs of the coefficient currents of $s$ through $t$-order $J$ agree with the prescribed balanced Laurent germs, modulo local errors with component node order at least
\[
K-\left\lfloor J/a\right\rfloor
\]
and modulo $t^{J+1}$. These local errors need not themselves extend to global component currents. Their node residue pairing with a propagated functional of level at most this order is zero. No assertion concerning the physical relative-weight dilation is included in this estimate.
\end{lemma}
\begin{proof}
Write $O_J=O/(t^{J+1})$, and let $Q$ denote the finite set of actual
source nodes.  We construct the prescribed jet using an algebraic coherent quotient and establish the necessary algebraization.  At one selected node let $A$ be the Noetherian
local ring of $C_O$, $\mathfrak m$ its maximal ideal, and
$\widehat A=B$ with maximal ideal $\widehat{\mathfrak m}=(t,s_+,s_-)$.  The
quotient $B/I_{J,K}$ is a quotient of
$\C[t,s_+,s_-]/(t^{J+1},s_+^K,s_-^K)$, hence is Artinian of
finite $\C$-dimension (at most $(J+1)K^2$).  In particular
$\widehat{\mathfrak m}^{n_0}\subset I_{J,K}$ for some $n_0$.  Completion
induces the canonical isomorphism
$A/\mathfrak m^{n_0}\simeq
 \widehat A/\widehat{\mathfrak m}^{n_0}$.
Take the inverse image of $I_{J,K}/\widehat{\mathfrak m}^{n_0}$
in $A/\mathfrak m^{n_0}$, and call the resulting ideal of $A$
$I_A$.  Thus $A/I_A\simeq B/I_{J,K}$, compatibly with completion.
The chosen balanced trivialization similarly identifies the finite
jet of the actual current sheaf with
$(B/I_{J,K})\otimes_{\C}V_{\le d}$; the ideal $I_A$ defines an
actual coherent quotient supported at the node.  There is
no claim that a general infinite formal germ is algebraic.

The inertia generator acts on the node coordinates by
$s_+\mapsto\epsilon s_+$ and $s_-\mapsto\epsilon^{-1}s_-$,
and on states by the prescribed inertia representation.  It preserves
$I_{J,K}$; the hypothesis that $\sigma_k(b)$ is
inertia-invariant says exactly that its class is fixed by this
stabilizer.  Strict deck transport sends this jet to each of the other
nodes in the same $H$-orbit, and the ideal, being symmetric in the two
oriented branches, is carried to the corresponding ideal there.
The finite-state transition is a $B$-linear change of frame of the
actual coherent sheaf; it cannot change the ideal defining its
jet quotient.  The complete opposite-branch balanced transition is
the constant, invertible
$\Sigma_d=e^{L(1)}(-1)^{L(0)+1}|_{V_{\le d}}$
from \eqref{eq:balancednodeframe}; it includes every positive
$L(1)$ correction and preserves $V_{\le d}$.  This also verifies
compatibility when deck transport compares opposite orientations.
Taking the direct sum of these local finite-length quotients over
$Q$ therefore gives an actual $H$-equivariant coherent quotient
$\mathscr F_{d,O}\twoheadrightarrow\mathscr J_{J,K}$,
with coherent kernel $\mathscr K_{J,K}$, and the prescribed jets give
$\bar\sigma\in H^0(C_O,\mathscr J_{J,K})^H$.
Finite sums of compatible monomial jets are handled by the same
quotient and linearity.

The external marked divisor $P$ is disjoint from $Q$, so
$\mathscr J_{J,K}(mP)$ is canonically identified with
$\mathscr J_{J,K}$ using the nowhere-vanishing canonical section
of $\mathcal O(mP)$ near $Q$.  The allowed actual admissible-cover
geometry ensures that $P$ meets every geometric irreducible
special-fiber component in positive degree.  Hence $\mathcal O(P)$
is relatively ample on this proper nodal family.  Relative Serre
vanishing \cite[Tag~02O1]{Stacks}, applied specifically to the
coherent \emph{kernel} $\mathscr K_{J,K}$, yields for $m\gg0$
\begin{equation}\label{eq:matchingjetserre}
 R^1p_*\bigl(\mathscr K_{J,K}(mP)\bigr)=0.
\end{equation}
The genuine coherent exact sequence
\[
 0\longrightarrow\mathscr K_{J,K}(mP)
 \longrightarrow\mathscr F_{d,O}(mP)
 \longrightarrow\mathscr J_{J,K}\longrightarrow0
\]
and affineness of $\operatorname{Spec}O$ now give
\[
 H^0(C_O,\mathscr F_{d,O}(mP))
 \twoheadrightarrow H^0(C_O,\mathscr J_{J,K}).
\]
Lift $\bar\sigma$ and apply $|H|^{-1}\sum_{h\in H}h$.
Since the chosen jet is already invariant, averaging leaves it
unchanged.  The resulting $s$ is an $H$-invariant global section of
the actual current sheaf, with the required jet on every node
orbit.  Here $m$ depends only on the fixed finite data, never on a
subsequently chosen Ward test tensor.

We check the coefficientwise estimate in the original balanced
frames.  Express the difference between the completed jet of $s$
and its prescribed lift as
\[
 E=t^{J+1}E_0+s_+^KE_++s_-^KE_-,
 \qquad E_0,E_+,E_-\in B\otimes V_{\le d}.
\]
It is enough to calculate on the monomials of these \emph{integral}
series, since a fixed $t$-coefficient and fixed node exponent receives
only finitely many contributions.  For
$t^js_+^{K+r}s_-^s$ ($j,r,s\ge0$), the substitutions in the two
punctured normalized branches give
\[
 t^{j+as}s_+^{K+r-s},\qquad
 t^{j+a(K+r)}s_-^{s-K-r}.
\]
Set $L=K-\lfloor J/a\rfloor>0$.  If the first term survives
modulo $t^{J+1}$, then $j+as\le J$ and hence
$K+r-s\ge K-\lfloor J/a\rfloor=L$.  The second term is zero
modulo $t^{J+1}$ because $a(K+r)\ge aK\ge J+1$.
The calculation for $s_-^KB$ is the same with the branches exchanged,
and $t^{J+1}B$ vanishes modulo $t^{J+1}$ on both branches.
Consequently, for every $0\le j\le J$, the coefficient of $t^j$
in the error is a local node germ vanishing to order at least $L$
in its balanced frame. This conclusion concerns germs; subtraction
of a prescribed local monomial does not assert that the difference
is a global component current.

It remains to relate these node germs to the actual global residue
operators, without forgetting state-coordinate factors. The
opposite-branch frame change is precisely the complete constant
$\Sigma_d$, so its action replaces each coefficient state $c$
by the finite sum $\Sigma_dc$ \emph{at the same node monomial}.
It cannot alter either the $t$-order or the node order above.
The root rotation multiplies $s_+^p$ or $s_-^p$ by a root of unity;
inertia selection discards forbidden congruence classes, and typed
deck pullback transports $c$ and its module sector together.
None decreases the order. By Lemma~\ref{lem:canonicalnodecurrent},
the coefficient restrictions $s_{j,\bullet}$ of the actual global
$s$ are genuine meromorphic currents on the normalized component,
with poles only at its external markings and nodes. Their node
germs, rather than their global differences from local monomials,
are the objects just estimated.

For $\lambda\in\mathscr D_{\bullet,k_0}$ with $k_0\le L$, pair
$s_{j,\bullet}$ with the actual propagated dual section of
Lemma~\ref{lem:partialequivariantresidue}. In a balanced root
frame, pairing that section with a basis current has the form
$f(s_\pm)ds_\pm$, with
$f\in s_\pm^{-k_0}\C[[s_\pm]]$ by
\eqref{eq:partialnodepole78}. Multiplication by a local error in
$s_\pm^L\C[[s_\pm]]\otimes V_{\le d}$ therefore has zero
node residue. The global residue theorem on the proper normalized
component gives
\[
 \lambda(R_\bullet(s_{j,\bullet})w_\bullet)
 =-\sum_{\text{node orbits}}\Res^{\rm base}
    \langle P\lambda,s_{j,\bullet}\rangle,
\]
and each node residue on the right may be replaced by its
prescribed Laurent monomial. The orbit factors are exactly
$1/N_i$ as in \eqref{eq:canonicalmarkedres}; the external output
uses its full restricted-contragredient field. Each residue is a
finite mode sum on a fixed algebraic tensor. This proves the
claimed residue consequence without globalizing a local error or
using any node-matching identity.

All these estimates are in the balanced \emph{relative component}
frames.  In particular, if a scaled external leg has a mode $b_n$
of weight shift $\Delta=\mathrm{wt}(b)-n-1$, then
$\mathsf D_t^{-1}b_n\mathsf D_t=t^{-e\Delta}b_n$;
its exponent can be negative.  Therefore our modulo-$t^{J+1}$
comparison cannot be conjugated to the physical frame without a
separate coefficientwise analysis.  No such physical-frame claim is
used here.
\end{proof}

For a fixed actual product $F$, use its normalized collision coefficients and write $\ell_\beta=\ell_{\beta,0}$; $\ell_{\beta,l}=0$ for $l>0$ by Proposition~\ref{prop:partialnodefiniteordinary}.  By
Theorem~\ref{thm:finitecollisionkernels} and the partial node spaces, each $\ell_\beta$ belongs to the algebraic tensor $\mathscr D_{\rm out}\otimes\mathscr D_{\rm in}$. Set $\ell_\gamma=0$ for an absent exponent $\gamma$ and define
\begin{equation}\label{eq:extractdepth716}
 J_\Sigma(\gamma)=\max\bigl(\{0\}\cup
 \{n\in\Z_{\ge0}:\gamma=\sigma+n,\ \sigma\in\Sigma\}\bigr).
\end{equation}
The set $\Sigma$ is the finite uniform exponent support of Proposition~\ref{prop:collisionuniformsupport}; in particular
\eqref{eq:extractdepth716} is well-defined even when
$\ell_\gamma=0$.  The shift $ak=e(k/N)\in\Z$ preserves exponent classes modulo $\Z$.

For homogeneous $u\in V_d$, $d\in\Z_{\ge0}$, and an allowed mode index $m$ define, on the inverse-inertia node module,
\begin{equation}\label{eq:oppositenodemode}
 C_m(u)=(-1)^d\sum_{j=0}^{d}\frac{1}{j!}
 \bigl(L(1)^j u\bigr)_{2d-j-m-2}.
\end{equation}
Thus $C_m(u)$ is an operator (not a state); the modes on the right act on the specified inverse-twisted module.

\begin{lemma}\label{lem:oppositemodes}
Let $u\in V_d$ be homogeneous with $gu=e^{2\pi ir/N}u$, and let $m\in r/N+\Z$.  Put $\alpha=m-d+1$. Then \eqref{eq:oppositenodemode} is the opposite-field mode appearing in the restricted contragredient action of
\cite[(3.16), Lemma~3.8]{DLM}.  For every nonzero summand,
\[
 \deg\bigl((L(1)^j u)_{2d-j-m-2}\bigr)=\alpha,
 \qquad 2d-j-m-2\in-r/N+\Z.
\]
Moreover, as operators on the inverse-twisted module,
\begin{equation}\label{eq:oppositederivative717}
 C_m(L(-1)u)=-m C_{m-1}(u).
\end{equation}
\end{lemma}
\begin{proof}
Let $M$ be the indicated inverse-$g$-twisted module, with all
modes on the right interpreted as operators on $M$.  Starting from
the original opposite-field expression (not a new field
axiom), expand a homogeneous $u\in V_d$:
\[
 Y_M^{\mathrm{op}}(u,z)
 =Y_M\bigl(e^{zL(1)}(-z^{-2})^{L(0)}u,z^{-1}\bigr)
 =(-1)^d\sum_{j=0}^d\frac{z^{j-2d}}{j!}
   \sum_{q\in-r/N+\Z}(L(1)^ju)_qz^{q+1}.
\]
The $j$-sum ends at $d$ because $L(1)^ju\in V_{d-j}$,
$V_{<0}=0$, and $gL(1)=L(1)g$.  The coefficient at
$z^{-m-1}$ occurs exactly when
$q+j-2d+1=-m-1$, that is,
$q=2d-j-m-2$.  This proves formula
\eqref{eq:oppositenodemode}, including the sign and the mode
index.  Since $m\in r/N+\Z$, each such $q$ lies in
$-r/N+\Z$, the actual mode coset for $M$; no two
sectors have been identified.  Its weight shift on $M$ is
\[
 (d-j)-q-1=m-d+1=\alpha,
\]
independently of $j$.  Each sum is finite on every algebraic
vector; its opposite-field coefficient is well defined in the
usual restricted-contragredient pairing.  More precisely, for
the actual restricted dual $M'$,
\[
 \langle u_m^{M'}w',w\rangle
   =\langle w',C_m^M(u)w\rangle,
 \qquad w'\in M',\ w\in M,
\]
which fixes the transpose direction and gives the DLM
contragredient action.  The transpose on $M'$ is well defined
because the homogeneous modes have fixed finite weight shift.

For the derivative, let $v_j=L(1)^ju$, and use the Virasoro
commutators to calculate (with $v_{-1}=0$)
\[
 L(1)^jL(-1)u=L(-1)v_j+j(2d-j+1)v_{j-1}
 \quad(j\ge0).
\]
Now $L(-1)u\in V_{d+1}$, so the \emph{left} opposite-mode
formula has $0\le j\le d+1$, not merely $j\le d$.
The $j=d+1$ term is retained in the reindexing below, even
though $v_{d+1}=0$.  Applying the ordinary state-mode derivative
$(L(-1)v_j)_q=-q(v_j)_{q-1}$ to the first term gives
\[
\begin{aligned}
 C_m(L(-1)u)
 &=(-1)^{d+1}\sum_{j=0}^{d+1}\frac1{j!}
 \Bigl[-(2d-j-m)(v_j)_{2d-j-m-1}\\
 &\hspace{42mm}{}+j(2d-j+1)(v_{j-1})_{2d-j-m}\Bigr].
\end{aligned}
\]
Put $i=j-1$ in the second sum.  For each $0\le i\le d$,
the common state $v_i$ and common mode index $2d-i-m-1$
then have total coefficient
\[
 \frac{(-1)^{d+1}}{i!}
 \bigl[-(2d-i-m)+(2d-i)\bigr]
 =\frac{(-1)^{d+1}m}{i!}.
\]
The first summand at $j=d+1$ is zero; the final second-sum term
at $j=d+1$ contributes to $i=d$ and has therefore \emph{not}
been omitted.  Since
$C_{m-1}(u)=(-1)^d\sum_{i=0}^d
(v_i)_{2d-i-m-1}/i!$, the claimed identity follows:
\[
 C_m(L(-1)u)=-mC_{m-1}(u).
\]
The calculation holds for fractional $m$ with the same algebraic
mode convention and is independent of any later node-matching theorem.
\end{proof}

\begin{theorem}\label{thm:matchedcoefficientmodes}
Put $h=g_Ag_B$, $U=\mathscr D_{\rm in}'$, and $e=aN$ as in \eqref{eq:nodejetlocalring}.  Let $u\in V_d\cap V_h^{(r)}$ be homogeneous, where $hu=e^{2\pi ir/N}u$, and let
\[
 k\in r+N\Z,\qquad m=d-1+\frac{k}{N},\qquad
 \alpha=\frac{k}{N}.
\]
For $\beta\notin\Sigma+\Z_{\ge0}$ set $\ell_\beta=0$. Then, as an identity in the algebraic tensor product $\mathscr D_{\rm out}\otimes\mathscr D_{\rm in}$,
\begin{equation}\label{eq:matchedcoefficientidentity}
 (u_m\otimes\id)\ell_\beta
 =(\id\otimes C_m(u))\ell_{\beta-e\alpha}
 \qquad(\beta\in\C).
\end{equation}
Furthermore,
\begin{equation}\label{eq:matchedcoefficientweights718}
 \ell_\beta\in
 (\mathscr D_{\rm out})_{\beta/e}\otimes
 (\mathscr D_{\rm in})_{\beta/e},
\end{equation}
where a missing weight space is zero.  The maps
\begin{equation}\label{eq:matchedcoefficientmap718}
 T_F|_{U_w}:U_w=(\mathscr D_{\rm in})_w^*
 \longrightarrow(\mathscr D_{\rm out})_w,
 \qquad x\longmapsto(\id\otimes x)\ell_{ew},
\end{equation}
assemble to a homomorphism of ordinary $h$-twisted $V$-modules $T_F:U\to\mathscr D_{\rm out}$.
\end{theorem}
\begin{proof}
All tensors in this proof are algebraic.  By
Theorem~\ref{thm:finitecollisionkernels} each coefficient is a finite-rank
bilinear form, so it is an element of
$\mathscr D_{\rm out}\otimes\mathscr D_{\rm in}$, and, by
Proposition~\ref{prop:partialnodefiniteordinary}, both node modules
are finite direct sums of ordinary twisted modules.  In particular,
they have finite-dimensional homogeneous subspaces, and the mode
operators in the assertion act on the actual algebraic node modules.
We work throughout in the balanced \emph{relative} frames, never
transporting a congruence by a possibly nonintegral physical dilation.

The comparison furnished by a global current is as follows.
On the $A,B$ input legs, the relative-weight dilation is
$\mathsf D_t|_{A_\mu\otimes B_\nu}
=t^{e((\mu-h_A)+(\nu-h_B))}\id$; it is the identity on the
outer $D',C$ legs.  For a mode $b_n$ of $b\in V_d$ on either
of the two scaled legs, put $\Delta=d-n-1$.  On a homogeneous vector
$v$ of relative weight $q$ one has the \emph{modewise} identity
\[
  \mathsf D_t^{-1}b_n\mathsf D_tv
   =t^{eq-e(q+\Delta)}b_nv=t^{-e\Delta}b_nv.
\]
Every such exponent is integral for the common ramification $e$;
there is no additional $h_A$ or $h_B$ factor, because their
contributions cancel in this conjugation.  An unscaled leg has no
such factor.  For nonprimary states the actual coordinate gauge
produces the finite list $L(1)^jb$ of weights $d-j$, and the same
calculation is applied separately to each resulting mode.
At an output $D'$ leg the full restricted-contragredient opposite field
is used instead of an untyped transpose; the root phases, stabilizer
$1/N_i$ and strictly ordered deck pullbacks are those in
\eqref{eq:canonicalmarkedres} and
Lemma~\ref{lem:matchingnodejetlift}.  This is the actual
physical-to-relative residue calculation of
\eqref{eq:collisiongauge76}--\eqref{eq:collisionrelativeward76}, and
it gives, in the coefficientwise relative-frame space,
\begin{equation}\label{eq:matchedrelativeward718}
 \widehat F(t;R_{\rm rel}(s)w)=0,\qquad
 R_{\rm rel}(s)=\mathsf D_t^{-1}R_{\rm phys}(s)\mathsf D_t,
 \qquad \widehat F(t;w)=c_tF(t;\mathsf D_tw)
\end{equation}
for every actual invariant global current $s$ and algebraic external
tensor $w$.  The conjugation formula alone gives no assertion
about integral $t$-orders.  The latter, in the same relative frame,
is supplied only by the independently proved coefficientwise
residue statement of Lemma~\ref{lem:matchingnodejetlift}.

For the finite-level coefficient comparison, fix $\beta$
and the allowed $m=d-1+k/N$, and set
\[
 \gamma_+=\beta,\quad \gamma_-=\beta-ak,\quad
 J=\max\{J_\Sigma(\gamma_+),J_\Sigma(\gamma_-)\}.
\]
For every $0\le j\le J$ retain the finite collection of
nonzero coefficients
$\ell_{\gamma_+-j},\ell_{\gamma_--j}$, where absent exponents mean
zero.  By Theorem~\ref{thm:finitecollisionkernels} each of these
coefficients admits a finite algebraic tensor decomposition with
factors in $\mathscr D_\bullet$.  By the union definition
\eqref{eq:partialnodedualspace}, choose $k_{\max}$ so that all
of these finitely many factors belong to
$\mathscr D_{\bullet,k_{\max}}$.  If a node mode is applied to a
factor in one of the finitely many comparisons, increase
$k_{\max}$ by the finite mode-level bound of
\eqref{eq:nodelevelbounds713}.  This choice uses neither the desired
matching identity nor any later factorization theorem, and is
independent of the subsequently chosen external test tensor.
Choose
\[
 K>\max\{|k|, k_{\max}+\lfloor J/a\rfloor\},\qquad
 aK\ge J+1,
 \quad L_K=K-\lfloor J/a\rfloor>k_{\max}.
\]
Lemma~\ref{lem:matchingnodejetlift} then produces an actual invariant
current $s$ with the integral balanced jet $\sigma_k(u)$. Write
$s_{j,\bullet}$ for its genuine coefficient currents on the
normalized components. For $0\le j\le J$ their node germs differ
from the prescribed principal monomials by local germs
$E_{j,\bullet}$ of order at least $L_K$. We make no assertion that
these differences are global currents. Pair $s_{j,\bullet}$ with
the propagated section of each factor
$\lambda\in\mathscr D_{\bullet,k_{\max}}$ of
$\ell_{\gamma_\pm-j}$. The residue consequence proved in
Lemma~\ref{lem:matchingnodejetlift} applies: the pairing with
$E_{j,\bullet}$ belongs to
$s_\pm^{L_K-k_{\max}}\C[[s_\pm]]ds_\pm$ and has zero node
residue. The global residue theorem for the actual
$s_{j,\bullet}$ therefore expresses its external Ward action
solely by the prescribed principal node residues. This holds for
each factor in the finite tensor decompositions, hence for the
selected coefficients themselves. For a term $t^q$ with
$q\ge J+1$ to contribute to
one of the two selected exponents $\gamma_\pm$, one would need
$\gamma_\pm-q=\sigma+n$ for some $\sigma\in\Sigma,n\ge0$;
then $q\le J_\Sigma(\gamma_\pm)\le J$, a contradiction.
This also shows local finiteness: for a fixed target exponent there
are only finitely many possible integer $q\ge0$; for each such $q$
the residue of the coefficient current on a fixed algebraic tensor
is a finite sum of actual modes.  Logarithms are absent by
\eqref{eq:partialnodelogfree714}; multiplication by integral powers
of $t$ only shifts an exponent inside its fixed $\C/\Z$ class.
In particular no infinite series of endomorphisms has been applied
to a completed tensor product.  This proves the \emph{finite-level}
vanishing of the local error residues for both selected exponents,
not an identity obtained by applying a Ward ideal to a local germ.

We identify the principal node residues and their signs directly.
For an algebraic $w_\bullet$ and
$\lambda\in\mathscr D_{\bullet,k_{\max}}$, apply the global residue
theorem to the actual complete one-free-point propagated dual section
$P\lambda$ of Lemma~\ref{lem:partialequivariantresidue}, paired with
the component current obtained from $s$.  Its external-marked
residue is $\lambda(R_\bullet(s_\bullet)w_\bullet)$; its remaining
residue is the \emph{normalized} node coefficient of
$P\lambda$, with the opposite output field on the inner side.
Their sum is zero.  On the plus branch the balanced density
$e_u=\iota_z(u)z^{d-1}dz$, $z=s_+^N$, satisfies
\[
 \frac1N\Res_{s_+}s_+^{N(d-1)+k}
   \mathcal P_{\lambda,u}(s_+^N;w_\bullet)
   Ns_+^{N-1}ds_+
 =\Res_z z^{d-1+k/N}\mathcal P_{\lambda,u}(z;w_\bullet)dz
 =(u_m\lambda)(w_\bullet),
\]
where the equality is the original normalized-root mode definition,
not a new action.  On the opposite branch the complete state/current
involution
$\Sigma=e^{L(1)}(-1)^{L(0)+1}$, together with reversal of
the node-residue orientation, gives $-C_m(u)$, where
$C_m(u)$ is precisely the finite sum in
\eqref{eq:oppositenodemode}.  This statement is checked mode by mode
by expanding $\Sigma u=\sum_j(-1)^{d+1}
L(1)^ju/j!$ and using the opposite local density
$z^{d-j-1}dz$ from \eqref{eq:balancednodeframe}; its mode label is
$2d-j-m-2$.  The node residue on each physical component is the
root residue divided by the stabilizer order $N$, already included
in \eqref{eq:normalizednoderesidue}.  The full state gauge includes
all $j>0$ terms; neither a second factor $N$ nor a new root phase
appears.  The two component-current residues sum with the signs
shown below, since each individual global residue theorem writes its
external-mark sum as minus its node residue.  All comparisons use
the same selected root, so a noncommuting deck element only changes
both typed source and target of this comparison, never their order.

For $k\ge0$, the \emph{integral} jet is $s_+^ku$ and its opposite
branch representative on the punctured annulus is
$t^{ak}s_-^{-k}\Sigma u$.  After the above marked-versus-node
residue identification and the finite error calculation, the global
relative Ward identity yields at the two selected exponent depths
\begin{equation}\label{eq:matchedpositive718}
 [t^{\gamma_+}]\left((u_m\otimes\id)\widehat\ell(t)
   -t^{ak}(\id\otimes C_m(u))\widehat\ell(t)\right)=0.
\end{equation}
The factor $t^{ak}$ belongs to the \emph{opposite branch}
contribution, not to a physical-frame conjugation.  Hence
\[
 (u_m\otimes\id)\ell_\beta
   =(\id\otimes C_m(u))\ell_{\beta-ak}
 \qquad(k\ge0).
\]
For $k<0$ the \emph{other} jet $s_-^{-k}u$ in the common
plus-branch trivialization is the integral one. Its plus-branch
expression is $t^{-ak}s_+^ku$, while its minus-branch expression
is $s_-^{-k}\Sigma u$. Thus the same principal node residues
are $t^{-ak}u_m$ and $-C_m(u)$, respectively; the full state
transition is applied exactly once. The finite-level Ward comparison
has a different positive power of $t$:
\[
 [t^{\gamma_-}]\left(t^{-ak}(u_m\otimes\id)\widehat\ell(t)
          -(\id\otimes C_m(u))\widehat\ell(t)\right)=0,
 \qquad \gamma_-=\beta-ak.
\]
This is legitimate with the same $J$ and $K$ because
$-ak>0$ and both target exponents were included when $J$ and
$k_{\max}$ were chosen.  Taking this coefficient gives
\[
 (u_m\otimes\id)\ell_\beta
   =(\id\otimes C_m(u))\ell_{\beta-ak}
 \qquad(k<0).
\]
Finally $ak=aN(k/N)=e\alpha$.  The calculation made no assumption
that either of the shifted coefficients is nonzero.  If both are
absent their evaluations vanish; if just one is absent, the same
finite-level Ward equality forces the mode of the other to vanish.
Exponents in an unrelated $\C/\Z$ class do not contribute, as $ak$
is integral.  Since $\beta$ and the test tensor were arbitrary,
finite-dimensional restricted dual separation proves
\eqref{eq:matchedcoefficientidentity} as an equality of algebraic
tensors for all $\beta\in\C$.

The inner $L(0)$ eigenvalue is determined independently by
Lemma~\ref{lem:collisionenergy},
\eqref{eq:partialnodeenergy} and
\eqref{eq:partialnodelogfree714}:
\[
 (\id\otimes L_{\rm in}(0))\ell_\beta
     =(\beta/e)\ell_\beta.
\]
Only now put $u=\omega$, $d=2$, $m=1$, $k=0$ in the
already proved mode relation.  Because $L(1)\omega=0$,
$C_1(\omega)=\omega_1=L_{\rm in}(0)$, and the outer
mode $\omega_1=L_{\rm out}(0)$.  Therefore
\[
 (L_{\rm out}(0)\otimes\id)\ell_\beta
   =(\id\otimes L_{\rm in}(0))\ell_\beta
   =(\beta/e)\ell_\beta.
\]
The ordinary finite semisimple decompositions of the two node
modules give the weight statement
\eqref{eq:matchedcoefficientweights718}, including the case in which
$\beta/e$ is not an occurring weight.

Let $x\in U_w=(\mathscr D_{\rm in})_w^*$, so that
$T_Fx=(\id\otimes x)\ell_{ew}$ is an actual element of the finite
weight space $(\mathscr D_{\rm out})_w$.  An algebraic vector
of $U=\mathscr D_{\rm in}'$ has finite weight support; hence this
defines an algebraic weight-preserving linear map on $U$.
Both $U$ and $\mathscr D_{\rm out}$ are $h$-twisted, since the
inner node has inverse inertia $h^{-1}$ and restricted
contragredience reverses it.  By Lemma~\ref{lem:oppositemodes},
$\deg u_m=d-m-1=-\alpha$, while
$\deg C_m(u)=\alpha$, and for $y\in\mathscr D_{\rm in}$,
$(u_m^Ux)(y)=x(C_m(u)y)$.  Consequently
\[
 \begin{aligned}
 u_m^{\rm out}T_Fx
 &= (\id\otimes x)(u_m\otimes\id)\ell_{ew}\\
 &= (\id\otimes xC_m(u))\ell_{e(w-\alpha)}\\
 &=T_F(u_m^Ux).
 \end{aligned}
\]
Every contraction is legitimate in a finite-dimensional algebraic
weight space.  For an absent exponent $ew$ or $e(w-\alpha)$ the
corresponding coefficient is zero and the same identity applies.
This establishes commutation with every allowed twisted mode and
therefore $T_F\in\Hom_V(U,\mathscr D_{\rm out})$, without appeal to
Schur--Hom decomposition or exact spanning.
\end{proof}

\begin{theorem}\label{thm:partialnodeexactspan}
Every actual product of two complete twisted intertwining operators on the fixed inherited cover is a finite linear combination of actual iterates on the positive product--iterate overlap.  The intermediate objects belong
to $\mathcal C_{g_Ag_B}$, the factors satisfy Definition~\ref{def:completeintertwiner}, and the equality is between the original physical-coordinate, restricted-dual weight sums with their specified branches.
\end{theorem}
\begin{proof}
It suffices initially to take ordinary simple external modules and one
actual product channel with a simple ordinary intermediate module.  A
finite direct sum of intermediate modules gives a finite sum of such
channels; the final extension to arbitrary objects will use only the
finite semisimplicity of Proposition~\ref{prop:regularmodules}.  Retain
throughout the collision roots, based paths, physical coordinates and
positive logarithms fixed before
\eqref{eq:collisionnormalizedkernel}.  Put
\[
 h=g_Ag_B,\quad U=\mathscr D_{\rm in}',\quad
 \mathcal I_S^{\rm in}=\cI\binom S{A\ B},\quad
 \mathcal I_S^{\rm out}=\cI\binom D{S\ C}
 \qquad(S\in\mathcal S_h).
\]
Both $U$ and $\mathscr D_{\rm out}$ belong to the $h$-sector and are
finite ordinary modules by Proposition~\ref{prop:partialnodefiniteordinary}.
By the actual evaluation identifications of
Proposition~\ref{prop:partialnodevertices},
\[
 \Hom_V(U,S)\xrightarrow{\sim}\mathcal I_S^{\rm in},
 \quad f\longmapsto\overline f\circ\mathcal Y^{\rm in},
 \qquad
 \Hom_V(S,\mathscr D_{\rm out})\xrightarrow{\sim}
 \mathcal I_S^{\rm out}.
\]
In the second arrow the operator associated to $g$ is specified,
without a new choice of root, by
$g(s)(d'\otimes c)=\langle d',\mathcal Y_g^{\rm out}(s,1)c\rangle$.
Finite semisimplicity, Schur's lemma and these two evaluation maps
give the finite decompositions
\[
 U\cong\bigoplus_{S\in\mathcal S_h}
        S\otimes(\mathcal I_S^{\rm in})^*,\qquad
 \mathscr D_{\rm out}\cong\bigoplus_{S\in\mathcal S_h}
        S\otimes\mathcal I_S^{\rm out}.
\]
The first multiplicity space is the dual of
$\Hom_V(U,S)=\mathcal I_S^{\rm in}$; the second is
$\Hom_V(S,\mathscr D_{\rm out})=\mathcal I_S^{\rm out}$.
Every space appearing here is finite-dimensional, and
$\Hom_V(S,S')=0$ for nonisomorphic simple $S,S'$; equal conformal
weights do not identify distinct simple modules.

We describe the Hom isomorphism explicitly. For $f\in\Hom_V(U,S)$ and
$g\in\Hom_V(S,\mathscr D_{\rm out})$, send their two associated
complete operators to $g\circ f$.  This defines
\begin{equation}\label{eq:exacthom719}
 \bigoplus_{S\in\mathcal S_h}
 \mathcal I_S^{\rm out}\otimes\mathcal I_S^{\rm in}
 \xrightarrow{\ \sim\ }
 \Hom_V(U,\mathscr D_{\rm out}),
 \qquad \mathcal Y_g^{\rm out}\otimes\mathcal Y_f^{\rm in}
 \longmapsto g\circ f.
\end{equation}
To construct its inverse, use the two displayed finite semisimple
decompositions, take the $S$-diagonal block of a homomorphism
$T:U\to\mathscr D_{\rm out}$, and regard it as a linear map
$(\mathcal I_S^{\rm in})^*\to\mathcal I_S^{\rm out}$.
Finite-dimensional evaluation identifies this linear map with a
unique tensor in
$\mathcal I_S^{\rm out}\otimes\mathcal I_S^{\rm in}$.
These operations are inverse, because on $S$-isotypic blocks a
homomorphism is $\id_S$ tensored with its multiplicity map.
In particular, multiple copies of $S$ are retained in that map.
This description also shows that changing bases changes the tensor
coordinates but not the associated homomorphism or operator pair.

Fix the original product $F$.  Theorem~\ref{thm:matchedcoefficientmodes}
provides an actual $h$-twisted module map
$T_F:U\to\mathscr D_{\rm out}$; no iterate factorization is used
in its construction.  Choose finite bases of the two operator spaces.
Under the inverse of \eqref{eq:exacthom719}, write
\[
 T_F\longleftrightarrow
  \sum_{S,\mu,\nu}c^S_{\mu\nu}
  \mathcal Y^{\rm out}_{S,\mu}\otimes
  \mathcal Y^{\rm in}_{S,\nu}.
\]
The sum is finite because the number of $S$ and both fusion spaces
are finite.  By Proposition~\ref{prop:partialnodevertices}, each
factor is an \emph{actual complete} intertwining operator satisfying all
four conditions of Definition~\ref{def:completeintertwiner}; it is
not merely a three-point Ward functional.  The constants
$c^S_{\mu\nu}$ are the coordinates of the one module map $T_F$:
they are independent of external test vectors, $q$ and $t$.
Their change under a new basis leaves the sum of composed maps,
and hence the eventual finite iterate, unchanged.

We prove the coefficient identity directly from the defining
pairing.  Let $a\in A_\alpha$, $b\in B_\delta$ be homogeneous,
$o=d'\otimes c\in T_{\rm out}$, and $w\in\C$.
The canonical inner evaluation
$j(a,b)\in\overline U$ of \eqref{eq:nodej715} has a well-defined
component $x_w=\pi_w^Uj(a,b)\in U_w$; for
$\lambda\in(\mathscr D_{\rm in})_w$ it satisfies
$\langle\lambda,x_w\rangle=\lambda(a\otimes b)$.
Theorem~\ref{thm:matchedcoefficientmodes} gives
$\ell_{ew}\in(\mathscr D_{\rm out})_w\otimes
(\mathscr D_{\rm in})_w$ and
$T_Fx_w=(\id\otimes x_w)\ell_{ew}$.  Therefore
\[
 \begin{aligned}
 \ell_{ew}(o;a,b)
 &=\bigl(T_F\pi_w^Uj(a,b)\bigr)(o)\\
 &=\sum_{S,\mu,\nu}c^S_{\mu\nu}
    \bigl(g_{S,\mu}f_{S,\nu}
       \pi_w^Uj(a,b)\bigr)(d'\otimes c)\\
 &=\sum_{S,\mu,\nu}c^S_{\mu\nu}
   \left\langle d',\mathcal Y^{\rm out}_{S,\mu}
     \bigl(\pi_w^S\mathcal Y^{\rm in}_{S,\nu}(a,1)b,1\bigr)c
   \right\rangle.
 \end{aligned}
\]
Here $f_{S,\nu}:U\to S$ and
$g_{S,\mu}:S\to\mathscr D_{\rm out}$ are the actual module maps
specified above.  The last equality uses
$\overline f_{S,\nu}\,j(a,b)
 =\mathcal Y^{\rm in}_{S,\nu}(a,1)b$ and
$\overline f_{S,\nu}\pi_w^U=
 \pi_w^S\overline f_{S,\nu}$, followed by the actual outer
node evaluation.  Thus the intermediate projection has exactly the
position displayed; the fixed-$w$ contraction involves only
finite-dimensional weight spaces.  In particular, no sum in
$\overline U$ is treated as an algebraic vector.  The $S$-summands
and all of their multiplicity coordinates remain distinct even if
they have the same $w$.  If $ew$ is outside the collision support,
$\ell_{ew}=0$ by definition, so this formula imposes the same zero
condition on the \emph{sum} of actual iterate coefficients.  If
$\beta\ne ew$ for every node weight $w$, then
$\ell_\beta=0$ by \eqref{eq:matchedcoefficientweights718}; this
also excludes missing exponents.  Lemma~\ref{lem:generation} says
that the homogeneous vectors $\pi_w^Uj(a,b)$ span $U$ as
$a,b,w$ vary, and hence gives uniqueness of the recovered map;
it is not needed to assert convergence of an infinite evaluation.

On the prescribed positive branch set $q=t^e$ with
$\log q=e\log t$ and $\log1=0$.  Ordinary $L(0)$-homogeneity of
each \emph{complete} inner operator gives, component by component,
\[
 \pi_w^S\mathcal Y^{\rm in}_{S,\nu}(a,q)b
 =\exp\bigl((w-\alpha-\delta)\log q\bigr)
  \pi_w^S\mathcal Y^{\rm in}_{S,\nu}(a,1)b,
\]
so that
\[
 q^{\alpha+\delta}\mathcal Y^{\rm in}_{S,\nu}(a,q)b
  =q^{L_S(0)}\mathcal Y^{\rm in}_{S,\nu}(a,1)b,
 \qquad q^{L_S(0)}\pi_w^S=q^w\pi_w^S.
\]
These are equalities in the completed output interpreted
\emph{weightwise}, not assertions that the sum lies in $S$.
No rationality or integrality of $w$ is required and there is no
additional Kummer phase along the positive branch.
Substitution of the already established coefficient pairing shows
that the normalized finite iterate has the \emph{formal} expansion
\[
 \begin{aligned}
 &q^{\alpha+\delta}
 \sum_{S,\mu,\nu}c^S_{\mu\nu}
 \left\langle d',\mathcal Y^{\rm out}_{S,\mu}
  (\mathcal Y^{\rm in}_{S,\nu}(a,q)b,1)c\right\rangle\\
 &\hspace{14mm}\sim
  \sum_{w}q^w\ell_{ew}(o;a,b)
  =\sum_{\beta\in\Sigma+\N}t^\beta\ell_\beta(o;a,b).
 \end{aligned}
\]
The notation $\sim$ at this stage means equality of all
homogeneous coefficients, before asserting an analytic equality.
The second equality uses \eqref{eq:collisioncoefficientkernels},
\eqref{eq:partialnodelogfree714},
\eqref{eq:matchedcoefficientweights718} and $q^w=t^{ew}$.
The union $\Sigma+\N$ consists of disjoint progressions, and the
fixed branch separates their exponent classes.  Thus no exponent
is lost or identified by an unlicensed reordering.

The analytic step follows independently of the preceding formal comparison.  For the original product $F$, its original
intermediate-grade contractions converge absolutely and locally
normally on $|z_1|>|z_2|>0$ by
Corollary~\ref{cor:smallconvergenceautomatic}: its original
operators are complete and its intermediate module is ordinary
(or a finite ordinary direct sum).  For each newly constructed
iterate the two factors above are complete and $S$ is simple
ordinary, so the \emph{same prior corollary}, applied separately to the
iterate chamber, gives absolute local normal convergence of its
original intermediate-weight sum on
$|z_2|>|z_1-z_2|>0$.  The finite combination remains normally
convergent.  In particular, when
$z_1=1+q$, $z_2=1$ and $0<q<1$, both sets of original sums
converge: $1+q>1>q$.  They use the selected roots and the same
positive real logarithms, not branches transported along unrelated
paths.

Moreover, Proposition~\ref{prop:collisionuniformsupport} constructs
for the \emph{original normalized product} an \emph{actually convergent}
finite-exponent-class Nilsson expansion near $t=0$, and
Proposition~\ref{prop:partialnodefiniteordinary} removes all its
positive logarithmic powers.  The iterate's expansion converges
on the same sufficiently small positive $q$-interval by the
preceding normal-convergence result.  Since every coefficient of
the two convergent generalized-power expansions has just been
identified, the two analytic functions are equal for all
sufficiently small $q>0$.  This uses both independent convergence
statements: coefficient equality by itself would not exclude a
nonzero analytic function flat to infinite order at the collision.
The common slice of the product--iterate overlap is the open set
\[
  \bigl\{q\in\C:0<|q|<1,\ |1+q|>1\bigr\};
\]
its component containing the small positive interval has the
common inherited rooted branch.  By the one-variable identity
theorem, the equality extends throughout that based component,
with the continuation of both complete-operator sums fixed from
the positive interval.  Other unrelated monodromy branches are
not identified.

Finally, let $a,b,c,d'$ be homogeneous and put
$\kappa=\wt(d')-\alpha-\delta-\wt(c)$.
For each product contraction, the two $L(0)$-commutators
contribute scaling exponents which add to $\kappa$.
For an iterate through a homogeneous intermediate vector of weight
$w$, the inner operator contributes $w-\alpha-\delta$ and the
outer operator contributes
$\wt(d')-w-\wt(c)$; their sum is again
$\kappa$, independent of $S$ and $w$.  Consequently both
\emph{normally convergent original sums} satisfy, on consistently
continued branches,
\[
 \Phi(z_1,z_2)=
 \exp(\kappa\log z_2)\,
 \Phi(z_1/z_2,1),
 \quad \log\bigl((z_1-z_2)/z_2\bigr)
   =\log(z_1-z_2)-\log z_2
\]
with the logarithms initially chosen on the positive-real overlap.
Scaling is verified termwise on every fixed intermediate weight
before invoking normal convergence; it does not commute an
uncontrolled infinite completion with an operator.
The equality on the $z_2=1$ slice therefore gives, for
$z_2>0$ and $0<z_1-z_2<z_2$, the physical equality
\[
 F(z_1,z_2;o;a,b)
 =\sum_{S,\mu,\nu}c^S_{\mu\nu}
  \left\langle d',\mathcal Y^{\rm out}_{S,\mu}
   \bigl(\mathcal Y^{\rm in}_{S,\nu}(a,z_1-z_2)b,z_2\bigr)c
  \right\rangle.
\]
The holomorphic identity theorem, on the jointly based connected
component of $|z_1|>|z_2|>|z_1-z_2|>0$, extends this equality to
the original product--iterate overlap with its specified rooted
continuation.  There is no sewing-factorization or associativity
assumption here.  Finite homogeneous decomposition extends the
identity to algebraic external vectors, while the finite simple
direct-sum decompositions of Proposition~\ref{prop:regularmodules}
extend it to arbitrary external objects; composition with module
morphisms preserves the concrete evaluation pairings.
\end{proof}

The finite choices in the preceding constructions are summarized below. Throughout, the source family is the actual admissible extension of the framed generic cover; component covers and roots are inherited from it.
\begin{center}
\begingroup\small
\setlength{\tabcolsep}{4pt}
\renewcommand{\arraystretch}{1.12}
\begin{tabular}{@{}P{0.17\linewidth}P{0.27\linewidth}P{0.49\linewidth}@{}}
\toprule
Result & Fixed data & Choice of bounds and comparison \\
\midrule
Lemma~\ref{lem:canonicalnodecurrent}
& One generic current, its finite state space and inherited node roots.
& Clear its poles by $t^{m(s)}$. For each tested coefficient the window of Lemma~\ref{lem:nodewindow} is finite; no bound uniform in the current is used. \\[3pt]
Lemma~\ref{lem:truncatedcomponentlift}
& Finite state cutoff $d$, parameter depth $J$ and node level $k$.
& First require $a_jk\ge J+1$ at every node; then choose the external pole bound by Serre vanishing. The congruence is in relative component frames. \\[3pt]
Lemma~\ref{lem:partialpropagationward}
& One normalized component, level $k$, finitely many states and algebraic tests.
& Enlarge the finite state and incoming pole bounds for these tests. Finite-jet uniqueness gives parameter holomorphy and the full OPE; compatible cutoffs retain each regular coefficient. \\[3pt]
Theorem~\ref{thm:matchedcoefficientmodes}
& One exponent $\beta$ and one mode, with root degree $k=N(m-d+1)$.
& Choose $J$ for $\beta,\beta-ak$ and a common tested level $k_{\max}$; then take $K>\max\{|k|,k_{\max}+\lfloor J/a\rfloor\}$ and $aK\ge J+1$. Choose the external pole bound last. Apply the physical relative-weight conjugation before coefficient extraction. \\
\bottomrule
\end{tabular}
\endgroup
\end{center}

\section{Associativity and the elementary structure maps}
\label{sec:associativity}
\subsection{Associativity from actual iterate factorization}
Fix homogeneous $A\in\cC_g$, $B\in\cC_h$, $C\in\cC_k$ and put $U=A\bo B,\; W=B\bo C,\; L=U\bo C,\; R=A\bo W$. We use the natural identifications
\begin{equation}\label{eq:assocHom81}
 \mathsf P_{ABC}^{D}\cong\Hom_V(R,D),\qquad
 \mathsf I_{ABC}^{D}\cong\Hom_V(L,D),
 \quad D\in\cC_{ghk},
\end{equation}
from \eqref{eq:pairHom42}.  For $f:R\to D$ and $v:L\to D$, denote their original restricted-dual scalar realizations by $\Phi_D^P(f)$ and $\Phi_D^I(v)$, respectively. These converge in their respective chambers by
Corollary~\ref{cor:smallconvergenceautomatic}. The universal $P(z)$-map and its Hom characterization are those of \cite[Section~4]{DH}, as realized here by Theorem~\ref{thm:binary}.

\begin{lemma}\label{lem:productepi}
The factorization of Theorem~\ref{thm:partialnodeexactspan} induces a unique natural epimorphism
\[
 b_{A,B,C}:L\longrightarrow R
\]
characterized on the positive product--iterate overlap by
\begin{equation}\label{eq:onewayscalar}
 \overline b_{A,B,C}\,
 \mathcal Y^{\can}_{A\bo B,C}
 (\mathcal Y^{\can}_{A,B}(a,z_1-z_2)b,z_2)c
 =\mathcal Y^{\can}_{A,B\bo C}(a,z_1)
 \mathcal Y^{\can}_{B,C}(b,z_2)c
\end{equation}
after evaluation against every restricted-dual vector.
\end{lemma}
\begin{proof}
Write $U=A\bo B$, $W=B\bo C$, $L=U\bo C$, and
$R=A\bo W$.  All four objects and every permitted output
$D\in\cC_{ghk}$ are finite direct sums of ordinary modules in
the indicated sectors.  The natural algebraic pair-space
identifications \eqref{eq:pairHom42} are obtained by applying the
binary universal property twice, before defining any three-fold
associator.  They represent finite sums of \emph{actual complete}
operator pairs, since the fusion spaces are finite-dimensional.
We regard a scalar realization as the family of its original
restricted-dual matrix coefficients for all homogeneous algebraic
inputs, on the inherited positive-rooted overlap.  The relevant
original intermediate-weight sums are absolutely and locally
normally convergent by Corollary~\ref{cor:smallconvergenceautomatic}.

We prove iterate faithfulness.  Let
$f\in\Hom_V(L,D)$ and suppose its entire iterate scalar family
vanishes on the positive overlap.  Fix homogeneous
$a\in A_\alpha$, $b\in B_\beta$, $c\in C_\gamma$ and
$d'\in D'_\delta$; put $\kappa=\delta-\alpha-\beta-\gamma$.
At $z_2=R_0>0$, $z_1=R_0(1+q)$, $0<q<1$, its original
convergent weight sum is
\begin{equation}\label{eq:iteratefaithfulness81}
 R_0^\kappa\sum_{\mu\in\operatorname{wt}(U)}
 q^{\mu-\alpha-\beta}
 \big\langle d',f\mathcal Y^{\can}_{U,C}
 (\pi^U_\mu\mathcal Y^{\can}_{A,B}(a,1)b,1)c\big\rangle.
\end{equation}
Here $\operatorname{wt}(M)$ is the weight set of an ordinary module
and $\pi^U_\mu:\overline U\to U_\mu$ is its finite-dimensional
weight projection.  This is an equality of scalar analytic
functions, not an application of an uncompleted operator to an
infinite sum.  For some positive integer $N$ and finitely many
lowest weights $\lambda_j$,
\[
 \operatorname{wt}(U)\subseteq
 \bigcup_{j=1}^r(\lambda_j+N^{-1}\Z_{\ge0}).
\]
Thus after $q=\xi^N$ and a fixed logarithm, the exponents
$N(\mu-\alpha-\beta)$ lie in finitely many cosets of $\Z$,
with a lower integer bound in each coset.  Group \emph{equal}
exponents, including those coming from different simple summands,
before taking Laurent coefficients.  Distinct cosets modulo $\Z$
are eigencharacters of continuation $\xi\mapsto e^{2\pi i}\xi$;
finite-dimensional Vandermonde projection onto the distinct
characters separates their convergent generalized Laurent series.
Within each coset the usual uniqueness of Laurent coefficients
separates every distinct $\mu$.  This remains valid when the
$\lambda_j$ are not rational and when their differences are
integers.  Hence the vanishing family implies
\[
 \big\langle d',f\mathcal Y^{\can}_{U,C}
 (\pi^U_\mu\mathcal Y^{\can}_{A,B}(a,1)b,1)c\big\rangle=0
 \qquad(\text{every }\mu,a,b,c,d').
\]
By Lemma~\ref{lem:generation}, the displayed inner projected
vectors span $U_\mu$ for each $\mu$, and the projected
canonical coefficients for $(U,C)$ span $L$.  The restricted
dual of $D$ separates $D$, so $f=0$.

For \emph{product faithfulness} take
$f\in\Hom_V(R,D)$ and $z_1=R_0>0$, $z_2=R_0\rho$.
Its original weight sum is
\[
 R_0^\kappa\sum_{\nu\in\operatorname{wt}(W)}
 \rho^{\nu-\beta-\gamma}
 \big\langle d',f\mathcal Y^{\can}_{A,W}
 (a,1)\pi^W_\nu\mathcal Y^{\can}_{B,C}(b,1)c\big\rangle.
\]
The actual positive-real product--iterate overlap on this slice
is $\tfrac12<\rho<1$: here
$z_1>z_2>z_1-z_2>0$.  If the scalar family vanishes on
that overlap, analytic uniqueness \emph{within the product
chamber}, where this same original weight sum converges for
$0<\rho<1$, extends its vanishing to $0<\rho<1$.
We may therefore analyze the convergent expansion at $\rho=0$.
The same finite-coset and monodromy--Laurent uniqueness argument
shows that each $\nu$-coefficient vanishes.  Lemma~\ref{lem:generation}
first spans the individual $W_\nu$, then the canonical
coefficients for $(A,W)$ span $R$; hence $f=0$.
No associator or inverse factorization enters either injectivity
argument.  In particular, coincident conformal weights and
multiplicities cause no loss of coefficients because the \emph{full}
$U_\mu$ or $W_\nu$ is tested.

For $f\in\Hom_V(R,D)$, the exact-spanning theorem applied to
each actual complete product pair in its finite expansion via
\eqref{eq:pairHom42} produces a finite sum of actual complete
iterate pairs.  The iterate identification in that same equation
assembles this finite sum into a map $F_D(f):L\to D$ whose
original normally convergent scalar family satisfies
\begin{equation}\label{eq:factorizationFD81}
 \Phi_D^I(F_D(f))=\Phi_D^P(f)
\end{equation}
on the specified positive overlap.  This is not merely an
identity of formal collision coefficients.  Iterate faithfulness
proves uniqueness of $F_D(f)$; it therefore proves its linearity.
For $u:D\to D'$, composition with $u$ preserves all typed
complete operators, original weight projections and the inherited
rooted branch.  Both $uF_D(f)$ and $F_{D'}(uf)$ consequently
have iterate realization $\Phi_{D'}^P(uf)$, whence
$F_{D'}(uf)=uF_D(f)$.  This proves naturality in $D$.

Take $D=R$ and define
$ b_{A,B,C}:=F_R(\id_R)\in\Hom_V(L,R)$.
With $u=f:R\to D$ in the established output naturality,
\[
 F_D(f)=f\circ b_{A,B,C}.
\]
The defining equation for $D=R$ is precisely
\eqref{eq:onewayscalar} after restricted-dual evaluation.
Iterate faithfulness gives uniqueness of $b$ among morphisms
satisfying this scalar equation.  The completion $\overline b$
is legitimate because $b$ preserves all ordinary weights and
acts componentwise on the corresponding products.

For a triple of morphisms
$\mathbf u=(u_A:A\to A',u_B:B\to B',u_C:C\to C')$,
let $L(\mathbf u)=(u_A\bo u_B)\bo u_C$ and
$R(\mathbf u)=u_A\bo(u_B\bo u_C)$, defined solely by the
binary representing property.  Its defining equation identifies
the two canonical operator pairs after substituting
$u_Aa,u_Bb,u_Cc$.  Therefore
$R(\mathbf u)b_{A,B,C}$ and
$b_{A',B',C'}L(\mathbf u)$, both maps $L\to R'$, have
identical original iterate realizations.  Iterate faithfulness
forces
\[
 R(\mathbf u)b_{A,B,C}=b_{A',B',C'}L(\mathbf u).
\]
All maps have sector $ghk$ at their respective targets; no
noncommuting sector labels have been interchanged.  The
finite direct-sum case follows by the biadditivity of the
binary representing functor, which follows from additivity of
complete-operator spaces and does not use associativity.

Lastly let $p:R\twoheadrightarrow Q:=\operatorname{coker}b$.
The category $\cC_{ghk}$ is finite semisimple abelian, so
$Q\in\cC_{ghk}$ and both scalar realizations with output $Q$
are defined.  As $pb=0$, the already established equation with
$f=p$ gives
\[
 \Phi_Q^P(p)=\Phi_Q^I(pb)=0.
\]
Product faithfulness implies $p=0$.  Since the cokernel
projection $p$ is epic, $Q=0$; thus
$ b_{A,B,C}:L\twoheadrightarrow R$ is an epimorphism.
No reverse exact spanning, associativity or categorical
coherence has been invoked.
\end{proof}

\begin{lemma}\label{lem:finiterotation}
Let $\mathcal A=\bigoplus_{g\in G}\mathcal A_g$ be a $G$-graded finite-length abelian category, with $G$ finite, a biadditive product $\bo$ satisfying $\mathcal A_g\bo\mathcal A_h\subseteq\mathcal A_{gh}$, and an action $T_s$ satisfying
$T_sT_t=T_{st}$, $T_1=\id$ and $T_s(\mathcal A_g)\subseteq\mathcal A_{sgs^{-1}}$. Suppose there are isomorphisms $X\bo Y\xrightarrow{\ \sim\ }T_gY\bo X\; (X\in\mathcal A_g)$ and epimorphisms
$ b_{X,Y,Z}:(X\bo Y)\bo Z\twoheadrightarrow
X\bo(Y\bo Z)$ for all homogeneous triples.
Then every $b_{X,Y,Z}$ is an isomorphism. Neither the pentagon nor either crossed hexagon is assumed.
\end{lemma}
\begin{proof}
For a homogeneous triple $\mathbf X_i=(X_i,Y_i,Z_i)$ with
labels $(g_i,h_i,k_i)$ define
\[
 \rho(\mathbf X_i)=(T_{g_ih_i}Z_i,X_i,Y_i),\qquad
 \mathbf X_{i+1}=\rho(\mathbf X_i),\qquad
 L_i=(X_i\bo Y_i)\bo Z_i,\quad R_i=X_i\bo(Y_i\bo Z_i).
\]
The new first label is $(g_ih_i)k_i(g_ih_i)^{-1}$,
so the total label of $\mathbf X_{i+1}$ is
$ (g_ih_i)k_i(g_ih_i)^{-1}g_ih_i=g_ih_ik_i$,
with the displayed noncommutative order.  Since
$X_i\bo Y_i\in\mathcal A_{g_ih_i}$, crossed
\emph{binary} commutativity supplies the precisely typed
isomorphism
\[
 C_i=c_{X_i\bo Y_i,Z_i}:
 L_i=(X_i\bo Y_i)\bo Z_i
 \xrightarrow{\ \sim\ }
 T_{g_ih_i}Z_i\bo(X_i\bo Y_i)=R_{i+1}.
\]
The parenthesized object $X_i\bo Y_i$ is used as one input of
a binary isomorphism; no reassociation has been performed.
Write $b_i=b_{X_i,Y_i,Z_i}:L_i\twoheadrightarrow R_i$.
For composition length $\ell$ we obtain
\begin{equation}\label{eq:lengthrotation}
 \ell(R_i)\le\ell(L_i)=\ell(R_{i+1}).
\end{equation}

We establish periodicity without an associator.  The three
entries of $\rho(\mathbf X)$ are a cyclic permutation of
those of $\mathbf X$, with one entry replaced by a $G$-translate.
Inductively, every entry in every $\mathbf X_i$ is isomorphic
to $T_sX_0$, $T_sY_0$ or $T_sZ_0$ for some $s\in G$:
indeed $T_sT_t=T_{st}$, and permutation introduces no new
object.  There are at most $3!|G|^3$ ordered triple
isomorphism classes of these forms, even for nonabelian $G$.
Rotation is invertible on such classes: explicitly its inverse is
\[
 \rho^{-1}(A,B,C)=(B,C,T_{bc}^{-1}A),
 \qquad b=\deg B,\ c=\deg C.
\]
The two formulas are inverse because $\deg B\deg C=bc$
and $T_{bc}T_{bc}^{-1}=\id$.  Hence the initial triple lies on
a finite cycle, not merely a preperiodic tail.  Choose $r>0$
and componentwise isomorphisms
$\mathbf X_r\simeq\mathbf X_0$.  Biadditivity (indeed
ordinary functoriality) of $\bo$ induces
$\eta:R_r\xrightarrow{\sim}R_0$, without any rebracketing.
Thus $\ell(R_r)=\ell(R_0)$, and summing
\eqref{eq:lengthrotation} for $i=0,\ldots,r-1$ yields
\[
 0\le\sum_{i=0}^{r-1}
  \bigl(\ell(L_i)-\ell(R_i)\bigr)
  =\sum_{i=0}^{r-1}
  \bigl(\ell(R_{i+1})-\ell(R_i)\bigr)=0.
\]
Every nonnegative summand is zero.  Since each $b_i$ is epic,
its kernel has length
$\ell(\ker b_i)=\ell(L_i)-\ell(R_i)=0$.
Thus every $b_i$ is monic and epic and hence an isomorphism
in the abelian category.

An explicit inverse also follows from
only the already given epis and binary isomorphisms.
Define
\[
 E_i=b_iC_i^{-1}:R_{i+1}\twoheadrightarrow R_i,
 \qquad
 Q=C_0^{-1}E_1\cdots E_{r-1}\eta^{-1}:
 R_0\twoheadrightarrow L_0,
\]
where the $E$-product is empty for $r=1$ and composition is
read right to left.  For $r=1$, $Q=C_0^{-1}\eta^{-1}$; for
$r=2$, the typed chain is
$R_0\xrightarrow{\eta^{-1}}R_2
 \xrightarrow{E_1}R_1\xrightarrow{C_0^{-1}}L_0$.
Every factor is epic or invertible; thus $Q$ is epic.
The map $H=b_0Q:R_0\twoheadrightarrow R_0$ is an epic
endomorphism of a finite-length object.  The exact sequence
$0\to\ker H\to R_0\xrightarrow{H}R_0\to0$
forces $\ell(\ker H)=0$, so $H$ is invertible.  Put
$c=QH^{-1}:R_0\to L_0$.  Then
\[
 b_0c=\id_{R_0},\qquad
 (cb_0)Q=cH=Q.
\]
As $Q$ is epic, $cb_0=\id_{L_0}$.
This proves the displayed formula is a genuine two-sided inverse;
its resulting morphism is independent of the choice of $r$,
$\eta$ and the binary crossed isomorphisms by uniqueness of
inverses.

In the present application each sector of $\cC_G$ is a
finite semisimple abelian category by
Proposition~\ref{prop:regularmodules}, and $G$ is finite.
The binary object assignment of Theorem~\ref{thm:binary}
is biadditive: on each variable a direct-sum input yields a
direct sum of complete-operator spaces, hence a direct sum of
the corresponding representing objects, uniquely by the binary
universal property.  It preserves sector products, and
$T_sT_t=T_{st}$ strictly as functors by the preliminary
convention.  The maps $c_{X,Y}$ of
Lemma~\ref{lem:units}, obtained from actual Du--Huang
skew-symmetry and positive-path $P(z)$-transport, are typed
binary isomorphisms, without using any triple associator.
Lemma~\ref{lem:productepi} provides the required epis on every
homogeneous triple.  Therefore the abstract argument applies;
for nonhomogeneous objects it applies componentwise to their
finite direct-sum decomposition.  No action tensorator coherence,
pentagon, hexagon or rigidity is used.
\end{proof}

\begin{theorem}\label{thm:smallassoc}
For $A\in\cC_g$, $B\in\cC_h$ and $C\in\cC_k$, put $L=(A\bo B)\bo C$ and $R=A\bo(B\bo C)$. The natural epimorphism $b_{A,B,C}:L\to R$ of
Lemma~\ref{lem:productepi} is an isomorphism.  Set $a_{A,B,C}:=b_{A,B,C},\; c_{A,B,C}:=b_{A,B,C}^{-1}:R\longrightarrow L$. For $D\in\cC_{ghk}$, the mutually inverse natural maps
\begin{equation}\label{eq:assocdouble83}
\begin{aligned}
 \mathsf P_{ABC}^{D}\longrightarrow\mathsf I_{ABC}^{D},
 &\qquad f\longmapsto f\circ b_{A,B,C},\\
 \mathsf I_{ABC}^{D}\longrightarrow\mathsf P_{ABC}^{D},
 &\qquad v\longmapsto v\circ c_{A,B,C}
\end{aligned}
\end{equation}
identify the original product and iterate weight sums:
\begin{equation}\label{eq:assocdoubleidentity83}
 \Phi_D^I(fb_{A,B,C})=\Phi_D^P(f),\qquad
 \Phi_D^P(vc_{A,B,C})=\Phi_D^I(v).
\end{equation}
The equalities hold normally on the positive product--iterate overlap with compatible logarithms, and after \emph{common} analytic continuation along a collision-free path with its framing retained.  In the common formal binary representatives,
$a_{A,B,C}$ also equals the positive-real $P(z)$-transport normalization \eqref{eq:positivenormalization}.
\end{theorem}
\begin{proof}
Fix $A\in\cC_g$, $B\in\cC_h$, $C\in\cC_k$, and put
$U=A\bo B$, $W=B\bo C$, $L=U\bo C$, $R=A\bo W$.
All four objects lie in their ordered sectors, and both $L,R$
lie in $\cC_{ghk}$.  Lemmas~\ref{lem:productepi}
and~\ref{lem:finiterotation} give $b:L\xrightarrow{\sim}R$
and $c=b^{-1}:R\xrightarrow{\sim}L$.  For every
$D\in\cC_{ghk}$, the independent binary representing
identifications \eqref{eq:assocHom81} are
$\mathsf P_{ABC}^D\simeq\Hom_V(R,D)$ and
$\mathsf I_{ABC}^D\simeq\Hom_V(L,D)$.  They are not
identified in advance.  Postcomposition with $b$ gives
$f:R\to D\mapsto fb:L\to D$, and postcomposition with $c$
gives $v:L\to D\mapsto vc:R\to D$; indeed
$(fb)c=f$ and $(vc)b=v$.  This proves the two inverse
maps \eqref{eq:assocdouble83} with their asserted types.

The first scalar identity is the original normally
convergent equality \eqref{eq:factorizationFD81}, since
$F_D(f)=fb$ by Lemma~\ref{lem:productepi}.  Applying it to
$f=vc:R\to D$ gives
\[
 \Phi_D^P(vc)=\Phi_D^I((vc)b)=\Phi_D^I(v).
\]
Thus the reverse identity is a consequence of $cb=\id_L$;
no reverse exact-spanning assertion is being used.  Both
realizations here mean families of matrix coefficients of the
original restricted-dual intermediate-weight sums,
not abstract continuation solutions.  Their injectivity, already
proved separately for product and iterate in
Lemma~\ref{lem:productepi}, shows that the correspondence is
uniquely determined by these scalar equalities.  No equality
of conformal weights of distinct intermediate simple modules
is needed.

Fix homogeneous $a\in A_\alpha$, $b_0\in B_\beta$,
$c_0\in C_\gamma$, $d'\in D'_\delta$, and put
$\kappa=\delta-\alpha-\beta-\gamma$.
For $R_0>0$ and $0<q<1$ set
$z_2=R_0$, $z_1=R_0(1+q)$, so that
$z_1-z_2=R_0q$.  By the $L(0)$-commutator of the actual
complete operators, an inner $B\bo C$ component of weight
$\nu$ contributes $R_0^{\nu-\beta-\gamma}$ at $z_2$
and the outer operator contributes
$[R_0(1+q)]^{\delta-\alpha-\nu}$.  Similarly an inner
$A\bo B$ component of weight $\mu$ contributes
$(R_0q)^{\mu-\alpha-\beta}$, and its outer operator contributes
$R_0^{\delta-\mu-\gamma}$.  Their products have the same
$R_0^\kappa$ factor.  Consequently the two original sums are
\begin{align}
 \Phi_D^P(f)
 &=R_0^\kappa
 \sum_{\nu\in\operatorname{wt}(B\bo C)}
 (1+q)^{\delta-\alpha-\nu}
 \big\langle d',f\,\mathcal Y^{\can}_{A,B\bo C}(a,1)
 \pi_\nu^{B\bo C}\mathcal Y^{\can}_{B,C}(b_0,1)c_0\big\rangle,
 \label{eq:assocproduct83}\\
 \Phi_D^I(v)
 &=R_0^\kappa
 \sum_{\mu\in\operatorname{wt}(A\bo B)}
 q^{\mu-\alpha-\beta}
 \big\langle d',v\,\mathcal Y^{\can}_{A\bo B,C}
 (\pi_\mu^{A\bo B}\mathcal Y^{\can}_{A,B}(a,1)b_0,1)c_0
 \big\rangle.
 \label{eq:associterate83}
\end{align}
Every fixed $\nu$ or $\mu$ term is a finite-dimensional
algebraic contraction; the vectors evaluated by the inner
complete operator before projection belong in its \emph{completed}
output and need not belong to its algebraic direct sum.  Finite
sums over homogeneous components handle nonhomogeneous algebraic
inputs.  All fractional and possibly nonrational conformal powers
are evaluated with the indicated real logarithms.

The common physical chamber is
\[
 \mathcal O=\{(z_1,z_2)\in\C^2:
 |z_1|>|z_2|>|z_1-z_2|>0\}.
\]
The biholomorphic change $(z_1,z_2)=(z_2(1+q),z_2)$
identifies it with
$\C^*\times\{0<|q|<1,\ |1+q|>1\}$.
For $z_2=R_0>0$ and $q\in(0,1)$ the inequalities hold;
indeed at $(3,2)$ one has $q=1/2$.
Corollary~\ref{cor:smallconvergenceautomatic}, applied to
the actual complete operator pairs obtained from the
binary representatives, supplies absolute and local normal
convergence of each unregrouped original weight sum on its
own product or iterate chamber, hence simultaneously on
$\mathcal O$.  No associativity assertion is used in this
convergence input.  The inherited rooted logarithms on the
positive slice satisfy
\[
 \log z_2=\log R_0,\quad
 \log(z_1-z_2)=\log R_0+\log q,\quad
 \log z_1=\log R_0+\log(1+q).
\]
On a simply connected sheet of $\mathcal O$ with these
continued logarithms, $L(0)$-homogeneity identifies the two
holomorphic scalar functions with $z_2^\kappa$ times their
respective functions of $q$ from the above original sums.
The scalar equality proved by Lemma~\ref{lem:productepi}
holds for all sufficiently small positive $q$ and all such
external tests.  The one-variable identity theorem first extends
it to the connected $q$-sheet; the common homogeneity factor
then extends it in $z_2$ on the same two-variable logarithmic
sheet.  Equivalently, the identity theorem can be applied
successively in $q$ and $z_2$.  This proves
\eqref{eq:assocdoubleidentity83} as an equality of the
original convergent grade sums throughout the designated
sheet, not just equality of formal collision expansions.

For $\mathbf u=(u_A:A\to A',u_B:B\to B',u_C:C\to C')$,
let $L(\mathbf u)=(u_A\bo u_B)\bo u_C$ and
$R(\mathbf u)=u_A\bo(u_B\bo u_C)$, using only binary
functoriality.  Lemma~\ref{lem:productepi} gives
\[
 R(\mathbf u)b_{A,B,C}=b_{A',B',C'}L(\mathbf u),\qquad
 L(\mathbf u)c_{A,B,C}=c_{A',B',C'}R(\mathbf u).
\]
For any output morphism $h:D\to D'$, one has
$h(fb)=(hf)b$ and $h(vc)=(hv)c$.  These identities give
naturality in all four variables; finite semisimple direct sums
extend it to nonhomogeneous objects without invoking any
coherence condition.  By Theorems~\ref{thm:fullsmallward}
and~\ref{thm:E0covariantcontinuation}, the two original
scalar realizations, on a common framed placement chart,
solve the same finite-rank Ward connection and possess
covariant analytic continuation along collision-free paths.
In fact, once their equality of holomorphic germs has been
established above, uniqueness of analytic continuation already
forces equality along \emph{one and the same lifted path}, with
its continued roots, deck frames and logarithms.  It imposes
no equality between functions reached along different paths;
nontrivial monodromy remains possible.  Output and input module
maps commute with the placements and their $L(-1)$ connection
operators, so these comparisons preserve naturality.

Use the fixed formal representative $U_{X,Y}$ and its complete
canonical operator $\mathcal Y^{\can}_{X,Y}$ of
Theorem~\ref{thm:binary} at every $P(t)$, $t>0$.  Although
the representing \emph{objects} are then literally the same,
the canonical $P(t)$-maps are the distinct evaluations
$I^t_{X,Y}=\operatorname{ev}_{\log t}(\mathcal Y^{\can}_{X,Y})$.
For the positive-real interval from $s$ to $t$,
\eqref{eq:Pztransport} is satisfied by
$\mathsf T_{s,t}^{X,Y}=\id_{U_{X,Y}}$, since its two
sides are precisely the evaluation at $s$ of this same
formal operator.  By the $P(s)$-universal property this
transport is \emph{uniquely} the identity.  It is not claimed
to be the identity for unrelated independently chosen
representing objects or for paths with winding.

At the actual positive configuration $(z_1,z_2,0)=(3,2,0)$,
put
$L_{3,2}=(A\bo_1B)\bo_2C$ and
$R_{3,2}=A\bo_3(B\bo_2C)$.
The two canonical scalar contractions are, respectively,
\[
 \left\langle d',\,
 \mathcal Y^{\can}_{A\bo_1B,C}
 (\mathcal Y^{\can}_{A,B}(a,1)b_0,2)c_0\right\rangle,
 \qquad
 \left\langle d',\,
 \mathcal Y^{\can}_{A,B\bo_2C}(a,3)
 \mathcal Y^{\can}_{B,C}(b_0,2)c_0\right\rangle,
\]
with their original intermediate-weight projections understood;
the two displayed expressions are not uncompleted products
of arbitrary formal infinite vectors.  Via the common formal
representatives, $L_{3,2}=L$ and $R_{3,2}=R$ as modules.
The standard product-to-iterate comparison, traditionally
written $\mathcal A_{3,2}^{-1}:L_{3,2}\to R_{3,2}$,
is characterized by identifying the two original analytic
scalar germs on the positive-rooted overlap in a neighbourhood
of $(3,2)$, for every algebraic input and restricted-dual test.
Its evaluation at $(3,2)$ identifies the displayed two
canonical contractions.  Here its \emph{existence} is supplied
by the previously constructed $b$, not by a pre-existing
associativity isomorphism.  Its uniqueness follows from the
iterate scalar faithfulness in Lemma~\ref{lem:productepi}
with output $R_{3,2}$: two maps $L_{3,2}\to R_{3,2}$
with that full-germ property have zero difference under the
original iterate realization.  Equality at only one numerical
input and one placement would not suffice for this inference.  Equation~\eqref{eq:onewayscalar} at
$(3,2)$ therefore gives
\begin{equation}\label{eq:assocnormalization83}
 \mathcal A_{3,2}^{-1}=b_{A,B,C}
 \quad\text{in the fixed common formal representatives}.
\end{equation}
This argument defines and identifies the comparison before
using the positive-real formula \eqref{eq:positivenormalization}.

Finally, the four arrows of that formula have the following
\emph{literal} domains and codomains:
\[
\begin{gathered}
 (A\bo_1B)\bo_1C
 \xrightarrow{\ \mathsf T_{1,2}^{A\bo_1B,C}\ }
 (A\bo_1B)\bo_2C
 \xrightarrow{\ \mathcal A_{3,2}^{-1}\ }
 A\bo_3(B\bo_2C)\\
 \xrightarrow{\ \id_A\bo_3\mathsf T_{2,1}^{B,C}\ }
 A\bo_3(B\bo_1C)
 \xrightarrow{\ \mathsf T_{3,1}^{A,B\bo_1C}\ }
 A\bo_1(B\bo_1C).
\end{gathered}
\]
All terms lie in $\cC_{ghk}$ with the prescribed ordered
sector product.  Each $\mathsf T$ is the identity on its fixed
formal representative, and by binary functoriality so is
$\id_A\bo_3\mathsf T_{2,1}^{B,C}$.  Their composite is
therefore $\mathcal A_{3,2}^{-1}=b_{A,B,C}$.
Thus the \emph{independently specified} positive-real
normalization \eqref{eq:positivenormalization} agrees with
$a_{A,B,C}=b_{A,B,C}$, without assuming a pentagon,
triangle, crossed hexagon or any other associator coherence.
\end{proof}

\subsection{The unit triangle and the monoidal group action}
\label{sec:unittriangle}
For $t>0$ write $\bo_t=\bo_{P(t)}$ and $\mathsf T_{s,t}^{X,Y}:X\bo_sY\xrightarrow{\sim}X\bo_tY$ for the transport of \eqref{eq:Pztransport} along $[s,t]\subset\R_{>0}$,
with $\log s,\log t\in\R$.  For $X,Y,Z\in\cC_G$ put $\mathcal A_{3,2}^{-1}:(X\bo_1Y)\bo_2Z \xrightarrow{\ \sim\ }X\bo_3(Y\bo_2Z)$. By Theorem~\ref{thm:smallassoc}, this is the inverse product-to-iterate map characterized by the original
weight sums.  The associator is fixed by
\begin{equation}\label{eq:positivenormalization}
 a_{X,Y,Z}=
 \mathsf T_{3,1}^{X,Y\bo_1Z}
 \circ(\id_X\bo_3\mathsf T_{2,1}^{Y,Z})
 \circ\mathcal A_{3,2}^{-1}
 \circ\mathsf T_{1,2}^{X\bo_1Y,Z}.
\end{equation}
The positive half-braid is \eqref{eq:c}.  The normalization agrees with the $P(z)$ convention of
\cite[Sections~12.2--12.3]{HLZ}; the twisted $P(z)$
transport used here is that of \cite[Section~4]{DH}.

\begin{proposition}\label{prop:unittriangle}
For $X,Y\in\cC_G$,
\begin{equation}\label{eq:triangle82}
 (\id_X\bo l_Y)\circ a_{X,V,Y}
 =r_X\bo\id_Y:(X\bo V)\bo Y\longrightarrow X\bo Y.
\end{equation}
Moreover, $l_V=r_V:V\bo V\to V$.
\end{proposition}
\begin{proof}
For $t>0$ define $l_Y^t:V\bo_tY\to Y$ and
$r_X^t:X\bo_tV\to X$ by the actual complete operators
$Y_Y^g$ and $\Gamma_X=\Omega_-Y_X^g$:
\[
 \overline{l_Y^t}I^t_{V,Y}(v\otimes y)=Y_Y(v,t)y,
 \qquad
 \overline{r_X^t}I^t_{X,V}(x\otimes v)=\Gamma_X(x,t)v.
\]
The first target is in the sector of $Y$, and the second in the
sector of $X$.  Skew symmetry uses the fixed $\Omega_-$ branch;
for the vacuum $\Gamma_X(x,t)\mathbf1=e^{tL_X(-1)}x$
in $\overline X$.  Both displayed operators satisfy the complete
intertwiner axioms, so the $P(t)$ universal property gives these
maps without invoking a triangle.  For $s,t>0$ apply the transport
defining equation \eqref{eq:Pztransport} before applying the unit
maps.  In the first case it gives, with both sides evaluated on
$I^s_{V,Y}(v\otimes y)$,
\[
 \overline{l_Y^t\mathsf T_{s,t}^{V,Y}}I^s_{V,Y}(v\otimes y)
 =Y_Y(v,s)y
 =\overline{l_Y^s}I^s_{V,Y}(v\otimes y).
\]
For $r_X$ the same calculation gives
$\Gamma_X(x,s)v$ on both sides, with the same $\Omega_-$
branch.  By uniqueness in \eqref{eq:Pzuniversal},
\begin{equation}\label{eq:unittransport82}
 l_Y^t\mathsf T_{s,t}^{V,Y}=l_Y^s,
 \qquad r_X^t\mathsf T_{s,t}^{X,V}=r_X^s.
\end{equation}
The transports are isomorphisms, and so are $l_Y^t,r_X^t$
by Lemma~\ref{lem:units}.

At $(z_1,z_2)=(3,2)$ put
\begin{equation}\label{eq:rawtriangle82}
\underbrace{(\id_X\bo_3 l_Y^2)\mathcal A_{3,2}^{-1}}_{L}
=\underbrace{\mathsf T_{2,3}^{X,Y}
 (r_X^1\bo_2\id_Y)}_{R}:
 (X\bo_1V)\bo_2Y\longrightarrow X\bo_3Y.
\end{equation}
On the left the intermediate arrows have targets
$X\bo_3(V\bo_2Y)$ and $X\bo_3Y$; on the right they have
targets $X\bo_2Y$ and $X\bo_3Y$.  No completed vector is
inserted as an argument of either of these module maps.

Fix algebraic homogeneous $x\in X_\alpha$, $y\in Y_\beta$,
and $w'\in(X\bo_3Y)'_\delta$.  Use the \emph{legal algebraic}
inner inputs $(x,\mathbf1)$ and the original restricted-dual
contraction of
$I^2_{X\bo_1V,Y}(I^1_{X,V}(x\otimes\mathbf1)\otimes y)$,
with its intermediate-weight projectors.  The complete
product--iterate identity from Theorem~\ref{thm:smallassoc}
and the vacuum identity give the $L$-contraction
\[
 \langle w',\mathcal Y^3_{X,Y}(x,3)y\rangle.
\]
For $R$, first use binary naturality for $r_X^1$ and then
transport.  Since
$\overline{r_X^1}I^1_{X,V}(x\otimes\mathbf1)
 =e^{L_X(-1)}x$, its contracted original product is the
intermediate-weight series
\[
 \sum_{n\ge0}\frac1{n!}
 \langle w',\mathcal Y^3_{X,Y}(L_X(-1)^nx,2)y\rangle.
\]
This expression does \emph{not} denote evaluation of a
$P(2)$-intertwining map on a general vector of $\overline X$.
To see precisely that it is legitimate, put
$x_{\le M}=\sum_{n=0}^M L_X(-1)^nx/n!\in X$.
Every $x_{\le M}$ is algebraic, and the weight-preserving
isomorphism $r_X^1$ identifies the projections of
$I^1_{X,V}(x\otimes\mathbf1)$ with those of
$\sum_n(r_X^1)^{-1}L_X(-1)^nx/n!$ in the completed
intermediate space.  For each $M$, the outer evaluation on
$x_{\le M}$ is defined in the ordinary sense; the original
product convergence identifies its restricted-dual scalar limit
with the displayed intermediate-weight series.  More directly,
for $\lambda=\delta-\alpha-\beta$, Lemma~\ref{lem:IOlogs}
gives a scalar $C$ such that
\[
 \langle w',\mathcal Y^3_{X,Y}(x,z)y\rangle=Cz^\lambda,
 \qquad
 \langle w',\mathcal Y^3_{X,Y}(L(-1)^nx,2)y\rangle
 =C\,n!\binom\lambda n2^{\lambda-n}.
\]
The second identity follows by differentiating the same
scalar germ $n$ times, not by applying an infinite exponential
operator.  The binomial series converges absolutely for
$|(3-2)/2|=1/2<1$, on the prescribed positive branch.  Therefore
both original scalar contractions equal $C3^\lambda$.
This also checks their analyticity and their agreement at the
actual placement $(3,2,0)$, rather than only formally at $z=0$.

Let
\[
 \phi=(L-R)((r_X^1)^{-1}\bo_2\id_Y):
 X\bo_2Y\longrightarrow X\bo_3Y.
\]
Its domain and codomain are genuine ordinary modules.  Set
$F_\phi(z)=\langle w',\overline\phi\,
\mathcal Y^{\mathrm{can}}_{X,Y}(x,z)y\rangle
=C_\phi z^\lambda$ for the corresponding complete formal
operator.  By the already established original contractions and
binary naturality, its \emph{finite truncated} Taylor
coefficients satisfy
\[
 0=\lim_{M\to\infty}
 \sum_{n=0}^M\frac1{n!}
 \langle w',\overline\phi\,
 \mathcal Y^{\mathrm{can}}_{X,Y}(L(-1)^nx,2)y\rangle
 =C_\phi\sum_{n\ge0}\binom\lambda n2^{\lambda-n}
 =C_\phi3^\lambda.
\]
The limit is justified exactly as above, for the $\phi$-composed
complete operator; composition with a module map preserves
weight projections and the $L(-1)$-derivative identity.  It follows
that $C_\phi=0$ for every homogeneous $x,y,w'$.  The restricted dual
separates the weight completion; thus the associated formal
operator is identically zero.  Faithfulness of the $P(2)$
representing map \eqref{eq:Pzuniversal} gives $\phi=0$,
and invertibility of $r_X^1$ gives \eqref{eq:rawtriangle82}.
In particular, no value at a single positive radius has been used
to infer vanishing of an arbitrary infinite collection of
unrelated coefficients: for each \emph{fixed homogeneous} test
there is just the single homogeneous scalar germ $C_\phi z^\lambda$.

The entire source of the left side of \eqref{eq:triangle82}
is $(X\bo_1V)\bo_1Y$.  Apply \eqref{eq:positivenormalization}
with the ordered triple $(X,V,Y)$; use transport naturality,
\eqref{eq:unittransport82}, and then the raw identity:
\begin{align*}
 &(\id_X\bo_1l_Y^1)a_{X,V,Y}\\
 &=(\id_X\bo_1l_Y^1)\mathsf T_{3,1}^{X,V\bo_1Y}
 (\id_X\bo_3\mathsf T_{2,1}^{V,Y})
 \mathcal A_{3,2}^{-1}\mathsf T_{1,2}^{X\bo_1V,Y}\\
 &=\mathsf T_{3,1}^{X,Y}(\id_X\bo_3l_Y^2)
 \mathcal A_{3,2}^{-1}\mathsf T_{1,2}^{X\bo_1V,Y}\\
 &=\mathsf T_{3,1}^{X,Y}\mathsf T_{2,3}^{X,Y}
 (r_X^1\bo_2\id_Y)\mathsf T_{1,2}^{X\bo_1V,Y}\\
 &=\mathsf T_{3,1}^{X,Y}\mathsf T_{2,3}^{X,Y}
 \mathsf T_{1,2}^{X,Y}(r_X^1\bo_1\id_Y)
 =r_X\bo\id_Y.
\end{align*}
The third equality is binary functoriality of the transport,
not an identification of differently parameterized maps.
For the vacuum module itself, ordinary VOA skew symmetry is
$Y_V(v,z)w=e^{zL(-1)}Y_V(w,-z)v$ on its fixed branch.
Hence $\Gamma_V(v,1)w=Y_V(v,1)w$ for all $v,w\in V$.
The two universal $P(1)$ maps \eqref{eq:Pzunits} therefore
represent the same complete operator; uniqueness gives
$l_V=r_V$, independently of the triangle.\end{proof}

\begin{proposition}\label{prop:action}
The pullback functors $T_k$ with the tensorators
\eqref{eq:J} and unit maps $j_k$ form a strong monoidal
$G$-action: for $k,l\in G$ they satisfy
\eqref{eq:Jcomp} and, for $X,Y,Z\in\cC_G$,
\begin{equation}\label{eq:actionassoc82}
 \begin{aligned}
 &T_k(a_{X,Y,Z})J^k_{X\bo Y,Z}(J^k_{X,Y}\bo\id_{T_kZ})\\
 &\qquad=J^k_{X,Y\bo Z}(\id_{T_kX}\bo J^k_{Y,Z})
 a_{T_kX,T_kY,T_kZ}.
 \end{aligned}
\end{equation}
The left and right unit diagrams commute and the crossed braiding satisfies \eqref{eq:equiv}.
\end{proposition}
\begin{proof}
For $X\in\cC_g$, the strict pullback is $T_kX\in\cC_{kgk^{-1}}$
with $Y_{T_kX}(u,z)=Y_X(k^{-1}u,z)$.  Direct substitution gives
\[
 Y_{T_kT_lX}(u,z)=Y_X(l^{-1}k^{-1}u,z)
 =Y_{T_{kl}X}(u,z),\quad T_1=\id.
\]
For a complete intertwining operator $\mathcal Y$ of type
$\binom Z{X\ Y}$ the same formal series, now of type
$\binom{T_kZ}{T_kX\ T_kY}$, obeys the four original Du--Huang
conditions: lower truncation and the $L(-1)$ identity are unchanged
because $k\omega=\omega$; one- and multiple-$V$ analytic identities
are those for $\mathcal Y$ with each inserted $v$ replaced by
$k^{-1}v$; the original branches, grades and projections are
unchanged.  This gives $J^{k,t}_{X,Y}$ at every fixed
$P(t)$ by its defining universal equation
\[
 \overline{J^{k,t}_{X,Y}}I^t_{T_kX,T_kY}
 =I^t_{X,Y}:T_kX\otimes T_kY\to\overline{T_k(X\bo_tY)}.
\]
(The equality on the right uses the unchanged underlying vector
spaces, not an identification of two differently twisted module
structures.)  Pulling back by $k^{-1}$ constructs its inverse.
The same construction at $t=1$ is \eqref{eq:J}.

Since $(kl)^{-1}=l^{-1}k^{-1}$, both
$J^{kl}_{X,Y}$ and $T_k(J^l_{X,Y})J^k_{T_lX,T_lY}$
are maps
$T_{kl}X\bo T_{kl}Y\to T_{kl}(X\bo Y)$
whose composition with the source canonical map is exactly the
pulled-back $I^1_{X,Y}$.  Its universal uniqueness gives
\[
 J^{kl}_{X,Y}=T_k(J^l_{X,Y})J^k_{T_lX,T_lY},
 \qquad J^1_{X,Y}=\id.
\]
For the unit arrow $j_k:V\to T_kV$, $j_k(v)=k^{-1}v$,
\[
 Y_{T_kV}(u,z)j_k(v)=k^{-1}Y_V(u,z)v
 =j_k(Y_V(u,z)v),
 \qquad j_{kl}=T_k(j_l)j_k,\quad j_1=\id.
\]
The last equality acts on $v$ as $l^{-1}k^{-1}v$;
its order cannot be interchanged when $G$ is nonabelian.
These are precisely the group-level tensorator and unit coherences,
proved without reference to any monoidal associativity coherence.

For $s,t>0$, the two well-typed maps
\[
 J^{k,t}_{X,Y}\mathsf T_{s,t}^{T_kX,T_kY},\quad
 T_k(\mathsf T_{s,t}^{X,Y})J^{k,s}_{X,Y}:
 T_kX\bo_sT_kY\to T_k(X\bo_tY)
\]
represent the same pulled-back evaluation at $s$ of the formal
operator belonging to the $P(t)$ representative.  More explicitly,
compose each with $I^s_{T_kX,T_kY}$ and use the defining
transport equation \eqref{eq:Pztransport}; both results are
$\operatorname{ev}_{\log s}(\mathcal Y^t_{X,Y})$
on the pulled-back modules.  Uniqueness in
\eqref{eq:Pzuniversal} gives
\begin{equation}\label{eq:actiontransport82}
 J^{k,t}_{X,Y}\mathsf T_{s,t}^{T_kX,T_kY}
 =T_k(\mathsf T_{s,t}^{X,Y})J^{k,s}_{X,Y}.
\end{equation}
The same proof applies to the specific upper-half-plane path
$-1\to1$ if all evaluation logarithms are \emph{continued along
that path}; this does not assert path independence under winding.

At the actual configuration $(3,2,0)$ consider the following two
morphisms, both with source
$((T_kX\bo_1T_kY)\bo_2T_kZ)$ and target
$T_k(X\bo_3(Y\bo_2Z))$:
\[
\begin{gathered}
 (T_kX\bo_1T_kY)\bo_2T_kZ
 \xrightarrow{J^{k,1}_{X,Y}\bo_2\id}
 T_k(X\bo_1Y)\bo_2T_kZ
 \xrightarrow{J^{k,2}_{X\bo_1Y,Z}}
 T_k((X\bo_1Y)\bo_2Z)\\
 \xrightarrow{T_k(\mathcal A_{3,2}^{-1})}
 T_k(X\bo_3(Y\bo_2Z)),\\[2pt]
 (T_kX\bo_1T_kY)\bo_2T_kZ
 \xrightarrow{\mathcal A_{3,2}^{-1}|_{T_kX,T_kY,T_kZ}}
 T_kX\bo_3(T_kY\bo_2T_kZ)\\
 \xrightarrow{\id\bo_3J^{k,2}_{Y,Z}}
 T_kX\bo_3T_k(Y\bo_2Z)
 \xrightarrow{J^{k,3}_{X,Y\bo_2Z}}
 T_k(X\bo_3(Y\bo_2Z)).
\end{gathered}
\]
To compare them, precompose with the actual nested canonical
$P(1),P(2)$-maps on \emph{algebraic} inputs and pair with an
arbitrary homogeneous restricted-dual output vector.  The first
chain becomes the $k$-pullback of the original canonical iterate
contraction; the second becomes the original canonical product
contraction for $T_kX,T_kY,T_kZ$ followed by the raw comparison.
Theorem~\ref{thm:smallassoc} equates these original normally
convergent contractions, and their pullbacks are obtained by
replacing every additional inserted state $v$ by $k^{-1}v$;
$L(0)$-projections, based branches and algebraic weight sums
are unchanged.  The scalar germ equality holds for every input
and every restricted-dual test, and
Lemma~\ref{lem:generation} separates the maps by the canonical
coefficient vectors.  Thus
\begin{equation}\label{eq:actionrawassoc82}
\begin{aligned}
 &T_k(\mathcal A_{3,2}^{-1})J^{k,2}_{X\bo_1Y,Z}
 (J^{k,1}_{X,Y}\bo_2\id_{T_kZ})\\
 &\quad=J^{k,3}_{X,Y\bo_2Z}
 (\id_{T_kX}\bo_3J^{k,2}_{Y,Z})
 \mathcal A_{3,2}^{-1}|_{T_kX,T_kY,T_kZ}.
\end{aligned}
\end{equation}
Insert \eqref{eq:actionrawassoc82} and
\eqref{eq:actiontransport82} into each of the four \emph{typed}
positive transport arrows of \eqref{eq:positivenormalization}.
Every transported outer/inner tensorator has matching $P(t)$
parameters, so the resulting two composites are exactly
\[
 T_k(a_{X,Y,Z})J^k_{X\bo Y,Z}(J^k_{X,Y}\bo\id_{T_kZ})
 \quad\text{and}\quad
 J^k_{X,Y\bo Z}(\id_{T_kX}\bo J^k_{Y,Z})a_{T_kX,T_kY,T_kZ}.
\]
They have source $(T_kX\bo T_kY)\bo T_kZ$ and target
$T_k(X\bo(Y\bo Z))$, proving \eqref{eq:actionassoc82}.
This uses the \emph{already proved} analytic associator comparison,
not a pentagon or an axiom imposed on $J^k$.

Both sides of the first equation below have source $V\bo T_kX$
and target $T_kX$, and both sides of the second have source
$T_kX\bo V$ and target $T_kX$:
\begin{equation}\label{eq:actionunit82}
\begin{aligned}
 T_k(l_X)J^k_{V,X}(j_k\bo\id_{T_kX})&=l_{T_kX},\\
 T_k(r_X)J^k_{X,V}(\id_{T_kX}\bo j_k)&=r_{T_kX}.
\end{aligned}
\end{equation}
Precompose with the corresponding canonical maps.  The left unit
then yields $Y_X(k^{-1}v,1)x$, equal to
$Y_{T_kX}(v,1)x$.  The right unit gives
$\Gamma_X(x,1)(k^{-1}v)=\Gamma_{T_kX}(x,1)v$
by the same $\Omega_-$ skew formula and
$kL(-1)=L(-1)k$.
The source $j_k$ already contributes the necessary $k^{-1}$;
there is no second pullback.  Binary universality proves both
equalities.  In particular, the maps $j_k$ are not silently
replaced by the identity $T_kV=V$.

For $X\in\cC_g$, $Y\in\cC_h$, both sides of the raw identity below
have source $T_kX\bo_1T_kY$ and target
$T_k(T_gY\bo_{P(-1)}X)$, of sector $kghk^{-1}$:
\begin{equation}\label{eq:actionrawbraid82}
 T_k(\mathcal R^1_{X,Y})J^{k,1}_{X,Y}
 =J^{k,-1}_{T_gY,X}\mathcal R^1_{T_kX,T_kY}.
\end{equation}
Indeed the raw skew operator defining \eqref{eq:Pzrawcrossed}
has the fixed upper-half-plane $\log(-1)=\pi i$ and the
actual translation $e^{L(-1)}$.  On pulling back, the state
insertions $v$ become $k^{-1}v$, and the twisting sector of
the reordered input is
$T_{kgk^{-1}}(T_kY)=T_k(T_gY)$, since
$T_kT_g=T_{kgk^{-1}}T_k$.  Specifically, after precomposition with
$I^1_{T_kX,T_kY}(x\otimes y)$ the left side is the
underlying-vector-space evaluation of
$e^{L(-1)}I^{-1}_{T_gY,X}(y\otimes x)$ in
$T_k(T_gY\bo_{P(-1)}X)$.  On the right, the raw commutor
first gives
$e^{L(-1)}I^{-1}_{T_k(T_gY),T_kX}(y\otimes x)$;
$J^{k,-1}_{T_gY,X}$ sends this to exactly the same
pulled-back evaluation.  These are the full evaluations of
$\Omega_+$ with the original $+\pi i$ logarithm, including
its fractional-power phases; there is no \emph{extra} phase
from pullback.  Uniqueness of the $P(1)$ representative proves
the raw equality.
Continue \eqref{eq:actiontransport82} along the same
upper-half-plane path from $-1$ to $1$ with its inherited
logarithms.  Compose it with \eqref{eq:actionrawbraid82} and
\eqref{eq:c} to obtain
\[
 T_k(c_{X,Y})J^k_{X,Y}
 =J^k_{T_gY,X}c_{T_kX,T_kY}.
\]
The two maps have source $T_kX\bo T_kY$ and target
$T_k(T_gY\bo X)$.  The assertion compares one and the same
lifted half-turn and makes no claim about arbitrary winding
paths.  Naturalities of $J^k,j_k$ and $c$ follow from the
same canonical evaluations and binary universality.  No
crossed hexagon, equivariantization or rigidity has been used.
\end{proof}

\section{Coherence}\label{sec:coherence}\label{sec:pentagon}
\subsection{Four-input convergence and the pentagon}
We establish a common convergence domain for the five four-input tree channels.

\begin{lemma}
\label{lem:criticaloneedge}
Let $L_0=\C((s))$, and let $M$ be an ordinary (in particular, grading-restricted and $L(0)$-semisimple) $g$-twisted $V$-module of finite order $N$, with restricted contragredient $M'$.  Let $\lambda_{-}$ be a three-point Ward
functional with an $M'$-input, and let $\lambda_{+}$ be an actual outside-component Ward functional with an $M$-input.  Assume
\begin{equation}\label{eq:oneedgewardhyp}
 \lambda_{\pm}\bigl(R_{\pm}(\sigma_{\pm})w_{\pm}\bigr)=0
 \quad\bigl(\sigma_{\pm}\in\mathscr J_{\pm},\
 w_{\pm}\text{ algebraic}\bigr),
\end{equation}
where $\mathscr J_{\pm}$ denotes the invariant meromorphic current space on the corresponding actual component, with arbitrary finite node-pole order.  The values of $\lambda_{+}$ may belong to
\begin{equation}\label{eq:oneedgemonomial}
 \mathscr R=\bigoplus_{\alpha\in\mathcal A}s^\alpha L_1,
 \qquad |\mathcal A|<\infty,\quad
 L_1=\C((s^{1/f})),\quad f\ge1,
\end{equation}
where $s^\alpha$ is a formal monomial and the summands are $L_1$-lines.  Fix compatible framed branches, local root coordinates, and a balanced degeneration of the original based cover in
which the two-leaf component is attached to the outside component along $M\otimes M'$.  Then the relative-grade contraction
\[
\Phi=\left(\lambda_{+}\otimes\lambda_{-}\right)
 \left(\Gamma_M(\rho)\right),\qquad
 \Gamma_M(\rho)=\sum_{\gamma\in\Lambda_M}
 \rho^\gamma\sum_{i=1}^{\dim M_{h_M+\gamma}}
 m_{\gamma,i}\otimes m_{\gamma,i}^{\vee},
 \quad\Lambda_M\subseteq N^{-1}\Z_{\ge0},
\]
with the physical frame and lowest-weight factors restored as in Lemma~\ref{lem:physicalcontraction}, is defined coefficientwise in the attaching parameter.  For every invariant generic current
$\sigma$ on the original based cover and every algebraic external state $w$, the following holds for a suitable integer $d\ge1$:
\begin{equation}\label{eq:oneedgefullward}
\begin{gathered}
 \Phi\bigl(R_{\mathrm{gen}}(\sigma)w\bigr)=0
 \quad\text{in }\displaystyle\bigoplus_{\beta\in\mathcal B}
 t^\beta\mathscr R((t^{1/d})),\\
 |\mathcal B|<\infty,\qquad
 \rho=t^{aN}\quad\text{after finite base change}.
\end{gathered}
\end{equation}
The identity descends to the original framed coefficient field and does not require convergence of the remaining $s$-series.
\end{lemma}
\begin{proof}
The coefficient module is used with its explicit support, not as a
completion of the node tensor product.  Use the ordinary grading
$M=\bigoplus_{p\in e\Lambda_M}M_{h_M+p/e}$ from the statement;
if a finite semisimple sum has several lowest-weight cosets, apply
the same argument separately to its finitely many summands.
Each $M_{h_M+p/e}$ and its restricted dual are finite dimensional,
so $P_p=\sum_i m_{p,i}\otimes m_{p,i}^{\vee}$ is an
\emph{algebraic} tensor.  The symbol $\Gamma_M(t^e)=\sum_p t^pP_p$
is never evaluated as a single vector of $M\otimes M'$: by definition
its contraction is the collection of $\mathscr R$-valued coefficients
obtained with $p$ fixed.  The only infinite sum remaining at this
stage is in the spectator variable, inside the Laurent field $L_1$.
Choose the monomial shifts in \eqref{eq:oneedgemonomial} modulo
$f^{-1}\Z$ without repetition; thus different summands are formally
independent, including after the common finite coefficient-field
extension used below.

Put $K_0=L_0((\rho))$.  By
\cite[Corollary~3.0.5 and Theorem~4.3.2]{ACV}, the balanced
$H$-cover over $K_0$ extends, after a finite extension of the $\rho$-adic DVR, to an admissible $H$-cover over its stable four-input boundary.  The projective nodal source and its inertia
charts follow from \cite[Lemma~2.2.1 and Theorem~4.3.2]{ACV}. Preserve the original marked sections, the auxiliary framing lift, and the remaining ratio $s$ during this extension.  Lemma~\ref{lem:pureramification}
then gives, after enlarging the coefficient field $L_1\subset L$,
\[
\begin{gathered}
 O=L[[t]],\qquad K=L((t)),\qquad U_+U_-=t^a,\\
 x_+=U_+^N,\qquad x_-=U_-^N,\qquad
 \rho=x_+x_-=t^{aN},\quad a,N\ge1.
\end{gathered}
\]
Here $U_+,U_-$ are \emph{node} root coordinates, not the $L^\times$-unit introduced below.  Neither an additional outside marking nor a remaining-ratio coordinate varies with $t$.

We justify the use of Lemma~\ref{lem:canonicalnodecurrent} in this five-point degeneration.  Let $\mathscr F_E$ be its coherent finite-state current sheaf, enlarged by the prescribed polar
divisor at the extra outside marking, and let $\sigma$ be one fixed generic invariant current.  To justify extension of an actual generic section, let
$p_*\mathscr F_E(D_{\mathrm{ext}})$ denote its coherent direct-image
module over the complete DVR $O$ (the source is proper and the marked
polar divisor is Cartier).  Its generic fibre contains the given
section $\sigma$ by coherent base change over the fraction field:
a fixed element of this localization is a finite sum of fractions with
powers of $t$ as denominators.  Hence a power $t^{m(\sigma)}$
clears its denominator.  Averaging over $H$ is legitimate in
characteristic zero and commutes with this localization.  Thus there
is an integer $m=m(\sigma)\ge0$ such that
\begin{equation}\label{eq:oneedgeextension}
 t^m\sigma\in
 H^0\bigl(\mathcal C,\mathscr F_E(D_{\mathrm{ext}})\bigr)^H
 \quad\text{over }O.
\end{equation}
The argument is the same coherent-lattice localization as in the proof of Lemma~\ref{lem:canonicalnodecurrent}; it is
\emph{not} an assertion that $m$ is uniform over all currents.
The extra section is disjoint from the node.  Consequently the comparison maps $\iota_{v,J}$ on punctured components and the finite Kummer charts at the old markings used in that proof persist
over $O/(t^{J+1})$.  The invariance of finite \'etale covers under nilpotent thickenings \cite[Tag~0BQB]{Stacks} and the fixed root frames
preserve these comparisons for every $J$.  At each source node the state transition is still $\Sigma=e^{L(1)}(-1)^{L(0)+1}$, with all $L(1)$-terms.  The normalized residue is the selected-lift
residue of \eqref{eq:normalizednode62}.  The statement about
``coefficients on the normalization'' must be interpreted after the
actual rooted comparison, not by taking coefficients in arbitrary
$t$-dependent local trivializations.  More precisely, fix $J$.  Delete
the node and the external marked divisors from each normalization and
use the fixed limiting marked coordinates there.  The nilpotent
finite-etale lifting comparison, pinned down by the chosen based lift
and compatible selected roots, identifies the actual pulled-back cover
on this punctured open modulo $t^{J+1}$ with the constant normalized
cover over $O/(t^{J+1})$.  Apply the full finite-state coordinate
transition to the current sheaf under this identification.  The
result is a single globally defined Laurent polynomial in $t$ of
degree at most $J$ on the fixed punctured normalization.
Its coefficient of $t^j$ patches across overlaps, hence is an actual
meromorphic component current; poles occur only at the prescribed
marks and node, and have finite order for each fixed $j$.  At the
node it is evaluated in the inherited balanced root frame, whose
finite $L(1)$ corrections have already been included.  These
comparisons are compatible when $J$ increases.  Without this
transport a coefficient of a global family section need not itself
be a global component section, since its untransported transition
functions may depend on $t$.

For precision, if $D_j$ bounds the poles of the coefficient of $t^j$ and $F$ is a field of definition for the outside framed component and $\mathscr F_E(D_j)$, then flat field extension gives
\begin{equation}\label{eq:oneedgebasechange}
 H^0\bigl(C_F,\mathscr F_E(D_j)\bigr)\otimes_F L
 \xrightarrow{\ \sim\ }
 H^0\bigl(C_L,\mathscr F_E(D_j)_L\bigr).
\end{equation}
Thus each coefficient is a finite $L$-linear combination of actual outside currents, and \eqref{eq:oneedgewardhyp} applies $L$-linearly, also in $\mathscr R\otimes_{L_1}L$.  This uses the
hypothesis for \emph{every algebraic $M$-state}; it does not apply any Ward identity to an uncompleted infinite vector of $M$.

We check the interchange of node residues and contraction.  Fix $\sigma$, one homogeneous $b\in V_d$ occurring in its finite-state part, and a monomial multiplying the standard node current $b\,x_+^{d-1}dx_+$ in its balanced expansion
\[
t^{-m+\ell}U_+^rU_-^s b,
 \qquad \ell,r,s\ge0,\quad k=r-s,\quad
 n=d-1+\frac{k}{N},\quad e=aN.
\]
Set $p=e\gamma\in\Z_{\ge0}$ for a homogeneous summand of $\Gamma_M(t^e)$.  The $+$-branch mode $b_n$ has degree $\delta=-ak$ in $t$-grading.  At coefficient $t^j$, Lemma~\ref{lem:nodewindow} yields
\begin{equation}\label{eq:oneedgewindow}
 \begin{gathered}
 p+as+\ell=j+m,\qquad
 a\max(0,k)\le p\le j+m-a\max(0,-k),\\
 |k|\le\left\lfloor\frac{j+m}{a}\right\rfloor,
 \qquad r=s+k\ge0.
 \end{gathered}
\end{equation}
Indeed, when $k\ge0$, the mode lowers relative degree by
$ak$; hence $p\ge ak$, and $s\ge0$ forces $p\le j+m$.
When $k<0$, integrality of the original node monomial requires
$s\ge-k$, whence $p\le j+m+ak$; there is no extra lower
bound beyond $p\ge0$.  These inequalities also imply
$ak\le j+m$ for $k\ge0$ and $-ak\le j+m$ for $k<0$.
For fixed $j$, $p,k$ run through finite sets and
$as+\ell=j+m-p$ leaves only finitely many nonnegative $(s,\ell)$;
$r=s+k$ is then determined.  The calculation holds for a completed
integral node section by inspecting its image in the finite quotient
needed for this one coefficient, rather than interchanging an
unrestricted formal sum with a residue.
For fixed $j$ this is a finite set of $(p,k,r,s,\ell)$, including the convention $M_{h_M+p/e}=0$ for absent grades. All operations below are therefore on finite-dimensional spaces, although the $s$-pole orders of their $\mathscr R$-valued
coefficients need not be uniformly bounded as $p$ varies.

Write $P_p\in M_{h_M+p/e}\otimes M_{h_M+p/e}'$ for the coevaluation tensor and $\Gamma=\sum_p t^pP_p$.  The two node residues of \eqref{eq:oneedgemonomial} act by
\begin{equation}\label{eq:oneedgecancel}
 \begin{aligned}
 (b_n\otimes\id)\Gamma
   &=t^{-\delta}(\id\otimes b_n^{\mathsf t})\Gamma,\\
 R_{+,\mathrm{node}}\Gamma
   &=(b_n\otimes\id)\Gamma,\\
 R_{-,\mathrm{node}}\Gamma
   &=-t^{-\delta}(\id\otimes b_n^{\mathsf t})\Gamma,\\
 (R_{+,\mathrm{node}}+R_{-,\mathrm{node}})\Gamma&=0.
 \end{aligned}
\end{equation}
Here $b_n^{\mathsf t}$ is the algebraic transpose on $M'$; the full opposite-mode formula, including the $\sum_{j=0}^dL(1)^jb/j!$ correction, is
\eqref{eq:node63opposite}, as in
\cite[Section~3]{DLM} and \cite[Section~2]{DH}.
For the direction of transposition, if
$b_n:M_p\to M_{p+\delta}$, its algebraic transpose is
$b_n^{\mathsf t}:M_{p+\delta}^*\to M_p^*$.
For any basis $e_i$ of $M_p$ and any basis $f_a$ of
$M_{p+\delta}$, the elementary identity
$\sum_i b_ne_i\otimes e^i=\sum_a f_a\otimes b_n^{\mathsf t}f^a$
gives $(b_n\otimes1)P_p=(1\otimes b_n^{\mathsf t})P_{p+\delta}$.
Consequently the coefficient on the right of the corresponding
identity for $\Gamma$ is shifted from $t^{p+\delta}$ to
$t^p$ precisely by $t^{-\delta}=t^{ak}$, not its inverse.
The actual opposite residue is not an untyped transpose: the
balanced opposite-node field, with its full finite
$e^{L(1)}(-1)^{L(0)+1}$ state transformation, is what realizes
$-b_n^{\mathsf t}$ after the displayed power of $t$.
Both sides use the same selected-root $1/N$ orbit normalization;
under a deck change the primal and restricted-dual arrows are
inverse transposes, so this elementary identity is unchanged.
Equation~\eqref{eq:oneedgecancel} is read coefficientwise using
\eqref{eq:oneedgewindow}; the two residue signs and the factor
$1/N$ are those of \eqref{eq:normalizednode62} and
\eqref{eq:canonicalpuremode}.  On induced source components the
simultaneous $H$-pullback and the orbit-normalization identity give the same cancellation.  Apply the two component Ward identities \eqref{eq:oneedgewardhyp} and sum: the node terms cancel by \eqref{eq:oneedgecancel}, leaving
\begin{equation}\label{eq:oneedgerelativeward}
 \bigl(\lambda_+\otimes\lambda_-\bigr)
 \bigl(R_{\mathrm{ext}}(\sigma)w\otimes\Gamma_M(t^e)\bigr)=0
 \quad\text{coefficientwise in }t.
\end{equation}

It remains to identify this with the physical contraction. Write the physical cherry displacement as $r_\mathrm{ch}=U\rho$, with $U\in L^\times$ independent of $t$; fix its logarithm and put
$(U\rho)^\zeta=U^\zeta\rho^\zeta$ on that branch. For homogeneous $a,b$ of weights $\lambda_A,\lambda_B$, and $v_\mu=P_\mu\mathcal Y(a,1)b\in M_\mu$, define the outside node-frame change by
\[
\Psi_{\mathrm{out}}^{x}(m)
 =\Psi_{\mathrm{out}}^{\mathrm{phys}}
 \bigl(U^{L_M(0)}m\bigr),\qquad
 x=(Z-\mathrm{anchor})/U,
\]
with all other marked local frames held fixed.  Then, weight by weight,
\begin{equation}\label{eq:oneedgephysical}
 \begin{aligned}
 &(U\rho)^{\mu-\lambda_A-\lambda_B}
 \Psi_{\mathrm{out}}^{\mathrm{phys}}(v_\mu)\\
 &\quad=\rho^{h_M}\Psi_{\mathrm{out}}^{x}
 \Bigl(\rho^{\mu-h_M}P_\mu
 \mathcal Y\bigl((U\rho)^{-L_A(0)}a,1\bigr)
 (U\rho)^{-L_B(0)}b\Bigr).
 \end{aligned}
\end{equation}
Explicitly, $P_\mu\mathcal Y((U\rho)^{-L_A(0)}a,1)
(U\rho)^{-L_B(0)}b=(U\rho)^{-\lambda_A-\lambda_B}v_\mu$,
whereas $\Psi_{\mathrm{out}}^x(v_\mu)
=U^\mu\Psi_{\mathrm{out}}^{\mathrm{phys}}(v_\mu)$.
The right side of \eqref{eq:oneedgephysical} is thus
$\rho^{\mu}U^{\mu-\lambda_A-\lambda_B}
\rho^{-\lambda_A-\lambda_B}
\Psi_{\mathrm{out}}^{\mathrm{phys}}(v_\mu)$, exactly its left side.
The factor $\rho^{h_M}$ is indispensable: it restores the removed
lowest-weight monomial.  This computation is performed at fixed
$\mu$, before either the spectator series or the node sum is formed.
The identity is an equality of individual finite homogeneous contractions; it does not invoke convergence of the $\mu$-sum. The full finite state-coordinate gauge conjugates every external
residue operator in \eqref{eq:oneedgerelativeward} to its physical counterpart.  Equation~\eqref{eq:oneedgephysical} restores the weights and the retained lowest-weight monomial; the global root scale is restored by the $D_R$ factor of
Lemma~\ref{lem:physicalcontraction}.  Consequently
\eqref{eq:oneedgerelativeward} is precisely
\eqref{eq:oneedgefullward}.  Finally the coefficient-field and
DVR extensions are injective (also on each formal monomial line); the zero identity descends to the original framed coefficient field.
\end{proof}

\begin{definition}\label{def:E1cherrydata}
Write the five four-input trees as
\[
\begin{aligned}
 T_1&=(((AB)C)E),& T_2&=((A(BC))E),&
 T_3&=(A((BC)E)),\\
 T_4&=(A(B(CE))),& T_5&=((AB)(CE)).
\end{aligned}
\]
For each tree use the following normalized physical coordinates
$(z_A,z_B,z_C,z_E)$, with actual positions
$R(z_A,z_B,z_C,z_E)$, $R\ne0$ and the indicated logarithm:
\begin{equation}\label{eq:criticalfivepositions}
\begin{array}{c|c}
 T_1&(1+u+uv,\ 1+u,\ 1,\ 0)\\
 T_2&(1+u,\ 1+uv,\ 1,\ 0)\\
 T_3&(1,\ u(1+v),\ u,\ 0)\\
 T_4&(1,\ u,\ uv,\ 0)\\
 T_5&(1+u,\ 1,\ v,\ 0).
\end{array}
\end{equation}
For its two-leaf cherry let $M_i$ be its ordinary output,
$\rho_i$ the collision ratio, $s_i$ the remaining ratio and
$r_i^{\rm ch}=U_i\rho_i$ its physical displacement.  The node inertia
and outside three-input channel are, respectively,
\[
\begin{array}{c|c|c|c|c|c}
 &\text{cherry}&g_{M_i}&(\rho_i,s_i)&U_i&\text{outside}\\ \hline
 T_1&AB&g_Ag_B&(v,u)&Ru&I\\
 T_2&BC&g_Bg_C&(v,u)&Ru&I\\
 T_3&BC&g_Bg_C&(v,u)&Ru&P\\
 T_4&CE&g_Cg_E&(v,u)&Ru&P\\
 T_5&AB&g_Ag_B&(u,v)&R&P
\end{array}
\]
Here $P$ and $I$ indicate the original product and iterate E0
channels.  These are coordinate definitions, not assertions that the
five channel functionals coincide.
\end{definition}

\begin{lemma}\label{lem:E1jointsourcecomparison}
Use Definition~\ref{def:fiveoutsetframing}.  For each of the five trees
\(T_i\) of Definition~\ref{def:E1cherrydata}, carry its specified
positive-chamber five-point cover, including its based lift, complete
lower-comb paths and external roots, to the actual cherry boundary.
The two normalization components are given the inherited paired node lift
and opposite node roots; choose their E0 and binary reference
\emph{expressions} by transporting these induced frames to the positive
Du--Huang preferred \((1,0,\infty)\) coordinates.  The three-point
operator at every vertex remains the original fixed operator at
\(z=1,\log1=0\).  In this one jointly rooted physical frame the
coefficient of each cherry grade \(\mu\) of the original three-operator
channel is precisely
\begin{equation}\label{eq:E1originalsourcecomparison}
 (U_i\rho_i)^{\mu-\alpha_i}
 \Psi_i\bigl(\pi_\mu^{M_i}\mathcal Y_i(p_i,1)q_i;s_i\bigr),
\end{equation}
with the same vectors and intertwining-operator labels as in
\eqref{eq:criticalfivefinitecontraction}.  The outside functional has
full-state \emph{formal} Ward over its spectator Laurent coefficient field,
and the binary functional has three-point Ward.  If the cherry parameter
is written $\rho_i=t^e$, then $e$ is chosen divisible by the node inertia
order, the two cherry marking root orders, and the finite ramification
orders of the actual rooted Hurwitz comparison; a further finite extension
of the spectator field accommodates its root orders.  All logarithms
remain those continued from the positive outset.  This claim concerns a
specified jointly induced rooted frame, not literal equality of arbitrary
independently preselected three-point frames.

For comparison with two independent reference frames let
\(A_\mu:M_\mu^{\rm act}\to M_\mu^{0}\) and
\(B_\mu:(M_\mu^{\rm act})'\to(M_\mu^0)'\) be the typed residual root/deck comparison arrows after the
explicit physical scale in \eqref{eq:oneedgephysical} has been extracted.
They preserve each individual homogeneous node grade; no nonlinear
jet change is concealed in these grade-wise arrows.  Put \(K_\mu=B_\mu A_\mu^t\).  If
\(\lambda_O^0,\lambda_B^0\) are the fixed original three-point
coefficients and \(G_O,G_B\) the external residual comparison arrows, the same
physical finite contraction in those references is
\begin{equation}\label{eq:E1referencecompensation}
 (U_i\rho_i)^{\mu-\alpha_i}
 \sum_j\lambda_O^0(f_j;G_Ow_O)
               \lambda_B^0(K_\mu f^j;G_Bw_B).
\end{equation}
Here \(f_j,f^j\) are a reference dual basis.  When the reference
\emph{pair} is induced from the same actual node, \(B_\mu=(A_\mu^{-1})^t\)
and \(K_\mu=1\).  With independently chosen reference frames, the compensation in
\eqref{eq:E1referencecompensation} is retained; the assertion does not
identify the two untransported scalar expressions.
\end{lemma}
\begin{proof}
The positive chamber of four ordered finite positions is convex.
Transporting the one based lift and the full lower-comb path system along
its positive paths is therefore independent of the chosen positive
homotopy.  For a cherry of consecutive leaves \(X,Y\), the lower boundary
stem represents the ordered word \(\gamma_X\gamma_Y\); the bubble
has meridians \((g_X,g_Y,(g_Xg_Y)^{-1})\), and the outside has the
collapsed meridian \(g_Xg_Y\) with the other original meridians in order.
Apply the actual stable reduction to this specified generic cover,
not to a separately reconstructed cover with only the same meridian tuple.
By the balanced local model of Lemma~\ref{lem:criticaloneedge}, the
normalizations inherit the same lifted boundary stem and opposite node
root parameters.  Every original external marking retains its already
selected lift and root.  The unramified outside based lift is obtained by
lifting the right-hand segment from the transported five-point basepoint
to the outside basepoint; it is not independently set equal to an E0 lift.
For a potentially disconnected normalization component, use the induced
\(H\)-cover from its actual monodromy-image subgroup, as in the
convention preceding \eqref{eq:localcurrent25}.

Write \(r=U_i\rho_i>0\) for the physical cherry difference and
\(\xi=(Z-Z_Y)/r\).  Since this affine change has positive slope, the
restriction of the original lower/direct paths is homotopic, with its
specified based stem, to the preferred three-point paths for
\((\xi_X,\xi_Y)=(1,0)\).  The original marking roots satisfy
\(r_X^{\rm phys}=r^{1/N_X}r_X^{\rm bin}\) and
\(r_Y^{\rm phys}=r^{1/N_Y}r_Y^{\rm bin}\) after transporting \emph{those
same} roots; the complex powers use the continued positive logarithm.
The comparison here uses the entire based path: along the
lower corridor the chosen cluster stem followed by either terminal leaf
stem is homotopic, relative to its end disk, to the original lower-comb
reference path.  Uniqueness of path lifting identifies the marking lifts;
the normalized local root is the particular root
\(r_X^{\rm phys}/r^{1/N_X}\), not an arbitrarily chosen root of the
same coarse coordinate.  On positive real changes of scale its root
phase is constant.  The initial calibration in
Definition~\ref{def:fiveoutsetframing} therefore identifies this
rooted/logged lower germ with the Du--Huang preferred one: at
\(\xi_X-\xi_Y=1\) its continued logarithm is exactly zero.  The
uniqueness assertion is based path lifting plus the prescribed
special-root torsor element, rather than a classification by monodromy
tuples alone.  Thus the bubble germ of the specified operator is its original
\(\mathcal Y_i(p_i,1)q_i\), not a new conformal block obtained merely
from the inertia tuple.  The output node root is the balanced opposite
of the outside node root, so their module maps are typed contragredient
inverses.  In particular the two actual node maps preserve the
restricted coevaluation, on each homogeneous grade including its
lowest-weight line.

To compare the outside with an auxiliary E0 outset, use the actual
right-hand based path and the given roots to form the typed
coordinate/deck arrow of Lemma~\ref{lem:E0formalframedtransport}.
It transports \emph{the original} E0 formal full-state Ward, including
all \(V_{\le D}\) cutoffs, into the actual outside component.  No
analytic continuation or five-point convergence is invoked.
For an outside physical local coordinate \(X=Z-Z_Y\) and its collar
coordinate \(x=X/U_i\), their node state frames are related by the
complete action of Lemma~\ref{lem:relativecoordinates}; restoring the
lowest-weight line gives \(U_i^{L_{M_i}(0)}\) on the node input.
By the original \(L(0)\)-commutator for \(\mathcal Y_i\), and not by
specifying a new functional, its homogeneous bubble coefficient carries
\(r^{\mu-\alpha_i}\).  The finite homogeneous identity
\eqref{eq:oneedgephysical} therefore identifies the resulting physical
contraction with \eqref{eq:E1originalsourcecomparison}.  This
identification is made before taking any infinite sum.

The two distinct scales can be seen explicitly. For
\(T_3=A((BC)E)\), take
\((z_A,z_B,z_C,z_E)=(1,u(1+v),u,0)\),
\(h=g_Bg_C\), \(r=uv\), \(U=u\),
\(\xi=(Z-u)/(uv)\), \(x=(Z-u)/u\).  The outside is the
original product E0 in the global \(Z\)-positions
\((1,u,0)\).  Its node insertion is transported by
\(u^{L_M(0)}\), and the original bubble coefficient is
\(m_\mu=\pi_\mu\mathcal Y_{BC}(b,1)c\).  On this induced frame
\((uv)^{\mu-\wt(b)-\wt(c)}F_P(m_\mu;u)\) is exactly the original
T3 coefficient.  In independently chosen reference frames, it is
\(F_P(u^{L_M(0)}J_Mm;J_{\rm ext}w;u)\) before restoring the
bubble and lowest-weight factors: the \(J\)'s are actual-to-reference
typed transport arrows, not numbers chosen to fit the weight exponent.
For \(T_5=((AB)(CE))\), take
\((z_A,z_B,z_C,z_E)=(1+u,1,v,0)\),
\(h=g_Ag_B\), \(r=u\), \(U=1\),
\(\xi=(Z-1)/u\), \(x=Z-1\); the outside is the original product
E0 at \((1,v,0)\).  The original coefficient is
\(u^{\mu-\wt(a)-\wt(b)}F_P(\pi_\mu\mathcal Y_{AB}(a,1)b;v)\),
with no additional spectator dilation.  Both identities are obtained
from the actual normalized coordinate changes and fixed operators.

The five relevant ordered contractions, all taken on the same
positive physical chamber, can therefore be checked without a
cross-tree factorization hypothesis: their cherry/outside data are
$(AB,I),(BC,I),(BC,P),(CE,P),(AB,P)$, respectively; the
ordered internal twist is $g_Ag_B$, $g_Bg_C$, $g_Bg_C$,
$g_Cg_E$, $g_Ag_B$.  In particular, the outside use of E0
Ward is on a genuine three-input source with the collapsed
meridian in its original position, not on a source identified solely
from a tuple of conjugacy classes.

For the other three trees the additional data are: \(T_1\) uses
\(AB\), \(h=g_Ag_B\), \(r=uv\),
\(\xi=(Z-1-u)/(uv)\), outside iterate \((1+u,1,0)\);
\(T_2\) uses \(BC\), \(h=g_Bg_C\), \(r=uv\),
\(\xi=(Z-1)/(uv)\), outside iterate \((1+u,1,0)\);
\(T_4\) uses \(CE\), \(h=g_Cg_E\), \(r=uv\),
\(\xi=Z/(uv)\), outside product \((1,u,0)\).
Their outside node collar scales are \(U=u\).  The actual
five-point right basepoints and the outside three-point right basepoints
are, in order \(T_1,\ldots,T_5\),
\[
 (5+2u+uv,4+u),\quad(5+u+uv,4+u),\quad
 (3+2u+uv,3+u),\quad(3+u+uv,3+u),\quad
 (4+u+v,3+v).
\]
In each case the specified comparison path is the right-hand real
segment between these two points, lifted from the original five-point
based lift.  For every cherry leaf \(X\), the comparison of its
original special root with the bubble root is exactly the
positive factor \(r^{1/N_X}\); for a non-cherry external leaf its
original centered physical root is inherited on the outside without
replacement.  Because the three root normalizations for the cherry
species \(AB,BC,CE\) all have positive scaling and the same initial
lower-comb marking path, this list gives a single root assignment to
each shared external marking, not five conflicting assignments.
The same \emph{transported} physical external roots are used wherever a
marking occurs in more than one tree.  In particular the fixed-zero
marking \(E\) in \(T_4\) is not re-rooted on its bubble.  The lower
stems for \(AB\) in \(T_1,T_5\) and for \(BC\) in \(T_2,T_3\)
are homotopic through the ordered positive chamber; this compares their
smooth rooted origins without asserting any equality of tree-channel
conformal blocks.

Only after the physical scale and the full opposite-node
coordinate change have been extracted do the residual based-root
arrows act on individual node grades.  On the positive affine
comparison charts these residual arrows consist of root rotations
and typed deck maps; they commute with $L(0)$ (and thus with each
homogeneous projection).  Any nonlinear finite-state coordinate
change is instead retained in the already specified external or
node-current gauges, not discarded by asserting grade preservation.

Finally choose a reference basis \(m_j\) with dual \(m^j\).  If
\(A_\mu\) takes an actual primal node basis to this reference, its
reference dual is \((A_\mu^{-1})^tm^j\); sending the actual dual by
\(B_\mu\) therefore inserts exactly \(B_\mu A_\mu^t=K_\mu\) in
one leg.  The remaining external and output factors act through the
specified typed \(G_O,G_B\), with inverse transpose at the output.
More explicitly, if $(A_\mu)_{ji}$ and $(B_\mu)_{li}$
are the matrices of the primal and dual changes in these bases, the
transported coevaluation is
$\sum_i(A_\mu m_i)\otimes(B_\mu m^i)
=\sum_j f_j\otimes(B_\mu A_\mu^t)f^j$.
This identifies the side on which the compensation acts and
explains why $K_\mu$ is generally nontrivial.
This proves \eqref{eq:E1referencecompensation}.  When the references
are inherited as one pair from the actual node, the maps are dual
inverse, proving \(K_\mu=1\).  Independent reference pairs are not
assumed to have this property.  The original coefficients and the
original \(\log1=0\) germs were never changed. 
\end{proof}

\begin{theorem}\label{thm:criticalfourward}
Let $A,B,C,E$ be ordinary finite-order twisted $V$-modules, $D$ an ordinary output module of the requisite ordered-product sector, and $a\in A$, $b\in B$, $c\in C$, $e\in E$, $d'\in D'$ homogeneous.  Write $T_i$ for the trees of Definition~\ref{def:E1cherrydata}.
Fix the single smooth rooted/logged five-input outset of Definition~\ref{def:fiveoutsetframing} and three actual typed twisted intertwining operators at the internal vertices of one of the five ordered binary trees $T_i$.  Every coefficient refers to the original $z=1$, $\log1=0$ intertwining operators in the physical external frames induced from that outset.  For its original two-edge homogeneous contraction $\Phi_{T_i}$,
on the anchored framed physical chart and for every canonical five-point current $\sigma$ and algebraic external state $w$,
\begin{equation}\label{eq:criticalfiveward}
 \Phi_{T_i}\bigl(R(\sigma)w\bigr)=0.
\end{equation}
The two internal grade sums are defined coefficientwise, with no convergence assumption.  For every derivation $\delta$ of the chosen framed base,
\begin{equation}\label{eq:criticalfivehorizontal}
 \partial_\delta\Phi_{T_i}(w)
   =\Phi_{T_i}(\nabla_\delta w),\qquad
 \Phi_{T_i}(\cJ_{\diff})=0.
\end{equation}
The equalities are in the formal generalized-power module determined by the chosen root coordinates and logarithms.
\end{theorem}
\begin{proof}
Put $g_A,g_B,g_C,g_E\in G$ for the input twists.  Fix simple ordinary intermediate modules and their types; finite sums of such choices are taken at the end.  Use the coordinates, cherry scales and outside channel types fixed in
Definition~\ref{def:E1cherrydata}; in particular the actual physical
positions are $R(z_A,z_B,z_C,z_E)$.  For the chosen cherry let $M_i$ be
its output, let $\alpha_i$ be the sum of its two input weights, and put
$r_i^{\rm ch}=U_i\rho_i$.
The symbols $P,I$ denote the original product and iterate channels.  The actual outside rooted source is identified, without changing its original IO coefficients, by Lemma~\ref{lem:E1jointsourcecomparison}; its coefficientwise full-state Ward is the \emph{formal} Theorem~\ref{thm:fullsmallward} transported by Lemma~\ref{lem:E0formalframedtransport}.  Analytic continuation is not used here.  After contracting the cherry, the three remaining physical input positions, divided by $R$, are $(1+u,1,0)$ for $T_1,T_2$,
$(1,u,0)$ for $T_3,T_4$, and $(1,v,0)$ for $T_5$. For the formal $\rho_i$-adic argument, hold $R\in\C^\times$ fixed with its chosen logarithm.  Then $U_i\in\C((s_i))^\times$ is independent of $\rho_i$;
the common $R$-dependence is restored by homogeneity.  The ordered inertia at the cherry equals the indicated $g_{M_i}$; the opposite node carries its inverse.  The induced cover components, selected lifts and original rooted/logged IO comparisons are those of Lemma~\ref{lem:E1jointsourcecomparison}, applied to the actual degeneration in Lemma~\ref{lem:criticaloneedge}.

Let $p_i,q_i$ be the homogeneous vectors at the two cherry leaves and $\alpha_i=\wt(p_i)+\wt(q_i)$.  Write $\mathcal Y_i$ for its actual intertwining operator and $\Psi_i(-;s_i)$ for the outside actual product or iterate,
including its prescribed physical scale and logarithmic branch. For each $\mu\in\wt(M_i)$, set
\begin{equation}\label{eq:criticalfivefinitecontraction}
 v_{i,\mu}=\pi^{M_i}_\mu
           \mathcal Y_i(p_i,1)q_i\in(M_i)_\mu,\qquad
 \Phi_{T_i}=\sum_{\mu\in\wt(M_i)}
 (U_i\rho_i)^{\mu-\alpha_i}\Psi_i(v_{i,\mu};s_i).
\end{equation}
The expression is understood weight by weight; the factor in
\eqref{eq:criticalfivefinitecontraction} follows from the
$L(0)$-commutator for ordinary intertwining operators
\cite[Definition~2.8]{DH}.  Every $v_{i,\mu}$ is algebraic and
lies in a finite-dimensional weight space.  The outside channel therefore satisfies its \emph{full-state formal} current Ward identity for each $v_{i,\mu}$ by Theorem~\ref{thm:fullsmallward} and Lemma~\ref{lem:E0formalframedtransport}, after the explicit actual source comparison of Lemma~\ref{lem:E1jointsourcecomparison}; the cherry satisfies the three-point identity by Proposition~\ref{prop:binaryward}.  In independent three-point references the comparison includes the compensator \eqref{eq:E1referencecompensation}, not an identification of unchanged scalar expressions. Neither Ward identity is applied to the uncompleted sum over $\mu$.

No full-state Ward identity is applied to a completed cherry
output: for each fixed $\mu$, $v_{i,\mu}$ is an algebraic
vector and the outside Ward statement applies to that vector in its
actual induced framed source.  The eventual use of the full
$\mu$-sum is postponed to the finite-window operation with the
same generic current.

To specify the formal coefficient ring, for either internal edge $f$ let $I_f$ be its descendant leaf set and $\mu_f=h_f+n_f/N_f$, $n_f\ge0$, its intermediate weight. Set $\lambda_A=\wt(a)$, $\lambda_B=\wt(b)$,
$\lambda_C=\wt(c)$, $\lambda_E=\wt(e)$ and $\lambda_D=\wt(d')$.  Lemma~\ref{lem:treesupport}, applied to the original three-operator channel, gives
\begin{equation}\label{eq:criticalfivesupport}
 \Phi_{T_i}
 =R^{\lambda_D-\lambda_A-\lambda_B-\lambda_C-\lambda_E}
 \prod_{f=1}^{2}q_f^{\,h_f-\sum_{j\in I_f}\lambda_j}
 \sum_{n_1,n_2\ge0}C_{n_1,n_2}
                   q_1^{n_1/N_1}q_2^{n_2/N_2},
\end{equation}
where $q_f\in\{u,v\}$ is the corresponding edge ratio, and every $C_{n_1,n_2}$ is a finite-dimensional homogeneous contraction. Absent weights contribute zero.  Choose a common root power $d$ with $N_f\mid d$, and let
$\Lambda_i$ be a finite set of representatives for the resulting exponent classes modulo $d^{-1}\Z^2$.  It embeds the finite sum of generalized-power series \emph{injectively} into the iterated module
\begin{equation}\label{eq:criticalfivecoefficientring}
 \bigoplus_{(\alpha,\beta)\in\Lambda_i}
 \rho_i^\alpha s_i^\beta
 \C((s_i^{1/d}))((\rho_i^{1/d})),\qquad
 |\Lambda_i|<\infty.
\end{equation}
The initial complex powers are formal monomial symbols, with fixed branches; no assertion that their exponents are rational is made.  For fixed $n_1,n_2$ the coefficient in
\eqref{eq:criticalfivesupport} is finite; for fixed $\rho_i$-degree
the remaining outside contraction is an element of the $s_i$-Laurent module, not necessarily a finite scalar sum. Thus \eqref{eq:criticalfivecoefficientring} does not presuppose analytic convergence or exchange the two infinite sums.

The two-ratio support is particularly transparent for these
five trees: their unshifted internal edge ratios are the independent
variables $u,v$.  After selecting the finitely many initial
exponent classes, every fixed $\rho_i$ exponent selects a finite
number of cherry weights; the remaining spectator coefficients form
an ordinary lower-bounded Laurent series in $s_i^{1/d}$, possibly
with a bound depending on that selected cherry exponent.
It would be invalid to infer this for a substituted nonadapted
ratio, or to evaluate the spectator series at a nonzero complex
number before the later convergence theorem.  The finite-root
base change is injective on each distinct formal monomial class,
and localization does not collapse these classes.

Fix a canonical generic five-point current $\sigma$.  Its polar divisor is finite; after the finite DVR extension of Lemma~\ref{lem:criticaloneedge}, a power $t^{m(\sigma)}\sigma$ extends to the original balanced admissible cover.
The extra outside marked section is disjoint from the chosen node; the coefficient restrictions are actual component currents by
\eqref{eq:oneedgeextension}--\eqref{eq:oneedgebasechange}.
The inner and outside Ward identities verified above meet exactly the hypotheses \eqref{eq:oneedgewardhyp} for this selected cut. Lemma~\ref{lem:criticaloneedge} therefore gives, in the localized
iterated module and with $\rho_i=t^{aN}$ after base change,
\begin{equation}\label{eq:criticalfivewardlocalized}
 \Phi_{T_i}(R(\sigma)w)=0.
\end{equation}
The use of the original current extension, the finite node-grade window \eqref{eq:oneedgewindow}, the complete opposite-mode residues and the frame conjugation
\eqref{eq:oneedgephysical} are necessary here: they identify
\eqref{eq:criticalfivewardlocalized} with the physical
contraction \eqref{eq:criticalfivefinitecontraction}, including $U_i^{\mu-\alpha_i}$ and the lowest-weight powers. This is not an appeal to sewing convergence.  The formal coefficient module of \eqref{eq:criticalfivesupport} is torsion-free
under localization in $u,v$ and injects into
\eqref{eq:criticalfivecoefficientring}.  Thus
\eqref{eq:criticalfivewardlocalized} proves
\eqref{eq:criticalfiveward} in the original formal module.
The same argument applies to every current in the finite-state marked seed and hence to its Ward submodule.

For the horizontal assertion one can compute the scalar
coefficient derivatives without referring to analytic sums.
For homogeneous external tests write
$\kappa=\lambda_D-\lambda_A-\lambda_B-\lambda_C-\lambda_E$
and let $a_{n_1},b_{n_2}$ be the two internal exponents in
\eqref{eq:criticalfivesupport}.  For a fixed finite pair,
$C_{n_1,n_2}R^\kappa u^{a_{n_1}}v^{b_{n_2}}$ obeys
\[
\partial_\delta(\cdot)
=\left(\kappa\frac{\delta R}{R}
+a_{n_1}\frac{\delta u}{u}
+b_{n_2}\frac{\delta v}{v}\right)
 C_{n_1,n_2}R^\kappa u^{a_{n_1}}v^{b_{n_2}}
\]
when the external algebraic input and framed output are held fixed.
For $\delta=\partial_u,\partial_v,\partial_R$ these give
respectively $a_{n_1}/u,b_{n_2}/v,\kappa/R$; the chain rule
for the actual positions in \eqref{eq:criticalfivepositions}
produces the corresponding sums of external $L(-1)$ insertions.
On an internal edge the two $L(-1)$ commutator contributions cancel
at each fixed pair of projected weights.  This is the
explicit coefficientwise content of Lemma~\ref{lem:channelderivative}
for these placements.  Differentiation is performed on the finite homogeneous contractions, not on an unverified double analytic series.

Finally, at each vertex the actual twisted intertwining operator satisfies the $L(-1)$-derivative identity
\cite[Definition~2.8(4)]{DH}.  On every fixed pair of
intermediate weights, the chain rule and the $L(-1)$ commutator are finite homogeneous identities.  Internal translation terms cancel after contraction.  By Lemma~\ref{lem:channelderivative},
for every base derivation $\delta$ the remaining terms equal $\Phi_{T_i}(\nabla_\delta w)$, proving the first equality of
\eqref{eq:criticalfivehorizontal}.  Its combination with
\eqref{eq:criticalfiveward} annihilates the differential closure
$\cJ_{\diff}$, proving the second equality.  The finite-rank meromorphic residue construction follows \cite[Sections~5--6]{FS}; the balanced-node step follows Lemma~\ref{lem:criticaloneedge}.
\end{proof}

\begin{proposition}
\label{prop:criticalfiveconvergence}
Fix the original three-operator channel for each $T_i$ in Theorem~\ref{thm:criticalfourward}, including its chosen ordinary simple internal modules and intertwining operators.  Put
\begin{equation}\label{eq:criticalfivedomains}
\begin{aligned}
 \Omega_1=\Omega_3
 &=\{(u,v)\in\C^2:|v|<1,\ |u|(1+|v|)<1\},\\
 \Omega_2=\Omega_4
 &=\{(u,v)\in\C^2:|u|<1,\ |v|<1\},\\
 \Omega_5
 &=\{(u,v)\in\C^2:|u|+|v|<1\},\qquad
 \Omega_i^*=\Omega_i\cap(\C^\times)^2.
\end{aligned}
\end{equation}
On every fixed branch of $\Omega_i^*$, the original intermediate-weight double series converges absolutely and normally:
\begin{equation}\label{eq:criticalfiveoriginalnormal}
 \sum_{\mu_1,\mu_2}
 \sup_{(u,v)\in K}
 \left|C_{\mu_1,\mu_2}
 R^{\lambda_D-\lambda_A-\lambda_B-\lambda_C-\lambda_E}
 u^{\mu_1-\sum_{j\in I_1}\lambda_j}
 v^{\mu_2-\sum_{j\in I_2}\lambda_j}\right|<\infty
 \quad(K\Subset B\subset\Omega_i^*).
\end{equation}
Here $I_j$ is the descendant-leaf set of the $j$th internal edge, $\mu_j$ its intermediate weight ($j=1,2$), and $\lambda_A,\lambda_B,\lambda_C,\lambda_E,\lambda_D$ are the fixed external weights of Theorem~\ref{thm:criticalfourward}.
The number $R\in\C^\times$ and its logarithm are fixed, $B$ is simply connected with fixed logarithms, and each $C_{\mu_1,\mu_2}$ denotes the finite contraction at the specified \emph{pair} of intermediate weights. The conclusion
also holds for a finite direct sum of ordinary internal modules and locally uniformly for $R$ on a fixed compact subset of its nonzero branch.  In the common positive framing the five domains
contain one physical neighborhood of $(z_A,z_B,z_C,z_E)=(7,6,4,0)$.
\end{proposition}
\begin{proof}
We identify the punctured configuration locus. For the positions of \eqref{eq:criticalfivepositions}, let $\Delta_i=R^{-6}\prod_{p<q}(z_p-z_q)$. Its nonmonomial factors, up to nonzero constants and repetition, belong respectively to
\begin{equation}\label{eq:criticalfiveunits}
\begin{array}{c|l}
 i&\text{nonmonomial factors of }\Delta_i\\\hline
 1&1+v,\ 1+u,\ 1+u+uv\\
 2&1-v,\ 1+u,\ 1+uv\\
 3&1+v,\ 1-u,\ 1-u-uv\\
 4&1-u,\ 1-v,\ 1-uv\\
 5&1+u,\ 1-v,\ 1+u-v.
\end{array}
\end{equation}
These factors follow by taking all six differences of the four
positions in \eqref{eq:criticalfivepositions}.  For example, the
six differences for $T_3$ are $R(1-u-uv)$, $R(1-u)$, $R$,
$Ruv$, $Ru(1+v)$, $Ru$, while those for $T_5$ are
$Ru$, $R(1+u-v)$, $R(1+u)$, $R(1-v)$, $R$, $Rv$.
For $T_1,T_3$ we have $|u|(1+|v|)<1$; for $T_5$,
$|u|+|v|<1$; and for $T_2,T_4$ both $|u|,|v|<1$.
The triangle inequality proves that each factor in the right column
is nonzero on its entire $\Omega_i$, not just near the origin.
Consequently
\[
\Delta_i\ne0\ \text{on }\Omega_i^*,\qquad
 \Omega_i\setminus\Omega_i^*
 =\Omega_i\cap\{uv=0\}.
\]
The inequalities \eqref{eq:criticalfivedomains} define connected complete Reinhardt domains containing $0$.  In particular, $\pi_1(\Omega_i^*)\cong\Z^2$ by
\eqref{eq:omegapi1}.  This verification is necessary:
regular singularity on an arbitrary local polydisk would not by itself yield convergence throughout $\Omega_i^*$.

Fix $i$, fix the stated original simple internal modules and the external
homogeneous tensor, and keep the rooted/logged branch inherited from
Definition~\ref{def:fiveoutsetframing}.  The substitution
\eqref{eq:criticalfivepositions}, with $R\ne0$ fixed, is algebraic on
$\Omega_i^*$ and pulls back the finite \emph{framed} \`etale cover of
Lemma~\ref{lem:framedcharts}.  Its permutations $P_u,P_v$ about the two
coordinate axes commute, since $\pi_1(\Omega_i^*)=\Z^2$.
For the two internal relative-grade denominators $N_1,N_2$, take $d>0$
large enough to kill these permutations and to make both relative grades
integral after pullback:
\begin{equation}\label{eq:criticalfiveKummer}
 N_1\mid d,\quad N_2\mid d,\qquad
 P_u^d=P_v^d=1,\qquad
 u=\xi^d,\quad v=\eta^d,\quad
 \widetilde\Omega_i
 =\{(\xi,\eta):(\xi^d,\eta^d)\in\Omega_i\}.
\end{equation}
Include in $d$ the finite orders of any additional marking-root torsors
present in the chosen framed chart.  The unpunctured $\widetilde\Omega_i$ is complete Reinhardt;
its punctured part $\widetilde\Omega_i^*$ is Reinhardt but
\emph{not} complete Reinhardt (it omits both axes).  Since the
coordinatewise radial deformation retract used in
\eqref{eq:omegapi1} applies to the complete unpunctured domain,
$\pi_1(\widetilde\Omega_i^*)=\Z^2$.  After this pullback, the monodromy representation of the
finite cover sends both generators of
$\pi_1(\widetilde\Omega_i^*)\cong\mathbb Z^2$ to the identity:
a lifted generator maps to $d$ times its original generator.  Thus
the finite cover is analytically a disjoint union of trivial degree-one
sheets over $\widetilde\Omega_i^*$, and the selected sheet is the
continuation of the original based, rooted frame (not an arbitrary
relabeling).  This kills the \emph{finite Hurwitz-cover} monodromy
only; it makes no assertion about the Ward connection monodromy.

We verify separately that a chosen sheet extends to the algebraic
smooth normal-crossing setting required by Proposition~\ref{prop:algebraictests}.
Let $Z_i^\circ\subset\mathbb A^2_{\xi,\eta}$ be the principal open obtained
by inverting the Kummer pullbacks of all nonmonomial factors of
\eqref{eq:criticalfiveunits}; thus $\widetilde\Omega_i\subset
(Z_i^\circ)^{\rm an}$ and the deleted configuration divisor is
$D^\circ=\{\xi\eta=0\}\cap Z_i^\circ$.  Form the finite integral
normalization $Y_i\to Z_i^\circ$ in the pullback of the finite
\`etale Hurwitz algebra over $Z_i^\circ\setminus D^\circ$.
Finiteness follows from excellence of this finite-type complex algebra.

For any $p\in\widetilde\Omega_i$, take a sufficiently small
analytic bidisk $P$ about $p$ and put
$P^*=P\setminus\{\xi\eta=0\}$.  The restriction of the selected
framed-cover sheet to $P^*$ is an actual \emph{degree-one} analytic
cover: its monodromy is trivial by \eqref{eq:criticalfiveKummer}
and $P^*$ is connected.  Let $Y_P^{\mathrm{an}}$ be the analytification
of the normalization over $P$.  Normalization commutes with the relevant
complex-analytic finite normalization, so $Y_P^{\mathrm{an}}$ is finite
and normal over $P$.  Take the analytic closure $W_P$ in
$Y_P^{\mathrm{an}}$ of the chosen connected degree-one sheet over
$P^*$.  This closure is a normal irreducible component: normal
analytic spaces have disjoint irreducible components locally, and
$P^*$ meets every irreducible component of a finite dominant cover
of the bidisk.  The finite map $W_P\to P$ is generically an
isomorphism, hence finite and bimeromorphic.  Because the bidisk $P$
is normal, the analytic finite-birational normality criterion gives
\[
 W_P\xrightarrow{\ \sim\ }P.
\]
In local rings the same assertion can be seen without taking a
limit of arbitrary multivalued generators: the rank-one fraction-field
factor of the finite normalization is the base fraction field, whose
integral closure over the regular, integrally closed ring
$\mathcal O_P$ is just $\mathcal O_P$.  Thus the chosen sheet has a
\emph{genuine} extension across both boundary axes and their
intersection, not merely a bounded continuation of its coordinate
functions on a punctured disk.  Equivalently its local direct-factor
idempotent extends in the finite normal algebra, and its factor is
$\mathcal O_P$.  It follows that the normalization is smooth at the
entire extended chosen sheet, the functions $\xi,\eta$ remain
independent local coordinates, and the reduced preimage of the two
axes is a simple normal-crossing divisor there.  Smoothness and the
normal-crossing condition are algebraic open conditions on the finite
algebraic normalization $Y_i$.  Hence an algebraic smooth SNC open
neighbourhood of this chosen analytic sheet exists, and its
analytification contains the entire sheet over
$\widetilde\Omega_i$.  In the resulting coordinates the boundary is
\begin{equation}\label{eq:criticalfiveNC}
 \widetilde D_i=\{\xi\eta=0\}\cap\widetilde\Omega_i.
\end{equation}
This degree-one direct-factor argument is essential: triviality of a
finite cover only on a punctured analytic neighbourhood, by itself,
would not justify the algebraic SNC assertion.

Pull back the algebraic finite-rank Ward quotient and its algebraic
finite-state quotient map to this actual smooth chart.  Its
integrable connection is regular singular along
\eqref{eq:criticalfiveNC} by Theorem~\ref{thm:canonicalregular}
(applied to the algebraic pullback, not to a postulated sewing solution).
The Deligne logarithmic extension
$(\overline{\cQ}_i,\overline\nabla_i)$ exists by
\cite[Chapter~II, Proposition~5.4]{Deligne}.  For each fixed external
algebraic tensor $w$, the image $q_w=[w]$ is an algebraic section of the
actual pulled-back quotient on the punctured smooth chart:
this follows from the pulled-back algebraic quotient map
$\mathcal O\otimes T\twoheadrightarrow\mathcal Q_i$, not from a
hypothetical factorization of the original double series.  Near any
point of the axis boundary, choose a local algebraic frame of the
Deligne logarithmic lattice.  Both the Ward test vector and the
change to this algebraic frame are regular on the complement of the
boundary, hence their coordinates are in the localization by
$\xi\eta$.  Since this is a finite-type algebraic localization,
each coordinate has a finite denominator $\xi^a\eta^b$ locally;
its poles have no other support.  Proposition~\ref{prop:algebraictests}
therefore gives precisely
$q_w\in\Gamma(\widetilde\Omega_i,
\overline{\mathcal Q}_i(*\widetilde D_i))$.
For the finite collection of state generators required to express an
arbitrary logarithmic frame, the minor-identity argument of
Lemma~\ref{lem:jointsupport} supplies a finite common pole bound
near the crossing.  It uses only the algebraic quotient, finite
rank and Noetherianity; it does not presume convergence of the two
intermediate-weight sums.
Theorem~\ref{thm:criticalfourward} supplies the original
coefficientwise formal horizontal functional, while
\eqref{eq:criticalfivesupport} supplies its adapted two-edge,
nonnegative relative-grade support and zero logarithmic degree.
Thus Theorems~\ref{thm:formalconv} and~\ref{thm:globalratio}, with
Proposition~\ref{prop:algebraictests}, apply here to the actual
finite-rank logarithmic connection and the specified algebraic test
section.  In particular the original formal evaluation has a
convergent analytic realization near the crossing, and the resulting
scalar $f_i$ admits continuation throughout
$\widetilde\Omega_i^*$ with only the finite Nilsson behaviour allowed
by the regular-singular connection.  Nothing in this application uses
an equality between different trees, associativity coherence, or
five-input sewing convergence.

We identify the analytic Taylor coefficients with the original
finite intermediate-weight contractions, before any summation or
reordering.  Write
\[
 \kappa_j=h_j-\sum_{a\in I_j}\lambda_a,\qquad
 \mu_j=h_j+\frac{n_j}{N_j},\quad n_j\in\Z_{\ge0},\quad j=1,2.
\]
The original IO $L(0)$-commutator and
\eqref{eq:criticalfivesupport}, with the continued logarithms
$\log u=d\log\xi$ and $\log v=d\log\eta$, give the \emph{formal}
identity
\begin{equation}\label{eq:criticalfiveTaylor}
 R^{-\lambda_D+\sum_a\lambda_a}
 \xi^{-d\kappa_1}\eta^{-d\kappa_2}
 \widehat\Phi_{T_i}(\xi^d,\eta^d)
 =\sum_{n_1,n_2\ge0}C_{n_1,n_2}
 \xi^{dn_1/N_1}\eta^{dn_2/N_2}.
\end{equation}
The exponents are nonnegative integers and the map from the
\emph{pair} of relative-grade indices to the exponent pair is injective.
The local analytic realization therefore has precisely these Taylor
coefficients by uniqueness of the formal germ; no coefficient is
introduced by analytically re-expanding another channel.

Let $H_i$ denote the continuation of the scalar obtained by removing
the displayed lowest-weight powers from $f_i$.  Initially $H_i$ is
the Taylor germ in \eqref{eq:criticalfiveTaylor}.  A small loop in
either coordinate about the origin returns this ordinary Taylor germ
to itself.  Every loop of $\widetilde\Omega_i^*$ is a product of these
two generators, since its fundamental group is $\Z^2$.
Continuation of \emph{this scalar} therefore has trivial monodromy,
and the monodromy theorem makes $H_i$ single-valued on
$\widetilde\Omega_i^*$.  This does not trivialize the monodromy of
other Ward solutions or of arbitrary pure braids.

At a generic point of either coordinate axis, the global finite
Nilsson expression from Theorem~\ref{thm:globalratio} is a finite sum
of powers and powers of logarithms times meromorphic test factors.
To use single-valuedness correctly, fix the tangential variable in a
small open set and group equal transverse exponents modulo $\mathbb Z$.
For one group, continuation once around the axis acts as the scalar
$e^{2\pi i\lambda}$ times the unipotent translation
$\log\xi\mapsto\log\xi+2\pi i$.  If the total scalar is
invariant, the highest nonzero logarithmic coefficient must vanish
unless its degree is zero, and the remaining exponent class must be
integral.  Apply this to the finitely many exponent classes, then to
the other axis.  Consequently $H_i$ has at worst integral Laurent
poles along either axis.  The conclusion refers to this chosen
scalar, not to all flat sections of the Ward connection.  A single-valued function need not be
holomorphic (for example $\xi^{-1}$), so a further argument is
required.  For each axis the meromorphic extension has, locally, a finite
Laurent principal part, whose coefficients are holomorphic in the
tangential variable away from the other axis (meromorphic at its
intersection).  On overlaps these coefficients agree by uniqueness
of the Laurent expansion; hence they continue along the connected
punctured axis.  They vanish on a nonempty punctured segment near
the crossing, since \eqref{eq:criticalfiveTaylor} is a nonnegative
Taylor germ on a whole bidisk, rather than merely an asymptotic
series.  The one-variable identity theorem now kills every principal
part along that entire punctured axis.  This remains valid even
without assuming a \emph{globally} uniform pole order at the
outset: the argument is local on the axis and then propagated through
overlapping neighborhoods.  The function is thus holomorphic at every smooth
point of both axes, and Hartogs (or the already known Taylor germ)
handles the crossing.  We have obtained
\[
 H_i\in\cO(\widetilde\Omega_i),\qquad
 H_i(\xi,\eta)=\sum_{n_1,n_2\ge0}C_{n_1,n_2}
 \xi^{dn_1/N_1}\eta^{dn_2/N_2}
\]
as the ordinary Taylor expansion on the entire complete
Reinhardt root domain, not merely as a local branch continuation.

We bound each original finite weight-pair contraction.
Set $k_j=d/N_j$.  On a closed polydisk
$\overline P_{\boldsymbol\sigma}\subset
 P_{\boldsymbol\varrho}\Subset\widetilde\Omega_i$, let
$M_{\boldsymbol\varrho}=\sup_{\overline P_{\boldsymbol\varrho}}|H_i|$.
Cauchy's double contour formula gives, separately for each
original weight pair,
\[
 C_{n_1,n_2}=\frac{1}{(2\pi i)^2}
 \int_{|\xi|=\varrho_1}\!\int_{|\eta|=\varrho_2}
 \frac{H_i(\xi,\eta)}{\xi^{k_1n_1+1}\eta^{k_2n_2+1}}
 \,d\eta\,d\xi.
\]
All non-sparse Taylor coefficients are zero, since they are zero in
the already identified convergent germ at the crossing; uniqueness
of analytic continuation retains these zeros throughout the
connected Reinhardt domain.  Therefore the bound is for the
original finite contraction $C_{n_1,n_2}$, not a
coefficient of a rearranged series, and gives
\begin{equation}\label{eq:criticalfiveCauchy}
 \sum_{n_1,n_2\ge0}\sup_{(\xi,\eta)\in
 \overline P_{\boldsymbol\sigma}}
 |C_{n_1,n_2}\xi^{dn_1/N_1}\eta^{dn_2/N_2}|
 \le M_{\boldsymbol\varrho}
 \prod_{j=1}^2
 \frac{1}{1-(\sigma_j/\varrho_j)^{d/N_j}}<\infty.
\end{equation}
This applies to the finite contraction at each pair before either
infinite sum is formed; no interchange of sums is assumed.
Complete Reinhardt geometry gives, for every compact
$\widetilde K\Subset\widetilde\Omega_i$, a finite family of such
nested origin-centred polydisks covering $\widetilde K$ (one may
choose their larger polydisks independently inside the domain).
Summing the corresponding positive majorants gives
\[
 \sum_{n_1,n_2\ge0}\sup_{(\xi,\eta)\in\widetilde K}
 |C_{n_1,n_2}\xi^{dn_1/N_1}\eta^{dn_2/N_2}|<\infty.
\]
For $K\Subset B\subset\Omega_i^*$ with fixed logarithms, lift $K$
to the root domain and multiply by the bounded external factor
$R^{\lambda_D-\sum_a\lambda_a}u^{\kappa_1}v^{\kappa_2}$.
This proves \eqref{eq:criticalfiveoriginalnormal} on the entire
specified branch, with the unchanged original coefficients and their
fixed lowest-weight exponents.  Finite sums over simple internal
choices and uniformity in $R$ on compact nonzero logarithmic branches
are immediate.

Finally, at $(z_A,z_B,z_C,z_E)=(7,6,4,0)$,
\begin{equation}\label{eq:criticalcommonpoint}
 \begin{array}{c|c|c|c}
  &R&u&v\\\hline
 T_1&4&1/2&1/2\\
 T_2&4&3/4&2/3\\
 T_3&7&4/7&1/2\\
 T_4&7&6/7&2/3\\
 T_5&6&1/6&2/3
 \end{array}
\end{equation}
and the three coupled quantities $|u|(1+|v|)$ for $T_1,T_3$ and $|u|+|v|$ for $T_5$ are respectively $3/4,6/7,5/6$, all strictly less than $1$. The two remaining pairs satisfy $|u|,|v|<1$.
These strict inequalities persist in one complex physical neighborhood. The full framing and all logarithms there are continued from the same original framed point by Lemma~\ref{lem:framedcharts}.
The proof of the convergence is internal to the Ward quotient, its regular-singular extension, and the Cauchy estimate; the untwisted sewing-factorization results of
\cite{GuiZhangSF} are used only as geometric comparison,
not as a twisted convergence theorem.
\end{proof}

\begin{proposition}\label{prop:criticalconstantmaps}
For $X\in\cC_g$ and $Y\in\cC_h$, fix the formal canonical operator $\mathcal Y^{\can}_{X,Y}$ of Theorem~\ref{thm:binary}.  The realizations
\begin{equation}\label{eq:constantPz94}
 X\bo_{P(z)}Y=X\bo Y,\qquad
 I^z_{X,Y}=\operatorname{ev}_{\log z}
       (\mathcal Y^{\can}_{X,Y})
 \quad(z>0,\ \log z\in\R)
\end{equation}
are compatible with
\begin{equation}\label{eq:constanttransport94}
 \mathsf T^{X,Y}_{s,t}=\id_{X\bo Y}
 \qquad(s,t>0).
\end{equation}
In these representatives, the associator $a_{X,Y,Z}$ of \eqref{eq:positivenormalization} is the fixed map $b_{X,Y,Z}$ of Theorem~\ref{thm:smallassoc}; its actual product--iterate identity holds whenever
\begin{equation}\label{eq:positiveoverlap94}
 z_1>z_2>z_1-z_2>0,
\end{equation}
with the real logarithms and the prescribed common framing. The crossed commutor $c_{X,Y}:X\bo Y\to T_gY\bo X$ satisfies the identity of formal twisted intertwining operators
\begin{equation}\label{eq:criticalrawc}
 \overline c_{X,Y}\,\mathcal Y^{\can}_{X,Y}(x,z)y
 =\Omega_+\bigl(\mathcal Y^{\can}_{T_gY,X}\bigr)(x,z)y
 =e^{zL(-1)}\mathcal Y^{\can}_{T_gY,X}
          (y,e^{\pi i}z)x,
 \quad \log(e^{\pi i}z)=\log z+\pi i.
\end{equation}
Here $y$ on the right is the underlying vector of $T_gY$, and \eqref{eq:criticalrawc} is understood in the completed output, coefficientwise in each $L(0)$-weight.
\end{proposition}
\begin{proof}
Fix $X\in\cC_g$ and $Y\in\cC_h$ and write
$U=U_{X,Y}\in\cC_{gh}$ and
$\mathcal Y=\mathcal Y^{\can}_{X,Y}\in\cI\binom U{X\ Y}$.
For $t>0$ put $I^t=\operatorname{ev}_{\log t}(\mathcal Y)$;
this is a $P(t)$-map in the precise image-space sense of
Definition~\ref{def:Pzproduct}.  For any $W\in\cC_{gh}$,
composition with the evaluation bijection gives
\[
\begin{aligned}
\Hom_V(U,W)&\xrightarrow{\ \sim\ }\cI\binom W{X\ Y}
\xrightarrow{\ \operatorname{ev}_{\log t}\ }\cI_{P(t),\log t}\binom W{X\ Y},\\
 f&\longmapsto\overline f\mathcal Y
\longmapsto\overline f I^t.
\end{aligned}
\]
The first arrow is the complete-operator representing property of
Theorem~\ref{thm:binary}; the second is bijective by
\eqref{eq:Pzinverse} and the \emph{definition} of the $P(t)$-map
space.  A module morphism preserves $L(0)$-weights, so evaluation
commutes with its componentwise extension $\overline f$.
The displayed bijection proves the $P(t)$-universal property for the
same object $U$ for every $t>0$, including all finite
semisimple direct sums and multiplicity factors.  It does not identify
independently selected representing objects without their unique
universal isomorphisms.

For the positive interval $\gamma:s\to t$, the logarithm at the
endpoint is $L_2=\log t$ and its continuation backwards is
$L_1^\gamma=\log s$.  Equation~\eqref{eq:Pztransport} therefore reads
\[
 \overline{\mathsf T_{s,t}^{X,Y}}\,I^s
 =\operatorname{ev}_{\log s}(\mathcal Y)=I^s.
\]
The identity of $U$ satisfies this equation, and the $P(s)$-universal
property makes it its unique solution.  Hence
$\mathsf T_{s,t}^{X,Y}=\id_U$ in these common representatives.
For a winding path the backwards logarithm may instead differ by
$2\pi i\mathbb Z$, so no loop-transport identity is asserted.
For the evaluation used here, if
$x\in X_\alpha$, $y\in Y_\beta$ and
$w'\in(W')_\delta$, Lemma~\ref{lem:IOlogs} gives
\[
 \langle w',\mathcal Z(x,z)y\rangle
 =\exp\bigl((\delta-\alpha-\beta)(\log z-L)\bigr)
 \langle w',\operatorname{ev}_L(\mathcal Z)(x\otimes y)\rangle
\]
for every $\mathcal Z\in\cI\binom W{X\ Y}$.
Thus one evaluation determines \emph{each output-weight coefficient},
not an arbitrary sum of different powers from one numerical value.
The restricted dual separates the product completion, and arbitrary
algebraic inputs are finite sums of homogeneous ones.

Every arrow in \eqref{eq:positivenormalization} has the following
source and target:
\[
\begin{aligned}
 (X\bo_1Y)\bo_1Z
 &\xrightarrow{\mathsf T_{1,2}^{X\bo_1Y,Z}}
 (X\bo_1Y)\bo_2Z
 \xrightarrow{\mathcal A_{3,2}^{-1}}X\bo_3(Y\bo_2Z)\\
 &\xrightarrow{\id_X\bo_3\mathsf T_{2,1}^{Y,Z}}
 X\bo_3(Y\bo_1Z)
 \xrightarrow{\mathsf T_{3,1}^{X,Y\bo_1Z}}
 X\bo_1(Y\bo_1Z).
\end{aligned}
\]
The ordered sectors of each intermediate tensor product are $ghk$.
At the actual $(z_1,z_2,0)=(3,2,0)$ placement,
Theorem~\ref{thm:smallassoc}, including its original scalar
realization and \eqref{eq:assocnormalization83}, identifies
$\mathcal A_{3,2}^{-1}=b_{X,Y,Z}$ in these common representatives;
its existence and uniqueness were obtained there without a pentagon.
The other three arrows are identity maps by the preceding
positive-transport calculation, giving $a_{X,Y,Z}=b_{X,Y,Z}$.

If $z_1>z_2>z_1-z_2>0$, then
$|z_1|>|z_2|>|z_1-z_2|>0$.  The two original contracted
weight sums, with their intermediate-weight projections exactly as in
\eqref{eq:assocproduct83} and \eqref{eq:associterate83}, are
therefore absolutely and locally normally convergent by
Corollary~\ref{cor:smallconvergenceautomatic}.
For homogeneous exterior inputs of weights
$\alpha,\beta,\gamma,\delta$ their common scale factor is
$z_2^{\delta-\alpha-\beta-\gamma}$, and the remaining variables are
$q=(z_1-z_2)/z_2\in(0,1)$ and $1+q$.
Theorem~\ref{thm:smallassoc} identifies these very sums for the
morphism $b$, with the real choices
$\log z_2$, $\log(z_1-z_2)=\log z_2+\log q$, and
$\log z_1=\log z_2+\log(1+q)$ and the inherited rooted framing.
Neither the constant representative nor its identity transports alone
would imply this analytic equality; the latter is the separate
content of the previously proved product--iterate comparison.
Only the continuation paths and coherence identities already established enter this argument.

The canonical target for the exchange is
$W=T_gY\bo_1X\in\cC_{(ghg^{-1})g}=\cC_{gh}$,
where $T_gY\in\cC_{ghg^{-1}}$ has the underlying vector space of $Y$
but the pulled-back vertex action.  Let
$W_{-1}=T_gY\bo_{P(-1)}X$ and write
$\mathcal Z_{-1}=\mathcal Y^{-1}_{T_gY,X}$ for its universal
formal operator, evaluated at $L_{-1}=\pi i$.
Let $\mathcal Z_1=\mathcal Y^{\can}_{T_gY,X}$ be the canonical
formal operator on $W$.  From \eqref{eq:Pzrawcrossed} at $z=1$,
with the source evaluation $\log1=0$, we have the \emph{typed}
equality of completed $W_{-1}$-valued maps
\[
 \overline{\mathcal R^1_{X,Y}}\,I^1_{X,Y}
   =\operatorname{ev}_{0}(\Omega_+\mathcal Z_{-1}).
\]
Fix the upper-half-plane path $\gamma:-1\to1$ of
\eqref{eq:c}.  Its backward continuation of $\log1=0$ is
$\log(-1)=\pi i$, so the actual transport defining equation is
\[
 \overline{\mathsf T_\gamma^{T_gY,X}}
 \operatorname{ev}_{\pi i}(\mathcal Z_{-1})
 =\operatorname{ev}_{\pi i}(\mathcal Z_1).
\]
Because transport is a $V$-module isomorphism and preserves every
$L(0)$-weight, evaluation commutes with its completion.
The injectivity of $\operatorname{ev}_{\pi i}$ from
\eqref{eq:Pzinverse} now gives the \emph{formal} identity
\[
 \overline{\mathsf T_\gamma^{T_gY,X}}\mathcal Z_{-1}
 =\mathcal Z_1.
\]
Furthermore, $\Omega_+$ commutes with a target module map:
its definition \eqref{eq:skewprelim} consists of applying that map
to the original operator and then $e^{zL(-1)}$, with which a
$V$-module map commutes.  Therefore, using only
$c_{X,Y}=\mathsf T_\gamma\mathcal R^1_{X,Y}$ and the preceding
three equations,
\[
 \overline c_{X,Y}\,\mathcal Y^{\can}_{X,Y}
 =\Omega_+\bigl(\overline{\mathsf T_\gamma}\mathcal Z_{-1}\bigr)
 =\Omega_+(\mathcal Z_1).
\]
The skew definition then gives \eqref{eq:criticalrawc}, with
$\log(e^{\pi i}z)=\log z+\pi i$; in particular each term of
exponent $\lambda$ is multiplied by $e^{\pi i\lambda}$, not
$e^{-\pi i\lambda}$.  The output is the genuine $gh$-twisted
module $W$, and $\Omega_+$ has the asserted type
$\cI\binom W{X\ Y}$ since
$T_{g^{-1}}(T_gY)=Y$ as labeled pullback modules.
No untyped identification $Y\cong T_gY$ is made.

Finally, if $\pi_\mu^W$ is a weight projection, the translation
on the completed output is defined by the finite formula
\[
 \pi_\mu^W e^{zL_W(-1)}H(z)
   =\sum_{j\ge0}\frac{z^j}{j!}
       L_W(-1)^j\pi_{\mu-j}^W H(z),
 \qquad H(z)=\mathcal Z_1(y,e^{\pi i}z)x.
\]
There are only finitely many nonzero terms for fixed $\mu$, since
$W$ is ordinary and lower bounded.  Thus this is an equality in
$\overline W\{z\}$, interpreted by its algebraic homogeneous
coefficients; it does not feed an arbitrary element of $\overline W$
into another intertwining operator.  The restricted dual of $W$
separates all these coefficients.  The complete-operator properties
of the skew transform are those of \eqref{eq:skewprelim} and the
cited Du--Huang skew theorem, with no crossed hexagon invoked.
The same upper-half-plane choice is the one fixed in
\eqref{eq:actionrawbraid82}; the identity
$T_kT_g=T_{kgk^{-1}}T_k$ then makes its conjugated sector labels
agree with the current commutor normalization.
\end{proof}

\begin{theorem}\label{thm:criticalpentagon}
For $A,B,C,E\in\cC_G$, the associator $a_{-,-,-}$ of Theorem~\ref{thm:smallassoc} satisfies
\begin{equation}\label{eq:criticalpentagon95}
 a_{A,B,C\bo E}\circ a_{A\bo B,C,E}
 =(\id_A\bo a_{B,C,E})\circ
   a_{A,B\bo C,E}\circ(a_{A,B,C}\bo\id_E)
\end{equation}
as morphisms $(((A\bo B)\bo C)\bo E)\to A\bo(B\bo(C\bo E))$.
\end{theorem}
\begin{proof}
It suffices first to take homogeneous objects
$A\in\cC_g$, $B\in\cC_h$, $C\in\cC_k$ and $E\in\cC_l$.
Write $M=A\bo B$, $Q=B\bo C$, $P=C\bo E$.  Put
\[
\begin{aligned}
 U_1&=((M\bo C)\bo E),&U_2&=((A\bo Q)\bo E),\\
 U_3&=A\bo(Q\bo E),&U_4&=A\bo(B\bo P),&
 U_5&=M\bo P.
\end{aligned}
\]
Every $U_i$ lies in $\cC_{ghkl}$: the sector at every
binary vertex is the ordered product of the sectors of
its children.  The five typed arrows are
\begin{equation}\label{eq:pentagonedges95}
\begin{aligned}
 e_{12}&=a_{A,B,C}\bo\id_E:U_1\longrightarrow U_2,&
 e_{23}&=a_{A,Q,E}:U_2\longrightarrow U_3,\\
 e_{34}&=\id_A\bo a_{B,C,E}:U_3\longrightarrow U_4,&
 e_{15}&=a_{M,C,E}:U_1\longrightarrow U_5,\\
 e_{54}&=a_{A,B,P}:U_5\longrightarrow U_4.
\end{aligned}
\end{equation}
Each $a$ is the fixed $P(1)$-normalized map of
Proposition~\ref{prop:criticalconstantmaps}, hence the
induced maps in \eqref{eq:pentagonedges95} are determined by
binary functoriality; no independent four-point comparison
map is chosen.

Fix homogeneous $a\in A_\alpha$, $b\in B_\beta$,
$c\in C_\gamma$, $e\in E_\epsilon$, and the jointly based
physical placement
\begin{equation}\label{eq:pentagonphysical95}
 (z_A,z_B,z_C,z_E)=(7,6,4,0).
\end{equation}
From \eqref{eq:criticalfivepositions} we obtain the five triples
$(R,u,v)$ of \eqref{eq:criticalcommonpoint}.
The triple associated with each edge is, in the actual physical
coordinate of the relevant three-input subtree,
\begin{equation}\label{eq:pentagonoverlaps95}
\begin{array}{c|c|c}
\text{edge}&\text{triple}&(z_1,z_2)\\ \hline
12&(A,B,C)&(3,2)\\
23&(A,Q,E)&(7,4)\\
34&(B,C,E)&(6,4)\\
15&(M,C,E)&(6,4)\\
54&(A,B,P)&(7,6).
\end{array}
\end{equation}
All five satisfy $z_1>z_2>z_1-z_2>0$, and their
collision ratios are respectively $1/2$, $3/4$, $1/2$, $1/2$, $1/6$.
Here an internal coefficient $q\in Q_\sigma$, $m\in M_\mu$,
or $p\in P_\tau$ is inserted at the physical position of the
rightmost leaf of that subtree; its original coordinate dilation
and lowest-weight powers are those of the corresponding complete
canonical operator.  Thus \eqref{eq:assocdoubleidentity83}
applies with the original real logarithms.  No path
around a discriminant divisor is taken.

Write
$I^t_{X,Y}(x,y)=\operatorname{ev}_{\log t}
 (\mathcal Y^{\can}_{X,Y})(x\otimes y)$, $t>0$, the completed
canonical binary output; a use of $I^t$ as an input to another
$I^s$ means only its \emph{individual homogeneous projections}
are inserted first, followed by the indicated original
intermediate-weight sum.  This is the meaning of each original
four-input channel $\mathscr Y_{T_i}\in\overline{U_i}$.
By Proposition~\ref{prop:criticalfiveconvergence}, for every
restricted-dual test the two-grade sums for all five $T_i$
converge absolutely and locally normally on one physical
neighborhood $\mathscr O$ of \eqref{eq:pentagonphysical95}.
Indeed the five ratios at the displayed point strictly satisfy
\eqref{eq:criticalfivedomains}; the same holds throughout a
small simply connected $\mathscr O$, on which the roots,
external coordinates and logarithms are the \emph{jointly
inherited} ones of Definition~\ref{def:fiveoutsetframing}.
To justify the completed channel, fix a final weight $\eta$.
Choose a basis of the finite-dimensional $(U_i)_\eta^*$;
the convergence theorem applied to these finitely many restricted
tests produces the coordinates of one well-defined vector in
$(U_i)_\eta$.  These vectors, over all $\eta$, form an element
of $\overline{U_i}$; no assertion of convergence in the algebraic
direct sum $U_i$ is made.  A module homomorphism $f$ preserves
weights, so $\pi_\eta\overline f=f\pi_\eta$ and
$\langle d',\overline f\mathscr Y_{T_i}\rangle
 =\langle f'd',\mathscr Y_{T_i}\rangle$ for $d'$ restricted.
In particular, the following desired five equations are
well-typed scalar holomorphic equalities:
\begin{equation}\label{eq:pentagonedgegerms95}
 \left\langle d',\overline{e_{ij}}\mathscr Y_{T_i}
 (a,b,c,e;\boldsymbol z)\right\rangle
 =\left\langle d',\mathscr Y_{T_j}
 (a,b,c,e;\boldsymbol z)\right\rangle,
 \quad d'\in U_j',\quad\boldsymbol z\in\mathscr O.
\end{equation}
We verify the edges separately, using only the three-input
analytic identity of Theorem~\ref{thm:smallassoc}.

Put $W=M\bo C$ and $W'=A\bo Q$, with
$a_0=a_{A,B,C}:W\to W'$.  For each numerical weight
$\nu\in\wt(W')$ and each fixed outer displacement $t=z_C-z_E$
near $4$, set
\[
 \lambda_\nu^{12}(w')=
 \left\langle d', I^t_{W',E}
 (\pi_\nu^{W'}w',e)\right\rangle,
 \qquad w'\in W'.
\]
The functional vanishes on every weight except $\nu$ and
therefore belongs to the \emph{restricted} dual $(W')_{\rm res}'$;
$\lambda_\nu^{12}\circ a_0\in (W)_{\rm res}'$, the
restricted dual of the source module $W$.
Naturality of the binary universal map gives, at every fixed
algebraic $w\in W$,
$\overline{a_0\bo\id_E}I^t_{W,E}(w,e)
 =I^t_{W',E}(a_0w,e)$.
For the triple $(A,B,C)$ the original complete analytic identity
at $(z_1,z_2)=(3,2)$ reads
\[
 \overline{a_0}\,I^2_{M,C}(I^1_{A,B}(a,b),c)
 =I^3_{A,Q}(a,I^2_{B,C}(b,c)),
\]
where each side is expanded in its own intermediate weights.
Pair this identity with $\lambda_\nu^{12}$, applying it first
to each algebraic projected coefficient of either expansion.
It equates the full original inside-grade contractions for the
$\nu$th outer input weight on the two sides of $12$.
The test may depend holomorphically on $t$; the three-input
identity is valid for each such restricted test and all
parameters in the common local overlap.  No unprojected
completed $W'$-vector has been fed into the outer operator.

Now $W=Q\bo E$, $W'=B\bo P$ and
$a_0=a_{B,C,E}:W\to W'$.  For $\nu\in\wt(W')$ define,
with $t=z_A-z_E$ near $7$,
\[
 \lambda_\nu^{34}(w')=
 \left\langle d',I^t_{A,W'}
 (a,\pi_\nu^{W'}w')\right\rangle.
\]
Again $\lambda_\nu^{34}\in(W')_{\rm res}'$ and
$\lambda_\nu^{34}a_0\in(W)_{\rm res}'$; these are
the restricted duals of the target and source respectively.
Using the binary naturality in the \emph{second} variable and
then the original $(B,C,E)$ identity at $(6,4)$ yields
\[
 \overline{a_0}\,I^4_{Q,E}(I^2_{B,C}(b,c),e)
 =I^6_{B,P}(b,I^4_{C,E}(c,e))
\]
against this restricted test after projecting the outer $W'$-
weight $\nu$.  Thus edge $34$ holds for the original fixed
$\nu$th outer coefficient.  The input order here differs from
edge $12$: the exterior $A$ is the first argument
of $I_{A,W'}$, which is why $\lambda_\nu^{34}$ is defined
with the projection on the second argument.

For $23$ fix an actual algebraic coefficient
$q_\sigma=\pi_\sigma^Q I^2_{B,C}(b,c)\in Q_\sigma$.
The $(A,Q,E)$ identity at $(7,4)$ gives
\[
 \overline{a_{A,Q,E}}\,
 I^4_{A\bo Q,E}(I^3_{A,Q}(a,q_\sigma),e)
 =I^7_{A,Q\bo E}(a,I^4_{Q,E}(q_\sigma,e))
\]
after the respective original middle-weight sums.
For $15$ instead fix
$m_\mu=\pi_\mu^M I^1_{A,B}(a,b)\in M_\mu$; the
$(M,C,E)$ identity at $(6,4)$ is
\[
 \overline{a_{M,C,E}}\,
 I^4_{M\bo C,E}(I^2_{M,C}(m_\mu,c),e)
 =I^6_{M,P}(m_\mu,I^4_{C,E}(c,e)).
\]
Finally, for $54$ fix
$p_\tau=\pi_\tau^P I^4_{C,E}(c,e)\in P_\tau$;
the $(A,B,P)$ identity at $(7,6)$ gives
\[
 \overline{a_{A,B,P}}\,
 I^6_{M,P}(I^1_{A,B}(a,b),p_\tau)
 =I^7_{A,B\bo P}(a,I^6_{B,P}(b,p_\tau)).
\]
Every displayed identity is between completed \emph{three-input}
channels whose scalar realizations are the original restricted-
dual weight sums.  Each fixed spectator vector is algebraic,
and the sectors of its module are respectively $hk$, $gh$,
and $kl$.  These are exactly the three allowed triples in
\eqref{eq:pentagonoverlaps95}; no associativity or pentagon
for four inputs is assumed.

Fix an edge $i\to j$, $d'\in U_j'$, and
$K\Subset\mathscr O$.  Its source scalar series is tested by
the restricted vector $e_{ij}'d'\in U_i'$, which is legitimate
because $e_{ij}$ preserves $L(0)$-weights.  The convergence
proposition gives the \emph{two separate} bounds
\begin{equation}\label{eq:pentagonMtest95}
 \sum_{\mu,\nu}\sup_{\boldsymbol z\in K}
    |F^{(i)}_{\mu,\nu}(\boldsymbol z)|<\infty,
 \qquad
 \sum_{\sigma,\tau}\sup_{\boldsymbol z\in K}
    |F^{(j)}_{\sigma,\tau}(\boldsymbol z)|<\infty.
\end{equation}
The dummy indices here deliberately refer to the actual
two internal modules of each tree; in general there is no
term-by-term identification between $(\mu,\nu)$ and
$(\sigma,\tau)$.  For edges $12$ and $34$, first sum the
triple associativity equality against $\lambda_\nu^{12}$ or
$\lambda_\nu^{34}$ for each fixed outer weight $\nu$;
for edges $23$, $15$, $54$, first fix the indicated algebraic
spectator coefficient and apply the corresponding whole
three-input analytic identity.  Only then restore the spectator
weights.  The two M-tests separately show that the tails of
both original finite double truncations tend to zero uniformly
on $K$.  Applying the fixed-coefficient identities before this
limit yields the equality of the two original scalar sums.
This argument permits Tonelli/Fubini \emph{within} each absolutely
summable series but never asserts a bijection of its weight-pair
indices with those of the other tree.  Analytic uniqueness then
gives \eqref{eq:pentagonedgegerms95} on the stated common
holomorphic neighborhood, with the original rooted paths and
real logarithms of Definition~\ref{def:fiveoutsetframing}.

Define $F_{\rm top}=e_{34}e_{23}e_{12}$ and
$F_{\rm bottom}=e_{54}e_{15}$, both in $\Hom_V(U_1,U_4)$.
All five edge identities hold on the same $\mathscr O$;
thus application of the weightwise extended morphisms in their
stated order gives
\begin{equation}\label{eq:pentagonroute95}
 \overline{F_{\rm top}}\mathscr Y_{T_1}
 =\mathscr Y_{T_4}
 =\overline{F_{\rm bottom}}\mathscr Y_{T_1}
 \qquad\text{weakly on }\mathscr O.
\end{equation}
These are genuinely the same final channel, not merely two
solutions of the same Ward differential equations.  There is no
transport along an independent loop, and hence no separate monodromy-cancellation identity is needed.

Put $H=F_{\rm top}-F_{\rm bottom}\in\Hom_V(U_1,U_4)$.
For every $d'\in U_4'$ and homogeneous external algebraic inputs,
\begin{equation}\label{eq:pentagonzeros95}
 \big\langle d',\overline H\,
  \mathscr Y_{T_1}(a,b,c,e;\boldsymbol z)\big\rangle=0
 \quad(\boldsymbol z\in\mathscr O).
\end{equation}
Continue this zero germ along a positive collision-free path
inside the connected $T_1$ domain to a neighborhood of its
$(u,v)=(0,0)$ corner.  Proposition~\ref{prop:criticalfiveconvergence}
identifies the resulting convergent generalized Laurent
expansion with the original intermediate-weight
coefficients.  Write $N=M\bo C$ and $U_1=N\bo E$.
For fixed homogeneous output weight $\eta$ put
\[
\begin{aligned}
 K_{\mu,\nu,\eta}(a,b,c,e)
  =\pi_\eta^{U_1}I^1_{N,E}\bigl(
     \pi_\nu^N I^1_{M,C}(
         \pi_\mu^M I^1_{A,B}(a,b),c),e\bigr)
 \in (U_1)_\eta .
\end{aligned}
\]
This is an algebraic vector: every displayed projection
has finite-dimensional image before the next $I^1$ is applied.
For the tree $T_1$ the contribution of this fixed pair of
intermediate weights is
\[
 R^{\eta-\alpha-\beta-\gamma-\epsilon}
 u^{\nu-\alpha-\beta-\gamma}
 v^{\mu-\alpha-\beta}
 K_{\mu,\nu,\eta}(a,b,c,e).
\]
The powers follow successively from the physical displacements
$Ruv$ (the $AB$ cherry), $Ru$ (the $ABC$ stem), and $R$
(the final stem), by the $L(0)$-commutator.  For fixed input and
output weights the pair of exponents distinguishes each
numerical pair $(\mu,\nu)$.  Distinct fractional congruence
classes remain distinct after the common Kummer root
substitution of the convergence theorem.  As $H$ preserves
weights, the vanishing in \eqref{eq:pentagonzeros95},
coefficient uniqueness for convergent generalized power series,
and restricted-dual separation give
\[
 H K_{\mu,\nu,\eta}(a,b,c,e)=0
 \quad\text{for all }\mu,\nu,\eta,a,b,c,e.
\]
If different simple summands have the same numerical weight,
$\pi_\mu^M$ and $\pi_\nu^N$ include all of them,
and no spurious individual separation is required.
Lemma~\ref{lem:generation}, applied successively to the three
binary vertices, states that these algebraic vectors span the
whole $U_1$, not merely its completion.  Hence $H=0$.
This is precisely \eqref{eq:criticalpentagon95}.
Finally finite additivity extends the identity from homogeneous
sectors and homogeneous test vectors to all the indicated finite
semisimple objects.  No later crossed hexagon, rigidity, or
coherence result was used.
\end{proof}

\subsection{The crossed hexagons}\label{sec:hexagons}
For $X\in\cC_g$, $Y\in\cC_h$, and $Z\in\cC_k$, we use the fixed branches of \eqref{eq:Pzrawcrossed} and
\eqref{eq:criticalrawc} to prove
\eqref{eq:hex1}--\eqref{eq:hex2}.  Throughout this subsection
$\pi_\eta^M$ denotes the $L_M(0)$-weight projector of the ordinary twisted module $M$; for a twisted intertwining operator, $\operatorname{ev}_L$ means evaluation with the specified logarithm $L$ (Definition~\ref{def:Pzproduct}).

\begin{lemma}\label{lem:literalskew}
Let $X\in\cC_g$, $Y\in\cC_h$, and put
\[
 N=T_gY\bo X,\quad
 \mathcal Z=\mathcal Y^{\can}_{T_gY,X},\quad
 \mathcal K=\overline c_{X,Y}\,\mathcal Y^{\can}_{X,Y}
             =\Omega_+\mathcal Z.
\]
For $x\in X_\alpha$, $y\in Y_\beta$ set $\lambda=\alpha+\beta$ and $Z_\nu=\pi_\nu^N\mathcal Z(y,1)x$. For $r>0$, with $\log r\in\R$, one has
\begin{equation}\label{eq:gradedskew}
 \pi_\eta^N\mathcal K(x,r)y
 =r^{\eta-\lambda}
 \sum_{\substack{j\ge0\\\eta-j\in\wt(N)}}
 \frac{e^{\pi i(\eta-j-\lambda)}}{j!}
 L_N(-1)^j Z_{\eta-j}.
\end{equation}
Here the sum is finite for each $\eta$, and the continuation along $\log r\mapsto\log r-\pi i$ satisfies
\begin{equation}\label{eq:skewinverse}
 \operatorname{ev}_{\log r-\pi i}(\mathcal K)(x\otimes y)
  =e^{-rL_N(-1)}\mathcal Z(y,r)x,
 \qquad \Omega_-\mathcal K=\mathcal Z.
\end{equation}
For a complete twisted intertwining operator $\mathcal O\in\cI\binom{W_3}{W_1\ W_2}$, $m\in W_1$, $b\in W_2$, the following are identities after each fixed-output-weight projection:
\begin{equation}\label{eq:skewtranslations96}
 \begin{aligned}
 \mathcal O(e^{tL_{W_1}(-1)}m,z)b
   &=\mathcal O(m,z+t)b,\\
 \mathcal O(m,z)e^{tL_{W_2}(-1)}b
   &=e^{tL_{W_3}(-1)}\mathcal O(m,z-t)b.
 \end{aligned}
\end{equation}
The series are Taylor expansions at $t=0$ with $|t|<|z|$ and compatible logarithms.  In each of the adjacent and block endpoints used below, the resulting identities of original three-input channel functions
hold as analytic germs on the indicated collision-free overlap.
\end{lemma}
\begin{proof}
By \eqref{eq:criticalrawc}, $\mathcal K=\Omega_+\mathcal Z$ with the indicated $V$-module types; in particular $T_{g^{-1}}(T_gY)=Y$.  By the $L(0)$ commutator formula for an ordinary twisted intertwining operator (\cite[Definition~2.8]{DH}),
\[
 \mathcal Z(y,r)x
   =\sum_{\nu\in\wt(N)}r^{\nu-\lambda}Z_\nu,
 \qquad L_N(0)Z_\nu=\nu Z_\nu,
\]
coefficientwise in $\overline N$.  Equation~(3.1) of
\cite[Section~3]{DH} gives, on the branch
$\log(e^{\pi i}r)=\log r+\pi i$,
\[
 \mathcal K(x,r)y
 =e^{rL_N(-1)}
  \sum_{\nu\in\wt(N)}e^{\pi i(\nu-\lambda)}
                     r^{\nu-\lambda}Z_\nu.
\]
As $L_N(-1)N_\nu\subseteq N_{\nu+1}$,
\[
 \pi_\eta^N e^{rL_N(-1)}Z_\nu
  =\begin{cases}
    \displaystyle\frac{r^{\eta-\nu}}{(\eta-\nu)!}
                L_N(-1)^{\eta-\nu}Z_\nu,
       &\eta-\nu\in\Z_{\ge0},\\
    0,&\text{otherwise}.
   \end{cases}
\]
For fixed $\eta$, only $0\le j\le\Re(\eta)-B_N$ can contribute, where $B_N\in\R$ is a lower bound for $\{\Re(\nu):\nu\in\wt(N)\}$; thus the projected sum is algebraically finite.  More precisely, every nonzero summand has
$0\le j\le\lfloor\Re(\eta)-B_N\rfloor$ and
$\nu=\eta-j$; if this upper bound is negative, the entire
projection is zero.  The projector $\pi^N_{\eta-j}$ is the
projector onto the whole numerical weight space, including
all simple summands and multiplicity copies.  The mode cosets are
therefore not identified or reordered.  Substitution yields
\eqref{eq:gradedskew}; there is no need to rearrange
an unprojected infinite series.

Continue the same expression along $L(t)=\log r-\pi it$, $0\le t\le1$. Every output weight $\eta$ acquires $e^{-\pi i(\eta-\lambda)}$ at $t=1$.  Hence
\[
 \begin{aligned}
 \pi_\eta^N\operatorname{ev}_{\log r-\pi i}
             (\mathcal K)(x\otimes y)
 &=r^{\eta-\lambda}
   \sum_{\substack{j\ge0\\\eta-j\in\wt(N)}}
   \frac{(-1)^j}{j!}L_N(-1)^jZ_{\eta-j}\\
 &=\pi_\eta^N e^{-rL_N(-1)}\mathcal Z(y,r)x.
 \end{aligned}
\]
It follows coefficientwise that the first identity in
\eqref{eq:skewinverse} holds.  Applying the defining
$\Omega_-$ transformation, with its prescribed $-\pi i$ substitution, gives
\[
 \begin{aligned}
 (\Omega_-\mathcal K)(y,r)x
 &=e^{rL_N(-1)}
     \operatorname{ev}_{\log r-\pi i}
                         (\mathcal K)(x\otimes y)\\
 &=\mathcal Z(y,r)x.
 \end{aligned}
\]
Equality for $r>0$ and each homogeneous input/output weight implies equality of the formal operators by uniqueness of their generalized Laurent coefficients. This is the inverse skew-symmetry assertion of
\cite[Section~3, Theorem~3.1]{DH}.
The lower path $r e^{-\pi it}$ is distinct from the
\emph{upper} transport $-1\to1$ defining $c$ in
\eqref{eq:Pzrawcrossed}: its backward logarithm at
$-1$ is $+\pi i$.  Thus \eqref{eq:criticalrawc} and the present inverse substitution use compatible, explicitly specified branches; no triviality of winding monodromy is asserted.

The $L(-1)$-derivative axiom and the translation commutator for $\mathcal O$ are
\[
 \partial_z\mathcal O(m,z)
  =\mathcal O(L_{W_1}(-1)m,z),\qquad
 [L_{W_3}(-1),\mathcal O(m,z)]
  =\mathcal O(L_{W_1}(-1)m,z),
\]
where the second bracket uses the right action of $L_{W_2}(-1)$.  Exponentiating these two identities coefficientwise gives \eqref{eq:skewtranslations96}. More explicitly, for a fixed output weight $\eta$,
\[
 \pi_\eta^{W_3}e^{tL_{W_3}(-1)}v
 =\sum_{\substack{j\ge0\\\eta-j\in\wt(W_3)}}
    \frac{t^j}{j!}L_{W_3}(-1)^j
           \pi_{\eta-j}^{W_3}v,
\]
which is finite because $W_3$ is ordinary and lower bounded. This projector formula does \emph{not} imply that an exponential
of a completed \emph{input} can be substituted into an IO: if the input
has unbounded weights, infinitely many input-weight summands can
contribute to one fixed output weight.  We use the following precise
coefficient convention instead.  For algebraic homogeneous $m,b$ and
$d'\in (W_3')_\eta$, write
$F(z)=\langle d',\mathcal O(m,z)b\rangle=Cz^{\eta-\wt(m)-\wt(b)}$
on the selected logarithmic branch, by Lemma~\ref{lem:IOlogs}.
For the first translation identity, the coefficient of $t^n$ after
pairing is
$\langle d',\mathcal O(L(-1)^nm,z)b\rangle/n!
=F^{(n)}(z)/n!$.  Hence its sum is precisely the
Taylor series of $F(z+t)$ and converges for $|t|<|z|$ on a compatible
branch.  For the second translation identity the commutator gives,
for the same algebraic inputs,
\[
 \big\langle d',\mathcal O(m,z)e^{tL_{W_2}(-1)}b\big\rangle
 =\big\langle e^{tL_{W_3'}(1)}d',
             \mathcal O(m,z-t)b\big\rangle.
\]
The left is initially a formal power series in $t$; on the right
$e^{tL_{W_3'}(1)}d'$ is a \emph{finite polynomial} because the
restricted dual is ordinary and lower bounded.  Each of its
finitely many scalar IO coefficients is a fixed complex power of
$z-t$, so the right side converges for $|t|<|z|$.  This proves the
analytic Taylor interpretation, including nonintegral weights,
without claiming a general IO action on completed inputs.
Linearity handles algebraic vectors with finitely many weights.

For the first adjacent endpoint, let $R>0$ and $0<r<R/2$.  The $-r$ substitution in
\eqref{eq:skewtranslations96} yields
\[
 \mathcal O(e^{-rL_{W_1}(-1)}m,R)b
       =\mathcal O(m,R-r)b.
\]
The terminal ratio $q=r/(R-r)$ has $0<q<1$, and the original iterate in the $(R-r,r)$ coordinates is normally convergent by Corollary~\ref{cor:smallconvergenceautomatic}. Its coefficientwise Taylor expansion is the one
obtained above.  On every compact subdomain of this chamber the scalar original intermediate-grade series (after pairing with $W_3'$) converges normally; Cauchy integration on a smaller polydisc identifies the Taylor
coefficients of the sum with those of the continued holomorphic germ.  In particular the finite output-grade identities suffice: no exchange between two uncontrolled infinite weight sums is used.

For the other adjacent endpoint let $R>2s>0$, $\mathcal B\in\cI\binom M{U\ W}$ and $\mathcal O\in\cI\binom D{A\ M}$ be the relevant canonical operators.  If $a\in A$, $u\in U$, $v\in W$, and $d'\in D'$ is homogeneous, then
\begin{equation}\label{eq:skewanchored96}
 \begin{aligned}
 &\big\langle e^{sL_{D'}(1)}d',
       \mathcal O(a,R-s)e^{-sL_M(-1)}
                             \mathcal B(u,s)v\big\rangle\\
 &\qquad=\big\langle e^{sL_{D'}(1)}d',
        e^{-sL_D(-1)}\mathcal O(a,R)
                             \mathcal B(u,s)v\big\rangle\\
 &\qquad=\big\langle d',
                   \mathcal O(a,R)\mathcal B(u,s)v\big\rangle.
 \end{aligned}
\end{equation}
The first equality is the second formula of
\eqref{eq:skewtranslations96}, with $t=-s$, but its completed
middle input requires a specific order of evaluation.  Set
$v_\mu(s)=\pi^M_\mu\mathcal B(u,s)v\in M_\mu$.
For each fixed $\mu$, first take the convergent Taylor
series in the translation parameter $-s$ and pair its result with
the finite algebraic vector $e^{sL_{D'}(1)}d'$.
The inequality $s<R-s$ (equivalently $R>2s$) is precisely the
radius condition needed for this inner Taylor series.  The
preceding translation calculation gives the \emph{exact} identity
\[
\left\langle e^{sL_{D'}(1)}d',
\mathcal O(a,R-s)e^{-sL_M(-1)}v_\mu(s)\right\rangle_{\mathrm{Taylor}}
=\langle d',\mathcal O(a,R)v_\mu(s)\rangle.
\]
Only after this fixed-$\mu$ equality is established do we
sum over $\mu$.  The right-hand series is the original product
channel at $(R,s)$, and is absolutely locally normally convergent
since $R>s>0$, by Corollary~\ref{cor:smallconvergenceautomatic}.
Thus the left-hand expression in \eqref{eq:skewanchored96} means
this ordered (translation-Taylor first, intermediate-weight second)
scalar contraction.  We do \emph{not} assert absolute convergence
of the unrestricted two-index $(\mu,j)$ series, nor an extension
of $\mathcal O$ to all of $\overline M$.  The outer sum of these
already evaluated Taylor coefficients is locally normal by
its termwise equality to the original product series.
For the second, the restricted contragredient relation
$\langle L_{D'}(1)d',w\rangle
  =\langle d',L_D(-1)w\rangle$
implies
$\langle e^{sL_{D'}(1)}d',e^{-sL_D(-1)}w\rangle
=\langle d',w\rangle$.
The exponential on $D'$ is finite because $L_{D'}(1)$ strictly lowers the ordinary weight.  On $d'\in D'$ the exponential is a finite algebraic sum;
the middle-input exponential is handled in the fixed-$\mu$
Taylor order just specified, not by the output-projector
formula alone.

Finally, the right side of \eqref{eq:skewanchored96} is the original three-input product in the chamber $R>s>0$; the left side is evaluated in the smaller safe overlap $R>2s$.  By Corollary~\ref{cor:smallconvergenceautomatic}, the
original intermediate-grade contractions on both sides are absolutely and normally convergent there, after pairing with the indicated finite restricted-dual vectors. Consequently the ordered scalar contraction just defined
agrees with the original normally convergent three-input product
and determines the same analytic germ.  This does not establish
normal convergence for an arbitrarily reordered double series
of translation orders and intermediate weights; only the
precisely displayed contraction is used in the endpoint
comparison.  Derivatives of the original channel are controlled
by Cauchy formulas on smaller compact subdomains.  The same argument applies to a block exchange after fixing the other internal output
grade: \eqref{eq:skewinverse} is then a finite projector identity, and the remaining actual grade sum is governed by the corresponding normal-convergence statement. This proves the asserted endpoint legality without
assuming either crossed hexagon.
\end{proof}

\begin{lemma}\label{lem:fullframes}
Let
\[
\begin{gathered}
 B_3=F_3(\C),\qquad B_0=F_2(\C^*),\\
 \vartheta:B_3\xrightarrow{\ \sim\ }B_0\times\mathbb A^1,\\
 (q_1,q_2,q_3)\longmapsto
 (q_1-q_3,q_2-q_3,p=q_3).
\end{gathered}
\]
Fix the based admissible-cover component, inertia labels, marked lifts and special roots of Lemma~\ref{lem:framedcharts}. The anchored three-input differential Ward quotient $\cQ_{\rm an}$ of Theorem~\ref{thm:finiteconn} pulls back along
$\vartheta$ to the corresponding quotient $\cQ_{\rm phys}$ on $B_3$ (and their finite framed covers).  For the output restricted contragredient $D'$ the state-frame identification is
\begin{equation}\label{eq:fullframegauge97}
 I_p=\id_{X\otimes Y\otimes Z}\otimes U_p,
 \qquad U_p=e^{pL_{D'}(1)},\qquad
 I_{p+q}=I_pI_q,
 \quad I_p^{-1}=I_{-p}.
\end{equation}
It intertwines the seed Ward actions and the differential closures.  In fixed physical input frames the common-translation connection is
\begin{equation}\label{eq:fulltranslation97}
 \nabla^{\rm phys}_{\partial_p}
 =\partial_p+L_X(-1)+L_Y(-1)+L_Z(-1).
\end{equation}
Consequently the resulting finite-rank flat quotient admits homotopy-invariant transport on the \emph{full} marked-lift and special-root atlas over every collision-free bounded path and
endpoint-fixed homotopy.  No trivial-monodromy assertion on the unframed configuration space is made.
\end{lemma}
\begin{proof}
Put $z_i=q_i-p$ and let $\tau_p(z)=z-p$. Then $z_i\ne0$, $z_1\ne z_2$, and $B_3\simeq B_0\times\mathbb A^1$ as displayed.  Translation $\tau_p$ fixes $\infty$ and takes the physical cover with
branch points $(q_1,q_2,p)$ to the anchored cover with branch points $(z_1,z_2,0)$.  The input local parameters agree: $(z-q_i)=(\tau_p(z)-z_i),\; i=1,2,3,\; (z_3=0)$. At $\infty$, let $s=z^{-1}=u^{N_\infty}$ be the
physical coarse/root coordinates.  The anchored coordinates are
\begin{equation}\label{eq:infinityroot97}
 w=(z-p)^{-1}=\frac{s}{1-ps},\qquad
 \rho_{p,N_\infty}(u)
 =u(1-pu^{N_\infty})^{-1/N_\infty},
 \quad \rho_{p,N_\infty}(u)^{N_\infty}=w.
\end{equation}
Here the $N_\infty$th root is the unique power series with constant term $1$; in particular it preserves the specified leading root.
The binomial expansion, in the actual special-root coordinate, is
\[
 \rho_{p,N}(u)=u+\frac pN u^{N+1}
       +\frac{N+1}{2N^2}p^2u^{2N+1}+O(p^3).
\]
Differentiating at $p=0$ gives exactly
$N^{-1}u^{N+1}\partial_u$; by
\eqref{eq:relativeautlie} its infinitesimal module action is
$-L_{D'}(1)$, not $-NL_{D'}(1)$.
Since $w(s,p)=s/(1-ps)$ satisfies
$w(w(s,p),q)=w(s,p+q)$, the same is true of the uniquely based
$N$th roots.  The action
\eqref{eq:relativeautlie} gives
\[
R_{D',h_{D'}}(\rho_{p,N_\infty})
 =e^{-pL_{D'}(1)},\qquad
 R_{D',h_{D'}}(\rho_{p,N_\infty})^{-1}=U_p.
\]
The exponentials are finite on every $d'\in D'$: $L_{D'}(1)$ strictly decreases ordinary weights, which are bounded below.  Formula \eqref{eq:relativefieldcov}, including its $V$-state transformation as in\ \cite[Sections~5--6]{FS},
now gives for every invariant finite-rank current $\sigma$ on the physical cover
\begin{equation}\label{eq:seedcovariance97}
 I_p\,R_{\rm phys}(\sigma)
 =R_{\rm an}(\tau_p^*\sigma)\,I_p
 \qquad\text{on }X\otimes Y\otimes Z\otimes D'.
\end{equation}
This is a conjugacy of residue operators in the chosen frames, not an identification of unframed conformal blocks.
At every finite input, the centered coordinate identity
$(z-q_i)=((z-p)-(q_i-p))$ makes the centered state action and
its differential Jacobian identical; at the output the
comparison is $R_{D',h_{D'}}(\rho_{p,N_\infty})^{-1}=U_p$.
For $b\in V_{\le E}$ the entire finite jet matrix, including all
$L(j)b$ terms for $j>0$, is included by
\eqref{eq:relativefieldcov}.  Both the input and output
normalized residues use the same inertia factors $1/N_i$,
and the based/deck pullback preserves their order.  Thus each
individual local residue on the left becomes its corresponding
anchored local residue on the right, giving
\eqref{eq:seedcovariance97} after summation.  This remains true
for $b=\mathbf1$ and $b=\omega$: the coarse translation is
M\"obius, its Schwarzian derivative is zero, and the remaining
finite Virasoro coordinate terms are exactly those already in
the full-state covariance law.  Lemma~\ref{lem:canonicalfinite}
ensures that the two sides act in the same finite state-current
bundle; no infinite state completion has been invoked.  Thus $I_p\cJ_{E,\rm phys}=\cJ_{E,\rm an}$. Also
$\rho_{p,N_\infty}\circ\rho_{q,N_\infty}
 =\rho_{p+q,N_\infty}$, so \eqref{eq:fullframegauge97}
follows.  The coarse transformation in
\eqref{eq:infinityroot97} is M\"obius and has zero
Schwarzian.  The \emph{root} map $\rho_{p,N}$ is generally not
M\"obius when $N>1$; for example
$\{\rho_{p,N},u\}=p(N^2-1)u^{N-2}+O(p^2)$.
This is not an additional Ward scalar: its full state-coordinate
correction is already built into \eqref{eq:relativefieldcov}, and
the infinitesimal action of the relevant positive Virasoro
subalgebra has no central cocycle.  The global translation stress
residue therefore has no unaccounted central term.

To verify the connection, use the global translation vector field $\partial_z$ and its stress-tensor current.  At a finite marking its physical local vector field is
$\partial_{z-q_i}$; at infinity, in $s=z^{-1}=u^{N_\infty}$,
it is $-s^2\partial_s=-N_\infty^{-1}u^{N_\infty+1}\partial_u$.
By \eqref{eq:relativeautlie}, the latter coordinate action is
$+L_{D'}(1)$; the actual \emph{outgoing} Ward residue has the
opposite sign because it acts in the restricted-contragredient
output leg.  Thus the global current contributes
$+\sum_iL_{M_i}(-1)-L_{D'}(1)$ to the seed Ward relation.
Hence the normalized seed relation on the anchored tensor is
\[
L_X(-1)+L_Y(-1)+L_Z(-1)-L_{D'}(1)\equiv0\pmod{\cJ_{E,\rm an}}.
\] There is no ramified scalar residue: for the finite constant
vector field and the infinity quadratic vector field, respectively, the vacuum term in \eqref{eq:stress} has powers $u^{-N_i-1}$ and $u^{N_\infty-1}$, neither of which is $u^{-1}$.  The $1/N_i$ inertia normalization
is that of Definition~\ref{def:canonicalmarked}. Write $D=L_X(-1)+L_Y(-1)+L_Z(-1)$. For a section of the algebraic external tensor, the anchored connection along the $p$-factor is $\partial_p$, whereas
\eqref{eq:Tconn} gives
$\nabla^{\rm phys}_{\partial_p}=\partial_p+D$. Since $[I_p,D]=0$ and $\partial_pI_p=L_{D'}(1)I_p$,
\begin{equation}\label{eq:connframe97}
 \nabla^{\rm an}_{\partial_p}I_p
       -I_p\nabla^{\rm phys}_{\partial_p}
 =\bigl(L_{D'}(1)-D\bigr)I_p
 \equiv0\pmod{\cJ_{E,\rm an}}.
\end{equation}
For each relative placement derivative $\delta(z_i)$ the gauge is independent of $z_i$;
\eqref{eq:seedcovariance97} and \eqref{eq:Tconn}
therefore give the analogous intertwining relation. Applying covariant derivatives repeatedly proves
\[
 I_p\cJ_{\diff,\rm phys}
      =\cJ_{\diff,\rm an},\qquad
 \cQ_{\rm phys}\xrightarrow[\ I_p\ ]{\ \sim\ }
      \vartheta^*\operatorname{pr}_{B_0}^*\cQ_{\rm an}
\]
with their integrable connections.  The inverse is $I_{-p}$; local freeness and finite rank follow from Theorem~\ref{thm:finiteconn}. Dualizing \eqref{eq:fulltranslation97}, a horizontal scalar
functional satisfies
\[
 \partial_p\Phi(w)=\Phi(Dw)
 \quad (w\in X\otimes Y\otimes Z\otimes D').
\]

We verify the framing.  Over the smooth Hurwitz chart, a lift of each marking is a finite deck choice.  Once its physical coarse coordinate $x_i$ is fixed, two special roots $u_i^{N_i}=x_i$ differ by precisely one
$\zeta_i\in\mu_{N_i}$: if $u_i'=u_i a(u_i^{N_i})$, then $(u_i')^{N_i}=u_i^{N_i}$ forces $a^{N_i}=1$. Thus the choices of marked lifts and \emph{leading} root phases form a finite \'etale cover of the smooth Hurwitz chart;
all higher root coefficients are uniquely determined. For a root transition $\zeta_i=e^{2\pi i r_i/N_i}$, Lemma~\ref{lem:relativecoordinates} prescribes on each simple ordinary block
\begin{equation}\label{eq:rootgauge97}
 S_{i,r_i}=R_{M_i,h_i}(\zeta_i u_i)
 =e^{-2\pi i r_i(L_{M_i}(0)-h_i)}.
\end{equation}
Because $[L(0),L(j)]=-jL(j)$ and $j\in\Z$, $S_{i,r_i}L(j)=L(j)S_{i,r_i}$ for $j=-1,1$. The associated field covariance is
\eqref{eq:relativefieldcov}; the factors
$S_{i,r_i}$ are \emph{not} discarded, nor is monodromy around a noncontractible loop declared trivial. A deck transition $k\in H$ replaces $(g_i,M_i,b)$ by $(kg_ik^{-1},T_kM_i,kb)$, where $Y_{T_kM_i}(kb,x)=Y_{M_i}(b,x),\; T_kT_l=T_{kl}$.
The cocycles can be checked on their actual types.
For fixed coarse coordinate, root phases satisfy
$S_{i,r}S_{i,s}=S_{i,r+s}$, with $r+s$ understood modulo $N_i$:
$S_{i,N_i}=\id$ follows since the relative weights belong to
$N_i^{-1}\mathbb Z_{\ge0}$.  For a deck label $k\in H$ the typed
identification satisfies
$T_kR_{M,h}(\rho)=R_{T_kM,h}(\rho)T_k$, because $k$
fixes $\omega$ and hence all Virasoro operators.  Therefore a
successive deck pair $k,l$ gives $T_kT_l=T_{kl}$, retaining the
order even when $kl\ne lk$.  More general changes of coarse
special coordinates use the \emph{full} group action
$R_{M,h}(\rho)$, not merely the leading root phase;
the coordinate-action composition law of
Lemma~\ref{lem:relativecoordinates} supplies their transition
cocycle, with the opposite-to-point-map convention used there.
Finally $U_{p+q}=U_pU_q$ and $[U_p,S_{i,r}]=0$ for the output
root gauge (the relative-weight phase commutes with $L(1)$).
Consequently the transition cocycles on the complete marked-root
atlas preserve both the seed Ward actions and the connection,
by \eqref{eq:relativefieldcov} and \eqref{eq:connframe97}.

Finally, for a bounded collision-free path and a fixed-endpoint homotopy, choose a compact neighborhood $K\Subset B_3$ of their images, and pull the finite marked-root cover back to $K$.
Fix an initial lifted point; lift each path and its homotopy from this point, without assuming a global section over $K$. The full framed quotient is the above finite-rank flat connection, pulled back from Theorem~\ref{thm:finiteconn};
its analytic parallel transport exists along these paths and is invariant under fixed-endpoint homotopy
\cite[Chapter~I, Section~2]{Deligne}.
The homotopy-lifting property of the finite \'etale cover identifies the \emph{terminal full frame} for two homotopic lifts starting at the same framed point. If a different endpoint root chart is subsequently chosen,
\eqref{eq:rootgauge97} gives the required (generally
nonidentity) change of coordinates.  The claim holds on the full framed atlas.
\end{proof}

\begin{proposition}\label{prop:braidendpoints}
Let $X\in\cC_g$, $Y\in\cC_h$, $Z\in\cC_k$ be ordinary twisted modules.  Put $e_t=e^{-\pi it}$ $(0\le t\le1)$ and use the lower-half-turn convention of \eqref{eq:skewinverse}.  The following four paths in $F_3(\C)$ are collision-free:
\begin{align}
 \gamma_{12}(t)&=\left(\frac{R+s+re_t}{2},
                         \frac{R+s-re_t}{2},0\right),
       &r&=R-s,\quad 0<r<s,\label{eq:adjacent1298}\\
 \gamma_{23}(t)&=\left(R,\frac{s(1+e_t)}2,
                         \frac{s(1-e_t)}2\right),
       &&R>2s>0,\label{eq:adjacent2398}\\
 \gamma_{X\mid YZ}(t)&=\left(\frac{R(1+e_t)}2,
                    \frac{R(1-e_t)}2+s,\frac{R(1-e_t)}2\right),
       &&R>s>0,\label{eq:blockX98}\\
 \gamma_{XY\mid Z}(t)&=\left(\frac{s(1+e_t)}2+r,
                    \frac{s(1+e_t)}2,\frac{s(1-e_t)}2\right),
       &&0<r<s.\label{eq:blockZ98}
\end{align}
The first two paths exchange adjacent inputs; the last two exchange $X$ with the block $YZ$, and the block $XY$ with $Z$, respectively.  Their endpoint identities are precisely the fixed $c_{X,Y}$, $c_{Y,Z}$, $c_{X,Y\bo Z}$ and
$c_{X\bo Y,Z}$, together with the specified tensorators $J^g_{Y,Z}$ and positive-real associators as needed to change the displayed bracketing.  These statements refer to the original normally convergent weight sums, with the
complete marked-root framing of Lemma~\ref{lem:fullframes}. The chronological crossing words of the two block paths are $\sigma_1\sigma_2$ and $\sigma_2\sigma_1$, respectively.
\end{proposition}
\begin{proof}
Write $(p_X,p_Y,p_Z)$ for a configuration in
\eqref{eq:adjacent1298}--\eqref{eq:blockZ98}.  The pairwise differences, with
$e=e_t$, are the following actual functions, rather than
schematic strands:
\begin{equation}\label{eq:braidactualdifferences}
\begin{array}{c|ccc|c}
 &p_X-p_Y&p_Y-p_Z&p_X-p_Z& (p_X,p_Y,p_Z)_{t=1}\\ \hline
\gamma_{12}&re&(R+s-re)/2&(R+s+re)/2&(s,R,0)\\
\gamma_{23}&R-s(1+e)/2&se&R-s(1-e)/2&(R,0,s)\\
\gamma_{X\mid YZ}&Re-s&s&Re&(0,R+s,R)\\
\gamma_{XY\mid Z}&r&se&r+se&(r,0,s).
\end{array}
\end{equation}
For $\gamma_{12}$, $|(R+s\pm re)/2|\ge (R+s-r)/2=s>r$;
for $\gamma_{23}$, $|R-s(1\pm e)/2|\ge R-s>s$.
The block paths have the respective inequalities
$|Re|=R>s$ and $|se|=s>r$, with remaining
nonzero differences bounded below by $R-s$ and $s-r$.
Thus all four curves remain in their \emph{specified original}
iterate/product chambers and avoid every diagonal, including
at both endpoints.  These bounds also hold uniformly on a
small complex neighborhood of each compact image.  For fixed
ordinary-module restricted-dual tests, the applicable original
one-edge intermediate-weight series are locally normally convergent
there by Corollary~\ref{cor:smallconvergenceautomatic}; this statement
is about the original scalar series, not about arbitrary reordering
of translation and weight sums.  Along
$\gamma_{12}$ and $\gamma_{23}$, respectively, one has
\begin{equation}\label{eq:adjacentbounds98}
 \begin{aligned}
 |p_X-p_Y|&=r<s\le |p_Y-p_Z|,\\
 |p_Y-p_Z|&=s<R-s\le |p_X-p_Z|.
 \end{aligned}
\end{equation}
Indeed the two unrotated anchors are circles of radii $r/2$ and $s/2$.  The inequalities are strict in the respective iterate and product chambers.  Both paths therefore admit normally convergent original
three-input contractions throughout their images, by Corollary~\ref{cor:smallconvergenceautomatic}. At $t=1$ their ordered positions are $(s,R,0)$ and $(R,0,s)$.  Apply the fixed-output-weight identity
\eqref{eq:skewinverse}; the remaining translations are
\eqref{eq:skewtranslations96} and
\eqref{eq:skewanchored96}.  In particular their output
factors are the actual finite restricted-dual transforms $e^{pL_{D'}(1)}$, not discarded scalar factors.  Their endpoint canonical types are $(T_gY:R,X:s,Z:0),\; (X:R,T_hZ:s,Y:0)$. The first $T_gY$ and second $T_hZ$ agree with the
pullback types in \eqref{eq:Pzrawcrossed}; the sign of each half-turn is the $-\pi i$ continuation in
\eqref{eq:skewinverse}, not an arbitrary choice of
$\arg(-r)$.

To see exactly how the adjacent endpoint maps are obtained, let
$\mathcal Z_{X,Y}=\mathcal Y^{\can}_{T_gY,X}$ and fix every
output weight $\eta$ before making a continuation.  The raw crossed
map and the independently checked inverse-skew identity give
\[
\operatorname{ev}_{\log r-\pi i}
 (\overline c_{X,Y}\mathcal Y^{\can}_{X,Y})(x\otimes y)
 =e^{-rL(-1)}\mathcal Z_{X,Y}(y,r)x.
\]
The left side is continued along the actual $r e_t$ separating
coordinate of $\gamma_{12}$.  Its right-hand side is compared with
the positive-real terminal operator only after the translated
output has been paired with the finite dual gauge
$e^{pL_{D'}(1)}$.  The identity in
Lemma~\ref{lem:literalskew}, equation~\eqref{eq:skewanchored96},
is precisely the relevant finite-$\eta$ comparison.  For
$\gamma_{23}$ apply the same inverse-skew calculation to the
inner pair $(Y,Z)$, with the outer original operator held fixed,
and perform the Taylor sum before its intermediate-weight sum.
This yields the typed output $(X,T_hZ,Y)$, not $(X,Z,Y)$.
For the first exchange the output is $(T_gY,X,Z)$, and the
relevant maps are exactly the pre-existing $c_{X,Y}$ and
$c_{Y,Z}$ (tensored with the spectator and transported by the
pre-existing associator as required).  The comparison is
between their \emph{actual matrix coefficients}, not merely a
permutation of the input labels.

For the block paths, the \emph{inner} relative coordinate has constant modulus, while the \emph{root} coordinate has constant modulus:
\begin{equation}\label{eq:blockbounds98}
 \begin{array}{c|c|c|c}
   &\text{inner displacement}&\text{root displacement}
   &\text{other differences}\\ \hline
 \gamma_{X\mid YZ}&p_Y-p_Z=s&p_X-p_Z=Re_t
     &|p_X-p_Y|=|Re_t-s|\ge R-s>0\\
 \gamma_{XY\mid Z}&p_X-p_Y=r&p_Y-p_Z=se_t
     &|p_X-p_Z|=|r+se_t|\ge s-r>0.
 \end{array}
\end{equation}
Thus $s<R$ and $r<s$ are exactly the iterate-chamber inequalities for the root/block products along these paths; no collision or loss of normal convergence occurs. Their terminal configurations are $(X:0,Y:R+s,Z:R),\; (X:r,Y:0,Z:s)$.
The first root half-turn uses the twist $g$ of $X$; the second uses the twist $gh$ of $X\bo Y$.  For the first one, with $M=Y\bo Z$, the finite-weight form of
\eqref{eq:skewinverse}, followed by the translation of
Lemma~\ref{lem:fullframes}, cancels $e^{-RL(-1)}$ against the terminal restricted-dual gauge $U_R=e^{RL_{D'}(1)}$.  The naturality defining the actual tensorator gives
\begin{equation}\label{eq:blocktensor98}
 \overline{(J^g_{Y,Z})^{-1}}\,
 T_g\!\left(\mathcal Y^{\can}_{Y,Z}(y,s)z\right)
 =\mathcal Y^{\can}_{T_gY,T_gZ}(y,s)z,
\end{equation}
interpreted after each internal weight projection.  The same calculation on a homogeneous coefficient
$\pi^M_\mu\mathcal Y^{\can}_{Y,Z}(y,s)z$ first transports it as
a \emph{vector of $T_gM$}; since $J^g$ is an $L(0)$-preserving
module map, its inverse then acts on that $\mu$-component, giving
$\pi^{T_gY\bo T_gZ}_\mu\mathcal Y^{\can}_{T_gY,T_gZ}(y,s)z$.
No infinite completed input is passed to $(J^g)^{-1}$ prior to its
weightwise extension.  The result is the iterate with terminal
positions $(T_gY:R+s,T_gZ:R,X:0)$; Theorem~\ref{thm:smallassoc}
identifies its original sum with the product, since
$R+s>R>s>0$.  At the other block endpoint the corresponding
finite projection is taken in $X\bo Y$, whose meridian is
$gh$ in the stipulated ordered product.  Consequently the
first entry after exchange is $T_{gh}Z$.  In the overlap
$s>r>s-r>0$ the existing inverse associator maps the terminal
product to the terminal iterate with no change of rooted branch.
For the second block, inverse skew and $U_s$ similarly give the product with positions $(T_{gh}Z:s,X:r,Y:0)$.  In the overlap
\[
 \frac{s}{2}<r<s,\qquad s>r>s-r>0,
\]
Theorem~\ref{thm:smallassoc} supplies its inverse associator and identifies the corresponding iterate. In both cases the cited calculations are on the actual finite projections; the remaining intermediate-grade sums
are locally normally convergent along the compact paths by Corollary~\ref{cor:smallconvergenceautomatic}. No equality of hexagon composites is used here.

Finally, the nonconstant differences of the first block whose real parts cross zero are $Re_t-s$ and $Re_t$; their zeros in real projection occur successively at $\cos(\pi t)=s/R$ and $\cos(\pi t)=0$.  For the second block they are
$se_t$ and $r+se_t$, with crossings successively at $\cos(\pi t)=0$ and $\cos(\pi t)=-r/s$. At the crossing times their imaginary parts are strictly negative.  Moreover \eqref{eq:blockbounds98} excludes
all collisions, and the inner strands never cross. A generic small perturbation gives the same positive braid diagrams, hence the stated chronological words.  Their
Hurwitz action is exactly \eqref{eq:meridianhurwitz}, for the same input-ordered stems used in Lemma~\ref{lem:canonicalnodecurrent}. In particular $(aba^{-1})a=ab$, so the output meridian retains the inverse of the total input product. The complete triple grades are carried by
\[
 \begin{aligned}
 (g,h,k)&\longmapsto(ghg^{-1},g,k)
       \longmapsto(ghg^{-1},gkg^{-1},g),\\
 (g,h,k)&\longmapsto(g,hkh^{-1},h)
       \longmapsto(ghkh^{-1}g^{-1},g,h).
 \end{aligned}
\]
These are the pullback labels of the four canonical commutors.  The analytic root lifts are obtained by lifting each path from its prescribed initial frame in the finite marked-root atlas, as in Lemma~\ref{lem:fullframes};
their terminal roots are compared with the deck transformations.  The skew and inverse-skew operations and their branch restrictions agree with
\cite[Section~3, Theorem~3.1 and its following corollary]{DH}.

We verify the \emph{endpoint-fixed} homotopies of the two \emph{complete}
hexagon routes, including their inserted positive-real motions.  Keep the
original names $(X,Y,Z)$ for the three moving markings, and use the global
coordinates
\begin{equation}\label{eq:hexconfigurationcoordinates}
 p=p_Z,\qquad w=p_Y-p_Z,\qquad
 q=\frac{p_X-p_Z}{p_Y-p_Z},\qquad
 (p_X,p_Y,p_Z)=(p+qw,p+w,p).
\end{equation}
They identify $F_3(\C)$ with
$\C_p\times\C_w^\times\times(\C\setminus\{0,1\})_q$;
in particular all three collision differences are $w,qw,(q-1)w$.

For the first hexagon take $S<R<2S$ and write $r=R-S$.
The block path starts at $(R,S,0)$ and ends at $(0,R+S,R)$;
its coordinates are $w=S$ and $q=(R/S)e_t$.  The adjacent
route consists, in temporal order, of $\gamma_{12}$, the
positive movement $p_Y:R\to R+S$, the exchange of the then
adjacent $X,Z$ at fixed $p_Y=R+S$, and the positive movement
$p_Z:S\to R$.  On its two rotating pieces, respectively,
\[
 q=\frac{a+be_t}{a-be_t},\quad
 a=\frac{R+S}{2},\ b=\frac{R-S}{2};
 \qquad
 q=\frac{S e_t}{R+S/2+(S/2)e_t}.
\]
Both have $\Im q<0$ for $0<t<1$, while the inserted
real motions avoid $q=0,1$.  In addition $\Re w>0$ on
both full routes.  Hence, after an arbitrarily small
negative imaginary perturbation of $q$ on their interior
real segments, both routes lie in the convex product
$\{\Re w>0\}\times\{\Im q<0\}$ at all interior
times, with their original endpoints fixed.

For the second hexagon take $S/2<r<S$.  Its block route
starts at $(S+r,S,0)$ and ends at $(r,0,S)$, with
$w=S e_t$ and $q=1+(r/S)e^{\pi it}$.  The adjacent route
moves $p_X:S+r\to3S$, exchanges $Y,Z$, moves
$p_X:3S\to S+r$, exchanges the currently adjacent
$X,Z$, and finally moves $(p_X,p_Z):(S,S+r)\to(r,S)$
within $0<p_X<p_Z$ (always keeping $p_Y=0$ on this last
piece).  On its two rotating pieces the coordinates are,
respectively,
\[
 (w,q)=\left(S e_t,\frac12+\frac52e^{\pi it}\right),
 \qquad
 (w,q)=\left(-\frac{2S+r-r e_t}{2},
              -\frac{2r e_t}{2S+r-r e_t}\right).
\]
Here $\Im w<0$ and $\Im q>0$ on both open rotating
arcs; the inserted real motions avoid $w=0$ and $q=0,1$.
The block route has the same inequalities.  Perturbing the
interior real segments of $w$ slightly downwards and those
of $q$ slightly upwards places both routes in
$\{\Im w<0\}\times\{\Im q>0\}$, again leaving
endpoints unchanged.

For reference, the positive-real pieces omitted from the rotating-arc
inequalities have the following exact real restrictions.  In the first
route, after the $XY$ exchange the motion of $Y$ from $R$ to
$R+S$ has $p=0$, $w=p_Y>0$, and $q=S/p_Y\in(0,1)$;
following the $XZ$ exchange, the motion of $Z$ from $S$ to $R$
has $p=p_Z$, $w=R+S-p_Z>0$, and
$q=-p_Z/w<0$.  In the second route, $X:S+r\to3S$
gives $w=S$, $q=p_X/S>1$; after the $YZ$ exchange,
$X:3S\to S+r$ gives $w=-S$ and $q=1-p_X/S<0$;
the terminal motion with $0<p_X<p_Z$ gives
$w=-p_Z<0$ and $q=1-p_X/p_Z\in(0,1)$.
Thus every real piece avoids the deleted loci, and each
enters the indicated half-plane by the specified small perturbation.
In each closed subarc on which a difference has nonzero imaginary
part, its sign is strict; at a joining point on the real axis its
real value is separated from $0$ and $1$, where applicable.
Compactness of the full finite union of subarcs therefore permits a
single $\varepsilon>0$ for which the entire straight perturbation
is collision-free, even arbitrarily close to its fixed endpoints.

To make the homotopies explicit, reparametrize each complete
route on $[0,1]$, and in the first case set
$\widetilde q(t)=q(t)-i\varepsilon t(1-t)$; in the second
set $\widetilde w(t)=w(t)-i\varepsilon t(1-t)$ and
$\widetilde q(t)=q(t)+i\varepsilon t(1-t)$.
For sufficiently small $\varepsilon>0$ these perturbations,
and their straight-line removal, avoid all deleted loci.
If subscripts $\mathrm{blk},\mathrm{adj}$ denote the two
perturbed routes, interpolate their three coordinates by
\begin{equation}\label{eq:hexexplicitisotopy}
 (p_\lambda,w_\lambda,q_\lambda)(t)
 =(1-\lambda)(p,w,q)_\mathrm{blk}(t)
       +\lambda(p,w,q)_\mathrm{adj}(t),
 \qquad 0\le\lambda\le1.
\end{equation}
Convexity keeps $w_\lambda\ne0$ and
$q_\lambda\notin\{0,1\}$ throughout; applying the inverse
of \eqref{eq:hexconfigurationcoordinates} gives the required
endpoint-fixed homotopies in the \emph{full} $F_3(\C)$, rather
than merely equality of braid diagrams.  For the differences
$(p_X-p_Y,p_X-p_Z,p_Y-p_Z)$ the continuous logarithm
increments from the initially positive-real branches are,
respectively,
\[
 (-\pi i,-\pi i,0)\quad\text{and}\quad
 (0,-\pi i,-\pi i).
\]
Indeed, on the \emph{block} curves the ordered separating
coordinates have the factorizations
\[
 (p_X-p_Y,p_X-p_Z,p_Y-p_Z)=
 \begin{cases}
 (Re_t-S,Re_t,S), &X\mid YZ,\\
 (r,r+Se_t,Se_t),&XY\mid Z.
 \end{cases}
\]
The nonconstant entries move strictly through the lower half-plane
from positive to negative real values when appropriate; the constant
positive entries acquire no argument.  Thus the displayed logarithm
increments are actual changes of continuously chosen logarithms,
not a deduction from the group-valued Hurwitz action.  The explicit
homotopy gives the same increments on the adjacent routes.
The moving third marking is retained in $p$; the fixed
output point at infinity is accounted for by the
contragredient gauge of Lemma~\ref{lem:fullframes}.
Starting both routes at the same complete framed
point, lift these fixed-endpoint homotopies in the finite
marked-root Hurwitz atlas of that lemma and in the continuous
logarithm covering.  This identifies their terminal typed
based lifts, marking lifts, special roots, deck transports,
and logarithms, not merely the terminal unframed configurations.
It does \emph{not} assert trivial monodromy for an arbitrary
pure braid.
\end{proof}

For a concrete check of the same convention, use right-to-left composition of permutations and take
\[
 g=(12),\qquad h=(23),\qquad k=(123)\quad\text{in }S_3.
\]
Then $gh=(123)$, $hk=(13)$ and $ghk=(132)$. The product and iterate nodes of Lemma~\ref{lem:canonicalnodecurrent} therefore have opposite labels $((13),(13))$ and $((123),(132))$, respectively; all four component triples have product one. The two chronological block exchanges above give
\[
\begin{aligned}
 ((12),(23),(123))
 &\xrightarrow{\sigma_1}((13),(12),(123))
 \xrightarrow{\sigma_2}((13),(132),(12)),\\
 ((12),(23),(123))
 &\xrightarrow{\sigma_2}((12),(132),(23))
 \xrightarrow{\sigma_1}((123),(12),(23)).
\end{aligned}
\]
Both terminal products are $(132)$, and the output meridian stays $(123)$. The labels change by conjugation, with the corresponding pullbacks $T_gT_h=T_{gh}$; the root phases are still those obtained by lifting the actual paths, not determined by this group calculation alone.

\begin{lemma}\label{lem:sourcecontinuation}
Let $U$ be either of the canonical three-input bracketings of $X,Y,Z$, let $N$ be an ordinary output twisted module, and let $F_1,F_2\in\Hom_V(U,N)$.  Set
\begin{equation}\label{eq:sourcefunctional99}
 \Lambda_i(d';x,y,z;\boldsymbol p)
 =\left\langle d',\overline{F_i}\,
   \mathscr Y_U(x,y,z;\boldsymbol p)\right\rangle,
 \quad d'\in N',\quad i=1,2,
\end{equation}
where $\mathscr Y_U$ is the original canonical product or iterate with its specified full initial frame.  Suppose two paths $\gamma_1,\gamma_2$ in the collision-free three-point configuration space start at that frame,
have the same full framed endpoint up to its specified root/deck gauge, and are endpoint-fixed homotopic after lifting from the same framed starting point.  Assume that the \emph{actual normally convergent} continuations of
\eqref{eq:sourcefunctional99}, including the endpoint
gauges, have identical scalar germs in two independent terminal configuration variables for every $x\otimes y\otimes z$ and $d'\in N'$.  Then $F_1=F_2$.
\end{lemma}
\begin{proof}
Use the canonical differential Ward quotient $\cQ$ of Definition~\ref{def:diff}, with ordered insertion modules $X,Y,Z,N'$ and the fixed original physical/root frames. It is finite-rank locally free, with flat connection, by
Theorem~\ref{thm:finiteconn}; its actual moving-frame identification is Lemma~\ref{lem:fullframes}.  For each $d'\in N'$, the transpose $F_i'd'\in U'$ is again
in the \emph{restricted} dual, because $F_i$ commutes with $L(0)$.  Postcomposition of the canonical operators by $F_i$ obeys the complete twisted-intertwining-operator axioms.  Theorem~\ref{thm:E0covariantcontinuation}, applied to these actual input/output modules and to the transported normally convergent original germs, implies $\Lambda_i(\cJ_{\diff})=0,\; \partial_\delta\Lambda_i(w) =\Lambda_i(\nabla_\delta w)$. Thus $\Lambda_1$ and $\Lambda_2$ are horizontal sections of the same dual $\cQ^\vee$ at the initial frame.
The three input positions may be unanchored by Lemma~\ref{lem:fullframes}; its output gauge and the root/deck gauges intertwine the two Ward connections.

Write $\mathsf P_{\gamma_i}$ for the induced invertible parallel transports on the framed dual $\cQ^\vee$. The homotopy-lifting property and flatness give
\begin{equation}\label{eq:sourcehomotopy99}
 \mathsf P_{\gamma_1}=\mathsf P_{\gamma_2}
 \quad\text{after the prescribed common terminal gauge}.
\end{equation}
This is an equality of transport maps between the same framed fibers, not an assertion of trivial monodromy. By the hypothesis on the terminal germs, the evaluations of $\mathsf P_{\gamma_1}\Lambda_1$ and
$\mathsf P_{\gamma_2}\Lambda_2$ on all algebraic external tensors coincide; hence these two terminal sections of $\cQ^\vee$ coincide.  Equation~\eqref{eq:sourcehomotopy99} and invertibility of parallel transport imply
$\Lambda_1=\Lambda_2$ on the original overlap.

This equality of horizontal functionals detects
$F_1-F_2$, rather than only its value on one completed vector.
It is important that endpoint agreement is a \emph{separate
hypothesis}: equality of the two transport maps alone does not
imply equality of two initially different horizontal sections.
Conversely, agreement of the terminal scalar germs for all
algebraic external tests identifies the two sections of the
finite-rank terminal Ward dual, since those evaluations separate
its underlying quotient.  Only then does invertible transport
return equality to the initial framed fiber.  Put $H=F_1-F_2$.  On any nested convergence chamber, the scalar function
\[
\big\langle d',\overline H\,
   \mathscr Y_U(x,y,z;\boldsymbol p)\big\rangle=0
\]
is an original normally convergent generalized Laurent expansion.  By uniqueness of its coefficients, for a source $U=(X\bo Y)\bo Z$ one obtains
\begin{equation}\label{eq:sourcecoeff99}
 \left\langle d',H\,
 \pi^U_\eta\mathcal Y^{\can}_{X\bo Y,Z}
 \left(\pi^{X\bo Y}_\mu
       \mathcal Y^{\can}_{X,Y}(x,1)y,1\right)z\right\rangle=0
 \quad(\mu,\eta\in\C).
\end{equation}
For the iterated bracketing $U=(X\bo Y)\bo Z$, let
$V=X\bo Y$, write $x\in X_\alpha$, $y\in Y_\beta$,
$z\in Z_\gamma$, and take the anchored positions
$(R(1+q),R,0)$ with $R>0$, $0<q<1$.
For homogeneous $d'\in N'_\eta$, the original scalar
contraction of $H$ has the normally convergent expansion
\[
 R^{\eta-\alpha-\beta-\gamma}
 \sum_{\mu\in\wt(V)}q^{\mu-\alpha-\beta}
 \left\langle d',H\pi^U_\eta
  \mathcal Y^{\can}_{V,Z}
  (\pi^V_\mu\mathcal Y^{\can}_{X,Y}(x,1)y,1)z\right\rangle.
\]
Distinct \emph{numerical} weights $\mu$ give distinct powers
of $q$ even when two irreducible summands have the same lowest
weight.  A finite set of exponent cosets modulo $\mathbb Z$
can be separated by continued logarithms and then by ordinary
Laurent uniqueness on each coset.  Since the scalar is zero on a
nonempty open set, the coefficient of every power vanishes.
For the product bracketing $U=X\bo(Y\bo Z)$ set $W=Y\bo Z$,
take the anchored positions $(R,R\rho,0)$ with $R>0$,
$0<\rho<1$, and obtain instead
\[
 R^{\eta-\alpha-\beta-\gamma}
 \sum_{\mu\in\wt(W)}\rho^{\mu-\beta-\gamma}
 \left\langle d',H\pi^U_\eta
 \mathcal Y^{\can}_{X,W}
 (x,1)\bigl(\pi^W_\mu
 \mathcal Y^{\can}_{Y,Z}(y,1)z\bigr)
 \right\rangle.
\]
The same independent exponent separation applies.  Thus in
both bracketings we recover every actual doubly projected
canonical coefficient, not merely an equality of completed vectors
at one numerical placement.  In each formula the projections
create genuine algebraic vectors, and the map $H$ preserves the
outer weight $\eta$.
 Only \emph{numerical} weights are separated: contributions from distinct simple summands of equal numerical weight
are kept together in $\pi_\mu$.  Lemma~\ref{lem:generation} shows that these doubly projected vectors span the relevant source $U$.  Since $N'$ separates $N$, \eqref{eq:sourcecoeff99} (and its product counterpart) gives $H=0$.
All coefficient extractions follow from the normal convergence supplied by Corollary~\ref{cor:smallconvergenceautomatic}; no formal series is identified with an analytic germ without this step.
\end{proof}

\begin{theorem}\label{thm:criticalhexagons}
Let $X\in\cC_g$, $Y\in\cC_h$, and $Z\in\cC_k$.  With $T_gT_h=T_{gh}$, the fixed maps $a,c,J$ satisfy both
\eqref{eq:hex1} and \eqref{eq:hex2}.  Their common sources and
codomains, respectively, are
\[
 \begin{aligned}
 L_1&=(X\bo Y)\bo Z,&
 N_1&=T_gY\bo(T_gZ\bo X),\\
 L_2&=X\bo(Y\bo Z),&
 N_2&=(T_{gh}Z\bo X)\bo Y.
 \end{aligned}
\]
The commutors are the fixed $\Omega_+$ maps of
\eqref{eq:criticalrawc}; their inverse endpoint calculations
use the $-\pi i$ continuation in \eqref{eq:skewinverse}. No coherence identity is assumed in constructing these maps.
\end{theorem}
\begin{proof}
We give the object types, the two independent endpoint calculations,
and the passage from terminal germs to maps separately.  In particular,
path homotopy will not be used to assert equality of two initially
different horizontal sections.

For the first hexagon the adjacent composite has the chain
\[
\begin{split}
 (X\bo Y)\bo Z&\xrightarrow{c_{X,Y}\bo\id_Z}
 (T_gY\bo X)\bo Z\\
 &\xrightarrow{a_{T_gY,X,Z}}T_gY\bo(X\bo Z)
 \xrightarrow{\id_{T_gY}\bo c_{X,Z}}
 T_gY\bo(T_gZ\bo X),
\end{split}
\]
whereas the block composite has the chain
\[
\begin{split}
 (X\bo Y)\bo Z&\xrightarrow{a_{X,Y,Z}}X\bo(Y\bo Z)
 \xrightarrow{c_{X,Y\bo Z}}T_g(Y\bo Z)\bo X\\
 &\xrightarrow{(J^g_{Y,Z})^{-1}\bo\id_X}
 (T_gY\bo T_gZ)\bo X
 \xrightarrow{a_{T_gY,T_gZ,X}}T_gY\bo(T_gZ\bo X).
\end{split}
\]
For the second hexagon the adjacent chain is
\[
\begin{split}
 X\bo(Y\bo Z)&\xrightarrow{\id_X\bo c_{Y,Z}}
 X\bo(T_hZ\bo Y)\xrightarrow{a^{-1}_{X,T_hZ,Y}}
 (X\bo T_hZ)\bo Y\\
 &\xrightarrow{c_{X,T_hZ}\bo\id_Y}
 (T_gT_hZ\bo X)\bo Y,
\end{split}
\]
and the block chain is
\[
 X\bo(Y\bo Z)\xrightarrow{a^{-1}_{X,Y,Z}}(X\bo Y)\bo Z
 \xrightarrow{c_{X\bo Y,Z}}T_{gh}Z\bo(X\bo Y)
 \xrightarrow{a^{-1}_{T_{gh}Z,X,Y}}(T_{gh}Z\bo X)\bo Y.
\]
All arrows are $V$-module maps in the sector $ghk$, since
$T_gY\in\cC_{ghg^{-1}}$, $T_gZ\in\cC_{gkg^{-1}}$ and
$T_gT_hZ=T_{gh}Z\in\cC_{ghkh^{-1}g^{-1}}$.
The indicated inverse tensorator, not $J^g_{Y,Z}$, has the
necessary source and target.  The maps $a$, $c$ and $J$ have
already been constructed by Theorem~\ref{thm:smallassoc},
\eqref{eq:criticalrawc} and Proposition~\ref{prop:action};
no hexagon is used in their definitions.  Fix algebraic
$x\in X$, $y\in Y$, $z\in Z$ and $d'$ in the restricted dual
of the displayed common target.  By finite linearity we may first
assume these vectors homogeneous.

Take independent $R,S>0$ with $S<R<2S$ and put $r=R-S$.
The initial labeled positions are $(X:R,Y:S,Z:0)$.
The block path $\gamma_{X\mid YZ}$ of \eqref{eq:blockX98}
ends at $(X:0,Y:R+S,Z:R)$, keeping the $YZ$ relative
distance $S$ and the root/block distance $R$.  For an
internal weight $\mu$ of $M=Y\bo Z$ put
\[
 v_\mu=\pi_\mu^M\mathcal Y^{\can}_{Y,Z}(y,S)z\in M_\mu.
\]
This is an \emph{algebraic} vector.  The inverse-skew formula
\eqref{eq:skewinverse} on this coefficient gives, in the prescribed
terminal frame, the contraction
\[
 \big\langle U_Rd',e^{-RL(-1)}
 \mathcal Y^{\can}_{T_gM,X}(T_gv_\mu,R)x\big\rangle,
 \qquad U_R=e^{RL_{N_1'}(1)}.
\]
The complete output translation is interpreted after a fixed
output-weight projection, as in Lemma~\ref{lem:literalskew}.
Restricted contragredience gives
$\langle U_Rd',e^{-RL(-1)}w\rangle=\langle d',w\rangle$.
For every $\mu$, the defining tensorator identity
\eqref{eq:blocktensor98} is
\[
 \overline{(J^g_{Y,Z})^{-1}}T_gv_\mu
 =\pi_\mu^{T_gY\bo T_gZ}
 \mathcal Y^{\can}_{T_gY,T_gZ}(y,S)z.
\]
Hence, after applying the specified inverse tensorator, the
root contraction is exactly the original iterate
of the pair of complete operators with inner displacement $S$
and outer displacement $R$.  Since $R>S>0$, the original
product--iterate identity of Theorem~\ref{thm:smallassoc},
with its positive-real logarithms and original restricted-dual
weight sums, identifies its analytic germ with
\begin{equation}\label{eq:hexterminal1}
 \mathcal F_1(R,S)
 =\left\langle d',
 \mathcal Y^{\can}_{T_gY,T_gZ\bo X}(y,R+S)
 \mathcal Y^{\can}_{T_gZ,X}(z,R)x\right\rangle.
\end{equation}\label{eq:hexfullgerm1}
The displayed operator pair is its original intermediate-weight
sum, not an operator evaluated on an arbitrary completed input.

The first lower-half exchange $\gamma_{12}$ ends at
$(T_gY:R,X:S,Z:0)$.  The relevant inner gap is $r=R-S$,
which satisfies $0<r<S$; hence the needed original
three-input associativity comparison is in its strict overlap.
Move the $T_gY$ marking from $R$ to $R+S$ while keeping
$X:S$, $Z:0$.  As $R+S>2S$, the second lower-half exchange
of $(X,Z)$ is an admissible instance of $\gamma_{23}$,
with the spectator $T_gY$ at $R+S$.  It ends at
$(T_gY:R+S,T_gZ:S,X:0)$.
Precisely here the output-translation comparison is, for a
fixed internal coefficient and then its normally convergent
sum,
\begin{equation}\label{eq:hexadjacenttranslation}
\begin{split}
 &\left\langle U_Sd',
 \mathcal Y^{\can}_{T_gY,T_gZ\bo X}(y,R)
 e^{-SL(-1)}\mathcal Y^{\can}_{T_gZ,X}(z,S)x\right\rangle\\
 &\qquad=\left\langle d',
 \mathcal Y^{\can}_{T_gY,T_gZ\bo X}(y,R+S)
 \mathcal Y^{\can}_{T_gZ,X}(z,S)x\right\rangle.
\end{split}
\end{equation}
The equality is \eqref{eq:skewanchored96} with its outer
parameter $R+S$ and shift $S$; the source in the left-hand
operator is therefore $R$, not an independently chosen
numerical value.  This choice matches the same
positive-real motions used in Proposition~\ref{prop:braidendpoints}.
Finally move the $T_gZ$ marking from $S$ to $R$, keeping
$T_gY:R+S$ and $X:0$ fixed.  Throughout,
$R+S>R\ge S>0$; this is the original product chamber,
and the terminal evaluation is exactly \eqref{eq:hexterminal1}.
Thus both first-hexagon routes give the same actual
terminal scalar analytic function for all $S<R<2S$, not
only for $R=3S/2$.

Here restrict further to $S/2<r=R-S<S$.
The initial positions are $(X:R,Y:S,Z:0)$, and the inverse
associator changes the canonical product bracketing to
$(X\bo Y)\bo Z$ on the strict overlap $S>r>S-r>0$.
The block path $\gamma_{XY\mid Z}$ keeps $X-Y=r$ and
ends at $(X:r,Y:0,Z:S)$.
Its skew transformation is by the ordered $gh$-twist;
its terminal product, including the inverse-skew phase and
restricted-dual output gauge, is
\[
 \left\langle d',
 \mathcal Y^{\can}_{T_{gh}Z,X\bo Y}(z,S)
 \mathcal Y^{\can}_{X,Y}(x,r)y\right\rangle.
\]
At each internal weight one uses the actual $T_{gh}Z$-typed
opposite field and then the output-frame cancellation;
no completed inner output is treated as algebraic input.
Applying the original inverse associator in
$S>r>S-r>0$ gives the common terminal iterate
\begin{equation}\label{eq:hexterminal2}
 \mathcal F_2(R,S)
 =\left\langle d',
 \mathcal Y^{\can}_{T_{gh}Z\bo X,Y}
 \left(\mathcal Y^{\can}_{T_{gh}Z,X}(z,S-r)x,r\right)y
 \right\rangle.
\end{equation}\label{eq:hexfullgerm2}

Start with $X:R,Y:S,Z:0$; move $X$ to $3S$ (the product
condition $3S>S$ remains strict), exchange $(Y,Z)$ along
$\gamma_{23}$, and move $X$ back to $R$ with
$(T_hZ:S,Y:0)$ fixed.  The inverse associator is now
applied at $(X:R,T_hZ:S,Y:0)$, where $0<r=R-S<S$.
The adjacent exchange of $(X,T_hZ)$ is an instance of
$\gamma_{12}$ with separation $r$ and external marking
$Y:0$; its terminal ordered labels are
$(T_gT_hZ:R,X:S,Y:0)$.  After the complete
inverse-skew and translation comparison, the original
terminal iterate at \emph{these} positions is
\[
 \left\langle d',
 \mathcal Y^{\can}_{T_{gh}Z\bo X,Y}
 \left(\mathcal Y^{\can}_{T_{gh}Z,X}(z,r)x,S\right)y
 \right\rangle.
\]
This is the actual endpoint before the last positive-real
motion.  Move the two positive positions along
\[
 b(\tau)=(1-\tau)R+\tau S,\qquad
 d(\tau)=(1-\tau)S+\tau r,\qquad 0\le\tau\le1.
\]
The inner displacement $b(\tau)-d(\tau)$ remains strictly
between $0$ and $d(\tau)$: both inequalities are affine in
$\tau$, and at the endpoints they read $0<r<S$ and
$0<S-r<r$ (the latter uses $r>S/2$).
Consequently the original \emph{iterate} remains locally
normally convergent, and its final value is
\eqref{eq:hexterminal2}.  The pullback is the strictly
ordered functor $T_gT_h=T_{gh}$, not $T_{hg}$.

These four calculations are statements about original complete
intertwining operators.  More precisely, fix a restricted-dual
vector, an internal numerical weight $\mu$, and an algebraic
input before applying the skew or tensorator formulas.
The projector $\pi_\mu$ makes the intermediate argument
algebraic, $J$ and every associator preserve its weights,
and $U_p=e^{pL_{N_i'}(1)}$ is finite on the restricted dual.
Translations of completed outputs are taken only after fixed
output-weight projection.  For a translated \emph{input},
we use the ordered Taylor interpretation of
Lemma~\ref{lem:literalskew}, not an undefined substitution
of a vector of $\overline M$ into another operator.

The resulting scalar identity is also an identity of
\emph{normally convergent original series}, rather than a
formal rearrangement.  Indeed, the compact real path pieces
above, and sufficiently small complex neighborhoods of them,
stay within the explicit strict original convergence margins
of \eqref{eq:adjacentbounds98}, \eqref{eq:blockbounds98}
and the positive-real product--iterate overlaps.
Corollary~\ref{cor:smallconvergenceautomatic} supplies normal
convergence of the original intermediate-weight sums.
On a smaller complex parameter polydisc, Cauchy's formula
allows derivatives of those sums to be computed termwise;
the complete IO $L(-1)$-identities and the full-frame
gauges identify the resulting \emph{analytic germs}.
Their Taylor continuation is thus justified at the level of
the normally convergent function and need not be interchanged
with an unrelated infinite intermediate-weight sum.  Where
two IO channels have different intermediate weight indices,
Theorem~\ref{thm:smallassoc} compares the two original
normally convergent sums as analytic functions; it does not
assert a coefficient-by-coefficient bijection of their indices.

For the first route the independent positive variables $(R,S)$
range over $\{S<R<2S\}$.  For the second they range over
$\{3S/2<R<2S\}$, where $r=R-S$ and $S/2<r<S$.
The two terminal coordinate maps are
\[
 (R,S)\longmapsto(R+S,R),\qquad
 (R,S)\longmapsto(S,r),
\]
with respective Jacobian determinants $-1$ and $-1$.
Their images contain nonempty \emph{real two-dimensional}
open sets inside the corresponding complex product or iterate
convergence chambers.  Both routes of each hexagon therefore
give one and the same \emph{specified rooted} scalar analytic
function \eqref{eq:hexterminal1} or \eqref{eq:hexterminal2}
on a real open set.  Local normal convergence gives complex
holomorphic germs about every such point; applying the
one-variable identity theorem successively in the two
independent complex parameters upgrades this to equality
of the full terminal two-variable germs.  The argument
uses the same continuous logarithms and rooted output gauge
on each route, with no assertion of trivial monodromy along
an unrelated path.

Proposition~\ref{prop:braidendpoints} proves that the block
and adjacent routes with the same initial frame are
endpoint-fixed homotopic in $F_3(\C)$, including all of the
positive-real pieces just used.  The lifted homotopy in the
full marked-root atlas of Lemma~\ref{lem:fullframes}
identifies their terminal roots, based and marked lifts,
noncommutative deck transports and continuous logarithms.
The three separating-coordinate increments are, in order,
$(-\pi i,-\pi i,0)$ for the first hexagon and
$(0,-\pi i,-\pi i)$ for the second; they cannot be suppressed.
Thus the equality of terminal analytic germs above holds
\emph{after the same specified terminal gauge}, not merely
on the same unframed positions.

The composites in the two type diagrams at the beginning
are actual $V$-module maps
$H_i^{\mathrm{adj}},H_i^{\mathrm{blk}}:L_i\to N_i$.
Postcomposing the fixed original canonical three-input
channel of $L_i$ with either map gives a horizontal
restricted-dual Ward functional in the same
finite-rank framed connection, by
Theorem~\ref{thm:E0covariantcontinuation} and
Lemma~\ref{lem:fullframes}.  Homotopy identifies the two
transport maps, but does \emph{not} by itself identify
the two horizontal functionals; the latter identification
is precisely the independent terminal-germ calculation above.
Lemma~\ref{lem:sourcecontinuation} now applies, once to
$(L_1,N_1)$ and once to $(L_2,N_2)$.
Its proof recovers every doubly projected original canonical
coefficient from generalized-Laurent uniqueness and uses
Lemma~\ref{lem:generation} for their \emph{algebraic}
span, together with restricted-dual separation.  It yields
\[
 H_1^{\mathrm{adj}}=H_1^{\mathrm{blk}},\qquad
 H_2^{\mathrm{adj}}=H_2^{\mathrm{blk}}.
\]
These are \eqref{eq:hex1} and \eqref{eq:hex2}.  Naturality
and finite additivity of the preceding canonical operators,
associators, commutors and tensorators extend both identities
to finite semisimple direct sums and then to arbitrary objects
of $\cC_G$.
\end{proof}

\needspace{6\baselineskip}
\section{Rigidity and the main theorem}\label{sec:rigidity}
\subsection{Contragredient pairings and equivariantization}

\begin{lemma}\label{lem:contragredientpairing}
Let $D(X)=X'$ be the restricted contragredient. There is a canonical weightwise identification $\vartheta_{k,X}:T_kDX\xrightarrow{\sim}D(T_kX)$, under which
\[
 D:\cC_G^{\mathrm{op}}\xrightarrow{\sim}\cC_G,
 \quad D^2\cong\id,\quad
 D(\cC_g)=\cC_{g^{-1}},\quad DT_k\cong T_kD.
\]
For each $X\in\cC_G$ there is a natural morphism $e_X:X\bo DX\to V$, nonzero when $X\ne0$, satisfying
\begin{align}
 e_Y(f\bo\id_{DY})&=e_X(\id_X\bo Df),
                     &&f:X\to Y,\label{eq:pairingnaturality}\\
 e_{T_kX}&=j_k^{-1}\,T_k(e_X)\,J^k_{X,DX}
 (\id_{T_kX}\bo\vartheta_{k,X}^{-1}),
                     &&k\in G,\label{eq:pairingcovariance}
\end{align}
where the second identity has source $T_kX\bo D(T_kX)$ and target $V$; all tensorators are those constructed in Proposition~\ref{prop:action}.
\end{lemma}
\begin{proof}
For a homogeneous $X\in\cC_g$, the restricted dual is
$DX=\bigoplus_{\lambda}X_\lambda^*$. The DLM contragredient
identity reads
\[
 \langle Y_{DX}(v,z)x',x\rangle
 =\big\langle x',Y_X(e^{zL(1)}(-z^{-2})^{L(0)}v,z^{-1})x\big\rangle,
\]
with the inherited logarithm. Its twist is $g^{-1}$. Ordinary,
finite-dimensional homogeneous spaces make $D^2X\to X$ the canonical
weightwise bidual isomorphism; $Df$ is the restricted transpose of
$f$, so $D(f_2f_1)=Df_1Df_2$. For every $k\in G$, $k\omega=\omega$
ensures that $k^{-1}$ commutes with $L(0),L(1)$ on the $V$ input.
Consequently the two displayed contragredient formulas for
$T_kDX$ and $D(T_kX)$ agree coefficientwise:
\[
 \langle Y_{T_kDX}(v,z)x',x\rangle
 =\big\langle x',Y_X(e^{zL(1)}(-z^{-2})^{L(0)}k^{-1}v,z^{-1})x\big\rangle
 =\langle Y_{D(T_kX)}(v,z)x',x\rangle.
\]
Thus $\vartheta_{k,X}$ is the identity on the \emph{underlying
restricted dual spaces}, not an identification of unrelated twisted
module structures. Its naturality and the coherence
$\vartheta_{kl,X}=\vartheta_{k,T_lX}\,T_k(\vartheta_{l,X})$
follow on each weight space from the strict identity
$T_kT_l=T_{kl}$. The construction commutes with finite direct sums,
including multiplicity spaces, and takes the zero object to zero.

The normalized invariant form $(\mathbf1,\mathbf1)=1$ on the simple,
self-contragredient CFT-type VOA $V$ is unique. Since each $k$ fixes
$\mathbf1$ and $\omega$, $(kv,kw)$ is another normalized invariant
form. Therefore $(kv,kw)=(v,w)$, and the identification
$\iota:V\xrightarrow{\sim}V'$ is $G$-equivariant. The unit comparison
$j_k:V\to T_kV$ is $v\mapsto k^{-1}v$, with inverse $j_k^{-1}$
represented by $k$; it is not the identity of twisted $V$-modules.

Let now $0\ne X\in\cC_g$. The actual $g$-twisted vertex operator
$Y_X^g\in\cI\binom X{V\ X}$ has the complete IO properties required
by the skew theorem. Hence the chosen lower-half transform
$\Gamma_X=\Omega_-Y_X^g\in\cI\binom X{X\ V}$ satisfies
\[
 \Gamma_X(x,z)\mathbf1=e^{zL_X(-1)}x.
\]
Indeed this is the vacuum-creation identity, since the other input
is $\mathbf1$ and the logarithm of $z$ is fixed. For
$\mathcal Y\in\cI\binom{W_3}{W_1\ W_2}$ Du--Huang's
contragredient transformation used here is
\[
 \big\langle(A_+\mathcal Y)(w_1,z)w_3',w_2\big\rangle
 =\big\langle w_3',\mathcal Y(
 e^{zL_{W_1}(1)}e^{\pi iL_{W_1}(0)}z^{-2L_{W_1}(0)}w_1,
 z^{-1})w_2\big\rangle.
\]
Theorem~3.3 of \cite{DH} gives its complete IO axioms and its
precise output type $\binom{T_{g_1}W_2'}{W_1\ W_3'}$, where
$W_1\in\cC_{g_1}$. In the present instance $W_1=X$, $W_2=V$,
$W_3=X$, so
\[
 A_X:=A_+(\Gamma_X)\in\cI\binom{T_gV'}{X\ DX},
 \qquad
 \mathcal E_X:=j_g^{-1}\iota^{-1}A_X
 \in\cI\binom V{X\ DX}.
\]
Here $\iota^{-1}$ between $T_gV'$ and $T_gV$ is the corresponding
pullback of the same module map; $j_g^{-1}$ has source $T_gV$.
The product sector is $gg^{-1}=1$, so the output really can be $V$.
The binary universal property supplies a unique
$e_X:X\bo DX\to V$ such that
$\overline{e_X}\mathcal Y^{\mathrm{can}}_{X,DX}=\mathcal E_X$.
No categorical evaluation or duality is used.

Choose a simple direct summand of $X$, a lowest $L(0)$-weight $h$
in that summand, $0\ne x\in X_h$ supported there, and
$x'\in(DX)_h$ with $\langle x',x\rangle\ne0$. Since
$L(1)x=0$ on this lowest-weight vector, the $V'$ vacuum
coefficient, calculated from the preceding \emph{definition} rather
than from a categorical evaluation, is
\[
 \begin{aligned}
 \big\langle A_X(x,z)x',\mathbf1\big\rangle
 &= e^{\pi ih}z^{-2h}
 \big\langle x',\Gamma_X(x,z^{-1})\mathbf1\big\rangle\\
 &=e^{\pi ih}z^{-2h}\langle x',x\rangle\ne0
 \qquad (z>0).
 \end{aligned}
\]
The first bracket uses the canonical pairing $V'\times V$; the
weight-$h$ projection of $e^{z^{-1}L(-1)}x$ is exactly $x$.
Thus $A_X$, $\mathcal E_X$, and, by their universal factorization,
$e_X$ are nonzero, even for nonintegral $h$. This calculation holds
for an arbitrary complex lowest weight on the selected branch.

For $f:X\to Y$ in the same sector, $f\Gamma_X(x,z)v
=\Gamma_Y(fx,z)v$ and $f$ commutes with all Virasoro modes.
The defining $A_+$ matrix coefficient therefore gives, for every
$x\in X$, $y'\in DY$ and $v\in V$,
\[
 \langle A_Y(fx,z)y',v\rangle
 =\langle y',f\Gamma_X(e^{zL(1)}e^{\pi iL(0)}
 z^{-2L(0)}x,z^{-1})v\rangle
 =\langle A_X(x,z)(Df)y',v\rangle.
\]
After the identical $V'\to V$ and $j_g^{-1}$ maps, the binary
universal property proves \eqref{eq:pairingnaturality}. This argument
applies to all morphisms, rather than only scalar morphisms of simple
objects.

For $k\in G$, on the common underlying vector space one has
$\Gamma_{T_kX}(x,z)v=\Gamma_X(x,z)k^{-1}v$.
Because $k$ commutes with all $L(n)$, the $A_+$ definition yields
\[
 \langle A_{T_kX}(x,z)x',v\rangle
 =\langle A_X(x,z)x',k^{-1}v\rangle.
\]
Through the invariant form this equality identifies the output
vector of $A_{T_kX}$ with $k$ times that of $A_X$; before
this identification their targets are different pullbacks of $V'$.
Applying the appropriate $j^{-1}$ maps then gives the
$V$-valued complete IO equality
$\mathcal E_{T_kX}=kg\,A_X$ on the common underlying spaces,
where $kg=(kgk^{-1})k$ is the ordered product.
On the other side, the defining tensorator identity for $J^k$
identifies the pullback of the original canonical binary IO with
the canonical IO on $(T_kX,T_kDX)$. Inserting
$\vartheta_{k,X}^{-1}:D(T_kX)\to T_kDX$, its composite with
$j_k^{-1}T_k(e_X)$ has the same $V$-valued complete IO coefficients
$kg\,A_X$. Uniqueness in the binary universal property proves the
fully typed \eqref{eq:pairingcovariance}. All equalities use the
same $+\pi i$ branch of $A_+$; there is no commutation of $k$ and $g$.

For $X=\bigoplus_{g\in G}X_g$, take
$DX=\bigoplus_gDX_g$, with $DX_g\in\cC_{g^{-1}}$. By finite
biadditivity of $\bo$, define $e_X$ as the direct sum of the
canonical homogeneous pairings on $X_g\bo DX_g$ and as zero on
$X_g\bo DX_h$ for $g\ne h$. These cross summands have degree
$gh^{-1}\ne1$, so any morphism from them to $V\in\cC_1$
vanishes. On each diagonal summand the preceding operator proof
applies, including the individual simple and multiplicity blocks;
the sum is finite, natural, and $G$-covariant because $T_k$
permutes degrees by conjugation. If $X\ne0$, at least one
homogeneous diagonal block has a nonzero pairing.
\end{proof}

\begin{lemma}\label{lem:equivariantization}
The equivariantization $\mathcal B:=\cC_G^G$ is a finite semisimple
$\C$-linear braided monoidal category with simple unit
$\mathbf u=(V,(j_k^{-1})_{k\in G})$. Its forgetful functor
$U:\mathcal B\to\cC_G$ is faithful, exact and strong monoidal.
Every simple $X=(\underline X,u)\in\mathcal B$ satisfies
\[
 \mathbf u\ \mid\ (D\underline X,u^D)\bo X,
 \qquad
 u_k^D:=D((u_k)^{-1})\,\vartheta_{k,\underline X}:
 T_kD\underline X\longrightarrow D\underline X,
\]
where $\mathbf u\mid W$ means that $\mathbf u$ is a direct summand
of $W$. No rigidity assertion is made here.
\end{lemma}
\begin{proof}
An equivariant object is $(X,u)$ with
$u_k:T_kX\xrightarrow{\sim}X$, $u_1=\id_X$, and
$u_{kl}=u_kT_k(u_l)$. The action tensorator has source
$T_kX\bo T_kY$ and target $T_k(X\bo Y)$, so the induced
structure is
\[
 u_k^{X\bo Y}=(u_k^X\bo u_k^Y)(J^k_{X,Y})^{-1},
 \qquad u_k^{\mathbf u}=j_k^{-1}.
\]
To check the cocycle, first apply the inverse of
$J^{kl}_{X,Y}=T_k(J^l_{X,Y})J^k_{T_lX,T_lY}$ and then the
naturalities of $J^k$ and $u_{kl}=u_kT_k(u_l)$. Both resulting
maps from $T_{kl}(X\bo Y)$ to $X\bo Y$ are
$\big(u_k^XT_k(u_l^X)\bo u_k^YT_k(u_l^Y)\big)(J^{kl}_{X,Y})^{-1}$.
The unit cocycle is $j_{kl}^{-1}=j_k^{-1}T_k(j_l^{-1})$.
Proposition~\ref{prop:action} gives precisely the associator and
unit compatibility needed for these maps to be equivariant
morphisms. Thus they define a monoidal category before rigidity.

For $(X,u)$ and $(Y,v)$, let
$X=\bigoplus_gX_g$ be the finite homogeneous decomposition.
The direct sum of the component maps
\[
 \beta_{X_g,Y}:=(v_g\bo\id_{X_g})c_{X_g,Y}:
 X_g\bo Y\longrightarrow Y\bo X_g
\]
defines $\beta_{X,Y}:X\bo Y\to Y\bo X$; an individual $X_g$
need not itself be an equivariant object. The maps are natural
isomorphisms because the crossed commutors and every $v_g$ are.
To check equivariance on the $T_kX_g$ component, put
$a=kgk^{-1}$. Crossed-braiding equivariance of
Proposition~\ref{prop:action} has the precisely typed form
\[
 T_k(c_{X_g,Y})J^k_{X_g,Y}
 =J^k_{T_gY,X_g}\,c_{T_kX_g,T_kY},
\]
and the equivariant cocycle gives the equality of maps
from $T_{kg}Y$ to $Y$,
\[
 v_a\,T_a(v_k)=v_{ak}=v_{kg}=v_k\,T_k(v_g).
\]
Naturality of $c$ and $J$ now identifies
$\beta_{X,Y}u_k^{X\bo Y}$ with
$u_k^{Y\bo X}T_k(\beta_{X,Y})$, component by component;
no identification of $X_g$ as an equivariant object was made.

With associators in the indicated directions, the ordinary hexagon identities are
\[
 \begin{aligned}
 \beta_{X,Y\bo Z}
 &=a^{-1}_{Y,Z,X}(\id_Y\bo\beta_{X,Z})
 a_{Y,X,Z}(\beta_{X,Y}\bo\id_Z)a^{-1}_{X,Y,Z},\\
 \beta_{X\bo Y,Z}
 &=a_{Z,X,Y}(\beta_{X,Z}\bo\id_Y)
 a^{-1}_{X,Z,Y}(\id_X\bo\beta_{Y,Z})a_{X,Y,Z}.
 \end{aligned}
\]
In the first equation for $X_g$, expand the left side as
\[
 \bigl(((v_g^Y\bo v_g^Z)(J^g_{Y,Z})^{-1})\bo\id_{X_g}\bigr)
 c_{X_g,Y\bo Z}
\]
and apply the first crossed hexagon
\eqref{eq:hex1}, then the naturality of $a$ with respect to
$v_g^Y,v_g^Z$. The resulting five maps are exactly the first
ordinary hexagon above. For the second equation for
$X_g,Y_h$, expand its left side as
$\big(v_{gh}^Z\bo\id_{X_g\bo Y_h}\big)c_{X_g\bo Y_h,Z}$.
Apply \eqref{eq:hex2} and use in that order
$v_{gh}^Z=v_g^Z T_g(v_h^Z)$, the naturality of $c$ for
$v_h^Z:T_hZ\to Z$, and the naturality of $a$.
This gives precisely the second displayed ordinary hexagon.
These computations are valid for every homogeneous pair, and finite
biadditivity extends them to arbitrary equivariant objects.
Consequently $\mathcal B$ is braided monoidal; neither proof uses
duals or the conclusion of Proposition~\ref{prop:rigidity}.

The underlying category $\cC_G$ is finite semisimple over $\C$.
Exactness of each pullback $T_k$ makes the kernel and cokernel of
an equivariant map stable under the induced equivariance
isomorphisms. If $p:X\twoheadrightarrow Y$ is equivariant, take an
underlying splitting $s:Y\to X$ and define
\begin{equation}\label{eq:section10average}
 \widetilde s=\frac1{|G|}\sum_{k\in G}
 u_k^X\,T_k(s)\,(u_k^Y)^{-1}.
\end{equation}
Every summand has source $Y$ and target $X$; because $p$ is
equivariant and $ps=\id_Y$, each summand satisfies
$p u_k^XT_k(s)(u_k^Y)^{-1}=\id_Y$. Conjugating the summands
by the equivariance at $l$ and reindexing $k\mapsto lk$ gives
$\widetilde s\,u_l^Y=u_l^XT_l(\widetilde s)$.
Thus $p\widetilde s=\id_Y$ and every equivariant epimorphism
splits; this uses $|G|^{-1}\in\C$ but no duality. Underlying finite
length and splitting show that $\mathcal B$ is semisimple and
idempotent complete.

For $M\in\cC_G$ set
$\operatorname{Ind}(M)=\bigoplus_{k\in G}T_kM$. Its equivariance
on $T_l\operatorname{Ind}(M)=\bigoplus_kT_{lk}M$ sends the
summand indexed by $k$ to the identical summand indexed by $lk$;
left multiplication gives the cocycle since
$T_lT_k=T_{lk}$. For equivariant $(X,u)$, the components
$ p_k=u_k^X:T_kX\to X$ and
$i_k=|G|^{-1}(u_k^X)^{-1}:X\to T_kX$
assemble into equivariant morphisms
$p:\operatorname{Ind}(UX)\to X$ and
$i:X\to\operatorname{Ind}(UX)$. Indeed on the summand labeled
$k$, equivariance of $p$ is precisely
$u_{lk}=u_lT_l(u_k)$, and equivariance of $i$ is its inverse;
moreover $pi=|G|^{-1}\sum_k\id_X=\id_X$.
Writing $S_1,\ldots,S_r$ for the finitely many simple
underlying types and decomposing $UX$ into a finite direct sum of
them shows that every simple of $\mathcal B$ is a direct summand
of one of the finitely many objects $\operatorname{Ind}(S_j)$.
Thus the number of simple isomorphism classes is finite.
Finite-dimensionality of the underlying Hom spaces gives the
same for $\operatorname{Hom}_{\mathcal B}$.
Since $U\mathbf u=V$ is simple, $\mathbf u$ is simple and
$\operatorname{End}_{\mathcal B}(\mathbf u)=\C$.
The forgetful functor is faithful, exact, and strong monoidal:
on the underlying tensor object its structure map is the
identity, and equivariance of that identity is exactly the
already verified formula for $u_k^{X\bo Y}$.

Now define the equivariant contragredient without presupposing it
to be a rigid dual. For $(X,u)$ put
\[
 u_k^{DX}:=D((u_k^X)^{-1})\,\vartheta_{k,X}:
 T_kDX\longrightarrow DX.
\]
Its type is correct since $(u_k^X)^{-1}:X\to T_kX$, so
$D((u_k^X)^{-1}):D(T_kX)\to DX$. To check the cocycle, use
\[
 (u_{kl}^X)^{-1}=T_k((u_l^X)^{-1})(u_k^X)^{-1}
\]
and contravariance, together with the naturality and coherence
of $\vartheta$ proved in Lemma~\ref{lem:contragredientpairing}.
The resulting common map from $T_{kl}DX$ to $DX$ is
\[
 D((u_k^X)^{-1})D(T_k((u_l^X)^{-1}))
 \vartheta_{k,T_lX}\,T_k(\vartheta_{l,X}),
\]
which is both $u_{kl}^{DX}$ and $u_k^{DX}T_k(u_l^{DX})$.
Hence $(DX,u^D)$ is a well-defined object of $\mathcal B$.

The pairing is a morphism in $\mathcal B$. Indeed its equivariance
condition, written with source $T_kX\bo T_kDX$ and target $V$, is
\[
 e_X(u_k^X\bo u_k^{DX})
 =j_k^{-1}T_k(e_X)J^k_{X,DX}.
\]
Apply \eqref{eq:pairingnaturality} with
$f=u_k^X:T_kX\to X$, and precompose its second input with
$D((u_k^X)^{-1})\vartheta_{k,X}$. The composite
$D(u_k^X)D((u_k^X)^{-1})\vartheta_{k,X}$ equals
$\vartheta_{k,X}$, so the left side becomes
$e_{T_kX}(\id_{T_kX}\bo\vartheta_{k,X})$.
Equation~\eqref{eq:pairingcovariance} makes this exactly the
right side. If $X\ne0$, the underlying $e_X$ is nonzero by
Lemma~\ref{lem:contragredientpairing}, and the faithful forgetful
functor detects this.

For a nonzero simple object $X\in\mathcal B$, the braiding
$\beta_{DX,X}:DX\bo X\xrightarrow{\sim}X\bo DX$
is an isomorphism in $\mathcal B$; hence
\[
 \widetilde e_X=e_X\beta_{DX,X}:
 DX\bo X\longrightarrow\mathbf u
\]
is a nonzero morphism in $\mathcal B$. The unit is simple, so
$\widetilde e_X$ is an epimorphism. Semisimplicity gives a section
$i:\mathbf u\to DX\bo X$ with
$\widetilde e_Xi=\id_{\mathbf u}$. Therefore
$\mathbf u\mid DX\bo X$. This is only a unit retract and makes
no claim that $DX$ is a categorical dual or that any triangle
identity holds; those are separate questions for
Proposition~\ref{prop:rigidity}.
\end{proof}

\subsection{Rigidity and the proof of the main theorem}
\begin{proposition}\label{prop:rigidity}
The $G$-graded monoidal category $\cC_G$ is rigid. If $X\in\cC_g$ is homogeneous, its left and right duals can be chosen in $\cC_{g^{-1}}$.
\end{proposition}
\begin{proof}
Write $\mathcal B=\cC_G^G$ and $\mathbf u=(V,(j_k^{-1})_{k\in G})$.
We apply the rigidity criterion of Etingof--Penneys in
$\mathcal B$, then transfer the resulting duals to $\cC_G$.
No duality on $\cC_G$ is used in the first step.

By Lemma~\ref{lem:equivariantization}, $\mathcal B$ is an additive,
idempotent-complete (hence Cauchy complete), $\C$-linear braided
monoidal category. It has finite-dimensional Hom spaces, finitely
many simple isomorphism classes, finite semisimple decompositions,
and $\End_{\mathcal B}(\mathbf u)=\C$. In particular each nonzero
simple object is indecomposable. These properties were established
without invoking rigidity.

Choose representatives $S_1,\dots,S_m$ and a nonzero simple
$S\in\mathcal B$. Set
\[
 (N_S)_{ji}=\dim_\C\Hom_{\mathcal B}(S_j,S\bo S_i)
 \quad(1\le i,j\le m).
\]
Semisimplicity and Schur's lemma identify this integer with the
multiplicity of $S_j$ in $S\bo S_i$; the matrix has finitely many
nonnegative integral entries. By biadditivity, for every $n\ge1$,
\[
 S^{\bo n}\simeq\bigoplus_{i=1}^m S_i^{\oplus a_i(n)},
 \qquad a(n+1)=N_Sa(n),
\]
where $a(n)$ is the column of multiplicities. No positivity,
invertibility or faithfulness of $N_S$ is required. Put
$A=\max\{1,\|N_S\|_1\}$, where
$\|N_S\|_1=\max_i\sum_j(N_S)_{ji}$, and put
$c=\|a(1)\|_1\ge1$. Since $a(n)$ has nonnegative entries,
\[
 \sum_i a_i(n)=\|a(n)\|_1\le cA^{n-1},\qquad
 \dim_\C\End_{\mathcal B}(S^{\bo n})
 =\sum_i a_i(n)^2\le c^2A^{2n-2}.
\]
To verify the strict factorial bound, take $n=2r$, with $r$
an integer larger than $2A^4$ and $\log_2(c^2)$.
Then $(2r)!\ge r^r>(2A^4)^r>c^2A^{4r-2}$;
the last inequality uses $2^r>c^2$ and $A\ge1$.
Thus $\dim\End(S^{\bo n})<n!$ for this $n$, including when
some tensor powers vanish. This is precisely the moderate-growth
condition of \cite[Theorem~1.1]{EtingofPenneys}.

For the same simple $S$, Lemma~\ref{lem:equivariantization}
constructs an actual equivariant object $DS$ and an equivariant
split epimorphism
\[
 \widetilde e_S:DS\bo S\longrightarrow\mathbf u,
 \qquad i_S:\mathbf u\longrightarrow DS\bo S,
 \qquad\widetilde e_Si_S=\id_{\mathbf u}.
\]
In particular $DS$ is an object of $\mathcal B$ (not yet
identified with a categorical dual), and the retract is in the
required order $Y\bo S$, with $Y=DS$. This is precisely
non-negligibility in the sense of Etingof--Penneys; no triangle
identity has entered the argument. Consequently every nonzero
simple $S$ is indecomposable, non-negligible, and of moderate
growth in the Cauchy complete braided category $\mathcal B$,
with finite-dimensional Hom spaces and scalar unit endomorphisms.
By \cite[Theorem~1.1]{EtingofPenneys}, it admits both a left and
a right categorical dual, including evaluation and coevaluation
maps satisfying the respective two triangle identities.

Duals extend to finite direct sums independently of pivotal or spherical structure. If $S_i$ has a left dual $L_i$,
let $X=\bigoplus_iS_i$, $L=\bigoplus_iL_i$. Define
$\operatorname{ev}:L\bo X\to\mathbf u$ to be the evaluation
on each $L_i\bo S_i$ and zero on $L_i\bo S_j$ for $i\ne j$;
define $\operatorname{coev}:\mathbf u\to X\bo L$ as the sum of
the coevaluations followed by the diagonal biproduct inclusions.
Both snake composites have zero off-diagonal components and their
$i$-th diagonal component is the corresponding triangle for $S_i$;
thus they are the identity of $\bigoplus_iS_i$ and
$\bigoplus_iL_i$, respectively. With right evaluations on
$S_i\bo R_i$ and right coevaluations into $R_i\bo S_i$,
the same componentwise calculation gives both right triangles.
The zero object has zero as either dual, its identity being the
zero morphism. Hence $\mathcal B$ is rigid.

For equivariant objects $A,B$, the equivariant tensor product
has underlying object $UA\bo UB$ and equivariance
$u_k^{A\bo B}=(u_k^A\bo u_k^B)(J^k_{UA,UB})^{-1}$.
Thus the tensorator
$\mu_{A,B}:UA\bo UB\xrightarrow{\sim}U(A\bo B)$
of $U$ is the identity on that underlying object; its unit
constraint $\mu_0:V\xrightarrow{\sim}U\mathbf u$ is also the
identity. The monoidal coherence follows from the associator
and unit compatibility already checked in
Lemma~\ref{lem:equivariantization}. In general, for a strong
monoidal functor $F$ with tensorator $\mu$ and unit $\mu_0$,
the images of left evaluation $\epsilon_A:L\bo A\to\mathbf u$
and coevaluation $\eta_A:\mathbf u\to A\bo L$ are
\[
 \mu_0^{-1}F(\epsilon_A)\mu_{L,A}:FL\bo FA\to\mathbf u,
 \qquad
 \mu_{A,L}^{-1}F(\eta_A)\mu_0:\mathbf u\to FA\bo FL.
\]
The first and second triangle diagrams are the images of the
original triangles, using the functor's associator/unit coherence.
For a right dual $R$, the corresponding maps are
$\mu_0^{-1}F(\epsilon_A^{\rm r})\mu_{A,R}$ and
$\mu_{R,A}^{-1}F(\eta_A^{\rm r})\mu_0$; their two triangles
follow from the same coherence. Braidedness, fullness and
essential surjectivity of $F$ are not used. In particular $U$
preserves both kinds of duality.

For an arbitrary, not necessarily homogeneous, $X\in\cC_G$,
use the actual equivariant object
\[
 \widetilde X=\operatorname{Ind}(X)
 =\bigoplus_{k\in G}T_kX\in\mathcal B.
\]
For each $l\in G$, its equivariance
$u_l:T_l\widetilde X=\bigoplus_kT_{lk}X\to\widetilde X$
identifies the summand indexed by $k$ on the source with the
summand indexed by $lk$ on the target. Left multiplication
$k\mapsto lk$ gives $u_lT_l(u_t)=u_{lt}$, by
$T_lT_t=T_{lt}$, so these are valid equivariant structure maps.
The direct sum is finite; put
$P=U\widetilde X=\bigoplus_{k\in G}T_kX$.
Its $k=1$ summand is $T_1X=X$, giving actual inclusion and
projection maps
\[
 i:X\longrightarrow P,\qquad r:P\longrightarrow X,
 \qquad ri=\id_X,
 \qquad p=ir\in\End_{\cC_G}(P),\quad p^2=p.
\]
They are maps in $\cC_G$ also for inhomogeneous $X$.
Because $\widetilde X$ is rigid in $\mathcal B$ and $U$
is strong monoidal, $P$ has both duals in $\cC_G$.
We construct both duals of its retract $X$, without
appealing to a duality functor on the whole of $\cC_G$.

Let $L=P^\vee$ be a chosen left dual of $P$, with
$\epsilon:L\bo P\to V$ and $\eta:V\to P\bo L$.
All associators $a$ and unitors $l,r$ below are those of
$\cC_G$. The two triangles for $P$ and $L$ are
\begin{align}
 r_P(\id_P\bo\epsilon)a_{P,L,P}
       (\eta\bo\id_P)l_P^{-1}&=\id_P,\label{eq:rigidity-leftP}\\
 l_L(\epsilon\bo\id_L)a^{-1}_{L,P,L}
       (\id_L\bo\eta)r_L^{-1}&=\id_L.
       \label{eq:rigidity-leftL}
\end{align}
Define, for $h:L\bo P\to V$, the map
\[
 \Phi^{-1}(h):=
 l_L(h\bo\id_L)a^{-1}_{L,P,L}
       (\id_L\bo\eta)r_L^{-1}:L\longrightarrow L,
 \qquad \Phi(t):=\epsilon(t\bo\id_P).
\]
The first and second displayed triangles, together with the
naturality of $a,l,r$, give respectively
$\Phi\Phi^{-1}(h)=h$ and $\Phi^{-1}\Phi(t)=t$.
For example, in the coherent strict notation the latter reads
$(\epsilon\bo\id_L)(t\bo\id_P\bo\id_L)
(\id_L\bo\eta)=t$, by the triangle for $L$;
the former is
$\epsilon(\Phi^{-1}(h)\bo\id_P)=h$
by the triangle for $P$. Thus $\Phi$ is a bijection, proved
using only the chosen dual of $P$. Define the left transpose
\[
 p^\vee:=\Phi^{-1}\bigl(\epsilon(\id_L\bo p)\bigr):L\to L.
\]
Its defining equation is the first of
\begin{equation}\label{eq:section10transpose}
 \epsilon(p^\vee\bo\id_P)=\epsilon(\id_L\bo p),
 \qquad
 (p\bo\id_L)\eta=(\id_P\bo p^\vee)\eta.
\end{equation}
For the second equation, the map
$\Theta(f):=r_P(\id_P\bo\epsilon)a_{P,L,P}
(f\bo\id_P)l_P^{-1}$ from
$\Hom(V,P\bo L)$ to $\End(P)$ is bijective,
with inverse $t\mapsto(t\bo\id_L)\eta$, by the same
triangles and coherence. Applying $\Theta$ to both sides of
the claimed equation gives $p$: for the first side use
\eqref{eq:rigidity-leftP}, and for the second side use the
defining equation for $p^\vee$ and
\eqref{eq:rigidity-leftP}. Injectivity of $\Theta$ proves the
second transpose identity. Also
$\Phi((p^\vee)^2)=\epsilon(\id_L\bo p^2)
=\Phi(p^\vee)$, hence $(p^\vee)^2=p^\vee$.

Split this idempotent in the finite semisimple category $\cC_G$:
\[
 p^\vee=js,\qquad j:X^\vee\to L,\quad
 s:L\to X^\vee,\quad sj=\id_{X^\vee}.
\]
Then $p^\vee j=j$ and $sp^\vee=s$. The maps
\[
 \epsilon_X=\epsilon(j\bo i):X^\vee\bo X\to V,
 \qquad
 \eta_X=(r\bo s)\eta:V\to X\bo X^\vee
\]
are both well typed. Their complete, non-strict triangle
identities to be verified are
\begin{align}
 r_X(\id_X\bo\epsilon_X)a_{X,X^\vee,X}
       (\eta_X\bo\id_X)l_X^{-1}&=\id_X,
       \label{eq:rigidity-leftX}\\
 l_{X^\vee}(\epsilon_X\bo\id_{X^\vee})
       a^{-1}_{X^\vee,X,X^\vee}
       (\id_{X^\vee}\bo\eta_X)r_{X^\vee}^{-1}
       &=\id_{X^\vee}.
       \label{eq:rigidity-leftXdual}
\end{align}
The calculation uses the parenthesizations determined by $a,l,r$ in
\eqref{eq:rigidity-leftX}--\eqref{eq:rigidity-leftXdual}.
Functoriality and naturality of these constraints give, for the
first left-hand side,
\begin{align*}
 &(\id_X\bo\epsilon_X)(\eta_X\bo\id_X)\\
 &\quad=r(\id_P\bo\epsilon)
       (\id_P\bo js\bo\id_P)(\eta\bo i)\\
 &\quad=r(\id_P\bo\epsilon)
       (\id_P\bo\id_L\bo p)(\eta\bo i)
       &&\text{by \eqref{eq:section10transpose}}\\
 &\quad=r(\id_P\bo\epsilon)(\eta\bo pi)
       =ri=\id_X
       &&\text{by \eqref{eq:rigidity-leftP} and $pi=i$.}
\end{align*}
For the second left-hand side one obtains
\begin{align*}
 &(\epsilon_X\bo\id_{X^\vee})
       (\id_{X^\vee}\bo\eta_X)\\
 &\quad=s(\epsilon\bo\id_L)
       (j\bo ir\bo\id_L)(\id_{X^\vee}\bo\eta)\\
 &\quad=s(\epsilon\bo\id_L)
       (p^\vee j\bo\id_P\bo\id_L)
       (\id_{X^\vee}\bo\eta)
       &&\text{by \eqref{eq:section10transpose}}\\
 &\quad=s(\epsilon\bo\id_L)(j\bo\eta)
       =sj=\id_{X^\vee}
       &&\text{by \eqref{eq:rigidity-leftL} and $p^\vee j=j$.}
\end{align*}
Thus \eqref{eq:rigidity-leftX} and
\eqref{eq:rigidity-leftXdual} hold, with the actual non-strict
associators and unitors, proving that $X^\vee$ is a left dual.

Independently choose a right dual $R={}^\vee P$, with
$\epsilon^{\rm r}:P\bo R\to V$ and
$\eta^{\rm r}:V\to R\bo P$. Its right triangle identities,
with the opposite parenthesizations to the left case, are
\begin{align}
 l_P(\epsilon^{\rm r}\bo\id_P)a^{-1}_{P,R,P}
       (\id_P\bo\eta^{\rm r})r_P^{-1}&=\id_P,
       \label{eq:rigidity-rightP}\\
 r_R(\id_R\bo\epsilon^{\rm r})a_{R,P,R}
       (\eta^{\rm r}\bo\id_R)l_R^{-1}&=\id_R.
       \label{eq:rigidity-rightR}
\end{align}
The map
$\Psi(q)=\epsilon^{\rm r}(\id_P\bo q)$ from
$\End(R)$ to $\Hom(P\bo R,V)$ has inverse
\[
 \Psi^{-1}(h)=r_R(\id_R\bo h)a_{R,P,R}
       (\eta^{\rm r}\bo\id_R)l_R^{-1},
 \qquad h:P\bo R\longrightarrow V.
\]
The two right triangles show $\Psi\Psi^{-1}=\id$ and
$\Psi^{-1}\Psi=\id$, by the right-hand analogues of the two
explicit snake computations above. Set
$\prescript{\vee}{}{p}=\Psi^{-1}
(\epsilon^{\rm r}(p\bo\id_R))$. Then
\begin{equation}\label{eq:rigidity-righttranspose}
 \epsilon^{\rm r}(\id_P\bo\prescript{\vee}{}{p})
       =\epsilon^{\rm r}(p\bo\id_R),
 \qquad
 (\prescript{\vee}{}{p}\bo\id_P)\eta^{\rm r}
       =(\id_R\bo p)\eta^{\rm r}.
\end{equation}
Indeed the latter identity follows by applying the bijection
\[
 f\longmapsto l_P(\epsilon^{\rm r}\bo\id_P)
 a^{-1}_{P,R,P}(\id_P\bo f)r_P^{-1}
 \quad\bigl(f:V\to R\bo P\bigr),
\]
whose inverse is $t\mapsto(\id_R\bo t)\eta^{\rm r}$;
both sides have image $p$, by
\eqref{eq:rigidity-rightP} and the defining transpose identity.
Uniqueness gives
$(\prescript{\vee}{}{p})^2=\prescript{\vee}{}{p}$,
since $p^2=p$. Choose a splitting
\[
 \prescript{\vee}{}{p}=j_{\rm r}s_{\rm r},\qquad
 j_{\rm r}:{}^\vee X\to R,\quad
 s_{\rm r}:R\to{}^\vee X,\quad
 s_{\rm r}j_{\rm r}=\id_{{}^\vee X}.
\]
It follows that $\prescript{\vee}{}{p}j_{\rm r}=j_{\rm r}$.
Put
\[
 \epsilon_X^{\rm r}=\epsilon^{\rm r}(i\bo j_{\rm r}):
 X\bo{}^\vee X\to V,
 \qquad
 \eta_X^{\rm r}=(s_{\rm r}\bo r)\eta^{\rm r}:
 V\to{}^\vee X\bo X.
\]
The two correctly typed right triangles, retaining the
associators and unitors, are
\begin{align}
 l_X(\epsilon_X^{\rm r}\bo\id_X)
 a^{-1}_{X,{}^\vee X,X}
 (\id_X\bo\eta_X^{\rm r})r_X^{-1}
 &=\id_X,\label{eq:rigidity-rightX}\\
 r_{{}^\vee X}(\id_{{}^\vee X}\bo\epsilon_X^{\rm r})
 a_{{}^\vee X,X,{}^\vee X}
 (\eta_X^{\rm r}\bo\id_{{}^\vee X})l_{{}^\vee X}^{-1}
 &=\id_{{}^\vee X}.\label{eq:rigidity-rightXdual}
\end{align}
Their component calculations, with the parenthesizations of
\eqref{eq:rigidity-rightX}--\eqref{eq:rigidity-rightXdual}
restored by coherence, are independent of the left calculations:
\begin{align*}
 &(\epsilon_X^{\rm r}\bo\id_X)
    (\id_X\bo\eta_X^{\rm r})\\
 &\quad=r(\epsilon^{\rm r}\bo\id_P)
 (i\bo j_{\rm r}s_{\rm r}\bo\id_P)
 (\id_X\bo\eta^{\rm r})\\
 &\quad=r(\epsilon^{\rm r}\bo\id_P)
 (pi\bo\id_R\bo\id_P)
 (\id_X\bo\eta^{\rm r})
 &&\text{by \eqref{eq:rigidity-righttranspose}}\\
 &\quad=r(\epsilon^{\rm r}\bo\id_P)
 (i\bo\eta^{\rm r})=ri=\id_X
 &&\text{by \eqref{eq:rigidity-rightP};}
\end{align*}
\begin{align*}
 &(\id_{{}^\vee X}\bo\epsilon_X^{\rm r})
    (\eta_X^{\rm r}\bo\id_{{}^\vee X})\\
 &\quad=s_{\rm r}(\id_R\bo\epsilon^{\rm r})
 (\id_R\bo ir\bo j_{\rm r})
 (\eta^{\rm r}\bo\id_{{}^\vee X})\\
 &\quad=s_{\rm r}(\id_R\bo\epsilon^{\rm r})
 (\prescript{\vee}{}{p}\bo\id_P\bo j_{\rm r})
 (\eta^{\rm r}\bo\id_{{}^\vee X})
 &&\text{by \eqref{eq:rigidity-righttranspose}}\\
 &\quad=s_{\rm r}\prescript{\vee}{}{p}
 (\id_R\bo\epsilon^{\rm r})
 (\eta^{\rm r}\bo j_{\rm r})
 =s_{\rm r}\prescript{\vee}{}{p}j_{\rm r}
 =s_{\rm r}j_{\rm r}=\id_{{}^\vee X}
 &&\text{by \eqref{eq:rigidity-rightR}.}
\end{align*}
This verifies \eqref{eq:rigidity-rightX} and
\eqref{eq:rigidity-rightXdual}. Therefore every object
of $\cC_G$ has both left and right duals.

Let $X\in\cC_g$ and decompose the constructed left dual
$X^\vee=\bigoplus_{h\in G}L_h$ according to the external
$G$-grading. Because $L_h\bo X\in\cC_{hg}$ and
$X\bo L_h\in\cC_{gh}$ whereas $V\in\cC_1$,
the evaluation $\epsilon_X$ vanishes on $L_h\bo X$ for
$h\ne g^{-1}$, and the coevaluation $\eta_X$ factors through
$X\bo L_{g^{-1}}$. The homogeneous projections and inclusions
are morphisms of $\cC_G$, and the associator preserves each
sector because its source and target have the same ordered product
of degrees. Hence restricting both maps to $L_{g^{-1}}$
retains the $X$-triangle, as all its nonzero terms already lie
in that sector. The second triangle restricted to $L_{g^{-1}}$
remains its identity, by naturality of the associator and
unitors with the corresponding inclusion and projection.
More strongly, on each $L_h$ with $h\ne g^{-1}$, the
$L$-triangle has a zero evaluation component, hence is the
zero endomorphism of $L_h$; since the full triangle is
$\id_{X^\vee}$, this forces $L_h=0$. Thus even the chosen
left dual has no non-inverse component.

For the right dual write ${}^\vee X=\bigoplus_hR_h$.
Now $X\bo R_h\in\cC_{gh}$ and
$R_h\bo X\in\cC_{hg}$. The right evaluation is zero on
$X\bo R_h$ for $h\ne g^{-1}$, and the right coevaluation
factors through $R_{g^{-1}}\bo X$. The right $X$-triangle
and right $R_{g^{-1}}$-triangle persist by the same
sector-preserving naturality; the second right triangle
vanishes on any $R_h$, $h\ne g^{-1}$, and forces it to be zero.
Consequently $X^\vee,{}^\vee X\in\cC_{g^{-1}}$.
For $X=0$, choosing both duals to be zero is compatible with
all triangles and with every sector. This completes the proof.
\end{proof}

\begin{proof}[Proof of Theorem~\ref{thm:main}]
By Proposition~\ref{prop:regularmodules}, every object of each
$\cC_g$ is a finite direct sum of ordinary irreducibles, there are
finitely many simple types, and $\cC_g\ne0$. Kernels and cokernels
of module maps, and the images of their idempotents, are again such
finite direct sums. Lemma~\ref{lem:Schur} therefore gives finite-dimensional
Hom spaces and scalar endomorphisms of simples. Since $G$ is finite,
the external direct sum $\cC_G$ is a finite semisimple, idempotent-complete
$\C$-linear abelian category. Its cross-sector Hom spaces vanish by
definition, and its homogeneous projections are morphisms.

Theorem~\ref{thm:binary} represents the entire operator space of
Definition~\ref{def:completeintertwiner}, naturally in both inputs and
every output in $\cC_{gh}$. With $P(1)$ and $\log1=0$ this gives
$\Hom_V(X\bo Y,Z)\cong\cI\binom Z{X\ Y}$ and
$X\bo Y\in\cC_{gh}$. Additivity of the four complete-operator
axioms gives a $\C$-bilinear, biadditive tensor bifunctor on the
external direct sum. Every exact sequence splits, so this functor
is exact in each variable. The unitors of Lemma~\ref{lem:units}
identify its unit with the simple adjoint module $V\in\cC_1$.
Every sector is nonzero, so the external $G$-grading is faithful;
none of these arguments commutes two group elements.

Theorem~\ref{thm:partialnodeweakmodule} constructs the weak node module
before Proposition~\ref{prop:partialnodefiniteordinary} uses twisted
regularity. Proposition~\ref{prop:partialnodevertices} and
Theorem~\ref{thm:matchedcoefficientmodes} give the complete node
operators and the all-mode kernel relation. Together with the original
convergence of Corollary~\ref{cor:smallconvergenceautomatic},
Theorem~\ref{thm:partialnodeexactspan} gives actual product--iterate
factorization. The natural morphisms $b:L\to R$, $c:R\to L$ of
Lemmas~\ref{lem:productepi} and \ref{lem:finiterotation} satisfy
$bc=\id_R$, $cb=\id_L$. Theorem~\ref{thm:smallassoc} identifies
$a=b$ with the canonical positive-real $P(z)$ normalization and
retains the original product--iterate identity. This construction
precedes and does not use the pentagon.

Proposition~\ref{prop:unittriangle} gives the triangle and $l_V=r_V$.
Theorem~\ref{thm:criticalpentagon}, using
Propositions~\ref{prop:criticalfiveconvergence} and
\ref{prop:criticalconstantmaps}, gives the pentagon for this same
associator. Proposition~\ref{prop:action} gives the strong monoidal
pullbacks $T_k$, their natural tensorators $J^k$, and unit maps $j_k$,
including their associator and unitor compatibility and
\eqref{eq:Jcomp}. Thus $T_kT_l=T_{kl}$ and $T_1=\id$ are strict
identities of underlying functors, with identity action compositors;
the tensorators themselves need not be identities. Pullback sends
$\cC_g$ to $\cC_{kgk^{-1}}$.

The crossed commutors are natural isomorphisms by
Lemma~\ref{lem:units} and Proposition~\ref{prop:action}. Their
upper-half-plane $\Omega_+$ normalization is
\eqref{eq:criticalrawc}, proved in
Proposition~\ref{prop:criticalconstantmaps}. Proposition~\ref{prop:action}
proves the equivariance identity \eqref{eq:equiv}, and
Theorem~\ref{thm:criticalhexagons} proves both
\eqref{eq:hex1} and \eqref{eq:hex2}, with their displayed ordered
sectors. The corresponding unit comparisons, for $X\in\cC_g$, are
\[
 r_Xc_{V,X}=l_X,\qquad
 l_X(j_g^{-1}\bo\id_X)c_{X,V}=r_X.
\]
Indeed the first is the complete-operator equality
$\Omega_+\Gamma_X=Y_X$. For the second, after identifying $T_gV$
with $V$ by $j_g^{-1}$, the represented operator is
$e^{zL(-1)}Y_X(gv,e^{\pi i}z)x$; twisted equivariance identifies
this with $e^{zL(-1)}Y_X(v,e^{-\pi i}z)x=\Gamma_X(x,z)v$.
The binary universal property gives both displayed morphism identities.
Hence the monoidal action, equivariance and two crossed hexagons are
exactly the data and diagrams in
\cite[Definition~8.24.1, (8.104)--(8.107)]{EGNO}, specialized to
identity action compositors.

Finally Lemma~\ref{lem:contragredientpairing} constructs a nonzero
pairing without presupposing a categorical dual.
Lemma~\ref{lem:equivariantization} constructs a finite semisimple
braided monoidal equivariantization, its strong monoidal forgetful
functor and the required unit retracts before rigidity is invoked.
Proposition~\ref{prop:rigidity} applies the Etingof--Penneys criterion
there and transfers both duals to $\cC_G$, including all four
triangle identities and
$X_g^\vee,{}^\vee X_g\in\cC_{g^{-1}}$. Thus $\cC_G$ is a fusion
category and the preceding data give the asserted faithfully
$G$-graded, $G$-crossed braided fusion structure.
\end{proof}

\subsection*{Use of AI tools}
ChatGPT (OpenAI) was used in mathematical exploration, literature searching, preparation of preliminary drafts, proof checking, language editing, and \LaTeX{} typesetting. The author assumes full responsibility for the mathematical arguments, references, and final text.

\end{document}